\documentclass[11pt]{mybooke}
\usepackage{amssymb,amscd,ascmac}
\usepackage{graphicx}
\usepackage{eurm}
\usepackage[mathscr]{eucal}
\usepackage{longtable}
\usepackage{minitoc}
\renewcommand{\mtcfont}{\small\rm}
\newenvironment{abst}{%
\vspace{2ex}\begin{center}\begin{minipage}{12.5cm}%
\setlength{\baselineskip}{.85\baselineskip}%
\begin{small}{\sc Abstract}:\ }%
{\end{small}\end{minipage}\end{center}}

\newcommand{\Sh}{\mathcal{S}h}

\newcommand{\getsfrom}{\ensuremath{\gets\kern-.45em\lower-.1ex\hbox%
{$\shortmid\,$}}}

\newcommand{\ovl}[1]{\overline{#1}}
\newcommand{\udl}[1]{\underline{#1}}
\newcommand{\til}[1]{\widetilde{#1}}

\newcommand{\FRAC}[2]{\leavevmode\kern-.1em%
\raise.5ex\hbox{\the\scriptfont0 #1}\kern-.1em%
/\kern-.15em\lower.25ex\hbox{\the\scriptfont0 #2}}

\newcommand{\bp}{{\lower -0.4ex \hbox{\bf{.}}}}
\newcommand{\bbp}{{\lower -0.4ex \hbox{\huge\bf{.}}}}
\newcommand{\bD}{\boldsymbol{D}}
\newcommand{\bcD}{\boldsymbol{\mathcal{D}}}
\newcommand{\D}{\mathbf{D}}

\newcommand{\Aut}{\mathop{\mathrm{Aut}}\nolimits}

\newcommand{\ch}{\mathop{\mathrm{char}}\nolimits}

\newcommand{\Coim}{\mathop{\mathrm{Coim}}\nolimits}

\newcommand{\CoLie}{\mathop{\mathrm{CoLie}}\nolimits}
\newcommand{\Cov}{\mathop{\mathrm{Cov}}\nolimits}

\newcommand{\Def}{\mathop{\mathrm{Def}}\nolimits}

\newcommand{\dlog}{\mathop{\mathrm{dlog}}\nolimits}
\newcommand{\End}{\mathop{\mathrm{End}}\nolimits}

\newcommand{\Ind}{\mathop{\mathrm{Ind}}\nolimits}

\newcommand{\Gal}{\mathop{\mathrm{Gal}}\nolimits}

\newcommand{\Gr}{\mathop{\mathrm{Gr}}\nolimits}
\newcommand{\Grass}{\mathop{\mathrm{Grass}}\nolimits}

\newcommand{\Hom}{\mathop{\mathrm{Hom}}\nolimits}

\newcommand{\id}{\mathop{\mathrm{id}}\nolimits}
\renewcommand{\Im}{\mathop{\mathrm{Im}}\nolimits}
\newcommand{\inv}{\mathop{\mathrm{inv}}\nolimits}
\newcommand{\Isom}{\mathop{\mathrm{Isom}}\nolimits}

\newcommand{\Ker}{\mathop{\mathrm{Ker}}\nolimits}

\newcommand{\Lie}{\mathop{\mathrm{Lie}}\nolimits}
\newcommand{\LOG}{\mathop{\mathrm{LOG}}\nolimits}
\newcommand{\Log}{\mathop{\mathrm{Log}}\nolimits}

\newcommand{\Nil}{\mathop{\mathrm{Nil}}\nolimits}

\newcommand{\Nr}{\mathop{\mathrm{Nr}}\nolimits}

\newcommand{\ord}{\mathop{\mathrm{ord}}\nolimits}

\newcommand{\qisog}{\mathop{\mathrm{qisog}}\nolimits}
\newcommand{\rank}{\mathop{\mathrm{rank}}\nolimits}
\newcommand{\rk}{\mathop{\mathrm{rk}}\nolimits}

\newcommand{\Res}{\mathop{\mathrm{Res}}\nolimits}

\newcommand{\Sets}{\mathop{\mathrm{Sets}}\nolimits}
\renewcommand{\sp}{\mathop{\mathrm{sp}}\nolimits}
\newcommand{\Spec}{\mathop{\mathrm{Spec}}\nolimits}

\newcommand{\Spf}{\mathop{\mathrm{Spf}}\nolimits}
\newcommand{\Spm}{\mathop{\mathrm{Spm}}\nolimits}
\newcommand{\Stab}{\mathop{\mathrm{Stab}}\nolimits}
\newcommand{\Supp}{\mathop{\mathrm{Supp}}\nolimits}

\newcommand{\Tr}{\mathop{\mathrm{Tr}}\nolimits}
\newcommand{\tr}{\mathop{\mathrm{tr}}\nolimits}

\def\dfrac{\displaystyle\frac}
\def\A{\mathbb{A}}
\def\C{\mathbb{C}}
\def\E{\mathbb{E}}
\def\F{\mathbb{F}}
\def\G{\mathbb{G}}

\def\N{\mathbb{N}}
\def\P{\mathbb{P}}
\def\Q{\mathbb{Q}}

\def\R{\mathbb{R}}
\def\Z{\mathbb{Z}}

\def\AA{\mathcal{A}}

\def\FF{\mathcal{F}}

\def\MM{\mathcal{M}}
\def\OO{\mathcal{O}}
\def\PP{\mathcal{P}}

\def\VV{\mathcal{V}}

\def\q{\mathfrak{q}}
\newcommand{\liminj}{\varinjlim}
\newcommand{\limproj}{\varprojlim}
\newcommand{\an}[1]{\langle#1\rangle}
\newcommand{\lan}[1]{\bigl\langle#1\bigr\rangle}
\newcommand{\Lan}[1]{\Bigl\langle#1\Bigr\rangle}

\newcommand{\ie}{{\it i.e. }}
\newcommand{\cf}{{\it cf. }}
\newcommand{\resp}{{\it resp. }}
\newcommand{\eg}{{\it e.g. }}
\newcommand{\loccit}{{\it loc.~cit. }}

\numberwithin{equation}{section}
\theoremstyle{plain}
\newtheorem{thm}{\bf Theorem}[subsection]
\newtheorem{pro}[thm]{Proposition}
\newtheorem{lem}[thm]{Lemma}
\newtheorem{klem}[thm]{Key lemma}
\newtheorem{cor}[thm]{Corollary}
\newtheorem{coj}[thm]{Conjecture}

\newtheorem{que}[thm]{Question}

\theoremstyle{definition}
\newtheorem{dfn}[thm]{Definition}

\newtheorem{para}[thm]{\bf\hspace{-.34em}}

\theoremstyle{remark}
\newtheorem{exa}[thm]{Example}
\newtheorem{exs}[thm]{Examples}
\newtheorem{rem}[thm]{Remark}
\newtheorem{rems}[thm]{Remarks}

\newcommand{\thmref}[1]{Theorem~\ref{#1}}

\newcommand{\defref}[1]{Definition~\ref{#1}}
\newcommand{\chapref}[1]{Chapter~\ref{#1}}

\newcommand{\secref}[1]{\S~\ref{#1}}
\newcommand{\subsecref}[1]{\S\S\ref{#1}}
\newcommand{\lemref}[1]{Lemma~\ref{#1}}

\newcommand{\examref}[1]{\ref{#1}}
\newcommand{\propref}[1]{Proposition~\ref{#1}}

\newtheorem{appthm}{\bf Theorem}[section]
\newtheorem{apppro}[appthm]{Proposition}
\newtheorem{applem}[appthm]{Lemma}
\newtheorem{appcor}[appthm]{Corollary}

\theoremstyle{definition}
\newtheorem{appdfn}[appthm]{Definition}

\newtheorem{appexa}[appthm]{Example}

\newtheorem{apppara}[appthm]{\bf\hspace{-.34em}}

\theoremstyle{remark}
\newtheorem{apprem}[appthm]{Remark}

\renewcommand{\labelenumi}{(\theenumi)}
\renewcommand{\theenumi}{\roman{enumi}}
\title{Period mappings and differential equations. From $\C$ to $\C_p$ \\
\hspace{0em}{\mdseries\sffamily\scshape\large
T\^{o}hoku-Hokkaid\^{o} lectures in Arithmetic Geometry}}
\author{Yves Andr\'e}

\makeindex
\begin{document}
\pagenumbering{roman}

\dominitoc
\nomtcrule
\cleardoublepage \thispagestyle{empty}%
\begin{center}
 \setlength{\baselineskip}{2.0\baselineskip}
 {\LARGE\bfseries Period mappings and differential equations.
 From $\C$ to $\C_p$}\\
\hspace{0em}{\mdseries\sffamily\scshape\large
T\^{o}hoku-Hokkaid\^{o} lectures in Arithmetic Geometry}
\end{center}
\bigskip
\vspace{24pt}
\begin{center}
 {\LARGE Yves Andr\'e}
\end{center}
\vspace{50pt}
\begin{center}
 \large (with appendices by F. Kato and N. Tsuzuki)
\end{center}
\vfill
\newpage\thispagestyle{empty}

\tableofcontents

\pagenumbering{arabic}

\markboth{INTRODUCTION}{}

\begin{center}
 {\large\bf Introduction.}
\end{center}

\addcontentsline{toc}{subsection}{Introduction.}
\subsection*{Purpose of this book.}

This book stems from lectures given at the Universities of T\={o}hoku
and Hokkaid\={o}.  

Its main purpose is to introduce the reader to $p$-adic analytic
geometry and to the theory of $p$-adic analytic functions and
differential equations, by focusing on the theme of {\it period
mappings}. Of course, this approach is not meant to replace more
systematic expositions of $p$-adic analysis or geometry found in a
number of good treatises, but to complement them.

 It is our general rule, in this book, to follow as closely as possible
the complex theory, and to go back and forth between the complex and the
$p$-adic worlds. We hope that this will make the text of interest both
to some complex geometers and to some arithmetic geometers.

In the course of Chapter I, this approach will eventually become a
kaleidoscope of half-correspondences and broken echoes. We hope that the
reader will then have gained enough hindsight and wariness about these
analogies, and will enjoy seeing unity being restored at a deeper level
in Chapters II and III.

\bigskip

We have chosen the theme of period mappings because of its
central role in the nineteen-century mathematics as a fertile place of
interaction between differential equations, group theory, algebraic and
differential geometry, topology (and even number theory).  In fact, it
was a guiding thread in the early harmonious development of these
branches of mathematics, from Gauss and Riemann to Klein and Poincar\'e
(\cf [{\bf Gray86}]).

\bigskip

The origins lie in Gauss' largely unpublished work on elliptic
functions on one hand (rediscovered and extended by Abel and by Jacobi),
and on the hypergeometric differential equation in the complex domain on
the other hand. Gauss knew the connection between the two topics: the
inverse of the indefinite integral
$\int\frac{dx}{\sqrt{(1-x^2)(1-k^2x^2)}}$ is a single-valued function
with two independent periods $\omega_1$ and $\omega_2$ which are
solutions of the 
hypergeometric differential equation with parameters
$(\frac{1}{2},\frac{1}{2};1)$ in the variable $k^2$.

Riemann's point of view of the `Riemann surface' of a multivalued
complex function has given a geometric framework for all of complex
analysis. He applied this idea with equal success to Jacobi's inversion
problem for more general indefinite algebraic integrals on one hand; and
to the elucidation of the paradoxical polymorphism of hypergeometric
functions which had puzzled Gauss and Kummer by introducing the concept
of monodromy, on the other hand. He studied in detail the monodromy of
the multivalued map $ \; k^2\mapsto \tau ={\omega_2/\omega_1}\;$, the
first and basic example of a `period mapping' , whose inverse is
single-valued. He also showed (and Schwarz rediscovered) that under
certain conditions, the quotient of two solutions of more general
hypergeometric equations maps the upper half-plane onto a curvilinear
triangle.

 Group-theoretic aspects of monodromy were studied by Jordan, and the
 geometrization of complex analysis around the concept of automorphic
 function was carried on by Klein, who once described his work as
 `blending Galois with Riemann'\footnote{quoted from [{\bf Gray86},
 p. 179].}.  The achievement fell to Poincar\'e. Starting from Fuchs'
 problem of finding all second-order differential equations for which a
 quotient of solutions $\tau$ admits a single-valued inverse, Poincar\'e
 founded the theory of `uniformizing differential equations' (in modern
 terminology), and eventually recognized that every Riemann surface of
 genus $>1$ should arise from the action of a discrete group of
 non-euclidean moves on the upper half-plane.

\bigskip
\bigskip

 What about non-archimedean analogues of this saga?

\bigskip

It is clear that the development of $p$-adic analysis and geometry took
a strikingly different route.

The very beginnings looked similar, indeed: J. Tate introduced rigid
geometry as a proper framework for the uniformization $\C_p^\times/q^\Z$
of elliptic curves with multiplicative reduction; B. Dwork developed
$p$-adic analysis starting from the hypergeometric differential equation
with parameters $(\frac{1}{2},\frac{1}{2};1)$, where he discovered the
essential notion of `Frobenius structure'. But it was immediately clear
that since solutions of $p$-adic differential equations do not converge
up to the next singularity, no faithful counterpart of complex monodromy
could take place.

From there on, different $p$-adic theories grew apart, with their own
languages, most of them relying on sophisticated parts of contemporary
algebraic geometry, and all claiming some analogy with ``the complex
case'':

\begin{itemize}
 \renewcommand{\labelitemi}{\normalfont\bfseries\textendash}
 \item the theory of differential equations matured slowly (with a strong
       orientation toward applications to exponential sums), struggling with
       the problems of singularities, without being able to tackle global
       problems,
 \item crystalline theory offered a global viewpoint on differential
       equations (oriented toward the cohomology of varieties in positive
       characteristic), but did not help to understand singularities,
 \item the theory of $p$-adic representations and $p$-adic Hodge theory
       developed independently of differential equations, as did the several
       avatars of rigid geometry, motivated by idiosyncratic problems.
\end{itemize}

 \bigskip

However, it is the author's opinion that the situation has somewhat
changed over the last years, that isolated branches are merging by fits
and starts\footnote{let us mention notably the incursion of rigid
geometry into the crystalline viewpoint and into the geometric theory of
finite coverings, the maturity of index theory and its applications to
algebraic geometry in positive characteristic, the new connections
between $p$-adic representations and differential equations on
annuli...}, and that a synthesis is gradually emerging. This encourages
to hope that, after many twists and turns, period mappings can indeed
become a unifying theme in the $p$-adic context. This book is intended
to be a contribution in this direction, by bringing closer the $p$-adic
theory to its complex precursor --- from periods and monodromy up to
triangle groups. Much remains to be done in order to achieve comparable
harmony and clarity.

\bigskip

\subsection*{Contents of this book.}

{\it Chapter I} is preparatory. It deals with the problems of analytic
continuation --- with emphasis on the case of solutions of differential
equations --- and periods of abelian integrals, in the $p$-adic context.

\medskip Multivalued complex-analytic functions can be handled in two
different, but essentially 
equivalent, ways:

$1)$ in a geometric way using Riemann surfaces, coverings and paths; 

$2)$ as limits of algebraic functions; this less orthodox way leads to a more
Galois-theoretic viewpoint on analytic continuation.

\medskip  Both ways are practicable in the $p$-adic context, but
eventually lead to completely different theories of analytic continuation. 

Following the first way demands to have at disposal $p$-adic spaces
which are locally arcwise connected.  Surprisingly enough in view of the
fractal nature of $p$-adic numbers, such a nice $p$-adic geometry does 
exist: it has been built by V. Berkovich (from the categorical
viewpoint, it is essentially equivalent to rigid geometry, in the sense
of Tate and Raynaud). In this framework, the monodromy of differential
equations can be analyzed in the usual way; but a new phenomenon occurs:
it is no longer true that there exists a basis of solutions around every
ordinary point. This approach is therefore limited to a rather special 
class of connections. We shall see in the sequel how the theory of
$p$-adic period mappings provides interesting examples in this class.
 
Following way $2)$, one encounters Dwork's notion of Frobenius
structure, which has often been considered as a plausible substitute for
monodromy in the $p$-adic context. We discuss so-called unit-root
crystals and overconvergence, and illustrate these notions by Dwork's
treatment of the $p$-adic Gamma function $\Gamma_p$, the Gross-Koblitz
formula, and by a detailed study of the $p$-adic hypergeometric function
$F(\frac{1}{2},\frac{1}{2},1; z)$.

\medskip

Abelian periods show themselves in two different ways:

$1)$ as integrals of algebraic differentials over loops, 

$2)$ when the integrand depends on a parameter $z$, as solutions of
certain linear differential equations in $z$ (Picard-Fuchs, or Gauss-Manin).

\medskip

Both ways are passable in the $p$-adic context, and again lead
to completely different theories. 

Following up the first way --- integrals being interpreted as Riemann sums
--- leads to P. Colmez' construction of abelian $p$-adic periods (in the
sense of Fontaine-Messing). These periods relate the Tate module 
($p$-primary torsion points) to the first De Rham cohomology module of
a given abelian variety over a $p$-adic field. They do not live in
$\C_p$, but in a certain $p$-adic ring $\matheur{B}^+_{\mathrm{dR}}$ (of which
$\C_p$ is a quotient), and it is not possible in general to express them 
function-theoretically, when the abelian variety moves in a
family. Nevertheless, we stress some particular cases where this can
actually be done, throwing a bridge between ways $1)$ and $2)$. 

The second way --- Dwork's viewpoint of periods as solutions of
Gauss-Manin connections --- will prevail at the modular level, when
$p$-divisible groups and filtered De Rham modules attached to families
of abelian varieties will be compared not directly, but via their moduli
spaces by the period mapping.

This way of looking at periods lacks arithmetic structure: namely, a
rational, or at least a $\ovl{\Q}$-structure, on the $\C_p$-space of
solutions, to convey some arithmetical meaning to the periods. We
discuss some cases where such a canonical structure exists (`$p$-adic
Betti lattices'), notably the case of abelian varieties with
supersingular reduction. We get here well-defined $p$-adic periods
(completely different from the Fontaine-Messing periods), and compute
them in terms of $\Gamma_p$-values in the case of elliptic curves with
complex multiplication ($p$-adic analogue of the Lerch-Chowla-Selberg
formula).

\bigskip

{\it Chapter II} is an introduction to the geometric theory of $p$-adic
period mappings, in the sense of Drinfeld-Rapoport-Zink.

We have tried to keep prerequisites at a minimum, and to emphasize as
much as possible the analogies between the complex and $p$-adic
contexts. Basic definitions about $p$-divisible groups and crystals are
recalled, and the proof of some basic results is sketched.

\medskip

The theory of period mappings attaches to a family of algebraic
varieties its periods, viewed abstractly as a moving point in a suitable
grassmannian. Due to constraints of Riemann-type, the period mapping
actually takes its values in an open subset of the grassmannian, the
period domain, which is a symmetric domain. It is multivalued, the
ambiguity being described by the projective monodromy of the Gauss-Manin
connection.

The younger $p$-adic theory is far less advanced: at present, there is a
theory of $p$-adic period mappings only for $p$-divisible groups or
closely related geometric objects. For want of wide range, it has
nevertheless gained richness and depth.

\medskip

In the presentation of Drinfeld-Rapoport-Zink, one starts by
constructing a moduli space $\MM$ for $p$-divisible groups which are
(quasi-)isogenous to a given one in characteristic $p$. The {\it
$p$-adic period mapping} $\PP$ then relates this moduli space to a
suitable grassmannian which parametrizes the Hodge filtration in the
Dieudonn\'e module. There are again constraints of Riemann-type, which
force the period mapping to take its values in an open subset of the
grassmannian, the period domain, which is a `symmetric domain'. The best
known example of such a period domain is the Drinfeld space
$\C_p\setminus \Q_p$, a $p$-adic analogue of the double-half-plane
$\C\setminus \R$.

In the situation where the $p$-divisible groups are algebraizable, \ie
come from $p$-primary torsion of abelian varieties parametri\-zed by a
certain Shimura variety $\Sh$, the moduli space $\MM$ provides a
uniformization of some tubular region in the Shimura variety $\Sh$. 

We show that, just as in the complex case, $\PP$ can be described in
terms of quotients of solutions of the associated Gauss-Manin connection
(this feature does not seem to appear in the literature, except in the
old special case of Dwork's period mapping for elliptic curves with 
ordinary reduction). This allows to give explicit formulas `\`a la
Dwork' for $\PP$ in many cases.

The most interesting cases, investigated by Drinfeld-Rapoport-Zink,
arise when $\Sh$ itself (more accurately: a whole connected component) is
uniformized by $\MM$: $\Sh$ is then a quotient of $\MM$ by an arithmetic
discrete group $\Gamma$.  We show that, up to replacing $\Sh$ by a finite
covering (to kill torsion in $\Gamma$), the solutions of the Gauss-Manin 
connection extend to global multivalued functions: in other words, the
Gauss-Manin connection has global monodromy in the 
sense of chapter I, and the projective monodromy group coincides with
$\Gamma$. 

We review the case of Shimura curves attached to quaternion algebras
over totally real fields (\v{C}erednik-Drinfeld-Boutot-Zink):
there is a global $p$-adic uniformization when $p$ divides the
discriminant of the quaternion algebra. We then construct, using the
$p$-adic Betti lattices of chapter I, a canonical $\ovl{\Q}$-space of
solutions of the Gauss-Manin connection which is stable under global
monodromy.

 \bigskip {\it Chapter III} explores the group-theoretic aspects of the
 theory of period mappings, in the $p$-adic context.

 The modern presentation of ramified coverings, uniformization, and of
 the Gauss-Riemann-Fuchs-Schwarz theory uses the notions of orbifold and
 uniformizing differential equation. We develop $p$-adic counterparts of
 these notions.

\medskip

The right notion of ramified covering to adopt is not obvious, because
\'etale coverings are only exceptionally topological coverings, in the
$p$-adic context. Ber\-ko\-vich's geometry provides a very convenient
framework for this kind of problems, since his spaces are locally arcwise
connected.  In the first section, we discuss the formalism of
fundamental groups attached to categories of (possibly infinite) \'etale
coverings satisfying simple axioms. Although the definition of these
topological groups is simple and natural, their topology itself may be
very complicated (not locally compact in general).

It turns out that there are too few topological coverings and too many
\'etale coverings, in general, to provide a workable theory of
monodromy. To remedy this, we introduce in the second section the notion
of {\it temperate \'etale covering}. Such coverings are essentially
built from (possibly infinite) topological coverings of finite \'etale
coverings. They are classified by temperate fundamental groups, which
seem to be the right $p$-adic equivalents, in many ways, of fundamental
groups of complex manifolds. These groups are not discrete in general,
but nevertheless often possess many infinite discrete quotients. We
present a large selection of examples.

\medskip Temperate \'etale coverings are well-suited to the definition
of $p$-adic orbifold charts and of the category of {\it $p$-adic
orbifolds} (they replace the unramified coverings, over $\C$). As in the
complex case, orbifolds and orbifold fundamental groups are the tools
for a theory of ramified coverings, developed in section 4.

Before turning to the $p$-adic analogue of the {\it uniformizing
differential equation} attached to an orbifold of dimension one (via
Schwarzian derivatives), we discuss local and global monodromy of
$p$-adic differential equations in section 3.  
 
We outline the Christol-Mebkhout theory of $p$-adic slopes of
differential equations over annuli, and the relation with Galois
representations of local fields of characteristic $p$.

We then define and study the {\it $p$-adic \'etale Riemann-Hilbert
functor}, which attaches a vector bundle with integrable connection to
any discrete representation of the \'etale fundamental group of a
$p$-adic manifold; connections in the image are characterized by the
fact that the \'etale sheaf of germs of solutions is locally
constant. This is a vast generalization of the phenomena of global
monodromy studied in Chapter I; for instance, the differential equation
$y'=y$ over the affine line belong to this class.

Uniformizing differential equations of orbifolds of dimension one also
belong to this class, and the representation actually factors through
the temperate fundamental group.  The case of a Shimura orbifold is of
special interest: the period mapping (complex or $p$-adic) is given by a
quotient of two solutions of a fuchsian differential equation defined
over a number field, which can be interpreted as uniformizing
differential equation either over $\C$ or over $\C_p$.

In Section 5, we examine the case of Schwarz orbifolds: the projective
line with $0,1,\infty$ as branched points (endowed with suitable
finite multiplicities). Over $\C$, the uniformizing differential
equations are of hypergeometric type, with projective monodromy group
identified with the orbifold fundamental group, namely with a
(cocompact) triangle group.

In the $p$-adic case, we define {\it $p$-adic triangle groups} to be
projective monodromy groups of those hypergeometric differential
equations which are in the image of the $p$-adic \'etale Riemann-Hilbert
functor: thus by definition, there exists a finite \'etale covering of
$\C_p\setminus \{0,1\}$ over which the hypergeometric function extends
to a global multivalued analytic function.

We give a purely geometric description of these discrete subgroups of
$PSL_2(\C_p)$ (without reference to differential equations; it is
perhaps here that we are closest to Fuchs and Schwarz).  From this, and
recent combinatorial work by F. Kato, it follows that infinite $p$-adic
triangle groups exist only for $p\leq 5$.

We then construct the $p$-adic analogues of Takeuchi's list of
arithmetic triangle groups, \ie the list of all ``arithmetic'' $p$-adic
triangle groups for every $p$, using the 
\v{C}erednik-Drinfeld-Boutot-Zink uniformization of Shimura curves and
$p$-adic period mapping. Special values of the corresponding
hypergeometric functions at CM points are expressed in terms of
$\Gamma_p$-values.

\medskip

\subsection*{On the style.} 

\medskip

The first chapter is rather down-to-earth. It has kept
something of the informal style of lecture notes, and this also holds to
a lesser extent for the second chapter; the level is inhomogeneous, the
pace may sometimes be brisk, proofs are often omitted and replaced by
references.  In contrast, the last and longest chapter is devoted to a
systematic exposition of new material.

We hope that our constant function-theoretic viewpoint brings some unity
to the exposition, throughout the chapters.  We have tried to keep them
(and even the sections) as logically independent as possible. Thus, for
example:

only subsections 1.5 and 5.3 of chapter I are needed in the sequel;

until section 7, chapter II is almost self-contained;

until subsection 4.7, chapter III is almost self-contained;

reading fragments of III should be enough to grasp the ins and outs
of $p$-adic triangle groups.
 
\medskip

Sections I.3 and III.6 are small pieces of `computational mathematics'
(without computer) intended as testing ground for notions developed in
the preceding sections.

\bigskip

{\bf Acknowledgments.}
It is a pleasure to acknowledge here our main sources of inspiration:

first of all, V. Berkovich's views on $p$-adic geometry: without this
extremely convenient way of seeing $p$-adic spaces, our systematic
to-ings and fro-ings between complex and $p$-adic worlds would have 
been impossible;

Rapoport-Zink's book on period spaces for $p$-divisible groups has been
a constant reference for Chapter II;

in Chapter III, our treatment has been much influenced by the work of
M. van der Put on $p$-adic discontinuous groups, by J. De Jong's study
of \'etale coverings, and no less by M. Yoshida's orbifold viewpoint on
the (complex) Gauss-Schwarz theory.

\medskip

Concerning this book, I am much indebted to F. Kato in many ways: not
only did he organize the T\={o}hoku and the Hokkaid\={o} lectures
series, asking me to speak on these matters, but he also suggested the
very project of writing these notes, and influenced a lot both the
content and the style by his numerous questions and comments; I have
been much encouraged in writing them by his constant interest and by his
own recent work on $p$-adic triangle groups.

\medskip

The T\={o}hoku lectures, as well as the birth of some ideas presented in
III (especially $p$-adic triangle groups), took place during my stay in
Japan thanks to a fellowship from the Japan Society for the promotion of
Science. I am very grateful to this institution, and to H. Shiga who
arranged the stay.

\medskip

I thank M. Matignon, F. Baldassarri and S. Bosch for invitations to
lecture on these topics at the Universities of Bordeaux, Padova, and
M\"unster respectively.

 \medskip

I am grateful to F. Kato and N. Tsuzuki for writing the appendices of
this book. I am also very grateful to S. Matsuda who kindly took care,
with great patience, of the nightmarish electronic harmonization of the
main text and the appendices, and drawing the figures.

 \medskip

Thanks are also due to P. Colmez and O. Gabber for several useful
comments, and to J.F. Boutot for sending me overseas his joint preprint
with T. Zink about uniformization of Shimura curves.

I also thank very much P. Bradley for his careful reading of a version
of chapter III and for kindly providing a long list of errata (of
course, I am responsible for any remaining error).

\chapter{Analytic aspects of $p$-adic periods.}\label{chap-ana-p-period}
\addtocontents{toc}{\protect\par\vskip2mm\hskip20mm(Analysis)\par\vskip5mm}
\minitoc
\newpage
\section{Analytic continuation; topological point of view.}\label{sec-1}

\begin{abst}
 Thanks to Berkovich's presentation of $p$-adic analytic geometry, it is
 possible to make sense of the familiar monodromy principle in the
 exotic world of $p$-adic manifolds. Its application is however more
 limited than in the classical case, because (1) sheaves of solutions of
 linear differential equations 
 are usually not locally constant, and (2) many spaces (for instance, 
 annuli) are simply connected.
\end{abst}

\subsection{The Monodromy principle.}
 
Let $S$ be a topological space, and $\mathcal{F}$ an abelian sheaf on $S$.
For an open set $U\subseteq S$ and a section $f\in\Gamma(U,\mathcal{F})$, the 
support of $f$ is the subset $\Supp(f)=\{u\in U\ |\ f_u\neq 0\}$, 
which is easily seen to be closed in $U$.

\begin{dfn}\label{dfn-unique-cont}
We say that the sheaf $\mathcal{F}$ satisfies the {\it principle of unique
continuation} if for any open set $U\subseteq S$ and any section 
$f\in\Gamma(U,\mathcal{F})$ the support $\Supp(f)$ is open in $U$.
\end{dfn}

Note that the principle of unique continuation is a local property.

\begin{lem}\label{lem-unique-cont}
If $\mathcal{F}$ satisfies the principle of unique continuation, then
any two sections $f,g\in\Gamma(U,\mathcal{F})$
on a connected open subset $U\subseteq S$ coincide if (and only if)
their germs $f_s$, $g_s$ at some point $s\in U$ coincide.
The converse also holds if $S$ is locally connected.
\end{lem} 

\begin{proof}
 The first assertion is clear. For the converse, let $u$ be a point adherent to the complement 
 of $\Supp(f)$. There exists a sufficiently small connected open
 neighborhood $V$ of $u$ contained in $U$. Since $f_s=0$ at some
 point $s\in V$, $f\equiv 0$ on $V$. Hence $u\not\in\Supp(f)$. This shows that $\Supp(f)$ is open.
\end{proof}

\begin{lem}\label{lem-separated}
 Let $p: F\rightarrow S$
 be the local homeomorphism canonically attached to
 $\mathcal{F} $: $\mathcal{F}(U)$ is the set of (continuous) sections of $p$ over $U$. 
 If $F$ is separated (\ie Hausdorff), then $\mathcal{F}$ satisfies the
 principle of unique continuation. The converse also holds if $S$ is
 separated.
\end{lem}

\begin{proof} If $F$ is separated, then for any two sections $f, g$ of
 $p$ over an open subset $U\subset S$, the set of $s\in U$ such that
 $f(s)=g(s)$ is closed.  Applying this to the zero-section $g=0$, we see
 that the support of $f$ is open.

 Conversely, let $x,y $ be two points of $F$. If $p(x)\neq p(y)$, there
 are disjoint open neighborhoods $V_x, V_y$ of $p(x)$ and $p(y)$
 respectively in the separated space $S$; then $p^{-1}V_x$ and
 $p^{-1}V_y$ are disjoint open neighborhoods of $x$ and $y$
 respectively. If $p(x)= p(y)$, there are open neighborhoods $U_x, U_y$
 of $x$ and $y$ respectively, an open subset $V\subset S$, and sections
 $f,g$ of $p$ over $V$, such that $f(V)=U_x,\;g(V)=U_y$. By the
 principle of unique continuation, $f\neq g$ defines an open subset
 $W\subset V$ containing $p(x)=p(y)$; then $f(W)$ and $g(W)$ are
 disjoint open neighborhoods of $x$ and $y$ respectively.
\end{proof}

\begin{dfn}
 We say that the sheaf $\mathcal{F}$ 
 satisfies the {\it monodromy principle} if it has the following property:

 let $\Gamma: [a,b]\times [0,1] \rightarrow S$ be any 
 continuous map with
 $\Gamma(\{a\}\times [0,1])=\{x_a\}, \Gamma(\{b\}\times
 [0,1])=\{x_b\}$. Let $f_{x_a}$ be an element of the stalk
 $\mathcal{F}_{x_a}$. Assume that for any
 $t\in [0,1]$, $f_{x_a}$ extends to a global section $f^{\Gamma_t}$ of
 $\Gamma_t^{\ast} \mathcal{F}$ on
 $[a,b]$. Then this extension is unique and
 $f^{\Gamma_t}(b)\in \mathcal{F}_{x_b}$ is independent of $t$.

 In this situation, the section $f^{\Gamma_t}$ of 
 $\Gamma_t^{\ast} \mathcal{F}$ is called the
 {\it continuation of $f_{x_a}$ along the path $\Gamma_t$}, and 
 $f^{\Gamma_t}(b)\in \mathcal{F}_{x_b}$ its {\it value at $x_b$}.
\end{dfn}

\begin{pro}
 If $\mathcal{F}$ satisfies the principle of
 unique continuation, then it satisfies the principle of monodromy. 
 The converse also holds if
 $S$ is locally arcwise connected.
\end{pro}

\begin{proof}
 Assume that $\mathcal{F}$ satisfies the principle of
 unique continuation. This guarantees the uniqueness of the 
 continuation $f^{\Gamma_t}$ of
 $f_{x_a}$ along the path $\Gamma_t$ for any $t$. Moreover, by a 
 special case of the proper base change theorem
 (applied to the first projection $[a,b]\times [0,1]\rightarrow 
 [a,b]$, \cf \cite[IV, 1.4]{cos}) $f^{\Gamma_t}$
 extends to a section of
 $\Gamma^{\ast}
 \mathcal{F}$ on a suitable subset of the form $[a,b]\times 
 (t-\epsilon,t+\epsilon)$. By unicity, these sections glue
 together to a global section of $\Gamma^{\ast} \mathcal{F}$. The restriction 
 of this sheaf to
 $\{b\}\times [0,1]$ is the constant sheaf with stalk
 $\mathcal{F}_{x_b}$. Therefore the value 
 $f^{\Gamma_t}(b)\in \mathcal{F}_{x_b}$ is independent of $t$.

 Conversely, let $U$ be an arcwise connected open subset of $S$, let
 $x_a, x_b$ be two points of
 $U$, and let $f$ be a section of
 $\mathcal{F}$ over $U$. Let $\gamma: [a,b]\rightarrow S$ be any path from 
 $x_a$ to $ x_b$. It
 follows from the unicity of continuation of $f_{x_a}$ along $\gamma$ 
 (requested in the monodromy
 principle) that if $f_{x_a}$ is
 $0$, so is
 $f_{x_b}$.
\end{proof}

\begin{exa}\label{lem-homotopy}
  (1) Any locally constant abelian sheaf $\mathcal{F}$ on a
 topological space $S$ satisfies the principle of unique 
 continuation, hence the principle of monodromy.
 
 In fact, for any $\Gamma: [a,b]\times [0,1]\rightarrow S$ as in
 \lemref{lem-separated}, the inverse image
 $\Gamma^{\ast}\mathcal{F}$ is locally constant, hence constant and 
 canonically isomorphic to the
 constant sheaf attached to
 $\mathcal{F}_{x_a}$. Therefore the extension $f^\gamma$ of any 
 $f_{x_a}\in \mathcal{F}_{x_a}$ along any
 path $\gamma$ exists, and the value at the other extremity 
 $x_b=\gamma(b)$ depends only on the
 \emph{homotopy class} of $\gamma$.

  (2) When $S$ is a complex manifold, the structure sheaf
 $\mathcal{O}_S$ satisfies the principle of unique continuation, hence
 the principle of monodromy.
\end{exa}

\begin{para}
 Let $S$ be a topological space, connected and locally arcwise 
 (or simply) connected,
 and $\mathscr{F}$ a locally constant abelian sheaf on $S$.
 We fix a point $s\in S$, and denote by $\pi_1(S,s)$ the fundamental
 group based at $s$.\footnote{to remove any ambiguity, let us say that we
 adopt Deligne's convention: the composition in $\pi_1(S,s)$
 is induced by the juxtaposition of loops in the reverse order. As 
 such, $\pi_1(S,s)$ acts \emph{on the right} on the
 pointed universal covering.}
 To any loop $\gamma\colon [0,1]\rightarrow S$ based at $s$, let us associate
 the so-called {\it monodromy} along $\gamma$, defined by the composite
 \begin{equation*}
  \mathscr{F}_s\stackrel{\sim}{\longrightarrow}(\gamma^{\ast}\mathscr{F})_0
   \stackrel{\sim}{\longrightarrow}\Gamma([0,1],\gamma^{\ast}\mathscr{F})
   \stackrel{\sim}{\longrightarrow}(\gamma^{\ast}\mathscr{F})_1
   \stackrel{\sim}{\longrightarrow}\mathscr{F}_s.
 \end{equation*}
 By \examref{lem-homotopy} (1), the monodromy
 $\mathscr{F}_s\xrightarrow{\sim}\mathscr{F}_s$
 depends only on the class of $\gamma$ in
 $\pi_1(S,s)$,
 hence gives rise to a left $\Z[\pi_1(S,s)]$-module structure on 
 $\mathscr{F}_s$.
 This construction yields an equivalence of the categories
 \begin{equation*}
  \left\{\textrm{locally constant sheaves on} \; S\right\}
   \stackrel{\sim}{\longrightarrow}
   \left\{\textrm{left\ }\Z[\pi_1(S,s)]\textrm{-modules}\right\}\;:
 \end{equation*}
  giving a locally constant sheaf amounts to giving its value at $s$
 together with the monodromy action.
\end{para}

\subsection{Rigid geometry and the problem of unique continuation.}
\label{subsection-rigid-and-unique-cont}
\begin{para}
 For $p$ a prime number, let $\Q_p$ denote as usual the
 completion of $\Q$ for the $p$-adic absolute value
 $|\;\;|_p$ : $|p^n\frac{a}{b}|_p = p^{-n}$ if the rational integers $a,
 b$ are prime 
 to $p$. This ultrametric absolute value extends in a unique way to each
 finite extension of $\Q_p$. These finite extensions are locally
 compact and totally disconnected. They are all complete, but ``the''
 algebraic closure $\ovl{\Q}_p$ of $\Q_p$ itself is not complete. Its
 completion, denoted by $\C_p$, turns out to be algebraically
 closed, and plays the role of $\C$ in $p$-adic analysis.

 \noindent In the sequel, $\matheur{D}(a,r^+)$
 (\resp $\matheur{D}(a,r^-)$) stands for the disk --- archimedean or not 
 --- of radius $r$ centered at $a$ with (\resp without) circumference.

 It might be surprising at first that geometries can be built upon
 $p$-adic numbers, whose ``fractal'' nature makes them hardly amenable to
 intuition as a continuum.  Nevertheless, Bourbaki's presentation of
 analytic geometry \cite{vdea} treats the real, complex and $p$-adic
 cases on equal footing. This approach is based on a local definition of
 analytic functions as sums of convergent power series. Its major
 drawback is that these analytic functions fail to satisfy the principle
 of unique continuation, essentially because two ultrametric disks are
 either concentric or disjoint (like drops of mercury)

 \begin{figure}[h]
 \begin{picture}(250,60)(0,0)
  \put(0,0){\includegraphics{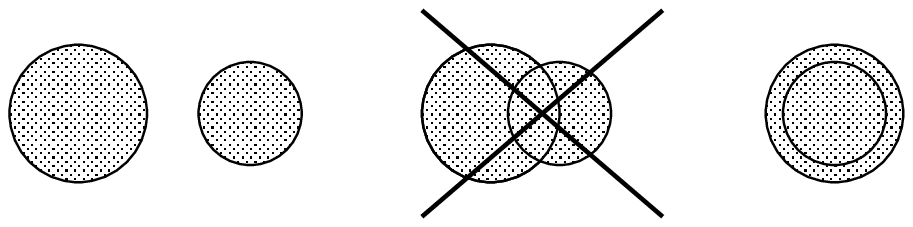}}
 \end{picture}
 \caption{}
 \label{fig1}
 \end{figure}

\end{para}

\begin{para}
 M. Krasner had the idea to remedy this by using a definition of
 analytic functions \emph{\`a la Runge}.
 He properly founded ultrametric analysis by introducing his \emph{analytic
 elements} defined as
 uniform limits of rational functions, a global notion which overcomes,
 to some extent, problems stemming from the disconnectedness of the $p$-adics.

 The next step is due to J.~Tate. In order to deal with more general spaces
 than just subsets of the line, he introduced and developed the
 so-called \emph{rigid analytic geometry} (as opposed to Bourbaki's
 ``wobbly'' analytic geometry), based on affinoid algebras (topological
 $K$-algebras isomorphic to quotients of rings of restricted formal
 power series, \ie whose coefficients tends to $0$) and a suitable
 Grothendieck topology \cite{ras}.

 In a connected rigid analytic variety $S$, one has the following 
 avatar of analytic continuation \cite[0.1.13]{crecrasp}:  assume that
 there is no admissible covering formed by two disjoint non-empty open
 subsets, and let $f$ be an analytic function 
 on $S$. If there is a connected open subset $U$ over which $f$ 
 vanishes, then $f=0$.
\end{para}

\begin{para}
 The first achievement of rigid geometry was Tate's representation of an
 elliptic curve with ``bad reduction at $p$'' as a rigid analytic quotient
 $\C_p^{\times}/q^{\Z}$ of the multiplicative group by the discrete
 group generated by $q$.  Here $q$ is given by the usual series
 $\frac{1}{j}+\ldots$ in the $j$-invariant of the elliptic 
 curve ($|j| >1$), interpreted $p$-adically.  The quotient is
 obtained by gluing two affinoid annuli of width $|q|^{-1/2}$
 along their boundary; the inverse image in $\C_p^{\times}$ of each of
 these annuli consists of countably many disjoint copies of it.  This is
 analogous to Jacobi's partial uniformization $\C^{\times}/q^{\Z}$
 of a complex elliptic curve.
\end{para}

\begin{para}
 This result was extended by D.~Mumford to
 curves of higher genus \cite{aacodcoclr}.  Here, the
 partial uniformization which serves as a complex model is the 
 Schottky uniformization of complex curves of genus
 $g>1$, which we briefly recall (\cf \cite{osats}).  Let
 $D$ be an open subset of ${\mathbb{P}}^1(\C)$, limited by $2g$ 
 disjoint circles $C_1,C'_1,\ldots,C_g,C'_g$.  One
 assumes the existence, for any $i=1\ldots, g,$ of an element 
 $\gamma_i\in PGL_2(\C)$, with two distinct fixed
 points, which sends $C_i$ to $C'_i$ and $D$ outside itself. Then the 
 subgroup $\Gamma$ generated by the
 $\gamma_i$'s is discrete and freely generated by them, the complement 
 $\Omega\subset {\mathbb{P}}^1(\C)$ of the
 topological closure of the set of fixed points of
 $\Gamma$ is an open dense subset ($=\bigcup_{\gamma \in \Gamma} 
 \gamma(D)$), and $\Omega/\Gamma$ is a projective
 smooth complex curve of genus $g$.

 When $\C$ is replaced by $\C_p$, the same construction applies. The
 rigid analytic quotient $\Omega/\Gamma$ exists and is called a Mumford
 curve. Among projective smooth curves of genus $g$ over $\C_p$, Mumford
 curves are characterized by the existence of a reduction over the
 residue field ${\ovl{\F}}_p$ such that every irreducible component is
 isomorphic to ${\mathbb{P}}^1$ and intersect the others at double
 points, \cf \cite{sgamc}.

 Let us mention that T. Ichikawa has proposed a unified theory of 
 the archimedean and non-archimedean
 Schottky-Mumford uniformizations, \cf \cite{sutorsamcoig}.
\end{para}

\subsection{Berkovich geometry and the principle of monodromy}

\begin{para}
 Rigid analytic spaces are endowed with a 
 Grothendieck generalized topology, and their
 structure sheaf is a sheaf with respect to this topology. Hence
 it cannot be said to satisfy the principle of
 unique continuation in the strict sense of \defref{dfn-unique-cont}.
 Moreover, there is no non-trivial path in such spaces.  Therefore,
 rigid geometry is not a suitable setting for discussing $p$-adic 
 analytic continuation in an intuitive way.

 In contrast, Berkovich's viewpoint on $p$-adic geometry does 
 not suffer from these drawbacks: Berkovich's
 analytic spaces are genuine locally ringed topological spaces, which 
 are even locally arcwise connected. We refer to
 \cite{spasalc} for a compact technical presentation.
\end{para}

\begin{para}
 The buildings blocks are the same: affinoid algebras $\matheur{A}$
 (called strictly affinoid algebras by Berkovich), \ie topological
 $K$-algebras isomorphic to quotients of rings of restricted power
 series. 
 However, instead of attaching to $\matheur{A}$ its maximal spectrum
 $\mathrm{Spm}(\matheur{A})$ as in rigid geometry, 
 Berkovich analytic geometry deals with the {\it affinoid space}
 $\matheur{M}(\matheur{A})$ of all bounded multiplicative seminorms
 on $\matheur{A}$.

 This ``spectrum'' contains ``more points'' than the maximal spectrum
 $\mathrm{Spm}(\matheur{A})$, namely something 
 like ``generic points'', which ``complete'' $\mathrm{Spm}(\matheur{A})$.
 There is a natural inclusion $\mathrm{Spm}(\matheur{A})\subseteq 
 \matheur{M}(\matheur{A})$ which identifies $\mathrm{Spm}(\matheur{A})$
 with the subspace of all seminorms $|\cdot|$ with
 $\matheur{A}/\Ker |\cdot| =\C_p$.
 This ``completion'', in fact, simplifies the topology; e.g.\ Berkovich
 analytic spaces are locally arcwise connected. 
\end{para}

\begin{para}
 {\itshape Berkovich's affinoids.\ }
 Let us be a little more precise about the definition of 
 {\it affinoid spaces} in Berkovich geometry.
 Let $\matheur{A}$ be an affinoid algebra over a complete subfield $K$ of 
 $\C_p$.

 \begin{enumerate}
  \item A point of $\matheur{M}(\matheur{A})$ is a bounded multiplicative
	seminorm on $\matheur{A}$.
  \item The topology of $\matheur{M}(\matheur{A})$ is the weakest one so
	that the mapping
	$\matheur{M}(\matheur{A})\ni\chi\mapsto\chi(f)\in \R_{\ge 0}$
	is continuous for any $f\in\matheur{A}$.
  \item The sheaf of rings $\OO_{\matheur{M}(\matheur{A})}$ is defined by 
	$\Gamma(U,\OO_{\matheur{M}(\matheur{A})})=\varprojlim\matheur{A}_V$, 
	where $V$ runs over finite unions $\bigcup V_i$ of affinoid
	domains contained in $U$, and
	$\matheur{A}_V= \Ker\bigl(\prod_i\matheur{A}_{V_i}
	\rightrightarrows
	\prod_{ij}\matheur{A}_{V_i\cap V_j}\bigr)$, $\matheur{A}_{V_i}$
	being the affinoid algebra attached to the affinoid domain $V_i$.
 \end{enumerate}

 The value of a ``function'' $f\in \matheur{A}$ at a point $\chi\in 
 \matheur{M}(\matheur{A})$ is its image in the field
 $\matheur{A}/\Ker\chi$. This field inherits the absolute value induced by 
 the seminorm $\chi$, and its completion is
 denoted by $\mathscr{H}(\chi)$.\index{000H@$\mathscr{H}(\;)$}
\end{para}

\begin{exa}\label{exa-berkovich-space}
Let $K\{t\}$ denote the ring of restricted power series in one
variable.
Let us assume, for simplicity, that $K$ is algebraically closed.
In rigid geometry, the maximal spectrum $\mathrm{Spm}(K\{t\})$ is
just the closed disk $\matheur{D}(0,1^{+})$ in $K$ of radius $1$ in the usual sense.
In $\matheur{M}(K\{t\})$, the following four kinds of points occur:
\begin{enumerate}
 \renewcommand{\labelenumi}{(\arabic{enumi})}
 \item classical points (\ie those coming from $\mathrm{Spm}(K\{t\})$):
       $x\in\matheur{D}(0,1^{+})$, $\chi_x(f)=|f(x)|_K$ for $f\in K\{t\}$,
 \item  generic points of disks: $\chi=\eta_{x,r}$ for $0<r\leq 1$ with 
       $r\in |K^{\times}|$, 
       $\chi(f)=|f|_{\matheur{D}(x,r^{+})}$ (the sup-norm), 
 \item the same, for $r\not\in |K^{\times}|$,
 \item generic points of infinite decreasing families
       $\{\matheur{D}_{\alpha}\}$
       of closed disks with radius $\leq 1$: 
       $\chi(f)=\mathrm{inf}|f|_{\matheur{D}_{\alpha}}$.
\end{enumerate}

\vspace{2ex}
\begin{figure}[h]
 \setlength{\unitlength}{1pt}
 \begin{picture}(150,60)(0,0)
  \put(25,10){\line(1,1){35}}
  \put(95,10){\line(-1,1){35}}
  \put(25,9){\makebox(0,0)[c]{$\bp$}}
  \put(95,9){\makebox(0,0)[c]{$\bp$}}
  \put(60,44){\makebox(0,0)[c]{$\bp$}}
  \put(45,29){\makebox(0,0)[c]{$\bp$}}
  \put(55,58){\makebox(0,0)[l]{$\eta_{x,|x-y|}$}}
  \put(17,0){$x$}
  \put(96,0){$y$}
  \put(20,35){$\matheur{\eta}_{x,r}$}
 \end{picture}
 \caption{Generic points and paths on a Berkovich space.}
 \label{pic1}
\end{figure}
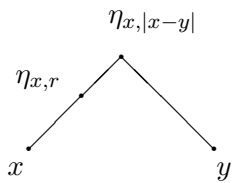
\vspace{1ex}

Therefore, two arbitrary distinct points in $\matheur{M}(K\{t\})$ can be 
connected by a unique path.
For example, two points of type (1), $x$ and $y$, are connected by a
path consisting of points of type (2) and (3) associated to the disks
$\matheur{D}(x,r^{+})$ and $\matheur{D}(y,r^{+})$ for $0<r\leq |x-y|$.
The complement of a point of type (1) or (4) is connected.
The affine line $\mathbb{A}^1$ is the union of the affinoid spaces
associated to the algebras of power series convergent in 
$\matheur{D}(0,r^{+})$, and the projective line $\mathbb{P}^1$ is the
Alexandroff compactification of $\mathbb{A}^1$.
\end{exa}

\begin{rem}
 Adding generic points to rigid spaces 
 to come up with Berkovich spaces simplifies the topology in contrast to 
 what happens in algebraic geometry, when
 generic points are added to varieties to produce schemes.
 Actually Berkovich's generic points are closed, unlike Grothendieck's
 ones.
\end{rem}

\begin{para}\label{para:strict-analytic}
 The construction of (strictly) analytic $K$-spaces by gluing affinoid
 spaces together is a little delicate (one is not gluing open
 subspaces); we refer to \cite{ecfnas}.  These spaces are \emph{locally 
 compact, locally countable at infinity, and locally arcwise connected.} 

 There is a fully faithful functor between Berkovich's Hausdorff 
 (strictly) analytic $K$-spaces and rigid analytic varieties: at the
 level of underlying sets, this functor sends a space $S$ to the subset of
 ``classical points'', \ie points $s$ for which $[\mathscr{H}(s):K]<\infty$.

 This functor establishes \emph{an equivalence of categories 
 between the category of paracompact (strictly)
 analytic $K$-spaces and quasi-separated rigid spaces over $K$ having 
 an admissible affinoid covering of finite type}.

 This equivalence respects the notion of dimension
 (topological dimension in the case of Berkovich analytic spaces), as
 well as the standard properties of local rings such as: reduced,
 normal, smooth..., and the properties: finite, \'etale...  of
 morphisms. Also, a paracompact analytic space is connected (in the
 usual sense) if and only if the corresponding rigid space does not
 admit an admissible covering by two disjoint non-empty open subsets.

 Furthermore, there is a canonical functor ``analytification''
 \begin{equation*}
  \left\{
  \begin{matrix}
   \text{separated schemes}\\
   \text{locally of finite type}\\
   \text{over $K$}
  \end{matrix}
  \right\}
  \rightarrow
  \left\{
  \begin{matrix}
   \text{paracompact (strictly)}\\
   \text{analytic $K$-spaces}
  \end{matrix}
  \right\}
 \end{equation*}
 \noindent
 and a canonical functor ``generic fiber''
 \begin{equation*}
  \left\{
  \begin{matrix}
   \text{separated formal schemes locally}\\
   \text{finitely presented over $\mathcal{O}_K$}
  \end{matrix}
  \right\}
  \rightarrow
  \left\{
  \begin{matrix}
   \text{paracompact (strictly)}\\
   \text{analytic $K$-spaces}
  \end{matrix}
  \right\}.
 \end{equation*}
\end{para}

\begin{para}\label{para:p-adic-manifold}
 {\itshape $p$-adic manifolds.} We shall be mainly 
 concerned with paracompact (strictly) $K$-analytic spaces $S$ which
 satisfy the following assumption:
 \begin{quote}
  any $s\in S$ has a neighborhood $U(s)$ which is isomorphic 
  to an affinoid subdomain of some space
  $V_s$ which admits locally an \'etale morphism to the affine space 
  $\mathbb{A}^{\dim S}$.
 \end{quote}
 For convenience, we shall call such a space $S$ a {\it
 $K$-manifold}\index{Kmanifold@$K$-manifold}, and simply 
 a {\it $p$-adic manifold}\index{padicmanifold@$p$-adic manifold} if $K=\C_p$.

 \noindent It is an important result of V. Berkovich \cite{spasalc}
 that \emph{$p$-adic manifolds are locally contractible}.  Therefore
 they have universal coverings. 
\end{para}

\begin{pro}  
 The structure sheaf of any $p$-adic 
 manifold $S$ satisfies the principle of unique
 continuation, and (equivalently) the principle of monodromy.
\end{pro}
\begin{proof} 
 Since the principle of unique continuation is local, we may replace $S$
 by the neighborhood $U(s)$ of any 
 point $s$ given in advance. Thus we may assume that $S$ 
 itself is affinoid: $S=\matheur{M}(\matheur{A})$,
 and connected, and we have to show that the homomorphism 
 $\iota_s:\;\matheur{A}\rightarrow \matheur{A}_s\simeq
 \underset{s\in V}{\liminj}\matheur{A}_V$ is injective ($V$ runs over
 the affinoid neighborhoods of $s$).  This follows from the
 fact that the homomorphism of completion at 
 $s : \matheur{A}\rightarrow  \widehat{\matheur{A}}_s =
 \underset{n}{\limproj}\matheur{A}/\mathcal{I}_s^n\;$ is injective and factors
 through $\iota_s$.
\end{proof}  

In the case of curves (\ie $p$-adic manifolds of dimension $1$),
analytic continuation is particularly intuitive, because there is a basis
of open subsets $\mathcal{U}$ with finite boundary such that two arbitrary
points in $\mathcal{U}$ are connected by a unique geometric path lying in
$\mathcal{U}$.


\subsection{Topological coverings and \'etale coverings.}
\label{sub-top-cov-et-cov}
\begin{para}\label{para-complex-vs-padic}
 Complex manifolds are locally contractible,
 and have universal coverings.  There is no need to distinguish between
 topological
 coverings and \'etale coverings (finite or infinite).

 In $p$-adic rigid geometry, the situation is more complicated. It is
 natural to call
 topological covering any morphism $f: Y \rightarrow X$ such that there is
 an admissible
 cover $(X_i)$ of $X$ and an admissible cover $(Y_{ij})$ of $f^{-1}(X_i)$
 with disjoint $Y_{ij}$
 isomorphic to $(X_i)$ via  $f$. Indeed, such topological coverings
 correspond to locally constant
 sheaves of sets on $X$. It is still true that \emph{topological coverings
 are \'etale, but the converse is wrong}, even if one restricts to
 finite surjective morphisms. Indeed, the Kummer covering
 $z \mapsto z^n$ of the punctured unit disk is an \'etale 
 covering, but not a topological covering, if $n>1$.
\end{para}

\begin{para}
 It is again more convenient to tackle these questions in the framework
 of Berkovich's geometry.  For instance, one sees immediately that a
 Kummer covering $z\mapsto z^n$ as above is not a topological covering
 because a classical point has $n$ preimages, while the ``generic
 point'' $\eta_{0,1}$ (corresponding to the sup-norm on the disk) is its
 own single preimage. 
 Topological coverings of a $p$-adic manifold $X$ are defined in the
 usual way; they correspond to locally constant sheaves of sets on $X$.
 They coincide \cite[3.3.4]{staagonf} with topological coverings of the
 rigid analytic variety associated to $X$.
\end{para}

\begin{para}
 There are only three one-dimensional simply connected complex manifolds up to
 isomorphism: the projective line ${\P}^1(\C)$, the affine line $\C$, and
 the disk $\matheur{D}$ ($\simeq \mathfrak{h}$, the complex upper half
 plane). In contrast, there are plenty of \emph{simply-connected}
 one-dimensional 
 $p$-adic manifolds, for instance: any annulus, the line deprived from a
 finite number of points, more generally any connected $p$-adic manifold
 homeomorphic to a \emph{subset of} $\;{\P}^1(\C_p)$, any smooth
 projective curve with ``good reduction''... (all these Berkovich spaces
 look like ``bushy trees'').
\end{para}

\begin{para}
 One defines the (discrete) topological fundamental group
 $\pi_1^{\mathrm{top}}(S,s)$ of
 a pointed $p$-adic manifold $(S,s)$ in the usual way. The general
 topological theory of
 coverings applies: $(S,s)$ is naturally isomorphic to the quotient of the
 pointed universal
 covering $({\til{S}},{\tilde{s}})$ by $\pi_1^{\mathrm{top}}(S,s)$, and
 $\pi_1^{\mathrm{top}}(S,s)$ classify the
 topological coverings of $(S,s)$.   In the one-dimensional case, the
 topological fundamental
 group is a discrete \emph{free} group (more precisely, it is isomorphic to
 the fundamental
 group of the dual graph of the so-called semistable reduction of $S$
 (\cite[5.3]{efgonas})). 

 For instance, if $S=\Omega/\Gamma$ is the
 uniformization of a Mumford curve
 (\cf~\ref{subsection-rigid-and-unique-cont}), $\Omega$ is the universal
 covering of $S$ and the Schottky group $\Gamma$ is isomorphic to the 
 topological fundamental group.

 On the other hand, one can define the (profinite) algebraic fundamental
 group $\pi_1^{\mathrm{alg}}(S,s)$ \`a la Grothendieck, classifying the finite
 \'etale coverings of $(S,s)$. In contrast to the complex situation, the
 natural map $\pi_1^{\mathrm{alg}}(S,s) \rightarrow
 \pi_1^{\mathrm{top}}(S,s)^{\wedge}$
 to the profinite completion of $\pi_1^{\mathrm{top}}(S,s)$ is
 \emph{generally not injective} (\eg for annuli).
\end{para}

\subsection{Connections with locally constant sheaves of solutions.}
\label{sub-connection}
\begin{para}
 Let us briefly recall the complex situation. Let $S$ be a complex
 connected manifold, $(\mathcal{M},\nabla)$ a vector bundle of rank $r$
 with integrable connection on $S$. The classical Cauchy theorem shows
 that for any $s \in S$, the solution space
 $(\mathcal{M}\otimes\mathcal{O}_{S,s})^{\nabla}$ at $s$ has dimension
 $r$. Analytic continuation along paths gives rise to a homomorphism
 $\pi_1^{\mathrm{top}}(S,{s}) \rightarrow
 \Aut_{\C}((\mathcal{M}\otimes\mathcal{O}_{S,s})^{\nabla})$ (the
 \emph{monodromy}). The sheaf of germs of solutions
 $\mathcal{M}^{\nabla}$ is locally constant on $S$: its pull-back over
 the universal covering $\til{S}$ of $S$ is constant. Conversely, any
 complex representation $V$ of $\pi_1^{\mathrm{top}}(S,{s})$ of dimension $r$
 gives rise naturally to a vector bundle $\mathcal{M}$ of rank $r$ with
 integrable connection $\nabla$:  $\mathcal{M}= (V \times
 \til{S})/\pi_1^{\mathrm{top}}(S,{s}), \; \nabla(V) = 0$.  This sets up an
 equivalence of categories:
 \begin{equation*}
  \left\{
  \begin{matrix}
   \text{finite dimensional}\\
   \text{representations of $\pi_1^{\mathrm{top}}(S,{s})$}
  \end{matrix}
  \right\}
  \simeq
  \left\{
  \begin{matrix}
   \text{$S$-vector bundles with}\\
   \text{integrable connection}
  \end{matrix}
  \right\}.
 \end{equation*}
\end{para}  
\begin{para}
 Let $(S,s)$ be now a pointed connected $p$-adic manifold. It is still
 true that any $\C_p$-linear representation $V$ of
 $\pi_1^{\mathrm{top}}(S,s)$ of dimension $r$ gives rise to a vector
 bundle $\mathcal{M}=\mathcal{M}_V$ a vector bundle of rank $r$ with
 integrable connection $\nabla =\nabla_V$ (same formula).
 The functor
 \begin{equation*}
  \left\{
  \begin{matrix}
   \text{finite dimensional}\\
   \text{representations of $\pi_1^{\mathrm{top}}(S,s)$}
  \end{matrix}
  \right\}
  \rightarrow
  \left\{
  \begin{matrix}
   \text{$S$-vector bundles with}\\
   \text{integrable connection}
  \end{matrix}
  \right\}
 \end{equation*}
 is still fully faithful, but \emph{no longer surjective}; its essential
 image consists of those connections whose sheaf of solutions is locally
 constant, \ie  becomes constant over $\til{S}$.  In fact, the
 classical ``\emph{Cauchy theorem}'' according to which the solution
 space $(\mathcal{M}\otimes\mathcal{O}_{S,s})^{\nabla}$ at $s$ has
 dimension $r$ holds for every classical point $s$ of $S$ --- which
 corresponds to a point of the associated rigid variety ---, but
 \emph{does not hold} for non-classical points $s$ of the Berkovich
 space $S$ in general.
\end{para}

\begin{para}\label{para-cauchy} 
 Any non-trivial connection over the projective line minus a few points
 gives an example when the $p$-adic analogue of the ``Cauchy theorem''
 does not hold: indeed, in this case, the topological fundamental group
 is trivial.  

 However, a $p$-adic variant of Cauchy's theorem in the neighborhood of a
 non-classical point $s$ may be restored as follows (Dwork's technique of
 generic points): it suffices to extend the scalars from $\C_p$
 to a complete algebraically closed extension of $\mathscr{H}(s)$.  This
 transforms $s$ into a classical point, and neighborhoods of $s$ after
 scalar extensions are ``smaller'' than before.
\end{para}

\begin{para}\label{para-cauchy-monodromy} 
 When ``Cauchy's theorem'' holds at every point of $S$, one can continue
 the local solutions along paths and get the monodromy representation
 just as in the complex situation. 

 This nice category of $p$-adic connections has not yet
 attracted much attention.
\end{para}

\begin{exa}\label{exa-tate-ell-curve}
 Let us consider the case when $S$ is a Tate elliptic curve:
 $S=\C_p^{\times}/q^{\Z}$, with $s =$ its origin. Then
 $\pi_1^{\mathrm{top}}(S,s)=q^{\Z}$, and the connections on $S$ which arise
 from representations of $q^{\Z}$ are those which become trivial over
 $\til{S}=\C_p^{\times}$.
 In this correspondence, the multivalued solutions $\vec y$ (\ie 
 horizontal sections) of such a connection are the solutions of the
 associated linear $q$-difference equation $\vec y(qt)=M(q) \vec y(t)$,
 where $t$ is the standard coordinate on $\C_p^{\times}$ and the
 matrix $M(q)$ is the image of $q$ in the monodromy representation.

 The simplest example is given by the obvious one-dimensional 
 representation $M(q)=q$. It corresponds to
 $\mathcal{M}=\mathcal{O}_S,\;\nabla(1)=\omega_{\mathrm{can}}$ (the canonical
 differential on $S$ induced by $dt/t$); here the $q$-difference
 equation is $t.dy=y.dt$ with obvious solution $y=t$. 

 Let us now look at the representation $M(q)=\sqrt{q}$.  The associated
 $q$-difference equation is $t.dy=\frac{1}{2}y.dt$. There is the 
 obvious solution $\sqrt t$, which leads to
 $\mathcal{M}=\mathcal{O}_S,\;\nabla(1)=\frac{1}{2}\omega_{\mathrm{can}}$. 
 Here we encounter an interesting paradox: $\sqrt{t}$ is not a multivalued
 analytic function on $S$ (\ie it is not an analytic function on the
 universal covering $\C_p^{\times}\;$.  In the complex situation, a
 similar paradox arose in the work of G. Birkhoff in his theory of
 $q$-difference equations, which was pointed out
 and analyzed by M.~van der Put and M.~Singer in the last chapter of
 their book \cite{gtode}.  The solution of the paradox is that the vector
 bundle $\mathcal{M}$ associated to the representation $q \rightarrow
 \sqrt{q}$ of $q^{\Z}$ (or to the $q$-difference equation
 $y(qt)=\sqrt{q}.t$) is in fact a non-trivial vector bundle of rank one,
 and the basic solution is not $\sqrt{t}$, but
 $\frac{\theta(t/\sqrt{q})}{\theta(t)}$, where
 $\theta(t)=\prod_{n>0}(1-q^nt)\prod_{n\leq 0}(1-q^n/t)$.

\end{exa}

\begin{para}\label{para-Simpson} 
 In the previous example, it is easily seen that rank-one vector bundles
 with connection on $S$ which arise from a representation of
 $\pi_1^{\mathrm{top}}(S,s)$ form a space of dimension one, while the
 space of all rank-one vector bundles with connection on $S$ has
 dimension two. 

 We next consider the case of a $p$-adic manifold $S$ which ``is'' a
 geometrically irreducible algebraic ${\C}_p$-curve. Let $\ovl{S}$ be
 its projective completion. It follows from the Van Kampen theorem,
 together with the fact that punctured disks are simply-connected, that
 $\pi_1^{\mathrm{top}}(S,s) \to \pi_1^{\mathrm{top}}(\ovl{S},s)$ is
 an isomorphism. It 
 follows that the vector bundles with connection attached to
 representations of $\pi_1^{\mathrm{top}}(S,s)$ automatically extend
 to $\ovl{S}$. Hence we may assume without loss of generality that $S$ is
 compact.

 Vector bundles with connection on $S$ are algebrizable, and one can
 use C.~Simpson's construction (\cite{morotfgoaspv1}) to define the {\it
 moduli space of connections} of rank $r$ over $S$, denoted by
 $M_{\mathrm{dR}}(S,r)$. On the other hand, we have seen that the topological
 fundamental group $\pi_1^{\mathrm{top}}(S,s)$ is free on $b_1(\Delta)$
generators, being isomorphic to the fundamental group of the dual graph
$\Delta$ of the semistable reduction of $S$. Simpson has also studied
the {\it moduli space of representations of dimension $r$} of such a
group. We denote it by $M_{\mathrm{B}}(S,r)$; in fact, it depends only on the
couple of integers $(b_1(\Delta),r)$.

 Let us assume that $S$ is of genus $g\geq 2$. Simpson shows that 
 $M_{\mathrm{dR}}(S,r)$ is algebraic irreducible of dimension
 $2(r^2(g-1)+1)$.   On the other hand, $M_{\mathrm{B}}(S,r)$ is an algebraic
 irreducible affine variety of dimension $(r^2(b_1(\Delta)-1)+1)$.  We
 note that this dimension is maximal when the Betti number $b_1(\Delta
 )$ takes its maximal value, namely $g$. This corresponds to the case 
 where $S$ is a Mumford curve (a curve with totally degenerate reduction).  

 It turns out that the functor which associates a vector
 bundle with connection to any finite-dimensional representation of the
 topological fundamental group induces an {\rm injective analytic map}
 of moduli spaces $\iota : M_{\mathrm{B}}(S,r)\to M_{\mathrm{dR}}(S,r)$.

[The map is induced by the functor $V \to (\mathcal{M}_V,{\nabla}_V)$,
and its injectivity follows from the faithfulness of this functor. The
difficulty in showing that $\iota$ is analytic is that
$M_{\mathrm{dR}}(S,r)$ is 
a priori a moduli space for algebraic connections, not for all analytic
connections; we shall not pursue here in this direction].

In the complex situation, the corresponding map $\iota$ would always be
an analytic isomorphism (Riemann-Hilbert-Simpson). In the $p$-adic
cases, we see that the connections which satisfy Cauchy's theorem at
every point (classical or not) form a stratum of maximal dimension half
of the dimension of the total moduli space.
\end{para}

\vspace{10pt}\par
We shall leave the closer analysis of this kind of $p$-adic 
differential equations with global monodromy until later chapters, where
we show how they arise in the context of period mappings. In the next 
section, we shall deal with a very different kind of differential 
equations, which play a distinguished role in $p$-adic analysis under 
the name of {\it unit-root F-crystals}.

\newpage
\section{Analytic continuation; algebraic approach.}

\begin{abst}
 An equivalent, but more algebraic, approach to analytic continuation
 consists in interpreting complex analytic multivalued functions as
 limits of algebraic functions.
 Since the concept of topological coverings and that of \'etale coverings
 do not coincide in the $p$-adic setting, the algebraic approach 
 leads in that case to a theory quite different 
 from that of the previous section.  It turns out to be well-suited to
 the function-theoretic study of so-called unit-root $F$-isocrystals, and induces
 us to revisit two fundamental notions due to Dwork: Frobenius structure 
 and overconvergence.
\end{abst}

\subsection{Limits of rational functions.}\label{subsec-limit-rational}

\begin{para}\label{para-complex-approx}
{\itshape Over $\C$.\ }
 Let us briefly survey the complex analytic situation.
 Let $S$ be a connected complex analytic curve, and $U$ an open set
 of $S$ (not necessarily connected) which is holomorphically convex; 
 \index{holomorphically convex@holomorphically convex}
 \ie $S\setminus U$ has no compact connected component.
 Then $\OO(S)$ is dense in $\OO(U)$ for the topology of uniform
 convergence on every compact set.
 We denote this situation by $\OO(U)=\widehat{\OO(S)^U}$
 (\cf \eg \cite[13]{qadsdr}).

 Moreover, if $S$ is the Riemann surface coming from an affine algebraic
 curve $S^{\mathrm{alg}}$, then we also have 
 $\OO(U)=\widehat{\OO(S^{\mathrm{alg}})^U}$, which generalizes the theorem
 of Runge on approximation by rational functions.
 Indeed, considering an embedding $S^{\mathrm{alg}}\hookrightarrow
 (\mathbb{A}^N)^{\mathrm{alg}}$, we may extend analytic functions on $S$ to 
 analytic functions on $(\mathbb{A}^N)^{\mathrm{alg}}$
 (due to the vanishing of the first cohomology group of the
 coherent ideal sheaf defining $S$).
 An analytic functions on $(\mathbb{A}^N)^{\mathrm{alg}}$ can be 
 approximated by polynomials, which we restrict to $S$.
\end{para}

\begin{para}\label{para-ultra-approx}
{\itshape Over the $p$-adics.\ } We have an analogous situation. Let $K$ 
be a complete subfield of $\C_p$, and let $S$ be an analytic curve 
coming from a smooth affine algebraic curve over $K$. We consider a 
closed immersion $S \hookrightarrow {\A}^N$. For any $r > 0$, $S_r := S
\cap \matheur{D}_{{\A}^N}(0,r^+)$ is an affinoid domain in $S$. 
The same argument of polynomial approximation used in 
\ref{para-complex-approx} shows that $\mathcal{O}(S_r)$ is the 
completion of $\mathcal{O}(S^{\mathrm{alg}})$. 

 For any compact $Z \subset S$ (\eg an affinoid domain), we define the
 topological ring $\mathcal{H}(Z)$\index{000H@$\mathcal{H}(\;)$}
 of \emph{analytic elements}
 \index{analytic elements@analytic elements} on $Z$ to
 be the completion of $\Gamma (Z, \mathcal{O}_S) = \liminj_{Z\subset
 U \text{: open}}\Gamma (U, \mathcal{O}_S) $ under the sup-norm.
 Whenever $Z$ is contained in a $S_r$, $\mathcal{H}(Z)$ is also the
 completion of $\mathcal{O}(S_r)_Z$, where the subscript $Z$ denotes the
 localization with respect to the set of elements which do not vanish on
 $Z$ (\cf\ \cite{staagonf} and \cite{rdldaecpecd} for a more precise
version of the Runge theorem).

 The holomorphic convexity condition on $Z$ is that for some $r$, $Z$ is
 the intersection of affinoid neighborhoods defined by inequalities of
 the form $|f_i| \leq 1$ in $S_r$. If this condition is satisfied,
 $\mathcal{H}(Z)$ is in fact the completion of $\mathcal{O}(S_r)$ itself.
\end{para}

\begin{para}
{\itshape The Krasner-Dwork viewpoint.\ }
For $S=\mathbb{A}^1$, and $K=\C_p$, we deduce
$\mathcal{H}(Z)=\widehat{\left(K[z]_Z\right)}$, which is nothing but the 
M.\ Krasner's original definition of analytic elements on $Z$.

Let $\widehat{K(z)}$ be the completion of $K(z)$ by means of Gauss
norm.
One can interpret its elements as the analytic elements on a generic
disk, and then specialize in the complement $Z\subseteq\mathbb{P}^1$ of
a finite union of disks $\matheur{D}(a_i,1^-)$, $(i=1,2,\ldots)$.
This viewpoint leads us to understand analytic continuation
as a \emph{specialization}. 
\end{para}

\subsection{Limits of algebraic functions; complex case.}

Let $S$ be a Riemann surface coming from a complex affine algebraic
curve $S^{\mathrm{alg}}$, smooth and connected.
Let $\til{S}$ be the universal covering of $S$
(here we tacitly fix a base point $s\in S$).
We endow $\OO(\til{S})$ with the topology of uniform convergence on
every compact; note that the group $\pi_1(S)$ acts (continuously) on
$\OO(\til{S})$.

On the other hand, let $\OO(S^{\mathrm{alg}})^{\mathrm{et}}$ be the
maximal unramified integral extension of $\OO(S^{\mathrm{alg}})$.
Elements in $\OO(S^{\mathrm{alg}})^{\mathrm{et}}$ can be regarded as
unramified algebraic functions on $S^{\mathrm{alg}}$; thus we have
$\OO(S^{\mathrm{alg}})^{\mathrm{et}}\subset\OO(\til{S})$.

\begin{pro}\label{pro-gabber}
$\OO(S^{\mathrm{alg}})^{\mathrm{et}}$ is dense in $\OO(\til{S})$.
\end{pro}

For example, when $S^{\mathrm{alg}}=\mathbb{P}^1\setminus\{0,\infty\}$,
the function $\log z\in\OO(\til{S})$ can be written 
as $\displaystyle \lim_{n\rightarrow \infty}n(z^{1/n}-1)$, 
uniformly on every compact
subset of $\til{S}\simeq\mathbb{A}^1$.

\noindent We can use this formula to compute $\log_{\gamma}1=\lim
n(\zeta_n-1)=2i\pi$, where $\gamma$ is the counter-clockwise loop around
$0$ with the base point $1$.
\begin{figure}[h]
 \begin{picture}(70,80)(-35,-40)
  \put(0,-1){\makebox(0,0)[c]{$\bp$}}
  \put(0,-8){\makebox(0,0)[c]{$0$}}
  \put(30,-1){\makebox(0,0)[c]{$\bp$}}
  \put(33,-8){\makebox(0,0)[c]{$1$}}
  \qbezier(0,20)(5,20)(15,15)
  \qbezier(30,0)(25,10)(15,15)
  \qbezier(0,-20)(5,-20)(15,-15)
  \qbezier(30,0)(25,-10)(15,-15)

  \qbezier(-20,16)(-19,18)(-18,20)
  \qbezier(-20,16)(-18,15)(-16,14)

  \qbezier(0,20)(-10,21)(-20,16)
  \qbezier(-30,0)(-30,10)(-20,16)
  \qbezier(0,-20)(-10,-21)(-20,-16)
  \qbezier(-30,0)(-30,-10)(-20,-16)
 \end{picture}
 \caption{}
\end{figure}

\begin{proof}[Proof of \ref{pro-gabber}] (after O. Gabber).
 Consider a countable covering of $\til{S}$
 ($\simeq \mathbb{P}^1(\C),\C$ or $\matheur{D}(0,1^-)$) by
 relatively compact contractible open subsets which are oriented
 manifolds $U_n$ (\eg disks), such that $\ovl{U}_n\subset U_{n+1}$. For
 any $n$, let us consider the set
 \begin{equation*}
 \mathcal{F}_n=\{\gamma \in \pi_1(S,s)\mid
 \gamma\neq 1\text{ and }\gamma\ovl{U}_n \cap \ovl{U}_n \neq \emptyset\}.
 \end{equation*}
 If this set were infinite, we could find a sequence of
 pairwise distinct elements $\gamma_m \in \pi_1(S,s)$ and a sequence of
 points $s_m\in \ovl{U}_n$ such that $\gamma_m s_m\in \ovl{U}_n$. By
 compacity of $\ovl{U}_n$, we might assume that the sequences $s_m$ and
 $\gamma_m s_m$ converge to points $s'$ and $s''$ respectively. Then the
 sequence $\gamma_m s'$ converges to $s''$, which contradicts the
 discreteness of the $\pi_1(S,s)$-orbits in $\til{S}$.  So
 $\mathscr{F}_n$ is finite.

 On the other hand, $S$ is topologically a surface of genus
 $g$ with $t\geq 1$ punctures, so $\pi_1(S)$ is generated by $2g+t$
 generators $\gamma_i$ subject to the relation
 $$[\gamma_1,\gamma_{g+1}]\cdots[\gamma_g,\gamma_{2g}]\gamma_{2g+1}\cdots
 \gamma_{2g+t}=1.$$ Thus $\pi_1(S)$ is free with $2g+t-1$
 generators.  Therefore it is residually finite (\cf\ \cite[p.\ 150, ex.\
 34]{a1}); so we can find a subgroup $\Gamma_n\subseteq\pi_1(S)$ of
 finite index which avoids $\mathscr{F}_n$.

 We immediately see that 
 $\bigl((\Gamma_n\setminus\{1\})\!\cdot\!\ovl{U}_n\bigr)\cap\ovl{U}_n =
 \emptyset$.
 The restriction to $\ovl{U}_n$ of the canonical projection $\til{S}
 \rightarrow S_n := \til{S}/\Gamma_n$ is thus an embedding.
 We identify $\ovl{U}_n$ with its image.

 The part of the exact sequence of cohomologies with compact support
 \begin{equation*}
 \cdots\longrightarrow\mathrm{H}^1(\ovl{U}_n,\Z)\longrightarrow
 \mathrm{H}^2_{\mathrm{c}}(S_n\setminus\ovl{U}_n,\Z)\longrightarrow
 \mathrm{H}^2_{\mathrm{c}}(S_n,\Z)\longrightarrow\cdots
 \end{equation*}
 which can also be written as
 \begin{equation*}
 \ldots \rightarrow
 \mathrm{H}_1(\ovl{U}_n,\Z)=0\rightarrow
 \mathrm{H}_0(S_n\setminus \ovl{U}_n,\Z)\rightarrow
 \mathrm{H}_0(S_n,\Z)=\Z\rightarrow \ldots
 \end{equation*}
 and from which we deduce that $S_n\setminus \ovl{U}_n$ is
 connected.  Because $\partial U_n\subset \ovl{U_{n+1}\setminus
 \ovl{U}_n}$, we have $\ovl{S_n\setminus \ovl{U}_n}=S_n\setminus
 U_n$, hence $S_n\setminus U_n$ is connected.  Since $S_n\setminus U_n$
 is not compact, $U_n$ is holomorphically convex in $S_n$.
 Riemann's existence theorem assures that $S_n$ is the 
 analytification of an affine algebraic curve $S^{\mathrm{alg}}_n$.
 Here we can apply \ref{para-complex-approx} to deduce
 $\OO(U_n)=\widehat{\OO(S^{\mathrm{alg}}_n)^{U_n}}$.
 Note that 
 $\OO(S^{\mathrm{alg}}_n)\subset\OO(S^{\mathrm{alg}})^{\mathrm{et}}$.
 Since every compact set of $\til{S}$ is contained in some $U_n$, we 
 conclude the desired equality 
 $\OO(\til{S})=\widehat{\OO(S^{\mathrm{alg}})^{\mathrm{et}}}$.
\end{proof}

\begin{para}
 It follows from \propref{pro-gabber} that if $\Omega/\Gamma$ is 
 a Schottky partial uniformization of $S$, then
 the intersection $\mathcal{O}(S^{\mathrm{alg}})^{\mathrm{et}}\cap
 \mathcal{O}(\Omega)$ is dense in $\mathcal{O}(\Omega)$.  A variant of this
 statement (with a similar proof) holds, in the $p$-adic case, for a 
 Mumford curve.  We leave it to the reader.
\end{para}

\subsection{Limits of algebraic functions; $p$-adic case.}
\label{subsection-remedy}
\begin{para}
 In the $p$-adic situation, we have seen that there are ``much more''
 \'etale coverings than topological coverings. For instance, the unit disk
 $\matheur{D}=\matheur{D}(0,1^+)$ is simply-connected, but admits many
 non-trivial finite \'etale coverings, \eg the Artin-Schreier covering
 $\matheur{D}\rightarrow \matheur{D}, \;y\mapsto z=y^p-y$.
 
\propref{pro-gabber} suggests to replace, in the $p$-adic case, the ring 
 $\mathcal{O}(\til{S})$, which is often too
 small, by some kind of completion of $\mathcal{O}(S)^{\mathrm{et}}= \varinjlim
 \mathcal{O}(S'),$ where $S'$ runs over the finite \'etale
 connected coverings of $S$.
 
 Let us now assume that $S$ is an affinoid curve with good 
 reduction (hence simply connected). Then there is a
 Gauss $p$-adic norm on $\mathcal{O}(S)$, which extends uniquely to 
 $\mathcal{O}(S)^{\mathrm{et}}$.
 The completion $\widehat{\mathcal{O}(S)^{\mathrm{et}}}$ is
 however pathological in several senses: for instance, 
 it is difficult to give a function-theoretic 
 meaning to its elements, and the continuous derivations of
 $\mathcal{O}(S)$ extends to $\mathcal{O}(S)^{\mathrm{et}}$ but not to
 $\widehat{\mathcal{O}(S)^{\mathrm{et}}}$ in a natural way.
\end{para}

\begin{para}\label{para:r}
 To bypass such problems, we are going to look for some convenient {\it subring} of
 $\widehat{\OO(S)^{\mathrm{et}}}$.
 For simplicity, let us limit ourselves to the following situation:
 \begin{itemize}
  \item $K = {\widehat{\Q}}^{\mathrm{ur}}_p \subset \C_p$
	is the completion of
	the maximal unramified algebraic extension of ${\Q}_p$,
  \item $\mathfrak{v} = {\widehat{\Z}}^{\mathrm{ur}}_p$ is its ring of
	integers (the Witt ring of ${\ovl{\F}}_p$),
  \item $S_0$ is the affine line over ${\ovl{\F}}_p$, deprived from
	finitely many points ${\ovl{\zeta}}_1\ldots,{\ovl{\zeta}}_{\nu}$;
	we choose a point $s_0$ on $S_0$,
  \item $\mathcal{R}$ is the $p$-adic completion of
	$\mathfrak{v}[z, \frac{1}{(z-\zeta_1)\ldots (z-\zeta_{\nu})}]$,
	where $\zeta_1\ldots \zeta_{\nu}$ are liftings of
	${\ovl{\zeta}}_1\ldots,{\ovl{\zeta}}_{\nu}$ in $\mathfrak{v}$; its
	residue ring is $\mathcal{O}(S_0)$,
  \item $S = \matheur{M}\bigl(\mathcal{R}[\frac{1}{p}]\bigr) =
	\matheur{D}(0,1^+) \setminus \bigcup 
	\matheur{D}(\zeta_i,1^-)$, the associated affinoid domain over $K$.
	There is a natural specialization map $\sp: S \rightarrow S_0$ from
	characteristic $0$ to characteristic $p$.
 \end{itemize}

 Let us consider the integral closure ${\mathcal{R}^{\mathrm{et}}}$ of
 $\mathcal{R}$ in ${\mathcal{O}(S)^{\mathrm{et}}}$ and its $p$-adic
 completion $\widehat{\mathcal{R}^{\mathrm{et}}} \subset
 \widehat{\OO(S)^{\mathrm{et}}}$.
 It is thus endowed with (a natural extension of) the $p$-adic
 valuation, and its residue ring is
 $\mathcal{O}(S_0)^{\mathrm{et}}$. Moreover, there is a natural
 structure of ``differential ring'' (better: a connection) on
 $\widehat{\mathcal{R}^{\mathrm{et}}}$ with $\mathfrak{v}$ as ring of
 constants. The ``remarkable equivalence of categories'' of
 A. Grothendieck \cite[18.1.2]{ega44} allows to identify
 $\Aut_{\mathrm{cont}}(\widehat{\mathcal{R}^{\mathrm{et}}}/\mathcal{R})$
 with the algebraic fundamental group $\pi_1^{\mathrm{alg}}(S_0,s_0)$.
\end{para}

\begin{para}\label{para-analytic-cont}
 {\itshape ``Analytic continuation'' as specialization.\ } Let
 $\matheur{D}(s_0,1^-) = \sp^{-1}\lbrace s_0\rbrace$ be the residue
 class of $s_0$ in $S$.  The morphism $\mathcal{R}\rightarrow
 \mathcal{O}(\matheur{D}(s_0,1^-))$ extends non-canonically to a
 continuous morphism $\widehat{\mathcal{R}^{\mathrm{et}}}
 [\frac{1}{p}]\rightarrow \mathcal{O}(\matheur{D}(s_0,1^-))$, determined
 \emph{only up to the action of } $\pi_1^{\mathrm{alg}}(S_0,s_0)$. This
 allows to interpret elements of
 $\widehat{\mathcal{R}^{\mathrm{et}}}[\frac{1}{p}]$ as certain
 multivalued locally analytic functions (in the ``wobbly'' sense). 

The main difference here with Krasner's analytic continuation
 (\cf~\subsecref{subsec-limit-rational}) is the ambiguity coming from
 $\pi_1^{\mathrm{alg}}(S_0,s_0)$.  This provides a kind of multivalued
 analytic continuation, which may be interpreted as analytic
 continuation in characteristic $0$ along an ``\'etale path'' in
 characteristic $p$ (figure \ref{pic4}). 
 
\begin{center}
 \begin{figure}[h]
  \begin{picture}(200,100)(0,0)
   \put(0,100){\makebox(0,0)[l]{$\matheur{D}(s_0,1^{-})$}}
   \put(0,0){\includegraphics{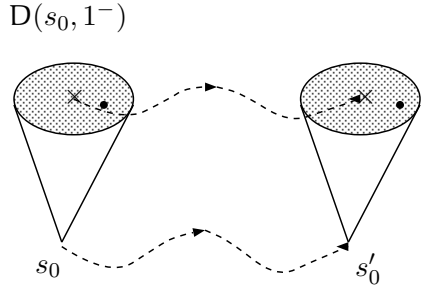}}
   \put(15,5){\makebox(0,0)[c]{$s_0$}}
   \put(136,5){\makebox(0,0)[c]{$s_0'$}}
   \put(25,70){\makebox(0,0)[c]{$\times$}}
   \put(135,70){\makebox(0,0)[c]{$\times$}}
   \put(36,66){\makebox(0,0)[c]{$\bbp$}}
   \put(148,66){\makebox(0,0)[c]{$\bbp$}}
  \end{picture}
  \caption{analytic continuation along an underground path}
  \label{pic4}
 \end{figure}
\end{center}

 This is very close to the Robba-Christol theory of ``algebraic
 elements'' \cite{feea}, in which the main player is the \emph{complete}
 subalgebra of $H^{\infty}(\matheur{D}_{\C_p}(0,1^-))$ (the algebra of bounded
 analytic functions) generated the algebraic functions analytic on
 $\matheur{D}_{\C_p}(0,1^-)$. Underlying this ``algebraic'' analytic
 continuation, 
 there is a combinatorics of \emph{automata} (whereas a combinatorics of
 graphs underlies the ``topological'' analytic continuation, as we have
 seen).\footnote{for these and further aspects of the Krasner, resp.
 Christol-Robba, analytic continuation, we refer to \cite{paeap} and to the mimeographed
 notes of the numerous conferences of the Groupe d'\'etude d'analyse
 ultram\'etrique, Paris, devoted to this subject (1973-1980).}
\end{para}

\begin{rem}
 The complex formula $\log z = \lim_n n\bigl(z^\frac{1}{n}-1\bigr)$ has the
 following $p$-adic analogue: $\log z = \lim_n p^{-n}(z^{p^n}-1)$ for $z
 \in \matheur{D}(1,1^-)$. However, the convergence is not uniform on
 $\matheur{D}(1,1^-)$, and the analytic function $\log z$ is not bounded
 on that disk.  Nevertheless $\log : \matheur{D}(1,1^-) \rightarrow
 {\A}^1$ defines an infinite Galois \'etale covering of the $p$-adic
 affine line (not at all a topological covering), with Galois group
 $\mu_{p^{\infty}}$, the $p$-primary torsion in $\C_p$.
\end{rem}

\subsection{$\widehat{\mathcal{R}^{\mathrm{et}}}$ and unit-root $F$-crystals.}
\label{sub-r-and-unit-root}
\begin{para}
 Let us assume for technical simplicity that $p\neq 2$. Let $\mathcal{M}$ be
 a free $\mathcal{R}$-module of finite rank $\mu$, endowed with a connection
 $\nabla: \mathcal{M} \rightarrow
 \Omega_\mathcal{R}\otimes_\mathcal{R}\mathcal{M}$, where
 $\Omega_\mathcal{R}$ is the rank-one module of continuous differentials
 (relative to $\mathfrak{v}$).

 In general, solutions make sense only very locally: typically, analytic
 solutions exist in disks of radius $p^\frac{-1}{p-1}$ (the radius of
 $p$-adic convergence of the exponential function, which is the basic
 example).  A very important criterion, due to B.~Dwork, for the
 convergence of analytic solutions in any open disk of radius $1$, is
 the existence of a so-called Frobenius structure.
\end{para}

\begin{para}\label{para:f-cris}
 {\itshape $F$-crystals.}  Let $\sigma$ be the Frobenius automorphism of
 $\mathfrak{v}$ lifting the $p$th-power map in $\ovl{\F}_p$. There
 are many $\sigma$-linear endomorphisms $\phi$ of $\mathcal{R}$ which reduce
 to the $p$th-power map of $\mathcal{O}(S_0)$ in characteristic
 $p$. For instance, we can take $\phi(z) = z^p$, so that
 $\phi\bigl((z-\zeta_1)\ldots (z-\zeta_{\nu})\bigr) =
 (z^p-\sigma(\zeta_1))\ldots
 (z^p-\sigma(\zeta_{\nu})) \equiv
 \bigl((z-\zeta_1)\ldots (z-\zeta_{\nu})\bigr)^p
 \pmod{p}$. One can show that any lifting $\phi$ extends uniquely to a 
 $\pi_1^{\mathrm{alg}}(S_0,s_0)$-equivariant endomorphism of
 $\widehat{\mathcal{R}^{\mathrm{et}}}$.

 A Frobenius structure on $(\mathcal{M},\nabla)$ is a rule $F$ which
 associates to every lifting $\phi$ a homomorphism $F(\phi):
 \phi^{\ast}\mathcal{M} \rightarrow \mathcal{M}$ (that is to say, a
 $\phi$-linear endomorphism of $\mathcal{M}$) such that:

\begin{enumerate}
 \item $F(\phi)\otimes \Q$ is an \emph{isomorphism},
 \item $F(\phi)$ is \emph{horizontal}, \ie compatible with
       the connections
       $\phi^{\ast}\nabla$ and $\nabla$ respectively,
 \item for any two liftings $\phi$ and $\phi^{\prime}$, the
       homomorphisms $F(\phi)$ and
       $F(\phi^{\prime})$ are related by $F(\phi^{\prime}) = F(\phi)\; \circ
       \;\chi(\phi^{\prime},\phi),$ where
       $\chi(\phi^{\prime},\phi): {\phi^{\prime}}^{\ast}\mathcal{M} \rightarrow
       \phi^{\ast}\mathcal{M} $ is
       the ``Taylor isomorphism'' given by the formulas
       \begin{align*}
	\chi(\phi^{\prime},\phi)({\phi^{\prime}}^{\ast}m) & =
	 \sum_{n\geq 0}{\phi}^{\ast}
	\biggl(\nabla \Bigl(\frac{d}{dz}\Bigr)^n(m)\biggr)
	 ({\phi^{\prime}}(z)-{\phi}(z))^n/n!\\
	& = \sum_{n\geq 0}{\phi}^{\ast} 
	 \biggl(\nabla \Bigl(z\frac{d}{dz}\Bigr)^n(m)\biggr)
	 \biggl(\log\frac{{\phi^{\prime}}(z)}{{\phi}(z)}\biggr)^n\big/n!.
       \end{align*}\label{item:3}
\end{enumerate} 

In virtue of (\ref{item:3}), a Frobenius structure amounts to the data of
the semi-linear horizontal endomorphism $F(\phi)$, for the standard
$\phi(z) = z^p$.

\noindent
The triple $(\mathcal{M},\nabla, F)$ is
called an $F$-{\it crystal}, \cf \cite{tdd}\footnote{the endomorphism
$F_{S_0}$ of $S_0$ given by the $p$th-power map, as well as $\sigma$, $\phi$
and $F$ are all called ``Frobenius'' in the usual jargon, without
causing too much confusion, it seems... }.  The very existence of a
horizontal isomorphism $\phi^{\ast}\mathcal{M}\otimes {\Q} \simeq
\mathcal{M} \otimes {\Q}$ implies that the radius of convergence
${\rho}$ of any local solution of $(\mathcal{M},\nabla)$ satisfies
$\min(1,{\rho})=\min(1,{\rho}^p)$, that is to say: ${\rho}\geq 1$ (since
${\rho}\neq 0$).
\end{para}

\begin{para}\label{para-unit-root}
 {\itshape Unit-root $F$-crystals.}  This is the case where $F(\phi)$
 (not only $F(\phi)\otimes {\Q}$) is an isomorphism (for one, or
 equivalently, for all $\phi$).  The name comes from their first
 appearance in Dwork's computation of the $p$-adic units among the
 reciprocal zeroes of the zeta-function of hypersurfaces in
 characteristic $p$.

 N. Katz has constructed a functor:
 \begin{equation*}
  \left(\begin{array}{l}
   \text{continuous ${\Z}_p$-representations}\\
   \text{of $\pi_1^{\mathrm{alg}}(S_0,s_0)$}
  \end{array}\right)
   \longrightarrow \text{(unit-root $F$-crystals)}
 \end{equation*}
 which associates to any free ${\Z}_p$-module $V$ of rank $r$
 with a continuous action of
 $\pi_1^{\mathrm{alg}}(S_0,s_0)$ a unit-root
 $F$-crystal $(\mathcal{U}_V,\nabla_V, F_V)$ over $\mathcal{R}$ of rank
 $r$.
 Let us recall its definition at the level of objects.
 For any $m\in \N$, we set $S_m=\Spec\;\mathcal{R}/p^m\mathcal{R}$.
 Starting from a representation $\rho$ of $\pi_1^{\mathrm{alg}}(S_0,s_0)$, let
 $G_n$ denote the image of $\rho$ in $GL(V/p^nV)$.  The homomorphism
 $\pi_1^{\mathrm{alg}}(S_0,s_0)\rightarrow G_n$ corresponds to an \'etale 
 covering $S_{n,0}\rightarrow S_0$, which has a
 unique \'etale lifting $\pi_{n,m}:\;S_{n,m}\rightarrow S_m$. The action 
 of $G_n$ on $S_{n,0}$ extends uniquely to
 $S_{n,m}$; on the other hand, the action of $\phi$ on $S_m$ extends 
 uniquely to $S_{n,m}$, and the $\phi$- and
 $G_n$-actions commute. The opposite action makes
 $\mathcal{O}_{S_{n,m}}$ into a right $G_n$-module.
 If we set $\mathcal{U}_n=
 \pi_{n,n_\ast}\mathcal{O}_{S_{n,n}}
 \otimes_{(\mathfrak{v}/p^n\mathfrak{v})[G_n]} V$, 
 we then have a compatible system of isomorphisms
 \begin{equation*}
  \Phi_n=\phi\otimes \id:\;\phi^{\ast}\mathcal{U}_n\rightarrow\mathcal{U}_n
 \end{equation*} 
 and a compatible system of connections
 \begin{equation*}
  \nabla_n=d\otimes \id:\;\mathcal{U}_n\rightarrow
 \pi_{n,n_\ast}\Omega_{S_{n,n}}\otimes_{\mathcal{R}/p^m\mathcal{R}}
 \mathcal{U}_n \simeq \Omega_\mathcal{R}\otimes_\mathcal{R} {\mathcal{U}}_n
 \end{equation*} 
 The unit-root crystal attached to $\rho$ is given 
 by $\mathcal{U}=\varprojlim_n \mathcal{U}_n,
 \nabla=\varprojlim_n \nabla_n,\Phi=\varprojlim_n \Phi_n$.

 The next statement summarizes results of Katz and R.~Crew \cite{fapr}.
\end{para}

\begin{pro}
 \begin{enumerate}
  \item The Katz functor is an equivalence of categories;
  \item any unit-root $F$-crystal $(\mathcal{U},\nabla, F)$ is
	solvable in $\widehat{\mathcal{R}^{\mathrm{et}}}$, \ie :
	$\mathcal{U}\otimes_\mathcal{R}
	\widehat{\mathcal{R}^{\mathrm{et}}}\simeq
	(\mathcal{U}\otimes_\mathcal{R}
	\widehat{\mathcal{R}^{\mathrm{et}}})^{\nabla}\otimes_{\mathfrak{v}}
	\widehat{\mathcal{R}^{\mathrm{et}}}$;
  \item the rule  $\;(\mathcal{U},\nabla, F) \mapsto
	(\mathcal{U}\otimes_\mathcal{R}
	\widehat{\mathcal{R}^{\mathrm{et}}})^{\nabla =0, F(\phi)=\id}$
	(with Galois action coming from that on
	$\widehat{\mathcal{R}^{\mathrm{et}}}$) provides an inverse of
	the Katz functor. 
  \end{enumerate}
\end{pro}
More precisely, we have the action of $\phi$ and $d/dz$ on 
$\widehat{\mathcal{R}^{\mathrm{et}}}$ commute with the action of
$\pi_1^{\mathrm{alg}}(S_0,s_0)$, and we have a canonical isomorphism
\begin{equation*}
 \mathcal{U}\otimes_\mathcal{R}\widehat{\mathcal{R}^{\mathrm{et}}}
  \simeq V\otimes_{\Z_p}\widehat{\mathcal{R}^{\mathrm{et}}}
\end{equation*}
compatible with $\phi, d/dz$ and $\pi_1^{\mathrm{alg}}(S_0,s_0)$
(diagonal action of $\phi, d/dz$ on the left hand side, diagonal action
of $\pi_1^{\mathrm{alg}}(S_0,s_0)$ on the right hand side), which allows
to reconstruct the representation $V$ from the unit-root $F$-crystal and
conversely. In fact, the connection as well as $V$ can be reconstructed
from $(\mathcal{U}, \Phi)$ alone.

This proposition may be compared with \ref{para-cauchy-monodromy}, though it
applies to a quite different type of $p$-adic connections. In
\ref{para-cauchy-monodromy} (for dimension 1), the analytic curve $S$ had
typically bad reduction and the main player was the \emph{discrete
fundamental group} $\pi_1(\Delta)$ of the dual graph of the semistable
reduction, together with the $\pi_1(\Delta)$-module
$\mathcal{O}(\til{S})$. Here the curve $S$ has good reduction
$S_0$ and the main player is the \emph{compact fundamental group}
$\pi_1^{\mathrm{alg}}(S_0)$, together with the
$\pi_1^{\mathrm{alg}}(S_0)$-module $\widehat{\mathcal{R}^{\mathrm{et}}}$.

Before presenting one of Dwork's classical unit-root $F$-crystals, let
us just mention that the above theory extends with little change to the
case when the base ring $\mathfrak{v}$ is a finite ramified extension of
${\widehat{\Z}}^{\mathrm{ur}}_p$: one has to fix an extension of $\sigma$ to
$\mathfrak{v}$, to replace ${\Z}_p$ by $\mathfrak{v}^{\sigma}$, and
to impose some mild nilpotence constraint on $\nabla$ if the
ramification index is $\geq p-1$ (also, it is customary to extend the
definition of Frobenius structure on replacing $\phi$ by some power).

\begin{exa}\label{exa-dwork-exp}
 {\itshape Dwork's exponential.} We denote by
 $\pi$ a fixed $(p-1)$th
 root of $-p$, and set $\mathfrak{v}={\widehat{\Z}}^{\mathrm{ur}}_p[\pi]$, with
 $\sigma(\pi) = \pi$. Let us
 consider the differential equation 
 \begin{equation*}\label{eq:ast}
  f^{\prime}(z)=-\pi f(z) \tag{$\ast$}
 \end{equation*}
 over $S=\matheur{D}(0,1^+)$, which has the analytic solution $f_a =
 e^{-\pi(z-a)}$ in any residue class $\matheur{D}(a,1^-)\subset
 \matheur{D}(0,1^+)$.  The change of variable $z\mapsto z^p$ leads to
 the differential equation

 \begin{equation*}\label{eq:astp}
   g'(z)=-p\pi z^{p-1}g(z) = {\pi}^p z^{p-1}g(z) \tag{$\ast_p$}
 \end{equation*} 

 \noindent
 Dwork's exponential function is 
 \begin{equation*}
  E_{\pi}(z) =
   e^{\pi(z-z^p)} \in {\Z}_p[\pi][[z]].
 \end{equation*}
 This is an invertible element of $\mathcal{R}$, the $\pi$-adic completion of
 $\mathfrak{v}[z]$. This function provides the unit-root Frobenius structure
 which relates (\ref{eq:ast}) and (\ref{eq:astp}).

 Let us reformulate this in the setting of \ref{para-unit-root}: the
 relevant $F$-crystal is $(\mathcal{R}, \nabla, F)$, with $\nabla(1) =
 \pi dz$, $F(\phi)(1) = E_{\pi}(z)$ for the standard $\phi$.  This is
 the unit-root $F$-crystal attached to the finite character of
 $\pi_1^{\mathrm{alg}}({\A}^1_{{\ovl{\F}}_p},0)$ corresponding to the
 Artin-Schreier covering $z=y^p-y.$ The solution $f_0$ of $(\ast)$
 belongs to $\widehat{\mathcal{R}^{\mathrm{et}}}$: indeed, $e^{-\pi z}=
 E_{\pi}(y) \in {\widehat{\mathfrak{v}[y]}} \subset
 \widehat{\mathcal{R}^{\mathrm{et}}}$;
 explicitly, $e^{-\pi z}=\lim_n(1-\pi p^n z)^{p^{-n}}$.

 Dwork has computed the value of his exponential for any $(p-1)$th
 root of unity $\zeta_{p-1}$ (\eg $\zeta_{p-1} =1$): this is

 \begin{equation*}
  E_{\pi}(\zeta_{p-1})=\zeta_p^{\zeta_{p-1}} 
 \end{equation*}

 where $\zeta_p$ is the unique $p$th root of unity $\equiv
 1+\pi\pmod{\pi}$ (\cf \cite[chap. 14]{cf1a2}).  In the spirit of
 \ref{para-analytic-cont}, this may be understood as follows:

 $\zeta_p^{\zeta_{p-1}}=$ \emph{the value at $(z=0)$ of the analytic
 continuation of $f_0$ along the ``wild loop''} corresponding to the
 path from $(y=0)$ to $(y=\zeta_{p-1})$ on the Artin-Schreier covering
 $z=y^p-y$ in characteristic $p$ (figure \ref{pic2}).  If one
 changes $\zeta_{p-1}$, this value $\zeta_p^{\zeta_{p-1}}$ is multiplied
 by some $p$th root of unity.

\begin{figure}[h]
 \begin{picture}(200,130)(0,0)
  \put(0,0){\includegraphics[height=4cm,clip]{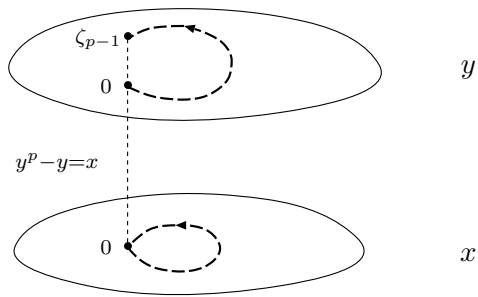}}
  \put(40,102){\makebox(0,0)[c]{$\scriptstyle{\zeta_{p-1}}$}}
  \put(48,99.5){$\bbp$}
  \put(43,84){\makebox(0,0)[c]{$\scriptstyle{0}$}}
  \put(48,81){$\bbp$}
  \put(40,55){\makebox(0,0)[r]{$\scriptstyle{y^p-y=x}$}}
  \put(43,23){\makebox(0,0)[c]{$\scriptstyle{0}$}}
  \put(48,20){$\bbp$}
  \put(180,90){\makebox(0,0)[c]{$y$}}
  \put(180,20){\makebox(0,0)[c]{$x$}}
 \end{picture}
 \caption{a wild underground loop}
 \label{pic2}
\end{figure}

In our special case, the discussion of \ref{para-analytic-cont} tells us
that $e^{-\pi z}$ \emph{admits an extension to any disk
$\matheur{D}(a,1^-) \subset \matheur{D}(0,1^+)$, analytic in that disk,
and well-defined up to multiplication by some $p$th root of unity}. It
is easy to find a formula for this extension: it must be proportional to
$e^{\pi (a-z) }$, and by evaluation at $a$, we find that it is $e^{\pi
(a-z)}E_{\pi}(b)$, where $b$ is any solution of the equation $b^p-b=a$.
Its $p$th power is, as expected, $e^{-p\pi z}$ itself.

Therefore, this multivalued exponential $e^{-\pi z}$ on $\matheur{D}(0,1^+)$
appears as a multivalued section of the logarithmic \'etale covering
$\frac{-1}{\pi}\log : \matheur{D}\bigl(1,p^{-\frac{1}{p-1}+}\bigr)
\rightarrow \matheur{D}(0,1^+)$.  This is just the \emph{opposite} way
from the complex situation, where the logarithm is a multivalued
section of the \'etale covering of $\C^{\times}$ given by the exponential.
\end{exa}

\subsection{$p$-adic chiaroscuro: overconvergence.}
\label{sub-chiaroscuro}
\begin{para}
 In fact, Dwork's exponential $E_{\pi}$ is more
 than just an element of the $\pi$-adic completion of $\mathfrak{v}[z]$: it
 is \emph{overconvergent}, \ie extends to an analytic function on a disk
 of radius $\rho > 1$ ($\rho = p^{(p-1)/p^2}$).

 Overconvergence is a fundamental notion in $p$-adic analysis: dimming
 the contours of an affinoid turns out to be the key to the finiteness
 properties of $p$-adic cohomology, as was recognized by Dwork, and
 subsequently developed by P. Monsky, G.  Washnitzer, P. Berthelot
 (rigid cohomology, $\mathcal{D}^{\dagger}$-modules)...
  \begin{figure}[h]
   \begin{picture}(80,80)(0,0)
    \put(0,0){\includegraphics{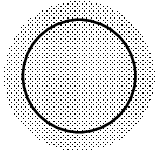}}
    \put(35,35){$\bbp$}
   \end{picture}
   \caption{}
  \end{figure}
\end{para}

\begin{para}
 Let again $S$ be $\matheur{D}(0,1^+)
 \setminus \bigcup \matheur{D}(\zeta_i,1^-)$.  For any
 $\rho >1$, let us consider the bigger affinoids
 $S_{\rho}=\matheur{D}(0,{\rho}^+)
 \setminus \bigcup \matheur{D}(\zeta_i,\frac{1}{\rho}^-)$.  The ring of
 \emph{overconvergent analytic functions} on $S$ is
 
 \begin{equation*}
  \mathcal{H}^{\dagger}(S)= \bigcup_{\rho >1}\mathcal{O}(S_{\rho}).
 \end{equation*}
 Its relevance to the algebraic viewpoint on analytic continuation comes
 from the following result \cite{utdppdfa}:
\end{para}
\begin{pro}
  Let $f \in \mathcal{O}(S)$ satisfy a \emph{monic} polynomial equation
 with coefficients in $\mathcal{H}^{\dagger}(S)$. Then $f \in
 \mathcal{H}^{\dagger}(S)$.
\end{pro}
\begin{para}
 A $F$-crystal $(\mathcal{M},\nabla, F)$ is {\it overconvergent} if
 $(\mathcal{M},\nabla)$ as well as $F(\phi)$ extends over some
 $S_{\rho}$.  This is the case in example \examref{exa-dwork-exp}.  For
 unit-root $F$-crystals, Crew has given the following characterization:

\begin{pro}
 The unit-root $F$-crystal $(\mathcal{U}_V,\nabla_V, F_V )$
 attached to a $p$-adic representation $V$ is overconvergent if and only
 if the images in $GL(V)$ of the inertia groups at the missing points
 ${\ovl{\zeta}}_1\ldots,{\ovl{\zeta}}_{\nu},\infty$ are finite.
\end{pro}

 \end{para}

\subsection{(Overconvergence and Frobenius) Dwork's derivation of
the $p$-adic Gamma function and exponential sums.}
\label{subsection-deriv-of-gamma}

\begin{para}
 We come back to the situation of
 \examref{exa-dwork-exp}.  For any $\alpha \in {\Z}_p$, let us consider
 $M_{\alpha}^{\dagger}:= z^{\alpha}e^{\pi z}\mathcal{H}^{\dagger}(S)$
 endowed with the derivation $z\frac{d}{dz}$. A simple computation shows
 that
 \begin{equation*}
  z^{\alpha}e^{\pi z}z^{k+1}\equiv
   -\frac{\alpha + k}{\pi}z^{\alpha}e^{\pi z}z^k 
 \end{equation*}
 in $M_{\alpha}^{\dagger}/z\frac{d}{dz}M_{\alpha}^{\dagger}$, from which
 one deduces that this cokernel has dimension $1$ over $K$ and is
 generated by $[z^{\alpha}e^{\pi z}]$.

 Following Dwork, let us introduce the operator $\psi$ defined by
 \begin{equation*}
  \psi\Bigl(\sum a_nz^n\Bigr)=\sum a_{pn}z^n.
 \end{equation*}
 This is a left inverse of the Frobenius operator induced by
 $\phi: z \mapsto z^p$.  It acts on $\mathcal{H}^{\dagger}(S)$ and
 commutes with $z\frac{d}{dz}$ up to multiplication by
 $p$: $\psi \circ z\frac{d}{dz}= pz\frac{d}{dz} \circ \psi$.
 The operator $\psi$ can also be applied to any element of
 $M_{\alpha}^{\dagger}$ for $\alpha \in \N$: one finds that
 \begin{equation*}
  \psi(z^{\alpha}e^{\pi z}f)=z^{\beta}e^{\pi z}\psi\bigl(z^{\alpha -
   p\beta}E_{\pi}(z)f\bigr) 
 \end{equation*}
 where $\beta$ is the {\it successor} of $\alpha \in {\Z}_p$, \ie  the
 unique $p$-adic integer $\beta$ such that $p\beta - \alpha \in {\Z}\cap
 [0,p[$ (note that since $\alpha - p\beta >-p\;$, the terms containing a
 negative power of $z$ disappear when applying $\psi$). This formula
 makes sense for any $\alpha \in {\Z}_p$, and we can see that $\psi$
 applies $M_{\alpha}^{\dagger}$ into $M_{\beta}^{\dagger}$, and commutes
 with $z\frac{d}{dz}$. Let us write the induced map of one-dimensional
 cokernels $[{\psi}]:M_{\alpha}^{\dagger} /
 z\frac{d}{dz}M_{\alpha}^{\dagger}\rightarrow
 M_{\beta}^{\dagger}/z\frac{d}{dz}M_{\beta}^{\dagger}$ in the form
 \begin{equation*}
  [{\psi}]([z^{\alpha}e^{\pi z}]) = 
   {\pi}^{p\beta  - \alpha}\Gamma_p(\alpha)[z^{\beta}e^{\pi z}].
 \end{equation*}
\end{para}

\begin{para}
 This function $\Gamma_p$ is then nothing but Morita's $p$-adic
 Gamma function, characterized by its continuity and the functional equation
 \begin{equation*}
  \Gamma_p(0)=1, \;\;\Gamma_p(\alpha +1)/\Gamma_p(\alpha) =
   \begin{cases}
    -\alpha & \text{ if $\alpha$ is a unit,}\\
    -1 & \text{ if $|\alpha|_p < 1$.}
   \end{cases}
 \end{equation*}
 Let us check this functional equation:
 \begin{itemize}
  \item for $\alpha = 0$, one has
	$[z^{k}e^{\pi z}]=0$ if $k>0$; hence
	$\Gamma_p(0)[e^{\pi z}]=[e^{\pi z}][\psi(E_{\pi}(z))]$ $=[e^{\pi z}]$;
  \item if $|\alpha|_p=1$,
	$\beta$ is also the successor of $\alpha +1$, and
	$\Gamma_p(\alpha +1)[z^{\beta}e^{\pi z}] = {\pi}^{\alpha - p\beta
	+1}[{\psi}]([z^{\alpha +1}e^{\pi z}])= {\pi}^{\alpha - p\beta
	+1}[{\psi}](-\frac{\alpha}{\pi}[z^{\alpha}e^{\pi
	z}])=-a\Gamma_p(\alpha)[z^{\beta}e^{\pi z}]$;
  \item finally, if $|\alpha|_p<1$, one has
	$\alpha = p\beta$, and $\beta +1$ is the successor of $\alpha
	+1$; one gets $\Gamma_p(\alpha
	+1)[z^{\beta +1}e^{\pi z}] = {\pi}^{-p+1}[{\psi}]([z^{\alpha +1}
	e^{\pi z}]) = 
	-\frac{1}{p}[{\psi}](-\frac{\alpha}{\pi}[z^{\alpha} e^{\pi z}])
	$ $=
	\frac{\beta}{\pi}\Gamma_p(\alpha)[z^{\beta}e^{\pi z}] =
	- \Gamma_p(\alpha)[z^{\beta +1}e^{\pi z}]$.
 \end{itemize}
\end{para} 

\begin{para}\label{para-cont-of-gamma}
 Let us now check the continuity of $\Gamma_p$ on $\Z_p$, or
 better, its analyticity on any disk
 $\matheur{D}(-k,\vert p\vert^+),\; k=0,1,\ldots,p-1$.
 Let us write $E_{\pi}(z)=e^{\pi(z-z^p)}=\sum_0^{\infty} e_n z^n$.

 For any $\alpha \in \Z_p\cap \matheur{D}(-k,1^-)$, one has
 \begin{align*}
  & \;\psi(z^{\alpha} e^{\pi z})\equiv \pi^k\Gamma_p(\alpha)z^{\beta} e^{\pi z}
  = z^{\beta} e^{\pi z}\psi(z^{-k}E_\pi(z)) \\
  = \;& z^{\beta} e^{\pi z}\sum_{n=0}^{\infty} e_{pn+k} z^n
  \equiv z^{\beta} e^{\pi z}\sum_{n=0}^{\infty} e_{pn+k}(-\pi)^{-n} (\beta)_n,
 \end{align*} 
 where $(\beta)_n=(\frac{\alpha + k}{p})_n$ is the Pochhammer symbol.
 Easy estimates now show that $\Gamma_p(\alpha)=\sum_0^{\infty}
 e_{pn+k}(-\pi)^{-n-k} (\frac{\alpha + k}{p})_n$ is analytic on 
 $\matheur{D}(-k,\vert p\vert^+)$.
\end{para}
\begin{para}\label{para-gross-koblitz}
 {\itshape Gross-Koblitz formula:}
 for $k=0,1,\ldots,p-2$, one has
 \begin{equation*}
  \Gamma_p\biggl(\frac{k}{p-1}\biggr)=-\pi^{-k}\sum_{\zeta_{p-1}}\;
   \zeta_{p-1}^{-k}\zeta_p^{\zeta_{p-1}}\;\;\in \Q_p[\pi]\cap
   {\ovl{\Q}}
 \end{equation*}
 where $\zeta_{p-1}$ runs over the $(p-1)$th roots of unity,
 and $\zeta_p^{\zeta_{p-1}}$ is as before the unique $p$th root of 
 unity $\equiv 1+\zeta_{p-1}\pi \pmod{\pi}$.

 \medskip
 Let us sketch the proof.  We choose $\alpha=\frac{k}{p-1}$, so that
 $\alpha=\beta$ and $\psi$ acts on the ind-Banach-space
 $M_\alpha^{\dagger}$ via the formula:
 $\psi(z^{-k}E_{\pi}(z)f) =
 (z^{\alpha}e^{\pi z})^{-1}\psi(z^{\alpha}e^{\pi z}f)$.
 Coming back to the definition of $\psi$, one observes that for any
 $g\in \mathcal{H}^{\dagger}(S)$, 
 \begin{equation*}
  \psi(z^{-k}g)(z)= \frac{1}{p}\sum_{t\in \phi^{-1}(z)}\;t^{-k}g(t).
 \end{equation*}
 Using the fact that the domain of analyticity of this function is 
 bigger than the domain of analyticity of $g$ itself, one shows that
 $\psi$ is a nuclear operator of $M_\alpha^{\dagger}$. In particular, it
 has a trace, which is the trace of the ``composition operator''
 $\Psi: g\in \mathcal{H}^{\dagger}(S)\mapsto \frac{1}{p}
 \sum_{t\in \phi^{-1}(z)}\;t^{-k}E_\pi(t)g(t)$.
 The computation of this trace is done by approximating $E_\pi$ by
 polynomials and studying the resulting action on the subspace of
 polynomials.  One finds $\Tr\Psi = \frac{1}{p-1}\sum_{\zeta_{p-1}}
 \zeta_{p-1}^{-k}E_\pi(\zeta_{p-1})$.

 At last, because $\psi \circ z\frac{d}{dz}= pz\frac{d}{dz}\circ\psi$,
 one has
 \begin{align*}
  & \pi^k\Gamma_p\biggl(\frac{k}{p-1}\biggr) =
 \Tr\biggl([\psi] \Big|
  M_{\alpha}^{\dag}\big/z\frac{d}{dz}M_\alpha^{\dag}\biggr)
 =(1-p)\Tr\bigl(\psi\big|{M_\alpha^{\dagger}}\bigr)\\
 = & -\sum_{\zeta_{p-1}}\zeta_{p-1}^{-k}E_\pi(\zeta_{p-1})
 = -\sum_{\zeta_{p-1}}\;\zeta_{p-1}^{-k}\zeta_p^{\zeta_{p-1}}.
 \end{align*}

 We refer to \cite{edpaase} for a detailed account. A more general form of the Gross-Koblitz,
proved along the same lines, shows that for any $1\leq k < p^r$, the product 
\begin{equation*}
  \prod_0^{r-1} \Gamma_p\biggl(\frac{p^i k}{p^r-1}\biggr)
 \end{equation*}
belongs to the cyclotomic field $\Q(\zeta_{p^r})$.
  
\end{para}
\begin{para}
 A more straightforward and elementary proof has been discovered
 by A. Robert. It goes as follows.

 As we have just seen, the right hand side of the Gross-Koblitz 
 formula can be written
 \begin{align*}
  &\;-\pi^{-k}\sum_{\zeta_{p-1}}\zeta_{p-1}^{-k}E_\pi(\zeta_{p-1})
  =-\pi^{-k}\sum_{n=0}^{\infty}
  \biggl(\sum_{\zeta_{p-1}}\zeta_{p-1}^{n-k}\biggr) e_n\\
  =&\;(1-p)\pi^{-k}\sum_{m=0}^{\infty} e_{(p-1)m+k}.
 \end{align*}

 Using the expansion of $\Gamma(\alpha)$ given in
 \ref{para-cont-of-gamma}, in the case of $\alpha=\frac{k}{p-1}$, we
 thus have to show that 
 \begin{equation*}
  \sum_{n=0}^{\infty}
   e_{pn+k}(-\pi)^{-n} \biggl(\frac{k}{p-1}\biggr)_n 
   = (1-p)\sum_{m=0}^{\infty} e_{(p-1)m+k}.
 \end{equation*}
 Denote the left hand side by $G_k$ (for any $k\in \N$).
 Due to the overconvergence of $E_\pi$, it is not difficult to see that 
 $\displaystyle\lim_{k\rightarrow \infty}G_k=0$.  One has Robert's formula:
 \begin{equation*}
  G_k-G_{p-1+k}=(1-p)e_k.
 \end{equation*}
 Summing up consecutive expressions, one gets a telescoping sum which yields
 the desired equality $\displaystyle  G_k= G_k-\lim_{m\rightarrow
 \infty}G_{(p-1)m+k} = (1-p)\sum_{m=0}^{\infty} e_{(p-1)m+k}$.
 It remains to prove Robert's formula. One first observes that 
 $z\frac{d}{dz}
 E_\pi=(\pi z-p\pi z^p)E_\pi(z)$, which yields the relation 
 $ne_n=\pi(e_{n-1}-pe_{n-p})$ for $n\geq p$, hence
 $\pi e_{p-1+m}=p\pi e_m + (m+p)e_{m+p}$ for $m\geq 0$. Then 
\begin{align*}
  &\; G_k-G_{p-1+k}\\
 =\;& e_k+ \sum_0^{\infty} e_{p(n+1)+k}(-\pi)^{-n-1}
 \biggl(\frac{k}{p-1}\biggr)_{n+1}\\
 &\phantom{(1-p)e_k+p(1-p)}
 -\sum_0^{\infty} e_{p-1+pn+k}(-\pi)^{-n} \biggl(\frac{k}{p-1}+1\biggr)_n\\
 =\;& e_k+ \sum_0^{\infty}
 \biggl[\frac{k}{p-1}e_{p(n+1)+k}+\pi e_{p-1+pn+k}\biggr]
 (-\pi)^{-n-1} \biggl(\frac{k}{p-1}+1\biggr)_{n}\\
 =\;& e_k + \sum_0^{\infty}
 \biggl[p\frac{(k+(n+1)(p-1))}{p-1}e_{p(n+1)+k}+p\pi e_{pn+k}\biggr]
 (-\pi)^{-n-1} \biggl(\frac{k}{p-1}+1\biggr)_{n}\\
 =\;&(1-p)e_k+p\sum_0^{\infty} \frac{(k+(n+1)(p-1))}{p-1}e_{p(n+1)+k}
 (-\pi)^{-n-1}\biggl(\frac{k}{p-1}+1\biggr)_{n}\\
  &\phantom{(1-p)e_k+p(1-p)}
 - p\sum_1^{\infty} \frac{k+n(p-1)}{p-1} e_{pn+k}(-\pi)^{-n}
 \biggl(\frac{k}{p-1}+1\biggr)_{n-1}\\
 =\;& (1-p)e_k.
\end{align*}
\end{para}
\begin{para}
 The right hand side of the Gross-Koblitz formula is a special case of an
 {\it exponential sum}, \ie an expression of the form
 \begin{equation*}
  S_r(\ovl{f},\ovl{g},\ovl{h})
   = \sum_{x_1,\ldots,x_d \in \F_{p^r},\;
   \ovl{h}(\udl{x})\neq 0}\;\chi\bigl(N_{\F_{p^r}/\F_{p}}
   \ovl{g}(\udl{x})\bigr)\exp\biggl(\frac{2\pi i}{p}
   \Tr_{\F_{p^r}/\F_{p}}\ovl{f}(\udl{x})\biggr) 
 \end{equation*}
 where $\ovl{h}$ is a polynomial with coefficients in
 $\F_{p^r}$, $\ovl{f}$ and $\ovl{g}$ are rational functions with
 coefficients in $\F_{p^r}$ with no pole where $\ovl{h}$
 vanishes, and where $\chi$ is a character of
 $\F_{p}^{\times}$.  This includes Gauss, Jacobi and Kloosterman
 sums as special cases.  In fact, it is known, after the classical works
 of Gauss, Artin, Weil, that counting solutions of systems of polynomial 
 equations in finite fields amounts to the computation of exponential
 sums $S_r(\ovl{f},\ovl{g},\ovl{h})$.  They are studied via their
 generating series, the so-called $L$-series:
 \begin{equation*}
  L(\ovl{f},\ovl{g},\ovl{h};t)=
   \exp\biggl(\sum_{r\geq 1}S_r(\ovl{f},\ovl{g},\ovl{h})\frac{t^r}{r}\biggr).
 \end{equation*}
 Dwork's methods, refined by P. Robba and others, allow to tackle the
 questions of the rationality of $L$, its degree, and of the functional
 equation relating $L(\ovl{f},\ovl{g},\ovl{h};t)$ to
 $L(\ovl{f},\ovl{g},\ovl{h};\frac{1}{p^rt})$.
 In the one-dimensional case, the solution is elementary and follows the
 pattern sketched in \ref{para-gross-koblitz}. Namely \cite{edpaase}:
 \begin{enumerate}
  \item One considers liftings $f,g,h$ of $\ovl{f},\ovl{g},\ovl{h}$ in 
	characteristic zero, and one introduces the
	affinoid set $S$ defined by $\vert h\vert =1$. One sets
	$F=g(z)^{1/p-1}\exp(\pi f(z))$. Then $F'/F$ is a rational function, 
	and one shows that the differential operator
	$\displaystyle \frac{d}{dz} + \frac{F'}{F}$ has an index in
	$\mathcal{H}^{\dagger}(S)$ (in a slightly generalized sense, and
	which can be computed thanks to the work of Robba). This is the 
	crucial point. Hence the cohomology spaces of the de
	Rham complex: $\Omega^0=F\mathcal{H}^{\dagger}(S)\rightarrow 
	\Omega^1=F\mathcal{H}^{\dagger}(S)dz$ are finite-dimensional over
	$\C_p$.
  \item One observes that $E:= F^{\phi-\id}$ (\ie $E(z)= 
	F(z^p)/F(z)$) is an element of $\mathcal{H}^{\dagger}(S)$. This
	allows to define two endomorphisms of the de Rham complex by setting:
	\begin{equation*}
	 \phi^0(F.g)=FEg^{\phi},\;\;\phi^1(Fgdz)
	  =FE\phi'g^{\phi}dz \; \text{ (with $\phi'(z)=pz^{p-1}$)};
	\end{equation*} 
	\begin{equation*}
	 \psi^0(F.g)(z)=F(z)\sum_{\phi(t)=z} \frac{g(t)}{E(t)}\;,\;\;
	  \psi^1(Fgdz)=\sum_{\phi(t)=z} \frac{g(t)}{E(t)\phi'(t)}.
	\end{equation*}
	One has $\psi^{\ast} \circ \phi^{\ast}=p\cdot\id$.  In particular, the
	operators $\mathrm{H}^i(\psi^{\ast}), \;i=0,1,$ are invertible
	on the cohomology spaces. The same argument as in
	\ref{para-gross-koblitz} shows that $\psi^i$ is nuclear on
	$\Omega^i$.  In particular, it has a trace.
  \item Polynomial approximation allows to prove the trace formula:
	\begin{equation*}
	 \tr(\psi^1)-\tr(\psi^0)=\tr(\mathrm{H}^1(\psi^{\ast}))
	  -\tr(\mathrm{H}^0(\psi^{\ast}))=S_r(\ovl{f},\ovl{g},\ovl{h}).
	\end{equation*}
	It follows that
	\begin{equation*}
	 L(\ovl{f},\ovl{g},\ovl{h};t)
	  =\frac{\det(\id-t\mathrm{H}^1(\psi^{\ast}))}
	  {\det(\id-t\mathrm{H}^0(\psi^{\ast}))}.
	\end{equation*}
	This is a rational function since the cohomology is
	finite-dimensional.  Moreover, since the
	$\mathrm{H}^i(\psi^{\ast})$
	are invertible, the computation of its degree amounts to the
	computation of the index of
	$\displaystyle \frac{d}{dz} + \frac{F'}{F}$ (the computation of
	the dimension of $\mathrm{H}^0$, $0 $ or $1$, being essentially
	trivial).
  \item The functional equation follows from a topological ``dual theory''
	in which the transpose of $\phi$ and $\psi$ play the roles of $\psi$
	and $\phi$ respectively.
 \end{enumerate}

 Let us remark at last that this method has an \emph{archimedean
 analogue} in the so-called thermodynamic formalism ---in dimension
 one.  The analogy is especially striking in  
 the presentation of D. Mayer \cite{cfart}:
 exponential sums correspond to ``partition functions'', $L$ to the 
 Ruelle zeta-function, Frobenius to the
 ``shift'', Dwork's operator $\phi^{\ast}$ to the ``transfer operator''
 (or Ising-Perron-Frobenius-Ruelle operator) and is given by a 
 ``composition operator'' (loc. cit. 7.2.2), its nuclearity is
 established by the same argument, and the trace formula has the same
 form (loc. cit. 7.17).

 In view of this close analogy, one could dream of an archimedean 
 proof of the rationality of Weil zeta
 functions, parallel to Dwork's $p$-adic proof...
\end{para}

\newpage
\section{The tale of $F(\frac{1}{2},\frac{1}{2},1;z)$.}\label{sec-tale}

\begin{abst}
This is a detailed concrete illustration of the somewhat abstract
non-archimedean notions discussed in the previous section. We refer to
\cite{hfml} (and \cite[9]{husemoeller87:_ellip}) for an excellent account of
the archimedean tale of $F(\frac{1}{2},\frac{1}{2},1;z)$.
\end{abst}

\subsection{Dwork's hypergeometric function.}\label{sub-dwork-hyper}

Let us consider the Legendre pencil of elliptic
curves with parameter $z\neq 0,1,\infty$, given in inhomogeneous coordinates by
\begin{equation*}
 y^2=x(x-1)(x-z).
\end{equation*}
As a scheme over ${\Z}[z,\frac{1}{{2z(1-z)}}]$, we denote it by $X$. The
first de Rham cohomology module $\mathrm{H}^1_{\mathrm{dR}}(X)$ is free
of rank $2$ over ${\Z}[z,\frac{1}{{2z(1-z)}}]$, and endowed with the
Gauss-Manin connection $\nabla$ (derivation with respect to the
parameter $z$); it is generated by the class $\omega$ of $\frac{dx}{y}$,
and $\nabla(\frac{d}{dz})(\omega)$. The canonical symplectic form
(cup-product) satisfies
$\an{\omega, \nabla(\frac{d}{dz})(\omega)}=\frac{2}{z(z-1)}$. The
Gauss-Manin connection is given by the hypergeometric differential
equation with parameters $(\frac{1}{2},\frac{1}{2},1)$
\begin{equation*}
\nabla(L_{\frac{1}{2},\frac{1}{2},1})(\omega)=0,\text{ with }
L_{\frac{1}{2},\frac{1}{2},1} =
z(1-z)\frac{d^2}{dz^2}+(1-2z)\frac{d}{dz}-\frac{1}{4}.
\end{equation*}
Let $p$ be an odd prime. The {\it Hasse invariant} is the polynomial
$h_p(z) \in {\Z}[\frac{1}{2}][z]$ obtained by truncating the
hypergeometric series $(-1)^\frac{p-1}{2}F(\frac{1}{2},\frac{1}{2},1;z)$
at order $\frac{p-1}{2}$. It enjoys the following well-known properties:
\begin{itemize}
 \item functional equations: $h_p(z)\equiv
       (-1)^\frac{p-1}{2}h_p(1-z)\equiv
       z^\frac{p-1}{2}h_p(\frac{1}{z})\pmod{p}$,
 \item for any $\ovl{\zeta} \in \ovl{\F}_p \setminus \lbrace0,1\rbrace$,
       one has $h_p(\ovl{\zeta})=0$ if and only if the elliptic
       curve $X_{\ovl{\zeta}}$ over $\ovl{\F}_p$ is {\it supersingular}, \ie
       has no geometric point of order $p$,
       \cf \cite[13.3]{husemoeller87:_ellip}.
\end{itemize}
The roots of $h_p \pmod{p}$ are distinct and lie in ${\F}_{p^2}$: we
denote them by $\ovl{\zeta}_1,\ldots,\ovl{\zeta}_{(p-1)/2}$, and choose
liftings ${\zeta}_1,\ldots,{\zeta}_{(p-1)/2}$ in
$\widehat{{\Z}_p^{\mathrm{ur}}}$.

Non-supersingular elliptic curves are called {\it ordinary}; they have
exactly $p$ points of order $p$.

Dwork's hypergeometric function is

\begin{equation*}
f_p(z)=
(-1)^\frac{p-1}{2}\frac{F(\frac{1}{2},\frac{1}{2},1;z)}
{F(\frac{1}{2},\frac{1}{2},1;z^p)}\in {\Z}\Bigl[\frac{1}{2}\Bigr][[z]].
\end{equation*}

One has $f_p(z)\equiv h_p(z)\pmod{p}$. Dwork discovered that although
the $p$-adic radius of convergence of this series is exactly $1$, $f_p$
\emph{does extend to a $p$-adic analytic function on ${\A}^1$ deprived
from the supersingular disks} $\matheur{D}(\zeta_j,1^-)$; we denote this
extension by the same symbol $f_p$.

In terms of this function, he obtained in 1958 his famous $p$-adic 
formula for the number of
rational points of an ordinary elliptic curve defined over ${\F}_{p^n}$:

\emph{for any $s_0\in {\F}_{p^n}, s_0\neq 0,1,\ovl{\zeta}_j,$
 the number of ${\F}_{p^n}$-points of $X_{s_0}$ is}

\begin{equation*}
1 - \prod_{k=0}^{n-1} f_p(\omega^{p^k}) +
 p^n \left(1 - \prod_{k=0}^{n-1} f_p(\omega^{p^k})^{-1}\right),
\end{equation*}
where $\omega$ denotes the unique $(p^n-1)$th root of unity $\equiv s_0
\pmod{p}$.

\subsection{The ordinary unit-root $F$-crystal.}\label{sub-ordinary-unit-root}

On completing $\widehat{{\Z}_p^{\mathrm{ur}}}[z,\frac{1}{{2z(1-z)}}]$
$p$-adically, $(\mathrm{H}^1_{\mathrm{dR}}(X), \nabla)$ gives rise to an
overconvergent $F$-crystal $(\mathcal{H}, \nabla, F)$. We do not discuss
here the construction of the Frobenius structure, which is a general
$p$-adic feature of Gauss-Manin connections. In fact, various analytic
or geometric constructions are available, but in our present case, the
Frobenius structure can be made quite explicit, \cf \cite{pc}.

\medskip
Let us introduce a few notations:
\begin{itemize}
 \item $\mathcal{R}$ = the $p$-adic completion of
       $\widehat{{\Z}_p^{\mathrm{ur}}}\bigl[z,\frac{1}{{h_p(z)}}\bigr]$ 
 \item $\mathcal{R}_{\mathrm{ord}}$ = the $p$-adic completion of
       $\widehat{{\Z}_p^{\mathrm{ur}}}\bigl[z,\frac{1}{{z(1-z)h_p(z)}}\bigr]$
 \item $S=\matheur{M}(\mathcal{R}
       [\frac{1}{p}]) = \matheur{D}(0,1^+) \setminus
       \bigl(\bigcup \matheur{D}(\zeta_j,1^-)\bigr)$
 \item $S_{\mathrm{ord}}=\matheur{M}(\mathcal{R}_{\mathrm{ord}}
       [\frac{1}{p}]) = S\setminus
       \bigl(\matheur{D}(0,1^-)\cup \matheur{D}(1,1^-)\bigr)$
       is the {\it ordinary locus} 
 \item $S_{\mathrm{nss}}= {\A}^1\setminus
       \bigl(\bigcup \matheur{D}(\zeta_j,1^-)\bigr)$
       is the {\it non-supersingular locus}.
\end{itemize}
\noindent The restriction of $(\mathcal{H},\nabla, F)$ to the ordinary
locus possesses a unique non-zero unit-root sub-$F$-crystal
$(\mathcal{U}, \nabla, F)$. This unit-root $F$-crystal extends over $S$
(and even over $S_{\mathrm{nss}}$ as an $F$-isocrystal with logarithmic
singularity at $\infty$ in the sense of \cite{atfff}) \cf also
\cite{ogus00:_ellip}; we use the same
symbol for the extension. It can be described along the following lines:
\begin{itemize}
 \item $\mathcal{U}$ is the unique rank-one horizontal
       submodule of $\mathcal{H}$.
 \item $\mathcal{U}\otimes_{\mathcal{R}_{\mathrm{ord}}}
       \widehat{\mathcal{R}_{\mathrm{ord}}^{\mathrm{et}}} = (\mathcal{H}
       \otimes_{\mathcal{R}_{\mathrm{ord}}}
       \widehat{\mathcal{R}_{\mathrm{ord}}^{\mathrm{et}}})^{\nabla}
       \otimes_{\widehat{{\Z}_p^{\mathrm{ur}}}}
       \widehat{\mathcal{R}_{\mathrm{ord}}^{\mathrm{et}}}.$ (Here
       $\widehat{\mathcal{R}_{\mathrm{ord}}^{\mathrm{et}}}$ is defined
       as in \ref{para:r}.)
 \item for any ordinary $s_0 \in \ovl{\F}_p$, the
       associated $p$-adic representation of
       $\pi_1^{\mathrm{alg}}(S_0,s_0)$ is
       $(\mathcal{H}\otimes_{\mathcal{R}_{\mathrm{ord}}}
       \widehat{\mathcal{R}_{\mathrm{ord}}^{\mathrm{et}}})^{\nabla =0,
       F(\phi)=\id} \simeq \mathrm{H}^1_{\mathrm{et}}(X_{s_0},{\Z}_p).$
 \item For any $u \in (\mathcal{U}\otimes_\mathcal{R}
       \widehat{\mathcal{R}^{\mathrm{et}}})^{\nabla}$, ``the'' image of
       $u$ in $\mathcal{O}(\matheur{D}(s_0,1^-))$ is a bounded solution of the
       differential operator $L_{\frac{1}{2},\frac{1}{2},1}$ on
       $\matheur{D}_{\C_p}(s_0,1^-)$; this also characterizes $\mathcal{U}$.
 \item $\mathcal{U}|_{\matheur{D}(0,1^-)}^{\nabla}$ has a
       canonical $\Z$-submodule, which can be identified with
       $\mathrm{H}^1((X_z)^{\mathrm{an}},{\Z})$ for any $z\in
       \matheur{D}(0,1^-)\setminus \lbrace 0\rbrace$; for a generator
       $u$, one has
       \begin{equation*}
	\an{\omega,u}=\sqrt{-1}F\Bigl(\frac{1}{2},\frac{1}{2},1;z\Bigr)
	 \text{ in } \mathcal{O}(\matheur{D}(0,1^-))
       \end{equation*}
       ($\sqrt{-1}$ appears as residue of $\frac{dx}{y}\big|_{z=0}$ at
       $x=0$), \cf \cite{pbl}.  Similarly,
       $\mathcal{U}|_{\matheur{D}(1,1^-)}^{\nabla}$ has a canonical
       $\Z$-submodule; for a generator $u$, one has
       \begin{equation*}
	\an{\omega,u} = \sqrt{-1}F\Bigl(\frac{1}{2},\frac{1}{2},1;1-z\Bigr)
	 \text{ in } \mathcal{O}(\matheur{D}(1,1^-)).
       \end{equation*}
\end{itemize}

\noindent The existence of the unit-root $F$-crystal $(\mathcal{U},
\nabla, F)$ over $\mathcal{R}$ then amounts to the two
function-theoretic facts:
\begin{equation*}
 \dlog F\Bigl(\frac{1}{2},\frac{1}{2},1;z\Bigr)\;\text{ and }\;
 f_p(z) =
 \Bigl(\sqrt{-1}F\Bigl(\frac{1}{2},\frac{1}{2},1;z\Bigr)\Bigr)^{1-\sigma}
\end{equation*}
both extend to units in $\mathcal{R}$. 

The inertia at any supersingular point maps \emph{onto} $\Aut
\mathrm{H}^1_{\mathrm{et}}(X_{s_0},{\Z}_p)\simeq {\Z}_p^{\times}$ (J.-I. 
Igusa, \cf \eg \cite{anopu}).  The inertia at $\infty$ acts as $\pm1$.
According to Crew's criterion, the $F$-crystal $(\mathcal{U},\nabla,F)$
is \emph{not} overconvergent.

We refer to \cite{opm} for a study in the same spirit of more general
hypergeometric equations and other differential equations ``coming from
geometry''.

\subsection{Analytic continuation: a pr\'ecis.}\label{sub-analytic-cont}

\begin{para}
 We have seen that Dwork's hypergeometric function extends to an analytic
 function on the whole of $S_{\mathrm{nss}}$. Some of its special values
 have been computed.  Let us mention \cite{thfwpp}, \cite{anjsasvogphf}
 \begin{enumerate}
 \item $f_p(1)=1$ (Koblitz),
 \item if $p \equiv 1\pmod{4}$,
       \[
       f_p(-1) =
       (-1)^\frac{p-1}{4}\frac{\Gamma_p(1/4)^2}{\Gamma_p(1/2)} =
       \frac{\Gamma_p(1/4)}{\Gamma_p(1/2)\Gamma_p(3/4)} \quad\text{ (Young).}
       \]
 \end{enumerate}
 (the condition $p \equiv 1\pmod{4}$ ensures that $-1$ is an
 ordinary modulus); $f_p(-1)$ is a Gauss integer (Van Hamme).
\end{para}

\begin{para}
 We have seen that the logarithmic derivative
 $\dlog F(\frac{1}{2},\frac{1}{2},1;z) $ also extends to an analytic
 function on $S_{\mathrm{nss}}$.  It satisfies the functional equations
 \begin{align*}
  \dlog F\Bigl(\frac{1}{2},\frac{1}{2},1;z\Bigr)
   & = -\dlog F\Bigl(\frac{1}{2},\frac{1}{2},1;1-z\Bigr)\\
   & = -\frac{1}{2z} - \frac{1}{z^2}
   \dlog F\Bigl(\frac{1}{2},\frac{1}{2},1;\frac{1}{z}\Bigr).
 \end{align*}
\end{para}

\begin{para}
 Let us turn to the more subtle case of $F(\frac{1}{2},\frac{1}{2},1;z)$
 itself. The general discussion of
 \ref{para-analytic-cont}/\ref{para-unit-root} applies to the ordinary
 unit-root $F$-crystal and tells us that $F(\frac{1}{2},\frac{1}{2},1;z)$
 \emph{admits an extension to any unit disk $\matheur{D}(s,1^-) \subset S$,
 analytic in that disk, and well-defined up to multiplication by a unit
 in ${\Z}_p$}. In other words, $\;F(\frac{1}{2},\frac{1}{2},1;z)$ is the
 specialization of an element (in fact, a unit) of
 $\widehat{\mathcal{R}^{\mathrm{et}}_{\mathrm{nss}}}$.
 It may be suggestive to denote
 such an element by $\;F(\frac{1}{2},\frac{1}{2},1;\eta)$, where $\eta$
 stands for the Berkovich generic point of $S$ corresponding to the
 sup-norm).

 \noindent In $\matheur{D}(1,1^-)$, this extension is
 $F(\frac{1}{2},\frac{1}{2},1;1-z)$ (up to $\Z_p^{\times}$). The
 refined structure of $F$-isocrystal with logarithmic singularity at
 $\infty$ allows to construct an extension of
 $F(\frac{1}{2},\frac{1}{2},1;z)^2$ to $\matheur{D}(\infty,1^-)$,
 analytic in that disk: it is $\frac{1}{z} F(\frac{1}{2},\frac{1}{2},1;\frac{1}{z})^2$ (up to $\Z_p^{\times}$); notice the
 square, which comes from the exponent $\pm 1/2$ at $\infty$.
\end{para}

\begin{para}
 The Krasner representation of $f_p(z)$ and $\dlog
 F(\frac{1}{2},\frac{1}{2},1;z)$ as uniform limits of rational functions
 without pole (nor zero) on $S$ can be made explicit: for $n>0$, let
 $g_n$ be the polynomial obtained by truncating
 $F(\frac{1}{2},\frac{1}{2},1;z)$ at order $p^n-1$; then $f_p(z)= \lim
 (-1)^{\frac{p-1}{2}}\frac{g_{n+1}(z)}{g_n(z^p)}$ and $\dlog
 F(\frac{1}{2},\frac{1}{2},1;z)=\lim\dlog g_n(z)$ \cite[3.4]{pc}.  On
 the other hand, $F(\frac{1}{2},\frac{1}{2},1;z)$ cannot be approximated
 by rational functions on $S$, but $L_{\frac{1}{2},\frac{1}{2},1}(g_n)$
 tends uniformly to $0$ on $\matheur{D}(0,1^+)$ \cite{sbdsdlaafh}.  One
 has

\begin{equation*}
\sqrt{-1}F\Bigl(\frac{1}{2},\frac{1}{2},1;z\Bigr) \equiv
 \bigl(h_p(z)\bigr)^{-\frac{1}{p-1}}
\pmod{p}
\end{equation*}
but we do not know any explicit representation of 
$F(\frac{1}{2},\frac{1}{2},1;z)$ as a $p$-adic limit of algebraic
functions on $S$.
\end{para}

\begin{para}
 Any solution of $L_{\frac{1}{2},\frac{1}{2},1}$ on an ordinary disk
 $\matheur{D}(s,1^-)$, which is not proportional to (the extended)
 $F(\frac{1}{2},\frac{1}{2},1;z)$, is not bounded in
 $\matheur{D}_{\C_p}(s,1^-)$,
 hence does not extend to neighboring disks in any reasonable
 \emph{analytic} sense. In spite of this obstruction to continuation,
 one can nevertheless ``jump from disk to disk'' (see figure \ref{fig7})
 \begin{figure}[h]
  \begin{picture}(300,100)(0,-10)
   \put(0,0){\includegraphics{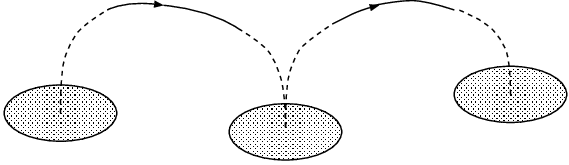}}
   \put(10,0){{\scriptsize $\matheur{D}(\zeta_{\mathrm{can}},1^-)$}}
  \end{picture}
  \caption{}
  \label{fig7}
 \end{figure}
 and extend them as locally analytic functions, in the following way.

 For each $\ovl{\zeta} \in S_0 \setminus \lbrace 0,1\rbrace$,
 there is a {\it canonical lifting} $\;\zeta_{\mathrm{can}} \in
 {\Z}_p^{\mathrm{ur}}$: namely, the modulus $\zeta_{\mathrm{can}}$ for which
 $X_{\zeta_{\mathrm{can}}}$ has complex multiplication by the quadratic
 order $\End(X_{\ovl{\zeta}})$.  Let us fix a branch of
 $F(\frac{1}{2},\frac{1}{2},1;z)$ in
 $\matheur{D}(\zeta_{\mathrm{can}},1^-)$, \ie an analytic solution
 $F_{\ovl{\zeta}}\in \widehat{{\Z}_p^{\mathrm{ur}}}[[z-{\zeta_{\mathrm{can}}}]]$
 of $L_{\frac{1}{2},\frac{1}{2},1}$ in
 $\matheur{D}(\zeta_{\mathrm{can}},1^-)$ which is a specialization of
 $F(\frac{1}{2},\frac{1}{2},1;\eta)\in
 (\widehat{\mathcal{R}^{\mathrm{et}}})^{\times}$ (any other branch is of
 the form $c.F_{\ovl{\zeta}}$ with $c \in {\Z}_p^{\times}$. Let
 $u_{\ovl{\zeta}} \in
 \mathcal{U}|_{\matheur{D}(\zeta_{\mathrm{can}},1^-)}^{\nabla}$ be
 defined by

 \begin{equation*}
 \an{\omega,u_{\ovl{\zeta}}}=\sqrt{-1}F_{\ovl{\zeta}}\;\text{ in }\;
 \mathcal{O}(\matheur{D}({\zeta_{\mathrm{can}}},1^-)).
 \end{equation*}
 On the other hand, let $v_{\ovl{\zeta}} \in
 \mathcal{M}|_{\matheur{D}(\zeta_{\mathrm{can}},1^-)}^{\nabla}$ be defined by
 \begin{equation*}
 \sqrt{-1} F_{\ovl{\zeta}}({\zeta_{\mathrm{can}}}).v_{\ovl{\zeta}}
  = \omega({\zeta_{\mathrm{can}}}),
  \;\text{ so that }\; \an{v_{\ovl{\zeta}},u_{\ovl{\zeta}}}=1.
 \end{equation*}
 Then $\an{\omega,v_{\ovl{\zeta}}}$ defines an
 unbounded solution of
 $L_{\frac{1}{2},\frac{1}{2},1}$ in
 $\matheur{D}_{\C_p}(\zeta_{\mathrm{can}},1^-)$.
\end{para}

\begin{para}
 Finally, in any supersingular disk $\matheur{D}_{{\C}_p}(\zeta_{i},1^-)$,
 there is no non-zero bounded solution of $L_{\frac{1}{2},\frac{1}{2},1}$.
\end{para}

\subsection{The Tate and Dwork-Serre-Tate parameters.}
\label{sub-serre-tate-para}

\begin{para}\label{para-tate-para}
 It turns out that the formal group
 ${\widehat{X}}_{(\widehat{\mathcal{R}_{\mathrm{ord}}^{\mathrm{et}}})}$
 of $X$ over $\widehat{\mathcal{R}_{\mathrm{ord}}^{\mathrm{et}}}$ is
 isomorphic to ${\widehat{\G}}_m$ (\cf~\cite[I]{fgade}). Any such
 isomorphism transforms the canonical differential
 $\omega_{\mathrm{can}}$ on ${\widehat{\G}}_m$ into
 $\Theta \frac{dx}{y}$, for a suitable element $\Theta \in
 \widehat{\mathcal{R}_{\mathrm{ord}}^{\mathrm{et}}}$ well-defined up to
 multiplication by an element of ${\Z}_p^{\times}$. It is possible to
 arrange normalizations so that the following relation hold (loc. cit.):
 
 \begin{equation*}
  \Theta
   = \biggl(\sqrt{-1}F\Bigl(\frac{1}{2},\frac{1}{2},1;\eta\Bigr)\biggr)^{-1}.
 \end{equation*}
 Let $\zeta\in \widehat{{\Z}_p^{\mathrm{ur}}}, \zeta \neq 0,1,$ be a point of
 $S$. By specialization, the element $\Theta(\zeta) \in
 \widehat{{\Z}_p^{\mathrm{ur}}}^{\times}$ (well-defined up to
 ${\Z}_p^{\times}$) links up $\omega_{\mathrm{can}}$ and $\frac{dx}{y}$ on
 $\widehat{X}_{\zeta}\simeq{\widehat{\G}}_m$; it is called the {\it Tate
 parameter} or {\it Tate constant} of $X_{\zeta}$.  It is well-defined
 up to multiplication by an element of ${\Z}_p^{\times}$.

 This is the $p$-adic analogue of the following familiar situation: let
 us consider the two-step uniformization of a complex elliptic curve
 \begin{equation*}
  \C \xrightarrow{\exp\bigl(\frac{2i\pi}{\omega_1}\cdot \;\,\bigr)}
   {\C}^{\times}
   \longrightarrow \C^{\times}/q^{\Z} \simeq
   {\C}/(\omega_1{\Z}+\omega_2{\Z})
 \end{equation*}
 with $q=\exp(2i\pi\frac{\omega_2}{\omega_1})$.  If the elliptic curve
 is $X_{\zeta}$ with $\zeta\in \matheur{D}(0,1^-)\setminus \lbrace
 0\rbrace$, and if $\omega_1$ and $\omega_2$ are fundamental periods of
 $\frac{dx}{y}$ ($\omega_1$ being the period attached to the vanishing
 cycle), then $\frac{2i\pi}{\omega_1}$ appears as the analogue of
 $\Theta(\zeta)$, and it is well-known that $\omega_1/2i\pi =
 iF(\frac{1}{2},\frac{1}{2},1;\zeta)$ (up to sign).
\end{para}

\begin{para}\label{para-dwork-para}
 We come back to the symplectic basis $v_{\ovl{\zeta}}, u_{\ovl{\zeta}}$
 and write $\frac{1}{{\sqrt{-1} F_{\ovl{\zeta}}}} \omega$ in the form
 $v_{\ovl{\zeta}} + \tau u_{\ovl{\zeta}}$, where

 \begin{equation*}
 \tau = -\frac{\an{\omega,v_{\ovl{\zeta}}}}{\an{\omega,u_{\ovl{\zeta}}}}
  \in
  (z-{\zeta_{\mathrm{can}}})\widehat{{\Q}_p^{\mathrm{ur}}}[[z-{\zeta_{\mathrm{can}}}]].
 \end{equation*}
 This defines an unbounded element of
 $\mathcal{O}(\matheur{D}_{{\C}_p}(\zeta_{\mathrm{can}},1^-))$ (notice that
 another choice of $u_{\ovl{\zeta}}$ multiplies $\tau$ by an element of
 $\Z_p^{\times}$).

 \noindent Applying $\nabla(d/dz)$ to
 $\frac{1}{{\sqrt{-1}F_{\ovl{\zeta}}}} \omega = v_{\ovl{\zeta}} + \tau 
 u_{\ovl{\zeta}}$ and using $\lan{\omega,\nabla(\frac{d}{dz})(\omega)} =
 \frac{2}{z(z-1)}$, one derives:
 \begin{equation*}
 \frac{d\tau}{dz} = \frac{2}{z(1-z) F_{\ovl{\zeta}}^2}.
 \end{equation*}
 The exponential of $\tau$ is called the {\it Dwork-Serre-Tate
 parameter}.  We shall recall later its meaning as a parameter for
 $p$-divisible groups. It satisfies
 the following remarkable integrality property \cite[Th.4]{pc}, \cite{tdd}
 \begin{equation*}
 q = e^{\tau} \in 1 + (z-{\zeta_{\mathrm{can}}})
 \widehat{{\Z}_p^{\mathrm{ur}}}[[z-{\zeta_{\mathrm{can}}}]].
 \end{equation*}
 This suggests the following question: does $q$ arise from an element of
 $\widehat{\mathcal{R}_{\mathrm{ord}}^{\mathrm{et}}}\;$? The answer is
 \emph{no}. This may be seen by considering $q\pmod{p}$: the 
 formula for $d\tau/dz = \dlog(q)$ shows that $q\pmod{p}$ is non-constant
 ($\dlog(q) \equiv \frac{2(h_p)^{2/p-1}}{z(z-1)}\pmod{p}$).  If $q$ comes
 from an element of $\widehat{\mathcal{R}_{\mathrm{ord}}^{\mathrm{et}}}$ which
 specializes to $1$ at every canonical modulus $\zeta_{\mathrm{can}}$, one has 
 $q\equiv 1\pmod{p}$, a contradiction.
\end{para}

\subsection{Complex counterpart: the supersingular locus.}
\label{sub-complex-supersing}

 Let us now revert things and try to understand the complex situation
 from the $p$-adic viewpoint!

\begin{figure}[h]
 \begin{picture}(150,180)(0,0)
  \put(0,0){\includegraphics[clip]{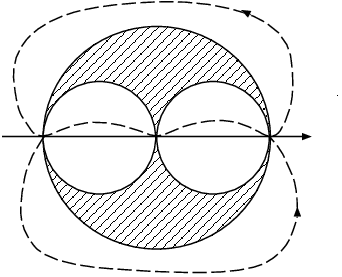}}
  \put(69,105){$\matheur{D}^{+}$}
  \put(69,30){$\matheur{D}^{-}$}
  \put(175,90){$\matheur{D}_{\infty}$}
  \put(42,82){$\matheur{D}_0$}
  \put(98,82){$\matheur{D}_1$}
  \put(0,58){$-1$}
  \put(55,58){$0$}
  \put(55,66.7){\makebox(0,0)[c]{$\bp$}}
  \put(80,58){$\FRAC{1}{2}$}
  \put(95,58){$1$}
  \put(95,66.7){\makebox(0,0)[c]{$\bp$}}
  \put(138,58){$2$}
  \put(145,108){$\gamma^{+}$}
  \put(150,40){$\gamma^{-}$}
 \end{picture}
 \caption{}
 \label{pic3}
\end{figure}

 Let $\matheur{D}^+$ and $\matheur{D}^-$ be the (complex-conjugate)
 connected components of $\matheur{D}(\frac{1}{2},\frac{3}{2}^-)\setminus
 \bigl(\matheur{D}(\frac{-1}{4},\frac{3}{4}^+) \cup
 \matheur{D}(\frac{5}{4},\frac{3}{4}^+)\bigr)$. Let $S$ 
 be the complement of $\matheur{D}^+ \cup \matheur{D}^-$ in the complex plane.
 Note that $S$ is closed and arcwise connected, and that
 $\pi_1^{\mathrm{top}}(S)$ 
 is a free group with two generators $\gamma^+,\gamma^-$ (see figure
 \ref{pic3}). Its 
 interior $S^{\circ}$ has three connected (simply-connected) components
 $\matheur{D}_0= \matheur{D}(\frac{-1}{4},\frac{3}{4}^-),\; \matheur{D}_1=
 \matheur{D}(\frac{5}{4},\frac{3}{4}^-)\;$ and $\matheur{D}_{\infty}$.

 Let $\mathcal{O}(S)$ denote the ring of continuous functions on $S$
 analytic in $S^{\circ}$. According to Mergelyan's theorem
 \cite{uatfoacv}, they are uniform limits of rational functions on every
 compact $K\subset S$ such that $\pi_0(\C\setminus K,0)$ is finite.
 This ring is however not stable under differentiation, and we consider
 its differential closure $\mathcal{R}$ in $\mathcal{O}(S^{\circ})$. We
 also consider the integral closure $\mathcal{R}^{\mathrm{et}}$ of
 $\mathcal{R}$ in $\mathcal{O}(S^{\circ})$. Every element of
 $\mathcal{R}^{\mathrm{et}}$ defines a multivalued locally analytic
 function on $S$, \ie an analytic germ which may be analytically
 continued along any path of $S$ not ending at $-1,\frac{1}{2}, 2$, in
 such a way that the germ at the other extremity is analytic.

 For instance, $F(\frac{1}{2},\frac{1}{2},1;z)$ may be viewed as an element
 of $\mathcal{R}^{\mathrm{et}}$: in fact,
 $F(\frac{1}{2},\frac{1}{2},1;z)^4 \in \mathcal{O}(S)$. An explicit
 continuous extension of $F(\frac{1}{2},\frac{1}{2},1;z)^4$ from
 $\matheur{D}_0$ to $S$ is given by $F(\frac{1}{2},\frac{1}{2},1;1-z)^4$
 in $\matheur{D}_1$, $\frac{1}{z^2 }
 F(\frac{1}{2},\frac{1}{2},1;\frac{1}{z})^4$ in $\matheur{D}_{\infty}$,
 and the values
 $\frac{{\Gamma(1/4)}^8}{64\pi^6},\;\frac{{\Gamma(1/4)}^8}{16\pi^6} ,
 \;\frac{{\Gamma(1/4)}^8}{64\pi^6}$ at $-1,\;\frac{1}{2},\;2$
 respectively \cite[p.189]{f2}.  On the other hand, note that $\dlog
 F(\frac{1}{2},\frac{1}{2},1;z)$ is in $\mathcal{R}$ but not
 in $\mathcal{O}(S)$.

The de Rham cohomology of the Legendre elliptic pencil gives
rise to an $\mathcal{R}[1/z(1-z)]$-module with connection
$(\mathcal{H},\nabla)$.  It admits a \emph{unique} non-zero horizontal
submodule $\mathcal{U}$ which extends to $S$; we use the same
symbol for the extension. It can be described along the following lines:

\begin{itemize}
 \item For any $u \in (\mathcal{U}\otimes_\mathcal{R}
       {\mathcal{R}^{\mathrm{et}}})^{\nabla}$, ``the'' image of $u$ in
       $\mathcal{O}(\matheur{D}_j) (j=0,1, \text{ or } \infty),$ is a
       bounded solution of the differential operator $L_{\frac{1}{2},\frac{1}{2},1}$.
 \item $\mathcal{U}\otimes_\mathcal{R}
       {\mathcal{R}^{\mathrm{et}}[\frac{1}{z(1-z)}]} =
       \bigl(\mathcal{H}\otimes_{\mathcal{R}[\frac{1}{z(1-z)}]}
       {\mathcal{R}^{\mathrm{et}}[\frac{1}{z(1-z)}]}\bigr)^{\nabla}\otimes_\C
       {\mathcal{R}^{\mathrm{et}}[\frac{1}{z(1-z)}]}.$
 \item $\bigl(\mathcal{H}\otimes_{\mathcal{R}[\frac{1}{z(1-z)}]}
       {\mathcal{R}^{\mathrm{et}}[\frac{1}{z(1-z)}]}\bigr)^{\nabla}$ has a
       canonical ${\Z}[i]$-submodule ($i=\sqrt{-1}$) which can be
       locally identified with the part
       $\mathrm{H}^1(X_z^{\mathrm{an}},{\Z}[i])_{\mathrm{isotriv}}$ of
       $\mathrm{H}^1(X_z^{\mathrm{an}},{\Z}[i])$ where
       $\pi_1^{\mathrm{top}}(S,z)$ acts through a finite group.
 \item $\mathcal{U}|_{\matheur{D}_0}^{\nabla}$ has a
       canonical $\Z$-submodule, which can be identified with the part
       of $\mathrm{H}^1(X_z^{\mathrm{an}},{\Z}[i])_{\mathrm{isotriv}}$ invariant
       under complex conjugation. One of the two generators $u$ satisfies
       \begin{equation*}
	\Lan{\omega,\frac{u}{2i\pi}}=
	 i F\Bigl(\frac{1}{2},\frac{1}{2},1;z\Bigr)
	 \;\text{ in }\; \mathcal{O}(\matheur{D}_0).
       \end{equation*}
       Similarly, $\mathcal{U}|_{\matheur{D}_1}^{\nabla}$ has a canonical
       $\Z$-submodule, and for one of the two generators $u$, one has
       \begin{equation*}
	\Lan{\omega,\frac{u}{2i\pi}} =
	 i F\Bigl(\frac{1}{2},\frac{1}{2},1;1-z\Bigr)\;\text{ in }\;
	 \mathcal{O}(\matheur{D}_1).
       \end{equation*}
       Note the occurrence of $2i\pi$ in these formulas. 
 \item The local monodromy $\gamma^+$ around $\matheur{D}^+$
       maps \emph{onto}
       \[\Aut \mathrm{H}^1(X_{s},{\Z}[i])_{\mathrm{isotriv}}\simeq
       {\Z}/4{\Z}.\]
       Same for $\gamma^-$.
       The local monodromy at $\infty$ acts as $\pm1$.
 \item In $\matheur{D}^+$ and in $\matheur{D}^-$, there is no non-zero
       bounded analytic solution
       of $L_{\frac{1}{2},\frac{1}{2},1}$.
\end{itemize}

This is quite similar to \subsecref{sub-ordinary-unit-root},
\subsecref{sub-analytic-cont}, $\matheur{D}^+\cup \matheur{D}^-$ playing
the role of the \emph{supersingular locus}. This picture is assuredly
very different from the traditional view of monodromy for
$F(\frac{1}{2},\frac{1}{2},1;z)$, and may shed some light upon the
divergences between the topological and algebraic approaches to analytic
continuation in the $p$-adic case (here of course, the coexistence of
the two pictures is explained by the fact that the principle of unique
continuation (\defref{dfn-unique-cont}) fails for the sheaf
of germs of continuous functions on $S$, analytic on $S^{\circ}$).

The tale of $F(\frac{1}{2},\frac{1}{2},1;z)$ is not finished: we have
not yet explored the islands of supersingularity. We shall reach them in
\subsecref{sub-p-adic-betti-lattice}.

\newpage
\section{Abelian periods as algebraic integrals.}

\begin{abst}
 We discuss periods of abelian varieties. Their $p$-adic counterparts 
 live naturally in Fontaine's ring $\matheur{B}_{\mathrm{dR}}$. We present 
 Colmez' construction of abelian $p$-adic periods, which relies on 
 $p$-adic integration and reflects as closely as possible the complex
 picture.
 We also deal with the concrete evaluation of elliptic $p$-adic periods.
\end{abst}

\subsection{Over $\C$.}\label{sub-complex-integ}

Let $A$ be an abelian variety over $\C$ of dimension $g$, 
and $\omega_1,\ldots,\omega_g$ a basis of invariant differential forms.
Let $\Lambda$ be the image of the map
\begin{equation*}
 \iota:\mathrm{H}_1(A(\C),\Z)\rightarrow \C^g\; ;\;\; \gamma\mapsto
 \left(\int_{\gamma}\omega_i\right)_i.
\end{equation*}
One has a canonical isomorphism 
\begin{equation*}
\C^g/\Lambda\stackrel{\sim}{\longrightarrow}A(\C).
\end{equation*}
Let us denote the projection $\C^g\rightarrow A(\C)$ by $\mathrm{pr}$.
For any differential one-form $\omega$ of the second kind, one can 
consider a primitive function $f_{\omega}$
of $\mathrm{pr}^{\ast}\omega$; it is meromorphic (univalued) on 
$\C^g$, because all residues of $\omega$ are $0$ by definition, and
it is unique up to addition of an arbitrary constant.
Then $f_{\omega}(z_1+z_2+z_3)-f_{\omega}(z_1+z_2)-f_{\omega}(z_1+z_3)
+f_{\omega}(z_1)$ defines a periodic function on $(\C^g)^3$, hence
induces a meromorphic function on $A(\C)^3$, which we denote by
$F^3_{\omega}$.
Note that the function $F^3_{\omega}$ can also be defined purely
algebraically by the following conditions:
\begin{itemize}
\item $F^3_{\omega}(z_1,0,z_3)=F^3_{\omega}(z_1,z_2,0)=0$.
\item $dF^3_{\omega}=m^{\ast}_{123}\omega-m^{\ast}_{12}\omega
-m^{\ast}_{13}\omega+m^{\ast}_1\omega,$
\end{itemize}
where $m_{123}$ is the addition $A^3\rightarrow A$ sending 
$(z_1,z_2,z_3)$ to $z_1+z_2+z_3$, etc.
For $\gamma\in\mathrm{H}_1(A(\C),\Z)$ the value
$f_{\omega}\bigl(i(\gamma)+a\bigr)-f_{\omega}(a)$ 
does not depend on $a\in\C^g$ ($a$ being chosen so that
$i(\gamma)+a$ and $a$ are not poles of $f_{\omega}$),
and one has the equality
\begin{equation*}
f_{\omega}\bigl(i(\gamma)+a\bigr)-f_{\omega}(a)=
\int_{\gamma}\omega,
\end{equation*}
which describes the period pairing 
\begin{equation*}
\mathrm{H}^1_{\mathrm{dR}}(A)
\times\mathrm{H}_1(A(\C),\Z)\longrightarrow\C
\end{equation*}
(or equivalently, the isomorphism
$\mathrm{H}^1_{\mathrm{dR}}(A)\otimes\C
\stackrel{\sim}{\rightarrow}
\mathrm{H}^1_{\mathrm{B}}(A(\C),\Z)\otimes\C$.)
Here $\mathrm{H}^1_{\mathrm{dR}}(A)$ denotes the first \emph{algebraic}
de Rham cohomology group of $A$, which coincides with the group of 
differential forms of the second kind modulo exact forms.

\subsection{Over the $p$-adics; prolegomena}\label{sub-prolegomena}

If we try to translate this into the $p$-adic setting, one has to face
at once the problem: what is \emph{integration over a loop}?

\begin{para}\label{para-padic-int-1}
 One tentative way is via Berkovich's theory, where loops do exist.  Let us for
 instance consider the Legendre elliptic curve $X_z$. Viewed as a
 $p$-adic space, for $z \in \matheur{D}(0,1^-)\setminus \lbrace 0 \rbrace$
 ($p\neq 2$), this is a Tate curve: $X_z^{\mathrm{an}} \simeq
 {\C}_p^{\times}/q^{\Z}$. The canonical differential
 $\omega_{\mathrm{can}}$ inherited from ${\C}_p^{\times}$ and
 $\frac{dx}{y}$ are proportional:

 \begin{equation*}
 \omega_{\mathrm{can}}=\Theta(z).\frac{dx}{y} \;\text{ with }\;
 1/\Theta(z)=\sqrt{-1}F\Bigl(\frac{1}{2},\frac{1}{2},1;z\Bigr).
 \end{equation*}

 The basic element $\gamma$ of $\mathrm{H}_1(X_z^{\mathrm{an}},{\Z})$
 can be identified with the generator $q$ of $q^{\Z}$, and it is natural
 to set $\int_{\gamma}\omega_{\mathrm{can}}=\int_{0}^q \frac{dt}{t}=
 \log(q)$ (choosing a branch of the $p$-adic logarithm). One then has

 \begin{equation*}
 \omega_2^{(p)}:=\int_{\gamma}^{(p)}\frac{dx}{y} =
 \sqrt{-1}F\Bigl(\frac{1}{2},\frac{1}{2},1;z\Bigr)\log(q).
 \end{equation*}
 This formula, as well as
 \begin{align*}
 & \sqrt{q} =
 \frac{z}{16}
 e^{\frac{F^{\ast}(\frac{1}{2},\frac{1}{2},1;z)}
         {F(\frac{1}{2},\frac{1}{2},1;z)}}\\
 \biggl( & \text{with }F^{\ast}\Bigl(\frac{1}{2},\frac{1}{2},1;z\Bigr)
  = 4\sum_{n>0}\binom{2n}{n}^2 \biggl(\sum_1^{2n} \frac{(-1)^{m-1}}{m}\biggr) 
 \Bigl(\frac{z}{16}\Bigr)^n\biggr)
 \end{align*}
 are the same as those encountered in the complex situation. This leaves
 the open problem: how can one construct ``the other period''
 $\omega_1^{(p)}$, \ie the $p$-adic analogue of
 $\omega_1=2\pi\sqrt{-1}F(\frac{1}{2},\frac{1}{2},1;z)$, since the
 corresponding loop is missing?
\end{para}

\begin{para}\label{para-padic-int-2}
 One simple tentative answer would be to replace the topological
 covering ${\C}_p^{\times}\rightarrow X_z^{\mathrm{an}}$ by other \'etale
 coverings. Natural candidates for this purpose are \'etale coverings of
 order $p^n$, especially those corresponding to torsion points of order
 $p^n$ which are close to the origin.  Let $\zeta_{p^n}$ be a $p^n$th
 root of unity in ${\C}_p^{\times}$ and let $x_n$ be its image in
 $X_z^{\mathrm{an}}$. In the corresponding complex case, taking
 $\zeta_{p^n}=e^{2i\pi/{p^n}}$ would give the right answer
 $p^n\int_{0}^{x_n}\frac{dx}{y}=\omega_1$.  In the $p$-adic case, we get
 instead $p^n\int_{0}^{x_n}\frac{dx}{y}=0$.  Indeed, already in the case
 of the multiplicative group, we have
 $p^n\int_{1}^{\zeta_{p^n}}\frac{dt}{t}=p^n \log(1+(\zeta_{p^n}-1))=0\,$  
 $p$-adically: in other words, $2i\pi$ ``\emph{is missing}''.
\end{para}

\begin{para}[{\itshape Riemann-Shnirelman sums.}]
 It is appropriate to evoke here Shnirel\-man's approach to integration over loops,
 indeed one of the earliest works in $p$-adic analysis. This
 is an adaptation to the $p$-adic case of the computation of integrals by
 Riemann sums
 \begin{align*}
  \frac{1}{2i\pi}\int_{C(a,r)} f(z)dz & = \int_0^1 f(a+re^{2i\pi\theta})
   re^{2i\pi \theta}d\theta
   & = \,\lim_{m\to \infty}\frac{1}{m}\sum_{\zeta^m=1}f(a+r\zeta) r\zeta.
 \end{align*}
 Now, let $a, r \in \C_p$ and let $f$ be a $\C_p$-valued function
 on the circumference
 $C(a,|r|)$.  Let us denote the limit
 \begin{equation*}
  \lim_{m\to \infty}\frac{1}{m}\sum_{\zeta^m=1} f(a+r\zeta).r\zeta 
 \end{equation*}
 symbolically
 by $\int_{C(a,r)} f(z)dz$. For better convergence, Shnirelman actually
 restricted the
 limit to those integers $m$ prime to $p$. In the case of an analytic
 function $f$ on $C(a,|r|)$, this restriction is unnecessary (one can even
 take $m=p^k$) and one has the following analogue
 of Cauchy's theorem of residues (\cf \cite[app.]{paascorw})

 \begin{lem}
  Assume that $f$ is a meromorphic function on $\matheur{D}(a,|r|^+)$,
  and that its poles $z_1,\ldots, z_{\nu}$ all lie in
  $\matheur{D}(a,|r|^-)$. Then $\int_{C(a,r)} f(z)dz$ exists and equals
  $\sum \Res_{z_i} f$.
 \end{lem}

 In the (trivial) special case $f = 1/z$, we get
 $\int_{C(0,1)}dz/z =1$ (not $0$ as
 in \ref{para-padic-int-2}!), but $2i\pi$ is still missing.
\end{para}

\begin{para}
 Actually, there is no way to remedy this if one remains in $\C_p$.  A
 deeper reason for that, due to Tate, is that while
 $\Gal({\ovl{\Q}}_p/{\Q}_p)$ acts on the inverse system of $p^n$th
 roots of unity ($n \geq 0$) through the cyclotomic character $\chi$,
 there is no element $(2i\pi)_p$ in $\C_p$ such that $g((2i\pi)_p)=\chi
 (g) (2i\pi)_p$ for every $g \in \Gal({\ovl{\Q}}_p/{\Q}_p)$.

 \noindent It turns out that $(2i\pi)_p$ is in a sense the only
 ``missing piece'': there is a good theory of $p$-adic periods (due to
 J.M. Fontaine, W. Messing \cite{ppapec}, G. Faltings) which lives in
 some ${\Q}_p$-algebra 
 isomorphic to $\C_p[[(2i\pi)_p]]$, and which we shall now describe in the
 case of abelian varieties.
\end{para}

\subsection{The Fontaine ring $\matheur{B}_{\mathrm{dR}}$.}
\begin{para}
 {\itshape The ring $\matheur{R}$.\ }\index{000R@$\matheur{R}$}
 Let $\mathcal{O}_{\C_p}$ be the ring of
 integers of $\C_p$, \ie $\{x \in \C_p$, $|x|_p \leq 1\}$.  Let us set
 \begin{equation*}
  \matheur{R}:=\underset{x\mapsto x^p}{\limproj}
   \OO_{\C_p},
 \end{equation*}
 \ie $\matheur{R}$ is the set of all series
 $\left(x^{(n)}\right)_{n\in\N}$
 such that $\left(x^{(n+1)}\right)^p=x^{(n)}$.
 This is in fact a ring of characteristic $p$ with 
 \begin{equation*}
  \left(\bigl(x^{(n)}\bigr)+\bigl(y^{(n)}\bigr)\right)^{(n)}=
   \lim_{m\rightarrow\infty}\left(x^{(n+m)}+y^{(n+m)}\right)^{p^m}
 \end{equation*}
 and
 \begin{equation*}
  \left(\bigl(x^{(n)}\bigr)\cdot\bigl(y^{(n)}\bigr)\right)^{(n)}=
   x^{(n)}\cdot y^{(n)}.
 \end{equation*}
 Let $\mathrm{W}(\matheur{R})$ be the Witt ring with coefficients in
 $\matheur{R}$. for $x\in \matheur{R}$, let $[x]$ denote the Teichm\"uller
 representative in $\mathrm{W}(\matheur{R})$.
\end{para}

\begin{para}
 {\itshape The ring $\matheur{B}^+_{\mathrm{dR}}$.\ }
 \index{000bdr+@$\matheur{B}^+_{\mathrm{dR}}$}
 To any element 
 $(x_0,x_1,\ldots,x_n,\ldots)=\sum p^n[x^{p^{-n}}_n]$ in
 $\mathrm{W}(\matheur{R})$, we associate
 \begin{equation*}
 \theta((x_n))=\sum^{\infty}_{n=0}p^nx^{(n)}_n.
 \end{equation*}
 This defines a surjective homomorphism
 $\theta\colon\mathrm{W}(\matheur{R})\rightarrow \OO_{\C_p}$,
 whose kernel is principal.
 It extends to a homomorphism $\theta:
 \mathrm{W}(\matheur{R})\left[\frac{1}{p}\right]\rightarrow\C_p$,
 and $\matheur{B}^+_{\mathrm{dR}}$ is defined as the $(\Ker\theta)$-adic 
 completion
 \begin{equation*}
 \matheur{B}^+_{\mathrm{dR}} :=\lim_{\longleftarrow}
 \mathrm{W}(\matheur{R})\textstyle{\left[\frac{1}{p}\right]}
 \Big/(\Ker\theta)^n.
 \end{equation*}
 By continuity, $\theta$ further extends to a homomorphism $\theta\colon
 \matheur{B}^+_{\mathrm{dR}}\rightarrow\C_p$.
 Then $\matheur{B}^+_{\mathrm{dR}}$ is a complete discrete valuation
 ring with maximal ideal $\Ker\theta$ and residue field $\C_p$.
 Moreover, the Galois group $\Gal(\ovl{\Q_p}/\Q_p)$ acts on 
 $\matheur{B}^+_{\mathrm{dR}}$ in such a way that $\theta$ is equivariant 
 with respect to this Galois action,
 \begin{equation*}
  \Gr^{\bp}\matheur{B}^+_{\mathrm{dR}}\simeq \bigoplus_{r\in \N}\C_p(r),
 \end{equation*}
 where $\Gr^{\bp}$ refers to the filtration by the powers of
 $\Ker\theta$, and where the ``twist'' $(r)$ indicates that the 
 $\Gal(\ovl{\Q_p}/\Q_p)$-action is twisted by the $r$th power of the
 cyclotomic character.

 It turns out that $\matheur{B}^+_{\mathrm{dR}}$ contains naturally a copy of
 ${\ovl{\Q}}_p$ (which $\theta$ maps isomorphically to ${\ovl{\Q}}_p \subset
 \C_p$). More precisely, P. Colmez \cite{lnasddb} has shown that
 $\matheur{B}^+_{\mathrm{dR}}$ is the separated \emph{completion} of
 ${\ovl{\Q}}_p$ with respect to the topology defined by taking
 $\Bigl(p^n\mathcal{O}^{(k)}_{{\ovl{\Q}}_p}\Bigr)_{n,k}$ as a basis of
 neighborhoods of $0$, where $\mathcal{O}^{(k)}_{{\ovl{\Q}}_p}$ denotes
 the subring of ${\ovl{\Z}}_p$ defined inductively as
 \begin{equation*}
 \mathcal{O}^{(0)}_{{\ovl{\Q}}_p}
  ={\ovl{\Z}}_p,\;\mathcal{O}^{(k)}_{{\ovl{\Q}}_p}
  =\Ker\Bigl(d: \mathcal{O}^{(k-1)}_{{\ovl{\Q}}_p}
  \rightarrow
  \Omega^1_{\mathcal{O}^{(k-1)}_{{\ovl{\Q}}_p}/{\Z}_p}\otimes
  {\ovl{\Z}}_p\Bigr).
 \end{equation*}
 It is easy to deduce from this description that
 $\matheur{B}^+_{\mathrm{dR}}$ contains $\widehat{{\Z}_p^{\mathrm{ur}}}$.
\end{para}

\begin{para}
 {\itshape Some remarkable elements of $\matheur{B}^+_{\mathrm{dR}}$.}
 For any $x \in \matheur{B}^+_{\mathrm{dR}}$ such that $|\theta(x)-1|_p< 1$,
 the series $\log(x) = - \sum_{n>0} \frac{(1-x)^n}{n}$ converges in
 $\matheur{B}^+_{\mathrm{dR}}$.  In particular, let $\underline {z} =
 (\ldots,z^{(1)},z^{(0)})$ be an element of $\matheur{R}$ such that $z^{(0)}\in
 {\ovl{\Z}}_p$ or $\widehat{{\Z}_p^{\mathrm{ur}}}$.  Then one can define the
 element \index{000Log@Log}
 $\Log \udl{z} = \log
 \bigl(\frac{z^{(0)}}{[\udl{z}]}\bigr)$\footnote{we follow Colmez' sign 
 convention; Fontaine's $\LOG$ is $\log - \Log$.
 }.
 Note that $\theta(\Log \underline{z})=0$. 

 In the special case where each
 $z^{(n)}=\zeta_{p^n}$ is a primitive $p^n$th root of unity,
 $\Log(\ldots, \zeta_{p^n}, \ldots,1)$ is the element $(2i\pi)_p$ 
 \index{0002ipp@$(2i\pi)_p$}
 we were looking for. This is a generator of $\Ker \theta$, and 
 $\Gal({\ovl{\Q}}_p/{\Q}_p)$ acts on it through the cyclotomic character;
 in other words, we have a $\Gal({\ovl{\Q}}_p/{\Q}_p)$-equivariant
 isomorphism

 \begin{equation*}
 \Gr^{\bp}\matheur{B}^+_{\mathrm{dR}}\;\simeq\; {\C}_p[(2i\pi)_p].
 \end{equation*}
 
 Note that another choice of $(\ldots, \zeta_{p^n},
 \ldots,1)$ changes $(2i\pi)_p$ by multiplication by a unit in ${\Z}_p$.
 Similarly,
 up to addition of an element of $(2i\pi)_p{\Z}_p$, $\Log\;\underline
 {z}\;$ depends only on
 $z^{(0)}$; it is sometimes simply denoted by $\Log\;z^{(0)}$.

 By choosing a double embedding of ${\ovl{\Q}}$ into $\C$ and
 ${\ovl{\Q}}_p$,
 one obtains a canonical element $(2i\pi)_p$, attached to the sequence
 $(\ldots, e^{2i\pi/{p^n}},\ldots,1)$. Similarly, if $z^{(0)}\in
 {\ovl{\Q}}$, then $\Log z^{(0)}$ is 
 well-defined up to addition of an element of $(2i\pi)_p{\Z}$.
\end{para}

\subsection{Colmez' construction of abelian $p$-adic periods.}

Let $A$ be an abelian variety defined over a $p$-adic local field
$K$. The Fontaine-Messing $p$-adic period pairing is a pairing

\begin{equation*}
\int^{(p)}\;:\;\mathrm{H}^1_{\mathrm{dR}}(A) \otimes T_p(A_{\ovl{K}})
 \longrightarrow \matheur{B}^+_{\mathrm{dR}}
\end{equation*}
where $T_p(A_{\ovl{K}}) = \varprojlim \Ker([p^n]: A_{\ovl{K}}\rightarrow
A_{\ovl{K}})$ is the Tate module (a ${\Z}_p$-module of rank $2\dim A$
with $\Gal({\ovl{K}}/K)$-action).

We present Colmez' construction, which is parallel to \ref{sub-complex-integ}.

Let $\omega$ be a differential form on $A$ of second kind, and
$F_{\omega}^3$ the function on $A^3$ determined as in
\ref{sub-complex-integ}.  Then:
\begin{pro}[{\cite[4.1]{ppdva}}]
There exists a locally meromorphic function $F_{\omega}$ on 
$A(\matheur{B}_{\mathrm{dR}})$, unique up to constant, such that:

\begin{enumerate}
 \renewcommand{\theenumi}{\arabic{enumi}}
 \item $dF_{\omega}=\omega$.
 \item $F_{\omega}(z_1+z_2+z_3)-F_{\omega}(z_1+z_2)-
       F_{\omega}(z_1+z_3)+F_{\omega}(z_1)=F^3_{\omega}(z_1,z_2,z_3)$.
 \item If $\omega=dF$, then $F_{\omega}=F$.
\end{enumerate}

\noindent
Moreover, if $\alpha\colon A_1\rightarrow A_2$ is a morphism of abelian 
varieties and $\omega$ is a differential form of second kind on
$A_2$, then $F_{\alpha^{\ast}\omega}=\alpha^{\ast}F_{\omega}$.
\end{pro}

Note that the function $F_{\omega}$ is \emph{not} multivalued; this
fact comes from the following lemma specific to the $p$-adic case:

\begin{lem}[{\cite[4.3]{ufdpplpdvaamc}}]
For any neighborhood $V$ of $0$ in $A(\matheur{B}_{\mathrm{dR}}^+)$, there
exists an open subgroup $U$ of $A(\matheur{B}_{\mathrm{dR}}^+)$ contained in 
$V$ such that $A(\matheur{B}_{\mathrm{dR}}^+)/U$ is a torsion group.
\end{lem}

Take a proper model $\mathscr{A}$ of $A$ over $\OO_{K}$.
Let $\gamma=(\cdots,u_n,\cdots,u_2,u_1=0)\in\mathrm{T}_p(A_{\ovl{K}})$ 
with each $u_n\in\mathscr{A}(\OO_{\C_p})$, and choose
$a_n\in\mathscr{A}(\matheur{B}^+_{\mathrm{dR}})$ so that neither
$a_n$ nor $a_n+_A\widehat{u}_n$ is close to a pole of $\omega$.
For suitable liftings 
$\widehat{u}_n\in\mathscr{A}(\matheur{B}^+_{\mathrm{dR}})$ of 
$u_n$ (\ie $\theta(\widehat{u}_n)=u_n$), the following holds. 

\begin{thm}[{\cite[5.2]{ufdpplpdvaamc}}]
The limit 
\begin{equation*}
 \oint^{(p)}_{\gamma}\omega \;:=\;
 \lim_{n\rightarrow\infty}p^n\bigl(F_{\omega}(a_n)-F_{\omega}
 (a_n+_A\widehat{u}_n)\bigr)
\end{equation*}
converges to an element in $\matheur{B}^+_{\mathrm{dR}}$, and it defines
a non degenerate bilinear pairing
\begin{equation*}
 \oint^{(p)}\;:\;
 \mathrm{H}^{1}_{\mathrm{dR}}(A)\otimes
 \mathrm{T}_p(A_{\ovl{K}})\longrightarrow\matheur{B}_{\mathrm{dR}}^+,
\end{equation*}
compatible with the Galois action and the filtrations.
\end{thm}

\begin{rem}
 On composing the period pairing $\int^{(p)}$ with $\theta$, one gets a
 bilinear map $\mathrm{H}^1_{\mathrm{dR}}(A) \otimes T_p(A_{\ovl{K}})
 \longrightarrow \C_p$ which sends $\Omega^1(A)\subset
 \mathrm{H}^1_{\mathrm{dR}}(A)$ to $0$.  This degenerate pairing
 describes ``half'' of the (Hodge)-Tate-Raynaud decomposition
 \begin{equation*}
  \mathrm{H}^1_{\mathrm{et}}(A_{\ovl{K}},{\Q}_p)\otimes {\C}_p \simeq
   \Omega^1(A)\otimes
   {\C}_p(-1) \oplus \mathrm{H}^1(\mathcal{O}(A))\otimes {\C}_p
 \end{equation*}
 of $\Gal({\ovl{K}}/K)$-modules.
\end{rem}

\subsection{Some computations of elliptic $p$-adic periods.}
\label{sub-computation}

\begin{para}{\itshape Tate elliptic curves.\ }
In this case, $T_p(A_{\ovl{K}})$ sits in an exact sequence 
\begin{equation*}
 0\rightarrow {\Z}_p(1)\rightarrow T_p(A_{\ovl{K}}) \rightarrow
  {\Z}_p \rightarrow 0.
\end{equation*}
Let us take for $\gamma_1$ the image of $(2i\pi)_p \in {\Z}_p(1)$ in
$T_p(A_{\ovl{K}})$, and for $\gamma_2$ any lifting of $1 \in {\Z}_p$.
For concreteness, let $A=X_z$ be in Legendre form as in
\ref{para-padic-int-1}.  Let $\eta \in \mathrm{H}^1_{\mathrm{dR}}$
satisfy $\an{\omega,\eta}=1$ (\eg $\frac{xdx}{4y}$).  We set
$\omega^{(p)}_i=\int^{(p)}_{\gamma_i}\omega,
\eta^{(p)}_i=\int^{(p)}_{\gamma_i}\eta$.  Then one has the ``Legendre
relation''
\begin{equation*}
\omega^{(p)}_1\eta^{(p)}_2-\eta^{(p)}_1\omega^{(p)}_2=(2i\pi)_p
\end{equation*}
and the formulas (\cf \cite{pbl}, \cite{get})
\begin{equation*}
\frac{\omega^{(p)}_1}{(2i\pi)_p} = \sqrt{-1}F\Bigl(\frac{1}{2},
\frac{1}{2},1;z\Bigr) = 1/\Theta(z), \;\; \omega^{(p)}_2 =
\sqrt{-1}F\Bigl(\frac{1}{2},\frac{1}{2},1;z\Bigr) \Log q
\end{equation*}
(compare with the complex case in \ref{para-tate-para}, and also
\cite[p.406]{f2}).  Note that the ``period'' we found in
\ref{para-padic-int-1} by integration along the Berkovich loop is
 \emph{not} a Fontaine-Messing period (there is a $\log$ instead of
 $\Log$).
\end{para}

\begin{para}{\itshape Elliptic curves with ordinary reduction.\ }
We now assume that $z\in {\Z}_p^{\mathrm{ur}}$ and that $A=X_z$ has good
ordinary reduction (with the notation of \subsecref{sub-dwork-hyper},
this means that $z(1-z)h_p(z)\neq 0$).  Again, $T_p(A_{{\ovl{\Q}}_p})$
sits in an exact sequence (\cite[A.2.4]{alraec})
\begin{equation*}
0\rightarrow {\Z}_p \gamma_1\rightarrow
T_p(A_{{\ovl{\Q}}_p}) \rightarrow T_p(A\;\mathrm{mod}\, p) \rightarrow 0.
\end{equation*}
 Then one still has the formula $\Theta(z)=
\frac{(2i\pi)_p}{\omega^{(p)}_1}$, where $\Theta(z)\in
\widehat{{\Z}_p^{\mathrm{ur}}}$ is the Tate constant discussed in
\ref{para-tate-para} (compare \cite{pbl}, \cite[4.3]{itoecwcm}).  In
particular, $\frac{\omega^{(p)}_1}{(2i\pi)_p}$ is given by the
evaluation at $z$ of ``the'' extension of
$F(\frac{1}{2},\frac{1}{2},1;?)$ discussed in \subsecref{sub-serre-tate-para}.

\noindent On the other hand, let $\gamma_2 \in T_p(A_{{\ovl{\Q}}_p})$ be
such that the Weil pairing $\an{\gamma_1,\gamma_2}\in {\Z}_p(1)={\Z}_p
(2i\pi)_p$ is the chosen generator $(2i\pi)_p$. Then one has the
``Legendre relation'', and $\omega^{(p)}_2 = \omega^{(p)}_1\cdot\frac{{\Log
\;q}}{{(2i\pi)_p}}$, where $q\in \widehat{{\Z}_p^{\mathrm{ur}}}$ denotes
here the so-called {\it Dwork(-Serre-Tate) parameter} of $A$ introduced
in \ref{para-dwork-para}.  If $q=1$, \ie when $A$ is the \emph{canonical
lifting} of $A\pmod{p}$, then one can choose $\gamma_2\;$ such the
corresponding value of ${\Log 1}$ is $0$, and then $\omega^{(p)}_2 =0.$
\end{para}

\subsection{Periods of CM elliptic curves and values of the Gamma function.}

\begin{para}
 The case of an elliptic curve $A$ with \emph{supersingular reduction} is
 more delicate.  If $A$ does not have complex multiplication (CM), one
 can show, using \cite[A.2.2]{alraec}, that the four basic periods
 $\omega^{(p)}_1,\omega^{(p)}_2,\eta^{(p)}_1,\eta^{(p)}_2$ are
 algebraically independent over ${\ovl{\Q}}_p$, so that one cannot
 expect ``formulas'' for the periods. 

  On the other hand, it is well-known that the complex
 periods of an elliptic curve with complex multiplication can be
 expressed in terms of special values of the $\Gamma$ function. So one
 may ask for a $p$-adic analogue involving $\Gamma_p$
 (\cf~\ref{subsection-deriv-of-gamma}).  Such a formula has been found by A.
 Ogus, actually not for the periods, but for the action of Frobenius.
\end{para}

\begin{para} {\itshape The Lerch-Chowla-Selberg-Ogus formulas.}
 \label{para-lerch-chowla}
 Let $A$ be an elliptic curve with complex multiplication by
 ${\Q}(\sqrt{-d})$ over ${\ovl{\Q}}\,$ ($-d$ denotes a fundamental
 discriminant).  Let $\epsilon = (\frac{-d}{}) $ be the quadratic
 character $({\Z}/d)^\times\to {\Z}/2$ induced by the embedding
${\Q}(\sqrt{-d})\subset {\Q}(\zeta_d) $, $h=$ the class number of
 ${\Q}(\sqrt{-d})$, and $w=$ the 
 number of roots of unity in ${\Q}(\sqrt{-d})$.  For any $u\in
 ({\Z}/d)^{\times}$, we denote by $\an{\frac{u}{d}}$ the unique rational
 number in $]0,1]$ such that $d\an{\frac{u}{d}}\equiv u$.

 Let $v$ be a place of ${\ovl{\Q}}$ of residue characteristic $p$, with
 associated embedding ${\ovl{\Q}}\subset {\ovl{\Q}}_p$.  To any lifting
 $\psi_v$ to $\Gal({\ovl{\Q}}_p/{\Q}_p)$ of the Frobenius element in
 $\Gal({\Q}_p^{\mathrm{ur}}/{\Q}_p)$, one attaches a $\psi_v$-linear
 endomorphism $\Psi_v$ of
 $\mathrm{H}^1_{\mathrm{dR}}(A/{\ovl{\Q}})\otimes {\ovl{\Q}}_p$
 \footnote{using its crystalline interpretation, \cf \cite[4]{fadrc1};
 we refer to \cite{chambert98:_cohom} for a nice recent survey of crystalline
 cohomology.}.
 On the other hand, it is well-known that $A$ has supersingular
 reduction at $v$ if and only if $\epsilon(p)=-1$ or $0$.
\end{para} 

\begin{thm}\label{thm-ogus-formula}
 There exists a basis $(\omega, \eta)$ of
 $\mathrm{H}^1_{\mathrm{dR}}(A/{\ovl{\Q}})$ of eigenvectors under the
 action of $\sqrt{-d}$ ($\omega$ being the class of a regular
 differential), and an element $\gamma \in \mathrm{H}_1(A(\C),{\Q})$,
 such that
 \begin{align*}
  \int_{\gamma} \omega
  & = \sqrt{2i\pi}\prod_{u\in ({\Z}/d)^{\times}}
  \Bigl(\Gamma\Lan{\frac{u}{d}}\Bigr)^{\epsilon(u)w/4h} \\
  \int_{\gamma} \eta 
  & = \sqrt{2i\pi}\prod_{u\in({\Z}/d)^{\times}}
  \Bigl(\Gamma\Lan{-\frac{u}{d}}\Bigr)^{\epsilon(u)w/4h}
 \end{align*}
 and such that for every place of ${\ovl{\Q}}$
 of residue characteristic $p$ satisfying $\epsilon (p)=-1$,
 \begin{align*}
  \Psi_v^{\ast} (\omega)
  & = p\prod_{u\in ({\Z}/d)^{\times}}
  \Bigl(\Gamma_p\Lan{p\frac{u}{d}}\Bigr)^{-\epsilon(u)w/4h}\eta \\
  \Psi_v^{\ast} (\eta)
  & = \prod_{u\in ({\Z}/d)^{\times}}
  \Bigl(\Gamma_p\Lan{-p\frac{u}{d}}\Bigr)^{-\epsilon(u)w/4h}\omega
 \end{align*}
 up to multiplication by some root of unity.
\end{thm}

\noindent This is a concatenation of \cite[3.15, 3.9]{ogus90:_chowl_selber},
taking into account the formula $2hd/w=- \sum_1^d \epsilon(u) u$,
\cf \cite[VII]{vudtdaz}).

\medskip\noindent {\it Remark.} The last two formulas extend to the case when
$\epsilon (p)=+1$ if one interchanges $\omega$ and $\eta$ in the right
hand sides (for $d=4$ and $p\equiv 1\pmod{4}$, this
is compatible with Young's formula in \subsecref{sub-analytic-cont}, since
$\Gamma_p(1/2)^4=1$, and with the formula
$F(\frac{1}{2},\frac{1}{2},1;-1)=\frac{1}{2}
\frac{\Gamma(1/4)}{\Gamma(1/2)\Gamma(3/4)}$.). 

 In that case, the expressions 
 \begin{equation*}
  \prod\Bigl(\Gamma_p\Lan{p\frac{u}{d}}\Bigr)^{-\epsilon(u)w/4h},\;
   \prod\Bigl(\Gamma_p\Lan{p\frac{u}{d}}\Bigr)^{-\epsilon(u)w/4h}
 \end{equation*}
 are algebraic numbers: indeed, let $r$ be the order of the subgroup
 $\an{p}$ of $({\Z}/d)^{\times}$ 
 generated by $p$, and $k= p^r-1/d$; then
\begin{equation*}
 \prod_{u\in
  ({\Z}/d)^{\times}}\Bigl(\Gamma_p\Lan{p\frac{u}{d}}\Bigr)^{-\epsilon(u)}=
  \prod_{w\in ({\Z}/d)^{\times}/\an{p}}\prod_1^{r}
  \Gamma_p\biggl(\frac{p^i wk}{p^r-1}\biggr)^{-\epsilon(w)}
\end{equation*}
and one concludes by the Gross-Koblitz formula.  
 
 \begin{para} {\itshape The ramified case.}\label{para-ramif} In the
  ramified case, \ie when
  $\epsilon (p)=0$, there is again an analogue of the
  last two formulas for
  $p\neq 2$. This relies on
  Coleman's computation of the Frobenius matrix of Fermat curves of degree
  divisible by $p$ --- which have arboreal reduction modulo $p$. The result
  takes the same form as in the case $\epsilon (p)=-1$: $\Psi$ is now
  attached to an element $\psi$ of degree one in the Weil group (\ie
  a lifting of the Frobenius element in
  $\Gal({\ovl{\Q}}_p/{\Q}_p)$), the expressions $\,{p\frac{u}{d}}\,$ have
  to be replaced by $\,{\psi\frac{u}{d}}$, and $\Gamma_p$ must be extended to 
  ${\Q}_p$ 
  \cite[6.5]{coleman90:_froben_fermat}. 
  
  Coleman's complicated expressions have been simplified by F.~Urfels in
  his thesis (Strasbourg, 1998; unpublished). If one passes to the
  Iwasawa logarithm $\log_p$ \index{Iwasawa logarithm@Iwasawa logarithm}
  \index{000logp@$\log_p$}
  (at the cost of rational powers of $p$), the result is that one should
  replace $\log_p \Gamma_p \lan{\psi\frac{u}{d}}$ in Ogus' formula by
  $G_p\bigl(\lan{\psi\frac{u}{d}}\bigr)-G_p\bigl(\lan{\frac{u}{d}}\bigr)$, 
  where $G_p$ denotes Diamond's LogGamma function
  \index{Diamond's LogGamma function@Diamond's LogGamma function}
  \index{000Gp@$G_p$}
  \begin{equation*}
   G_p(x)= \lim_{m\to \infty}\frac{1}{p^m}\sum_{n=0,\ldots, p^m-1}
    (x+n)\log_p(x+n) -(x+n) .
  \end{equation*}
 
\noindent {\it Example.} We take $p=3$. Let $\Q[\sqrt{-3n}]$ be an
 imaginary quadratic field with fundamental discriminant $-3n$, $n$
 prime to $3$. Let  $A$  be an elliptic curve with 
complex multiplication by some order in $\Q[\sqrt{-3n}]$.

Let $\psi$ be an element of the Weil group which acts by $-1$ on
$\Q_3/\Z_3\cong \mu_{3^\infty}$. One has a formula $\Psi_v^{\ast}
(\omega) = \kappa\eta$ with

\begin{equation*}
 \log_3 \kappa =-(\epsilon(u)w/4h)\sum_{u\in ({\Z}/3n)^{\times}}
 \Bigl(G_3\Bigl(\Lan{\psi\frac{u}{3n}}\Bigr) - 
 G_3\Bigl(\Lan{\frac{u}{3n}}\Bigr)\Bigr)
\end{equation*}
 We notice that for $u\in
 ({\Z}/3n)^{\times}$,
 $\epsilon\Bigl(3n\Lan{\psi\frac{u}{3n}}\Bigr)=-\epsilon (u)$, 
 whence
\begin{equation*}
  \log_3 \kappa = (w/2h) \sum_{u\in ({\Z}/3n)^{\times}}
 \epsilon(u)G_3\Bigl(\Lan{\frac{u}{3n}}\Bigr).
\end{equation*}

 On the other hand, for $p=3$,
the Teichm\"uller character is nothing but the Legendre symbol
$(\frac{-3}{{\,}})$. Using this, the latter expression can be rewritten,
by the Ferrero-Greenberg formula (\cf \cite[chap. 17]{cf1a2}), as
\begin{align*}
 \log_3\, \kappa = (w/2h) . (L'_3(0,\epsilon) +
 L_3(0,\epsilon)\log_3n)\\ =  (w/2h ) .\sum_{v\in ({\Z}/n)^{\times}} 
 \Bigl(\frac{n}{{v}}\Bigr)\log_3\Gamma_3 \Lan{\frac{v}{n}}.  
\end{align*} 
 When $(\frac{n}{{3}})=1$, then
$L_3(0,\epsilon)= 0$ and the Gross-Koblitz formula shows that
$L'_3(0,\epsilon)= \sum 
  \bigl(\frac{n}{{v}}\bigr)\log_3\Gamma_3 \lan{ \frac{v}{n}} $ is the
 Iwasawa logarithm of an algebraic number (a Gauss sum, {\it \cf
 loc.~cit.}), hence $\kappa\in \ovl{\Q}$.

 \noindent  A particularly simple case is $n=1$: one finds $\log_3
 \kappa=0$, thus $\kappa
 \in 3^\Q.\mu_\infty \subset \ovl{\Q}$.

 \medskip
 \noindent For $n=8$, a contrario, one has $(\frac{n}{{3}})=-11$: one finds  
 \begin{equation*} \log_3 \kappa = 
  (1/2).\Bigl[\log_3\Gamma_3 \Lan{\frac{1}{8}}-\log_3\Gamma_3
   \Lan{\frac{3}{8}}-\log_3\Gamma_3
   \Lan{\frac{5}{8}}+\log_3\Gamma_3 \Lan{\frac{7}{8}}\Bigr]=0
 \end{equation*}
 by the functional equation of
 $\Gamma_3$, so that again $\kappa
 \in 3^\Q.\mu_\infty \subset \ovl{\Q}$.
\end{para}

\begin{para} {\itshape Colmez' product formula.} There is a natural
 extension of $|\quad|_p$ on ${\ovl{\Q}}_p$ to
 $\matheur{B}^+_{\mathrm{dR}}$ (however not as an absolute value), such
 that $\vert (2i\pi)_p \vert_p = p^{-\frac{1}{p-1}}$
 \cf \cite{ufdpplpdvaamc}. 
 
 Colmez has remarked that the logarithm of the
 product $|2i\pi| \prod_p |(2i\pi)_p|_p $ $=$ $ 2\pi\prod_p
 p^{-\frac{1}{p-1}}$ 
 is formally equal to $\log (2\pi) + \zeta^{\prime}(1)/\zeta(1)$, a
 divergent sum which can be renormalized using the functional equation of
 $\zeta$: setting $\zeta^{\prime}(1)/\zeta(1)$ $=$ 
 $-\zeta^{\prime}(0)/\zeta(0)=-\log(2\pi)$, the renormalized product is
 \begin{equation*}
  |2i\pi| \sideset{}{'}\prod_p |(2i\pi)_p|_p = 1.
 \end{equation*}
 He has given an amazing generalization of this product formula to
 periods of CM elliptic curves and many other CM abelian varieties
 (loc.~cit.).
 \end{para}

\newpage
\section{Periods as solutions of the Gauss-Manin connection.}

\begin{abst}
 The variation of the periods in a family of complex algebraic varieties
 is controlled by the Picard-Fuchs differential equation. We discuss the
 p-adic situation, present Dwork's general viewpoint on p-adic
 periods. We then discuss the question of the existence of an arithmetic
 structure on the space of solutions, which would give an intrinsic
 meaning to the period modulo $\ovl{\Q}^\times$. We present such an
 arithmetic structure (analogous to the Betti lattices in the complex
 situation) in the case of abelian varieties with either multiplicative
 reduction or supersingular reduction. In the latter case, we relate the
 periods, in the presence of complex multiplication, to special values
 of the $p$-adic Gamma function.
 \end{abst}

\bigskip
\subsection{Stokes.}

We come back to the basic question of the meaning of \emph{integration
over a loop}, already discussed in \subsecref{sub-prolegomena}.  A
different approach is based on the Stokes lemma: integrating \emph{exact}
differentials on a loop gives zero. This is the approach favored by
K. Aomoto in the complex case and by B. Dwork in the $p$-adic case.

\begin{para}
 In order to see how this idea can be implemented, let us consider the
 Hankel expression for the gamma function
 \begin{equation*}
  \frac{1}{\Gamma(1-\alpha)} = \frac{1}{2i\pi}
   \int_{\gamma} x^{\alpha}e^x\frac{dx}{x}
 \end{equation*}
 where $\gamma$ is the following loop based at $-\infty$.

 \begin{figure}[h]
  \begin{picture}(200,40)(0,0)
   \put(0,0){\includegraphics{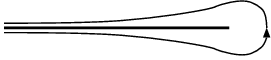}}
   \put(108,12.5){$\bbp$}
   \put(117,10){$\scriptstyle{0}$}
   \put(135,15){$\gamma$}
  \end{picture}
  \caption{}
 \end{figure}

 \noindent Stokes' lemma suggests to attach to this integral the
 following complex
 \begin{equation*}
 x^{\alpha}e^x\mathcal{O} \xrightarrow{xd/dx}  x^{\alpha}e^x\mathcal{O}
 \end{equation*}
 where $\mathcal{O}$ is a suitable ring of analytic functions.  From the
 point of view of index theory (Malgrange-Ramis), a natural choice is
 \begin{equation*}
 \mathcal{O}= {\C}[[x]]_{-1,1^-}=\biggl\{
 \sum_{n\geq 0} a_n x^n \biggm| \exists \kappa >0, \exists r \in ]0,1[,\;
 |a_n| \leq \kappa r^n/n!\biggr\}
 \end{equation*}
 identified with the ring of entire functions of {\it exponential order}
 $O(e^{r|x|})$ for some $r<1$\footnote{we use the traditional
 notation for rings of Gevrey series; the value $(-1,1^-)$ of the index
 is ``characteristic'', \ie an extremal value for which $\mathrm{H}^1$
 is of dimension one.}.  Using the formula $x^{\alpha}e^x x^k=
 \frac{-1}{\alpha + k}x^{\alpha}e^x x^{k+1} + xd/dx(\frac{1}{\alpha +
 k}x^{\alpha}e^x x^k)$, a simple computation then shows that
 $\mathrm{H}^0=0$ and that $\mathrm{H}^1$ is of dimension $1$; it is
 generated by $[\frac{1}{2i\pi}x^{\alpha}e^x]$, whose integration along
 $\gamma$ gives $\frac{1}{\Gamma(1-\alpha)}$.

 \par \noindent The $p$-adic analogue
 \begin{equation*}
 \C_p[[x]]_{-1,1^-}=\biggl\{
 \sum_{n\geq 0} a_n x^n \biggm|\exists \kappa >0, \exists r \in ]0,1[,\;
 |a_n|\leq \kappa r^n/| n!|_p \biggr\}
 \end{equation*}
 is nothing but the ring of \emph{overconvergent} analytic functions on
 $\matheur{D}(0,|\pi|_p)$ (here $\pi$ is Dwork's constant).  By change
 of variable $z=\pi x$, we thus recover the complex
 $M_{\alpha}^{\dagger} \xrightarrow{zd/dz} M_{\alpha}^{\dagger}$ of
 \ref{subsection-deriv-of-gamma}, where $\Gamma_p(\alpha)$ $\bigl(=\pm
 \frac{1}{{\Gamma_p(1-\alpha)}}\bigr)$ appeared.
\end{para}

\begin{para}
 This approach inspired by Stokes' lemma is well-suited for studying
 many kinds of hypergeometric functions (confluent or not). Basically,
 the integral will satisfy difference equations with respect to the
 exponents ($\alpha$ in the previous example), and differential
 algebraically on parameters. In the complex situation, one of the main
 problems is the construction of nice loops (such as $\gamma$); this is
 the object of a so-called ``topological intersection theory''
 generalizing the usual Betti homology and period pairing. In the
 $p$-adic case, the focus has been more on the construction and
 properties of the Frobenius structure, in particular its analyticity
 with respect to the exponents (``Boyarsky principle'' \cite{otbp}).
\end{para}

\begin{para}
 When no exponential is involved and when the exponents are rational,
 the integrand is algebraic; if it depends algebraically on parameters,
  the period integral is a solution of a {\it Gauss-Manin connection}. 

 Let us say a few words about the classical case of
 $F(a,b,c;z)$.  We assume that $a,b,c \in {\Z}_p$, and that $c-a,c-b,b,a$
 all lie outside $\Z$.  We set
 \begin{equation*}
 f_{a,b,c;z}=
 x^b(1-x)^{c-b}(1-zx)^{-a}
 \end{equation*}
 where, for simplicity, $z$ is limited to $\vert z(1-z)\vert = 1$.  Let
 $Z$ be the affinoid $\matheur{D}(0,1^+)\setminus
 (\matheur{D}(0,1^-)\cup \matheur{D}(1,1^-)\cup \matheur{D}(1/z,1^-))$.
 The relevant complex is
 \begin{equation*}
 f_{a,b,c;z}\mathcal{H}^{\dagger}(Z)
 \xrightarrow{xd/dx} f_{a,b,c;z}.\mathcal{H}^{\dagger}(Z).
 \end{equation*}
 Its $\mathrm{H}^1=\mathrm{H}^1_{a,b,c}$ is of dimension $2$; it is
 generated by $[f_{a,b,c;z}]$ and $[\frac{f_{a,b,c;z}}{1-x}]$. The
 Gauss-Manin connection is
 \begin{equation*}
 \nabla\Bigl(\frac{d}{dz}\Bigr)
  \begin{pmatrix}
   [f_{a,b,c;z}] \\ {}
   [\frac{f_{a,b,c;z}}{1-x}]
  \end{pmatrix}
 =
 \begin{pmatrix}
  -\frac{c}{z} & \frac{c-b}{z} \\ 
  \frac{c-a}{1-z} & \frac{a+b-c}{1-z}
 \end{pmatrix}
 \begin{pmatrix}
  [f_{a,b,c;z}] \\ {}
  [\frac{f_{a,b,c;z}}{1-x}]
 \end{pmatrix}
 \end{equation*}
  \cf \cite[3.1, 1.2]{lopde}.  The Frobenius structure relates
 $\mathrm{H}^1_{a,b,c}$ to
 $\mathrm{H}^1_{a^{\prime},b^{\prime},c^{\prime}}$, where
 $a^{\prime},b^{\prime},c^{\prime}$ are the successors of $a,b,c$
 respectively (\cf~\ref{subsection-deriv-of-gamma}).  For all these
 hypergeometric series, there is a story similar to that of
 $F(\frac{1}{2},\frac{1}{2},1;z)$.  In this context, it is fruitful to
 combine the Boyarsky principle and the contiguity relations.
\end{para}

\begin{para}
 In this sketch of Dwork's approach, ``periods''
 are just analytic solutions of the Gauss-Manin connection. One can ask
 more: namely, one can ask for a dual theory of $p$-adic cycles and a
 period pairing as in the complex case\footnote{Dwork has developed a
 ``dual theory'' and a pairing given by residues, which play somehow the
 role of ``$p$-adic cycles'' (\cf \cite[5]{uinacdd}); however, this dual
 theory is an avatar of de Rham cohomology with proper supports and
 does not enjoy the properties of discreteness and horizontality that we
 are looking for.}.  At best, we can expect these ``{\it $p$-adic Betti
 lattices}'' to be locally horizontal, functorial with respect to the
 endomorphisms of the geometric fibres, and defined over $\Z$ or at least
 over some number field.

 It turns out that this can be done in some
 cases, \eg in the cases of abelian varieties with multiplicative
 reduction, and of abelian
 varieties with supersingular reduction (\cite{pbl}, \cite{tdmeigdvpdg}). Let
us outline the results.
\end{para}

\subsection{$p$-adic Betti lattices for abelian varieties with 
multiplicative reduction.}\label{sub-p-adic-betti-lattice}

\begin{para}\label{para-multi-uniformize}
 {\itshape Over $\C$.}
 We first recall the ``multiplicative uniformization'' of complex
 abelian varieties.  Let $A$ be an abelian variety of dimension $g$ over
 $\C$.  One has the analytic representation
 \begin{equation*}
  A(\C)=T(\C)/M,
 \end{equation*}
 where $T$ is a torus of dimension $g$, and $M$ is a lattice of rank
 $g$.  Set $M' :=\mathrm{Hom}(T,\G_m)\,(=M_{A^{\vee}})$.  Then we
 have the exact sequence
 \begin{equation*}
  0\longrightarrow 2i\pi {M'}^{\vee}\longrightarrow\Lambda\longrightarrow
   M\longrightarrow 0,
 \end{equation*}
 where $\Lambda$ is the period lattice (of rank $2g$), and this sequence
 splits by choosing a branch of $\log$.  The inclusion $M\hookrightarrow
 T$ defines a pairing $M\times M' \rightarrow\C^{\times}$.  A
 polarization, on the other hand, induces $M\rightarrow M'$.  These
 data give rise to the pairing ({\it ``multiplicative period''})
\begin{equation*}
 q=(q_{ij})\colon M\otimes M\longrightarrow\C^{\times}
\end{equation*}
and $-\log|q|$ is a scalar product on $M_{\R}$.
\end{para}

\begin{para}\label{para-degenerate}
 {\itshape Over the $p$-adics.}
 Let $K$ be a $p$-adic local field
 and $A$ an abelian variety over $K$ having (split) multiplicative
 reduction.  One has a similar analytic representation \cite{aacodavocr}
 \begin{equation*}
  A(K)=T(K)/M
 \end{equation*}
 Set again $M' :=\Hom(T,\G_m)$.  We have the following identifications
 \begin{equation*}
  M\simeq \mathrm{H}_1(A^{\mathrm{rig}},\Z)\quad\mathrm{and}\quad
   M'\simeq \mathrm{H}_1(A^{\vee\,\mathrm{rig}},\Z),
 \end{equation*}
 where $A^{\mathrm{rig}}$ denotes the associated rigid analytic variety
 over $K$ to $A$.  Then, just as in the complex case, one can 
 construct a natural exact sequence:
 \begin{equation*}
  0\longrightarrow (2i\pi)_p{M'}^{\vee}\longrightarrow\Lambda
   \longrightarrow M\longrightarrow 0,
 \end{equation*}
 split by the choice of a branch of the $p$-adic logarithm $\log_p$, 
 and a non-degenerate pairing
 \begin{equation*}
 \int^p\colon\mathrm{H}^1_{\mathrm{dR}}(A)\times\Lambda\longrightarrow
 K[(2i\pi)_p].
 \end{equation*}
 (The construction involves one-motives and Frobenius \cite{pbl}).
 Restricted to $(2i\pi)_p{M'}^{\vee}$ which is in a 
 canonical way a subgroup of $T_p(A_{\ovl{K}})$, this pairing coincides
 with the Fontaine-Messing pairing.
\end{para}

\begin{para}
 {\itshape Degenerating abelian pencils.}  Let $\mathcal{A}\rightarrow
 {\A}^1\setminus \{0,\zeta_1,\ldots,\zeta_r\}$ be a pencil of polarized
 abelian varieties defined over a number field $k$. We consider an
 element $\omega$ of the relative $\mathrm{H}^1_{\mathrm{dR}}$ and its
 associated Picard-Fuchs differential equation $\mathcal{L}\omega=0$. We
 assume that the connected component of identity of the special fibre at
 $0$ of the Neron model of $\mathcal{A}$ is a split torus $T$ over $k$.
 The dual abelian pencil then has the same property, with a torus
 $T^{\prime}$.  We set $M=\Hom(T',{\G}_m)$, $M'=\Hom(T,{\G}_m)$.

 Let us fix $k\hookrightarrow \C$.  The constant subsheaf of
 $R_1(f^{\mathrm{an}}_{\C})_{\ast}{\Z}$ in the neighborhood of $0$
 identifies with $2i\pi M^{\prime \vee}$: its fibre at any $z (\neq
 0,\zeta_i)$ is the lattice $2i\pi M_z^{\prime \vee}$ attached to
 $\mathcal{A}_z$ (\ref{para-multi-uniformize}).  Similarly for $M$.  On
 the other hand, the choice of a branch of $\log$ identifies
 $\mathrm{H}_{1\mathrm{B}}(\mathcal{A}_z, {\Z})\simeq
 (R_1(f^{\mathrm{an}}_\C)_{\ast}{\Z})_z$ with $\Lambda = 2i\pi
 M_z^{\prime \vee}\oplus M_z$.  Then for any $\gamma \in 2i\pi M^{\prime
 \vee}$, $y(z)= \frac{1}{2i\pi}\int_{\gamma_z}\omega_z$ is a solution of
 $\mathcal{L}$ in $k[[z]]$.

 More generally, for any $\gamma^{\ast} \in \Lambda$, $\frac{1}{2i\pi}
 \int_{\gamma^{\ast}_z} \omega_z $ is a solution of the Picard-Fuchs
 differential equation of the form $y(z) \frac{\log az^n}{2i\pi} +
 y^{\ast}(z)$, with $y^{\ast}(z)\in k[[z]]$, $a\in k$, $n\in \Z$.  (In
 the case of the Legendre pencil, we have already met this situation
 with $y(z)=iF(\frac{1}{2},\frac{1}{2},1;z), y^{\ast}(z)= 2
 F^{\ast}(\frac{1}{2},\frac{1}{2},1;z), a= 2^{-8}, n = 2$).

 Let now fix an embedding $k\hookrightarrow K \subset \C_p$. For any $z
 \neq 0$ close enough to $0$, $M$ and $M^{\prime}$ identify respectively
 with the lattices $M_z$, $M_z^{\prime}$ attached to $\mathcal{A}_z$
 (\ref{para-degenerate}), and the choice of a determination of the
 $p$-adic logarithm identifies $\Lambda$ with $\Lambda_z$. It turns out
 that the $p$-adic pairing $\int^p$ of \ref{para-degenerate} extends to
 a horizontal pairing over a punctured disk around $0$, which is given
 by the ``same'' formula as in the complex case
 \begin{equation*}
 \frac{1}{(2i\pi)_p} \int^p_{\gamma^{\ast}_z}
 \omega_z = y(z)\frac{\log_p az^n}{(2i\pi)_p} + y^{\ast}(z),
 \end{equation*}
 where $y(z), y^{\ast}(z)$ are the \emph{same formal series as above,
 evaluated $p$-adically} (loc. cit.). The computation of
 \ref{para-padic-int-1} can be considered as a special case.

 \noindent This is also closely related to the work of T. Ichikawa on
 ``universal periods'' for Mumford curves
 \cite{sutorsamcoig}.
\end{para}

\begin{para}
 {\itshape Relation to the Fontaine-Messing periods.} Since the restriction
 of $\int^p$ to $(2i\pi)_p M_z^{\prime \vee}\subset
 T_p(\mathcal{A}_{z,\ovl{K}})$ is the Fontaine-Messing pairing, we
 obtain $(2i\pi)_p y(z)$ as a Fontaine-Messing period of
 $\mathcal{A}_z$. More generally, $y(z)\Log_p az^n +(2i\pi)_p
 y^{\ast}(z)$ appears as a period ($\Log$ instead of $\log$). We conclude
 that in this degenerating case, the Fontaine-Messing periods behave
 relatively well with respect to the Gauss-Manin connection. A similar
 observation applies to the pairing composed with $\theta$ (Hodge-Tate
 periods).
\end{para}

\subsection{$p$-adic Betti lattices for abelian varieties with
supersingular reduction.}\label{sub-p-adic-betti-lattice-supersing}

\begin{para}\label{para-ss-abelian}
 Any supersingular abelian variety over $\ovl{\F}_p$ is isogenous to
 a power of a supersingular elliptic curve. We refer to
 \cite[4]{waterhouse69:_abelian} for a detailed study of supersingular
 elliptic curves over finite fields (including the classification of
 isogeny classes and isomorphism classes, refining Deuring's classical
 work), and to \cite{li_o98:_modul} for the higher dimensional case. Let
 us simply recall a few basic facts:
 \begin{itemize}
  \item there is exactly one ${\ovl{\F}_p}$-isogeny class of
	supersingular elliptic curves.
  \item If $A_0$ is any supersingular elliptic curve over
	${\F}_{p^n}$, then $\bcD:= \End(A_0)_{\ovl{\F}_p}$ 
	\index{000D@$\bcD$}
	is a maximal order in a quaternion algebra over $\Q$, ramified
	exactly at $p,\infty$; we set $\bD=\bcD_{\Q}$.\index{000D@$\bD$}
  \item $A_0$ is ${\ovl{\F}_p}$-isomorphic (but not necessarily
	${\F}_{p^n}$-isomorphic) to an elliptic curve defined over
	${\F}_{p^2}$.
  \item If $n=1$, $\End(A_0)\otimes {\Q}={\Q}[Fr_{p}]$, an
	imaginary quadratic field; if $p\neq 2$, this is
	${\Q}[\sqrt{-p}]$. More precisely, the order $\End(A_0)$ can be either 
	${\Z}[Fr_{p}]$ or the maximal order in $\End(A_0)\otimes
	{\Q}$ (which coincide if $p\equiv 3\pmod{4}$). If $p=2$, there is
	another possibility, namely $\End(A_0)={\Z}[Fr_2]\simeq
	{\Z}[\sqrt{-1}]$.
  \item There exists a supersingular elliptic
	curve over ${\F}_p$ whose Frobenius endomorphism satisfies
	$(Fr_p)^2=-p.\id\,$ and whose endomorphism algebra is the maximal
	order in ${\Q}[Fr_{p}]$. Such elliptic curves belong to a single
	${\F}_p$-isogeny class.
  \item If $p\neq 2,3$, there is exactly one ${\F}_p$-isogeny class of
	supersingular elliptic curves over ${\F}_p$ (but several
	${\ovl{\F}_p}$-isomorphism classes in general).
 \end{itemize}
\end{para}

\begin{para}\label{para-supersing}
 Let $A_0$ be a supersingular abelian variety of 
 dimension $g$ over $\F_{p^n}$. By functoriality, the elements of 
 $\End\bigl((A_0)_{\ovl{\F}_p}\bigr)\otimes \Q\simeq M_g(\bD)$ act
 linearly on the crystalline cohomology 
 ${\mathrm{H}^{\ast}_{\mathrm{crys}}\bigl((A_0)_{\ovl{\F}_p}/\ovl{\Q}_p\bigr)}=
 \mathrm{H}^{\ast}_{\mathrm{crys}}\bigl((A_0)_{\ovl{\F}_p}/
 \widehat{\Z_p^{\mathrm{ur}}}\bigr)\otimes_{\widehat{\Z_p^{\mathrm{ur}}}}
 \ovl{\Q}_p$. 
 In degree one, this provides an embedding of $M_g(\bD)$ into the
 endomorphism ring of 
 $\mathrm{H}^1_{\mathrm{crys}}\bigl((A_0)_{\ovl{\F}_p}/\ovl{\Q}_p\bigr)$.

 On the other hand, one has a canonical isomorphism
 \begin{equation*}
  \bigwedge^2
  \mathrm{H}^1_{\mathrm{crys}}\bigl((A_0)_{\ovl{\F}_p}/\ovl{\Q}_p\bigr) \simeq
  \mathrm{H}^2_{\mathrm{crys}}\bigl((A_0)_{\ovl{\F}_p}/\ovl{\Q}_p\bigr) 
 \end{equation*} 
 and a canonical $\Q$-structure in $
  \mathrm{H}^2_{\mathrm{crys}}\bigl((A_0)_{\ovl{\F}_p}/\ovl{\Q}_p\bigr)
  \simeq 
  \ovl{\Q}_p^{g(2g-1)}$ coming from the fact that the whole cohomology
 in degree $2$ is generated by algebraic cycles (since
 $(A_0)_{\ovl{\F}_p}$ is isogenous to the $g$th power of a
 supersingular elliptic curve, one reduces easily to the case $g=2$, in
 which case this fact is well-known).

\medskip Let $F$ be either $\ovl{\Q}$ or a splitting number field for $\bD$: 
 $\bD\otimes_{\Q} F\simeq
 M_2(F)$.  In the latter case, this amounts to saying that $F$ is 
 totally imaginary, and for any place $v$ above $p$,
 $[F_v:\Q_p]$ is even.  We fix an embedding of $F$ in $\ovl{\Q}_p$.  Via 
 this embedding, the $M_g(\bD)$-action extends to a
 $F$-linear $M_{2g}(F)$-action on
 $\mathrm{H}^1_{\mathrm{crys}}\bigl((A_0)_{\ovl{\F}_p}/\ovl{\Q}_p\bigr)$.

 \noindent We consider an irreducible cyclic sub-$M_{2g}(F)$-module 
 $M_{2g}(F).u\subset \mathrm{H}^1_{\mathrm{crys}}\bigl((A_0)_{\ovl{\F}_p}/\ovl{\Q}_p\bigr)$.
 Obviously, such submodules exist, and have $F$-dimension $2g$.  We may and shall choose $u$ in such a 
 way that its exterior square of this $F$-space coincides with the canonical $F$-structure on $
\mathrm{H}^2_{\mathrm{crys}}\bigl((A_0)_{\ovl{\F}_p}/\ovl{\Q}_p\bigr) $ (this can be achieved by replacing $u$ by a
 suitable multiple).

 \begin{pro}\label{betti-la} Up to a homothety by a factor in 
 $\sqrt{F^\times}$, the normalized
 $M_{2g}(F)$-submodule $\;M_{2g}(F)\cdot u\subset
 \mathrm{H}^1_{\mathrm{crys}}\bigl((A_0)_{\ovl{\F}_p}/\ovl{\Q}_p\bigr)$
 depends only on $F\subset \ovl{\Q}_p$, not on the choice
 of $u$. 
  In particular, for $F=\ovl{\Q}$, this defines a canonical $\ovl{\Q}$-structure in $
\mathrm{H}^1_{\mathrm{crys}}\bigl((A_0)_{\ovl{\F}_p}/\ovl{\Q}_p\bigr)$, stable under 
$\End\bigl((A_0)_{\ovl{\F}_p}\bigr)$.
\end{pro}
 
 \begin{proof} Indeed, two such $M_{2g}(F)$-submodules are related 
 by some $h\in GL_{2g}(\ovl{\Q}_p)$, such that
 $\bigwedge^2(h)\in GL_{g(2g-1)}(F)$. Now 
 $M_{2g}(F)\cdot hu = h(M_{2g}(F)\cdot u)$ implies
 that $h$ normalizes
 $GL_{2g}(F)$. It follows that the image of $h$ in 
 $PGL_{2g}(\ovl{\Q}_p)$ lies in $PGL_{2g}(F)$.

  \noindent (When $g>1$, the proposition is not surprising, since
  $\bigwedge^2 V$ is a faithful representation of $PSL(V)$
  for any space $V$ of dimension $>1$).
 \end{proof}

\end{para}
\begin{rems}
 \begin{enumerate}
  \renewcommand{\theenumi}{\alph{enumi}}
  \item A minimal choice for $F$ is a splitting
	quadratic field $\Q(\sqrt{-d})\,$ (the splitting property amounts to
	saying that $d>0$ and $p$ ramifies  or remains prime in
	$\Q(\sqrt{-d})$. A natural choice is $d=p$).
  \item Let $\mathcal{O}_F$ be the ring of integers of $F$. Then
	$\End((A_0)_{\ovl{\F}_p})\otimes \mathcal{O}_F$ is an order in
	$M_{2g}(F)$, and there is a full
	$(\End((A_0)_{\ovl{\F}_p})\otimes \mathcal{O}_F)$-lattice in 
	$(M_{2g}(F)\cdot u)$. As an $\mathcal{O}_F$-module, it is
	projective of rank $2g$.  If $F=\ovl{\Q}$, \ie
	$\mathcal{O}_F=\ovl{\Z}$, then this module is 
	{\it canonically determined} by the
	normalization of its top exterior power, taking into account the 
	canonical isomorphism
	\begin{equation*}
	 \bigwedge^{2g}
	  \mathrm{H}^1_{\mathrm{crys}}\bigl((A_0)_{\ovl{\F}_p}/\ovl{\Z}_p\bigr)
	  \simeq
	  \mathrm{H}^{2g}_{\mathrm{crys}}
	  \bigl((A_0)_{\ovl{\F}_p}/\ovl{\Z}_p\bigr)
	  \simeq \ovl{\Z}_p,
	\end{equation*} 
	the latter isomorphism coming from the trace map. 
 
	Hence there is a 
	{\it canonical $\ovl{\Z}$-structure} in the
	crystalline cohomology of supersingular abelian 
	varieties (not used in the sequel).
 \end{enumerate}
\end{rems}

\begin{para}
 Let $A$ be an abelian variety defined over a 
 subfield $k\subset \C_p$, which has good
 {\it supersingular reduction} $A_0$ over the residue field of $k$. 
 There is a canonical embedding
 $\End(A)\hookrightarrow  \End((A_0)_{\ovl{\F}_p})$. 

The canonical $\ovl{\Q}$-subspace of $\; \mathrm{H}^1_{\mathrm{crys}}((A_0)_{\ovl{\F}_p}/\ovl{\Q}_p)\;$ defined in
\ref{betti-la} will be denoted by 
\begin{equation*}
  \mathrm{H}^1_{\mathrm{B}}(A_{\C_p}, \ovl{\Q}) \subset 
  \mathrm{H}^1_{\mathrm{crys}}\bigl((A_0)_{\ovl{\F}_p}/\ovl{\Q}_p\bigr) 
 \end{equation*}
because, as we shall see, it shares many properties with the Betti space
 $\mathrm{H}^1_{\mathrm{B}}(A_\C, \ovl{\Q})= \mathrm{H}^1_{\mathrm{B}}(A_\C,\Z)\otimes \ovl{\Q}$ of a
complex
 abelian variety $A_\C$). Its dual will be denoted by $ \mathrm{H}_{1,\mathrm{B}}(A_{\C_p}, \ovl{\Q}) $.

Similarly, we shall write $\mathrm{H}^1_{\mathrm{B}}(A_{\C_p}, F)$ for the $F$-subspace described \ref{betti-la}
(well-defined up to a homothety in $\sqrt{F^\times}$), and $ \mathrm{H}_{1,\mathrm{B}}(A_{\C_p},F) $ for its $F$-dual.
 
\medskip  Assume that $k\subset \ovl{\Q}_p$.  There is a functorial 
 isomorphism \cite{fadrc1}
 \begin{equation*}
  \mathrm{H}^1_{\mathrm{dR}}(A)\otimes_k\ovl{\Q}_p \simeq
 \mathrm{H}^1_{\mathrm{crys}}\bigl((A_0)_{\ovl{\F}_p}/\ovl{\Q}_p\bigr),
 \end{equation*} 
 from which one derives a canonical isomorphism
 \begin{equation*}
  \mathrm{H}^1_{\mathrm{dR}}(A)\otimes_k\ovl{\Q}_p \simeq 
  \mathrm{H}^1_{\mathrm{B}}(A_{\C_p},\ovl{\Q})\otimes_{\ovl{\Q}}\ovl{\Q}_p,
 \end{equation*} 
 which is functorial with respect to homomorphisms 
 of abelian varieties with supersingular reduction

 \noindent This isomorphism, which is far from being tautological, as we shall see, 
 can be translated into a $\ovl{\Q}_p$-valued pairing between 
 $\mathrm{H}_{1\mathrm{B}}(A_{\C_p}, F)$ and $\mathrm{H}^1_{\mathrm{dR}}(A)$. 
This pairing is conveniently expressed in the form of a {\it ``period 
 matrix''} $\Omega=(\omega_{ij})$ with
 entries in $\ovl{\Q}_p$, depending on the choice of a basis 
 $(\gamma_i)_{i=1,\ldots ,2g}$ of $\mathrm{H}_{1\mathrm{B}}(A_{\C_p}, F)$
 and a basis $(\omega_j)_{j=1,\ldots ,2g}$ of 
 $\mathrm{H}^1_{\mathrm{dR}}(A)$ (if $A$ is principally polarized, we
 choose symplectic bases). 

\noindent Warning: these $p$-adic periods attached to abelian varieties with supersingular reduction have little to
do with Fontaine-Messing periods  (a motivic interpretation of these periods is proposed in
\cite{tdmeigdvpdg}). 

\end{para}

\begin{para}
 {\itshape Horizontality.} If we let $A$ move in a family $A_z$,
 $\mathrm{H}^1_{\mathrm{crys}}\bigl((A_0)_{\ovl{\F}_p}/\ovl{\Q}_p\bigr)$ may be
 identified with a space of horizontal sections of the de Rham
 cohomology localized in the disk parameterizing the liftings $A_z$ of
 $(A_0)_{\ovl{\F}_p}$ \cite{fadrc1}. Thus the period matrix $\Omega(z)$ is a fundamental
 matrix of solutions of the Gauss-Manin connection.
\end{para}

In the sequel of this subsection, we assume for simplicity that 
 $g=1$, \ie $A$ is an elliptic curve.

\begin{para} {\itshape CM periods, $\,\Gamma_p$-values, and transcendence.}

 \medskip\noindent $\bullet$  Let $A_0$ be as before a
 supersingular elliptic curve over ${\F}_{p^n}$. Let $A$ be a CM-lifting of
 $A_0$, \ie an elliptic curve defined over some number field $k \subset
 \ovl{\Q}\subset \ovl{\Q}_p$ with residue field ${\F}_{p^n}$, such that
 $E= \End(A) \otimes \Q$ is
 quadratic. Note that $E$ is ipso facto a subfield of $\bD$. We denote by
 $c$ the conductor of the
 order $\End(A)$ (in particular, $ce \in \End(A)$).

 \medskip\noindent $\bullet$ Let $(\gamma_1,\gamma_2)$ be any symplectic
 basis of $\mathrm{H}_{1\mathrm{B}}(A,F)$. We can write
 $ce.\gamma_1=a\gamma_1+b\gamma_2$ with $a,b \in F$, so that
 \begin{equation*}
 \tau := \frac{\omega_{12}}{\omega_{11}}=
 \frac{-a+c\sqrt{-d}}{b}\in {\P}^1(FE).
 \end{equation*}
 In particular, $\tau \in {\P}^1({\ovl{\Q}})$.
 We conjecture the following converse, which is a $p$-adic analogue of
 T.~Schneider's theorem \cite{auei} (in the complex case):

 \begin{coj}\label{coj-complex-mult}
  Assume that $A$ is a lifting of $A_0$ defined over $\ovl{\Q}$, and
  that $\tau = \frac{\omega_{12}}{\omega_{11}} \in {\P}^1({\ovl{\Q}})$.
  Then $A$ has complex multiplication.
 \end{coj}

\medskip\noindent $\bullet$ 
 Let $E= {\Q}[\sqrt{-d}]$ and 
 $E'= {\Q}[\sqrt{-d'}]$ be two distinct subfields of
 $\bD\;$ ($d, d'$
 squarefree; we do not exclude the case $d=d'$). According to
 \cite[II 1.4]{adadq}, $E\cap \bcD=\OO_E,\;E'\cap \bcD=\OO_{E'}$. The
 elements  $e = \sqrt{-d}
 \in E$ and $e'=\sqrt{-d'} \in  E'$ generate the
 algebra $\bD$. We remark
 that
 $ee'+e'e$ commutes to $e$ and $e'$, hence is an
 integer $m \in \Z$ (an
 even integer if $d$ or $d'\equiv 3\pmod{4}$, which is divisible by
 $4$ if both $d,
 d'\equiv 3\pmod{4}$). We have
 $(ee'-e'e)^2=m^2-4dd'$. Since
 $ee'-e'e$ anticommutes with $e$, it is not an integer,
 hence $\vert m\vert <
 2\sqrt{dd'}$. On the other hand, the images of $e$ and $e'$
 in $\bcD\otimes
 {\F}_p$ commute ($ee'-e'e$ acts trivially on regular
 differentials in
 characteristic $p$), therefore $p$ divides $4dd'-m^2$. 
 
 \medskip\noindent $\bullet$ On the other hand, there exists a unique symplectic basis
 of eigenvectors for $e$ in $\mathrm{H}^1_{\mathrm{dR}}(A)$ of the form
 $(\omega_1,\omega_2+\sigma \omega_1)$, with
 $\sigma \in k$. We get the relation
 \begin{equation*}
  \omega_{22}+\sigma\omega_{12} = {\til{\tau}}
   (\omega_{21} + \sigma\omega_{11}), \text{ with } {\til{\tau}} =
   \frac{-a-c\sqrt{-d}}{b}.
 \end{equation*} 
 Let now $E'={\Q}[\sqrt{-d'}]\neq E$ be another subfield of $\bD$,
 and assume that $(\gamma_1,\gamma_2)$ is a basis of eigenvectors for
 $E'$: $e'.\gamma_1=\sqrt{-d'}\gamma_1$. We get
 $(ee^{\prime}+e^{\prime}e)\gamma_1=\frac{2a\sqrt{-d^{\prime}}}{c}\gamma_1
 = m\gamma_1$, so that
 \begin{equation*}
 \frac{\til{\tau}}{\tau} =
 \frac{m+2 \sqrt{dd^{\prime}}}{m-2\sqrt{dd^{\prime}}}
 \in {\Q}[\sqrt{dd^{\prime}}]
 \end{equation*}
  the sign of the square root being chosen in such a way that $p \mid
 {2\sqrt{dd'} -m}$ (in the unramified case).

 \medskip\noindent $\bullet$ The previous relations, together with $\det\Omega =
 1$, show that the transcendence degree over $\ovl{\Q}$ of the entries
 of $\Omega$ is at most one in the CM case. 

 \medskip\noindent $\bullet$ We now describe the
 $p$-adic number $\omega_{11}$ modulo $\ovl{\Q}^{\times}$, keeping the
 notation of \ref{para-lerch-chowla} and \thmref{thm-ogus-formula}
 ($\epsilon$ is the Dirichlet 
 character, $w$ the number of roots of unity, $h$ the class number).

 \begin{thm}\label{thm-cm-period}
  In case $p$ does not ramify in the field of
  complex multiplications ${\Q}[\sqrt{-d}] = \End(A)\otimes{\Q}$, one has
  \begin{equation*}
   \omega_{11}\; \sim \prod_{u\in ({\Z}/d)^{\times}}
    \Bigl(\Gamma_p\Lan{p\frac{u}{d}}\Bigr)^{-\epsilon(u)w/8h} \;\text{ in }\;
    \ovl{\Q}_p^{\times}/\ovl{\Q}^{\times}.
  \end{equation*}
   \end{thm}

\medskip \noindent {\it Examples.} If $A$ is the Legendre curve $X_{1/2}\;$
  ($y^2=x(x-1)(x-1/2)$)
  and $p\equiv 3\pmod{4}$, we find $\omega_{11} \sim \Gamma_p(1/4) .$
  Similarly, if $\End(A)
  \otimes {\Q} = {\Q}(\sqrt{-3})$ and $p\equiv 2\pmod{3}$, then
  $\omega_{11} \sim
  \Gamma_p(1/3)^{3/2}.$

 \begin{proof}
  We may assume that the number field $k$ ($\subset \ovl{\Q}_p$) is
  Galois over $\Q$. We denote by $k'=k\cap {\Q}_p$ the subfield of $k$
  fixed by the local Galois group $\Gal({\ovl{\Q}}_p/{\Q}_p)$. Let us
  consider the Weil restriction $B:= \mathcal{R}_{k/ k'}(A)$
  (\cf \cite[7.6]{bosch90:_nm}). This is an abelian variety
  over $k' \subset {\Q}_p$, with good, supersingular reduction at $p$
  (like $A$), and we have $B_k\simeq \prod_{\tau \in \Gal(k/k')}A^{\tau}$.

  \noindent On the other hand, let $\psi$ be any element of degree one
  in the Weil group (\ie any lifting of the Frobenius element in
  $\Gal({\ovl{\Q}}_p/{\Q}_p)$), and let 
  ${\Psi}_{A^{\tau}}$ (\resp  ${\Psi}_B$) be the corresponding
  semilinear endomorphism of $\mathrm{H}^1_{\mathrm{dR}}(A^{\tau}/k)\otimes
  {\ovl{\Q}}_p$ (\resp $\mathrm{H}^1_{\mathrm{dR}}(B/k')\otimes {\ovl{\Q}}_p$).
  If ${\Phi}_B$ denotes the canonical linear action induced by
  ${Fr_{{\F}_p}}$ on
  $\mathrm{H}^1_{\mathrm{crys}}(B_0/{\Q}_p)\otimes{\ovl{\Q}_p}\simeq
  \mathrm{H}^1_{\mathrm{dR}}(B)_{\ovl{\Q}_p}$, one has the formula ${\Psi}_B =
  {\Phi}_B \circ (\id \otimes \psi)$ (\cf \cite[4]{fadrc1}).  Moreover, ${\Psi}_B=\bigoplus {\Psi}_{A^{\tau}}$ with
  respect to the decomposition
  $\mathrm{H}^1_{\mathrm{dR}}(B,{\ovl{\Q}_p})\simeq
  \bigoplus_{\tau} \mathrm{H}^1_{\mathrm{dR}}(A^{\tau},{\ovl{\Q}_p})$. 
 
  \noindent If $(\omega_1,\omega_2)$ is any symplectic basis of
  eigenvectors for ${\Q}[\sqrt{-d}]$ in
  $\mathrm{H}^1_{\mathrm{dR}}(A)_{\ovl{\Q}}$,
  with $\omega_1 \in \Omega^1(A)_{\ovl{\Q}}$, we
  can write  ${\Psi}_A(\omega_1)= \kappa.\omega_2$,
  where $\kappa$ modulo $\ovl{\Q}^{\times}$ is given by Ogus' formula in
  \thmref{thm-ogus-formula}
  \begin{equation*}
   \kappa\; \sim \prod_{u\in ({\Z}/d)^{\times}}
    \Bigl(\Gamma_p\Lan{p\frac{u}{d}}\Bigr)^{-\epsilon(u)w/4h} \;\text{ in }\;
    \ovl{\Q}_p^{\times}/\ovl{\Q}^{\times}.
  \end{equation*}
   Because the conjugates $A^{\tau}$ are isogenous to each other over $k$ (and
  since isogenies preserve
  the regular differential forms), we derive that the constant $\kappa$, resp.
  the ``period''
  $\omega_{11}$, is the same modulo $\ovl{\Q}^{\times}$ for each of them.
  Therefore, for any
  $\omega \in \Omega^1(B) \subset
  \mathrm{H}^1_{\mathrm{crys}}(B_0/{\Q}_p)$, we can write
  ${\Phi}_B(\omega)={\Psi}_B(\omega)= \kappa.\eta$ for some $\eta \in
  \mathrm{H}^1_{\mathrm{dR}}(B)_{\ovl{\Q}}$ whose pairing with any
  element of $\mathrm{H}_{1\mathrm{B}}(B,{\ovl{\Q}})$ belongs to
  ${\omega_{21}}.{\ovl{\Q}}$. Since $\omega_{11}^{-1}\omega\in
  \mathrm{H}^1_{\mathrm{B}}(B,{\ovl{\Q}}),\,\omega_{21}^{-1}
  \eta\in \mathrm{H}^1_{\mathrm{B}}(B,{\ovl{\Q}})$, and since
  $\Phi_B$ respects
  $\mathrm{H}_{1\mathrm{B}}(B,{\ovl{\Q}})$, we get
  $\kappa \sim \frac{\omega_{11}}{\omega_{21}} \sim \omega_{11}^2$,
  whence the result. 
 \end{proof}
 
\begin{coj}\label{coj-no-complex-multi}
Assume that $A$ is a lifting of $A_0$ defined over $\ovl{\Q}$, and that
$p$ does not ramify $\End(A)\otimes \Q$ (which is obviously the case if
$A$ does not have complex multiplication). Then
$\mathrm{H}^1_{\mathrm{dR}}(B)_{\ovl{\Q}}\neq
\mathrm{H}^1_{1\mathrm{B}}(B,{\ovl{\Q}})$.
\end{coj} 

\noindent In the absence of complex multiplication, this follows from the previous conjecture.
In the presence of complex multiplication, it amounts to the transcendence of $\omega_{11}$. 

 \noindent Taking into account \ref{thm-cm-period}, the conjecture would
 imply, for instance, the transcendence of the
 adic numbers $\Gamma_3(1/4)$, $\Gamma_7(1/4)$, $\Gamma_2(1/3)$,
 $\Gamma_5(1/3)$ (note, in contrast,
 that $\Gamma_5(1/4)$ and $\Gamma_7(1/3)$ are algebraic numbers,
 according to the Gross-Koblitz formula (\ref{para-gross-koblitz})).

 \noindent $\bullet$ Actually, it may happen, in the ramified CM case, that
 $\mathrm{H}^1_{\mathrm{dR}}(A)_{\ovl{\Q}}=
 \mathrm{H}^1_{1\mathrm{B}}(A,{\ovl{\Q}})$. 

 \noindent Let us first notice that the argument given in the proof of theorem
 \ref{thm-cm-period} still works in the ramified case, and allows to
 conclude that $\kappa \sim  \omega_{11}^2$. The computation of $\kappa
 \pmod{{\ovl{\Q}}^\times}$ (for $p$ odd) according to
 Ogus-Coleman-Urfels was explained in \ref{para-ramif}. In the example
 of an elliptic curve with complex multiplication by $\Q[\sqrt{-3}]$, we
 have seen that $\kappa\in \ovl{\Q}$, hence $\omega_{11}\in \ovl{\Q}$
 and $\mathrm{H}^1_{\mathrm{dR}}(A)_{\ovl{\Q}}= 
 \mathrm{H}^1_{1\mathrm{B}}(A,{\ovl{\Q}})$. 

This is also the case for an elliptic curve with complex multiplication
 by $\Q[\sqrt{-3n}]$ 
if the Legendre symbol $(\frac{n}{3})$ is $1$, but probably not in general if
$(\frac{n}{3})=-1$ (although this happens for $n=8$, after the last example in
\ref{para-ramif}).
\end{para}

\begin{para}\label{para-supersingular-disk}
 {\it $L_{\frac{1}{2},\frac{1}{2},1}$ in a supersingular disk.} For
 concreteness, we consider the special case of the supersingular
 Legendre elliptic curve with parameter $z = 1/2$ in characteristic
 $p\equiv 3\pmod{4}$. A basis of solutions of the hypergeometric
 differential operator 
 $L_{\frac{1}{2},\frac{1}{2},1}= HGDO(\frac{1}{2},\frac{1}{2},1)$ in the
 supersingular disk $\matheur{D}(\FRAC{1}{2}, 1^-)$ is given by
 \begin{equation*}
 F\Bigl(\frac{1}{4},\frac{1}{4},\frac{1}{2};(1-2z)^2\Bigr),\;\;
 (1-2z)F\Bigl(\frac{3}{4},\frac{3}{4},\frac{3}{2};(1-2z)^2\Bigr).
 \end{equation*}

\begin{itemize}
 \item Let us consider the symplectic basis
       $\omega_1=[\frac{dx}{2y}],\omega_2=[\frac{(2x-1)dx}{4y}]$ of
       $\mathcal{M}$ (de Rham cohomology). At $z=1/2$, this is a basis of
       eigenvectors for the action of $E':=\End(X_{1/2}) \otimes {\Q}
       ={\Q}(\sqrt{-1})$. The Gauss-Manin connection satisfies
       \begin{equation*}
	\omega_2=2z(z-1)\nabla\Bigl(\frac{d}{dz}\Bigr)\omega_1
	 + \frac{4z-5}{6}\omega_1.
       \end{equation*}
       The fundamental solution matrix $Y$ of
       the Gauss-Manin connection (expressed in the basis
       $(\omega_1,\omega_2)$), normalized by $Y(1/2)=\id$, is then given by
       \begin{align*}
	& y_{11} =
	F\Bigl(\frac{1}{4},\frac{1}{4},\frac{1}{2};(1-2z)^2\Bigr) +
	\frac{1}{2}(1-2z)F\Bigl(\frac{3}{4},\frac{3}{4},\frac{3}{2};
	(1-2z)^2\Bigr),\\
	& y_{12} = (1-2z)F\Bigl(\frac{3}{4},\frac{3}{4},\frac{3}{2};
	(1-2z)^2\Bigr),\\
	& y_{2i} = 2z(z-1) \frac{dy_{1i}}{dz}+\frac{4z-5}{6} y_{1i},
	\quad i=1,2.
       \end{align*}

 \item We choose a symplectic basis $(\gamma_1,\gamma_2)$ of
       eigenvectors for $E'$ in $\mathrm{H}^1_{1\mathrm{B}}(X_{1/2},F)$
       (note that ${\Q}(\sqrt{-1})$ is always contained in $F$). We have
       $\mathrm{H}^1_{1\mathrm{B}}(X_{1/2},F) =
       \mathrm{H}^1_{1\mathrm{B}}(X_z,F)$ for any point $z$ in the
       supersingular disk $\matheur{D}(\FRAC{1}{2}, 1^-)$.
       Because the period matrix $\Omega(z)$ is another solution matrix
       of the Gauss-Manin connection, we have the relation
       $Y(z)=\Omega(z).\Omega(1/2)^{-1}$.
       
       Let us consider a CM-point $\zeta$ in the supersingular
       disk $\matheur{D}(\FRAC{1}{2}, 1^-)$: $\End(X_{\zeta}) =
       E ={\Q}(\sqrt{-d})$ (and $E\neq E'$ in $\bD$).

       \begin{figure}[h]
	\begin{picture}(150,120)(0,10)
	 \put(50,0){\includegraphics{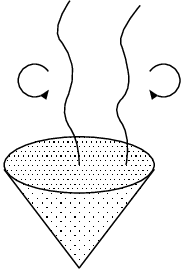}}
	 \put(55,118){$X_{\FRAC{1}{2}}$}
	 \put(120,118){$X_{\zeta}$}
	 \put(0,40){$\matheur{D}(\FRAC{1}{2}, 1^-)$}
	 \put(22,95){$\scriptstyle\Q(\sqrt{-1})$}
	 \put(142,95){$\scriptstyle\Q(\sqrt{-d})$}
	 \put(75,50){$\scriptstyle\FRAC{1}{2}$}
	 \put(114,51){$\scriptstyle\zeta$}
	\end{picture}
	\caption{}
       \end{figure} 
              
       \noindent Let $\sigma$ and $m$ be as before (in the present case, the
       integer $m$ is even:
       $m=2n$). Our period relations can be easily expressed in terms of
       $Y(\zeta)$. Apart from the
       obvious relation $\det Y(\zeta)=1$, we find a relation
       \footnote{F. Beukers has also found these
       relations independently, by a different
       method.}
       \begin{equation*}
	\left.\frac{{y_{11}(y_{22}+\sigma y_{12})}}{{y_{12}(y_{21}+\sigma
	 y_{11})}}\right|_{z=\zeta} = \frac{\til{\tau}}{\tau} =
	\frac{n+\sqrt{d}}{n-\sqrt{d}}.
       \end{equation*}

       \noindent Changing the viewpoint, one can consider $\zeta$ as
       fixed, and vary the prime $p\neq 2$ (or more precisely the place
       of ${\Q}(\zeta,\sqrt{-1},\sqrt{-d}$)).  This relation
       between values of $p$-adic hypergeometric functions at
       $\zeta$ holds whenever $\vert \zeta - 1/2 \vert_p < 1$; it
       depends on $p$ via the integer $n$ ($\vert n\vert < \sqrt{d}, \;
       p \mid d - n^2$). A relation of the same kind holds at the places
       at infinity (derived along similar lines, using the usual Betti
       lattices).  As a specific example, we can take $p=3, \zeta=$ a
       primitive $6$th root of unity.  Then $X_{\zeta}$ has complex
       multiplication by ${\Z}[\zeta]$, $\;\sigma =
       \frac{-\zeta^2}{2(1+\zeta)}\;, \;n=0$ and
       $\frac{\til{\tau}}{\tau}=-1$.
 \item On the other hand, one can combine \thmref{thm-cm-period} and
       the formula $Y(\zeta)=\Omega(\zeta).\Omega(1/2)^{-1}$ in order to
       express the value $\bmod{\,\ovl{\Q}^{\times}}$ of the $p$-adic
       hypergeometric functions $y_{ij}$ in terms of $\Gamma_p$ ({\it
       supersingular avatar of Young's formulas}).
       
       \noindent For instance, if $p=7, \zeta = 2(\sqrt{2} -1)$ (complex
       multiplication by ${\Z}[\sqrt{-2}]$), one has $\omega_{11}(\zeta) \sim
       (\Gamma_7(1/8)\Gamma_7(3/8))^{1/2},\,\omega_{11}(1/2) \sim
       \Gamma_7(1/4) $, from which
       one derives the 7-adic evaluation
       \begin{align*}
	F\left(\frac{3}{4},\frac{3}{4},\frac{3}{2};5-4\sqrt{2}\right) & \sim
	 \bigl(\Gamma_7(1/8)\Gamma_7(3/8)\bigr)^{1/2}\Gamma_7(1/4)^{-1}.
       \end{align*}
\end{itemize}
\end{para}

\subsection{Conclusion and vista.}
We have seen many instances of the following scenario: some familiar
notion or phenomenon from the complex realm shows itself in two
different aspects; each of these aspects has a natural $p$-adic
counterpart; these counterparts are not complementary aspects of the
same $p$-adic entity, but belong to totally different theories.  

In II, we shall see that the theory of $p$-adic period mappings is a privileged
field where this semantic splitting process integrates into an
harmonious picture, where the complementarity between the $p$-adic
counterparts is restored at a deeper level.

The $p$-adic period mappings which we shall deal with relate some
deformation spaces of $p$-divisible groups to certain grassmannians. It
will turn out that they can be described by quotients of analytic
solutions of the Gauss-Manin connection, just as in the complex case. In
particular, at the {\it modular level}, the Fontaine-Messing periods
(which links up directly the \'etale and crystalline representations of
$p$-divisible groups) go offstage or remain at the background.

On the other hand, we shall see that the differential equations
studied in section 1 (\resp 2) arise in the context of deformations of
supersingular (resp. ordinary) $p$-divisible groups.

\chapter{Introduction to the theory of $p$-adic period mappings.}
\label{chap:period-mappings}
\addtocontents{toc}{\protect\par\vskip2mm\hskip20mm(Geometry)\par\vskip5mm}
\minitoc
\newpage
\section{A survey of moduli and period mappings over $\C$.}\label{sec:1}

\begin{abst}
Review of moduli, Shimura varieties, 
Hodge structures, period mappings, period domains, Gauss-Manin connections.
\end{abst}
   
\subsection{Moduli spaces.}\label{sub:1.1}

\begin{para}\label{para:1.1.1} The term ``moduli'' was
 coined by Riemann to describe the continuous essential
 parameters of smooth complex compact curves of a given genus $g$. The use
 of this term has been extended to many classification problems in 
 analytic or algebraic
 geometry: typically, a classification problem consists of a discrete 
 part (the dimension, and other
 numerical invariants...) and a continuous part (the moduli). 

 One can distinguish between local and global moduli spaces: 

\begin{enumerate}
 \item local moduli are the parameters of the most general small (or
       infinitesimal) deformation of a given analytic manifold or
       algebraic variety. Let for instance $X$ be a compact analytic
       manifold.  Then a famous theorem of Kuranishi asserts that there
       is a local deformation $\pi: \udl{X}\rightarrow S$ with base
       point $s$ ($X={\udl{X}}_s$) such that any other deformation
       comes, locally around $s$, from $\pi$ by base change; moreover
       the tangent space of $S$ at $s$ is isomorphic to
       $\mathrm{H}^1(X,(\Omega^1_X)^\vee)$. For curves, the dimension of
       this space is $3g-3$ if $g\geq 2$ (the number of Riemann's
       complex moduli).
 \item Global moduli spaces, when they exist, provide the
       solution of the continuous part of moduli problems. For instance,
       Douady has constructed the global moduli space for closed analytic
       subspaces of a given compact complex manifold. Although global moduli
       spaces were already studied in the nineteenth century (in the case of
       curves), the concept was first put on firm foundations by Grothendieck
       in terms of representable functors.
\end{enumerate}

The existence of a global moduli space is often obstructed by the
presence of automorphisms of the objects under classification. One
possibility is to rule out these automorphisms by imposing some extra
structure (rigidification). When all automorphism groups are finite,
another possibility is to look for moduli orbifolds rather than moduli
spaces
\end{para}

\begin{para}\label{para:1.1.2}
 The most classical moduli problems concern complex elliptic
 curves. These are classified by their $j$-invariant. However the affine
 $j$-line ${\A}^1$ is not a fine moduli space: indeed, in any family of
 elliptic curves parametrized by the whole affine line, all the fibers
 are isomorphic.

 One can rigidify the problem by introducing a ``level $N$ structure''
 \index{levelNstructure@level $N$ structure} on 
 elliptic curves $E$, \ie fixing an identification $\Ker[N]_A\cong
 ({\Z}/N{\Z})^2$ such that the Weil pairing is given by
 $(\xi,\eta)\mapsto \exp(2\pi i(\xi\wedge\eta))$ (here $[N]_A$ denotes the
 multiplication by $N$ on the elliptic curve $A$).

For $N\geq 3$, this rules out automorphisms of $A$, and the
corresponding moduli functor is representable. For $N=3$, the moduli
space is ${\A}^1\setminus \{\mu^3=1\}$; the universal elliptic curve
with level $3$ structure is the Dixon elliptic curve
\begin{equation*}
 x^3+y^3+3\mu xy = 1 \text{ (in affine coordinates $x,y$)}
\end{equation*}
with $\Ker[3]_A=$ the set of flexes.

For $N=2$, the automorphism $[-1]_A$ remains; there is an elliptic curve
with level $2$ structure over ${\P}^1\setminus \{0,1,\infty\}$, namely
the Legendre elliptic curve
\begin{equation*}
 y^2=x(x-1)(x-\lambda) \text{ (in affine coordinates $x,y$)}
\end{equation*} 
with
$\Ker[2]_A=\{(0,0),(1,0),(\lambda, 0),(\infty,\infty)\}$, but it is not
universal: another one is given by $\lambda y^2=x(x-1)(x-\lambda)$.
\end{para}

\begin{para}\label{para:1.1.3}
 These moduli problems admit two natural generalizations: moduli
 problems for curves of higher genus, and moduli problems for polarized
 abelian varieties (polarizations ensure that the automorphism groups
 are finite).

 There is a moduli orbifold for smooth compact complex curves of genus
 $g>0$, denoted by $\mathcal{M}_g$, which can be described as
 follows. The set of complex structures on compact connected oriented
 surface $S_g$ of genus $g$ is in a natural way a contractible complex
 manifold of dimension $3g-3$ (\resp the Poincar\'e upper half plane
 $\mathfrak{h}$ if $g=1$), the Teichm\"uller space $\mathcal{T}_g$. The
 mapping class group $\Gamma_g$, \ie the component of connected
 components of the group of orientation-preserving self-homeomorphisms of
 $S_g$, acts on right on $\mathcal{T}_g$, and $\mathcal{M}_g$ is the
 quotient $\mathcal{T}_g/\Gamma_g$. For $g=1$, one recovers
 $\mathcal{M}_1=\mathfrak{h}/SL(2,{\Z})$.

There is a moduli orbifold for principally polarized abelian varieties of
dimension
$g>0$, denoted by $\mathcal{A}_g$, which can be described as follows. For any
abelian variety $A$, one has an exact sequence 
\begin{equation*}
 0\rightarrow L\rightarrow \Lie A\rightarrow A\rightarrow 0
\end{equation*}
and a principal polarization corresponds to a positive
definite hermitian form $H$ on $\Lie A$ whose imaginary part induces a
perfect alternate form on $L$. Choosing a $\C$-basis
$\omega_1,\ldots,\omega_g$ of $\Omega^1(A)\cong (\Lie A)^\vee$ and a
symplectic $\Z$-basis $\gamma_1^+,\gamma_1^-,\ldots,\gamma_g^+,
\gamma_g^-$ of $L$ gives rise to a period matrices 
$\Omega^+=(\int_{\gamma^+_j}\omega_i)$, $\Omega^-=(\int_{\gamma^-_j}\omega_i)$,
and $\tau=(\Omega^-)^{-1}\Omega^+$ does not depend on the
$\omega_i$'s.  By the properties of $H$, $\tau$ belongs to the Siegel
upper half space $\mathfrak{h}_g$ of symmetric matrices with positive
definite imaginary part (Riemann relations); it is well-defined up to
right multiplication by a matrix in $Sp(2g,{\Z})$ (which is the same as
the standard left action on $\mathfrak{h}_g$ of the transposed
matrix). This construction identifies $\mathcal{A}_g$ with
$\mathfrak{h}_g/Sp(2g,{\Z})$.

The construction which associates to a curve its jacobian variety with
the principal polarization given by the theta divisor gives rise to an
immersion $\mathcal{M}_g\hookrightarrow\mathcal{A}_g$, the Torelli map.
\end{para}

\begin{para}\label{para:1.1.4}
 Following G.~Shimura, one also studies refined moduli problems for
 ``decorated'' (principally polarized) abelian varieties, by prescribing
 in addition that the endomorphism algebra contains a given simple
 $\Q$-algebra, and imposing a level $N$ structure. This gives rise to
 moduli orbifolds (\resp moduli spaces if $N\geq 3$) of ``PEL type''.

 More precisely, let $B$ be a simple finite-dimensional $\Q$-algebra,
 with a positive involution $\ast$. Let $V$ be a $B$-module of finite
 type, endowed with an alternate $\Q$-bilinear form such that
 $\an{bv,w}=\an{v,b^\ast w}$.  One denotes by $G$ the algebraic $\Q$-group of
 $B$-linear symplectic similitudes of $V$.

 Let $\mathcal{B}$ be a maximal order in $B$ stable under $\ast$, and
let $L$ be a lattice in $V$, stable under $\mathcal{B}$ and autodual for
$\an{\;,\;}$ (this will represent the structure of $\mathrm{H}^1(A,{\Z})$ for
the abelian varieties $A$ under consideration). We choose in addition a
lagrangian subspace $F^1_0\subset V_{\C}$ stable under $B$ (which
will be $\Omega^1(A_0)\cong (\Lie A_0)^\vee$ for one particular abelian
variety $A_0$ under consideration).

There is a moduli orbifold (\resp moduli space if $N\geq 3$)
for principally $\ast$-polarized\footnote{a principal
$\ast$-polarization is a symmetric $\mathcal{B}$-linear isomorphism
between the abelian variety and its dual, endowed with the transposed
action of $\mathcal{B}$ twisted by $\ast$} abelian varieties $A$,
together with a given action $\iota: \mathcal{B}\rightarrow \End A$
such that $\det(\iota(b)| \Omega^1(A))=\det(b| F^1_0)$ for all
$b\in \mathcal{B}$ (Shimura type condition), and a given
$\mathcal{B}/N\mathcal{B}$-linear symplectic similitude
$A[N]=\Ker[N]_A\rightarrow \Hom(L,{\Z}/N {\Z})$.

It amounts to the same to consider abelian varieties $A$ up to isogeny,
with $B$-action (with Shimura type condition), together with a
$\Q$-homogeneous principal $\ast$-polarization and a class of $B$-linear
symplectic similitudes $\mathrm{H}^1(A,{\hat {\Z}}\otimes
\Q)\rightarrow {\hat {\Z}}\otimes V$ modulo the group $\{g \in
G({\hat {\Z}}\otimes {\Q}) \mid (g-1)({\hat {\Z}}\otimes L)\subset {\hat
{\Z}}\otimes L\}$.

These moduli orbifolds are called Shimura orbifolds of ``PEL
type''. \index{PELtype@PEL type} They are algebraic, and defined over
 the number field 
$E={\Q}[\tr(b| F^1_0)]$ (the so-called reflex field). Their
connected components are defined over finite abelian extensions of $E$,
\cf \cite{tds}.
\end{para}

\subsection{Period mappings.}\label{sub:1.2} 

\begin{para}\label{para:1.2.1}
 In the case of curves, the construction of period mappings belongs to
 the nineteenth
 century. The idea is to study the deformations of a curve by looking at
 the variation of its periods.  Let $f:{\udl{X}}\rightarrow S$ be an
 algebraic family of projective smooth (connected) curves of genus $g$
 over $\C$. Let $\tilde S$ be the universal covering of $S$. Then the
 homology groups $\mathrm{H}_1(-,{\Z})$ of the fibers form a constant
 local system on $\tilde S$. Let us choose a symplectic basis
 $\gamma_1^+,\gamma_1^-,\ldots ,\gamma_g^+, \gamma_g^-$. Choosing,
 locally, an auxiliary $\C$-basis $\omega_1,\ldots,\omega_g$ of
 $f_\ast\Omega^1_{\udl{X} /S}$ gives rise as above to a period matrices
 $\Omega^+=(\int_{\gamma^+_j}\omega_i)$, 
 $\Omega^-=(\int_{\gamma^-_j}\omega_i)$,  
 and $\mathcal{P}=(\Omega^-)^{-1}\Omega^+$ does not depend on the
 $\omega_i$'s and is well-defined on $\tilde S$.

 We get a commutative square 

 \begin{equation*}
  \begin{CD}
   \til{S} @>\mathcal{P}>> \mathfrak{h}_g\\
   @VVV @VVV\\
   S @>\mathcal{P}>> \mathfrak{h}_g/Sp(2g,\Z)
  \end{CD}
 \end{equation*}
 Applying this construction to the universal case
 $S=\mathcal{M}_g$ ($\mathcal{T}_g=\til{\mathcal{M}}_g$ in the orbifold
 sense), the bottom line of this commutative square is nothing but an
 analytic interpretation of the Torelli map.
\end{para}

\begin{para}\label{para:1.2.2}
 An important feature of the period mapping $\mathcal{P}$ is that it is
 a {\it quotient of solutions of a fuchsian differential system}, with
 rational exponents at infinity. This comes out as follows. The vector
 bundle $f_\ast\Omega^1_{\udl{X} /S}$ is a subbundle (of rank $g$) of a
 vector bundle $\mathcal{H}$ (of rank $2g$) whose fibers are the De Rham
 cohomology spaces $\mathrm{H}^1_{\mathrm{DR}}(\udl{X}_s)$. This bundle
 $\mathcal{H}$ carries an integrable connection, the Gauss-Manin
 connection
 \begin{equation*}
  \nabla_{\mathrm{GM}}:\;\mathcal{H}\rightarrow \Omega^1_S\otimes \mathcal{H}.
 \end{equation*}
 The local system of germs of horizontal analytic sections is
 $R^1f_\ast{\C}$ (with fibers
 $\mathrm{H}^1(\udl{X}_s,{\C})$). Completing, locally, the basis
 $\omega_1,\ldots,\omega_g$ of sections of $f_\ast\Omega^1_{\udl{X}/S}$
 into a basis $\omega_1,\ldots,\omega_g,\eta_1,\ldots,\eta_{g}$ of
 sections $\mathcal{H}$ allows to write the connection in form of a
 differential system, and a full {\it solution matrix} is given by the
 ``full $2g\times 2g$-period matrix''
 \index{full 2g 2g-period matrix@full $2g\times 2g$-period matrix} 
 $\begin{pmatrix} \Omega^+ &
  \Omega^- \cr N^+ & N^-  \end{pmatrix}$, while
 $\mathcal{P}=(\Omega^-)^{-1}\Omega^+$.

 In the case of the Legendre elliptic curve over ${\P}^1\setminus
 \{0,1,\infty\}$ (\cf I.\ref{sec-tale}), one finds that $\mathcal{P}$
 is given by the 
 quotient $\tau$ of two solutions of the hypergeometric differential
 equation with parameters $\bigl(\frac{1}{2},\frac{1}{2},1\bigr)$.
\end{para}

\begin{para}\label{para:1.2.3}
 The theory of period mappings has been generalized to any projective
 smooth morphism $f:\udl{X}\rightarrow S$ by P. Griffiths
 \cite{poioam3}.  We assume that $S$ and the fibers of $f$ are connected for
 simplicity.  The idea is the following: by abstracting the structure of
 the periods of the fibers of $f$, one arrives at the notion of
 polarized Hodge structure; one then construct a map from the universal
 covering $\tilde S$ to a classifying space $\mathcal{D}$ of polarized
 Hodge structures.

 For a fixed $n$, the real cohomology spaces
 $V_s^{\R}:= \mathrm{H}^n(\udl{X}_s,{\R})$) form a local system $V^{\R}$
 on $S$. We denote by $V^{\C}$ its complexification. The vector bundle
 $\mathcal{O}_S\otimes V^{\C}$ is the analytification of an algebraic
 vector bundle $\mathcal{H}$ on $S$, and $V^{\C}$ is the local system of
 germs of analytic solutions of an algebraic fuchsian integrable
 connection $\nabla_{\mathrm{GM}}$ on $\mathcal{H}$, the Gauss-Manin
 connection.

 On the other hand, the spaces $V_s^{\R}$ carry a natural algebraic
 representation of ${\C}^\times$ (Hodge structure), which amounts to a
 decomposition $V_s^{\C}=\bigoplus_{p+q=n,\;p\geq 0}
 \;V^{p,q},\;\;\overline{V^{p,q}}=V^{q,p}$. When $s$ varies, the
 filtrations
 \begin{equation*}
  F_s^p:=\bigoplus_{p'+q=n,\;p'\geq p} V^{p',q}
 \end{equation*}
 are the fibers of a decreasing filtration $F_s^p$ of $\mathcal{H}$ by
 vector subbundles (the Hodge filtration). The Gauss-Manin connection
 fails to preserve the Hodge filtration by just one notch:
 $\nabla_{\mathrm{GM}}(F^p)\subset \Omega^1_S\otimes F^{p-1}$.

 The last piece of data is the polarization: a bilinear form $Q$ on
 $V_s^{\R}$, $Q(v,v')=(-1)^nQ(v',v)$, such that $(-i)^nQ(v,\ovl{v})$ is
 hermitian on $V_s^{\C}$, definite (of sign $(-1)^p$) on each
 $V^{p,q}=F_s^p\cap \overline F_s^q$, and such that $F_s^p$ is the
 orthogonal of $F_s^{n-p-1}$ in $V_s^{\C}$.
 
 Let $G_{\R}$ denote the real group of similitudes of $Q$. The
 space of isotropic flags of type $0\subset \ldots \subset F_s^p\subset
 F_s^{p-1}\subset \ldots \subset V_s^{\C} \;$ is a homogeneous space
 \begin{equation*}
  \mathcal{D}^\vee = P\backslash G_{\C}
 \end{equation*}
 for a suitable parabolic subgroup $P$.

 The flags satisfying the above positivity condition (on each $F_s^p\cap
 \ovl{F}_s^q$) are classified by an analytic open submanifold, the
 {\it period domain},
 \begin{equation*}
  \mathcal{D}\subset \mathcal{D}^\vee ;\;\; \mathcal{D}
   \cong K\backslash G_\R^{\mathrm{ad}}({\R})^{0}
 \end{equation*}
 \index{000GRad@$G_\R^{\mathrm{ad}}$}
 for a suitable compact subgroup $K\;$ ($\mathcal{D}^\vee$ is called the
 compact dual of $\mathcal{D}$).

Over $\tilde S$, the local system $V^{\C}$ becomes constant, and the
construction which attaches to any $s\in \tilde S$ the classifying point
of the corresponding flag $0\subset \ldots \subset F_s^p\subset
F_s^{p-1}\subset \ldots \subset V_s^{\C}\;$ gives rise to a holomorphic
mapping, the {\it period mapping}

 \begin{equation*}
  \tilde{S}\xrightarrow{\mathcal{P}} \mathcal{D}.
 \end{equation*}
 Via Pl\"ucker coordinates, $\mathcal{P}$ is again
 given by quotients of solutions of the Gauss-Manin connection
 $\nabla_{\mathrm{GM}}$.  The projective monodromy group of
 $\nabla_{\mathrm{GM}}$ is a 
 subgroup $\Gamma \subset G_{\R}^{\mathrm{ad}}$ (well-defined up to
 conjugation), and we get a commutative square

\begin{equation*}
 \begin{CD}
  \tilde{S}@>\mathcal{P}>> \mathcal{D} @. \subset \;\mathcal{D}^{\vee}\\
  @VVV @VVV @.\\
  S@>\mathcal{P}>> \mathcal{D}/\Gamma.
 \end{CD}
\end{equation*}
\end{para}

\begin{para}\label{para:1.2.4} In the simple case of a
 family of elliptic curves and $n=1$, $G_{\R}=GL_{2,\R}$, the flag space
 is just $\mathcal{D}^\vee={\P}^1_{\C}$, $K=PSO(2)$ (the isotropy group
 of $i\in {\P}^1_{\C}$ in $PSL_2({\R})=G_\R^{\mathrm{ad}}({\R})^0$), and
 $\mathcal{D}\cong K\backslash G_\R^{\mathrm{ad}}({\R})^0$ is the
 Poincar\'e upper half plane $\mathfrak{h}$.

 More generally, for a polarized abelian
 scheme, $G_{\R}$ is the group of symplectic similitudes, and
 $\mathcal{D}^\vee$ is the grassmannian of lagrangian subspaces of
 ${\C}^{2g}$, and $\mathcal{D}$ is the Siegel upper half space.

 \bigskip
 When one deals with refined moduli problems for abelian schemes with PEL
 decoration as in \ref{para:1.1.4}, one should work with the group $G_{\R}$ of
 $B$-{\it linear} symplectic similitudes of $V_{\R}$ (notation as in
 \ref{para:1.1.4}). In that case, the corresponding period domain
 $\mathcal{D}$ is a symmetric domain.

 More precisely, the datum of our $B$-stable lagrangian subspace
 $F^1_0\subset V_{\C}$ (attached to our particular $A_0$ with
 $\mathrm{H}^1(A_0,{\C})=V_{\C}$) amounts to the datum of an algebraic
 $\R$-homomorphism $h_0: {\C}^\times\rightarrow G_{\R}\subset
 GL(V_{\R})$ such that $h_0$ acts by the characters $z$ and $\ovl{z}$ on
 $V_{\C}$ ($F^1$ is the subspace where $h_0$ acts through $z$), and such
 that the symmetric form $\an{v,h(i) w}$ is positive definite on
 $V_{\R}$. For any abelian variety $A$ with PEL structure of the correct
 type, together with a fixed isomorphism $\mathrm{H}^1(A,{\C})\cong
 V_{\C}$, the corresponding subspace
 $F^1\mathrm{H}^1(A,{\C})=\Omega^1(A)\subset V_{\C}$ is an eigenspace
 for some conjugate of $h_0$. The set of conjugates of $h_0$ under
 $G({\R})$ is a finite union of copies of a symmetric domain
 $\mathcal{D}$.  \par For any component $S$ of the {\it Shimura
 orbifold}, we thus get a commutative square
 
 \begin{equation*}
  \begin{CD}
   \tilde{S}@>\mathcal{P}>> \mathcal{D} @. \subset\; \mathcal{D}^{\vee}\\
   @V\mathcal{Q}VV @VV\mathcal{Q}V @.\\
   S@>\mathcal{P}>> \mathcal{D}/\Gamma
  \end{CD}
 \end{equation*}
 \par\noindent
 where the horizontal maps are {\it isomorphisms}, $\mathcal{Q}$ denotes
 the quotient maps, and $\Gamma$ is a congruence subgroup of level $N$
 in the semi-simple group $G^{\mathrm{ad}}$ (a conjugate of the standard
 one).

 \bigskip
 The complex manifold $\tilde S$ admits the following
 modular description. Let $U$ be the oriented real Lie group
 $\Hom(L,(\R/\Z)^\times)$ (with $\mathcal{B}$-action and polarization),
 and let $\hat U$ be the associated infinitesimal Lie group. The
 automorphism group $J$ of the decorated $\hat U$ coincides with
 $G(\R)^0$. Then $\tilde S$ is a moduli space for {\it marked} decorated
 abelian varieties, where the marking is an isomorphism $\rho: \hat
 A\rightarrow \hat U$ of associated infinitesimal decorated oriented
 real Lie groups (note that, by taking duals of Lie algebras, $\rho$
 amounts to an $B$-linear symplectic isomorphism
 $\mathrm{H}^1(A,\R)\cong V^{\R}$). Note that the standard right action
 of $J$ on $\tilde S$ is obtained by functoriality.

\newpage
\section{Preliminaries on $p$-divisible groups.}\label{sec:2}

\begin{abst}
 Definitions and basic theorems about $p$-divisible groups,
 quasi-isogenies, liftings and deformations.
\end{abst}
\end{para}

\subsection{$p$-divisible groups and quasi-isogenies.}\label{sub:2.1}

\begin{para}\label{para:2.1.1} We begin with four definitions. Let $S$ be a scheme. Let
$\Lambda$ and $\Lambda'$ be commutative group schemes over $S$.  

A
homomorphism $f:\Lambda \rightarrow \Lambda'$ is called an {\it isogeny}
if it is an (f.p.p.f.) epimorphism with finite locally free kernel.

Let $p$ be a prime number. $\Lambda$ is a {\it $p$-divisible group}
(or Barsotti-Tate group) if $\displaystyle
\Lambda=\underset{n}{\liminj}\Ker [p^n]$ and for all $n$, the multiplication
$[p^n]$ by $p^n$ on $\Lambda$ is an isogeny.

\par\noindent If $S$ is connected, the rank of $\Ker [p^n]$ is then of
the form $p^{hn}$, where $h$ is an integer called the height of $\Lambda$.  

From the fact that $[p]$ is an isogeny, it follows that for 
$p$-divisible groups, $\Hom_S(\Lambda,\Lambda')$ is a torsion-free 
${\Z}_p$-module. 

A {\it quasi-isogeny} of $p$-divisible groups $\Lambda, \Lambda'$ is a
global section $\rho$ of the Zariski sheaf $ \udl{\Hom}_S(\Lambda,
\Lambda')\otimes_{\Z}{\Q}$ such that there exists locally an integer $n$
for which $p^n\rho$ is an isogeny. By abuse, one writes $\rho: \Lambda
\rightarrow \Lambda'$ as for homomorphisms. We denote by
$\qisog_S(\Lambda, \Lambda')$\index{000qisog@qisog} the ${\Q}_p$-space of
quasi-isogenies.

The objects of the category of {\it $p$-divisible groups on $S$ up
to isogeny} are the $p$-divisible groups over $S$; morphisms between
$\Lambda $ and $\Lambda'$ are global sections of $ \udl{\Hom}_S(\Lambda,
\Lambda')\otimes_{\Z}{\Q}$. This is a $\Q_p$-linear category.
\end{para}

\begin{para}\label{para:2.1.2} The most important example of a
 $p$-divisible group is $\Lambda =A[p^\infty]$, the (inductive system
 of) $p$-primary torsion of an abelian scheme $A$ over $S$
 (\cite{pg}); in this case, the height 
$h$ is twice the relative dimension of $A$. 
\end{para}

\begin{para}\label{para:2.1.3}
 When $p$ is locally nilpotent on $S$, any $p$-divisible group $\Lambda$
 over $S$ is formally smooth; the completion $\hat\Lambda$ of a
 $p$-divisible group $\Lambda$ along the zero section is a formal Lie
 group \cite[VI,3.1]{grothendieck74:_group_barsot_tate_dieud}.  However
 ${\hat\Lambda}$ is not necessarily 
 itself a $p$-divisible group; it is if and only if the separable rank of
 the fibers of $\Ker[p]$ is a locally constant function on $S$ (in which
 case $\Lambda$ is an extension
 \begin{equation*}
  1\rightarrow \hat{\Lambda}\rightarrow \Lambda\rightarrow
   \Lambda^{\mathrm{et}}\rightarrow 1 
 \end{equation*} 
 of an ind-etale $p$-divisible group $\Lambda^{\mathrm{et}}$ by the
 infinitesimal $p$-divisible group $\hat\Lambda$, \cf
 \cite[III,7.4]{grothendieck74:_group_barsot_tate_dieud}).
 An example of an ind-etale (\resp infinitesimal)
 $p$-divisible group of height $h=1\;$ is
 $\displaystyle {\Q}_p/{\Z}_p=\underset{n}{\liminj} {\Z}/p^n{\Z}$ (\resp
 $\hat {\G}_m$).
\end{para}

\subsection{Two theorems on $p$-divisible groups.}\label{sub:2.2} 

\begin{para}\label{para:2.2.1}
 A considerable amount of work has been
 devoted to $p$-divisible groups and their 
 applications to the $p$-adic study of abelian varieties and their local
 moduli. Among the pioneers, let us mention: Barsotti, Tate, Serre,
 Grothendieck, Lubin, Messing...

 In the sequel, we shall outline the main concepts and results of the
 Rapoport-Zink theory of period mappings for $p$-divisible groups,
 emphasizing the analogy with the complex case. We refer to their book
 \cite{psfpg} for details and proofs. We shall also emphasize the
 relation to differential equations, which does not appear in
 \cite{psfpg}.

 The whole theory relies upon three basic theorems on $p$-divisible
 groups and their Dieudonn\'e modules, which we shall recall below:

 \begin{enumerate}
  \item the rigidity theorem for $p$-divisible groups up to isogeny,
  \item the Serre-Tate theorem, and
  \item the Grothendieck-Messing theorem,
 \end{enumerate}

 \noindent
 which deal with infinitesimal deformations of $p$-divisible groups.
\end{para}

\begin{para}\label{para:2.2.2} Let $\mathfrak{v}$ be a complete discrete
valued ring $\mathfrak{v}$ of mixed characteristic $(0,p)$. Let 
$\Nil_\mathfrak{v}$ denote the category of locally noetherian 
$\mathfrak{v}$-schemes $S$ on which $p$ is locally nilpotent.  For any
$S$ in $\Nil_\mathfrak{v}$, the ideal of definition 
$\mathcal{J}=\mathcal{J}_{S_{\mathrm{red}}}$ of the closed subscheme $S_{\mathrm{red}}$ 
(of characteristic $p$) is locally nilpotent.
\end{para}

\begin{thm}[Rigidity theorem for $p$-divisible groups up to isogeny]
 \label{thm:2.2.3}
 
 Every homomorphism
 $\;\ovl{\rho}:\;\Lambda\times_S S_{\mathrm{red}} \rightarrow \Lambda'\times_S
 S_{\mathrm{red}}$ of $p$-divisible groups up to isogeny admits a unique lifting
 $\rho: \Lambda \rightarrow \Lambda'$. Moreover, $\rho$ is a
 quasi-isogeny if $\ovl{\rho}$ is.
\end{thm}

\begin{proof}
 (following an argument of V. Drinfeld). We may assume that $S$ is
 affine, and that $\mathcal{J}$ is nilpotent. Then the connected part
 ${\hat\Lambda'}$ is a formal Lie group, hence there is a power $p^n$ of
 $p$ such that for any affine $S$-scheme $S'$, ${\hat\Lambda}'(S')$ is
 killed by $[p^n]$ (in fact if $\mathcal{J}^{r+1}=0$ on $S$, one can take
 $n=r^2$).

 We observe that $\Ker(\Lambda'(S')\rightarrow
 \Lambda'(S'_{\mathrm{red}}))= {\hat\Lambda}'(S')$. We have
 \begin{align*}
 &\Ker(\Hom(\Lambda(S') ,
 \Lambda'(S'))\rightarrow \Hom(\Lambda(
 S'_{\mathrm{red}}) , \Lambda'(S'_{\mathrm{red}}) )\\
 =&
 \Hom (\Lambda(S') ,\Ker(\Lambda'(S')
 \rightarrow \Lambda'(S'_{\mathrm{red}}))) 
 \end{align*}
 which is zero because $\Ker(\Lambda'(S')\rightarrow
 \Lambda'(S'_{\mathrm{red}})$ is killed by $[p^n]$ while $\Lambda$ is
 $p$-divisible.

 \noindent This implies the injectivity of
 $\qisog_S(\Lambda , \Lambda') \rightarrow 
 \qisog_{S_{\mathrm{red}}}(\Lambda\times_S S_{\mathrm{red}} ,
 \Lambda'\times_S S_{\mathrm{red}}) $ (taking into account the fact that 
 $\Hom(\Lambda( S'_{\mathrm{red}}), \Lambda'(S'_{\mathrm{red}}) )$ is
 torsion-free).

 \noindent For the surjectivity, let us first show that for any
 ${\ovl{f}}\in \Hom(\Lambda(S'_{\mathrm{red}}),
 \Lambda'(S'_{\mathrm{red}}))$, there is 
 a lifting $g\in \Hom(\Lambda(S'), \Lambda'(S'))$ of $p^n {\ovl{f}}$. For
 any $x\in \Lambda( S')$, let $\ovl{x}$ denote its image in $\Lambda(
 S'_{\mathrm{red}})$. Since $\Lambda'$ is formally smooth,
 ${\ovl{f}}({\ovl{x}})$ 
 admits a lifting $y\in \Lambda'(S')$.  Since
 $\Ker(\Lambda'(S')\rightarrow \Lambda'(S'_{\mathrm{red}}))) $ is killed
 by $p^n$, 
 $g(x):=[p^n] y$ is well-defined. This construction provides the desired
 lifting $g$.

 It remains to show that $g$ is an isogeny if
 $\ovl{f}$ is. Let $\ovl{f}'$ be a quasi-inverse of $\ovl{f}$
 ($\ovl{f}'\ovl{f}=[p^m]$), and $g'$ be a lifting of $p^n\ovl{f}'$. Then
 $g'g=[p^{2n+m}]$ by unicity of liftings. Thus $g$ is an epimorphism, and
 the subscheme $\Ker g$ of $\Ker [p^{2n+m}]$ is finite over $S$. On the
 other hand, $\Lambda$ is flat over $S$ and the fibers of $g'$ are flat
 (being isogenies). The criterium of flatness fiber by fiber
 \cite[11.3.10]{ega44} implies that $g'$ is flat, and we conclude that
 it is an isogeny.
\end{proof}

Let us observe that there is variant of both the statement and the
proof of \ref{thm:2.2.3} for abelian schemes, instead of $p$-divisible groups. 

 The next theorem asserts that deforming an abelian scheme of
 characteristic 
 $p$ is equivalent to deforming its $p$-divisible group. For any $S $ in
 $\Nil_\mathfrak{v}$, let $\Def_{p\text{-div}}(S)$ be the category of
 triples $(\ovl{A},\Lambda, \epsilon)$ consisting of an abelian scheme
 over $S_{\mathrm{red}}$, a $p$-divisible group $\Lambda$ over $S$ and an
 isomorphism $\epsilon: \ovl{A}[p^\infty]\cong \Lambda\times_S S_{\mathrm{red}}$.

\begin{thm}[The Serre-Tate theorem]\label{thm:2.2.4}
 The functor
 \begin{equation*}
  A\;\mapsto \;(A_{\mathrm{red}},\;\Lambda=A[p^\infty],
   \;\epsilon: A_{\mathrm{red}}[p^\infty]\cong \Lambda\times_S S_{\mathrm{red}})
 \end{equation*}
 induces an equivalence of categories between the
 category of abelian schemes over $S$ and the category
 $\Def_{p\text{-div}}(S)$.
\end{thm}

We again follow Drinfeld's argument \cite{copsr}, \cf also \cite{slm}.
We may assume that $S$ is affine and that
$\mathcal{J}^{n+1}=0$. We begin with the full faithfulness: given a
homomorphism $f^\infty: A[p^\infty]\rightarrow B[p^\infty]$ such that
$f^\infty$ $\pmod{\mathcal{J}}$ comes from a homomorphism $\ovl{f} :
A_{\mathrm{red}}\rightarrow B_{\mathrm{red}}$ of abelian schemes over $S_{\mathrm{red}}$, we have
to show that there exists a unique homomorphism $f: A\rightarrow B$
compatible with $f^\infty$ and $\ovl{f}$. The unicity follows from the
injectivity of $\Hom(A,B)\rightarrow \Hom(A_{\mathrm{red}},B_{\mathrm{red}})\;$
(\ref{thm:2.2.3}). Furthermore, \ref{thm:2.2.3} shows the existence of a
homomorphism $g$ 
compatible with $\;p^nf^\infty$ and $p^n\ovl{f}$. To show the existence
of an $f$ such that $g=p^nf=f\circ [p^n]$, the point is to show that $g$
kills $A[p^n]$, which can be seen on $g^\infty$.

For the essential surjectivity, we have to construct an abelian scheme
$A$ over $S$ compatible with a given datum $(\ovl{A},\Lambda, \epsilon)$
in $\Def_{p\text{-div}}(S)$.  By unicity of liftings, it suffices to do so
locally on $S$.  hence we may assume that $S$ is affine, and choose an
abelian scheme $A'$ over $S$ and an isogeny $A'_{\mathrm{red}}\rightarrow
\ovl{A}$. By the above argument, $p^n$ times the given isogeny
$A'_{\mathrm{red}}\rightarrow \ovl{A}$ lifts to a (unique) isogeny
$e:\;A'[p^\infty ]\rightarrow \Lambda$ of $p$-divisible groups over
$S$. Its kernel is a finite locally free subgroup scheme of $A'$, and we
can form the abelian scheme $A'':=A'/\Ker e$. The isogeny $A''[p^\infty
]\rightarrow \Lambda$ induced by $e$ is an isomorphism
$\pmod{\mathcal{J}}$, hence is an isomorphism.

\begin{para}\label{para:2.2.5} In the sequel, we shall have to deal
 occasionally with $p$-divisible groups over formal schemes rather 
than over schemes. Our formal schemes $\mathfrak{X}$ will be adic,
locally noetherian; hence there is a largest ideal of definition
$\mathcal{J}\subset \mathcal{O}_{\mathfrak{X}}$, and
$\mathfrak{X}=\liminj \;{\mathfrak{X}}_n\;,\;{\mathfrak{X}}_n =
\Spec(\mathcal{O}_{\mathfrak{X}}/\mathcal{J}^{n+1})\;\;
({\mathfrak{X}}_{\mathrm{red}}={\mathfrak{X}}_0)\;$.

A $p$-divisible group $\Lambda$ over $\mathfrak{X}$ will be an inductive
system of $p$-divisible groups $\Lambda_n$ over $\mathfrak{X}_n$ such
that $\Lambda_{n+1}\times_{{\mathfrak{X}}_{n+1}}{\mathfrak{X}}_n\cong
\Lambda_n$.

\newpage
\section{A stroll in the crystalline world.}\label{sec:3}

\begin{abst}
 Review of Dieudonn\'e modules, crystals, Grothendieck-Messing theory,
 convergent isocrystals.
\end{abst}
\end{para}

\subsection{From Dieudonn\'e modules to crystals.}\label{sub:3.1}

\begin{para}\label{para:3.1.1}
 At first, Dieudonn\'e theory presents itself as an analogue of Lie
 theory for formal Lie groups over a perfect field $k$ of characteristic
 $p>0$. In practice, its scope is limited to the case of {\it
 commutative} formal Lie groups (in which case the classical counterpart
 is not very substantial). The theory extends to the case of
 $p$-divisible groups over $k$.

 Let $W$ be the Witt ring of $k$, with its Frobenius automorphism
 $\sigma$. The theory associates to any 
 $p$-divisible group $\Lambda$ over $k$ a {\it Dieudonn\'e module}, \ie a
 finitely generated free module over $W$ 
 endowed with a $\sigma$-linear Frobenius $F$ and a $\sigma^{-1}$-linear
 Verschiebung 
 $V$ satisfying $FV=VF=p.\id$ (we work with the contravariant theory: 
 ${\D}(\Lambda)=\Hom(\Lambda , \mathcal{W}_k)$, 
 \index{000D@$\D$}
 where $\mathcal{W}_k\;$
 is the Witt scheme of $k$). This provides an 
 {\it anti-equivalence of categories between $p$-divisible groups over
 $k$ and Dieudonn\'e modules}. The rank of 
 ${\D}(\Lambda)$ is the height $h$ of $\Lambda$. 
\end{para}

\begin{exs}
  $\bullet \;\;{\D}(\Q_p/\Z_p)=W,$ with $F=\sigma,
 \;V=p\sigma^{-1}$. 
 \par\noindent $\bullet\;\; {\D}(\hat\G_m)=W(-1)$, \ie $W$ as an
 underlying module, with $F=p\sigma, 
 \;V=\sigma^{-1}$.
 \par\noindent $\bullet$ When
 $\Lambda=A[p^\infty]$ for an abelian variety $A$ over $k$,
 ${\D}(\Lambda)$ may be identified with 
 $\mathrm{H}^1_{\mathrm{cris}}(A/W)$. 
\end{exs}

 However, the structure of Dieudonn\'e modules is rather complicated in
 general unless one inverts $p$. 

\begin{para}\label{para:3.1.2}
 When $k$ is algebraically closed, Dieudonn\'e modules $\otimes \Q$ were
 completely classified by Dieudonn\'e (this provides a classification of
 $p$-divisible groups over $k$ up to isogeny). They are classified by
 their {\it slopes} $\lambda\in [0,1]\cap \Q$ and multiplicities
 $m_\lambda\in \Z_{>0}$ (these data are better recorded in the form of a
 Newton polygon). This comes as follows: each Dieudonn\'e module $\D$
 has a canonical increasing finite filtration by Dieudonn\'e submodules
 $\D_\lambda$, such that the associated graded $W$-module is free, and
 with the following property: if $\lambda=a/b$ in irreducible form,
 $\Gr_\lambda\otimes \Q$ admits a basis of $m_\lambda$ vectors $x$
 satisfying $F^bx=p^ax$.

 For instance, the slope of $W(-1)$ is $1$; the slope of the Dieudonn\'e
 module attached to 
 (the $p$-divisible group of) a supersingular elliptic curve over $k$ is
 $1/2$, with multiplicity $2$.
\end{para}

\begin{para}\label{para:3.1.3} The problem of generalizing Dieudonn\'e
 theory to $p$-divisible groups over more general bases $S$ (over which
$p$ is nilpotent) has been tackled and advertised by Grothendieck
 \cite{grothendieck74:_group_barsot_tate_dieud}, and further studied by many geometers.

 What should be the right substitutes for the $W$-modules ${\D}(\Lambda)$
 (\resp for the $W[\frac{1}{p}]$-spaces ${\D}(\Lambda)_{\Q}$)?

 Grothendieck's proposal was to define ${\D}(\Lambda)$ as a {\it
 $F$-crystal} on the crystalline site of $S$. 

 The right substitute for $W[\frac{1}{p}]$-spaces was proposed later by
 P. Berthelot and A. Ogus under the name of {\it convergent
 $F$-isocrystal}. Let us describe these notions.
\end{para}

\subsection{The Dieudonn\'e crystal of a $p$-divisible group.}\label{sub:3.2}

\begin{para}\label{para:3.2.1} Let $\mathfrak{v}$ be a finite extension of
 $W$ of degree $e$. We assume that $e<p$, and that $\varpi$ is a 
 uniformizing parameter of $\mathfrak{v}$ such that $\sigma$ extends to 
 $\mathfrak{v}$ by setting $\sigma(\varpi)=\varpi$.

Let $S$ be in $\Nil_\mathfrak{v}$. An
$S$-divided power thickening
$T_0\hookrightarrow T$ is given by
\begin{enumerate}
 \item an $S$-scheme $T_0$,
 \item a $\mathfrak{v}$-scheme $T$ on which $p$ is locally nilpotent,
 \item a closed immersion $T_0\hookrightarrow T$ over $\mathfrak{v}$,
 \item a collection of maps ``$\frac{x^n}{n!}$'' from the ideal of
 definition of $T_0$ in $T$ to $\mathcal{O}_T$, which 
 satisfy the formal properties of the divided powers. 
\end{enumerate}

 $S$-divided power thickenings form a category
 $(S/\mathfrak{v})_{\mathrm{cris}}$ (a fppf site). A sheaf on
 $(S/\mathfrak{v})_{\mathrm{cris}}$ is the data, for every
 $T_0\hookrightarrow T$, of a fppf sheaf $\mathcal{E}_{T_0,T}$ on $T$,
 and for every morphism $f: (T_0\hookrightarrow T)\rightarrow
 (T_0'\hookrightarrow T')$ in $(S/\mathfrak{v})_{\mathrm{cris}}$, of a
 homomorphism $f^\ast\mathcal{E}_{T_0,T}\rightarrow
 \mathcal{E}_{T_0',T'}$ satisfying the obvious transitivity
 condition. Example: $\mathcal{O}_{S_0/\mathfrak{v}}$, defined by
 $(\mathcal{O}_{S/\mathfrak{v}})_{T_0,T}=\mathcal{O}_T$. Another example:
 any fppf sheaf $\mathcal{F}$ on $S_0$ gives rise to a sheaf
 $\mathcal{F}_{\mathrm{cris}}$ on $(S/\mathfrak{v})_{\mathrm{cris}}$ by
 $(\mathcal{F}_{\mathrm{cris}})_{T_0,T}=\mathcal{F}_{T_0}$.

 A {\it crystal} on $(S/\mathfrak{v})_{\mathrm{cris}}$ (or,
 abusively, on $S$) is a sheaf $\mathcal{E}$ of
 $\mathcal{O}_{S/\mathfrak{v}}$-modules on
 $(S/\mathfrak{v})_{\mathrm{cris}}$ such that the homomorphisms
 $f^\ast\mathcal{F}_{T_0,T}\rightarrow \mathcal{F}_{T_0',T'}$ are
 isomorphisms (in Grothendieck's terms: ``crystals grow and are
 rigid''). One denotes by $\mathcal{E}_T$ the $\mathcal{O}_T$-module
 obtained by evaluating the crystal on $(\id: T\hookrightarrow T) \in
 (S/\mathfrak{v})_{\mathrm{cris}}$.

 Because of this rigidity, and using the canonical divided powers on the
 ideal $\varpi\mathcal{O}_S$ (since $e<p,\;\frac{\varpi^n}{n!}\in
 \mathfrak{v}$ for any $n$), one gets an equivalence of categories
 between crystals on $S$ and crystals on the scheme $S_0=S/\varpi$ of
 characteristic $p$.
\end{para}

\begin{para}\label{para:3.2.2}
 The Dieudonn\'e crystal attached to a $p$-divisible group $\Lambda$
 over $S$ is the crystalline $\mathcal{E}xt$-sheaf
 $\mathcal{E}xt^1(\Lambda_{\mathrm{cris}},
 \mathcal{O}_{S/\mathfrak{v}})$. It can be shown that ${\D}(\Lambda)$ is
 a finite locally free crystal on $S$, of rank the height of
 $\Lambda$. It depends functorially in $\Lambda$ (in a contravariant
 way).

 An important, though formal, consequence of the crystalline local
 character of this definition is that {\it the formation of
 ${\D}(\Lambda)$ commutes with base change in}
 $(S/\mathfrak{v})_{\mathrm{cris}}$, and for any $S$-divided power
 thickening $(T_0\hookrightarrow T)$, there is a canonical isomorphism
 of $\mathcal{O}_T$-modules ${\D}(\Lambda)_{T_0,T}\cong
 {\D}(\Lambda_{T_0})_{T_0,T}.$

 It follows that the datum of ${\D}(\Lambda)$ amounts to the datum of
 ${\D}(\Lambda\times_S S_0)$: the Dieudonn\'e crystal depends only on the
 $p$-divisible group modulo $\varpi$.
\end{para}

\begin{para}\label{para:3.2.3}
 The pull-back of a crystal $\mathcal{E}$ on $S_0$ via the absolute
 Frobenius $F_{S_0}$ is denoted by $\mathcal{E}^{(p)}$. An $F$-crystal
 structure on $\mathcal{E}$ is the datum of a morphism of crystals
 $F:\;\mathcal{E}^{(p)}\rightarrow \mathcal{E}$ which admits an inverse
 up to a power of $p$.

 If $\Lambda$ is a $p$-divisible group over $S_0$, ${\D}(\Lambda)$ is
 endowed by functoriality with a structure of $F$-crystal. The classical
 Dieudonn\'e module is obtained for
 $S_0=T_0=\Spec(k),\;T=\Spf(W)=\liminj\;\Spec(W/p^n)$.

 Although we shall not use this fact, let us mention that the
 ``Dieudonn\'e functor'' ${\D}$ from $p$-divisible groups on $S_0$ to
 $F$-crystals is fully faithful under further mild assumptions on $S_0$
 (and even an equivalence with values in the category of ``Dieudonn\'e
 crystals'', \ie finite locally free $F$-crystals for which $F$ admits an
 inverse up to $p$ --- the Verschiebung); without any further assumption on
 $S_0$, it induces an equivalence of categories between ind-etale
 $p$-divisible groups on $S_0$ and unit-root $F$-crystals (\ie finite
 locally free $F$-crystals for which $F$ is an isomorphism), 
 \cf \cite{tddc1}, \cite{jong95:_cryst_dieud}.
\end{para}

\subsection{The Hodge filtration.}\label{sub:3.3}

\begin{para}\label{para:3.3.1}
 Let $\Lambda^\vee$ be the (Serre) dual of $\Lambda$ (the $\liminj$ of
 the Cartier duals of the $\Lambda[p^n]$'s). There is a canonical,
 locally split, exact sequence of $\mathcal{O}_S$-modules
 \begin{equation*}
  0\rightarrow F^1\rightarrow {\D}(\Lambda)_S\rightarrow
   \Lie \Lambda^\vee\rightarrow 0
 \end{equation*}
 where $F^1=\CoLie \Lambda=\omega_\Lambda$ may be identified with the
 module of invariant differentials on the formal Lie group
 $\hat\Lambda$. In case $\Lambda$ is the $p$-divisible group attached to
 an abelian $S$-scheme $\udl{A}$, this exact sequence reduces to the
 Hodge exact sequence
 \begin{equation*}
  0\rightarrow F^1\rightarrow \mathrm{H}^1_{\mathrm{DR}}(\udl{A}/S)\rightarrow
   \Lie \udl{A}^\vee\rightarrow 0.
 \end{equation*}
\end{para}

\begin{para}\label{para:3.3.2}
 Now let $S_0\hookrightarrow S$ be a {\it nilpotent} divided power
 thickening (this means that the products
 ``$\frac{x_1^{n_1}}{{n_1}!}$''$\ldots$``$\frac{x_1^{n_m}}{{n_m}!}$'' vanish
 for $n_1+\ldots+n_m\gg 0$); this is for instance the case if
 $S_0=S/\varpi$ as before and $e<p-1$.

 The next theorem asserts that deformations of $p$-divisible
 groups are controlled by the variation of the Hodge filtration in the
 Dieudonn\'e module \cite{messing72:_barsot_tate}. For any $p$-divisible group
 $\ovl{\Lambda}$ on $S_0$, and any lifting $\Lambda$ over $S$, we have
 seen that ${\D}(\Lambda)_S$ identifies with ${\D}(\ovl{\Lambda})_S$,
 which is a lifting of ${\D}(\ovl{\Lambda})_{S_0}$ to $S$.
\end{para}

\begin{thm}[The Grothendieck-Messing theorem]\label{thm:3.3.3}
 The functor
 \begin{equation*}
  \Lambda\;\mapsto (\Lambda\times_S S_0,\;\omega_\Lambda)
 \end{equation*}
 induces an equivalence of categories between the category of
 $p$-divisible groups on $S$, and the category of pairs consisting of a
 $p$-divisible group $\ovl{\Lambda}$ over $S_0$ (Zariski-locally liftable
 to $S$) and a locally direct factor vector subbundle of
 ${\D}(\ovl{\Lambda})_S$ which lifts $\;\omega_{\ovl{\Lambda}}\subset
{\D}(\ovl{\Lambda})_{S_0}$.
\end{thm}

 We have already seen the faithfulness in \ref{thm:2.2.3}.  We just
 indicate here 
 the principle of proof of essential surjectivity and fullness. Working
 with nilpotent divided power thickenings allows one to construct, by
 the exponential method, a formally smooth group scheme
 $\;\E(\ovl{\Lambda})\;$ over $\;S\;$ together with a canonical
 isomorphism
 ${\D}(\ovl{\Lambda})_{S_0}=\Lie (\E(\ovl{\Lambda})\times_S{S_0})$, in
 such a way that the canonical exact sequence

 \begin{equation*}
  0\rightarrow \omega_{\ovl{\Lambda}}\rightarrow
   {\D}(\ovl{\Lambda})_{S_0}\rightarrow 
   \Lie \ovl{\Lambda}^\vee\rightarrow 0
 \end{equation*}
 lifts to an exact sequence of formally smooth group schemes
 \begin{equation*}
  0\rightarrow \omega_{\ovl{\Lambda}}\rightarrow
   \E(\ovl{\Lambda})\times_S{S_0}\rightarrow \ovl{\Lambda}^\vee\rightarrow 0
 \end{equation*}
 Via the exponential, it is equivalent to give a local summand $F^1$ of
 ${\D}(\ovl{\Lambda})_S$ which lifts
 $\;\omega_{\ovl{\Lambda}}$, or to give a vector sub-group scheme $V$ of
 $\E(\ovl{\Lambda})$ which lifts $\;\omega_{\ovl{\Lambda}}$. It then
 turns out that $\E(\ovl{\Lambda})/V$ is a 
 $p$-divisible group $\Lambda^\vee$, whose Serre dual is the looked for 
 $p$-divisible group $\Lambda$ over $S$: 
 $\Lambda\times_S S_0=\ovl{\Lambda}$ by duality, and $\omega_{\Lambda}=F^1$.

 The formation of $\E(\ovl{\Lambda})$ is functorial (contravariant in 
 $\ovl{\Lambda}$). If we have two $p$-divisible groups $\ovl{\Lambda}, 
 \ovl{\Lambda}'$ and a morphism $\ovl{f}:\;\ovl{\Lambda}=\Lambda\times_S 
 S_0\rightarrow \ovl{\Lambda}'=\Lambda'\times_S S_0$ such that 
 $\D(f)$ respects the $F^1$'s, then the homomorphism 
 $\E(\ovl{f}): \E(\ovl{\Lambda})\rightarrow \E(\ovl{\Lambda}')$ 
 preserves the vector subgroup schemes, hence induces a homomorphism 
 $f: \E(\ovl{\Lambda})/V\cong \Lambda\rightarrow 
 \E(\ovl{\Lambda}')/V'\cong\Lambda'$ which lifts $\ovl{f}$.

 The exponential construction becomes simpler when the ideal of $S_0$ in
 $S$ is of square zero (and endowed with the trivial divided
 powers). This is the case which will be used in the sequel.

\subsection{Crystals and connections.}\label{sub:3.4} 

\begin{para}\label{para:3.4.1}
 We have already encountered the terms ``$F$-crystal'', ``unit-root
 $F$-crystal'' in I.\ref{para:f-cris}, \ref{para-unit-root} in the more
 down-to-earth context of 
 connections. This comes as follows. Assume for simplicity that $S_0$ is
 affine of characteristic $p$ and lifts to a formally smooth $p$-adic
 affine formal scheme $\mathcal{S}$ over $\mathfrak{v}=W$.  Then there
 is an equivalence of categories between finite locally free crystals on
 $(S_0/W)_{\mathrm{cris}}$ and finite locally free
 $\mathcal{O}(\mathcal{S})$-modules with integrable topologically
 nilpotent connection (this interpretation of crystals extends to much
 more general situations, \cf \cite[2.2.2]{jong95:_cryst_dieud}).

 \noindent Using this equivalence of categories, the Katz (covariant) functor
 \begin{equation*}
  \begin{pmatrix}\text{continuous }{\Z}_p\text{-representations} \\
   \text{of } \pi_1^{\mathrm{alg}}(S_0,s_0)
  \end{pmatrix}
  \longrightarrow \text{(unit-root $F$-crystals)}
 \end{equation*}
 is the composite of the following (anti)equivalences of categories
 \begin{align*}
  \begin{pmatrix}\text{continuous } {\Z}_p\text{-representations} \\
   \text{of } \pi_1^{\mathrm{alg}}(S_0,s_0)
  \end{pmatrix}& \longrightarrow
  \text{(ind-etale $p$-divisible groups over $S_0$)}\\
  & \overset{\D}{\longrightarrow}
  \text{(unit-root $F$-crystals over $S_0$)}.
 \end{align*} 
\end{para}

\begin{para}\label{para:3.4.2}
 Let us briefly recall the construction of the connection attached to a
 finite locally free crystal $\mathcal{E}$ on
 $(S_0/W)_{\mathrm{cris}}$. It depends on the technique of formally
 adding divided powers to an ideal (the pd-hull). Let $\Delta_0$ be the
 diagonal in $S_0\times S_0$, and let $\hat\Delta$ be the $p$-completion
 of the pd-hull of the diagonal in $\mathcal{S}\hat\times\mathcal{S}$,
 endowed with the two projections $p_1,p_2$ to $\mathcal{S}$. Then by
 rigidity of crystals, we have for any $n$ an isomorphism
 \begin{equation*}
  p_2^\ast(\mathcal{E}_{S_0,\mathcal{S}/p^n})\;\cong
   \;\mathcal{E}_{\Delta_0,\hat\Delta/p^n}\;\cong
   \;p_1^\ast(\mathcal{E}_{S_0,\mathcal{S}/p^n})
 \end{equation*}
 and at the limit $n\rightarrow \infty$, an isomorphism of finite
 locally free $\mathcal{O}({\hat\Delta})$-modules  
 \begin{equation*}
  \epsilon:\; p_2^\ast(\mathcal{E}_{S_0,\mathcal{S}})\;\cong
   \;p_1^\ast(\mathcal{E}_{S_0,\mathcal{S}})
 \end{equation*}
 satisfying a cocycle condition. Such an isomorphism is the ``Taylor series'' 
 \begin{equation*}
  \epsilon(1\otimes e) = \sum_{n_1,\ldots,n_d\geq
   0}\biggl(\nabla\Bigl(\frac{\partial}{\partial t_1}\Bigr)^{n_1}
   \ldots\nabla\Bigl(\frac{\partial}{\partial t_d}\Bigr)^{n_d}\biggr)(e)
   \otimes \prod_i \frac{\text{``$(1\otimes t_i-t_i\otimes 1)^{n_i}$''}}{n_i!}
 \end{equation*}
 of an integrable connection $\nabla$ on the finite
 locally free $\mathcal{O}(\mathcal{S})$-module
 $\;\mathcal{E}_{S_0,\mathcal{S}}$ (here $t_1,\ldots t_d$ denote local
 coordinates on $\mathcal{S}$).
\end{para}

\begin{para}\label{para:3.4.3}
 In general, the Taylor series converges only in polydiscs of radius
 $| p |^{1/p-1}$, due to the presence of the factorials. One
 says that $\nabla$ is {\it convergent} if its Taylor series converges
 in (open) unit polydiscs. A well-known argument due to Dwork shows that
 $F$-crystals give rise to convergent connections $\nabla$ (\cite{tdd}).

 In order to put the discussion of convergence on proper foundations, we
 need the notion of tube in analytic geometry.
\end{para}

\subsection{Interlude : analytic spaces associated with formal schemes
 and tubes.}\label{sub:3.5}

\begin{para}\label{para:3.5.1} Let us first review the Raynaud-Berthelot
construction of the analytic space attached to a  
$\mathfrak{v}$-formal scheme $\mathfrak{X}$. It will be convenient for
later purpose to formulate the construction in the frame of Berkovich
spaces rather than rigid geometry (we refer to \cite[1.6]{ecfnas} for the
translation).

As in \ref{para:2.2.5}, we assume that $\mathfrak{X}$ is adic, locally
noetherian; moreover, we assume that 
${\mathfrak{X}}_{\mathrm{red}}=\Spec(\mathcal{O}_{\mathfrak{X}}/\mathcal{J})$ is a
separated $k$-scheme locally of finite type (in particular,
$p\in \mathcal{J}$). However, we {\it do not} assume that $\mathfrak{X}$ is
$p$-adic, \ie that the topology of $\mathcal{O}_{\mathfrak{X}}$
is $p$-adic.

 We describe the construction in the affine case ${\mathfrak{X}} =
 \Spf(\mathcal{A})$; the general case follows by gluing. One chooses
 generators $f_1,\ldots, f_r$ of 
 $J=\Gamma({\mathfrak{X}},\mathcal{J})$, and defines, for any $n>0$,
 \begin{equation*}
 \mathcal{B}_n=\mathcal{A}\{T_1,\ldots,T_r\}/(f_1^n-\varpi T_1,\ldots,
 f_r^n-\varpi T_r)
 \end{equation*}
 where $\mathcal{A}\{T_1,\ldots,T_r\}$ stands for the $\varpi$-adic
 completion of $\mathcal{A}[T_1,\ldots,T_r]$. One has
 $\mathcal{B}_n/\varpi\mathcal{B}_n\cong
 (\mathcal{A}/(\varpi,f_1^n,\ldots,f_r^n))[T_1,\ldots,T_r]$, a
 $k$-algebra of finite type; hence $\mathcal{B}_n$ is a
 $\mathfrak{v}$-algebra topologically of finite type (\ie quotient of an
 algebra of the form $\mathcal{A}\{T_1,\ldots,T_{r_n}\}$), so that
 $\mathcal{B}_n[\frac{1}{p}]$ is an affinoid algebra. Let
 $\matheur{M}(\mathcal{B}_n)$ be its Berkovich spectrum. The
 homomorphism $\mathcal{B}_{n+1}\rightarrow \mathcal{B}_n$ sending $T_i$
 to $f_iT_i$ identifies $\matheur{M}(\mathcal{B}_n)$ with an affinoid
 subdomain of $\matheur{M}(\mathcal{B}_{n+1})$ \cite[0.2.6]{crecrasp}.

 The analytic space $\mathfrak{X}^{\mathrm{an}}$ is defined as the union
 $\bigcup_n \matheur{M}(\mathcal{B}_n)$. It does 
 not depend on the choice of $f_1,\ldots, f_r$. This is a {\it
 paracompact (strictly) analytic space} over the fraction 
 field $K$ of $\mathfrak{v}$  (one can replace the nested ``ball-like''
 $\matheur{M}(\mathcal{B}_n)$ by nested ``annulus-like'' 
 spaces in order to obtain locally finite coverings). 
\end{para}

\begin{exs}\label{para:3.5.2}
 (i) if ${\mathfrak{X}}=\Spf(\mathfrak{v}\{t\})$,
 ${\mathfrak{X}}^{\mathrm{an}}$ is the closed unit disc.
 
 (ii) if ${\mathfrak{X}}=\Spf({\mathfrak{v}[[t]]}),
 \;\mathcal{J}=(\varpi, t)$, then $\matheur{M}(\mathcal{B}_n)$ is the 
 closed disc $\matheur{D}_K(0,|\varpi|^{1/n})$, and
 ${\mathfrak{X}}^{\mathrm{an}}$ is the open unit disc. 
\end{exs}

\begin{para}\label{para:3.5.3} Let $S_0$ be a closed subscheme of
$\mathfrak{X}_{\mathrm{red}}$, and let ${\mathfrak{X}}_{S_0}$ be the formal
completion of $\mathfrak{X}$ along $S_0$. The {\it tube of $S_0$ in
$\mathfrak{X}$} is $\mathfrak{X}_{S_0}^{\mathrm{an}}$. It is denoted by
$]S_0[_{\mathfrak{X}}$ (and considered as a subspace of
$\mathfrak{X}^{\mathrm{an}}$). It depends only on $S_{0,\mathrm{red}}$.  

\medskip
 \noindent {\it Example}: if $\mathfrak{X}$ is as in example (i) above
 (so that $\mathfrak{X}_{\mathrm{red}}=\A^1_k$), and $S_0$ is the point $0$,
 then $\mathfrak{X}_{S_0}$ is as in example (ii) above:
 $]0[_{\mathfrak{X}}\;=\matheur{D}(0,1^-)\subset \matheur{D}(0,1^+)$. 

 \noindent More generally, the polydiscs mentioned in \ref{para:3.4.3}
 are the tubes of closed points.
\end{para}

\subsection{Convergent isocrystals.}\label{sub:3.6}     

\begin{para}\label{para:3.6.1} Let $S_0$ be a separated $k$-scheme of finite
 type, and $S_0\hookrightarrow \mathfrak{Y}$ a closed immersion into a
 flat $p$-adic formal $\mathfrak{v}$-scheme that is formally smooth in a
 neighborhood of $S_0$. Let $E$ be a vector bundle on the tube
 $]S_0[_{\mathfrak{Y}}$. An integrable connection $\nabla$ on $E$ is
 {\it convergent} if its Taylor series
 \begin{equation*}
  \epsilon(1\otimes e) =
   \sum_{n_1,\ldots,n_d\geq 0}
   \biggl(\nabla\Bigl(\frac{\partial }{\partial t_1}\Bigr)^{n_1}\ldots
   \nabla\Bigl(\frac{\partial }{\partial t_d}\Bigr)^{n_d}\biggr)(e)
   \otimes \;\prod_i\frac{\text{``$(1\otimes t_i-t_i\otimes 1)^{n_i}$''}}{n_i!}
 \end{equation*}
 is induced by an isomorphism on the tube of the
 diagonal $]S_0[_{\mathfrak{Y}\times \mathfrak{Y}}$
 \begin{equation*}
  \epsilon:\; p_2^\ast(E)\;\cong \;p_1^\ast(E)
 \end{equation*} 
 It amounts to the same to say that for any affine $U$ in $\mathfrak{Y}$
 with local coordinates
 $t_1,\ldots, t_d$, any section $e\in \Gamma(U^{\mathrm{an}},E)$ and any $\eta<1$,
 one has 
 \begin{equation*}
  \Bigl\|\nabla\Bigl(\frac{\partial }{\partial t_1}\Bigr)^{n_1}
   \ldots\nabla\Bigl(\frac{\partial }{\partial t_d}\Bigr)^{n_d}(e)\Bigr\|.
   \eta^{n_1+\ldots + n_d}\rightarrow 0
 \end{equation*}
 where $\|\;\;\|$ denotes a Banach norm on $\Gamma(U^{\mathrm{an}},E)$.

 According to Berthelot, {\it the category of vector bundles
with convergent connection depends only on $S_0$} --- even only on
$S_{0,\mathrm{red}}$ --- (up to canonical equivalence), and is functorial in
$S_0/K$; this is the category of {\it convergent isocrystals on $S_0/K$}
(or, abusively, on $S_0$).

Since any separated $k$-scheme locally of finite type can be locally
embedded into a $\mathfrak{X}$ as above, one can define the category of
convergent isocrystals on such a scheme by gluing \cite[2.3.2]{crecrasp}.
\end{para}

\begin{para}\label{para:3.6.2} There is a natural functor from the category
 of $F$-isocrystals on $S_0/\mathfrak{v}$ to the category of 
convergent isocrystals on $S_0/K$ \cite[2.4]{crecrasp}.  The composite of
this functor with the Dieudonn\'e functor is denoted by $\D(-)_\Q$. It
factors through a $\Q_p$-linear functor, still denoted by $\D(-)_\Q:\;$
\begin{equation*}
\Bigl\{\begin{matrix}
	\text{$p$-divisible groups over $S_0$}\\
	\text{up to isogeny}
       \end{matrix}\Bigr\}
 \longrightarrow
 \{\text{convergent isocrystals on $S_0$}\}
\end{equation*}
\end{para}

\begin{para}\label{para:3.6.3} The advantage of convergent isocrystals over
 ``crystals up to isogeny'' lies in their strong 
``continuity property'', which allows to remove the restriction $e<p$ on
ramification (as we shall do henceforth), and work conveniently with
singular schemes $S_0$ (while divided powers present certain
``pathologies'' in the singular case). In our applications, $S_0$ will
be a countable (non-disjoint) union of irreducible projective varieties,
and $\mathfrak{v}$ may be highly ramified.
\end{para}

\begin{para}\label{para:3.6.4} If $S_0$ is the reduced underlying scheme 
 ${\mathfrak{X}}_{\mathrm{red}}$ of a formal scheme $\mathfrak{X}$ as in
 \ref{para:3.5.1}, there 
 is a notion of {\it evaluation of $E$ on $\mathfrak{X}^{\mathrm{an}}$}
 \cite[2.3.2]{crecrasp}: this is a vector bundle on
 ${\mathfrak{X}}^{\mathrm{an}}$ with an integrable connection.

 By functoriality of the construction, it will be enough to deal with the
 affine case ${\mathfrak{X}}=\Spf{\mathcal{A}}$ and to assume that there
 is a closed immersion $S_0=(\Spf{\mathcal{A}})_{\mathrm{red}}\hookrightarrow
 \mathfrak{Y}$ into a formally smooth $p$-adic formal
 $\mathfrak{v}$-scheme. One can construct inductively a compatible family
 of $\mathfrak{v}$-morphisms
 $\Spec(\mathcal{O}_{\mathfrak{X}}/\mathcal{J}^n)\rightarrow \mathfrak{Y}$,
 whence a $\mathfrak{v}$-morphism  $\mathfrak{X} \rightarrow
 \mathfrak{Y}$.  Let $\upsilon :\; ]S_0[_{\mathfrak{X}} 
 ={\mathfrak{X}}^{\mathrm{an}}\rightarrow\; ]S_0[_{\mathfrak{Y}}$ be the
 associated analytic map. The looked for evaluation of $E$ on
 ${\mathfrak{X}}^{\mathrm{an}}$ 
 is the pull-back by $\upsilon$ of the vector bundle with integrable
 connection on $]S_0[_{\mathfrak{Y}}$ defined by $E$.
\end{para}

\begin{para}\label{para:3.6.5} Let now $S$ be a separated scheme in
 $\Nil_\mathfrak{v}$ (it may also be viewed in  $\Nil_W$).
 Let $\mathcal{E}$ be a finite locally free $F$-crystal on $S/W$. This
 provides a vector bundle $\mathcal{E}_S$ on $S$, which 
 only depends on the inverse image $\mathcal{E}\times_W k$ of
 $\mathcal{E}$ on the $k$-scheme $S/p$ (\cf \ref{para:3.2.1}). On the other 
 hand, the convergent isocrystal $E$ on $S/p$ attached to $F$-crystal
 $\mathcal{E}\times_W k$ (\cf \ref{para:3.6.2}) depends only on
 $S_{\mathrm{red}}$.
\end{para}

\begin{para}\label{para:3.6.6} Let us apply this to Dieudonn\'e
 crystals. Let $\mathfrak{X}$ be a formal scheme as in \ref{para:3.5.1},
 and let 
 $\Lambda=\liminj \Lambda_n$ be a $p$-divisible on
 $\mathfrak{X}=\liminj {\mathfrak{X}}_n$. This gives rise to a vector bundle 
 $\displaystyle \D(\Lambda)_{\mathfrak{X}}= \liminj
 \D(\Lambda_n)_{{\mathfrak{X}}_n}$ on $\mathfrak{X}$. On the other 
 hand, we have the vector bundle $\displaystyle
 \D(\Lambda_0)_{\mathfrak{X}}=$ on $\mathfrak{X}$ and a natural
 homomorphism $\D(\Lambda_0)_{\mathfrak{X}}\rightarrow
 \D(\Lambda)_{\mathfrak{X}}$ obtained by pull-back 
 ${\mathfrak{X}}_{\mathrm{red}}={\mathfrak{X}}_0\rightarrow {\mathfrak{X}}$.

When $e\geq p$, this is not an isomorphism in general. However, {\it the
associated morphism of analytic vector bundles
$(\D(\Lambda_0)_{\mathfrak{X}})^{\mathrm{an}}\rightarrow
(\D(\Lambda)_{\mathfrak{X}})^{\mathrm{an}}$ is an isomorphism.  Moreover,
$(\D(\Lambda_0)_{\mathfrak{X}})^{\mathrm{an}}$ is the analytic vector bundle
underlying the evaluation of $\D(\Lambda_0)_\Q$ on
${\mathfrak{X}}^{\mathrm{an}}$}.

Indeed, working locally, and taking in account the very construction of
${\mathfrak{X}}^{\mathrm{an}}$, we reduce to the case when ${\mathfrak{X}}$ is
$p$-adic. Then $\mathfrak{X}/p$ is a $k$-scheme. Using the canonical
divided powers on $(p)$ and the rigidity of crystals,
$\D(\Lambda\times_{\mathfrak{X}}\mathfrak{X}/p)_{\mathfrak{X}}\cong
\D(\Lambda)_{\mathfrak{X}}$. It is thus enough to see that
$\D(\Lambda_0)^{\mathrm{an}}\cong\D
(\Lambda\times_{\mathfrak{X}}\mathfrak{X}/p)_{\mathfrak{X}}^{\mathrm{an}}$. This
follows from the fact that the convergent isocrystals $\D(\Lambda_0)_\Q$
and $\D(\Lambda\times_{\mathfrak{X}}\mathfrak{X}/p)_\Q$ are ``the
same'', since ${\mathfrak{X}}_0=({\mathfrak{X}}/p)_{\mathrm{red}}\;$: indeed, it
is clear that $(\D(\Lambda_0)_{\mathfrak{X}})^{\mathrm{an}}$ (\resp
$\D(\Lambda\times_{\mathfrak{X}}\mathfrak{X}/p)_{\mathfrak{X}}^{\mathrm{an}}$)
coincides with the evaluation of the convergent isocrystal
$\D(\Lambda_0)_\Q$ (\resp
$\D(\Lambda\times_{\mathfrak{X}}\mathfrak{X}/p)_\Q$) on
${\mathfrak{X}}^{\mathrm{an}}$ (\cf also \cite[6.4]{efgonas}, in a slightly less
general context).
\end{para}

\begin{para}\label{para:3.6.7} Let us draw from this two consequences:
 \begin{enumerate}
  \item taking into account the fact that the functor $\D(-)_\Q$ factors
	through the category of $p$-divisible
	groups up to isogeny, we see that if we are given a quasi-isogeny
	$\Lambda_0\rightarrow \Lambda_0'$, then there is a
	canonical isomorphism of vector bundles
	$\bigl(\D(\Lambda'_0)_{\mathfrak{X}}\bigr)^{\mathrm{an}}\rightarrow
	\bigl(\D(\Lambda)_{\mathfrak{X}}\bigr)^{\mathrm{an}}$.
  \item Through the isomorphism
	$\bigl(\D(\Lambda_0)_{\mathfrak{X}}\bigr)^{\mathrm{an}}\rightarrow
	\bigl(\D(\Lambda)_{\mathfrak{X}}\bigr)^{\mathrm{an}}$ (or equivalently,
	through the isomorphism
	$\D\bigl(\Lambda\times_{\mathfrak{X}}
	\mathfrak{X}/p\bigr)_{\mathfrak{X}}^{\mathrm{an}}
	\rightarrow
	\bigl(\D(\Lambda)_{\mathfrak{X}}\bigr)^{\mathrm{an}}$), the 
	analytic vector bundle
	$\bigl(\D(\Lambda)_{\mathfrak{X}}\bigr)^{\mathrm{an}}$ is naturally
	endowed with an integrable connection, which deserves to be
	called the {\it Gauss-Manin} connection attached to $\Lambda$.
 \end{enumerate}

 Indeed, when $\Lambda=\udl{A}[p^\infty]$ for an abelian scheme
 $\udl{A}$ over $\mathfrak{X}$,
 $(\D(\Lambda)_{\mathfrak{X}})^{\mathrm{an}}$ may then be identified
 with the De Rham cohomology bundle
 $\mathrm{H}^1_{\mathrm{DR}}(\udl{A}^{\mathrm{an}}/
 {\mathfrak{X}}^{\mathrm{an}})$
 together with its Gauss-Manin connection.  Here, one can use the
 comparison theorem between crystalline and De Rham cohomology, or the
 theory of universal vectorial extensions and Grothendieck's
 $\natural$-structures, \cf \cite{mazur74:_univer}, \cite{tddc2}.

 In the first alternative, one uses the
 divided powers on $(p)$ and construct a canonical isomorphism of vector
 bundles with connection 
 $\D(\Lambda\times_{\mathfrak{X}}\mathfrak{X}/p)_{\mathfrak{X}}\cong
 \mathrm{H}^1_{\mathrm{DR}}(\udl{A}/{\mathfrak{X}})$.

 In the second alternative, one uses Grothendieck's construction of a
 canonical connection on the universal vectorial extension $\E(\Lambda)$
 of $\Lambda^\vee$ --- \cf \ref{para:3.3.2}; the induced connection on the Lie
 algebra $\D(\Lambda)_{\mathfrak{X}}$ turns out to be the Gauss-Manin
 connection on $ \mathrm{H}^1_{\mathrm{DR}}(\udl{A}/{\mathfrak{X}})$.

 The Gauss-Manin connection is especially significant when
 ${\mathfrak{X}}^{\mathrm{an}}$ is smooth (even though $\mathfrak{X}$
 itself may not be formally smooth).

\newpage
\section{Moduli problems for $p$-divisible groups.}\label{sec:4}

\begin{abst}
 Review of moduli spaces of $p$-divisible groups quasi-isogenous to a
 fixed one modulo $p$, with examples.
\end{abst}
\end{para}

\subsection{Moduli problem for $p$-divisible groups quasi-isogenous to a
 fixed one modulo $p$.}\label{sub:4.1}

\begin{para}\label{para:4.1.1} Here, we work over the discrete valuation
 ring $\mathfrak{v}=W= {\widehat{\Z_p^{\mathrm{ur}}}}$ (the completion of the
maximal unramified extension of $\Z_p$), with uniformizing parameter
$\varpi=p$. Recall that $\Nil_\mathfrak{v}$ denote the category of
locally noetherian $\mathfrak{v}$-schemes $S$ on which $p$ is locally
nilpotent.

Let $\ovl{\Lambda}$ be a fixed $p$-divisible group over
$k={\ovl{\F}}_p$. The following theorem is due to M.  Rapoport and
T. Zink \cite[2.16, 2.32]{psfpg}.
\end{para}

\begin{thm}\label{thm:4.1.2} The functor
\begin{align*}
 \Nil_\mathfrak{v}&\rightarrow \Sets\\
 S&\mapsto \Bigl\{(\Lambda,\rho)\Bigm|
 \begin{array}{l}
  \Lambda \text{: $p$-divisible group on $S$},\\
  \rho\in \qisog(\Lambda\times_S S_{\mathrm{red}},
 \ovl{\Lambda}\times_{{\ovl{\F}}_p}S_{\mathrm{red}})
 \end{array}
 \Bigr\}\Big/\cong
\end{align*}
is representable by a formal scheme $\mathfrak{M}$
over $\mathfrak{v}$; $\mathfrak{M}_{\mathrm{red}}$ is locally of
finite type over ${\ovl{\F}}_p$ and its irreducible components are
projective ${\ovl{\F}}_p$-varieties.
\end{thm}

Of course, the rigidity of $p$-divisible groups up to isogeny (\ref{thm:2.2.3}) is
crucial here. The main difficulty in proving \ref{thm:4.1.2} is to control the
powers of $p$ which are needed in the process of lifting isogenies.

The separated formal scheme $\mathfrak{M}$ is not formally smooth in
general, but has been conjectured to be flat over 
$\mathfrak{v}$.  The group $J$ of self-quasi-isogenies of $\ovl{\Lambda}$
 acts (on the right) on $\mathfrak{M}$ by  $(\Lambda,\rho).j= 
(\Lambda,\rho\circ j$). This group $J$ is actually the group of
 ${\Q}_p$-points of an algebraic group over ${\Q}_p$. 

 We denote by $\mathcal{M}$ the Berkovich analytic space over
 $\C_p$ attached to $\mathfrak{M} : 
 \mathcal{M}={\mathfrak{M}}^{\mathrm{an}}\hat\otimes_K\;\C_p$
 (\cf \ref{sub:3.5}). This is a paracompact (strictly) analytic space. We shall
 see later that it is smooth, and a $p$-adic manifold in the sense of
 I.\ref{para:p-adic-manifold}.

\subsection{Examples.}\label{sub:4.2}
 We take $\ovl{\Lambda} ={\ovl{A}}[p^\infty]$, the $p$-primary
 torsion of an {\it elliptic curve ${\ovl{A}}$ over ${\ovl{\F}}_p$} ($h=2$).

\begin{para}\label{para:4.2.1}
 {\it The ordinary case}, \ie ${\ovl{A}}[p]({\ovl{\F}}_p)\neq 0$. 

 Then $\ovl{\Lambda}\cong\hat {\G}_m \oplus {\Q}_p/{\Z}_p$, and
 $J=({\G_m}_{{\Q}_p})^2$. For any $(\Lambda,\rho)$ in $\mathfrak{M}(S)$,
 the separable rank of $\Ker[p]$ is constant (due to the quasi-isogeny
 $\rho$). Therefore $\Lambda$ is an extension 
 \begin{equation*}
  \hat{\Lambda}\rightarrow \Lambda\rightarrow \Lambda^{\mathrm{et}}.
 \end{equation*}
 The quasi-isogeny $\rho$ respects this extension,
 hence splits into two parts $({\hat\rho},\rho^{\mathrm{et}})$. We have
 ${\hat\Lambda}\cong \hat {\G}_m ,\;\Lambda^{\mathrm{et}}\cong
 {\Q}_p/{\Z}_p$. Up to isomorphism, the pair
 $({\hat\rho},\rho^{\mathrm{et}})$ thus amounts to an element of
 $({\Q}_p^\times/{\Z}_p^\times)^2\cong {\Z}^2$. Therefore
 \begin{equation*}
  \mathfrak{M}=\coprod_{\Z^2}\;\mathfrak{M}^0
 \end{equation*}
 is a disjoint sum of copies of a formal (group) scheme
 $\mathfrak{M}^0$ which parametrizes extensions of ${\Q}_p/{\Z}_p$ by $
 \hat {\G}_m$. It is known that
 \begin{equation*}
  \mathfrak{M}^0 =\hat {\G}_m
 \end{equation*} 
 This identification is given by the
 following recipe. Recall that if $S=\Spec(R),\;R$ local artinian with
 residue field ${\ovl{\F}}_p$ and radical $\mathfrak{m}$
 ($\mathfrak{m}^{n+1}=0$), then $\hat {\G}_m(S)=1+\mathfrak{m}$; in
 particular $\hat {\G}_m(S)$ is killed by $[p^n]$. Let us choose any lift
 $q_n$ of $1/p^n\in {\Q}_p/{\Z}_p$ in $\Lambda$. Then $[p^n]q_n$ is a
 well-defined element $q$ in $ {\hat\Lambda}(S)\cong 1+\mathfrak{m}$,
 which is unchanged if $n$ is replaced by a bigger integer (but another
 choice of the isomorphism ${\hat\Lambda}\cong \hat {\G}_m$ would change
 $q$ into $q^a,\; a\in {\Z}_p^\times$).

 \noindent As a formal $W$-scheme, $\mathfrak{M}^0=
 \Spf(\mathfrak{v}[[q-1]])$, with the $(p,q-1)$-adic topology. 
 The associated analytic space $\mathcal{M}^0$ over $\C_p$ is the open
 unit disc $| q-1| <1$.
 
 The parameter $q$ is thus the local modulus for deformations of
 ${\ovl{A}}[p^\infty]$.  
 On the other hand, if
 $p\neq 2$, we have an algebraic local modulus for deformations of
 ${\ovl{A}}$ (say, with level two structure): the Legendre parameter
 $z=\lambda$ around the Legendre parameter $\zeta_{\mathrm{can}}$ of the
 canonical lifting $A_{\mathrm{can}}$ of ${\ovl{A}}$ (\cf 
 I.\ref{para-dwork-para}).
 The $p$-divisible group $A_{\mathrm{can}}[p^\infty]$ splits: its parameter is
 $q=1$. Therefore, in this special case, the Serre-Tate theorem \ref{thm:2.2.4}
 asserts that $q\in
 1+(z-\zeta_{\mathrm{can}}){\widehat{\Z}^{\mathrm{ur}}}[[z-\zeta_{\mathrm{can}}]]$. It has been
 proved by W. Messing and N. Katz that $q$ is the Dwork-Serre-Tate
 parameter discussed in I.\ref{sub-serre-tate-para}, 
 (\cite{q}, \cite{slm}) --- which explains the terminology.
\end{para}

\begin{para}\label{para:4.2.2} {\it The supersingular case}, \ie
 ${\ovl{A}}[p]({\ovl{\F}}_p)= 0$. 
 
 In this case ${\ovl{\Lambda}}
 =\hat{\ovl{\Lambda}}$ is a so-called Lubin-Tate formal group,
 $\;\End\;{\ovl{\Lambda}}=\mathcal{B}_p$, ``the'' maximal order in 
 the (non-split) quaternion algebra $B_p$ over ${\Q}_p$, and $J=B_p^\times$.

 Here again, there is a discrete invariant: the ``height'' of the
 quasi-isogeny $\rho$, which amounts to the 
 $p$-adic valuation of the norm of its image in $J=B_p^\times$.
 Therefore 
 \begin{equation*}
  \mathfrak{M}=\coprod_{\Z}\;\mathfrak{M}^0
 \end{equation*}
 is a disjoint sum of copies of a formal scheme $\mathfrak{M}^0$
 which parametrizes deformations of the Lubin-Tate group
 $\hat{\ovl{\Lambda}}$. According to Lubin-Tate, $\mathfrak{M}^0\cong
 \Spf(\mathfrak{v}[[t]])$, with the $(p,t)$-adic topology. The
 associated analytic space $\mathcal{M}^0$ over $\C_p$ is the open unit
 disc $|t|<1$.

 Assume for simplicity that ${\ovl{A}}$ is defined over ${\F}_p$. Then
 one can describe the universal formal deformation $\Lambda$ over
 $\Spf{\Z_p[[t]]}$ through its logarithm $\sum b_n(t)X^{p^{2n}}$, which
 is given recursively by
 \begin{equation*}
  b_0(t)=1,\;b_n(t)=\frac{1}{p}\sum_{0\leq m < n}b_m(t)t_{n-m}^{p^{2m}}
 \end{equation*}

The so-called canonical lifting, corresponding to $t=0$, has formal complex
multiplication by $\Z_{p^2}$, \cf \cite{gross94:_equiv_lubin_tate}, \cite{trapmlsasht}.
\end{para}

\subsection{Decorated variants.}\label{sub:4.3} 

\begin{para}\label{para:4.3.1} Just as in the complex case
 \ref{para:1.1.4}, it is useful to consider variants or the moduli
 problem \ref{sub:4.1} for ``decorated'' $p$-divisible groups, with prescribed
 endomorphisms and/or polarization.

Instead of
the $\Q$-algebra $B$ of \ref{para:1.1.4}, we shall consider a finite-dimensional semi-simple $\Q_p$-algebra $B_p$, since
the category of $p$-divisible groups up to isogeny is
$\Q_p$-linear. Let $V_p$ be a $B_p$-module
of finite type, and let $G_{\Q_p}$ be the group of 
$B_p$-linear endomorphisms of $V_p$.
\par\noindent  Let $\mathcal{B}_p$ be a maximal order in $B_p$, and let
$L_p$ be a lattice in $V_p$, stable under $\mathcal{B}_p$.

\end{para}

\begin{para}\label{para:4.3.2} We assume that $\mathcal{B}_p$ acts on our
 $p$-divisible group $\ovl{\Lambda}$. One then redefines $J$ to 
be the group of self-quasi-isogenies of $\ovl{\Lambda}$ which respect
the $\mathcal{B}_p$-action.

For technical reasons, one picks up a lifting $\tilde{\ovl{\Lambda}}$ of
$\ovl{\Lambda}$ with $\mathcal{B}_p$-action to some finite extension of
$W$, and looks at the $\mathcal{B}_p$-module
$F^1_0:=\omega_{\tilde{\ovl{\Lambda}}}$. We denote by $\mathfrak{v}$ the
finite extension $W[\tr(b| F^1_0)]$ of $W$.

One strengthens the moduli problem \ref{sub:4.1} by imposing to our pairs
 $(\Lambda,\rho)$ that $\Lambda$ is endowed with a
 $\mathcal{B}_p$-action $\iota:\mathcal{B}_p\rightarrow \End \Lambda$,
 $\rho$ respects the $\mathcal{B}_p$-action, and moreover a ``Shimura
 type condition'' $\det(\iota(b)| \omega(\Lambda))=\det(b|
 F^1_0)$.

 It follows easily from \ref{thm:4.1.2} \cite[3.25]{psfpg} that this
 moduli problem is {\it representable by a formal scheme, still denoted
 by $\mathfrak{M}$, over $\mathfrak{v}$. It is acted on by
 $J$. $\mathfrak{M}_{\mathrm{red}}$ is locally of finite type over ${\ovl{\F}}_p$
 and its irreducible components are projective
 ${\ovl{\F}}_p$-varieties. }
\end{para}

\begin{para}\label{para:4.3.3} Let us say a word about the polarized variant
 (here, it is safer to assume $p\neq 2$). One assumes 
 that $\mathcal{B}_p$ is endowed with an involution $\ast$, and that
 $L_p$ is autodual for an alternate $\Z_p$-bilinear form such that
 $\an{bv,w}=\an{v,b^\ast w}$. Of course, here, $G_{\Q_p}$ denotes the algebraic
 $\Q_p$-group of $B_p$-linear symplectic similitudes of $V_p$. One
 assumes that $\ovl{\Lambda}$ is endowed with a $\ast$-polarization
 \footnote{a $\ast$-polarization is a symmetric
 $\mathcal{B}_p$-linear quasi-isogeny between the $p$-divisible and its
 Serre dual, endowed with the transposed action of $\mathcal{B}_p$
 twisted by $\ast$}, and modify $J$ to respect the polarization up to a
 multiple in $\Q_p$.
 
 \bigskip
 The corresponding moduli problem is again representable by a formal
 scheme with the same properties. This formal scheme has been conjectured to be flat by Rapoport-Zink; although this has been verified in many cases (\cf \eg \cite{gortz}), there are some counter-examples  (\cf  \cite{pappas}).
\end{para}

\begin{exa}[fake elliptic curves at a critical prime $p$]\label{exa:4.3.4}
 Let us fix a indefinite quaternion algebra $B$
 over $\Q$ and a maximal order $\mathcal{B}$ in $B$ (after Eichler, they
 are all conjugate). A prime $p$ is called {\it critical} if it divides
 the discriminant of $B$, \ie if $B_p=B\otimes_\Q \Q_p$ is the
 (non-split) quaternion algebra over $\Q_p$. Note that $B_p$ then
 contains copies of the unramified quadratic extension on $\Q_p$; we
 denote by $\Q_{p^2}$ one of them. The maximal order (unique up to
 conjugation) can be written $\mathcal{B}_p=\mathcal{B}\otimes
 \Z_p=\Z_{p^2}[\Pi],\;\Pi^2=p,\;\Pi.a=a^\sigma\Pi$ for $a\in
 \Z_{p^2}$. We fix a critical prime $p$.

 We consider abelian surfaces with quaternionic multiplication by
 $\mathcal{B}$. They are sometimes called ``{\it fake elliptic curves}''
 \index{fake elliptic curves@fake elliptic curves}
 because they share many features with elliptic curves, 
 \cf \cite[III, 1]{updcds}.

 We fix a fake elliptic curve $\ovl{A}$ over $\ovl{\F}_p$, and set
 $\Lambda=\ovl{A}[p^\infty]$. This is a special formal
 $\mathcal{B}_p$-modules of height $4$ in the following sense.

 Let $S$ be in $\Nil_\mathfrak{v}$ (or else an inductive limit of such
 schemes --- as a formal scheme), and let $\Lambda$ be a $p$-divisible
 group over $S$. After Drinfeld, one says that $\Lambda$ is a {\it
 special formal $\mathcal{B}_p$-modules of height $4$} if

\begin{enumerate}
 \item $\Lambda$ is infinitesimal ($\Lambda=\hat\Lambda$) of height
       $h=4$,
 \item $\mathcal{B}_p$ acts on $\Lambda$,
 \item via this action, $\omega_\Lambda$ is a locally free
       $(\Z_{p^2}\otimes_{\Z_p}\mathcal{O}_S)$-module of rank $1$.
\end{enumerate}

\noindent We take $V_p=B_p$ where $B_p$ acts by left
multiplication. In this situation, $G_{\Q_p}\cong B_p^\times $. When
$S=\Spec(\ovl{\F}_p)$, there is a unique class of
$\mathcal{B}_p$-isogeny of special formal $\mathcal{B}_p$-modules of
height $4$. The group $J$ of self-quasi-isogenies of a special formal
$\mathcal{B}_p$-module of height $4$ is $J=GL_2(\Q_p)$. This can be
easily read off the Dieudonn\'e module $\otimes\Q$.

In this example, there is again a discrete invariant: the ``height'' of 
the quasi-isogeny $\rho$, which amounts to the 
$p$-adic valuation of the determinant of its image in $J$.
Therefore 
 \begin{equation*}
 \mathfrak{M}=\coprod_{\Z}\;\mathfrak{M}^0
\end{equation*}
 is a disjoint sum of copies of a formal scheme $\mathfrak{M}^0$
 which parametrizes deformations of the special formal module
 $\hat{\ovl{\Lambda}}$. This formal scheme has been described by
 Drinfeld, \cf \ref{para:6.3.4}.
\end{exa}

\begin{para}\label{para:4.3.5} This example and examples \ref{sub:4.2}
 were the sources of the general theory. 

\noindent Let us at once suggest a vague analogy with the
archimedean situation --- especially \ref{para:1.1.4}, \ref{para:1.2.4} ---
 which will become more 
and more precise in the next sections: one can see the $p$-divisible
group $\ovl{\Lambda}$ in characteristic $p$ (\resp $J$) as an analogue
of the infinitesimal oriented real Lie group $\hat U$ (\resp $G(\R)^0$).
    
\newpage
\section{$p$-adic period domains.}\label{sec:5}

\begin{abst}
 Review of $p$-adic symmetric domains and period domains via stability theory.
\end{abst}
\end{para}

\subsection{$p$-adic flag spaces.}\label{sub:5.1}

\begin{para}\label{para:5.1.1} We begin with the $p$-adic analogue of
 the flag space $\mathcal{D}^\vee$ of \ref{para:1.2.3},
 \ref{para:1.2.4}. In \ref{para:1.2.4}, this flag space parametrized the
 lagrangian subspaces $F^1$ in the $\mathrm{H}_{\mathrm{DR}}^1$ of abelian
 varieties of a certain type, marked by a fixed isomorphism
 $\mathrm{H}_{\mathrm{DR}}^1\cong V^\C$.

In the $p$-adic case, $\mathcal{D}^\vee$ will parametrize the Hodge
filtration $F^1$ of Dieudonn\'e modules of $p$-divisible groups of a
certain type (\ref{para:3.3.1}).
\end{para}

\begin{para}\label{para:5.1.2} More precisely, let us consider a moduli
 problem for $p$-divisible groups as in \ref{sub:4.1}, possibly
 decorated as in \ref{sub:4.3} (from which we follow the notation). The
 decoration reflects on the 
 Dieudonn\'e module $\D(\ovl{\Lambda})_\Q=\D(\ovl{\Lambda})\otimes \Q$:
 one has an isomorphism of $B_p$-module (with symplectic structure, in
 the polarized case):
 \begin{equation*}
  \D(\ovl{\Lambda})_\Q\cong
   V_p\otimes_{\Q_p}\widehat{\Q_p^{\mathrm{ur}}}.
 \end{equation*}
 We fix such an isomorphism, and consider that
 $G_{\Q_p}\otimes_{\Q_p}\widehat{\Q_p^{\mathrm{ur}}}$ acts on
 $\D(\ovl{\Lambda})_\Q$. But $\D(\ovl{\Lambda})_\Q$ has an extra
 structure: the $\sigma$-linear Frobenius automorphism. The group $J$ is
 nothing but the {\it subgroup} of $G_{\Q_p}(\widehat{\Q_p^{\mathrm{ur}}})$ of
 elements commuting with Frobenius.
\end{para}

\begin{para}\label{para:5.1.3} Let $\mathfrak{v}$ be the finite extension of
 $W={\Z_p^{\mathrm{ur}}}$ considered in \ref{para:4.3.2} and let $K$ be its 
 fraction field. We denote here by $F^1_0\subset V_p^{\C_p}=V_p\otimes
 \C_p$ the $\C_p$-span of the space
 $\omega_{\tilde{\ovl{\Lambda}}}\subset V_p\otimes_{\Q_p}K$ introduced in
 \ref{para:4.3.2}. We set $G_{\C_p}=G_{\Q_p}\otimes_{\Q_p}\C_p$.

 \noindent The flag variety $\mathcal{D}^\vee$ of all
 $G_{\C_p}$-conjugates of $F^1_0$ in $V_p^{\C_p}=\D(\ovl{\Lambda})
 \otimes_W\C_p$ is a homogeneous space 
 \begin{equation*}
  \mathcal{D}^\vee = P\backslash G_{\C_p}
 \end{equation*}
 for a suitable parabolic subgroup $P$. There is an
 ample line bundle $\mathcal{L}$ on $\mathcal{D}^\vee$
 which is homogeneous under the derived group $G_{\C_p}^{\mathrm{der}}$;
 \index{000GCpder@$G_{\C_p}^{\mathrm{der}}$}
 we denote by $\mathcal{D}^\vee \subset \P(W)$ an 
 associated $G_{\C_p}^{\mathrm{der}}$-equivariant projective embedding (if
 $G_{\C_p}=GL_{2g},\;\mathcal{D}^{\vee} = \Grass(g,2g)$, one 
 can just take the $SL_{2g}$-equivariant Pl\"ucker embedding).

 Note that, there is an obvious right action of $J$ on
 $\mathcal{D}^{\vee}$, which factors through the adjoint group
 $J^{\mathrm{ad}}:=\Im\bigl(J\rightarrow G_{\C_p}^{\mathrm{ad}}\bigr)$. 
\end{para}

\begin{para}\label{para:5.1.4} In our three basic examples (\ref{para:4.2.1}, \ref{para:4.2.2}, \ref{exa:4.3.4}), we
have $G_{\C_p}=GL_2$, hence $\mathcal{D}^\vee=\P^1_{\C_p}$.    
\end{para}

\subsection{Symmetric domains via stability theory.}\label{sub:5.2}

\begin{para}\label{para:5.2.1} In the $p$-adic case, one cannot define symmetric spaces by a positivity condition as in \ref{para:1.2.3}, \ref{para:1.2.4}.
One uses instead the Hilbert-Mumford notion of semi-stability, following an idea of M. van der Put and H. Voskuil
\cite{ssatsagoalf}.

\end{para}

\begin{para}\label{para:5.2.2}
 We set $J'=J\cap G_{\C_p}^{\mathrm{der}}$. This is
 again the group of $\Q_p$-points of an algebraic group over $\Q_p$. Let
 $T\cong ({\G_m}_{\Q_p})^j$ be a maximal $\Q_p$-split torus in $J'$. The
 {\it semi-stable locus} of $\mathcal{D}^\vee$ with respect to $T$-action
 is the open subset of points $x\in \mathcal{D}^\vee$ such that there is
 a $T$-invariant function of $\P(W)$ which does not vanish on $x$.

 One can check whether $x$ is semi-stable with respect to $T$ using the
 Hilbert criterion: let $\lambda$ be any cocharacter (= one-parameter
 subgroup) of $T$. By properness of $\mathcal{D}^\vee$, the map $z\in
 \G_m\mapsto x.\lambda(z)\in \mathcal{D}^\vee$ extends to $\A^1$, \ie the
 limit $\displaystyle x_0=\lim_{z\rightarrow 0} x.\lambda(z)$ exists. It
 is clearly fixed by $\G_m$, thus the fiber $\mathcal{L}_{x_0 }$
 corresponds to a character $z\mapsto z^{-\mu(x,\lambda)}$ of
 $\G_m$. Then $x$ is semi-stable if and only if the Mumford invariant
 $\mu(x,\lambda)$ is $\geq 0$ for every $\lambda$.
\end{para}

\begin{para}\label{para:5.2.3} One defines the {\it symmetric space
 $\mathcal{D}$ to be the intersection of the semi-stable loci 
 of $\mathcal{D}^\vee$ with respect to all maximal $\Q_p$-split tori $T$
 in $J'$.} This is a $p$-adic manifold in the sense of
 I.\ref{para:strict-analytic}.  It is
 also called the {\it $p$-adic period domain} (associated with the moduli
 problem).

 Since $J$ acts on the set of tori $T$ by conjugation, it is clear that
 $\mathcal{D}$ is stable under the $J$-action on $\mathcal{D}^{\vee}$.
\end{para}

\subsection{Examples.}\label{sub:5.3}

\begin{para}\label{para:5.3.1} In our three basic examples, we have
$G_{\C_p}=GL_2,\;\mathcal{D}^\vee=\P^1$.

In example \ref{para:4.2.1}, $J'\subset J=(\G_m)^2$ is the one-dimensional torus
$T$ with elements $(z^{-1},z)$. Its action on a point $x$ is by
$z^2$-scaling; the limit point $x_0$ is always $0$, except if
$x=\infty$. The only non-semi-stable point is $\infty$, hence
$\mathcal{D}=\A^1$.

In example \ref{para:4.2.2}, $J'=SL_1(B_p)\subset J=B_p^\times$. There is no
 non-trivial $\Q_p$-split torus in $J'$, hence
 $\mathcal{D}=\mathcal{D}^{\vee}=\P^1$.

\end{para}

\begin{para}\label{para:5.3.2} We now turn to example \ref{exa:4.3.4}. In this
 example, $J=GL_2(\Q_p)$ and $J'=SL_2(\Q_p)$. All maximal 
 $\Q_p$-split tori are conjugate $g^{-1}Tg$ of the above $T$. Therefore,
 $\mathcal{D}=\P^1\setminus \bigcup \infty.g$ 
 is the {\it Drinfeld space} over $\C_p$ 
 \begin{equation*}
  \mathcal{D}=\Omega_{\C_p} :=\P^1_{\C_p}\setminus \P^1(\Q_p)
 \end{equation*} 

 To be consistent, we consider the right action of
 $PGL_2$-action on $\P^1$, which is deduced from the customary left 
 action by the rule $x.g=g^{-1}.x$. 

 The Drinfeld space is actually defined over $\Q_p$, \ie is a
 $\Q_p$-analytic manifold $\Omega$. 
 It is closely related to the {\it Bruhat-Tits tree}
 $\mathcal{T}$ of $J^{\mathrm{ad}}=PGL_2(\Q_p)$. The vertices of this
 tree correspond to closed discs in 
 $\Q_p$; two such vertices are connected by an edge if for the
 corresponding discs $\matheur{D}',\matheur{D}''$, one has
 $\matheur{D}'\subset \matheur{D}''$ and the  
 radius of $\matheur{D}''$ is $p$ times the radius of $\matheur{D}'$
 (each vertex has $p+1$ 
 neighbors). The group $PGL_2(\Q_p)$ acts on $\mathcal{T}$ (on 
 the right, say).
 \vspace*{5mm}
 \begin{figure}[h]
  \begin{picture}(120,120)(0,0)
   \put(0,0){\includegraphics[scale=0.5,keepaspectratio,clip]{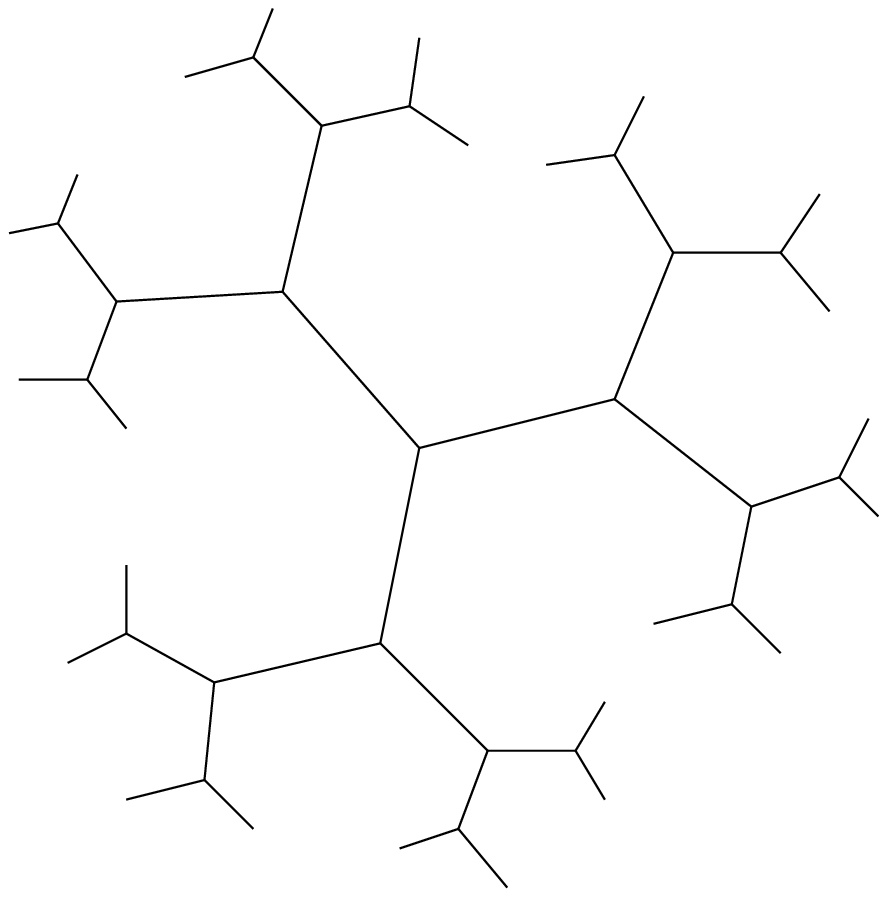}}
  \end{picture}\vspace{-1ex}
  \caption{$\mathcal{T}$ for $p=2$}
  \label{fig-bruhat-tits}
 \end{figure}

 For any disc $\matheur{D}$ centered in $\Q_p$, let $\eta_\matheur{D}$
 denote its Berkovich 
 ``generic point'' (which is a point in Berkovich's affine line $\A^1$
 over $\Q_p$). Then the geometric realization of 
 $\mathcal{T}$ may be identified with the closed subset 
 \begin{equation*}
  |\mathcal{T}|= \{\eta_{\matheur{D}(a,r^+)} \mid a\in \Q_p,\;0<r<\infty\}
  \,\subset \,\Omega
 \end{equation*}
 (a $PGL_2(\Q_p)$-equivariant embedding), and the set of vertices with
 the subset 
\begin{equation*}
 \mathrm{ver}\, \mathcal{T}=
  \{\eta_{\matheur{D}(a,r^+)} \mid a\in \Q_p,\;r\in p^\Z\}\, \subset
  \,|\mathcal{T}| 
\end{equation*}
 According to Berkovich, $|\mathcal{T}|$ is a retract of $\Omega$
 \cite{staagonf}; in particular, $\Omega$ is simply connected.

 On the other hand, $\Omega$ is the generic fiber of a $p$-adic formal
 scheme $\hat\Omega$ whose reduction modulo $p$ is an 
 infinite tree of projective lines, with dual graph $\mathcal{T}$ (\cf
 \cite[I]{updcds}).

 The Drinfeld space $\Omega$ is a $p$-adic analogue of the Poincar\'e
 double half-plane $\P^1_{\C}\setminus \P^1(\R)$. Note that the latter
 space could also be defined using stability, in the same way (with
 $J'=SL_2(\R)$).  Another similarity: according to \cite{staagonf},
 there is a $PGL_2(\Q_p)$-invariant metric on $\Omega$ (which 
 extends the standard metric on $|\mathcal{T}|$). 
\end{para}

\newpage
\section{The $p$-adic period mapping and the Gauss-Manin connection.}
\label{sec:6}

\begin{abst}
 Rapoport-Zink's period mapping as classifying map for the Hodge
 filtration in the Dieudonn\'e module. Basic properties: equivariance,
 \'etaleness. Its relation to the Gauss-Manin connection. Examples and
 formulas.
\end{abst}

\subsection{Construction.}\label{sub:6.1} 

\begin{para}\label{para:6.1.1} The construction of the $p$-adic period
 mapping is parallel to the complex case: one 
 attaches to a point $s\in \mathcal{M}$ representing a pair $(\Lambda,
 \rho)$ the point $\mathcal{P}(s)$ of $\mathcal{D}^\vee$ 
which parametrizes the notch $F^1$ of the Hodge filtration of the
 Dieudonn\'e module of $\Lambda$. The idea of the construction is
 already present in Grothendieck's talk
 \cite{grothendieck70:_group_barsot_tate}.
\end{para}

\begin{para}\label{para:6.1.2} More precisely, notation being as in
 \ref{para:4.1.1}--\ref{para:4.3.1}, let us consider the universal object 
 $(\udl{\Lambda},\udl{\rho})$ on the moduli formal scheme
 $\mathfrak{M}$, and the Dieudonn\'e crystal $\D(\udl{\Lambda})$. Its
 evaluation $\D(\udl{\Lambda})_{\mathfrak{M}}$ on $\mathfrak{M}$ is a
 vector bundle on $\mathfrak{M}$.
\end{para}

\begin{lem}\label{lem:6.1.3} The quasi-isogeny $\udl{\rho}$ induces a
 trivialization of the associated analytic vector bundle on
$\mathcal{M}={\mathfrak{M}}^{\mathrm{an}}_{\C_p} $:
\begin{equation*}
 (\D(\udl{\Lambda})_{\mathfrak{M}})^{\mathrm{an}}\cong
 V_p^{\C_p}\otimes_{\C_p}\mathcal{O}_\mathcal{M}
\end{equation*}

\end{lem}

\begin{proof}
 $\rho$ is a quasi-isogeny between $\udl{\Lambda}
 \times_{\mathfrak{M}}{\mathfrak{M}}_{\mathrm{red}}$ and
 $\ovl{\Lambda}\times_{\ovl{\F}_p}{\mathfrak{M}}_{\mathrm{red}}$. According to
 \ref{para:3.6.7} (i), by the rigidity of convergent isocrystals, it induces a
 canonical isomorphism $(\D(\udl{\Lambda})_{\mathfrak{M}})^{\mathrm{an}}\cong
 \bigl(\D(\ovl{\Lambda}\times_{\ovl{\F}_p}
 {\mathfrak{M}}_{\mathrm{red}})_{\mathfrak{M}}\bigr)^{\mathrm{an}}.$
 On the other hand, since the formation of $\D(\Lambda)$ commutes with
 base change (\ref{para:3.2.2}), we have
 \begin{equation*}
 \bigl(\D(\ovl{\Lambda}\times_{\ovl{\F}_p}
 {\mathfrak{M}}_{\mathrm{red}})_{\mathfrak{M}}\bigr)^{\mathrm{an}} =
 \D(\ovl{\Lambda})\otimes_W 
 \mathcal{O}_\mathcal{M}=V_p\otimes_{\Q_p}\mathcal{O}_\mathcal{M}.
 \end{equation*}
\end{proof}

\begin{para}\label{para:6.1.4}
 On the other hand, we have a locally direct summand
 $F^1=(\omega_{\udl{\Lambda}})^{\mathrm{an}}$ in
 $(\D(\udl{\Lambda})_{\mathfrak{M}})^{\mathrm{an}}\cong
 V_p\otimes_{\Q_p}\mathcal{O}_\mathcal{M}$ (\ref{para:3.3.1}). This
 determines an 
 $\mathcal{M}$-point of the flag space $\mathcal{D}^\vee$, \ie an
 analytic mapping 
 \begin{equation*}
  \mathcal{P}:\mathcal{M}\;\longrightarrow \;\mathcal{D}^\vee
 \end{equation*}

 This is called the {\it period mapping} for the moduli problem
 (decorated or not). It is clear from the construction that $\mathcal{P}$ is
 $J$-{\it equivariant}.

 If $\mathcal{M}^0$ is any connected component of $\mathcal{M}$, it is
 stable under $J'\subset J$, hence the induced period mapping
 $\mathcal{P}:\;\mathcal{M}^0\;\longrightarrow \;\mathcal{D}^\vee\;$ is
 $J'$-equivariant (the $J'$-action factors through ${J'}^{\mathrm{ad}}$).
\end{para}

\subsection{Properties.}\label{sub:6.2}  

\begin{thm}\label{thm:6.2.1} The period mapping $\mathcal{P}$ is etale. Its image lies in the period domain $\mathcal{D}$.

\end{thm}
The first assertion reduces, via an infinitesimal criterion of
etaleness, to the ``essential surjectivity part'' of the
Grothendieck-Messing theorem (\ref{thm:3.3.3}), in the case of an ideal
of square zero, \cf \cite[5.17]{psfpg}.  It corresponds to what is called
the ``local Torelli property'' in the complex case.

The second assertion is a
consequence of a theorem of B. Totaro on the equivalence of semi-stability and ``weak admissibility'' \cite{tpipht}.
The second assertion may be viewed as a
$p$-adic analogue of the ``Riemann relations'' (\ref{para:1.1.3}).

\begin{cor}\label{cor:6.2.2} $\mathcal{M}$ is smooth, and any component
 $\mathcal{M}^0$ is a connected $p$-adic manifold in the sense of
 I.\ref{para:p-adic-manifold}.
\end{cor}
\begin{para}\label{para:6.2.3}
 It is very likely that $\mathcal{P}$, viewed as a mapping
 $\mathcal{M}\rightarrow \mathcal{D}$, is surjective.  This is closely
 related to Fontaine's conjecture on the realizability of filtered
 Dieudonn\'e modules by $p$-divisible groups. Although recent work by
 C. Breuil \cite{breuil99:_repres}, \cite{breuil00:_group} settles the
 Fontaine conjecture at least for 
 finite residue fields of characteristic $\neq 2$ (and in other cases),
 it is not clear to the author to which extent the surjectivity
 conjecture is now established.

 On the other hand, the fibers of the period mapping are relatively
 well-understood, \cf \cite[5.37]{psfpg}. 
\end{para}

\begin{para}\label{para:6.2.4}
 The argument in the proof of \ref{lem:6.1.3} shows a little more. Recall from
 \ref{para:3.6.7} (ii) that the analytic vector bundle
 $(\D(\udl{\Lambda})_{\mathfrak{M}})^{\mathrm{an}}$ on $\mathcal{M}$ comes
 together with an integrable connection $\nabla_{\mathrm{GM}}$ (the Gauss-Manin
 connection).
\end{para}

 \begin{lem}\label{lem:6.2.5} The trivialization 
  \begin{equation*}
   (\D(\udl{\Lambda})_{\mathfrak{M}})^{\mathrm{an}}\cong
   V_p^{\C_p}\otimes_{\C_p}\mathcal{O}_\mathcal{M}
  \end{equation*}
  of \ref{lem:6.1.3} induces a trivialization of the Gauss-Manin connection:  
  \begin{equation*}
   \bigl(
    (\D(\udl{\Lambda})_{\mathfrak{M}})^{\mathrm{an}}
    \bigr)^{\nabla_{\mathrm{GM}}} 
    \cong V_p^{\C_p}.
  \end{equation*}

  The point is that the convergent isocrystal
  $\D(\Lambda\times_{\mathfrak{M}}{\mathfrak{M}}_{\mathrm{red}})$ is constant, 
  since it comes by base change, from $\D(\ovl{\Lambda})_\Q$.
 \end{lem}
 This useful lemma will allow us to give explicit expression of the
 period mapping in terms of quotients of solutions of differential
 equations, as in the complex case.

\subsection{Examples and formulas.}\label{sub:6.3} 

\begin{para}\label{para:6.3.1} We start with example \ref{para:4.1.1}
 ($p$-divisible groups of elliptic curves with ordinary reduction).
 The extension
 \begin{equation*}
  1\rightarrow {\hat\G_m}\rightarrow \Lambda\rightarrow \Q_p/\Z_p\rightarrow 1
 \end{equation*}
 parametrized by $q\in 1+pW$ gives rise, by contravariance, to an
 extension of Dieudonn\'e modules (\cf \ref{para:3.1.1}) 
 \begin{equation*}
  0\rightarrow W\rightarrow \D(\Lambda)\rightarrow W(-1)\rightarrow 0
 \end{equation*}
 which splits canonically due to the rigidity of Dieudonn\'e crystals:
 $\D(\Lambda)=\D(\ovl{\Lambda})=W\oplus W(-1)$.  We denote by
 $e_0,e_{-1}$ the corresponding canonical basis of $\D(\Lambda)$. The
 image $\tau=\mathcal{P}(q)\in W[\frac{1}{p}]$ of $q$ under the period
 mapping is, by definition, the unique element such that $\tau
 e_0+e_{-1}$ generates $F^1$. It turns out that
 \begin{equation*}
  \tau=\log q,
 \end{equation*} 
 thus
 \begin{equation*}
  \mathcal{M}^0=\matheur{D}(1,1^-)
  \xrightarrow{\mathcal{P}=\log} \A^1 =\mathcal{D}.
 \end{equation*}

 Let us indicate a sketch of proof of this remarkable formula (the
 formula goes back to Messing-Katz \cite{q}, \cite{slm}).  The argument
 is more transparent if 
 we think in terms of one-motives instead of elliptic curves: that is,
 rather than viewing $\Lambda$ as the $p$-divisible attached to some
 elliptic curve over $W$, one remarks that it is also the $p$-divisible
 group attached to the one-motive $M=[\Z
 \xrightarrow{1\mapsto q} \G_m]$. Then $\D(\Lambda)$ can be identified
 with $H_{\mathrm{DR}}(M)$, the dual of the Lie algebra of $G^\natural$, where
 $[\Z \rightarrow G^\natural]$ denotes the universal vectorial extension
 of $M$ (\cf \cite{tdh3}); this extension splits canonically
 ($G^\natural=\G_a\times \G_m$), and this splitting agrees with the
 above splitting of $\D(\Lambda)$. We are thus reduced to prove the
 following purely algebraic fact: let us view $M=[\Z
 \xrightarrow{1\mapsto q} \G_m]$ as a
 one-motive over $K((q-1))$; then in terms of the canonical basis
 $e_0,e_{-1}$ of $H_{\mathrm{DR}}(M)$, $F^1\subset
 H_{\mathrm{DR}}(M)$ is generated by the 
 vector $(\log q)e_0+e_{-1}$. This well-known fact is easily proved by
 transcendental means, replacing $K$ by $\C$.
\end{para}

\begin{para}\label{para:6.3.2} Since the basis $e_0,e_{-1}$ is horizontal
 under the Gauss-Manin connection $\nabla$ on $\matheur{D}(1,1^-)\;$ 
 (\ref{lem:6.2.5})\footnote{which is induced by the usual Gauss-Manin
 connection attached to the Legendre elliptic pencil, \cf
 \ref{para:3.6.6} (ii).}, $\tau$ is a quotient 
 of solutions of the Gauss-Manin connection. In terms of the Legendre
 parameter $z=\lambda$ (around the parameter $\zeta_{\mathrm{can}}$ of the
 canonical lifting $A_{\mathrm{can}}$ of the fixed ordinary elliptic curve
 $\ovl{A}$), the equality $\tau = \log q$ means that $q\in
 1+(z-\zeta_{\mathrm{can}})W[[z-\zeta_{\mathrm{can}}]]$ is the
 Dwork-Serre-Tate parameter 
 discussed in I.\ref{sub-serre-tate-para} (for $p\neq 2$).

 This allows to give explicit formulas for $\mathcal{P}$ as a function
 of the algebraic parameter $z$ rather than $q$. For instance, let us
 assume that $p\equiv 1\;(4)$. The disc $| z-\frac{1}{2}| <1$ is
 an ordinary disc. The canonical lifting is $\zeta_{\mathrm{can}}=\frac{1}{2}$
 (complex multiplication by $\Z[\sqrt{-1}]$).  The Gauss-Manin
 connection is represented by the hypergeometric differential operator
 with parameter $(\frac{1}{2},\frac{1}{2},1)$.  Let us recall the
 setting, which we have already met in I.\ref{para-supersingular-disk}.

 We consider the symplectic basis $\omega_1=\bigl[\frac{dx}{2y}\bigr],\;
 \omega_2=\bigl[\frac{(2x-1)dx}{4y}\bigr]$ of relative De Rham
 cohomology.  It has 
 the virtue to be, at $z=\frac{1}{2}$, a basis of eigenvectors for the
 complex multiplication. The Gauss-Manin connection satisfies
 \begin{equation*}
  \omega_2=2z(z-1)\nabla\Bigl(\frac{d}{dz}\Bigr)\omega_1 +
  \frac{4z-5}{6}\;\omega_1.
 \end{equation*} 
 The first row of the
 solution matrix $Y$ normalized by $Y\bigl(\frac{1}{2}\bigr)=\id$ is given by
 \begin{align*}
  y_{11} &=
  F\Bigl(\frac{1}{4},\frac{1}{4},\frac{1}{2};(1-2z)^2\Bigr) +
  \frac{1}{2}(1-2z)F\Bigl(\frac{3}{4},\frac{3}{4},\frac{3}{2};(1-2z)^2\Bigr),\\
  y_{12} &=
  (1-2z)F\Bigl(\frac{3}{4},\frac{3}{4},\frac{3}{2};(1-2z)^2\Bigr).
 \end{align*}
 We have to relate ${\omega_1}_{|_{1/2}},\;
 {\omega_2}_{|_{1/2}}$ to the symplectic basis $e_0,e_{-1}$. Because
 $z=1/2$ gives the canonical lifting, the canonical splitting of the
 Dieudonn\'e module is given by the complex multiplication: $e_0,e_{-1}$
 is a basis of eigenvectors. In fact, taking into account the fact that
 we deal with symplectic bases, we find
 \begin{equation*}
  {\omega_1}_{|_{1/2}}=\frac{\Theta}{2}e_{-1},\quad
   {\omega_2}_{|_{1/2}}=\frac{2}{\Theta}e_0,
 \end{equation*}
 where $\Theta\in W^\times$ denotes the Tate constant
 (for $z=\frac{1}{2}$) discussed in I.\ref{para-tate-para},
 \ref{sub-computation}. (the ambiguity --- by a 
 factor in $\Z_p^\times$ --- is removed since the basis $e_0,e_{-1}$ has 
 been fixed).  The period mapping $\tau=\tau(z)$ is given by the quotient
 of the entries in the first row of the matrix
 $Y.\begin{pmatrix}
     {\Theta/2}&0\\
     0&{2/\Theta}
    \end{pmatrix}$, that is to say:
 \begin{equation*}
  \tau(z)=\Bigl(\frac{\Theta^2}{4}\Bigr)\frac{(1-2z)F(\frac{3}{4},\frac{3}{4},
   \frac{3}{2};(1-2z)^2)}{F(\frac{1}{4},\frac{1}{4},
   \frac{1}{2};(1-2z)^2)+\frac{1}{2}(1-2z)F(\frac{3}{4},
   \frac{3}{4},\frac{3}{2};(1-2z)^2)}
 \end{equation*}
 
 Note that since the target of the period mapping is $\A^1$, the
 denominator does not vanish on $\matheur{D}(\frac{1}{2},1^-)$. Note
 also that, in 
 analytic terms, the etaleness of $\mathcal{P}$ (\ref{thm:6.2.1}) reflects the
 fact that $\frac{d\tau}{dz}$ does not vanish, being the quotient of
 linearly independent solutions of a linear differential equation of
 rank two.
\end{para}

\begin{para}\label{para:6.3.3}
 We now turn to example \ref{thm:4.1.2} ($p$-divisible groups of
 elliptic curves with supersingular reduction).  In this case, the period
 mapping 
 \begin{equation*}
  \mathcal{M}^0\cong \matheur{D}(0,1^-) \xrightarrow{\mathcal{P}}
  \P^1=\mathcal{D} 
 \end{equation*}
 was thoroughly investigated by Gross-Hopkins
 \cite{gross94:_equiv_lubin_tate}, \cite{trapmlsasht}. 
 In the special case considered in \ref{thm:4.1.2}. ii)
 (the fixed supersingular elliptic curve $\ovl{A}$ being defined over
 $\ovl{\F}_p$), they obtained the following formula
 \begin{equation*}
  \mathcal{P}(t)= \lim_{n\rightarrow \infty}
  p^{n}b_{2n}(t)/p^{n+1}b_{2n+1}(t) \in \Q_p[[t]]
 \end{equation*}
 and a closed formula was proposed by J. K. Yu in \cite[11]{Yu95}.

 \noindent In terms of the algebraic Legendre parameter
 $z=\lambda$ in a supersingular disc, it is still true, via
 \ref{lem:6.2.5}, that 
 $\tau=\mathcal{P}(z)$ is a quotient of solutions of the Gauss-Manin
 connection. This allows to give explicit formulas for $\mathcal{P}$ as a
 function of the algebraic parameter $z$. For instance, let us assume
 that $p\equiv 3\;(4)$. The disc $| z-\frac{1}{2}| <1$ is a
 supersingular disc. In terms of a symplectic basis of eigenvectors for
 the action of $[\sqrt{-1}]$ on $V_p^{\C_p}$, we get the following
 formula for the Gross-Hopkins period mapping
 \begin{equation*}
  \tau(z)=\kappa\;\frac{(1-2z)
   F\bigl(\frac{3}{4},\frac{3}{4},\frac{3}{2};(1-2z)^2\bigr)}
   {F\bigl(\frac{1}{4},\frac{1}{4},\frac{1}{2};(1-2z)^2\bigr)
   +\frac{1}{2}(1-2z)F\bigl(\frac{3}{4},\frac{3}{4},\frac{3}{2};(1-2z)^2\bigr)}
 \end{equation*}
 where $\kappa$ is a constant depending on the choice of the basis (a
 natural choice would be to take this basis inside the $p$-adic Betti
 lattice discussed in I.\ref{sub-p-adic-betti-lattice-supersing}; the
 corresponding constant $\kappa$ is then the one occuring in the proof
 of I.\thmref{thm-cm-period}.

 Note that since the target of the period mapping is $\P^1$, the
 denominator vanishes somewhere on $\matheur{D}(\frac{1}{2},1^-)$. 
 \par\noindent Note also that the period mapping in the complex
 situation\footnote{for the Legendre elliptic
 pencil, and computed with help of a basis of cycles --- with coefficients
 in the Gauss integers --- which are eigenvectors for complex
 multiplication at z=1/2.} is given, up to a constant factor, by
 {\it the same hypergeometric formula, viewed as a complex analytic function}. 
\end{para}

\begin{para}\label{para:6.3.4}
 At last, let us turn to example \ref{exa:4.3.4} ($p$-divisible groups of fake
 elliptic curves at a critical prime $p$). In this case, the period mapping 
 \begin{equation*}
  \mathcal{M}^0 \xrightarrow{\mathcal{P}}\Omega_{\C_p} =\mathcal{D}
 \end{equation*}
 is an {\it isomorphism}. This is a reinterpretation
 of Drinfeld's work, \cf \cite{psfpg} (and \cite[II8]{updcds}). Actually,
 Drinfeld establishes an isomorphism at the level of formal schemes 
 \begin{equation*}
  {\mathfrak{M}}^0 \cong \hat{\Omega}\,
  \hat\otimes_{\Z_p}\widehat{\Z_p^{\mathrm{ur}}}
 \end{equation*}
 which is $J^1$-equivariant (where $J^1$ denotes the
 image in $J^{\mathrm{ad}}=PGL_2(\Q_p)$ of the elements of $J$ whose
 determinant is a $p$-adic unit).
\end{para}

\begin{para}\label{para:6.3.5} Recall that the $p$-adic manifold
 $\Omega_{\C_p}$ is simply connected (in the usual topological
 sense). However, the above isomorphism and the modular property of
 $\mathfrak{M}$, allowed Drinfeld to construct a tower of finite etale
 connected Galois coverings of $\Omega_{\C_p}$ (or even of
 $\Omega_{\widehat{\Q_p^{\mathrm{ur}}}}$), as follows.

 One considers the universal special formal module $\udl{\Lambda}$ of
 height $4$ over
 $\mathfrak{M}^0\cong\;\hat\Omega\;\hat\otimes\widehat{\Z_p^{\mathrm{ur}}}$. For
 any $n\geq 1$, $\Ker[p^n]$ is a finite locally free formal group scheme
 of rank $p^{4n}$ over
 $\hat\Omega\;\hat\otimes\widehat{\Z_p^{\mathrm{ur}}}$. This gives rise to a
 finite etale covering of the analytic space
 $\Omega_{\widehat{\Q_p^{\mathrm{ur}}}}$ (fibered in rank-one
 $\mathcal{B}_p/p^n\mathcal{B}_p$-modules).

 To get a connected Galois covering, one looks at the $\Pi$-action
 (recall from \ref{exa:4.3.4} that $\mathcal{B}_p=\Z_{p^2}[\Pi],\;\Pi^2=p$): one
 has $\Ker[p^n]=\Ker[\Pi^{2n}]$, and the complement of $\Ker[\Pi^{2n-1}]$
 in $\Ker[\Pi^{2n}]$ provides a finite etale Galois covering $\Sigma^n$
 of $\Omega_{\widehat{\Q_p^{\mathrm{ur}}}}$ with group
 $(\mathcal{B}_p/p^n\mathcal{B}_p)^\times$. When $n$ grows, these form a
 projective system of Galois coverings (with Galois group the profinite
 completion of $\mathcal{B}_p^\times$). This tower is equivariant with
 respect to the $GL_2(\Q_p)$-action.

 Over $\C_p$, $\Sigma^n$ decomposes into finitely many copies of
 a connected $p$-adic manifold $\Upsilon^n$ which is a finite etale
 Galois covering of $\Omega_{\C_p}$ with group
 $SL_1(\mathcal{B}_p/p^n\mathcal{B}_p)$.  This can be seen either by
 using various maximal commutative $\Q_p$-subalgebras of $B_p$ and
 reasoning as in \cite[1.4.6]{varshavsky98:_shimur}, or by global means, using the
 \v{C}erednik-Drinfeld theorem (\ref{para:7.4.7}), or else by the
 properties of the so-called determinant map \cite[II]{tpuosc}.

Up to now, there is no other method for constructing non-abelian finite
etale coverings of $\Omega_{\C_p}$. The coverings $\Upsilon^n$ remain
rather mysterious: are they simply-connected? Is any connected finite
etale covering of $\Omega_{\C_p}$ a quotient of some $\Upsilon^n$? Does
there exist an infinite etale Galois covering of analytic manifolds
$\Upsilon\rightarrow \Omega_{\C_p}$ with Galois group
$SL_1(\mathcal{B}_p)$?

Drinfeld's construction works for any $GL_n$ over any local field $L$,
replacing $\Omega$ by the complement ${}_L\Omega^{(n-1)}$ of
$L$-rational hyperplanes in $\P^{n-1}$.

\subsection{De Jong's viewpoint on $p$-adic period mappings.}\label{sub:6.4}

 J. de Jong \cite{efgonas} has proposed an interpretation of the
 period mapping which clarifies a lot the nature of its fibers. We
 briefly present this interpretation under the following mild simplifying
 assumption:

\begin{itemize}
 \item[$\ast$] that there is a family of 
 affinoid domains $(X_i)$ of $\mathcal{M}^0$ such that the
 $\mathcal{P}(X_i)$ (which are finite unions of affinoid domains in
 $\mathcal{D}$ since
 $\mathcal{P}$ is etale) form an admissible covering of $\mathcal{D}$. 
\end{itemize}

This assumption is fulfilled in our three basic examples (but not
always, \cf \cite[5.53]{psfpg}).

 Under $(\ast)$, De Jong shows (\loccit\footnote{of course
 $p=\ell$ throughout \S\S  6 and 7 of \loccit }, intr. and \ref{sub:7.2}) that
 there is an etale local system of 
 $\Q_p$-spaces $\mathcal{V}$ on $\mathcal{D}$ such that
 $\mathcal{M}^0$ is a component of the space of
 $\Z_p$-lattices in $\mathcal{V}$. Here, an ``etale local system
 of $\Q_p$-spaces'' is given by the data
 \begin{equation*}
  \mathcal{V}=(\{U_i\hookrightarrow \mathcal{D}\}, \mathcal{V}_i, \phi_{ij})
 \end{equation*}
 where $U_i\hookrightarrow \mathcal{D}$ are open immersions,
 $\mathcal{V}_i=\underset{n}{\limproj}
 \;\mathcal{V}_i/p^n\mathcal{V}_i$ are etale local systems of
 $\Z_p$-lattices (each $\mathcal{V}_i/p^n\mathcal{V}_i$ being a finite
 locally free sheaf of $\Z/p^n\Z$-modules on the etale site of $U_i$),
 and $ \phi_{ij}: (\mathcal{V}_i\otimes \Q_p)_{|
 U_i\times_\mathcal{D}U_j} \rightarrow (\mathcal{V}_j\otimes \Q_p)_{|
 U_i\times_\mathcal{D}U_j} $ are isomorphisms satisfying the usual
 cocycle condition (\cf \cite[4]{efgonas}).

 Of course, in the ``decorated case'', $\mathcal{V}$ inherits the
 relevant decoration.

 It follows from this interpretation that in example \ref{thm:4.1.2}, the
 fibers of the period mapping are in natural bijection with the set of
 vertices of the Bruhat-Tits tree $\mathcal{T}$.

 Similarly, in example \ref{exa:4.3.4}, one can see that the fibers of
 $\mathcal{P}$ are in bijection with the 
 one-point set $SL_1(B_p)/SL_1(\mathcal{B}_p)$, \cf \ref{para:7.4.4}.
\end{para}

\newpage
\section{$p$-adic uniformization of Shimura varieties.}\label{sec:7}

\begin{abst}
 Application of the theory of $p$-adic period mappings: Rapoport-Zink's
 $p$-adic uniformization theorem. The case of global uniformization.
 \v{C}erednik-Drinfeld uniformization of Shimura curves and
 other examples.
\end{abst}

\subsection{$p$-integral models of Shimura varieties.}\label{sub:7.1}

\begin{para}\label{para:7.1.1} Let $\Sh$ be a Shimura variety of PEL type as
 in \ref{para:1.1.4}. We recall the setting of \ref{para:1.1.4}: $B$ is
 a simple finite-dimensional $\Q$-algebra, with a positive involution
 $\ast$, $V$ is a $B$-module 
 of finite type, endowed with an alternate $\Q$-bilinear form such that
 $\an{bv,w}=\an{v,b^\ast w}$; 
 $G$ denotes the $\Q$-group of
 $B$-linear symplectic similitudes of $V$; $\mathcal{B}$ is a maximal
 order in $B$ stable under $\ast$, and 
 $L$ is a lattice in $V$, stable under $\mathcal{B}$ and autodual for
 $\an{\;,\;}$.

 The (non-connected) Shimura variety $\Sh$ is defined over the reflex
 field $E$ (which is a number field in $\C$).  There are variants of
 $\Sh$ in which the principal congruence subgroup of level $N$ of
 $G(\hat\Z)$ is replaced by any open compact subgroup $C$ of $G(\A_f)$,
 where $\A_f=\hat\Z\otimes \Q$ is the ring of finite adeles of $\Q$. The
 complex points of $\Sh$ admit the following adelic description
\begin{equation*}
\Sh(\C)=C\backslash \bigl((G(\R)^{\mathrm{ad}}.h_0) \times G(\A_f)\bigr) /G(\Q)
\end{equation*}
(\cf \ref{para:1.2.4} about the one-parameter group $h_0$). 
\end{para}

\begin{para}\label{para:7.1.2} We set $B_p=B\otimes_{\Q}
 \Q_p,\;\mathcal{B}_p=\mathcal{B}\otimes_\Z \Z_p,$ and so on. We denote
 by $C_p$ the stabilizer of $L_p$ in $G(\Q_p)$

 {\it We assume that $C$ is of the form $C^pC_p$} for some open compact
 subgroup $C^p$ of $G(\A_f^p)$, where $\A_f^p=(\prod_{\ell\neq
 p} \Z_\ell)\otimes \Q$.

 In the case of a congruence subgroup of level $N$, this condition means
 that {\it the level is prime to $p$}. 

 We fix an embedding $v: E\hookrightarrow
 \C_p$, denote by $E_v$ the $v$-completion of $E$, by $\mathfrak{v}$ the
 discrete valuation ring
 $\widehat{\Z^{\mathrm{ur}}}.\mathcal{O}_{E_v}=
 \widehat{\Z^{\mathrm{ur}}}.\mathcal{O}_E$ 
 (compositum in $\C_p$), and by $K$ the fraction field of $\mathfrak{v}$.
\end{para}

\begin{para}\label{para:7.1.3} Let $S$ be a scheme. The objects of the the
 category of $S$-abelian schemes {\it up to prime to $p$ 
isogeny} are the abelian schemes over $S$; morphisms between $A $ and
$A'$ are global sections of $ \udl{\Hom}_S(A, A')\otimes_{\Z}\Z_{(p)}$,
where $\Z_{(p)}$ is the ring of rational numbers with denominator prime
to $p$. This is a $\Z_{(p)}$-linear category. Its relation to
$p$-divisible groups is given by the following elementary lemma:
\end{para}

\begin{lem}\label{lem:7.1.4} The faithful functor $A\mapsto A[p^\infty]$
\begin{equation*}
 \{\text{$S$-abelian schemes up to prime to $p$
  isogeny}\}\rightarrow \{\text{$p$-divisible groups over $S$}\}
\end{equation*}
is conservative, \ie reflects isomorphisms.

\end{lem}
The following result (which we state in a rather vague manner) is due to
R. Kottwitz \cite{possvoff} (Compare with \ref{para:1.2.4}.).

\begin{pro}\label{pro:7.1.5}
 For sufficiently small $C^p$, there is a moduli scheme $\Sh$ over
 $\mathcal{O}_{E_v}$ for isomorphism classes of abelian 
 varieties $A$ up to prime to $p$ isogeny, with $\mathcal{B}$-action 
 (with Shimura type condition), together with a $\Q$-homogeneous
 principal $\ast$-polarization and a class of $B$-linear symplectic
 similitudes $\mathrm{H}^1_{\mathrm{et}}(A,\A_f^p)\rightarrow \A_f^p\otimes V$
 modulo the group $C^p$. The generic fiber $\Sh
 \otimes_{\mathcal{O}_{E_v}}E_v$ is a finite sum of copies of
 $\Sh\otimes_E {E_v}$.
\end{pro}

In the case of a congruence subgroup of level $N$, this condition that
$C^p$ is sufficiently small means that $N$ is sufficiently large (and
kept prime to $p$).

\subsection{The uniformization theorem.}\label{sub:7.2}

\begin{para}\label{para:7.2.1}
 One fixes a point of $\Sh(\ovl{\F}_p)$, hence a
 $\ovl{\F}_p$-abelian variety $\ovl{A}$ up to prime to $p$ isogeny with
 $\mathcal{B}$-action, principal $\ast$-polarization, and level
 structure. We consider the $p$-divisible group
 $\ovl{\Lambda}={\ovl{A}}[p^\infty]$ (with $\mathcal{B}_p$-action,
 principal 
 $\ast$-polarization, level structure). We denote by $\ovl{G}$ the
 $\Q$-group of self-quasi-isogenies of ${\ovl{A}}$ respecting the
 additional structure. It is clear that $\ovl{G}(\Q_p)\subset J$; on the
 other hand $\ovl{G}(\A_f^p)\subset G(\A_f^p)$.

 We also consider an auxiliary lifting $\tilde{\ovl{A}}$ of $\ovl{A}$ (with
 $\mathcal{B}$-action) over some finite extension $\mathcal{O}_{E'}$ of
 $\mathcal{O}_E$ in $\mathfrak{v}$, and set
 $\tilde{\ovl{\Lambda}}=\tilde{\ovl{A}}_\mathfrak{v}[p^\infty]$.  We can
 then identify $\D(\ovl{\Lambda})\otimes_W K$ with
 $\mathrm{H}^1_{\mathrm{DR}}(\tilde{\ovl{A}})\otimes_{\mathcal{O}_{E'}}
 K$, and the submodule $\omega_{\tilde{\ovl{\Lambda}}}$ with
 $F^1\mathrm{H}^1_{\mathrm{DR}}(\tilde{\ovl{A}})\otimes_{\mathcal{O}_{E'}}
 K$.  Thus if the complex lagrangian space 
 $F^1_0$ of \ref{para:1.1.4} is taken to be
 $F^1\mathrm{H}^1_{\mathrm{DR}}(\tilde{\ovl{A}})\otimes_{\mathcal{O}_{E'}}
 \C$ (for some complex embedding of $E'$ 
 extending the natural complex embedding of $E$), the complex and
 $p$-adic Shimura type conditions ``agree'' (\ref{para:1.1.4}, \ref{para:4.3.2}). We can
 consider the corresponding decorated moduli problem for $p$-divisible
 groups, and the formal moduli scheme $\mathfrak{M}$ over $\mathfrak{v}$
 as in \ref{sub:4.3}.

 \noindent Moreover, the flags spaces $\mathcal{D}^\vee$ considered in 
 the complex and $p$-adic situations agree, \ie come from the same flag 
 space defined over the reflex field $E$.

\end{para}

\begin{para}\label{para:7.2.2}
 One defines a morphism of functors on $\Nil_\mathfrak{v}$: 
 \begin{equation*}
  {\mathfrak{M}} \rightarrow \Sh
 \end{equation*}
 as follows \cite[6.14]{psfpg}.
 Let $S$ be in $\Nil_\mathfrak{v}$ and let us consider a pair
 $(\Lambda,\rho)\in \mathfrak{M}(S)$. By rigidity of $p$-divisible
 groups up to isogeny (\ref{thm:2.2.3}), the
 $\mathcal{B}_p$-quasi-isogeny $\rho 
 \in \qisog(\Lambda\times_S S_{\mathrm{red}},
 \ovl{\Lambda}\times_{{\ovl{\F}}_p}S_{\mathrm{red}})$ lifts to a unique 
 $\mathcal{B}_p$-quasi-isogeny $\tilde{\rho}\in \qisog(\Lambda,
 \tilde{\ovl{\Lambda}}\times_\mathfrak{v}S)$. This can be algebraized in
 a weak sense: there is an $S$-abelian scheme $A_{(\Lambda,\rho)}$ up to
 prime to $p$ isogeny, with $\mathcal{B}$-action and with $p$-divisible
 group $\Lambda $, and a quasi-isogeny $A_{(\Lambda,\rho)}\rightarrow
 \tilde{\ovl{A}}\times_{\mathcal{O}_{E'}}S$ which induces $\tilde{\rho}$
 at the level of $p$-divisible groups. Moreover $A_{(\Lambda,\rho)}$ is
 unique in the category of $S$-abelian schemes up to prime to $p$ isogeny
 with $\mathcal{B}$-action (\ref{lem:7.1.4}), and its formation is
 functorial in 
 $(\Lambda,\rho)$. One then defines the required morphism of functors on
 $\mathfrak{M}(S)$ by setting $(\Lambda,\rho)\mapsto A_{(\Lambda,\rho)}$
 (endowed with the polarization and level structure inherited from $\ovl{A}$).

 This defines a morphism of formal schemes, hence a morphism
 of associated analytic spaces over $\C_p$:
 \begin{equation*}
  \mathcal{Q}:\;\mathcal{M}\;\rightarrow\; \Sh_{\C_p}^{\mathrm{an}}.
 \end{equation*}
\end{para}

\begin{para}\label{para:7.2.3}
 It turns out that the decorated abelian varieties (up to prime to $p$
 isogeny) which are {\it isogenous} to $\ovl{A}$ form in a natural way a
 {\it countable union $Z=\bigcup Z_n$ of projective irreducible
 subvarieties $Z_n$ of $\;\Sh_{\ovl{\F}_p}$}, each of them
 meeting only finitely many members of the family (\cite[6.23, 6.34]{psfpg}).
 One can slightly generalize the notion of tubes (\ref{sub:3.5}) and
 define the tube $]Z[\subset \Sh_{\C_p}^{\mathrm{an}}$. We denote by $S$ the
 connected component of $]Z[$ which contains the modular point $s$ of our
 chosen lifting $\tilde{\ovl{A}}$.
\end{para}

\begin{thm}\label{thm:7.2.4}
 For sufficiently small $C^p$, there is a connected component
 $\mathcal{M}^0$ of $\mathcal{M}$ such that the restriction
 $\mathcal{Q}:\;\mathcal{M}^0\rightarrow S$ is a topological covering of
 $S$. More precisely, $S$ is the (right) quotient of $\mathcal{M}^0$ by
 some torsion-free discrete subgroup of $J$.
\end{thm}

\cf \cite[6.23, 6.31]{psfpg}. Roughly speaking, the proof of \ref{thm:7.2.4}
involves two steps: using the Serre-Tate theorem \ref{thm:2.2.4}, one
shows that 
the morphism of $\mathfrak{v}$-formal schemes $\mathfrak{M}$
$\rightarrow \widehat{\Sh}$ defined above is
formally etale. One concludes by a careful study of this morphism on
geometric points.

In our examples \ref{sub:4.2}, $S$ is an open unit
disc, and the map $\mathcal{Q}$ is an isomorphism (in \ref{para:4.2.1}, it is
induced by $q\mapsto z=\lambda$).

\begin{para}\label{para:7.2.5} Let us consider the universal decorated
 abelian scheme $\udl{A}$ over $\Sh$, and the vector bundle
 $\mathcal{H}=\mathrm{H}^1_{\mathrm{DR}}(\udl{A}/\Sh)$ with its
 Gauss-Manin connection 
 $\nabla_{\mathrm{GM}}:\mathcal{H}\rightarrow \Omega^1_{\Sh}\otimes
 \mathcal{H}$. 
\end{para}

\begin{cor}\label{cor:7.2.6}  Viewed as a $p$-adic connection and
 restricted to $S$, the Gauss-Manin connection 
 ${\nabla_{\mathrm{GM}}}_{| S}$ comes from a representation of
 $\pi_1^{\mathrm{top}}(S,s)$ on $V^{\C_p}\;$ (hence it satisfies
 Cauchy's theorem, \cf I.\ref{sub-connection}).
\end{cor}

\begin{proof}
 Let $\tilde{S}$ be the universal covering of $S$. According to
 I.\ref{sub-connection}), we have to show that the pull-back of 
 $\nabla_{\mathrm{GM}}$ on $\tilde{S}$ is trivial ($\cong V^{\C_p}\otimes
 \mathcal{O}_{\tilde{S}})$. By \ref{thm:7.2.4} and \ref{lem:6.2.5}, this
 already holds on the quotient $\mathcal{M}^0$ of $\tilde{S}$.
\end{proof}

 \noindent {\it Remark.} Using the Weil descent data of \cite[6.21]{psfpg},
 one can show that the representation is defined over a
 finite unramified extension $E'_v$ of $E_v$, \ie comes from a
 representation of $\pi_1^{\mathrm{top}}(S,s)$ on $V\otimes E'_v$.

\begin{para}\label{para:7.2.7}
 The assumption $C=C^pC_p$ is necessary to define the morphism
 $\mathcal{Q}: \mathcal{M}^0\rightarrow S$. The assumption ``$C^p$ small
 enough'' is much less important. Theorem \ref{thm:7.2.4}. holds without it,
 except that the discrete subgroup of $J$ is then no longer torsion-free
 [the point is that there is always a normal subgroup ${C'}^p$ of finite
 index which is small enough, one applies the theorem to ${C'}^p$, and one
 passes to the quotient by ${C'}^p/C^p$]. However, to give a modular
 interpretation of the (ramified) quotient of $\mathcal{M}^0$ in this
 case, one needs the notion of $p$-adic orbifold (to be discussed in
 III); moreover the Gauss-Manin connection may acquire logarithmic
 singularities along the branched locus.

\subsection{Global uniformization}\label{sub:7.3}

We say that there is {\it global uniformization}
in case $S$ is (the analytification of) a whole connected component of
the Shimura variety $\Sh$. This occurs when the reduction mod $p$ of
the decorated abelian varieties $A$ parametrized by $\Sh$ form a 
{\it single isogeny class}. 

In the situation of global uniformization, the group $\ovl{G}$ of
self-quasi-isogenies of ${\ovl{A}}$ (with decoration) is an
{\it inner form of $G$}.

In view of our analogy \ref{para:4.3.5}, this situation can be expected to be
closest to the complex case (in the notation of \ref{para:1.2.4}, the $\hat A$'s
are all isomorphic to $\hat U$). Indeed, we then get a commutative
diagram of $p$-adic manifolds (similar to \ref{para:1.2.4}):
\begin{equation*}
 \begin{CD}
  \mathcal{M}^0@>\mathcal{P}>> \mathcal{D} @. \subset\;\mathcal{D}^{\vee}\\
  @V\mathcal{Q}VV @VV\mathcal{Q}V @.\\
  S@>\mathcal{P}>> \mathcal{D}/\Gamma.
 \end{CD}
\end{equation*}
 where the horizontal maps are etale, $\mathcal{Q}$ denotes the quotient
 maps, and $\Gamma$ is an arithmetic subgroup in the semi-simple group
 $\ovl{G}^{\mathrm{ad}}$. In all known examples, $\mathcal{P}$ is
 actually an isomorphism and $\mathcal{M}^0=\tilde S$ 
 \cite[6]{psfpg}, \cite{tpuosc}, \cite{varshavsky98:_shimur},
 \cite{varshavsky98:_shimur_ii}.

Here, the archimedean---non archimedean correspondence goes beyond a mere
analogy: indeed, here and in \ref{para:1.2.4}, what is denoted by $S$ (\resp
$\mathcal{D}^\vee $) is the --- complex or $p$-adic --- analytification of
the same algebraic variety.

\medskip
\noindent {\it Remark.} To avoid any ambiguity, let us point out
that we consider here only phenomena of global uniformization which can
be interpreted in terms of moduli spaces of $p$-divisible groups. This
does not account for all cases of global $p$-adic uniformization of
Shimura varieties. For instance, it is known that the elliptic modular
curve $X_0(p)$ is a Mumford curve over $\Q_p$, whereas the $p$-adic
uniformization of \ref{thm:7.2.4} (or variants with level $p$ structure)
is local.
\end{para}

\subsection{Example: the \v{C}erednik-Drinfeld uniformization of
Shimura curves at a critical prime.}\label{sub:7.4}

\begin{para}\label{para:7.4.1}
 The first historical example of global uniformization in the above
 sense is the \v{C}erednik uniformization of ``Shimura
 curves'' (modular curves for fake elliptic curves) at a critical prime
 $p$, \cf \cite{uoacbdsopwcq}.  \v{C}erednik's method was
 of group-theoretic nature and did not involve formal groups. A modular
 proof was subsequently proposed by Drinfeld \cite{copsr} (\cf also
 \cite[III]{updcds}). For examples, and related topics, we refer to
 \cite{van89:_les_shimur}, \cite{van92:_discr_mumfor}, \cite{van92:_unifor}.
\end{para}

\begin{para}\label{para:7.4.2}
 We fix an indefinite quaternion algebra $B$ over $\Q$, a maximal order
 $\mathcal{B}$ in $B$, and a positive involution $\ast$ of
 $\mathcal{B}$. We consider the standard representation $V=B$ of $B$, so
 that $G=GL_B(V)\cong B^\times$.

 We refer to the foundational paper \cite{shimura67:_const} for the
 properties of the (connected) Shimura curve ${\mathcal{X}^+(N)}$
 of level $N$ attached to $B$. It parametrizes principally 
 $\ast$-polarized abelian surfaces with multiplication by $\mathcal{B}$
 (fake elliptic curves, \cf \ref{exa:4.3.4}), with level $N$ structure.
 We recall the following properties of ${\mathcal{X}^+(N)}$, which 
 actually characterize it (\loccit 3.2):

\begin{enumerate}
\renewcommand{\theenumi}{\alph{enumi}}
 \item ${\mathcal{X}^+(N)}$ is a projective geometrically connected smooth
       curve defined over the cyclotomic field $\Q(\zeta_N)$,
 \item ${\mathcal{X}^+(N)}(\C)= \mathfrak{h}/\Gamma^+(N)$, that is to
       say: ${\mathcal{X}^+(N)}$ is the quotient of the Poincar\'e upper
       half plane by the image $\Gamma^+(N)=PSL_1(\mathcal{B})(N)$ of
       the group $\{g\in (1+N\mathcal{B})^\times, \Nr(g)>0\}$ in
       $PSL_2(\R)$, where $\Nr$ stands for the reduced norm,
 \item for any imaginary quadratic field $K$ in $B$ such that
       $\mathcal{O}_K\subset \mathcal{B}$, let $z$ be its fixed point in
       $\mathfrak{h}$. Then the extension of $\Q(\zeta_N)$ generated by the
       image of $z$ in ${\mathcal{X}^+(N)}$ is the class field over $K$ with
       conductor $N\mathcal{B}\cap \mathcal{O}_K$.
\end{enumerate}

 The commutative diagram of \ref{para:1.2.4}, involving the complex period 
 mapping, specializes to the following 

 \begin{equation*}
  \begin{CD}
   \tilde{S} @>\mathcal{P}>> \mathcal{D}=\mathfrak{h} @.
   \subset \; \mathcal{D}^{\vee}=\P^1\\
   @V\mathcal{Q}VV @VV\mathcal{Q}V\\
   S={\mathcal{X}^+(N)}_\C^{\mathrm{an}}@>\mathcal{P}>>\mathcal{D}/\Gamma^+(N)
  \end{CD}
 \end{equation*}
 where the horizontal maps are isomorphisms (and the
 vertical maps are topological coverings if $\Gamma^+_p(N)$ is
 torsion-free\footnote{if this group is not torsion-free $\tilde S$ 
 is the universal covering in the sense of orbifolds.}).

 There is a useful variant $\mathcal{X}^\ast$ of
 $\mathcal{X}^+={\mathcal{X}^+(1)}$, already considered by Shimura
 (\cite[3.13]{shimura67:_const}, in which 
 $\Gamma^+$ is replaced by the image $\Gamma^\ast$ in $PSL_2(\R)$ of the
 group $\{g\in B^\times \mid g\mathcal{B}=\mathcal{B}g,\,\Nr(g)>0\}$:
 \begin{equation*}
  \mathcal{X}^\ast(\C)= \mathfrak{h}/\Gamma^\ast.
 \end{equation*}
 It is defined over $\Q$. More generally, there is a
 Shimura curve $\mathcal{X}_\Gamma$ attached to any 
 congruence subgroup $\Gamma$ in $B^{\times^{+}}/\Q^\times$, where
 $B^{\times^{+}}$ denotes the group of elements of positive 
 reduced norm.
\end{para}

\begin{para}\label{para:7.4.3}
 It is known that ${\mathcal{X}^+(N)}$ has
 good reduction at any prime $p$ which does not divide the discriminant 
 $d(B)$ of $B$ nor the level $N$ (Y. Morita). If instead $p$ is critical
 (but still {\it does not divide $N$}), as we shall 
 assume henceforth, \v{C}erednik's theorem shows that
 ${\mathcal{X}^+(N)}_{\C_p}$ is a {\it Mumford curve} uniformized by 
 the Drinfeld space $\Omega_{\C_p}=\P^1_{\C_p}\setminus \P^1(\Q_p)$ (at
 least if $N$ is large enough); in fact, $\Sh$ becomes a 
 Mumford curve over some finite unramified extension of $\Q_p$ (note that
 $\Q(\zeta_N)\subset {\Q_p^{\mathrm{ur}}}$).  

 To be more explicit, let $\ovl{B}$ be the definite quaternion
 $\Q$-algebra which is ramified at the same primes as $B$ except $p$:
 $\ovl{B}_p\cong M_2(\Q_p)$. Let $\ovl{\mathcal{B}}$ be a maximal order in
 $\ovl{B}$. Let $\Gamma^+_p(N)$ be the image in $PGL_2(\Q_p)$ of the
 group $\{g\in (1+N\ovl{\mathcal{B}}[\frac{1}{p}])^\times \mid \ord_p(\Nr(g))
 \;{\rm even}\}$. Let $\Gamma^\ast_p$ be the image in $PGL_2(\Q_p)$ of
 the normalizer $\{g\in {\ovl{B}}^\times \mid g\ovl{\mathcal{B}}[\frac{1}{p}] =
 \ovl{\mathcal{B}}[\frac{1}{p}]g\}$.
 These are discrete subgroups, and one has the {\it \v{C}erednik
 uniformization}:
 \begin{align*}
  {\mathcal{X}^+(N)}_{\C_p}^{\mathrm{an}}&\cong\Omega_{\C_p}/\Gamma^+_p(N),\\
  {\mathcal{X}^{\ast}}_{\C_p}^{\mathrm{an}}&\cong\Omega_{\C_p}/\Gamma^\ast_p.
 \end{align*} 
 For recent applications to $p$-adic $L$-functions, Heegner points,
 etc... \cf \cite{bertolini98:_heegn_l_cered}.
\end{para}

\begin{para}\label{para:7.4.4}
 The case of $\mathcal{X}^{\ast}$ is the one originally considered in 
 \cite{uoacbdsopwcq}, except that \v{C}erednik uses the
 adelic language. In order to make the translation, it is useful to keep
 in mind the following few facts (\cf \cite[pp. 40, 99 and passim]{adadq}).

 For critical $p$, the unique maximal order
 $\mathcal{B}_p$ of $B_p$ is $\{b\in B_p \mid \Nr(b)\in
 \Z_p\}$. Hence, we have an equality $SL_1(\mathcal{B}_p) = SL_1(B_p)$
 and a chain of compact normal subgroups
 \begin{equation*}
  PSL_1(\mathcal{B}_p) = PSL_1(B_p)\subset PGL_1(\mathcal{B}_p) \subset
   PGL_1(B_p)
 \end{equation*}
 with $PGL_1(B_p)/PSL_1(\mathcal{B}_p)\cong \Q_p^\times/(\Q_p^\times)^2$
 (which is of type $(2,2)$ if $p\neq 2$, of type $(2,2,2)$ if $p=2$). 

 On the other hand, the normalizer of $\mathcal{B}_\ell$ in
 $B_{\ell}^{\times}$ is $B_{\ell}^{\times}$ if $\ell$ is critical, and
 $\mathcal{B}_\ell^\times.\Q_{\ell}^{\times}$ otherwise. It follows that  

 \begin{align*}
  \Gamma^\ast&=\{g\in B^{\times^+}\mid \forall \ell
  \nmid d(B), \; g\in \mathcal{B}_\ell^\times.\Q_\ell^\times\}/\Q^\times\\
  \Gamma^\ast_p&=\{g\in \ovl{B}^\times\mid \forall
   \ell \nmid d(B), \; g\in
   \ovl{\mathcal{B}}_\ell^\times.\Q_\ell^\times\}/\Q^{\times}\\
  \Gamma^+(N)&= \{g\in B^{\times^+}\mid \forall \ell, \; g\in
   (1+N\mathcal{B}_\ell)^\times.\Q_\ell^{\times}\}/\Q^{\times}\\
  \Gamma^+_p(N)&= \Bigl\{g\in \ovl{B}^\times\Bigm|
   \begin{array}{l}
    \ord_p(\Nr(g))\;{\rm even};\\
    \forall \ell \neq p,\; g\in
    (1+N\ovl{\mathcal{B}}_{\ell})^\times.\Q_\ell^\times\\ 
   \end{array}
  \Bigr\}\Big/\Q^{\times}
 \end{align*}
 If $N$ is prime to $d(B)$, this can also be rewritten
 \begin{align*}
 \Gamma^+(N)&= \Bigl\{g\in B^{\times^+}\Bigm|
  \begin{array}{l}
   \forall\ell\mid d(B),\, \ord_\ell(\Nr(g))\;{\rm even};\\
   \forall\ell\nmid d(B),\, g\in (1+N\mathcal{B}_\ell)^\times.\Q_\ell^\times\\
  \end{array}
  \Bigr\}\Big/\Q^\times \\
  \Gamma^+_p(N)&= \Bigl\{g\in \ovl{B}^\times\Bigm|
  \begin{array}{l}
   \forall\ell\mid d(B),\, \ord_\ell(\Nr(g))\;{\rm even};\\
   \forall\ell\nmid d(B),\,
    g\in (1+N\ovl{\mathcal{B}}_{\ell})^{\times}.\Q_{\ell}^{\times}\\
  \end{array}
  \Bigr\}\Big/\Q^{\times}
 \end{align*}
 In particular, for $N=1$, we see that there are canonical bijections 
 \begin{equation*}
 \Gamma^{\ast}/\Gamma^+=\Gamma^{\ast}_p/\Gamma^+_p =
  (\Z/2\Z)^{\sharp\{\ell| d(B)\}}
 \end{equation*} 
\end{para}

\begin{para}\label{para:7.4.5} We now give some indications on Drinfeld's
 modular approach to \v{C}erednik's uniformization.  
\end{para}

\begin{lem}\label{lem:7.4.6}
 \begin{enumerate}
  \item Fake elliptic curves have potentially good reduction at every
	prime; the reduction is supersingular reduction at the critical
	prime $p$. In particular, over $\ovl{\F}_p$,
	they form a unique isogeny class.
  \item Fake elliptic curves have a unique principal $\ast$-polarization.
 \end{enumerate}
\end{lem}

\begin{proof}
 (i): due to the presence of many endomorphisms, the potential good
 reduction follows easily from the semi-stable reduction theorem for
 abelian schemes; also, the non-split quaternion algebra $B_p$ acts on
 the $p$-divisible group up to isogeny, and this forces the slopes to be
 $1/2$; for details and (ii), \cf \cite[III]{updcds}.
\end{proof}

 It follows from (i) that any fake elliptic curve over
 $\ovl{\F}_p$ is isogenous to the square of a supersingular elliptic
 curve $\ovl{E}$.  However, while there are only finitely many
 supersingular elliptic curves up to isomorphism, fake elliptic curves
 over $\ovl{\F}_p$ have continuous moduli! The moduli spaces are chains
 of $\P^1$'s.\footnote{Another curious feature of the higher
 dimensional case is that the product of $n > 1$ supersingular elliptic
 curves over $\ovl{\F}_p$ form a single isomorphism class (Deligne,
 Ogus)}

 Let us introduce the quaternion $\Q$-algebra $D$ ramified only
 at $p$ and $\infty$: $D\cong \End^\Q(\ovl{E})$. 
 Let $\ovl{A}$ be a fake elliptic curve over $\ovl{\F}_p$. It follows from
 \ref{lem:7.4.6} (i) that $\End^\Q(\ovl{A})\cong M_2(D)$, and $\End_{B}^\Q
 (\ovl{A})\cong \ovl{B}$.  

 In particular, the inner form $\ovl{G}$ of $G$ considered in \ref{sub:7.3} is
 $\ovl{B}^\times$, and $\ovl{G}(\Q_p)=GL_2(\Q_p)=J$.  

 It follows from \ref{lem:7.4.6} (ii) that we may disregard polarizations. The
 $p$-divisible groups attached to fake elliptic curves with
 $\mathcal{B}$-multiplication are special formal $\mathcal{B}_p$-modules
 of height $4$; recall from \ref{exa:4.3.4} the corresponding formal
 moduli space 
 $\displaystyle \mathfrak{M} = \coprod_\Z\;\mathfrak{M}^0$. The period
 mapping is an isomorphism (\ref{para:6.3.4}, \ref{sub:6.4}), and the
 commutative diagram of \ref{sub:7.3} specializes to

\begin{equation*}
 \begin{CD}
  \tilde{S}=\mathcal{M}^0 @>\mathcal{P}>> \mathcal{D}=\Omega_{\C_p} @.
  \subset\; \mathcal{D}^{\vee}=\P^1\\
  @V\mathcal{Q}VV  @VV\mathcal{Q}V @.\\
  S={\mathcal{X}^+(N)}_{\C_p}^{\mathrm{an}} @>\mathcal{P}>>
  \mathcal{D}/\Gamma^+_p(N)
 \end{CD}
\end{equation*}
 where the horizontal maps are isomorphisms (and the
 vertical maps are topological coverings if $\Gamma^+_p(N)$ is
 torsion-free), in complete analogy with \ref{sub:7.2}.
 
\begin{para}\label{para:7.4.7} Let us end this section by considering the
 case where {\it $p$ divides the level $N$}: $N=p^n.M,\; p\nmid M$. 

 In this case, $\Gamma^+_p(N)=\Gamma^+_p(M)$. The finite etale covering
 ${\mathcal{X}^+(N)}_{\C_p}^{\mathrm{an}}\rightarrow
 {\mathcal{X}^+(M)}_{\C_p}^{\mathrm{an}}$ is not a topological covering anymore
 (even if $\Gamma^+_p(M)$ is torsion-free). As was shown by Drinfeld, it
 is hidden in the finite etale covering $\Upsilon^n$ of $\Omega_{\C_p}$
 (\ref{para:6.3.5}).

 Let us assume for simplicity that $M>2$, so that $-1\notin
 \Gamma^+(M)$. This implies that $\Gamma^+(N)/\Gamma^+(M)\cong
 SL_1(\mathcal{B}/p^n\mathcal{B})$. Then $\Upsilon^n$ is a topological
 covering of ${\mathcal{X}^+(N)}_{\C_p}^{\mathrm{an}}$, and more precisely, one
 has the commutative diagram

 \begin{equation*}
 \begin{CD}
  {\mathcal{X}^+(N)}_{\C_p}^{\mathrm{an}}@>\sim>> \Upsilon^n/\Gamma^+_p(M)\\
  @VVV @VVV\\
  {\mathcal{X}^+(M)}_{\C_p}^{\mathrm{an}} @>\mathcal{P}>>
  \Omega_{\C_p}/\Gamma^+_p(M)
 \end{CD}
 \end{equation*}
 where the vertical maps are finite etale Galois
 covering maps with group $SL_1(\mathcal{B}/p^n\mathcal{B})$, \cf 
\cite[III.5.5]{updcds}, \cite{varshavsky98:_shimur}.  Similarly, there is such a $p$-adic global
 uniformization for any Shimura curve $\mathcal{X}_\Gamma$ (\cf
 \ref{para:7.4.1}).
\end{para}

\subsection{Other examples of global $p$-adic uniformization.}\label{sub:7.5}

\begin{para}\label{para:7.5.1}
 Shimura's paper \cite{shimura67:_const} deals not only with quaternion algebra
 over $\Q$, but also with quaternion algebras $B/E$ over a totally real
 number field $E$. When $B$ {\it splits at a single real place,} he
 constructs curves ${\mathcal{X}^+(N)}$ and ${\mathcal{X}^\ast}$
 satisfying properties similar to (a),(b),(c) of \ref{para:7.4.2}.

 In this context, $N$ may be an arbitrary ideal of $\mathcal{O}_E$,
 $\Q(\zeta_N)$ has to be replaced by the ray class field with conductor
 $N.\infty$, and the condition $N(g)>0$ means that the reduced norm of
 $g$ is totally positive, or, what amounts to the same, is positive at
 the split real place $\infty_0$.\index{000infinity0@$\infty_0$}
 
 The novelty, however, is that these Shimura curves have {\it no
 direct interpretation} as moduli spaces for
 decorated abelian varieties, since $B$ is not in the Albert list of
 endomorphism algebras of abelian varieties when $[E:\Q]>1$.

 To overcome this problem, Shimura introduces the following remarkable
 trick. Let $K$ be a totally imaginary quadratic extension of $E$. Then
 $B^\bullet =B\otimes_E K$ is in Albert's list (we endow it with the
 involution which is the tensor product of $\ast$ on $B$ and complex
 conjugation on $K$). It turns out that the Shimura varieties of PEL
 type for abelian varieties of dimension $g=4[E:\Q]$ ($\ast$-polarized,
 with level structure) with complex multiplication by $B^\bullet$ are
 twisted forms of the same curve. The Shimura curve attached to $B$ is
 obtained by Weil descent from the Shimura curves attached to
 $B\otimes_E K$, for sufficiently many quadratic extensions $K$.
\end{para}

\begin{para}\label{para:7.5.2}
 \v{C}erednik's original $p$-adic uniformization theorem
 already applies to Shimura curves attached to quaternion algebras over
 a totally real number field $E$. His method has been further developed
 in \cite{varshavsky98:_shimur}, \cite{varshavsky98:_shimur_ii}. On the
 other hand, Drinfeld's modular 
 approach has been extended by Boutot-Zink \cite{tpuosc}, using Shimura's
 trick. Let us sketch the principle of the construction.

 Let $G$ be the multiplicative group of $B^\times$, viewed as an
 algebraic group over $\Q$. Let $\Gamma$ be a congruence subgroup of
 $G^{\mathrm{ad}}(\Q)$ (alternatively, one can consider an open compact
 subgroup 
 $C$ of $G(\A_f)$). Let us first assume that $\Gamma$ is {\it maximal at
 $p$} (alternatively $C$ is of the form $C^pC_p$ for some open compact
 subgroup $C^p$ of $G(\A_f^p)$, \cf \ref{para:7.1.2}). The (non-connected) Shimura
 curve $\Sh$ given by
 \begin{equation*}
  \Sh_\C^{\mathrm{an}}= C\backslash \bigl((\P_\C^1\setminus\P^1(\R)) \times
   G(\A_f)\bigr) /G(\Q)
 \end{equation*}
 is defined over the reflex field, which is $E$ (embedded into $\R$ via
 $\infty_0$). Over some finite abelian extension $E'/E$, it decomposes as
 a disjoint union of geometrically connected Shimura curves
 $\mathcal{X}_\Gamma$ ($\mathfrak{h}/\Gamma$ over $\C\;$; one possible
 $\Gamma$ is the image of $G(\Q)\cap C$ in $\G^{\mathrm{ad}}(\Q)$).

 Let $v=v_0$ be a {\it critical} finite place of $E$, \ie such that
 $B_v$ is a non-split quaternion algebra. Let $p$ be the residue
 characteristic of $v_0$ and let $v_1,\ldots, v_m$ be the other places
 of $E$ above $p$.

 Let $K$ be a totally imaginary quadratic extension of $E$, such that
 every place $v_i$ splits in $K$. We denote by $w_i, \ovl{w_i}$ the two
 places of $K$ above $v_i$. The quaternion algebra $B^\bullet =
 B\otimes_E K$ is ramified at $w=w_0$: $B^\bullet_w\cong B_v$.

 We fix an extension $\infty_0: K\hookrightarrow \C$ of the real
 embedding $\infty_0$ of $E$, and a double embedding $\C \hookleftarrow
 \ovl{\Q} \hookrightarrow \C_p$ whose restriction to $K$ is compatible
 with $(\infty_0,w)$. The embeddings $K \hookrightarrow \C_p$ which
 factor through some $w_i$ but not through some $\ovl{w_i}$ form a CM
 type $\Psi$: $\Psi\coprod \ovl{\Psi}=\Hom(K,\C_p)$ (identified with
 $\Hom(K,\C)$ through our chosen double embedding).

There is a natural $\Q$-group $G^\bullet$ such that
 $G^\bullet(\Q)=\{b\in B^\bullet| bb^\ast\in E\}$, acting on
 $V^\bullet=V\otimes_E K$. One can attach to $C$ an open compact
 subgroup $C^\bullet$ of $G^\bullet(\A_f)$, maximal at $p$, and a
 Shimura curve $\Sh^\bullet$ such that $\Sh$ is open and closed in
 $\Sh^\bullet$ (\cite[3.11]{tpuosc}). Alternatively, one can attach to
 $\Gamma$ a congruence subgroup $\Gamma^\bullet$ in $G^{\bullet ad}(\Q)$
 and a corresponding connected Shimura curve
 $\mathcal{X}_{\Gamma^\bullet}$ which is a twisted form of
 $\mathcal{X}_{\Gamma}$ (defined over some abelian extension of
 $K$). The point is that this new Shimura curve $\Sh^\bullet$ is of PEL
 type. More precisely, for $C$ small enough, it is a moduli space for
 $\ast$-polarized abelian varieties of dimension $g=4[E:\Q]$ with
 multiplication by
 $\mathcal{B}^\bullet=\mathcal{B}\otimes_{\mathcal{O}_E}\mathcal{O}_K$,
 Shimura type $\Psi$, and level structure. As in \ref{pro:7.1.5}, it
 admits a model $\Sh^\bullet$ over $\mathcal{O}_{E_v}$.

 Let $\udl{A} \rightarrow \Sh^\bullet$ be the universal abelian
scheme up to prime to $p$ isogeny.  Due to the $\mathcal{O}_K$-action,
the $p$-divisible group of $\udl{A} $ splits as $
\udl{A}[p^\infty]\cong \prod_i \Lambda_{w_i}\times \Lambda_{\ovl{w_i}}$,
where $ \Lambda_{w_i}$ and $\Lambda_{\ovl{w_i}}$ are $p$-divisible groups
of height $4[E_{v_i}:\Q_p]$ with formal multiplication by
$E_{v_i}$. Moreover, $ \Lambda_{w}= \Lambda_{w_0}$ is a special formal
$\mathcal{B}_v$-module of height $4.[E_{v}:\Q_p]$ (defined as in 
\ref{exa:4.3.4}, 
replacing $\Q_p$ by $E_v$), and $ \Lambda_{\ovl{w}}$ is the Serre dual of
$\Lambda_{w}$; the other $ \Lambda_{w_i},\;i>0,$ are etale, and the $
\Lambda_{\ovl{w_i}}$ are their respective Serre duals. Hence the Newton
polygon of the Dieudonn\'e module of $\udl{A}[p^\infty]$ has slopes
$0,\frac{1}{2},1$ (only $\frac{1}{2}$ if $m=0$).

\vspace*{2mm}
\begin{figure}[h]
 \begin{picture}(120,70)(-20,-10)
  \put(-10,0){\line(1,0){110}}
  \put(0,-10){\line(0,1){70}}
  \bezier{10}(60,15)(60,7.5)(60,0)
  \bezier{30}(90,45)(90,22.5)(90,0)
  \put(15,5){\makebox(0,0)[b]{{\tiny $0$}}}
  \put(40,15){\makebox(0,0)[c]{{\tiny $1/2$}}}
  \put(68,32){\makebox(0,0)[c]{{\tiny $1$}}}
  \thicklines
  \put(0,0){\line(1,0){30}}
  \put(30,0){\line(2,1){30}}
  \put(60,15){\line(1,1){30}}
 \end{picture}
 \caption{}
\end{figure}

Let $\mathfrak{v}$ be the discrete valuation ring
$\widehat{\Z^{\mathrm{ur}}}.\mathcal{O}_{E_v}$ as in
\ref{para:7.1.2}. The moduli problem 
for $\mathcal{B}_v$-module of height $4.[E_{v}:\Q_p]$ is representable
by the $p$-adic formal scheme $\displaystyle \coprod_\Z
\;({}_{E_v}\hat\Omega)_{\mathfrak{v}}$, where ${}_{E_v}\hat\Omega$ is
the analogue for (and over) $\mathcal{O}_{E_v}$ of the formal scheme
$\hat\Omega$, with associated analytic space
$(_{E_v}\Omega)_{\C_p}=\P^1_{\C_p}\setminus \P^1(E_v)$, \cf \ref{para:6.3.5}.
\end{para}

\begin{para}\label{para:7.5.3}
 This is the modular way taken by Boutot and Zink to reprove (and
 strengthen) \v{C}erednik's theorem.  As was remarked by
 Varshavsky \cite[3.13]{varshavsky98:_shimur_ii}, one obtains a simpler adelic
 formulation if one works with a twisted version ${}^{\mathrm{ad}}\Sh$ of $\Sh$
 (corresponding to a different normalization of $h_0$ in \ref{para:1.2.4}); this
 does not change the adelic description of the associated complex space
 \begin{equation*}
  {}^{\mathrm{ad}}\Sh_\C^{\mathrm{an}}= C\backslash
   \bigl((\P_\C^1\setminus\P^1(\R)) \times G(\A_f)\bigr)/G(\Q),
 \end{equation*}
 but makes the analogy with the following description of the associated
 $p$-adic space more striking.

 Let $\ovl{B}$ be the quaternion algebra over $E$ obtained from $B$ by
 interchanging the invariants at $\infty_0, v$, and let $\ovl{G}$ be the
 $\Q$-algebraic group attached to $\ovl{B}^\times$. Let us identify
 $\ovl{G}(\A^p_f)$ with $G(\A^p_f)$ via an anti-isomorphism
 $\ovl{B}\otimes_E \A^p_f \rightarrow B\otimes_E \A^p_f ,$ and set
 $\ovl{C}=GL_2(E_v)C^p$ viewed as a subgroup of $\ovl{G}(\A_f)$.  Then
 \begin{align*}
  {}^{\mathrm{ad}}\Sh_{\C_p}^{\mathrm{an}} & =C^p \backslash
  \bigl((_{E_v}\Omega)_{\C_p}\times G(\A_f^p)\bigr)/{\ovl{G}}(\Q)\\
  & =\ovl{C}\backslash \bigl((\P^1_{\C_p}\setminus\P^1(E_v))\times
  \ovl{G}(\A_f)\bigr)/{\ovl{G}}(\Q)
 \end{align*}
 where ${\ovl{G}}(\Q)$ acts on $(_{E_v}\Omega)_{\C_p}$ through
 ${\ovl{G}}^{\mathrm{ad}}(\Q_p)\cong PGL_2(E_v)$.  

 Moreover, one can drop the assumption that $C$ is maximal at $p$, on
 replacing $(_{E_v}\Omega)_{\C_p}$ by some finite etale covering as in
 \ref{para:7.4.7}. \cite[3]{tpuosc}, \cite[5]{varshavsky98:_shimur_ii}.
\end{para}

\begin{para}\label{para:7.5.4}
 Rapoport, Zink and Boutot on one hand (by modular methods \`a la
 Drinfeld), and Varshavsky on the other hand (by group-theoretic methods
 \`a la \v{C}erednik), have pointed out similar global
 $p$-adic uniformizations in much more general cases.

 Let us only mention a remarkable example: Mumford's {\it fake
 projective plane} \cite{aaswkak9pq0}. This is a smooth projective complex
 surface $S$ with the same Betti numbers as $\mathbb{P}^2$, but not
 isomorphic to $\mathbb{P}^2$. Actually, it is a surface of general
 type. Mumford's construction is of diadic nature. It is a quotient of
 the two-dimensional generalization $\Omega^{(2)}$ of $\Omega$ by an
 explicit arithmetic group. This $\Omega^{(2)}$ is the complement of the
 lines defined over $\Q_2$ in the projective plane.

 Using the last theorem of \cite{psfpg} or \cite{varshavsky98:_shimur_ii}, one can see
 that $S$ is actually a Shimura surface of PEL type, parametrizing
 polarized abelian varieties of dimension $9$, with multiplication by an
 order in a simple $\Q$-algebra $B$ of dimension $18$ with center
 $\Q(\sqrt{-7})$. The invariants of $B$ at the two primes of
 $\Q(\sqrt{-7})$ above $2$ are $1/3$ and $2/3$ (\cf a recent
 preprint by F. Kato for more detail). The associated group $G$ is a
 group of unitary similitudes of signature $(1,2)$.

 The ``diadic ball'' $\Omega^{(2)}$ (or rather its formal avatar) has an
 modular interpretation for certain $p$-divisible groups (with
 multiplication by $B_2$). These $p$-divisible groups split as a product
 of three $p$-divisible groups of slope $1/3$ and height $3$, and three
 $p$-divisible groups of slope $2/3$ and height $3$.

 The associated Gauss-Manin connection splits into six factors of rank
 $3$ (in two variables). It resembles the Appell-Lauricella connection
 studied in \cite{pdrefapdfhdpv}, \cite{calip}, \cite{fde}, but is different.
\end{para}

\subsection{Application of the theory of $p$-adic Betti lattices.}
\label{sub:7.6}

\begin{para}\label{para:7.6.1}
 We consider more closely the Gauss-Manin connection
 $(\mathcal{H},\nabla_{\mathrm{GM}})$ in the case of Shimura curves.

 We begin with the case $E=\Q$ as in \ref{sub:7.4}.  Then $\mathcal{H}$
 is a vector bundle of rank $4$, and the $B$-action on $\mathcal{H}$
 obtained by functoriality of De Rham cohomology commutes with
 $\nabla_{\mathrm{GM}}$.

 Let us choose an imaginary quadratic field $F=\Q(\sqrt{-d})$ inside
 $B$. Then $(\mathcal{H},\nabla_{\mathrm{GM}})\otimes F$ splits into two parts:
 $(\mathcal{H},\nabla_{\mathrm{GM}})_+\oplus
 (\mathcal{H},\nabla_{\mathrm{GM}})_{-}$ (on 
 which $F\subset B$ acts through identity and complex conjugation
 respectively).

 On the other hand, let us fix embeddings $F\subset \ovl{\Q} \subset \C$,
 and let us fix a base point $s\in \mathcal{X}_{\Gamma}(\ovl{\Q})$ of our
 Shimura curve.  This is the moduli point of a decorated fake elliptic
 curve $A/\ovl{\Q}$. We have $V^\C=\mathrm{H}^1_B(A_\C,\C)\cong
 (\mathcal{H}_\C^{\mathrm{an}})^{\nabla_{\mathrm{GM}}}$, and this again splits
 into $\pm$ parts. The Betti cohomology $F$-space
 $\mathrm{H}^1_B(A_\C,F)$ already splits into two parts:
 $\mathrm{H}^1_B(A_\C,F)_+\oplus \mathrm{H}^1_B(A_\C,F)_{-}$, both of
 rank two and stable under the action of $B\otimes F \cong M_2(F)$.

 The flag space $\mathcal{D}^\vee$ is actually defined over the reflex
 field $\Q$. Over $F$, it becomes isomorphic to $\P^1$, or more
 precisely, to $\P(\mathrm{H}^1_B(A_\C,F)_+)$. The period mapping describes the
 ``slope'' of the Hodge line $F^1\cap V^\C_+$ with respect to a basis of
 $\mathrm{H}^1_B(A_\C,F)_+\subset
 ((\mathcal{H}_+)_\C^{\mathrm{an}})^{\nabla_{\mathrm{GM}}}$, and is 
 thus given by a quotient of solutions of the (partial) Gauss-Manin
 connection $\nabla_{\mathrm{GM}+}$.

 Assume for simplicity that $\Gamma\subset B^\times/\Q^\times$ is
 torsion-free. Then $\Gamma$ may be identified with the projective
 monodromy group, \ie with the quotient by the homotheties of the image
 of the monodromy representation
 \begin{equation*}
  \pi_1((\mathcal{X}_{\Gamma })_\C,s)\rightarrow GL(\mathrm{H}^1_B(A_\C,F)_+).
 \end{equation*}

\end{para}

\begin{para}\label{para:7.6.2} Let us now turn to the $p$-adic side. We fix an
embedding $\ovl{\Q} 
\subset \C_p$. We shall limit ourselves to the case of a prime $p$ of
 {\it supersingular reduction} for $A$. Note that this includes the case
 of critical primes (\ref{lem:7.4.6}).

 Let $D$ be the definite quaternion algebra over $\Q$ which is ramified
 only at $p$, and let $\ovl{B}$ be the definite quaternion algebra over
 $\Q$ which is ramified at the same primes as $B$ except $p$. Then,
 whether $p$ is critical or not, the group $J$ (of self-quasi-isogenies
 of $\ovl{A}[p^\infty]$ commuting with $B$) is $\ovl{B}_p^\times$.

 By the local criterion for splitting fields \cite[III.3.5]{adadq}, $F$
 is also a splitting field for the quaternion algebra $D$ (and $\ovl{B}$
 as well). We can therefore apply the theory of
 I.\ref{sub-p-adic-betti-lattice-supersing}, and get an 
 $F$-structure $\mathrm{H}^1_B(A_{\C_p},F)$ inside
 $\mathrm{H}^1_{\mathrm{cris}}(\ovl{A},/\C_p)\cong
 V^{\C_p}=\D(\ovl{\Lambda})\otimes \C_p.$ By construction, it is stable
 under the action of $M_2(D)\otimes F\cong M_4(F)\cong B\otimes \ovl{B}
 \otimes F$.

 Using the embedding $F\otimes 1\otimes F\subset B\otimes \ovl{B} \otimes
 F$, one splits this $F$-structure $\mathrm{H}^1_B(A_{\C_p},F)$ into two parts
 $\mathrm{H}^1_B(A_{\C_p},F)_+\oplus \mathrm{H}^1_B(A_{\C_p},F)_{-}$ (on
 which $F\otimes 
 1\otimes 1\subset B\otimes \ovl{B} \otimes F$ acts through identity and
 complex conjugation respectively), both of rank two and stable under
 the action of $\ovl{B}\otimes F \subset M_4(F)$, hence under $J$.  \par
 The $p$-adic flag space $\mathcal{D}^\vee$ is again defined over the
 reflex field $\Q$. Over $F$, it becomes isomorphic to
 $\P(\mathrm{H}^1_B(A_{\C_p},F)_+)$ (note that the scaling ambiguity of
 $\mathrm{H}^1_B(A_{\C_p},F)_+$ by a factor in $\sqrt{F^\times}$ disappears
 here). The $J$-equivariant period mapping $\mathcal{P}$ describes the
 ``slope'' of the Hodge line $F^1\cap V^{\C_p}_+$ with respect to a
 fixed basis $(e_0,\,e_1)$ of $\mathrm{H}^1_B(A_{\C_p},F)_+\subset
 ((\mathcal{H}_+)_{\C_p}^{\mathrm{an}})^{\nabla_{\mathrm{GM}}}$, and is
 thus given by a quotient of solutions of $\nabla_{\mathrm{GM}+}$.

 Note that if $B$ is split at $p$, then $\mathcal{P}$ is of
 Gross-Hopkins type, and one has only a local $p$-adic uniformization of
 $\mathcal{X}_{\Gamma}$ around $s$ (in a supersingular disk).

 If instead $p$ is critical, then $\mathcal{P}$ is of Drinfeld type, and
 the group $\Gamma_p \subset J^1$ appearing in the global $p$-adic
 uniformization of $\mathcal{X}_{\Gamma}$ may be identified with the
 quotient by the homotheties of the image of the monodromy representation
 (in the sense of \ref{cor:7.2.6})
 \begin{equation*}
  \pi_1^{\mathrm{top}}((\mathcal{X}_{\Gamma })_{\C_p},s)\rightarrow
   GL(\mathrm{H}^1_B(A_{\C_p},F)_+).
 \end{equation*}

 A slight refinement provides a rank two
 $\mathcal{O}_F[\frac{1}{p}]$-lattice
 $\mathrm{H}^1_B(A_{\C_p},\mathcal{O}_F[\frac{1}{p}])_+$ in
 $((\mathcal{H}_+)_{\C_p}^{\mathrm{an}})^{\nabla_{\mathrm{GM}}}$ stable
 under $p$-adic monodromy.
\end{para}

\begin{para}\label{para:7.6.3}
 One can go one step further and see that the phenomenon encountered in
 \ref{para:6.3.3} shows up here again: ``the'' quotient $\tau$ of solutions of
 solutions of ${\nabla_{\mathrm{GM}}}$ which expresses locally the
 period mapping 
 $\mathcal{P}$ is essentially given by the same formula in the complex
 and the $p$-adic cases.

More precisely, assume that $A$ has complex multiplication by an order
in $M_2(F)$ (so that $A$ is isogenous to the square of an elliptic curve
with complex multiplication by $\mathcal{O}_{\Q(\sqrt{-d})}$ with
[fundamental] discriminant $-d$ ), and identify $\End_B(A)$ with $F$. We
denote by $\epsilon$ the Dirichlet character, $w$ the number of roots of
unity, $h$ the class number attached to $\mathcal{O}_{\Q(\sqrt{-d})}$.

\noindent We may choose our symplectic basis $e_0,e_1$ of
$\mathrm{H}^1_B(A_{\C},F)_+$ (\resp $\mathrm{H}^1_B(A_{\C_p},F)_+$) in
 such a way that 
$[\sqrt{-d}]^\ast e_0=-\sqrt{-d} \; e_0,\,$ $[\sqrt{-d}]^\ast
e_1=\sqrt{-d} \;e_1$ (here $[\sqrt{-d}]$ is viewed as an element of
$\End_B(A)=F$).

\noindent Let $\omega$ be a section of $\mathcal{H}_+$ in a Zariski
neighborhood $U$ of $s$, which corresponds to a relative differential of
the first kind (no pole) on the universal abelian scheme. Then $\omega$
is a cyclic vector for $\mathcal{H}_+$ over $U$ with respect to
$\nabla_{\mathrm{GM}}$, so that $\nabla_{\mathrm{GM}}$ amounts to a
 concrete linear 
differential equation of order two over $U$.

\noindent We calculate $\tau$ with respect the basis $e_0, e_1$
 of $\mathrm{H}^1_B(A_{\C},F)_+$ (\resp $\mathrm{H}^1_B(A_{\C_p},F)_+$)
 in a small 
 neighborhood of the base point $s$. Denoting by $z$ an algebraic local
 parameter at $s$ in $U$, and extending $e_0,\, e_1$ by horizontality, we
 then have the proportionality
\begin{equation*}
 \omega \sim \tau(z) e_0 + e_1
\end{equation*}
\end{para}

\begin{thm}\label{thm:7.6.4}
 There is a quotient $y=y_1/y_2$ of solutions of $\nabla_{\mathrm{GM}}$ in
 $\ovl{\Q}[[z]]$ such that $\tau$ takes the form
 $\;\tau(z)=\kappa^{-1}\;y(z)\;$ for a suitable constant $\kappa$.

\noindent In the complex case, $\;\kappa\sim\displaystyle
 \prod_{u\in (\Z/d)^\times}
 \Bigl(\Gamma\Bigl\langle\frac{u}{d}\Bigr\rangle\Bigr)^{\epsilon(u)w/2h}\;\in
 \C^\times/\ovl{\Q}^\times$.

\noindent In the $p$-adic case, assuming moreover that $p\nmid d $, 
\begin{equation*}
 \kappa\sim \prod_{u\in (\Z/d)^\times}
  \Bigl(\Gamma_p\Bigl\langle\frac{pu}{d}\Bigr\rangle\Bigr)^{-\epsilon(u)w/4h}
  \in \C_p^\times/\ovl{\Q}^{\times}.
\end{equation*}
\end{thm}
\begin{proof}
 Let us write $\omega_{| z=0}=\omega_{11}e_1$. Then
 $\tau(z)=\kappa^{-1} y(z)$, with $\kappa = \omega_{11}^2/2\pi i$ in
 the complex case, $\kappa = \omega_{11}^2$ in the $p$-adic case (the
 factor $2\pi i$ which arises in the complex case is the factor of
 proportionality of the symplectic forms in De Rham and in Betti
 cohomology respectively). The evaluation of $\omega_{11}$ in the
 complex case is given by the Lerch-Chowla-Selberg formula (\cf
 I.\ref{thm-ogus-formula}). In the $p$-adic case, it is given by
 I.\ref{thm-cm-period}.
\end{proof}

\begin{para}\label{para:7.6.5}
 This discussion generalizes without much complication to the case of
 Shimura curves attached to a quaternion algebra $D$ over a totally real
 number field $E$ of degree $>1$. We concentrate on the case of a
 critical place $v$ and assume, for simplicity, that $E$ is Galois over
 $\Q$ and that $v$ is the unique place of $E$ with the same residue
 characteristic $p$ (\ie $m=0$ in the notation of \ref{para:7.5.2}).

 Let $F$ be a totally imaginary quadratic extension of $E$ contained in
 $B$. Then $F$ is a splitting field for $B,\ovl{B}$, and $D$. Let us fix
 a double embedding $\C \hookleftarrow \ovl{\Q} \hookrightarrow \C_p$ as
 in \ref{para:7.5.2}, and an extension $FK\hookrightarrow \ovl{\Q}$ of $\infty_0$.
 For $C$ small enough, $\Sh^\bullet$ carries a universal abelian scheme
 of relative dimension $g=[FK:\Q]$. After tensoring $\otimes_K
 \ovl{\Q}$, its Gauss-Manin connection splits into pieces of rank two,
 indexed by the embeddings of $FK$ into $\ovl{\Q}$. We select the piece
 $(\mathcal{H},\nabla_{\mathrm{GM}})_+$ corresponding to $\infty_0$.

 Let us fix a base point $s\in \Sh(\ovl{\Q})\subset \Sh^\bullet(\ovl{\Q})
 $ of our Shimura curve.  This is the moduli point of a decorated
 abelian variety $A^\bullet/\ovl{\Q}$ (of dimension $g=[FK:\Q]$ and
 multiplication by $\mathcal{B}^\bullet$). We have $V^\bullet=V\otimes_E
 K \cong \mathrm{H}^1_B(A^\bullet_\C,\Q)$.  Note that
 $\mathrm{H}^1_B(A_\C,F)\cong 
 V\otimes_\Q FK$ splits into pieces of $FK$-rank two, indexed by the
 embeddings of $F$ into $\ovl{\Q}$. We select the piece
 $\mathrm{H}^1_B(A^\bullet_\C,F)_+$ corresponding to $\infty_0$.  We
 then have a 
 natural embedding $\mathrm{H}^1_B(A^\bullet_\C,F)_+\subset
 ((\mathcal{H}_+)_{\C}^{\mathrm{an}})^{\nabla_{\mathrm{GM}}}$, and
 $\mathrm{H}^1_B(A^\bullet_\C,F)_+$ is stable under monodromy.

 On the other hand, it follows from \ref{lem:7.4.6}, that $A^\bullet$ has
 supersingular reduction $\ovl{A}^\bullet$ at $v$. We can construct the
 $F$-structure $\mathrm{H}^1_B(A^\bullet_{\C_p},F)$ inside
 $\mathrm{H}^1_{\mathrm{cris}}(\ovl{A}^\bullet,/\C_p)\cong
 V^{\C_p}=\D(\ovl{\Lambda}^\bullet)\otimes \C_p.$ 
 By construction, it is stable under the action of $M_g(D)\otimes F\cong
 M_{2g}(F)$, which contains naturally $B\otimes_\Q \ovl{B}\otimes_\Q F$.

 Using the embedding $F\otimes 1\otimes F\subset B\otimes \ovl{B} \otimes
 F$, one splits this $F$-structure $\mathrm{H}^1_B(A_{\C_p},F)$ into parts of
 $FK$-rank two, indexed by the embeddings of $F$ into $\ovl{\Q}$. We
 select the piece $\mathrm{H}^1_B(A^\bullet_{\C_p},F)_+$ corresponding to
 $\infty_0$.  We then have a natural embedding
 $\mathrm{H}^1_B(A^\bullet_{\C_p},F)_+\subset
 ((\mathcal{H}_+)_{\C_p}^{\mathrm{an}})^{\nabla_{\mathrm{GM}}}$, and
 $\mathrm{H}^1_B(A^\bullet_{\C_p},F)_+$ is stable under $p$-adic monodromy.
\end{para}

\subsection{Conclusion.}\label{sub:7.7}

 Roughly speaking, the period mapping $\mathcal{P}$ attaches to
each member of an algebraic family of algebraic complex varieties (with
additional structure) a point in some flag space, which encodes the
Hodge filtration in the cohomology of this variety. This mapping can be
expressed in terms of quotients of solutions of the fuchsian
differential equation (Picard-Fuchs/Gauss-Manin) which controls the
variation of the cohomology.

We have described the theory of period mappings for $p$-divisible groups
in somewhat similar terms.

\bigskip
 Sometimes, it is possible to go through the looking-glass of sheer
 analogies. This occurs when we 
 restrict our attention to $p$-divisible groups attached to polarized
 abelian varieties with prescribed endomorphisms and fixed ``Shimura
 type''. In this situation, we have three algebraic objects defined over
 a number field (the reflex field):

\begin{enumerate}
 \renewcommand{\theenumi}{\alph{enumi}}
 \item the Shimura variety $\Sh$ parametrizing certain
       ``decorated'' abelian varieties,
 \item the Gauss-Manin connection $\nabla_{\mathrm{GM}}$ describing
       the variation of these abelian varieties over $\Sh$, and
 \item the flag variety $\mathcal{D}^\vee$.
\end{enumerate}

 These objects are transcendentally related, both in the complex and
 $p$-adic sense.

 The period mapping may be viewed as
 multivalued and defined on $p$-adic open domains in $\Sh$; it is given
 as in the complex case by quotients of solutions of the Gauss-Manin
 connection (as a $p$-adic connection).

 In most situations, these open domains are small and $\mathcal{P}$ is
 single-valued.

 However, it happens in some remarkable cases that $\mathcal{P}$ is a
 global multivalued function. This occurs when the abelian varieties
 under consideration form a {\it single isogeny class modulo $p$}. In
 this situation, the complex and $p$-adic transcendental relations
 between $\Sh, \nabla_{\mathrm{GM}}$ are completely similar: $\Sh$ is
 uniformized 
 by an analytic space which is etale over some open subset $\mathcal{D}$
 of $\mathcal{D}^{\vee}$, and the group of deck transformations $\Gamma$
 is isomorphic to the projective global monodromy group of
 $\nabla_{\mathrm{GM}}$.  In III, we shall analyze in detail the
 group-theoretic 
 aspects of this situation.

\chapter{$p$-adic orbifolds and monodromy.}
\addtocontents{toc}{\protect\par\vskip2mm\hskip20mm(Group theory)\par\vskip5mm}
\minitoc
\newpage

\section{\'Etale coverings and fundamental groups in $p$-adic geometry.}
\label{sec::1}

\markboth{\thechapter. $p$-ADIC ORBIFOLDS AND MONODROMY.}%
{\thesection. \'ETALE COVERINGS AND FUNDAMENTAL GROUPS.}

\begin{abst}
 In this section, inspired by De Jong's work \cite{efgonas}, we discuss various
 kinds of finite or infinite \'etale coverings in non-archimedean
 analytic geometry, and the associated fundamental groups.
 The temperate fundamental group of a $p$-adic manifold encapsulates
 information on the bad reduction properties of all its finite etale
 coverings.
\end{abst}

\subsection{More on the topological fundamental group.}\label{sub::1.1} 

\begin{para}\label{para::1.1.1}
 We shall work throughout in the context of Berkovich analytic geometry,
 over a fixed complete ultrametric field $(K,|\;|)$.

 What has been said about Berkovich spaces in I.\ref{sec-1} remains true
 in this broader context (\ie without any further restriction on
 $K$). 

 More precisely, we shall work with {\it paracompact strictly
 $K$-analytic spaces} in the sense of \cite{ecfnas}. We recall that these
 form a category equivalent to the category of quasi-separated rigid
 spaces over $K$ having an admissible covering of finite type.

 The main
 advantage of these spaces over rigid spaces lies in their topology: this
 is a topology in the usual sense (not a Grothendieck topology);
 moreover, each point has a fundamental system of open neighborhoods
 which are locally compact, countable at infinity and arcwise connected.

 We shall deal mostly with what we have called, for convenience,
 (analytic) {\it K-manifolds}, \ie 
 paracompact strictly $K$-analytic spaces which satisfy the following
 condition (I.\ref{para:p-adic-manifold}):
 
\begin{quote}
 {\it any $s\in S$ has a neighborhood $U(s)$ which is isomorphic 
 to an affinoid subdomain of some smooth space (in the sense of \cite{ecfnas}
 \ie a space which admits locally an \'etale morphism to the affine space 
  $\A^{\dim S}$).}
\end{quote}

\begin{rems}\label{rems:1.1.1.5}
 a) Any affinoid domain of a $K$-manifold is
 rigid-analytically smooth, \ie satisfies the jacobian criterion (\cf
 \cite[2.2]{ecfnas}; note also that ``smooth'' in the sense of Berkovich
 amounts to ``rigid-analytically smooth and having no boundary''); a
 fortiori, an affinoid $K$-manifold is rigid-analytically
 smooth.  Conversely, a rigid-analytically smooth (strictly) $K$-affinoid
 space $S$ is a $K$-manifold, \cf \cite[rem. 9.7]{spasalc}.

 b) This notion is slightly stronger than the notion of spaces
 locally embeddable in a smooth space considered in \cite{spasalc}. The
 main difference is that we impose that $S$ is a {\it good space} in the
 sense of \cite{ecfnas}, \ie that any point has an affinoid
 neighborhood, which is not automatic for non-classical points (this
 will be useful in \secref{sec::3}).

 c) It follows immediately from remark a) that {\it a
 $K$-manifold is nothing but a good and rigid-analytically smooth
 (paracompact strictly $K$-analytic) space.}

 d) Let $\phi:S'\to S$ be an \'etale morphism of paracompact
 strictly $K$-analytic spaces.  If $S$ is a $K$-manifold, so is $S'$.

 e) For an algebraic $K$-variety $X$, one has $X$ smooth
 $\iff$ $X^{\mathrm{an}}$ smooth $\iff$  $X^{\mathrm{an}}$ is a
 $K$-manifold.  A morphism $f$ is \'etale $\iff$ $f^{an}$ is \'etale.
\end{rems}

 \medskip 
 It follows from \cite[cor. 9.5]{spasalc} that $K$-manifolds are locally
 arcwise-connected and {\it locally simply-connected},
 hence subject to the familiar theory of topological coverings (the
 difficulty is to show that non-classical points admit contractible 
 neighborhoods).   
  In particular, any pointed $K$-manifold $(S,s)$ admits a {\it universal
 (pointed) topological covering} $(\tilde S,\tilde s)$, such that
 $S \cong \tilde S/\pi_1^{\mathrm{top}}(S,s)$; 
 $\tilde S$ is a (pointed) $K$-manifold, and  $\pi_1^{\mathrm{top}}(S,s)$
 is a discrete group of automorphisms of $\tilde S$.

 Let us also recall that in dimension one, $\pi_1^{\mathrm{top}}(S,s)$
 is a free group isomorphic to the fundamental group of the dual graph
 of the semistable reduction of $S$ (\cf \cite[5.3]{efgonas}). In higher
 dimension, there is no such restriction on the discrete groups
 $\pi_1^{\mathrm{top}}(S,s)$: actually, any differentiable manifold is
 homotopy equivalent to some $K$-manifold \cite[6.1.8]{staagonf}.
\end{para}

\begin{para}\label{para::1.1.2}
 Berkovich's analysis of the homotopy of analytic spaces of higher
 dimension makes use of De Jong's alteration theorem. This allows him to
 reduce the problem to the case of the generic fiber of so-called
 polystable fibrations of formal schemes over $K^0$ (the ring of
 integers in $K$). He associates to such fibrations a simplicial set
 which encodes the combinatorics of incidence of the intersection strata
 in the special fiber, and shows that the geometric realization of this
 simplicial set is homotopy equivalent to the generic fiber (viewed as
 an analytic $K$-manifold) \cite[5.4]{spasalc}.

 In particular, in the case of (the generic fiber of) an algebraic
 polystable fibration, the 
 simplicial set is finite, hence the topological fundamental group is
 finitely generated. Using De Jong's theorem and the following result,
 one can deduce that 
 {\it the topological fundamental group of an algebraic $K$-manifold is
 finitely generated}.
\end{para}

\begin{pro}\label{pro::1.1.3}
 Let $\ovl{S}$ be a $K$-manifold, and let $Z$ be a Zariski-closed,
 nowhere dense, reduced analytic subset. Then any topological covering
 of $S:=\ovl{S}\setminus Z$ extends uniquely to a topological covering of
 $\ovl{S}$. Equivalently (since $S$ and $\ovl{S}$ are locally
 simply-connected), the natural homomorphism $\pi_1^{\mathrm{top}}(S,s)\to
 \pi_1^{\mathrm{top}}(\ovl{S},s)$ is an isomorphism.
\end{pro}
\begin{proof}
 Extensions of topological coverings from $S$ to $\ovl{S}$ are unique if
 they exist, since $\ovl{S}$ is normal (note that $S$ is connected, 
 like $\ovl{S}$, by normality of $\ovl{S}$ \cite[3.3.16]{staagonf}).
 For the existence, we argue by induction on the dimension of $Z$
 (starting from the trivial case where $Z$ is empty).  This reduces the
 question to extending topological coverings from $S$ to $\ovl{S}\setminus
 Z^{\mathrm{sing}}$ (where $Z^{\mathrm{sing}}$ denotes the non-smooth
 locus of $Z$, which 
 is closed and of dimension less than the dimension of $Z$, \cf 
 \cite[2.2]{ecfnas}); equivalently, this reduces the question to the case
 where $Z$ is (rigid-analytically) smooth.

 By unicity, the question of existence of extensions is local, \ie is a
 question at the neighborhood in $S$ of each point $z$ of $Z$; we may
 assume that $\ovl{S}$ and $Z$ are (strictly) affinoid,
 rigid-analytically smooth and equidimensional ($Z$ of codimension $r>0$
 in $\ovl{S}$). The statement is clear if $z$ has a neighborhood $U(z)$
 which is a finite union of affinoid domains $U_i$'s such that
 $U_i\setminus (U_i\cap Z)$ is simply-connected. This is indeed the
 case, by Kiehl's theorem on the existence of tubular neighborhoods
 \cite{defeaidnf}: translated into Berkovich geometry, it asserts that (for
 $S$ affinoid and $Z$ a Zariski-closed subset, both rigid-analytically
 smooth and equidimensional) $Z$ has a neighborhood which is a finite
 union of affinoid domains $U_i$'s isomorphic to $(U_i \cap Z)\times
 \matheur{D}^r$, where $\matheur{D}^r$ is an $r$-dimensional polydisk,
 and $ (U_i\cap Z)$ corresponds to $(U_i \cap Z)\times \{0\}.$ Since
 $U_i\cap Z$ is a rigid-analytically smooth affinoid domain, it is
 locally simply connected according to Berkovich; on the other hand, in the
 non-archimedean situation, $\matheur{D}^r\setminus \{0\}$ is
 simply-connected. Whence the result.
\end{proof}

\subsection{\'Etale versus topological coverings (again).}\label{sub::1.2}  

\begin{para}\label{para::1.2.1}
 Let $S$ be a connected $K$-manifold. Let $S' \xrightarrow{f} S$ be a
 finite \'etale morphism ($S'$ is then a $K$-manifold).  Its degree $n$ may be
 defined as the rank of the locally free $\mathcal{O}_S$-module
 $f_\ast\mathcal{O}_{S'}$. Assume that $K$ is algebraically closed. Then
 one has the formula 

 \begin{equation*}
  n=\sum_{y\in f^{-1}(x)}[\mathscr{H}(y):\mathscr{H}(x)]
 \end{equation*}
 where $\mathscr{H}(x)$ stands for the residue field of
 the point $x$ (if $S=\matheur{M}(\mathcal{A})$ is affinoid and
 $x$ corresponds to the bounded multiplicative seminorm $|\;|_x$, this is
 just $\mathscr{H}(x)= 
 \widehat{Q(\mathcal{A}/\Ker|\;|_x)}$).

 Then $f$ is a {\it topological covering if and only if all fibers have
 the same cardinality}, which 
 amounts to the algebraic notion of being ``completely decomposed'': $
 \mathscr{H}(y)=\mathscr{H}(f(y))$ for any $y\in S'$, 
 \cf \cite[3.2.7]{staagonf}, \cite[6.3.1]{ecfnas}.

 This criterion shows at once that the Kummer \'etale covering
 $z\mapsto z^n$ (with $\mathrm{char}(K) \nmid n$) of the 
 punctured disk is {\it not} a topological covering. 
\end{para}

\begin{para}\label{para::1.2.2}   A {\it geometric point} $\ovl{s}$ of $S$ is a
 point with value in some complete algebraically closed extension 
 $(\Omega,|\;|)$ of
 $(K,|\;|)$,  \ie a morphism ${\ovl{s}} : \matheur{M}(\Omega) \to S$ in
 the category of analytic spaces over 
 $K$ \cite[p.48]{staagonf}.  For any finite \'etale covering there are
 exactly $n$ liftings ${\ovl{s}'}:\matheur{M}(\Omega)\to S'$ of $\ovl{s}$,
 \ie $n$ geometric points of $S'$ above $S$. We usually denote by $s$ the
 point of $S$ which is the image of $\ovl{s}$.  
\end{para}

\begin{dfn} [Berkovich, De Jong] \label{def::1.2.3}
 A morphism $f:S'\to S$ is a {\it covering} (\resp {\it \'etale 
 covering}, \resp {\it topological covering}) if $S$ is covered by open
 subsets $U$ such that $\coprod V_j = f^{-1}U \to U$ and the restriction 
 of $f$ to every $V_j$ is finite (\resp \'etale finite, \resp an
 isomorphism).
\end{dfn}

\noindent There is an obvious notion of morphism
 between coverings of $S$.  The category of \'etale coverings, \resp
 finite \'etale coverings, \resp topological coverings, will be denoted
 by $\Cov^{\mathrm{et}}_S$, \resp $\Cov^{\mathrm{alg}}_S$, \resp
 $\Cov^{\mathrm{top}}_S$.
 
 \medskip There are inclusions (fully faithful embeddings) 
 \begin{equation*}
  \Cov^{\mathrm{alg}}_S \hookrightarrow  \Cov^{\mathrm{et}}_S,
   \quad \Cov^{\mathrm{top}}_S \hookrightarrow  \Cov^{\mathrm{et}}_S.
 \end{equation*}

\begin{rems}
 (i) If $S'\to S$ is an \'etale covering, then $S'$ is a $K$-manifold.

 (ii) It is worth noticing that the notions of \'etale and of 
 topological coverings are local on the base $S$, but not the
 notion of algebraic coverings (except if $S$ is compact).

 (iii) Even if $S$ is covered by finitely many $U$'s as in the
 definition, it does not follow that $S'$ is a disjoint union of
 finite algebraic coverings of $S$. For instance, a Tate elliptic
 curve (with universal covering $\tilde S\cong \G_m$) has a
 covering by two open subsets $U_1,U_2$ such that $\tilde S_{|
 U_i}$ is a disjoint union of copies of $S$.

  (iv) It is an immediate but important fact that these
 categories are stable under taking fiber product (over $S$). 
 It is also stable under taking finite disjoint unions.

  (v) One should pay attention to the fact that unlike
 $\Cov^{\mathrm{top}}_S$, $\Cov^{\mathrm{et}}_S$ is not stable 
 under arbitrary disjoint unions, nor under composition (in the sense
 that the composition of an \'etale covering $S'/S$ and an \'etale
 covering $S''/S'$ may not be an \'etale covering
 $S''/S$). However, we have:
\end{rems}

\begin{lem}\label{lem::1.2.4}
 Any morphism which is the composition of a covering (\resp \'etale
 covering) followed or preceded by a finite (\resp finite \'etale)
 morphism is a covering (\resp \'etale covering).
\end{lem}
\begin{proof}
 Let us consider the case of a covering $g:S''\to S'$
 followed by a finite covering $f:S'\to S$ (the reverse situation is
 immediate from the definition). Let $s$ be a point of $S$. Then for any
 $s'$ in the finite set $f^{-1}(s)$, there exists an open neighborhood
 $U_{s'}\subset S'$ of $s'$ such that $g^{-1}(U_{s'})$ is a disjoint
 union of spaces $U_{s',i}$, each being finite over $U_{s'}$ via $g$. We
 may assume that the $U_{s'}$ are pairwise disjoint. Since $f$ is finite,
 it is in particular a closed continuous map, hence 
 $S'$ admits a basis of open neighborhoods of $f^{-1}(\bar{s})$ of
 the form $f^{-1}(V)$. 
 In particular,
 there is an open neighborhood $V_s$ of $s$ such that $\displaystyle
 \coprod_{s'\in f^{-1}(s)} U_{s'}$ contains $f^{-1}(V_s)$. We may
 replace each $U_{s'}$ by its intersection with $f^{-1}(V_s)$ (which is a
 union of connected components of $f^{-1}(V_s)$). It is then clear that
 $g\circ f$ induces a finite morphism from $U_{s',i}$ to $V_s$. Hence
 $g\circ f$ is a covering.  The \'etale variant is similar.
\end{proof}

\begin{exs}\label{exs::1.2.5}
  (i) The prototype of an \'etale covering is given by the logarithm map
	\begin{equation*}
	 \log:\; \matheur{D}_{\C_p}(1,1^-)\longrightarrow
	 (\A^1_{\C_p})^{\mathrm{an}}.
	\end{equation*}
	For any $m\in \N$, the inverse image of the open disk
	\begin{equation*}
	 \matheur{D}_m=\matheur{D}\left(0,
	 \bigl(|p|^{-m+\frac{1}{p-1}}\bigr)^-\right)
	 \subset (\A^1_{\C_p})^{\mathrm{an}}
	\end{equation*}
	is a disjoint union
	\begin{equation*}
	 \coprod_{\zeta \in \mu_{p^\infty}}\;\Bigl\{q\in
	  \matheur{D}(1,1^-) \Bigm| \zeta q^{p^m}\in 
	  \matheur{D}\left(1, \bigl(|p|^{\frac{1}{p-1}}\bigr)^-\right)\Bigr\}
	\end{equation*}
	and each of these components is isomorphic to $\matheur{D}_m$ if $m=0$,
	and finite \'etale onto $\matheur{D}_m$ in general (we are using
	here the 
	factorization $q\mapsto \zeta q^{p^m}=t\mapsto \log t \mapsto
	p^{-m}\log t =\log q$ of the logarithm).

  (ii) This example may be interpreted as a Dwork period
	mapping (II.\ref{para:6.3.1}). In fact, all examples of period
	mappings considered in II.\ref{sub:6.3} are \'etale coverings
	of the period domain 
	$\mathcal{D}$. For instance, the Gross-Hopkins period mapping  
	(II.\ref{para:6.3.2}) is an \'etale covering of $
	(\P^1_{\C_p})^{\mathrm{an}}$, \cf \cite[7.2]{efgonas}. 

	\noindent Any period mapping is \'etale, but it may happen in
	some cases that it is surjective without being an \'etale
	covering, \cf \cite[5.41]{psfpg} (this somewhat mysterious
	phenomenon is due to the existence of bijective
	non-quasi-compact \'etale morphisms).

  (iii) Any (possibly infinite) topological covering of a
	finite \'etale covering is an \'etale covering after \ref{lem::1.2.4}.
\end{exs}

The following two lemmas are due to V. Berkovich and J. De Jong (and not
stated in the most general form).
 
\begin{lem}\label{lem::1.2.6}
 \begin{enumerate}
  \item An \'etale presheaf $\mathcal{F}$ on $S$ which is representable
	by an \'etale covering $S'$ is a sheaf.  For any subdomain $U
	\subset S$, $\mathcal{F}_{|U}$ is representable by $S'\times_S
	U$.
  \item Let $\mathcal{F}$ be an \'etale sheaf. Let $U_i$ an
	open covering of $S$. If for every $i$, $\mathcal{F}_{|U_i}$ is
	representable by an \'etale covering, so is $\mathcal{F}$.
 \end{enumerate}
\end{lem}
 
 For (i), \cf \cite[4.1.3, 4.1.4, 4.1.5]{ecfnas}.  Here we take advantage
 of working with {\it strict} analytic spaces.  Gluing --- using 
 \cite[1.3.3]{ecfnas} --- gives (ii) (\cf also \cite[2.3]{efgonas}).

\begin{lem}\label{lem::1.2.7} Let $S'\to S$ be an \'etale covering, and
 let $R\subset S'\times_S S'$ be a union of connected components which
 is an equivalence relation on $S'$ over $S$. Then $S'/R$ (viewed as an
 \'etale sheaf on $S$) is representable by an \'etale covering $S''\to S$. 
\end{lem}

\cf \cite[2.4]{efgonas}: by the previous lemma, the question is local on
 $S$; the proof is by reduction to the case of a finite \'etale covering
 of affinoid spaces. 

\medskip In the situation of \ref{lem::1.2.7}, we say that the $S''\to
 S$ is a {\it quotient \'etale covering of} $S'\to S$.

\subsection{Existence of \'etale paths.}\label{sub::1.3}  
 
\begin{para}\label{para::1.3.1} Let us fix a geometric point $\ovl{s}$ of
 $S$. We consider the fiber functor 
 \begin{equation*}
  F^{\mathrm{et}}_{S,\ovl{s}}=F^{\mathrm{et}}_{\ovl{s}} :
   \Cov^{\mathrm{et}}_S \longrightarrow \Sets,\;\; S' \mapsto
   \{\text{geometric points $\ovl{s}'$ of $S'$ above $\ovl{s}$}\}
 \end{equation*}
and its restrictions

\begin{equation*}
 F^{\mathrm{alg}}_{S,\ovl{s}}=F^{\mathrm{alg}}_{\ovl{s}}:
 \Cov^{\mathrm{alg}}_S\longrightarrow \Sets,
\end{equation*}
(considered in \cite{sga1}) and
\begin{equation*}
 F^{\mathrm{top}}_{S,\ovl{s}}=F^{\mathrm{top}}_{\ovl{s}}:
  \Cov^{\mathrm{top}}_S\longrightarrow \Sets.
\end{equation*} 

 An {\it \'etale path} (\resp {\it algebraic path}) from
 $\ovl{s}$ to another geometric point $\ovl{t}$ of $S$ is an isomorphism
 of functors $F^{\mathrm{et}}_{\ovl{s}} \xrightarrow{\sim}
 F^{\mathrm{et}}_{\ovl{t}}$ (\resp $F^{\mathrm{alg}}_{\ovl{s}}
 \xrightarrow{\sim} F^{\mathrm{alg}}_{\ovl{t}}$). 

 The set of \'etale paths is topologized by taking as fundamental open
 neighborhoods of an \'etale path $\alpha$ the set 
 $\Stab_{S',\ovl{s}'}\circ \alpha$, where $\Stab_{S',\ovl{s}'}$ runs among
 the stabilizers in $\Aut(F_{\bar{s}}^{\mathrm{et}})$ of arbitrary geometric
 points $\ovl{s}'$ above 
 $\ovl{s}$ in arbitrary \'etale coverings $S'/S$.

 If we do the same for algebraic paths, we get the profinite topology
 considered in \cite{sga1}. 

 On the other hand, an isomorphism of functors
 $F^{\mathrm{top}}_{\ovl{s}} \xrightarrow{\sim}
 F^{\mathrm{top}}_{\ovl{t}}$ amounts to a path up to homotopy in the
 usual sense between $s$ and $t$ in the arcwise connected space $S$
 (this is compatible with composition of paths, which is juxtaposition
 in the reverse order according to our convention); they form a
 discrete set --- in fact, a principal homogeneous space under 
 $\pi^{\mathrm{top}}(S,s)$ acting on the left. 

 \medskip The inclusions $ \Cov^{\mathrm{alg}}_S \hookrightarrow
 \Cov^{\mathrm{et}}_S ,\;\Cov^{\mathrm{top}}_S \hookrightarrow
 \Cov^{\mathrm{et}}_S$ induce continuous maps from the space of \'etale
 paths between $\ovl{s}$ and $\ovl{t}$ to the space of algebraic paths
 between $\ovl{s}$ and $\ovl{t}$ (\resp to the discrete set of paths up
 to homotopy between $s$ and $t$).

 \medskip\noindent If $f: S\to X$ is any morphism of $K$-manifolds, and
 if $\alpha$ is an \'etale path from $\ovl{s}$ to $\ovl{t}$, there is an 
 obvious notion of push-forward $f_\ast\alpha$, which is an \'etale path
 from $f\circ \ovl{s}$ to $f\circ \ovl{t}$ ; it is compatible with
 composition of \'etale paths (as composition of isomorphisms of functors).
\end{para}

\begin{para}\label{para::1.3.2} It follows from \cite[exp. V]{sga1} that
 algebraic paths exist. Using the nice topology of Berkovich spaces,
 A.J. De Jong was able to prove \cite[2.9]{efgonas}: 
\end{para}

\begin{klem}\label{lem::1.3.3}  Etale paths exist. More precisely,
 any Berkovich path $\gamma$ (up to homotopy) between $s$ and 
 $t$ lifts to an \'etale path $\ovl{\gamma}$ between $\ovl{s}$ and $\ovl{t}$.
\end{klem}

\begin{proof}
 Up to introducing intermediate geometric points, we may assume
 that $S$ is affinoid (hence compact). For any {\it finite} open cover
 $\mathcal{U}=(U_i)$ of $S$, let us denote by
 $\Cov^{\mathrm{et}}_{S,\mathcal{U}}$ the full subcategory of
 $\Cov^{\mathrm{et}}_{S}$ of \'etale coverings $f:S'\to S$ such that
 $f^{-1}(U_i)= \coprod V_{ij}$ and $V_{ij}$ is \'etale finite over $U_i$
 via $f$.  We thus have
 \begin{equation*}
 \Cov^{\mathrm{et}}_{S}=\underset{\mathcal{U}}{\liminj}\;
  \Cov^{\mathrm{et}}_{S,\mathcal{U}}, \text{ and }
  \Isom(F^{\mathrm{et}}_{\ovl{s}}, F^{\mathrm{et}}_{\ovl{t}})
  =\underset{\mathcal{U}}{\limproj}\;
  \Isom(F^{\mathrm{et}}_{\ovl{s},\mathcal{U}},
  F^{\mathrm{et}}_{\ovl{t},\mathcal{U}})
 \end{equation*}

 \medskip\noindent where $F^{\mathrm{et}}_{\ovl{s},\mathcal{U}},
 F^{\mathrm{et}}_{\ovl{t},\mathcal{U}}$ denote the restrictions of
 $F^{\mathrm{et}}_{\ovl{s} }, F^{\mathrm{et}}_{\ovl{t} }$ to
 $\Cov^{\mathrm{et}}_{S,\mathcal{U}}$.

 We have to show that $\Isom(F^{\mathrm{et}}_{\ovl{s}},
 F^{\mathrm{et}}_{\ovl{t}}) \neq \emptyset$. The idea is to construct 
 non-empty compact subsets $K_\mathcal{U}$ of
 $\Isom(F^{\mathrm{et}}_{\ovl{s},\mathcal{U}},
 F^{\mathrm{et}}_{\ovl{t},\mathcal{U}})$ (topologized as above) such that
 if $\mathcal{U}'$ refines $\mathcal{U}$, $K_{\mathcal{U}'}$ maps to
 $K_\mathcal{U} $ under the natural map
 $\Isom(F^{\mathrm{et}}_{\ovl{s},\mathcal{U}'},
 F^{\mathrm{et}}_{\ovl{t},\mathcal{U}'})\to
 \Isom(F^{\mathrm{et}}_{\ovl{s},\mathcal{U}},
 F^{\mathrm{et}}_{\ovl{t},\mathcal{U}})$; it will follow that
 $\Isom(F^{\mathrm{et}}_{\ovl{s}}, F^{\mathrm{et}}_{\ovl{t}}) \supset
 \underset{\mathcal{U}}{\limproj} K_\mathcal{U}\neq \emptyset$.

 In fact, it suffices to do so for a cofinal system of $\mathcal{U}$'s,
 which are chosen as follows. Since $S$ is {\it arcwise connected}, there
 is a continuous embedding $I=[0,1]\hookrightarrow S$ with $0\mapsto
 s,\;1\mapsto t$.
 One considers finite open coverings $\mathcal{U} =(U_1,\ldots, U_n)$
 such that
 \begin{equation*}
 U_1\cap I = [0,r_1[,\;\ldots, U_i\cap I=]r_i,t_i[,\;\ldots, U_m\cap
  I=]r_m,1],\;U_{m+1}\cap 
  I=\emptyset,\ldots
 \end{equation*}
 \noindent for some points $r_i<t_{i-1}<r_{i+1}<t_i$ in $I \;(i=2,\ldots,
 m-1)$ and some $ m\leq n$.

 Given $\mathcal{U}$, one chooses $s_i$ in $U_i\cap U_{i+1}\cap I$ (which
 amounts to a real number between $r_{i+1}$ and $t_i$), and a geometric
 point $\ovl{s}_i$ above $s_i\; (i=1,\ldots, m-1)$; one completes this
 collection of geometric points by setting $\ovl{s}_0=\ovl{s},\;
 \ovl{s}_m=\ovl{t}$.  Let us look at the diagram

 \begin{equation*}
 \Cov^{\mathrm{et}}_{S,\mathcal{U}}
 \xrightarrow{\sigma_i}
 \Ind \Cov^{\mathrm{alg}}_{\mathcal{U}_{i+1}}
 \begin{array}{c}
  \xrightarrow{F^{\mathrm{alg}}_{\ovl{s}_i, U_{i+1}}}\\
  \xrightarrow{F^{\mathrm{alg}}_{\ovl{s}_{i+1}, U_{i+1}}}\\
 \end{array}
 \Sets.
 \end{equation*}               
 One has $F^{\mathrm{alg}}_{\ovl{s}_j, U_{i+1}}\circ \sigma_i =
 F^{\mathrm{et}}_{\ovl{s}_i, \mathcal{U} },\;\;j=i,i+1$, whence
 continuous maps 
 \begin{equation*}
 \Isom(F^{\mathrm{alg}}_{\ovl{s}_i,
 U_{i+1}},F^{\mathrm{alg}}_{\ovl{s}_{i+1}, U_{i+1}})\to
 \Isom(F^{\mathrm{et}}_{\ovl{s}_i, \mathcal{U}
 },F^{\mathrm{et}}_{\ovl{s}_{i+1}, \mathcal{U}}).
 \end{equation*} 
 Since algebraic paths exist, $\Isom(F^{\mathrm{alg}}_{\ovl{s}_i,
 U_{i+1}},F^{\mathrm{alg}}_{\ovl{s}_{i+1}, U_{i+1}})$ is a non-empty
 compact set for any $i=0,\ldots, m$.  On composing such isomorphisms for
 $i=0,\ldots, m$, one gets a continuous map
 \begin{equation*}
 \sigma_{\mathcal{U},(\ovl{s}_i)}: \prod_{i=0}^m
  \Isom(F^{\mathrm{alg}}_{\ovl{s}_{i},
  U_{i+1}},F^{\mathrm{alg}}_{\ovl{s}_{i+1}, U_{i+1}})\to
  \Isom(F^{\mathrm{et}}_{\ovl{s},\mathcal{U}},
  F^{\mathrm{et}}_{\ovl{t},\mathcal{U}})
 \end{equation*} 
 whose image is a non-empty compact set $K_\mathcal{U}$.

 In fact, while the map $\sigma_{\mathcal{U},(\ovl{s}_i)}$ depends on the
 choice of the geometric points $\ovl{s}_i$, its image $K_\mathcal{U}$
 does not. Indeed, for any other choice $\ovl{s}'_i$, there is an
 algebraic path
 \begin{equation*}
  \gamma_i:\;F^{\mathrm{alg}}_{\ovl{s}_i, U_{i+1}\cap U_i} \xrightarrow{\sim}
   F^{\mathrm{alg}}_{\ovl{s}'_i, U_{i+1}\cap
   U_i},\;\;\gamma_0=\id,\;\gamma_m=\id.
 \end{equation*} 
 We denote by the same letter the induced algebraic
 paths on $U_i$ and on $U_{i+1}$. Then 
 $\sigma_{\mathcal{U},(\ovl{s}_i)}(\alpha_0,\ldots,\alpha_m) 
 =
 \sigma_{\mathcal{U},(\ovl{s}'_i)}(\gamma_1\alpha_0\gamma_0^{-1},\ldots,
 \gamma_{m}\alpha_m\gamma_{m-1}^{-1})$, whence the independence.

 From this observation, it is straightforward to conclude that if
 $\mathcal{U}'$ refines $\mathcal{U}$, $K_{\mathcal{U}'}$ maps to 
 $K_\mathcal{U} $ under the natural map
 $\Isom(F^{\mathrm{et}}_{\ovl{s},\mathcal{U}'},
 F^{\mathrm{et}}_{\ovl{t},\mathcal{U}'})\to
 \Isom(F^{\mathrm{et}}_{\ovl{s},\mathcal{U}},
 F^{\mathrm{et}}_{\ovl{t},\mathcal{U}})$. This establishes the first
 assertion.

 For the second one (the surjectivity of
 $\Isom(F_{\ovl{s}}^{\mathrm{et}}), F_{\ovl{t}}^{\mathrm{et}}) \to
 \Isom (F_{\ovl{s}}^{\mathrm{top}}, F_{\ovl{t}}^{\mathrm{top}})$), 
 one may argue as follows: let $\ovl{\tilde s}$ be a
 lifting of $\ovl{s}$ to the universal topological covering $\pi:\tilde
 S \to S$. Then $\gamma.\bar{\tilde s}$ is a lifting of $\ovl{t}$
 to the universal topological covering.  Let $\alpha$ be an \'etale path
 between $\bar{\tilde s}$ and $\gamma.\bar{\tilde s}$. Then
 $\ovl{\gamma}:=\pi_\ast(\alpha)$ is an \'etale path between $\ovl{s}$
 and $\ovl{t}$ which maps to the ordinary path $\gamma$ up to homotopy,
 viewed as an isomorphism $F_{\ovl{s}}^{\mathrm{top}}\to
 F_{\ovl{t}}^{\mathrm{top}}$.
\end{proof}

\subsection{The formalism of fundamental groups.}\label{sub::1.4}  
 
\begin{para}\label{para::1.4.1} Let $\Cov^\bullet_S$ be a full
 subcategory of $\Cov^{\mathrm{et}}_S$ which is {\it stable under taking
 connected components, fiber products} (over $S$) {\it and quotients}
 (in the sense of \ref{lem::1.2.7}.). Examples:  $\Cov^{\mathrm{et}}_S$,
 $\Cov^{\mathrm{top}}_S$, $\Cov^{\mathrm{alg}}_S$.

 We denote by $F^\bullet_{S,\ovl{s}}=F^\bullet_{\ovl{s}}$ the
 restriction of $F^{\mathrm{et}}_{\ovl{s}}$ to $\Cov^\bullet_S$, and set 

 \begin{equation*}
  \pi_1^{\bullet}(S,\ovl{s}) = \Aut F^\bullet_{\ovl{s}}
 \end{equation*}
 topologized by considering as fundamental
 open neighborhoods of $1$ the stabilizers $\Stab^\bullet_{S',\ovl{s}'}$
 in $\pi_1^{\bullet}(S,\ovl{s})$ of arbitrary geometric points $\ovl{s}'$
 above $\ovl{s}$ in arbitrary \'etale coverings $S'\to S$ in
 $\Cov^\bullet_S\;$ (here we see $\pi_1^{\bullet}(S,\ovl{s})$ as acting
 on the left on   $F^\bullet_{\ovl{s}}(S')$).

 That the topology is indeed compatible with the group law is due to
 the fact that the system of subgroups $\Stab^\bullet_{S',\ovl{s}'}$ is
 stable under intersection (since $\Cov^\bullet_S$ is stable under
 taking fiber products) and under conjugation (which amounts to a change
 of $\ovl{s}'$).

 \medskip In the case $\Cov^\bullet_S = \Cov^{\mathrm{et}}_S,\;
 \pi_1^{\mathrm{et}}(S,\ovl{s})$ is the \'etale fundamental group
 introduced and studied by De Jong \cite{efgonas}.

 In the case $\Cov^\bullet_S = \Cov^{\mathrm{top}}_S,\;
 \pi_1^{\mathrm{top}}(S,\ovl{s})=\pi_1^{\mathrm{top}}(S, s)$ is (an
 incarnation of) the topological fundamental group.

 In the case $\Cov^\bullet_S = \Cov^{\mathrm{alg}}_S,\;
 \pi_1^{\mathrm{alg}}(S,\ovl{s})$ is a profinite group, which will be
 called the {\it algebraic fundamental group}. The theory is embodied in
 \cite[V]{sga1}.
 If $S$ is the analytification of an algebraic smooth
 $K$-variety $S^{\mathrm{alg}}$ and if $\mathrm{char}(K)=0$, then this
 group coincides with Grothendieck's algebraic fundamental 
 group. This follows from the {\it Gabber-L\"utkebohmert theorem},
 according to which the GAGA functor is an equivalence between the
 categories of finite \'etale coverings of $S^{\mathrm{alg}}$ and $S$
 respectively \cite{repfapf}.
\end{para}

\begin{lem}\label{lem::1.4.2} The natural continuous map 
 \begin{equation*}
  \phi: \pi_1^\bullet(S,\ovl{s})\to
  \underset{\Stab^\bullet_{S',\ovl{s}'}}{\limproj} 
  \pi_1^\bullet(S,\ovl{s})/\Stab^\bullet_{S',\ovl{s}'}
 \end{equation*}
 is a homeomorphism. In particular,
 $\pi_1^\bullet(S,\ovl{s})$ is a separated pro-discrete (hence totally
 discontinuous) topological space.
\end{lem}

\begin{proof}
 $\phi$ sends any neighborhood of $\id\in
 \pi_1^\bullet(S,\ovl{s})$ to a neighborhood of $\phi(\id)\in
 im(\phi)$. Since $\phi$ is equivariant with respect to the continuous
 action of $\pi_1^\bullet(S,\ovl{s})$ by left translation, it is open
 onto its image. Hence it suffices to show that it is bijective. Let
 $\underline{\gamma}:=(\gamma_{S',\ovl{s}'})$ be a compatible family of
 elements of $\pi_1^\bullet(S,\ovl{s})/\Stab^\bullet_{S',\ovl{s}'}$, \ie
 represent an element of $\limproj\;
 \pi_1^\bullet(S,\ovl{s})/\Stab^\bullet_{S',\ovl{s}'}$.

 \noindent For any $S'$ in $\Cov^{\bullet}_S$ and any $\ovl{s}'\in
 F^\bullet_{\ovl{s}}(S')$, $\gamma_{S',\ovl{s}'}.\ovl{s}'$ is another
 well-defined element of $F^\bullet(S')$, which we denote by
 $\gamma.\ovl{s}'$. From the compatibility of the
 $\gamma_{S',\ovl{s}'}$'s, it follows that $\gamma$ defines an
 endomorphism of the functor $F^\bullet_{\ovl{s}}$.

 On the other hand, one has a homeomorphism 
 \begin{equation*}
 \inv:\; \limproj\; \pi_1^\bullet(S,\ovl{s}) /\Stab^\bullet_{S',\ovl{s}'}
  \to \limproj\;
  \Stab^\bullet_{S',\ovl{s}'}\backslash\pi_1^\bullet(S,\ovl{s})
 \end{equation*}

 \noindent induced by $g\mapsto g^{-1}$. If we copy the construction
 $\underline{\gamma}\mapsto \gamma$ starting from
 $\inv(\underline{\gamma})$ and using right actions instead of left
 actions, the resulting endomorphism of $F^\bullet_{\ovl{s}}$ is nothing
 but the inverse of $\gamma$. Hence $\gamma$ is an element of $
 \pi_1^\bullet(S,\ovl{s})$, and it is clear that
 $\underline{\gamma}\mapsto \gamma$ is inverse to the natural map
 $\displaystyle\pi_1^\bullet(S,\ovl{s})\to \limproj\; 
 \pi_1^\bullet(S,\ovl{s})/\Stab^\bullet_{S',\ovl{s}'}$.
\end{proof}

\begin{para}\label{para::1.4.3} It is clear from the definitions that 
 $F^\bullet_{\ovl{s}}$ may be enriched to a functor
 \begin{equation*}
  \Cov^\bullet_S \longrightarrow
   \text{$\pi_1^\bullet(S,\ovl{s})$-$\Sets$},\;\;S'\mapsto
   F^\bullet_{\ovl{s}}(S')
 \end{equation*}
 where $\pi_1^\bullet(S,\ovl{s})$-$\Sets$ is the category of discrete
 sets endowed with a continuous 
 left action of $\pi_1^\bullet(S,\ovl{s})$. 
\end{para}

\begin{pro}\label{pro::1.4.4} Up to isomorphism (unique up to composition
 by inner automorphisms), $\pi_1^\bullet(S,\ovl{s})$ does not depend on
 $\ovl{s}$.
\end{pro}
\begin{proof}
 The space
 $\Isom(F^\bullet_{\ovl{s}},F^\bullet_{\ovl{t}})$, non-empty due to
 \ref{lem::1.3.3}, is formally a
 $\pi_1^{\mathrm{et}}(S,\ovl{t})-\pi_1^{\mathrm{et}}(S,\ovl{s})$-bitorsor,
 and the result follows.
\end{proof}

\begin{thm}\label{thm::1.4.5} The enriched functor $F^\bullet_{\ovl{s}}$
\begin{equation*}
 \Cov^\bullet_S \longrightarrow \text{$\pi_1^\bullet(S,\ovl{s})$-$\Sets$}
\end{equation*}
is fully faithful, and extends to an equivalence of categories  

\begin{equation*}
\{\text{disjoint unions of objects of }\Cov^\bullet_S\} \longrightarrow
\text{$\pi_1^\bullet(S,\ovl{s})$-$\Sets$}.
\end{equation*}
Connected coverings correspond to $\pi_1^\bullet(S,\ovl{s})$-orbits.
\end{thm}

\begin{proof}
 (\cf \cite[2.10]{efgonas} in the \'etale case). If
 $f:S'\to S$ is an \'etale covering in $\Cov^\bullet_S$ and
 $\ovl{s}',\;\ovl{s}''$ are elements of $F_{\ovl{s}}(S')$ in the same
 connected component of $S'$, proposition \ref{lem::1.3.3} ensures the existence of
 an \'etale path $\alpha$ from $\ovl{s}'$ to $\ovl{s}''$. Then
 $f_\ast(\alpha)\in \pi_1^{\mathrm{et}}(S,\ovl{s})$, and
 $f_\ast(\alpha).\ovl{s}'=\ovl{s}''$; of course, we may replace here
 $f_\ast(\alpha)$ by its image in $\pi_1^\bullet(S,\ovl{s})$. This shows
 that $\pi_1^\bullet(S,\ovl{s})$-orbits in $F^\bullet_{\ovl{s}}(S')$
 correspond bijectively to connected components of $S'$. In particular

 \begin{equation*}
 F^\bullet_{\ovl{s}}(S') = 
  \pi_1^\bullet(S,\ovl{s})/{\Stab^\bullet_{S',\ovl{s}'}}
 \end{equation*}
 if $S'$ is connected.

 This implies that the enriched functor $F^\bullet_{\ovl{s}}$ is fully
 faithful, because if $S'$ is connected, a morphism between \'etale
 coverings $S'\to S''$ of $S$ may be identified via its graph with a
 connected component $T$ of the \'etale covering $S'\times_S S''$ (which
 belongs to $\Cov^\bullet_S$ by assumption) such that
 $F^\bullet_{\ovl{s}}(T)\xrightarrow{\sim}
 F^\bullet_{\ovl{s}}(S')$ (this can be checked locally on $S$, and this
 reduces to a well-known fact about finite algebraic coverings).

 It remains to show that any $\pi_1^\bullet(S,\ovl{s})$-orbit comes from
 a (connected) \'etale covering. By definition, such an orbit may be
 written $ \pi_1^\bullet(S,\ovl{s})/H$, where $H$ is a subgroup
 containing the stabilizer ${\Stab^\bullet_{S',\ovl{s}'}}$ of some
 geometric point $\ovl{s}'$ above $\ovl{s}$ in some connected \'etale
 covering $S'$ of $S$. As we have seen, the connected components of
 $S'\times_S S'$ correspond bijectively to
 $\pi_1^\bullet(S,\ovl{s})$-orbits in $F^\bullet_{\ovl{s}}(S'\times_S
 S')= F^\bullet_{\ovl{s}}(S')^2.$ The union $R\subset S'\times_S S'$ of
 connected components which corresponds to the union of points
 $(\ovl{s}',\ovl{s}'.h),\;h\in H$, is an equivalence relation on $S'$
 over $S$. By lemma \ref{lem::1.2.7}, the quotient $S''=S'/R$ is an
 \'etale covering 
 of $S$ which belongs to $\Cov^\bullet_S$ by assumption. If $\ovl{s}''$
 denotes the image of $\ovl{s}'$ in $S''$, it is clear that
 $H=\Stab^\bullet_{S'',\ovl{s}''}$.
\end{proof}

\begin{para}\label{para::1.4.6} A connected \'etale covering $S'\to S$ is
 {\it Galois} if $S$ is the quotient $S'/(\Aut_S S')$; $\Aut_S S'$ is
 called the {\it Galois group} of the covering. By \cite[4.1.9]{ecfnas},
 an \'etale Galois covering of $S$ with group $G$ is the same as a
 connected principal homogeneous space under the constant 
 $S$-group $G_S$\footnote{{\it loc.~cit.}  4.1.4, which is used in the proof,
 is available since we are dealing with strictly analytic spaces, \cf
 {\it loc.~cit.} 4.1.5.}.

\noindent Through the dictionary
 of \ref{thm::1.4.5}, a Galois \'etale covering in $\Cov^{\bullet}_S$ 
 corresponds to a surjective continuous 
 homomorphism 

 \begin{equation*}
  \pi_1^\bullet(S,\ovl{s})\to G,
 \end{equation*}
 with $G$ discrete ($\cong \Aut_S S'$). In other words, a (pointed) covering
 $(S',\ovl{s}')$ in $\;\Cov^{\bullet}_S\;$ is Galois if and only if 
 ${\Stab^\bullet_{S',\ovl{s}'}}\subset
 \pi_1^\bullet(S,\ovl{s})$ is normal. 

 Example: the logarithm (\ref{exs::1.2.5}. (i)) is a Galois \'etale covering of
 $(\A^1_{\C_p})^{\mathrm{an}}$ with group $\mu_{p^\infty}\cong \Q_p/\Z_p$. 

The following is also clear:

\end{para}

\begin{cor}\label{cor::1.4.7} $\pi_1^\bullet(S,\ovl{s})$ is a
 pro-discrete group if and only if in $\Cov^\bullet_S$ any 
 connected covering is dominated by a Galois covering.
\end{cor}

\par\noindent {\it Remark}. This property holds in
$\Cov^{\mathrm{top}}_S$ and $\Cov^{\mathrm{alg}}_S$, but not in
$\Cov^{\mathrm{et}}_S$ in general. Let for instance $S=\P^1_{\C_p}$, and
let us consider the Galois covering given by the 
Gross-Hopkins period mapping (II.\ref{para:6.3.3}). According to 
\cite[7.4]{efgonas} 
(using the viewpoint sketched in II.\ref{sub:6.4}), the corresponding
discrete quotient of $\pi_1^{\mathrm{et}}(\P^1,\ovl{s})$ is the composition 

\begin{equation*}
 \pi_1^{\mathrm{et}}(\P^1,\ovl{s})\to SL_2(\Q_p) \to SL_2(\Q_p)/SL_2(\Z_p)
\end{equation*}
where the first map is a continuous surjective homomorphism,
and the quotient topology on $SL_2(\Q_p)$ is non-discrete ({\it
loc.~cit.} p. 116). Since the biggest normal subgroup of $SL_2(\Q_p)$
contained in $SL_2(\Z_p)$ is $\{\pm \id\}$, we see that its preimage in
$\pi_1^{\mathrm{et}}(\P^1,\ovl{s})$ cannot be open. Thus
$\pi_1^{\mathrm{et}}(\P^1,\ovl{s})$ is {\it not a pro-discrete group}.

\begin{cor}\label{cor::1.4.8} Let $\Cov^{\bullet\bullet}_S$ be a full
 subcategory of $\Cov^{\bullet}_S$ stable under taking connected
 components, fiber products (over $S$) and quotients. Then the
 corresponding continuous homomorphism
 $\pi_1^\bullet(S,\ovl{s})\to\pi_1^{\bullet\bullet}(S,\ovl{s})$ has
 dense image (it is neither strict\footnote{recall that a continuous
 homomorphism between topological groups is strict if the induced
 continuous bijective homomorphism between the coimage and the image is
 an isomorphism, i.e. bicontinuous; the coimage is the quotient of the
 source by the kernel, with the quotient topology.} nor surjective in
 general).  \hfill\break In particular, if
 $\;\Cov^{\mathrm{alg}}_S\subset \Cov^{\bullet}_S,\;$ then
 $\;\pi_1^{\mathrm{alg}}(S,\ovl{s})$ is the profinite
 completion\footnote{the profinite completion of a topological group is
 the inverse limit of its quotients by {\it open} normal subgroups of
 finite index.} of $\;\pi_1^\bullet(S,\ovl{s})$. 
\end{cor}

\begin{proof}
 By lemma \ref{lem::1.4.2}, it suffices to prove that for any $S'$ in
 $\Cov^{\bullet\bullet}_S$ and any $\ovl{s}'\in
 F^\bullet_{\ovl{s}}(S')=F^{\bullet\bullet}_{\ovl{s}}(S')$, the natural
 map $\pi_1^\bullet(S,\ovl{s})\to
 \pi_1^{\bullet\bullet}(S,\ovl{s})/\Stab^{\bullet\bullet}_{S',\ovl{s}'}$
 is surjective. We may replace $S'$ by the connected component of
 $s'$. According to the previous theorem,
 $\pi_1^{\bullet\bullet}(S,\ovl{s})/\Stab^{\bullet\bullet}_{S',\ovl{s}'}$
 and $\pi_1^{\bullet}(S,\ovl{s})/\Stab^{\bullet}_{S',\ovl{s}'}$ are then
 identified to
 $F^\bullet_{\ovl{s}}(S')=F^{\bullet\bullet}_{\ovl{s}}(S')$, whence the
 first assertion. The second assertion is clear, since
 $\;\pi_1^{\mathrm{alg}}(S,\ovl{s})$ is profinite.
 
 To see that the image of
 $\;\pi_1^\bullet(S,\ovl{s})\to\pi_1^{\bullet\bullet}(S,\ovl{s})$ is
 neither strict nor surjective in general, it suffices to take $S$ such that
 $\pi_1^{\mathrm{top}}(S,s)$ is infinite (\eg a Tate curve). Then
 $\pi_1^{\mathrm{et}}(S,\ovl{s})\to\pi_1^{\mathrm{alg}}(S,\ovl{s})$ is
 neither surjective nor strict. This follows from the commutative square 
 of topological groups
 \begin{equation*}
 \begin{CD}
  \pi_1^{\mathrm{et}}(S,\ovl{s}) @>>> \pi_1^{\mathrm{alg}}(S,\ovl{t})\\
  @VVV @VVV\\
  \pi_1^{\mathrm{top}}(S,s) @>>> \widehat{\pi_1^{\mathrm{top}}(S,t)}
 \end{CD}
 \end{equation*} 
 where the roof denotes the profinite completion; the coimage (\resp
 image) of the ``south-east map''
 $\pi_1^{\mathrm{top}} (S,\ovl{s}) \to  \widehat{\pi_1^{\mathrm{top}}(S,
 t)}$ is infinite discrete (\resp compact).
\end{proof}

\medskip In this context, the following general lemma is useful:

\begin{lem}\label{lem::1.4.9} Let $G, H$ be separated topological groups
 such that any neighborhood of $1$ contains an open 
 subgroup, and let $f: G\to H$ be a continuous homomorphism.
 \begin{enumerate}
  \item If the induced functor 
	$f^\ast:\,\text{$H$-$\Sets$}\to G$-$\Sets$ is essentially
	surjective, then $f$ is injective.
  \item If $f^\ast:\,H$-$\Sets \to G$-$\Sets$ is an equivalence, then $f$ is
	injective, strict and has dense image.
  \item If $f^\ast:\,H$-$\Sets\to G$-$\Sets$ is an equivalence and $G=
	\limproj G/U$ (the limit running over open subgroups of $G$),
	then $f$ is an isomorphism.
  \item If $G$ is complete and if the normal open subgroups are cofinal
	among the open subgroups $U$, then $G=\limproj G/U$.
 \end{enumerate}
\end{lem}
  
\begin{proof}
 (i) otherwise, let $g\neq 1$ be in the kernel of $f$. Because $G$ is
 separated, there exists a neighborhood $U$ of $1$ not containing $g$,
 and by hypothesis, one may assume that $U$ is an open subgroup of
 $G$. By assumption the (connected) $G$-set $G/U$ corresponds via
 $f^\ast$ to a (connected) $H$-set $X$. Since $f(g)=1$, $g$ acts 
 trivially on $f^\ast(X)$, in contradiction with the fact that $g\notin U$.

 (ii) To show that $f$ has dense image, it suffices to show that for $h\in
 H$ and any neighborhood $V$ of $1$, $f(G).h \cap V\neq \emptyset$. By
 hypothesis, one may assume that $V$ is an open 
 subgroup of $H$. Thus one has to show that $H=G.V$ for 
 any open subgroup $V$ of $H$, \ie $G$ acts transitively on $H/V$. If the
 action is not transitive, $f^{\ast}(H/V)$ can be 
 written as a disjoint union of $G$-sets, which contradicts the
 the assumption that $f^{\ast}$ is an equivalence. 

 \noindent It remains to show that $f$ is strict. But we have seen that
 $f^\ast$ induces an equivalence on transitive sets. 
 Therefore $f^{-1}$ induces a bijection between open subgroups of $H$
 and open subgroups of $G$, whence the result.

 (iii) Indeed one has a commutative diagram
 \begin{equation*}
 \begin{CD}
  G @>\subset>> H\\
  @V\wr VV @VVV \\
 \limproj G/U @>\sim>> \limproj H/V
 \end{CD}
 \end{equation*}
 and the right vertical map is injective since $H$ is separated.

 (iv) follows from \cite[III.7.2, prop. 2]{bourbaki74:_topol}.
\end{proof}

\begin{para}\label{para::1.4.10} Let $f:T\to S$ be any morphism of connected
 $K$-manifolds. We have the natural functor 

\begin{equation*}
 \Cov^{\mathrm{et}}_S\to \Cov^{\mathrm{et}}_T:\;\; S'\mapsto S'\times_S T.
\end{equation*}
Let us {\it assume} that it induces a functor 

\begin{equation*}
 f^{\bullet\ast}: \Cov^\bullet_S\to \Cov^\bullet_T.
\end{equation*}
Since
$F^\ast_{S,\ovl{s}}=F^\ast_{T,\ovl{t}}\circ f^{\bullet\ast}$ if
$\ovl{s}=f(\ovl{t})$, this gives rise to a continuous homomorphism

\begin{equation*}
f^\bullet_\ast :\pi_1^\bullet(T,\,\ovl{t})\to
 \pi_1^\bullet(S,\,\ovl{s}=f(\ovl{t})).
\end{equation*}
\end{para}

\begin{cor}\label{cor::1.4.11} $f^\bullet_\ast$ is an isomorphism of
 topological groups if and only if $f^{\bullet\ast}$ 
 is an equivalence of categories. 
\end{cor}

\begin{proof}
 If $f^\bullet_\ast$ is an isomorphism, it induces an equivalence 
 $\pi_1^\bullet(S,\,\ovl{s})$-$\Sets\cong
 \pi_1^\bullet(T,\,\ovl{t})$-$\Sets$, and it follows from \ref{thm::1.4.5} that
 $f^{\bullet\ast}$ is an equivalence of categories. The other
 implication is immediate.
\end{proof}

\begin{cor}\label{cor::1.4.12} Assume that $f$ is an \'etale covering in
 $\Cov^\bullet_S$, thus corresponding to the 
 $\pi_1^\bullet(S,\,\ovl{s})$-orbit
 $\pi_1^\bullet(S,\,\ovl{s})/\Stab^\bullet_{T,\ovl{t}}$ through the 
 dictionary of \ref{thm::1.4.5}.
 \begin{enumerate}
  \renewcommand{\theenumi}{\alph{enumi}}
  \item Then $\Stab^\bullet_{T,\ovl{t}}$ is the closure of the image of
	$f^\bullet_\ast$.  In particular, if $f$ is a Galois covering,
	its Galois group is
	$\;\pi_1^\bullet(S,\,\ovl{s})/(\Im\,f^\bullet_\ast)^-$.
  \item Assume moreover that for any covering $T'\to T$ in
	$\Cov^\bullet_T$, the composition with $f$ gives rise to a
	covering $T'\to S$ in $\Cov^\bullet_S$. Then $f^\bullet_\ast$
	is open and injective (hence $\;\pi_1^\bullet(T,\,\ovl{t})$ may
	be identified with the open and closed subgroup
	$\Stab^\bullet_{T,\ovl{t}}$ of $\pi_1^\bullet(S,\,\ovl{s})$).  
 \end{enumerate}
\end{cor}
\begin{proof}
 (a) Looking at the natural identifications 
 \begin{equation*}
 \pi_1^\bullet(S,\,\ovl{s})/\Stab^\bullet_{T,\ovl{t}}\cong
 F^\bullet_{\ovl{s}}(T)\cong F^\bullet_{\ovl{t},\ovl{t}}(T\times_S \,T)
 \end{equation*}
 of $\pi_1^\bullet(S,\,\ovl{s})$-orbits. The connected component of
 $T\times_S \,T$ containing $(\ovl{t},\ovl{t})$ is the 
 diagonal $T\subset T\times_S \,T$, therefore the
 $\pi_1^\bullet(T,\,\ovl{t})$-orbit is reduced to the point
 $(\ovl{t},\ovl{t})$. This shows that $\Im f^\bullet_\ast$ lies in
 $\Stab^\bullet_{T,\ovl{t}}$.

 To show that it is dense, we have to show that for any connected $S'$ in
 $\Cov^\bullet_S$ and any geometric point $\ovl{s}'$ of $S'$ above
 $\ovl{s}$, the map $\pi_1^\bullet(T,\,\ovl{t})\to
 \pi_1^\bullet(S,\,\ovl{s})/\Stab^\bullet_{T,\ovl{t}}$ induced by
 $f^\bullet_\ast $ is surjective. This follows from \ref{thm::1.4.5} as
 in \ref{cor::1.4.8}.

 (b) Under the extra assumption, one has functors $\Cov^\bullet_T\to
 \Cov^\bullet_S\to \Cov^\bullet_T: T'/T\mapsto T'/S ;\;S'/S\mapsto
 (S'\times_S \,T)/T$.

 \noindent For any open subgroup $H\subset  \pi_1^\bullet(T,\,\ovl{t})$
 (corresponding as in \ref{thm::1.4.5} to a pointed connected covering
 $(T',\ovl{t}')\;$), $\pi_1^\bullet(T,\,\ovl{t})/H\cong
 F^\bullet_{\ovl{t}' }(T' )$ is a subset of $
 F^\bullet_{\ovl{t}',\ovl{t}}(T'\times_S\,T)\cong 
 F^\bullet_{\ovl{s}}(T')\cong
 \pi_1^\bullet(S,\,\ovl{s})/\Stab^\bullet_{T',\ovl{t}'}$.
 In other words, $f^\bullet_\ast $ induces an open injective homomorphism 
 \begin{equation*}
 \pi_1^\bullet(T,\,\ovl{t})/H \hookrightarrow
 \pi_1^\bullet(S,\,\ovl{s})/f^\bullet_\ast(H).
 \end{equation*}
 Indeed, $f^\bullet_\ast(H)=\Stab^\bullet_{T',\ovl{t}'} $ is open in
 $\pi_1^\bullet(S,\,\ovl{s}) $: such groups correspond to connected
 coverings of $S$ which dominate $T$ (\ie which factorize over $f$).
 Since such coverings are cofinal in 
 $\Cov^\bullet_S$, we can pass to the limit, using \ref{lem::1.4.2}: 

 \begin{equation*}
 \pi_1^\bullet(T,\,\ovl{t}) = \limproj \pi_1^\bullet(T,\,\ovl{t})/H
 \hookrightarrow \limproj \pi_1^\bullet(S,\,\ovl{s}) /f^\bullet_\ast(H)
 \cong \pi_1^\bullet(S,\,\ovl{s})
 \end{equation*}
 and conclude that $f^\bullet_\ast$ is an open injective homomorphism.
\end{proof}

\subsection{Example: Ramero's locally algebraic \'etale fundamental
group.}\label{sub::1.5}

\begin{para}\label{para::1.5.1}
 When $S$ is countable at infinity, L. Ramero has introduced an
 interesting category $\Cov^{\mathrm{loc.alg}}_S$ of \'etale coverings, called
 `locally algebraic', which generalize the logarithmic covering of
 $\A^1$, and which are intimately connected with the local study of
 meromorphic differential equations around singularities \cite[4]{ramero98}.

 An \'etale covering is locally algebraic if the preimage of every {\it
 compact subdomain} $W\subset S$ is a disjoint union $\coprod W_i$ of
 subdomains $W_i$ which are finite \'etale over $W$. Since $S$ is (by
 assumption) a union of an increasing countable family of compact
 subdomains $S_j,\; j\geq 1$, it suffices to check the condition for 
 $W=S_j,\forall j.$ 

 \noindent This category is stable under taking connected components,
 fiber products (over $S$) and quotients. Moreover, it is stable by
 disjoint union and `composition'. Assuming that
 $s\in S_1$, one has a natural homomorphism  
 \begin{equation*}
 \liminj\pi_1^{\mathrm{loc.alg}}(S_j,\,\ovl{s}) = \liminj
  \pi_1^{\mathrm{alg}}(S_j,\,\ovl{s})\to \pi_1^{\mathrm{loc.alg}}(S,\,\ovl{s}).
 \end{equation*}
\end{para}

\begin{pro}\label{pro::1.5.2} This is an isomorphism:
 $\pi_1^{\mathrm{loc.alg}}(S,\,\ovl{s})\cong \liminj
 \pi_1^{\mathrm{alg}}(S_j,\,\ovl{s})$.
\end{pro}
   
\begin{proof}
 It is clear that an automorphism of $F^{\mathrm{loc.alg}}_{S,\ovl{s}}$
 amounts to a compatible family of automorphisms of
 $F^{\mathrm{alg}}_{S_j,\ovl{s}}$, hence the homomorphism is
 bijective. Moreover, the induced map $ \liminj
 \pi_1^{\mathrm{alg}}(S_j,\,\ovl{s})$-sets $\to
 \pi_1^{\mathrm{loc.alg}}(S,\,\ovl{s})$-sets is an equivalence: both
 sides are equivalent to $\Cov^{\mathrm{loc.alg}}_S$, \cf \ref{thm::1.4.5} and
 \cite[4.1.6]{ramero98}. Therefore open subgroups correspond to stabilizers as
 usual, and the bijection $\pi_1^{\mathrm{loc.alg}}(S,\,\ovl{s})\cong
 \liminj \pi_1^{\mathrm{alg}}(S_j,\,\ovl{s})$ is bicontinuous.
\end{proof}

\medskip\noindent It is not known whether
$\pi_1^{\mathrm{loc.alg}}(S,\,\ovl{s})$ is a pro-discrete group (\cf in
this direction {\it loc.~cit.} 4.3). 

\medskip\noindent {\it Example.} The logarithm of any Lubin-Tate formal
group over a $p$-adic field gives rise to a locally algebraic 
\'etale covering of the analytic affine line, \cf {\it loc.~cit.} 6.1.1.

On the other hand, the Gross-Hopkins covering of $\P^1$ is not locally
algebraic, although there is a finite affinoid 
cover $(V_i)$ of $\P^1$ such that the restriction of the covering to any
$V_i$ is a disjoint union of finite \'etale coverings.

\newpage
\section{Temperate fundamental groups.}\label{sec::2}

\begin{abst}
We introduce the notion of temperate \'etale
covering. Such coverings are essentially built from (possibly infinite)
topological coverings of finite \'etale coverings; they are classified
by temperate fundamental groups, which seem to be the right
non-archimedean equivalents of fundamental groups of complex
manifolds. These groups are not discrete, but nevertheless often possess
many infinite discrete quotients.
\end{abst}

\subsection{Temperate \'etale coverings and temperate fundamental
groups.}\label{sub::2.1}

\medskip In the archimedean context, the categories
$\Cov^{\mathrm{et}}_S$ and $\Cov^{\mathrm{top}}_S$ coincide; in the
non-archimedean situation, in contrast, $\Cov^{\mathrm{top}}_S$ is ``too
small'' (for instance, it consists of disjoint unions of copies of $S$
if $S$ has dimension one and tree-like reduction, \eg if $S$ is an
annulus), while  $\Cov^{\mathrm{et}}_S$ is ``too big'' (for $S=\P^1$, it
contains many non-trivial objects).

We introduce an intermediate category, which seems to be a closer
analogue of the category of unramified coverings in the 
archimedean context.
 
\begin{dfn}\label{def::2.1.1} An \'etale covering $S'\to S$ is {\rm
 temperate} if it is a quotient\footnote{in the sense of
 \ref{lem::1.2.7}.} of a 
 composite \'etale covering $T'\to T\to S$, where $T'\to T$ is a
 (possibly infinite) topological covering, and $T\to S$ is a finite
 \'etale covering.
\end{dfn}

Hence we have the commutative square 
\begin{equation*}
 \begin{matrix}
  &&T'&&\\
  & \swarrow &&\searrow &\\
  T&&&&S'\\
  &\searrow && \swarrow &\\
  &&S&&
 \end{matrix}
\end{equation*}
with $T'\in \Cov^{\mathrm{top}}_{T},\;T\in \Cov^{\mathrm{alg}}_{S},
S'\in \Cov^{\mathrm{et}}_S$. 

\medskip\noindent {\it Remark}. The notion of temperate covering is
stable under base change as in \ref{para::1.4.10}. 
 
\begin{lem}\label{lem::2.1.2} Let $T\to S$ be a finite Galois \'etale
 covering, and let $\tilde T\to T$ be the universal topological covering
 of $T$. Then the composite $\tilde T \to S$ is a Galois temperate
 covering. 
\end{lem}

\begin{proof}
 From \ref{lem::1.2.4}, we know that it is an \'etale covering, automatically
 temperate. It remains to show that it is Galois. By the universal
 property of $\tilde T$, any automorphism of $T$ lifts to $\tilde T$,
 hence there is an exact sequence of discrete groups:
 \begin{equation*}
 1\to \bigl(\pi_1^{\mathrm{top}}(T,t) \cong\bigr)\; \Aut_T \tilde T \to
 \Aut_S \tilde T \to \Aut_S\,  T\to 1. 
 \end{equation*}
 Since $\tilde T/(\Aut_T \tilde T) =T$ and $T/(\Aut_S T)=S,$ it
 follows that $ \tilde T/(\Aut_S \tilde T)=S$.
\end{proof}

We see from the lemma that in the definition of a temperate
covering, we may assume that $T'\to S$ and $T'\to T$ are Galois.  

\begin{rems}
 (i) It is {\it not true} that any quotient of a tower $T'\to T\to S$
 as in definition \ref{def::2.1.1} is again of this form. As a
 counter-example, let us consider a Tate elliptic curve $E$ over
 $K=\C_p$, and set $T= E\setminus E[2]$ (2-torsion points
 removed). Then $\tilde T $ is $\G_m^{\mathrm{an}}$ deprived from
 the countable subgroup $\pm (\sqrt q)^\Z$, and $T\cong \tilde
 T/q^\Z $. The inversion $x\mapsto 1/x$ on $\tilde T$ induces the
 inversion on $E$, and $\G_m/\{\id, \inv\}=\A^1 \;$ (via $x\mapsto
 x+1/x),\;E/\{\id,\inv\}=\P^1$. Let us consider the infinite
 temperate covering 

 \begin{equation*}
  S'=\tilde{T}/\{\id, \inv\} \to S= T/\{\id, \inv\}
   =\P^1\setminus\{0,1,\lambda,\infty\}.
 \end{equation*}
 Then $S'$ is not a topological covering of any finite \'etale
 covering $S''$ of $S$. Indeed, otherwise $\tilde T$ would be the
 universal topological covering of some connected component $T'$
 of $T\times_S S''$. The tower $\tilde T\to T' \to T$ would then
 extend to a tower of topological coverings $\G_m^{\mathrm{an}}\to
 E'\to E$, where $E'$ is another Tate curve isogenous to $E$ (we
 use the fact that the topological fundamental group of $T'$ and
 $E$ is the same). Then $S''\subset E'/\{\id, \inv\}\cong
 \P^1$. But any open subset in $\P^1_{\C_p}$ is simply-connected,
 whence $\tilde T= S''$, a contradiction.
 
 \noindent This example also shows that {\it a temperate covering may
 become a topological covering after finite \'etale base change, without
 being itself a topological covering}.

 (ii) The notion of temperate covering is not local on $S$. In fact,
 the Gross-Hopkins \'etale covering of
 $(\P^1_{\C_p})^{\mathrm{an}}$ is not temperate, although there is
 a finite open cover $(U_i)$ of $(\P^1_{\C_p})^{\mathrm{an}}$ such
 that it becomes a disjoint union of finite \'etale coverings over
 each $U_i$ (by compactness).
\end{rems}

\begin{para}\label{para::2.1.3} We denote the category of temperate
 (\'etale) coverings of $S$ by $\Cov^{\mathrm{temp}}_S$. 
 From the fact that topological and finite \'etale coverings
 respectively are stable under taking connected components and fiber
 products, it is not difficult to see that $\Cov^{\mathrm{temp}}_S$ is
 stable under taking connected components and fiber products (over
 $S$). By definition, it is stable under taking quotients. 
 
\medskip According to the results of 1.4, there is thus defined a
 separated totally discontinuous group, the {\it temperate fundamental group}
 $\pi_1^{\mathrm{temp}}(S,\ovl{s})$, which does not depend on the
 geometric point $\ovl{s}$ up to isomorphism. Moreover:
\end{para}

\begin{lem}\label{lem::2.1.4} $\pi_1^{\mathrm{temp}}(S,\ovl{s})$ is a
 pro-discrete group (hence complete, 
\cf \cite[III. 57.7]{bourbaki74:_topol}).  A basis of open neighborhoods of
$1$ is given by normal closed subgroups $H$ such that
 $\pi_1^{\mathrm{temp}}(S,\ovl{s})/H$ is the Galois group of the
 universal topological covering of some finite \'etale Galois covering of $S$.
\end{lem}
\begin{proof}
 This follows from \ref{cor::1.4.7}. and \ref{lem::2.1.2} (connected
 temperate coverings admit a Galois closure).
\end{proof}

 \medskip 
 The inclusion $\; \Cov^{\mathrm{top}}_S \hookrightarrow
 \Cov^{\mathrm{temp}}_S$ gives rise to a {\it surjective} homomorphism 
 \begin{equation*}
  \pi_1^{\mathrm{temp}}(S,\ovl{s}) \to \pi_1^{\mathrm{top}}(S,s) 
 \end{equation*} 
 (\cf \ref{cor::1.4.8}, taking
 into account the discreteness of $\pi_1^{\mathrm{top}}(S,s) $), which
 shows that $\pi_1^{\mathrm{temp}}(S,\ovl{s})$, although not a discrete 
 group, has many {\it infinite discrete quotients} in general; we say
 that $\pi_1^{\mathrm{temp}}(S,\ovl{s})$ is ``lacunary''.

 \medskip On the other hand, the inclusion $\; \Cov^{\mathrm{alg}}_S
 \hookrightarrow  \Cov^{\mathrm{temp}}_S$ gives rise to a 
 continuous homomorphism with dense image
 \begin{equation*}
  \pi_1^{\mathrm{temp}}(S,\ovl{s})\to \pi_1^{\mathrm{alg}}(S,\ovl{s}),
 \end{equation*}
 and $ \pi_1^{\mathrm{alg}}(S,\ovl{s}) $ may be identified with the
 profinite completion of $\pi_1^{\mathrm{temp}}(S,\ovl{s})$. This shows
 that $\pi_1^{\mathrm{temp}}(S,\ovl{s})$ is not ``too small'' (as opposed
 to $\pi_1^{\mathrm{top}}(S,s) $), and the following two propositions
 show that it is not ``too big'' as well (as opposed to
 $\pi_1^{\mathrm{et}}(S,\ovl{s}) $): 
 
\begin{pro}\label{pro::2.1.5} If $\dim\,S=1$, the homomorphism
 $\pi_1^{\mathrm{temp}}(S,\ovl{s})\to \pi_1^{\mathrm{alg}}(S,\ovl{s}) $
 is injective. 
\end{pro}

\begin{proof}
 By \ref{lem::2.1.2} and \ref{lem::2.1.4}, it suffices to see that if
 $T\to S$ is a finite 
 Galois \'etale covering, the Galois group of the temperate covering
 $\tilde T\to S$ is residually finite. But this group contains
 $\pi_1^{\mathrm{top}}(T,t)$ as a normal subgroup of finite index. The
 point is that in dimension one, topological fundamental groups are free,
 hence residually finite.
\end{proof}

\medskip The lacunary subgroup
$\pi_1^{\mathrm{temp}}(\P^1_{\C_p}\setminus \{0,1,\infty\})$ of the
profinite free group on two generators
$\pi_1^{\mathrm{alg}}(\P^1_{\C_p}\setminus \{0,1,\infty\})$ is
especially interesting and mysterious, and will be the 
object of our investigations in later sections\footnote{recall, on the
other hand, that $\pi_1^{\mathrm{top}}(\P^1_{\C_p}\setminus
\{0,1,\infty\})$ is trivial.}.

\begin{pro}\label{pro::2.1.6} If $S$ is algebraic and $K$ of
 characteristic zero, the Galois group of any Galois temperate 
 covering of $S$ is finitely generated. Moreover, 
 if $K$ has only countably many finite 
 extensions in a fixed algebraic closure $\ovl{K}$ (\eg if $K$ is a 
 $p$-adic field), then 
 $\pi_1^{\mathrm{temp}}(S,\ovl{s})$ and all its quotients by closed
 normal subgroups are `polish', \ie metrizable of countable type and
 complete.
\end{pro}

\begin{proof}
 Let us consider the first assertion. By \ref{lem::2.1.2}, it suffices to
 treat the case of the universal topological covering of a finite Galois
 \'etale covering.  By the already mentioned Gabber-L\"utkebohmert
 theorem, that finite \'etale covering is algebraic. The Galois group
 admits as a normal subgroup of finite index the topological fundamental
 group of that finite \'etale covering, which is finitely generated 
 (\cf \ref{para::1.1.2}); whence the result.

 \noindent Since $\pi_1^{\mathrm{temp}}(S,\ovl{s})$ is separated, its
 metrizability is equivalent to the existence of a countable fundamental
 system of neighborhood of $1$ \cite[IX.3.1 prop. 1]{bourbaki74:_topol}.
 By \ref{lem::2.1.4}, it 
 suffices to see that there are countably many finite \'etale coverings
 of $S$, 
 which follows, by d\'evissage, from the assumption that
 $\Gal(\ovl{K}/K)$ has countably many open subgroups, and from the
 topological finite generation of the algebraic fundamental group over
 $\ovl{K}$.
 This also shows that
 $\pi_1^{\mathrm{temp}}(S,\ovl{s})$ is not only a pro-discrete group
 (which already implies that it is complete), but even a countable
 inverse limit of discrete countable groups. It follows that the
 metrizable group $\pi_1^{\mathrm{temp}}(S,\ovl{s})$ admits a 
 countable dense subset, hence its topology has a countable basis 
 \cite[IX.2.8 prop. 12]{bourbaki74:_topol}.

 \noindent It follows that any closed subgroup of
 $\pi_1^{\mathrm{temp}}(S,\ovl{s})$ and any quotient of
 $\pi_1^{\mathrm{temp}}(S,\ovl{s})$ by a 
 closed normal subgroup is metrizable and complete 
 \cite[IX.3.1 prop. 4]{bourbaki74:_topol}, 
 and of countable type (for a quotient: because it 
 has a countable dense subset).
\end{proof}

Notice however that the first assertion does not imply that
$\pi_1^{\mathrm{temp}}(S,\ovl{s})$ is topologically finitely generated. 
On the other hand, we shall show below that it is not locally compact in
general (\ref{pro::2.3.12}).

\begin{pro}\label{pro::2.1.7} Let $f: T\to S$ be an \'etale finite Galois
 covering of a connected $K$-manifold with group $G$, and $\ovl{t}$ be a
 geometric point of $T$. Then there is an exact sequence 
 \begin{equation*}
  1\to \pi_1^{\mathrm{temp}}(T,\ovl{t})\to
 \pi_1^{\mathrm{temp}}(S,f(\ovl{t}))\to G\to 1. 
 \end{equation*}
\end{pro}

\begin{proof}
 It is clear from the definition \ref{def::2.1.1} that any temperate covering
 of $T$ gives rise, by composition with $f$, to a 
 temperate covering of $S$. Thus the proposition follows from
 \ref{cor::1.4.12}.b.
\end{proof}

\medskip\noindent {\it Remark.} We shall see below (\ref{para::2.3.3})
that on the 
other hand, an \'etale covering of $S$ which is a finite \'etale
covering of a topological covering is not necessarily temperate.

\medskip The following result is an important criterion to recognize
topological, \resp temperate, \resp locally algebraic coverings: 

\begin{thm}\label{thm::2.1.8} Let $S'$ be a Galois \'etale covering of
 $S$ with group $G$. Then
 \begin{enumerate}
  \renewcommand{\theenumi}{\alph{enumi}}
  \item if $G$ is torsion-free,
	$S'$ is a topological covering of $S$;
  \item if $G$ is virtually torsion-free, $S'$ is a temperate covering
	of $S$; more precisely, $S'$ is a topological covering of some
	finite \'etale Galois covering of $S$;
  \item assume that $S$ is countable at infinity; then $S'$ is a locally
	algebraic covering of $S$ if and only if $G$ is the union of an
	increasing sequence of finite subgroups.
 \end{enumerate}
\end{thm}

\begin{proof}
 (a) (\cf also \cite{efgonas}) Let
 $\rho:\;\pi_1^{\mathrm{et}}(S,\ovl{s}) \to G$ be the 
 continuous homomorphism corresponding $S'\to S$. Let $\ovl{t}$ 
 be another geometric point of $S$, and let $\ovl{\gamma}$ be an \'etale
 path from $\ovl{t}$ to $\ovl{s}$.  Under the composed continuous homomorphism
 \begin{equation*}
  \pi_1^{\mathrm{et}}(t,\ovl{t})\to \pi_1^{\mathrm{et}}(S,\ovl{t})
   \xrightarrow{\mathrm{ad}(\ovl{\gamma})}\pi_1^{\mathrm{et}}(S,\ovl{s})\to G  
 \end{equation*}
 the compact (profinite) group $\pi_1^{\mathrm{et}}(t,\ovl{t})\cong
 \Gal(\mathscr{H}(t)^{sep}/\mathscr{H}(t))$ is sent to 
 $1$ if $G$ is torsion-free. This means that for any $\ovl{t}'\in
 F_{\ovl{t}}(S')$,  
 $\mathscr{H}(t')=\mathscr{H}(t)$. Let $U$ be an open neighborhood of $t$
 (which we may assume to be contractible thanks to Berkovich's theorem),
 such that $S'\times_S U\cong \coprod V_j$, with $V_j$ finite \'etale
 over $U$. Since $\mathscr{H}(t')=\mathscr{H}(t)$ for any point $t'$ in
 $V_j$ above $t$, it follows from \ref{para::1.2.1} that
 $V_j$ is a topological covering of
 the contractible open set $U$, hence a disjoint union of finitely many
 copies of $U$. This holds 
 for any point $t\in S$, therefore $S'\to S$ is a topological covering. 

 (b) If $G$ is virtually torsion-free, let $N$ be a torsion-free normal
 subgroup of finite index, and let $S'\to T\to S$ 
 be the corresponding intermediate Galois covering. By $a)$, $S'\to T$ is
 a topological covering, hence $S'\to S$ is a temperate \'etale covering.

 (c) Assume that $S'$ corresponds to a principal homogeneous space under
 $G_S$, classified by an element of $\mathrm{H}^1(S,G)$. Let us write
 $S=\bigcup S_j$, with $S_j$ compact. Assume that $G=\bigcup G_k$, with
 $G_k$ finite. Then $\mathrm{H}^1(S_j,G)= \liminj \mathrm{H}^1(S_j, G_k)$, therefore the
 restriction of $S'$ to any $S_j$ is a disjoint union of finite \'etale
 coverings of $S_j$, hence 
 $S'\to S$ is locally algebraic. Conversely, if $S'\to S$ is Galois
 locally algebraic, its group $G$ is a discrete image of $\liminj
 \pi_1^{\mathrm{alg}}(S_j,\ovl{s})$, each
 $G_j:=\Im(\pi_1^{\mathrm{alg}}(S_j,\ovl{s}))$ is finite, and $G=\bigcup
 G_j$.
\end{proof}  

\begin{cor}\label{cor::2.1.9} In dimension one, $
 \pi_1^{\mathrm{temp}}(S,\ovl{s})$ is the pro-VTFD-completion of
 $\pi_1^{\mathrm{et}}(S,\ovl{s})$. 
\end{cor}

By pro-VTFD-completion of a topological group $\Gamma$, we mean the
inverse limit of its virtually torsion-free discrete quotients. This is
a filtered limit (if $\Gamma/H'$ and $\Gamma/H''$ are VTFD, so is
$\Gamma/(H'\cap H'')\;$). Pro-VTFD-completions are functorial (due to
the fact that any subgroup of a VTF group is VTF).

\begin{proof}
 By \ref{lem::2.1.2}, the coverings of the form $\tilde T\to T \to S$, with $T$
 finite \'etale Galois over $S$, are cofinal among temperate
 coverings. In dimension one, we have seen that the Galois groups of such
 Galois coverings are virtually torsion-free, whence the result (using
 \ref{para::1.3.2}).
\end{proof}

\medskip\noindent {\it Examples.} (i) It follows from \ref{cor::2.1.9}
that the logarithmic (locally algebraic) covering
\begin{equation*}
\log:\; \matheur{D}_{\C_p}(1,1^-)\longrightarrow  \A^1_{\C_p} 
\end{equation*}
of the affine line is {\it not }temperate, since its Galois group
$\Q_p/\Z_p$ has no non-trivial virtually torsion-free 
quotient. This also follows from \ref{pro::2.1.5}, since
$\pi_1^{\mathrm{alg}}(\A^1_{\C_p})=0$. 

\medskip\noindent {\it Remark.} It is easy to draw from \ref{thm::2.1.8} that a
connected \'etale covering which is both temperate 
and locally algebraic is (finite) algebraic.

\medskip \noindent (ii) It also follows from \ref{thm::2.1.8} that there is no
\'etale Galois covering of the punctured disk with Galois group
$\Z$. This fact has noteworthy consequences on the local monodromy of
$p$-adic differential equations (\cf next sections).

{\small
\medskip\noindent {\it Remark.} In \cite{poaptg}, we have limited
ourselves to temperate coverings for which the tower $T'\to T\to S$ 
occuring in definition \ref{def::2.1.1} is Galois with VTF group. The
fundamental attached to the category of such coverings is 
then the pro-VTFD-completion of $\pi_1^{\mathrm{et}}(S,\ovl{s})$. The
definition \ref{def::2.1.1} seems however more natural (and is much more 
general in higher dimension). In \cite{poaptg}, another topological
group is introduced, namely the coimage of 
$\pi_1^{\mathrm{et}}(S,\ovl{s})\to \pi_1^{\mathrm{top}}(S,\ovl{s})\times
\pi_1^{\mathrm{alg}}(S,\ovl{s})$, denoted by
$\pi_1^{\mathrm{red}}(S,\ovl{s})$. In dimension one, there is a continuous
injective homomorphism with dense image $\pi_1^{\mathrm{red}}(S,\ovl{s})\to
\pi_1^{\mathrm{temp}}(S,\ovl{s})$, but it not clear whether this is an
isomorphism (the question is the openness of 
$\pi_1^{\mathrm{et}}(S,\ovl{s})\to 
\pi_1^{\mathrm{temp}}(S,\ovl{s})$).
}

%
%
%

\begin{thm}\label{thm::2.1.10} Assume that $K$ is algebraically closed
 of characteristic zero. Let $\ovl{S}$ be  
 a (paracompact strictly) {\rm normal} connected $K$-analytic space
 $\ovl{S}$, and let $Z$ be a closed analytic subset such 
 that $S:=\ovl{S}\setminus Z$ is smooth.
 Assume that for any point $s$ of $\ovl{S}$, there is a covering (in the
 sense of \ref{def::1.2.3}) from a $K$-manifold 
 onto a neighborhood $U(s)$ of $s$ which is \'etale above $U(s)\cap S$. 

 \noindent
 Then any temperate \'etale covering $S'\to S $ extends to
 a (possibly ramified) covering  
 $\ovl{S}'\to \ovl{S}$, in the sense of \ref{def::1.2.3}: any point $x$
 has an open neighborhood $U(x)$ such that $S'\times_S U(x)$ 
 is a disjoint union of finite (ramified) coverings $V_j$ of
 $U(x)$. Moreover, one may assume that the $V_j$ are normal; 
 such an extension
 $\ovl{S}'\to \ovl{S}$ is then unique. 

 Furthermore, if $Z$ has codimension $\geq 2$ and $\ovl{S}$
 is smooth, any temperate
 \'etale covering $S'\to S $ extends to a unique temperate covering
 $\ovl{S}'\to \ovl{S}$.
\end{thm}

\begin{proof}
 We start with the statement about unicity. Since it is local on
 $\ovl{S}$, it is enough to establish it in the case of finite coverings
 (by definition \ref{def::1.2.3}). Let $\ovl{S}'$ and $\ovl{S}''$ be two
 connected 
 finite coverings of $\ovl{S}$ which coincide over $S$. The diagonal
 $\Delta$ is a connected component of $S'\times_S S''$. Its closure
 $\ovl{\Delta} $ in the normalization of $\ovl{S}'\times_{\ovl{S}}
 \ovl{S}''$ is again a connected component (using \cite[3.3.16]{spasalc}).
 The projections $\ovl{\Delta} \to \ovl{S}',\;\ovl{\Delta} \to
 \ovl{S}''$ are finite, and their restrictions to $\Delta$ are
 isomorphisms; hence they are themselves isomorphisms if $\ovl{S}'$ and
 $\ovl{S}''$ are normal.

 Let us turn to the existence of the extension (with the normality
 property). By unicity, this is a local question; we may and shall
 assume that $Z$ is a hypersurface, and that 
 $\ovl{S}$ is (strictly) affinoid, and that there exists a finite
 covering $\ovl{T}\to \ovl{S}$ by a 
 rigid-analytically smooth affinoid $\ovl{T}$, which restricts to a finite
 \'etale covering $T\to S$. 

 Let us first show that any quotient of a temperate covering $S'\to S$
 which has the extension property inherits this extension property: by
 unicity, the question is again local on $\ovl{S}$; by definition of
 coverings, this reduces the question to the case when $S'\to S$ is
 finite. By glueing (\cf \cite[1.3.3]{ecfnas}), it is again enough to treat
 the affinoid (normal) case. The argument is similar to the above one
 for unicity and left to the reader.
 
 This remark allows us to assume that $S'\to S$ is a topological covering of a
 finite \'etale covering of $S$. 
 Let us first treat the case of a finite \'etale covering. In this
 case, the extension property is proved in Gabber-L\"utkebohmert 
 \cite{repfapf}.
\end{proof}

The crucial case is the case of a punctured disk $S=\matheur{D}^{\ast}$
which is worth stating separately:

\begin{lem}[Gabber, L\"utkebohmert]\label{lem::2.1.11} Any finite
 \'etale covering of $\matheur{D}^\ast$ becomes a finite disjoint union
 of Kummer coverings on a smaller concentric punctured disk.
\end{lem}
\noindent (Here the assumption that $K$ is of characteristic zero is
essential). From there, one draws the extension property, in the product
case, when $\ovl{S}$ is a polydisk and $Z$ a union of axes (\cf {\it
loc.~cit.}, 3.3). By Kiehl's theorem on the existence of tubular
neighborhoods \cite{kiehl67:_rham_kohom_mannig_koerp}, the statement
also holds if $Z$ is a strict 
normal crossing divisor in $\ovl{S}$. In the general case, one can
invoke embedded resolution for $Z$, \cf \cite{schoutens99:_embed} (note that 
\cite[\S 4]{repfapf} proposes an alternative strategy which uses embedded
resolution only in dimension $2$).
 
It remains to deal at last with the case of a topological covering
$S'\to S$.  In order to avoid problems of singularities, we consider a finite 
etale covering $T\to S$ as above, $T$ being a rigid-analytically 
smooth affinoid.  It will be more convenient to replace $S'$ by 
$S'\times_S T$ (of which $S'$ is a quotient) so that the new $S'$ is 
a topological covering of $T$ 
itself a finite \'etale covering of $S$.
By \ref{pro::1.1.3}, $S'$ extends to a topological covering $\ovl{S}'
\to\ovl{T}$, and composition with $\ovl{T}\to \ovl{S}$ provides the 
desired extension of $S'/S$ to a covering $\ovl{S}'/\ovl{S}$.
 
The last statement of the theorem follows from the Nagata-Zariski purity
theorem (which becomes easy in the case of a finite morphism of
one-dimensional affinoids).

\medskip\noindent {\it Remark.} This extension property of temperate
coverings does not hold for arbitrary \'etale coverings. 

\noindent Indeed, one can show that the logarithmic covering of the
affine line does not extend to a (ramified) covering 
of the projective line.

\medskip The previous theorem is completed by the following result,
which is proven as in the complex case \cite[exp. XII, 5.3]{sga1}:

\begin{pro}\label{pro::2.1.12} Let $\ovl{S}$ be a normal analytic space,
 and let $Z$ be a closed analytic subset such that $S:=\ovl{S}\setminus
 Z$ is dense in $\ovl{S}$. Then the functor which attaches to any finite
 normal covering of $\ovl{S}$ its restriction to $S$ is fully faithful.
\end{pro}

It follows that in \ref{thm::2.1.10}, the extension of coverings from temperate
coverings of $S$ to `normal' coverings of $\ovl{S}$ `is' a fully
faithful functor.

\subsection{Interlude: tangential base-points.}\label{sub::2.2}

\begin{para}\label{para::2.2.1} Instead of basing fundamental groups at
 geometric points, it is sometimes more convenient\footnote{for instance
 in the study of Galois actions on fundamental groups.} to base them at 
 tangent vectors at points at infinity. 

 P. Deligne has explained the construction for affine algebraic curves in 
 characteristic zero, for topological, algebraic (= profinite) and other
 kinds of fundamental groups \cite{lgfdldpmtp}. In the profinite case, these
 tangential base-points amount to genuine geometric generic points,
 corresponding to embeddings of the function field of the curve into the
 field of Puiseux series at the chosen point at infinity (with respect
 to a fixed local parameter $z$ at that point; however, the embedding
 depends only on the tangent vector $\partial/\partial z$).

 If we try to imitate the construction in our non-archimedean setting, we
 face at once the following difficulty: let us assume for simplicity
 that the curve is the affine line; since there is no way to extend the
 Gauss norm on $K(z)$ to the field of Puiseux series with arbitrary coefficients in $K$, one does not get geometric (generic) points in the sense of \ref{para::1.2.2}.  We
 do not know how to overcome this difficulty for (analytic) \'etale coverings in
 general (and we doubt that this is possible). For temperate \'etale
 coverings, however, we can circumvent the difficulty as follows.
\end{para}

\begin{para}\label{para::2.2.2} We assume that $K$ is algebraically
 closed of characteristic $0$. Let $\ovl{S}$ be a $K$-manifold of
 dimension one, and let $S$ be the complement of a {\it classical} point
 $0$ in $\ovl{S}$ (so that $0$ has a basis of open neighborhoods
 consisting in disks centered at $0$). As above, it follows from
 \ref{lem::2.1.11} 
 that {\it the restriction of any temperate covering of $S$ to a
 punctured disk around $0$ of sufficiently small radius splits into a
 disjoint union of Kummer coverings of bounded degree}. If we replace
 $\ovl{S}$ by the its tangent space at $0$ and $0$ by the zero tangent
 vector, the same construction yields an equivalence of categories, with
 an obvious inverse given by `extension'. We thus get a {\it
 specialization functor}:
 \begin{equation*}
 \Cov^{\mathrm{temp}}_S\to  \Cov^{\mathrm{temp}}_{T_{\ovl{S},0}^0}.
 \end{equation*}         
 Any geometric point $\vec t$ of $T_{\ovl{S},0}^0$, 
 \eg any non-zero tangent vector at $0$, then gives rise by composition
 to a fiber functor 
 \begin{equation*}
  F^{\mathrm{temp}}_{\vec t}\;:\; \Cov^{\mathrm{temp}}_S\to  \Sets.
 \end{equation*}         
\end{para}

\begin{pro}\label{pro::2.2.3} This fiber functor is isomorphic to
 $F^{\mathrm{temp}}_{\ovl{t}}$ for any geometric point $\ovl{t}$.
\end{pro}
\begin{proof}
 We follow the line of proof of \ref{lem::1.3.3}. We may assume that
 $\ovl{S}$ is affinoid (hence compact). Let us consider finite open covers
 $\ovl{\mathcal{U}}=(\ovl{U}_i)_{i\geq 1}$ of $S$ {\it such that
 $\ovl{U}_1$ is a disk centered at $0$ and $\;0 \notin \ovl{U}_i\;$ 
 for $\;i>1$.}

 Let us denote by $ \mathcal{U}=(U_i)$ its trace on $S$, and by
 $\Cov^{\mathrm{temp}}_{S,\mathcal{U}}$ the full subcategory of 
 $\Cov^{\mathrm{temp}}_{S}$ of temperate coverings $f:S'\to S$ such that
 $f^{-1}U_i= \coprod V_{ij}$ and $V_{ij}$ is \'etale finite over 
 $U_i$ via $f$. 
 
 The point is that {\it by \ref{thm::2.1.10} (or \ref{para::2.2.2})}, we
 have
 $\displaystyle
 \Cov^{\mathrm{temp}}_{S}=\underset{\mathcal{U}}{\liminj} 
 \Cov^{\mathrm{temp}}_{S,\mathcal{U}}, $ where $\mathcal{U}$ runs among
 the open covers of $S$ satisfying the above properties. Also
 \begin{equation*}
  \Isom(F^{\mathrm{temp}}_{\vec t}, F^{\mathrm{temp}}_{\ovl{t}})
   =\underset{\mathcal{U}}{\limproj} \Isom(F^{\mathrm{temp}}_{\vec
   t,\mathcal{U}}, F^{\mathrm{temp}}_{\ovl{t},\mathcal{U}})
 \end{equation*}
 where $F^{\mathrm{temp}}_{\vec t,\mathcal{U}},
 F^{\mathrm{temp}}_{\ovl{t},\mathcal{U}}$ denote the restrictions of 
 $F^{\mathrm{temp}}_{\vec t }, F^{\mathrm{temp}}_{\ovl{t} }$ to
 $\Cov^{\mathrm{temp}}_{S,\mathcal{U}}$.

 To show that $\Isom(F^{\mathrm{temp}}_{\vec t},
 F^{\mathrm{temp}}_{\ovl{t}}) \neq \emptyset$, we construct 
 non-empty compact subsets $K_\mathcal{U}$ of
 $\Isom(F^{\mathrm{temp}}_{\vec t,\mathcal{U}},
 F^{\mathrm{temp}}_{\ovl{t},\mathcal{U}})$ (topologized 
 as before) such that if $\mathcal{U}'$ refines $\mathcal{U}$,
 $K_{\mathcal{U}'}$ maps to $K_\mathcal{U} $ under the natural map 
 $\Isom(F^{\mathrm{temp}}_{\vec t,\mathcal{U}'},
 F^{\mathrm{temp}}_{\ovl{t},\mathcal{U}'})\to
 \Isom(F^{\mathrm{temp}}_{\vec t,\mathcal{U}},
 F^{\mathrm{temp}}_{\ovl{t},\mathcal{U}})$; it will follow that
 $\Isom(F^{\mathrm{temp}}_{\vec t}, F^{\mathrm{temp}}_{\ovl{t}}) \supset
 \limproj_{\mathcal{U}}\, K_\mathcal{U}\neq 
 \emptyset$. 

 In fact, it suffices to do so for a cofinal system of
 $\mathcal{U}$'s. We now use a continuous embedding 
 $I=[0,1]\hookrightarrow \ovl{S}$ with $0\mapsto 0\in \ovl{S},\;1\mapsto t$. We
 consider open covers $\mathcal{U}  =(U_1,\ldots, U_n)$ as above, and
 such that
 \begin{equation*}
  U_1\cap I = 0,r_1[,\;\ldots, U_i\cap I=]r_i,t_i[,\;\ldots, U_m\cap
   I=]r_m,1],\;U_{m+1}\cap I=\emptyset,\ldots
 \end{equation*}
 for some points $r_i<t_{i-1}<r_{i+1}<t_i$ in $I \;(i=2,\ldots, m-1)$ and
 some $ m\leq n$.
 
 Given $\mathcal{U}$, one chooses $s_i$ in $U_i\cap U_{i+1}\cap I$ (which
 amounts to a real number between $r_{i+1}$ 
 and $t_i$), and a geometric point $\ovl{s}_i$ above $s_i\; (i=1,\ldots,
 m-1)$; one completes this collection by setting 
 $\ovl{s}_0=\vec t,\; \ovl{s}_m=\ovl{t}$. 

 Let us look at the diagram 
 \begin{equation*}
  \Cov^{\mathrm{temp}}_{S,\mathcal{U}}
   \xrightarrow{\sigma_i}
   \Ind\Cov^{\mathrm{alg}}_{\mathcal{U}_{i+1}}
   \begin{array}{c}
    \xrightarrow{F^{\mathrm{alg}}_{\ovl{s}_i, U_{i+1}}}\\
    \xrightarrow{F^{\mathrm{alg}}_{\ovl{s}_{i+1}, U_{i+1}}}\\
   \end{array}
   \longrightarrow \Sets.
 \end{equation*}               
 One has $F^{\mathrm{alg}}_{\ovl{s}_j, U_{i+1}}\circ \sigma_i =
 F^{\mathrm{temp}}_{\ovl{s}_i, \mathcal{U} },\;\;j=i,i+1$, whence
 continuous maps 
 \begin{equation*}
  \Isom(F^{\mathrm{alg}}_{\ovl{s}_i,
   U_{i+1}},F^{\mathrm{alg}}_{\ovl{s}_{i+1}, 
   U_{i+1}})\to \Isom(F^{\mathrm{temp}}_{\ovl{s}_i,
   \mathcal{U}},F^{\mathrm{temp}}_{\ovl{s}_{i+1}, \mathcal{U}}).
 \end{equation*} 
 But $\Isom(F^{\mathrm{alg}}_{\ovl{s}_i,
 U_{i+1}},F^{\mathrm{alg}}_{\ovl{s}_{i+1}, U_{i+1}})$ is a non-empty
 compact set for any $i=0,\ldots, m$: for $i=0\;$ ($s_0=\vec t$), one
 uses the analogue of \ref{pro::2.2.3} for the category of finite \'etale
 coverings, which follows from Deligne's discussion, {\it loc.~cit} (or
 even directly from \cite[V]{sga1}, interpreting
 $F^{\mathrm{alg}}_{\vec t, U_{i+1}}$ as a fiber functor at a geometric
 generic point). The rest of the proof goes exactly as in \ref{lem::1.3.3}.
\end{proof}

\begin{para}\label{para::2.2.4} This proposition allows to construct
 fundamental groups $\pi_1^\bullet(S,\,\vec t\,)$ based at a tangential
base-point $\vec t\,$ for any full subcategory $\Cov^\bullet_S\to
\Cov^{\mathrm{temp}}_S$ stable under taking connected components, fiber
products and quotients.  The theory of 1.4 extends without change to this
variant, as far as one restricts to curves over an algebraically closed
field of characteristic zero. For $\;\Cov^{\mathrm{top}}_S\to
\Cov^{\mathrm{temp}}_S$, the interpretation of tangential base-points is
not as transparent as in the complex-analytic context.
\end{para}

\subsection{Description of some temperate fundamental groups.}\label{sub::2.3}

 In this subsection, we assume that $K$ is {\it
 algebraically closed of characteristic $0$} (typically $K=\C_p$).

\begin{para}\label{para::2.3.1} It follows from \ref{pro::2.1.5}
 (together with \ref{para::1.4.1}) 
 that $\P^1$ and $\A^1$ have trivial temperate fundamental
 groups\footnote{in contrast, recall that $\A^1$ has non-trivial
 locally-algebraic fundamental group, and that $\P^1$ has trivial
 locally-algebraic fundamental group but non-trivial \'etale fundamental
 group.}.

Let us next examine $\G_m=\P^1\setminus\{0,\infty\}$. We have
$\pi_1^{\mathrm{top}}(\G_m,1)=0$; on the other hand, any connected finite
\'etale covering is Kummer $\G_m\to \G_m,\; x\mapsto x^n$. It follows that
\begin{equation*}
\pi_1^{\mathrm{temp}}(\G_m,1)=\pi_1^{\mathrm{alg}}(\G_m,1)=\widehat{\Z}(1):=
 \underset{n}{\limproj}\,\mu_n \cong \widehat{\Z}=\prod_\ell \Z_\ell.
\end{equation*}

\noindent Similarly, $\displaystyle\pi_1^{\mathrm{temp}}(\G_m^n,1)\cong
 \widehat{\Z}^n$. 
\end{para}

\begin{para}\label{para::2.3.2} {\it Elliptic curves.} Let us first consider
 the case of an {\it elliptic curve with good reduction} $S\;$ (so that
 $\pi_1^{\mathrm{top}}(S,1)=\{1\}$).  Then any finite
 \'etale covering is given by an isogeny
 $S'\to S$ and $S'$ has good reduction. Hence any temperate covering is
 an isogeny, and 
 \begin{equation*}
  \pi_1^{\mathrm{temp}}(S,1)=\pi_1^{\mathrm{alg}}(S,1)\cong \widehat{\Z}\times
   \widehat{\Z}.
 \end{equation*} 
 Let now $S=\G_m/q^\Z$ be instead a {\it Tate elliptic curve}. Then the
 homomorphisms 
 \begin{equation*}
  \pi_1^{\mathrm{temp}}(\G_m,1)\to \pi_1^{\mathrm{temp}}(S,1)\to q^\Z
 =\pi_1^{\mathrm{top}}(S,1)  
 \end{equation*}
 actually give rise to a canonical exact sequence
 (\ref{cor::1.4.12}. (b) may be applied here)
 \begin{equation*}
  1\to \widehat{\Z}(1)\to  \pi_1^{\mathrm{temp}}(S,1)\to  \Z \to 1,
\end{equation*}
 and since $ \pi_1^{\mathrm{temp}}(S,1)\subset
 \pi_1^{\mathrm{alg}}(S,1)\cong \widehat{\Z}\times \widehat{\Z}$ is
 abelian, to a non-canonical isomorphism 
 \begin{equation*}
  \pi_1^{\mathrm{temp}}(S,1)\cong \widehat{\Z}\times \Z.
 \end{equation*} 
 This is the simplest manifestation of the lacunary character of
 temperate fundamental groups.

 \medskip\noindent (The case of a Tate curve deprived from its origin is
 more difficult).
\end{para}

\begin{para}
 \label{para::2.3.3} {\it Mumford curves.} Let now $S$ be a Mumford
 curve of genus $g>1$:

 \noindent $S=\Omega/\Gamma,\;
 \;\Omega\cong \tilde S,\;$ and $\pi_1^{\mathrm{top}}(S,s)\cong \Gamma$,
 a Schottky group, free on $g$ generators.

 \noindent We have homomorphisms 
 \begin{equation*}
 \pi_1^{\mathrm{temp}}(\Omega,\bar{\tilde s})\to
  \pi_1^{\mathrm{temp}}(S,\ovl{s})\to  \pi_1^{\mathrm{top}}(S,s)
 \end{equation*}
 which give rise to an exact sequence  
 \begin{equation*}
 1\to H \to \pi_1^{\mathrm{temp}}(S,\ovl{s})\to 
 \pi_1^{\mathrm{top}}(S,s) \to 1 
 \end{equation*}
 with $H:= \Im(\pi_1^{\mathrm{temp}}(\Omega,\bar{\tilde
 s})\to \pi_1^{\mathrm{temp}}(S,\ovl{s}))^-\subset
 \pi_1^{\mathrm{temp}}(S,\ovl{s}) $. As an abstract group,
 $\pi_1^{\mathrm{temp}}(S,\ovl{s})$ is thus the semidirect product
 $H.\Gamma$.
 
 \noindent Passing to the profinite completion, we get an exact sequence
 \begin{equation*}
  \widehat{H} \to \pi_1^{\mathrm{alg}}(S,\ovl{s})\to 
 \widehat{\Gamma} \to 1, 
 \end{equation*}
 and the homomorphism
 $\pi_1^{\mathrm{alg}}(\Omega,\bar{\tilde s})\to
 \pi_1^{\mathrm{alg}}(S,\ovl{s})$ factors through 
 $\widehat{H}$. 

 \medskip This homomorphism is very far from being injective. The point
 is that the \'etale covering of $S$ obtained from an arbitrary finite
 Galois \'etale covering $\Omega'\to \Omega$ is not Galois (or, what
 amounts to the same by \ref{thm::2.1.8}, not temperate). The finite
 Galois \'etale coverings $\Omega'\to \Omega$ which give rise to a
 Galois covering of $S$ are called ``equivariant'' in \cite[1.3]{ecoamc}; 
 by \ref{thm::2.1.8}, they 
 correspond to a topological Galois covering of some Galois \'etale
 covering $S'\to S$.

 Moreover, there are examples \cite[2.7]{ecoamc} (one may choose for $\Omega$
 the Drinfeld half-plane) where the degrees of $\Omega'/
 \Omega$ and $S'/ S$ coincide, but $\Omega'$ is not the universal
 topological covering of $S'$ ; it follows that $\Omega'$ 
 is not simply-connected (in fact, its topological fundamental group is
 free of infinite type), and that the injection 
 $\pi_1^{\mathrm{temp}}(\Omega,\bar{\tilde s})\to
 \pi_1^{\mathrm{alg}}(\Omega,\bar{\tilde s})$ {\it is not a bijection in
 general}.

 Let us return to the non-injectivity of
 $\pi_1^{\mathrm{alg}}(\Omega,\bar{\tilde s})\to
 \pi_1^{\mathrm{alg}}(S,\ovl{s})$. Actually, the 
 following stronger result holds: {\it $\pi_1^{\mathrm{alg}}(S,\ovl{s})$
 is metrizable, while $\pi_1^{\mathrm{alg}}(\Omega,\bar{\tilde s})$ is 
 not}.  

 \noindent Recall that a profinite group is metrizable if and only if the
 set of its open subgroups is countable, 
 \cf \cite[1.3]{serre94:_cohom}. This is the case for
 $\pi_1^{\mathrm{alg}}(S,\ovl{s})$, which is even topologically finitely
 generated, but 
 not for $\pi_1^{\mathrm{alg}}(\Omega,\bar{\tilde s})$. In fact, one can
 show that for any 
 $n\geq 1$, the set of open subgroups $H\subset
 \pi_1^{\mathrm{alg}}(\Omega,\bar{\tilde s})$ endowed with an isomorphism
 $\pi_1^{\mathrm{alg}}(\Omega,\bar{\tilde s})/H\to \Z/n\Z$ is already
 uncountable. Indeed, this set may be identified with
 $\mathcal{O}(\Omega)^\times/(\mathcal{O}(\Omega)^\times)^n$, or else
 with the set of $\Z/n\Z$-valued 
 currents\footnote{\ie $\Z/n\Z$-valued functions on the edges of the
 tree, such that the sum of its values at the edges starting from an 
 arbitrary vertex is zero.} on the tree associated to
 $\Omega$ (\cf \cite[I.8.9, V.2.3]{garea}). Since $g\geq 2$, this 
 tree has infinitely many vertices connected to more than 
 two edges, and it follows that the $\Z/n\Z$-valued currents do not form
 a countable set.  

 A similar picture holds for Mumford curves deprived from a
 finite set of points.
\end{para}

\begin{para}\label{para::2.3.4} {\it Affinoid curves with good reduction.} In
 this case, one has 
 \begin{equation*}
  \pi_1^{\mathrm{temp}}(S,\ovl{s})=\pi_1^{\mathrm{alg}}(S,\ovl{s}) .
 \end{equation*} 
 We may assume that $s$ is the maximal point of $S$ corresponding to the
 sup-norm. The assertion then follows from De Jong's observation 
 \cite[7.5]{efgonas} that the map
 $\pi^{\mathrm{et}}_1(s,\ovl{s})=\pi^{\mathrm{alg}}_1(s,\ovl{s})\to
 \pi^{\mathrm{alg}}_1(S,\ovl{s})$ is already surjective. We shall only
 need the case of the closed unit disk. The proof goes 
 as follows (assuming that $S\subset \P^1$ as in {\it loc.~cit.}).

 One has to show that any finite \'etale covering $f:\;S'\to S$ of
 degree $n$ such that $S'_s= \{t_1,\ldots, t_n\}$ splits.  The
 assumption implies that $f$ is a local isomorphism around each $t_i$
 (\cf \ref{para::1.2.1}). Hence one can find an affinoid neighborhood $V$ of $s$ such
 that $S'_{| V}= W_1\coprod \ldots \coprod W_n$, and $W_i\cong
 V$. Because $V$ contains the maximal point $s$ which is the boundary of
 $S$ in $\P^1$ (in the sense of \cite[3.1]{staagonf}, we get an admissible
 covering (for the Grothendieck topology) of $\P^1= S \cup
 (\P^1\setminus(S\setminus V))$ by closed analytic subdomains. One then
 glues $S'$ and $n$ copies of $\P^1\setminus(S\setminus V)$ together via
 the isomorphisms $W_i\to V$ in order to get a finite \'etale covering of
 $\P^1$ which restricts to $S'$ over $S$. Such a covering is trivial.
 
 \medskip Concerning the algebraic fundamental group of the closed unit
 disk over $\C_p$, one has the following result, 
 essentially due to M. Raynaud (not used in the sequel):
\end{para}

\begin{thm}\label{thm::2.3.5} $\quad$ A finite group is a quotient of the
 profinite group $\;$ 
 $\pi_1^{\mathrm{alg}}\bigl(\matheur{D}_{\C_p}(0,1^+),0\bigr)$
 if and only if it is generated by its $p$-Sylow subgroups.
\end{thm}

 The necessity of the condition means that any finite Galois
 \'etale covering of $\matheur{D}_{\C_p}(0,1^+)$ of order prime to $p$
 is trivial, which is proven in \cite[6.3.3]{ecfnas} and
 \cite[2.11]{repfapf}.  The 
 converse is much harder and follows from Raynaud's solution of the
 Abhyankar conjecture on finite \'etale coverings of $\A^1_{\ovl{\F}_p}$,
 \cf \cite{rdldaecpecd}. Of course, we may replace $0$ by any geometric
point or tangential base-point (\cf \ref{pro::2.2.3}).

\medskip The standard example of a non-trivial finite Galois \'etale
 covering of $\matheur{D}(0,1^+)$ is the Artin-Schreier covering
 $\matheur{D}(0,1^+) \xrightarrow{z\mapsto z^p-z} \matheur{D}(0,1^+)$
 (with group $\Z/p$); it splits on any smaller concentric disk.
 A less standard 
 example, with non-solvable Galois group $A_5$, is given by $\matheur{D}(1,1^+)
 \setminus \matheur{D}(0,1^-) \xrightarrow{z\mapsto z^3-z^{-2}}
 \matheur{D}(0,1^+)$ for $p=3$.
 
\begin{para}\label{para::2.3.6} {\it Punctured disks.} Let us fix $r\in  {|
 K^\times|}$, and set $\; \matheur{D}_r=\matheur{D}_K(0,r^+),\;\;
 \matheur{D}_r^\ast=D_r\setminus \{0\}$.  We let $\ovl{s}$ be any
 geometric point or tangential base-point of $\matheur{D}^\ast$.
\end{para}

\begin{pro}\label{pro::2.3.7} There are canonical isomorphisms 
\begin{equation*}
\pi_1^{\mathrm{temp}}(\matheur{D}_r^\ast,\ovl{s})=
 \pi_1^{\mathrm{alg}}(\matheur{D}_r^\ast,\ovl{s})=
 \pi_1^{\mathrm{alg}}(\matheur{D}_r,\ovl{s})\times \widehat{\Z}(1).
\end{equation*}
\end{pro}

\begin{proof}
 By homothety, we may assume that $r=1$ (and drop the index
 $r$). These canonical isomorphisms are induced by the 
 obvious fully faithful functors (\cf \ref{cor::1.4.8}):
 \begin{equation*}
 \Cov^{\mathrm{kum}}_{\matheur{D}^\ast} \hookrightarrow
 \Cov^{\mathrm{alg}}_{\matheur{D}^\ast},\;\; 
 \Cov^{\mathrm{alg}}_{\matheur{D}} \hookrightarrow 
 \Cov^{\mathrm{alg}}_{\matheur{D}^\ast}\hookrightarrow
 \Cov^{\mathrm{temp}}_{\matheur{D}^\ast} 
 \end{equation*}
 where $\Cov^{\mathrm{kum}}_{\matheur{D}^\ast}$ denotes the category of
 finite disjoint unions of Kummer coverings of $\matheur{D}^\ast$; it is
 stable under taking connected 
 components, fiber products (the fiber product of a Kummer covering of
 degree $m$ and a Kummer covering of degree $n$ is a finite disjoint
 union of Kummer coverings of degree $\mathrm{lcm}(m,n)$) and quotients. The
 fundamental group $\pi_1^{\mathrm{kum}}(\matheur{D}^\ast,\ovl{s})$ is
 clearly $\Z(1)$.

 Assume for a while that $s$ is the maximal point of $\matheur{D}$ (the
 statement does not depend on the choice of $s$). Then as we have seen in
 \ref{para::2.3.4}, 
 the homomorphism $\pi_1^{\mathrm{temp}}(s,\ovl{s})\to
 \pi_1^{\mathrm{alg}}(\matheur{D},\ovl{s})$ is surjective. Since it
 factors through 
 $\pi_1^{\mathrm{temp}}(\matheur{D}^\ast,\ovl{s})$, it follows that
 $\pi_1^{\mathrm{temp}}(\matheur{D}^\ast,\ovl{s})\to\pi_1^{\mathrm{temp}}
 (\matheur{D},\ovl{s}) 
 = \pi_1^{\mathrm{alg}}(\matheur{D} ,\ovl{s})$ is surjective.
 
 Let us now take for $s$ a tangential base-point at $0$. Let us consider 
 the maps $\matheur{D}_r\to \matheur{D}$ (for $r\in {|K^\times|},\;r<1$) 
 and the corresponding homomorphism (\cf \ref{para::1.4.10}):
 \begin{equation*}
 \pi_1^{\mathrm{temp}}(\matheur{D}_r^\ast,\ovl{s})\to
 \pi_1^{\mathrm{temp}}(\matheur{D}^\ast,\ovl{s}).
 \end{equation*}
 It follows from \ref{lem::2.1.11} that 
 \begin{equation*}
 \underset{r\to 0}\limproj\,
  \Im(\pi_1^{\mathrm{temp}}(\matheur{D}_r^\ast,\ovl{s})\to
  \pi_1^{\mathrm{temp}}(\matheur{D}^\ast,\ovl{s})) = \Z(1),
 \end{equation*}
 the composition 

 \begin{equation*}
 \Z(1)\to \pi_1^{\mathrm{temp}}(\matheur{D}^\ast,\ovl{s})\to\Z(1)
 \end{equation*}
 being identity. It remains to show that the sequence 

 \begin{equation*}
 1\to \Z(1)\to \pi_1^{\mathrm{temp}}(\matheur{D}^\ast,\ovl{s}) \to
 \pi_1^{\mathrm{temp}}(\matheur{D} ,\ovl{s})=
 \pi_1^{\mathrm{alg}}(\matheur{D} ,\ovl{s})\to 1 
 \end{equation*}
 is exact. The problem is only in the middle; since
 $\pi_1^{\mathrm{temp}}(\matheur{D}^\ast,\ovl{s})$ is a pro-discrete group, it
 amounts to the following fact, which is a very special case of
 \ref{lem::2.1.11}: a 
 temperate Galois covering of $\matheur{D}^\ast$ extends to
 $\matheur{D}$ if and only if its 
 restriction to a sufficiently small punctured disk centered at $0$ splits.
\end{proof}

\medskip
On the other hand, the consideration of the logarithmic \'etale covering
restricted to a punctured disk centered at infinity shows that
$\pi_1^{\mathrm{loc.alg}}(\matheur{D}_r^\ast,\ovl{s})\neq
\pi_1^{\mathrm{alg}}(\matheur{D}_r^\ast,\ovl{s})$, and {\it a fortiori},
\begin{equation*}
 \pi_1^{\mathrm{et}}(\matheur{D}_r^\ast,\ovl{s})\neq
 \pi_1^{\mathrm{temp}}(\matheur{D}_r^\ast,\ovl{s}).
\end{equation*}
  
\begin{para}\label{para::2.3.8} {\it $\A^1$ minus a few points.}
 Let $S=\A^1\setminus\{\zeta_1,\ldots,\zeta_m\}$. For any $i=1,\ldots,
 m$, let $\vec t_i$ be a tangential base-point at $\zeta_i$, and let
 $\alpha_i$ be a temperate path from $\vec t_i$ to a fixed geometric
 point $\ovl{s} $ (or tangential base-point) of $ S$. Then we have a
 canonical homomorphism corresponding to Kummer coverings induced by
 $z-\zeta_i\mapsto (z-\zeta_i)^n$ (\cf \ref{para::2.3.6})

 \begin{equation*}
  \widehat{\Z}(1)\hookrightarrow \pi_1^{\mathrm{temp}}(S,\vec
  t_i)\to\widehat{\Z}(1)
 \end{equation*}
 whose composition is identity. 
 By composition with $\mathrm{ad}(\alpha_i)$, we get a monomorphism 
 \begin{equation*}
  \widehat{\Z}(1)\hookrightarrow \pi_1^{\mathrm{temp}}(S,\ovl{s}) 
 \end{equation*}
 which depends of course on $\alpha_i$.  Let us denote by $\gamma_i$ the
 image in $\pi_1^{\mathrm{temp}}(S,\ovl{s})$ of a topological generator
 of $\widehat{\Z}(1)\cong \widehat{\Z}$. The `local monodromies' $\gamma_i$ are
 not at all canonical, but the {\it smallest closed normal subgroup}
 $\an{\gamma_i}^-$ of $\pi_1^{\mathrm{temp}}(S,\ovl{s}) $ which contains
 them is independent of all choices. In fact we have\footnote{this is
 the correct version of cor. 5.3 in \cite{poaptg}, which was wrongly
 stated.}: 
\end{para}

\begin{pro}\label{pro::2.3.9}
 $\an{\gamma_i}^-=
 \pi_1^{\mathrm{temp}}(\A^1\setminus\{\zeta_1,\ldots, \zeta_m\},\ovl{s})$.
\end{pro}

\begin{proof}
 This amounts to saying that a temperate (Galois, if one wishes)
 covering splits if its restriction to sufficiently small punctured disks
 around any of the $\zeta_i$ splits. But the latter condition implies
 that the covering extends to a temperate covering of $\A^1$ 
 (\cf \ref{thm::2.1.10}), which automatically splits (\cf \ref{para::2.3.1}).
\end{proof}

\medskip By analogy with the complex situation, \ref{pro::2.3.9}
suggests two questions:

\begin{que}\label{que::2.3.10} Is it possible to choose the $\gamma_i$'s
 in such a way that they are topological generators of
 $\pi_1^{\mathrm{temp}}(S,\ovl{s})$ ? 
\end{que}

Following the above construction, one can introduce also introduce a
local monodromy $\gamma_\infty$ at $\infty$. 

\begin{que}\label{que::2.3.11} Is it possible to choose the $\gamma_i$'s
 and $\gamma_\infty$ in such a way that their product (in suitable
 order) is identity?
\end{que}

We do not know the answer to \ref{que::2.3.10}, but shall answer 
\ref{que::2.3.11} in \ref{sec::6}, when we have more tools at disposal.

For $m\geq 2$, $\pi_1^{\mathrm{temp}}(S,\ovl{s})$ is indeed a rather
complicated topological group, far from being discrete and from being compact:

\begin{pro}\label{pro::2.3.12} For $m\geq 2$,
 $\pi_1^{\mathrm{temp}}(\A^1\setminus\{\zeta_1,\ldots,
 \zeta_m\},\ovl{s}) $ is not locally compact.
\end{pro}

\begin{proof}
 Assume that $\pi_1^{\mathrm{temp}}(S,\ovl{s}) $ is
 locally compact. Then because $\pi_1^{\mathrm{temp}}(S,\ovl{s})$ is a
 prodiscrete group, there would exist a {\it compact} (necessarily
 totally discontinuous, hence profinite) open normal subgroup $H$. Let
 $S'\to S$ be the corresponding Galois temperate covering with group
 $G=\pi_1^{\mathrm{temp}}(S,\ovl{s})/H$. Up to replacing $H$ by a small
 subgroup, we may and shall assume that $S'=\tilde T$ is the topological
 universal covering of a finite Galois \'etale covering $T$ of $S$. Via
 \ref{pro::2.1.7}, we then have an exact sequence
 \begin{equation*}
 1\to H\to \pi_1^{\mathrm{temp}}(T,\ovl{t})\to \pi_1^{\mathrm{top}}(T,
 \ovl{t} )\to 1. 
 \end{equation*}
 Let $T'$ be any finite \'etale covering of $T$. Then the universal
 topological covering $\tilde T'$ is a temperate covering of $T$ which
 dominates $S'=\tilde T$. Due to the compactness of $H$, the map $\tilde 
 T' \to \tilde T$ must be finite. 

 We now choose $T'$ in order to contradict this property. Let $\ovl{X}$
 be a projective smooth analytic curve over $K$ such that $\rk
 \pi_1^{\mathrm{top}}(\ovl{X}) > \rk \pi_1^{\mathrm{top}}(T)$. Assume
 that $\ovl{X}$ is endowed with a finite morphism to $\P^1$ which is
 unramified above $S$ (since $m\geq 2$, this just amounts to saying that
 $\ovl{X}$ is defined over $\ovl{\Q}$, according to a well-known theorem
 of Belyi). We set $X=\ovl{X}\times_{\P^1} S \subset \ovl{X}$, and choose
 for $T'$ a connected component of $ X\times_S T$. Note that since $T$ is
 Galois over $S$, $T'$ is Galois over $X$. Let $Y$ be a connected
 component of $X\times_S \tilde T$ lying above $X$. Then $\tilde T' $ is
 a topological covering of $Y$, and is finite over $Y$ (if and) only if 
 it coincides with $Y$.  We then have a commutative square of Galois 
 temperate coverings
 \begin{equation*}
  \begin{matrix}
   && Y = \tilde T'&&\\
   &\swarrow &&\searrow &\\
   \tilde X &&&& T' \\
   &\searrow &&
   \swarrow &\\
   &&X&&
  \end{matrix}
 \end{equation*}
 Now $\Aut_X Y$ admits the free group $\Aut_{T'} Y$ as a normal subgroup
 of finite index, and admits the free group $\Aut_X \tilde X$ as a
 quotient. Since $\rank \Aut_{T'} Y \leq \rank \pi_1^{\mathrm{top}}(T)<
 \rank \pi_1^{\mathrm{top}}(\ovl{X})=\rank\,\pi_1^{\mathrm{top}}(X)$, this is
 impossible. Therefore $\pi_1^{\mathrm{temp}}(S) $ is not locally compact.  
\end{proof}

 \medskip\noindent {\it Remark.} In \cite{andre:_une_gal}, temperate fundamental
 groups of algebraic $p$-adic manifolds will be used to give a geometric 
 description of the local Galois group $\Gal(\ovl{\Q}_p/\Q_p)$. 
 
\newpage

\section{Local and global monodromy of $p$-adic differential
 equations.}\label{sec::3}

\markboth{\thechapter. $p$-ADIC ORBIFOLDS AND MONODROMY.}%
{\thesection. LOCAL AND GLOBAL MONODROMY.}

\begin{abst}
We briefly review and compare the theory of singularities of ordinary
linear differential equations in the formal, complex and $p$-adic
contexts. We then outline the Christol-Mebkhout theory of $p$-adic
slopes of differential modules over annuli. There is a close relation
with Galois representations of local fields of characteristic $p$,
expressed by Crew's $p$-adic local monodromy conjecture, which has been
recently proved \cite{andre01:_filtr_hasse_arf}. Complications occuring in the $p$-adic
context reflect the fact that the \'etale fundamental groups of annuli
are non-abelian.  \hfill\break We then define and study the
non-archimedean Riemann-Hilbert functor, which attaches a vector bundle
with integrable connection to any continuous representation with
discrete coimage of the \'etale fundamental group of a non-archimedean
manifold.  Connections in the image are characterized by the fact that
the \'etale sheaf of germs of solutions is locally constant.
\end{abst}

\subsection{Introduction.}\label{sub::3.1}

\begin{para}\label{para::3.1.1} Let us briefly review some well-known facts
 about singularities of ordinary linear differential equations. For
 concreteness, we consider a differential equation   

 \begin{equation}
  \alpha_\mu(z) \partial^\mu y + \cdots + \alpha_0 (z)
   y =0,\;\; \partial=z\frac{d}{dz}, \tag{$\ast$}
 \end{equation}
 with polynomial coefficients $\alpha_i(z)$, and concentrate on
 phenomena at infinity.
\end{para}

\begin{para}\label{para::3.1.2}  We begin the discussion with the formal
 setting. One defines the formal Newton polygon at $\infty$,
 $\widehat{NP}_\infty$ to be the convex hull of $\{(X\leq i,Y\geq \deg_z
 \alpha_0 -\deg_z \alpha_i )\}$ in the real 
 plane. Its `height'\footnote{the 
 height is the difference between the $Y$-coordinate of the highest 
 vertex and the $Y$-coordinate of the lowest vertex.}
 is called the (formal) irregularity of $(\ast)$ at
 $\infty$. The equation is called regular at $\infty$ if the
 irregularity is $0$, \ie if the only finite slope is $0$.

\medskip It is known that $(\ast)$ has a basis of `formal' solutions of
 the form 

\begin{equation*}
 \hat y=\sum {\hat u_i}\Bigl(\frac{1}{z}\Bigr) z^{-{\hat e}_i}(\log^{k_i}
  z) e^{P_i(z^{1/\mu !})},
\end{equation*}
 where the $P_i$'s are polynomials, and the ${\hat u_i}$ are just formal
 power series in $\frac{1}{z}$. The degree of the `ramified polynomials'
 $P_i(z^{1/\mu !}) $ actually corresponds to the slopes of $\widehat 
{NP}_\infty$. The ${\hat e}_i$'s are called the formal or Turrittin
 exponents; they can be computed algebraically.

 \medskip
 \noindent {\it Example.} Let us illustrate these notions with the {\it
 Bessel equation}
 \begin{equation*}
  (B)_{\nu}:\;\;\;\;(\partial^2  + (z^2  - \nu^2))y=0 
 \end{equation*}
 A basis of `formal solutions' at $\infty$ is given by
 \begin{equation*}
  \hat y_\pm= z^{-1/2}e^{\pm\sqrt{-1}z}.
   {}_2F_0\Bigl(\frac{\nu +1}{2},\frac{-\nu +1}{2}; \frac{\mp
   \sqrt{-1}}{2z}\Bigr), \quad {\hat e}_i=1/2,
 \end{equation*}
 and $\widehat{NP}_\infty$ is
 
 \begin{figure}[h]
  \begin{picture}(100,100)(0,0)
   \put(0,0){\line(1,0){40}}
   \put(40,0){\line(1,1){40}}
   \put(40,0){\makebox(0,0)[c]{{\small $\bullet$}}}
   \put(60,0){\makebox(0,0)[c]{{\small $\bullet$}}}
   \put(80,0){\makebox(0,0)[c]{{\small $\bullet$}}}
   \put(80,20){\makebox(0,0)[c]{{\small $\bullet$}}}
   \put(80,40){\makebox(0,0)[c]{{\small $\bullet$}}}
   \put(80,40){\line(0,1){40}}
  \end{picture}
  \caption{}
 \end{figure}
\end{para}

\begin{para}\label{para::3.1.3} {\it Complex-analytic setting}.  Solutions at
 any non-singular point $s$ (\ie a point in $\C^\times$ which is not a
 root of $\alpha_\mu$) converge up to the next singularity.

 Let ${\mathcal{C}}= \mathcal{C}(]r',r[) =\{z, \;r'<|z| < r\}$ be an
 annulus which contains no singularity. Using the fact the
 $\pi_1(\mathcal{C},s)\cong \Z$, it is classical to deduce that $(\ast)$
 has a basis of `analytic' solutions of the form

 \begin{equation*}
  y=\sum {u_i}\Bigl(\frac{1}{z}\Bigr) z^{-{e}_i} (\log^{k_i} z),
 \end{equation*}
 where the ${u_i}$ are Laurent series in $\frac{1}{z}$
 which converge in $\mathcal{C}$.
 The ${ e}_i$'s are called the analytic exponents.\footnote{in terms of
 differential systems, this amounts to the fact that a system
 $\partial Y = A Y$ on $\mathcal{C}$ has a fundamental solution matrix
 $Y$ of the form 
 $U.z^{-E}$, where $U$ is analytic in $\mathcal{C}$ and $E$ is constant (Fuchs
 normal form).}

 \medskip
 The relation between the $y$'s (analytic) and
 the $\hat y$'s (formal) is not straightforward. If $r=\infty$, the $\hat
 y$'s occur as asymptotic expansions of the $y$'s on suitable sectors,
 and the theory of multisummability/acceleration provides a kind of
 canonical inverse process (again on suitable sectors). If $\infty$ is a
 regular singularity, the $\hat u_i$ converge so that with $u_i=\hat u_i,
 e_i=\hat e_i$, the solutions $y$ and $\hat y$ coincide.

 \medskip
 \noindent {\it Example.} For the Bessel equation (say with non-integral
 parameter $\nu$ to simplify), one has a basis of solutions on
 $\mathcal{C}=\C^\times$ given by
  
\begin{equation*}
 y_\pm= z^{\pm\nu}\sum_0^\infty \;\frac{\;(-1)^n z^{2n }}{4^n n! (\pm\nu
  +1)_n}. 
\end{equation*}
 Moreover $e_i = \pm \nu$ (analytic exponents at $\infty = -$ analytic
 exponents at $0$); here we are `lucky', there are just two
 singularities, one of them is regular. In general, the computation of
 analytic exponents involves infinite determinants. Even in the
 so-called fuchsian case where there are only regular singularities, the
 computation of analytic exponents on an annulus which encloses several
 regular singularities is transcendental; e.g. the exponents of the
 differential equation of order two with solution
 ${}_2F_{1}(\frac{1}{2},\frac{1}{2},1;z^2)$, on an 
 annulus enclosing $0$ and $1$ but not $-1$ and
 $\infty$, are transcendental numbers:
 $\frac{\log(3\pm 2\sqrt 2)}{2\pi i}$\footnote{An independent
 solution is given by $i. \,{}_2F_{1}(\frac{1}{2}, \frac{1}{2}, 1;
 1-z^2)$, and in this basis, the local monodromy matrix 
 around $1$ (\resp $0$) is
 $\bigl(\begin{smallmatrix}1&0\\ -2&1\end{smallmatrix}\bigr)$
 (\resp $\bigl(\begin{smallmatrix}1&4\\0&1\end{smallmatrix}\bigr)$); the
 eigenvalues of the 
 monodromy around $\mathcal{C}$ are the eigenvalues of 
 $\bigl(\begin{smallmatrix}1&0\\-2&1\end{smallmatrix}\bigr)
  \bigl(\begin{smallmatrix}1&4\\0&1\end{smallmatrix}\bigr)
  =\bigl(\begin{smallmatrix}-7&4\\-2&1\end{smallmatrix}\bigr)$, 
 \ie $3\pm 2\sqrt 2$.}.
\end{para}

\begin{para}\label{para::3.1.4}  {\it $p$-adic-analytic setting}.  The first
 problem one encounters when dealing with $p$-adic differential
 equations (\cf I.\ref{sub-r-and-unit-root}) is that the solutions at
 non-singular points have small radii of convergence in general, and do
 not converge up to the next singularity.
 
 \medskip
 \noindent {\it Example.}
 $\frac{d}{dz}y=y,\;  e^z=\sum \frac{z^n}{n!} $ converges only in
 $\matheur{D}(0,|\pi|^-)$, $\pi^{p-1}=-p$. 

 \medskip
 This problem leads one to normalize somehow the situation by postulating 

 \medskip\noindent {\it Dwork's condition: the radius of convergence of
 solutions of $(\ast)$ at a generic point\footnote{a generic point of
 modulus $\rho$ in the sense of Dwork is a point $t_\rho$ in a complete
 extension $\Omega$ of $K$ with $| t_\rho|=\rho$ and such that the disk
 $\matheur{D}_\Omega(t_\rho,\rho^-)$ does not contain any element of $K$.}
 $\;t_1$ of modulus $1$ is $\geq 1$}.

 \medskip
 \noindent {\it Example.} One `normalizes' Bessel's equation by the
 change of variable $ z\mapsto 2\pi\sqrt{-1} z$. The normalized 
 Bessel differential operator is then   
 \begin{equation*}
  B_{\nu,\pi} = \partial^2  - (4\pi^2 z^2 +\nu^2)  
 \end{equation*}
 which is solvable in the $p$-adic generic disk whenever $\nu \in \Z_p$.
\end{para}

\begin{para}\label{para::3.1.5} Let us now discuss problems at $p$-adic
 singularities, concentrating again on phenomena at infinity.

 \medskip\noindent $\bullet$ The first one is the problem of `small
 divisors': if some formal exponent $\hat e_i\in \Z_p$ (or rather some
 difference $\hat e_i -\hat e_j$) is a Liouville number, \ie very
 closely approximated by infinitely many rational integers, the formal
 power series $\hat u_i$ occuring in the formal solutions $\hat y$ may
 diverge, even in the regular case. Otherwise, the $\hat u_i$ 
 converge, even in the irregular case \cite{baldassarri82:_differ}.

 \medskip
 \noindent {\it Example.} if $\nu \in \Z_p$ is well-approximable by
 negative integers, $\sum_0^\infty \frac{(-1)^n z^{2n}}{4^n n! (\nu
 +1)_n}$ may diverge, although $0$ is a regular singularity of the
 Bessel equation.

\medskip\noindent 
 $\bullet$ Even if Liouville numbers do not show up (\eg if everything
 is defined over $\ovl{\Q}$) so that the $\hat u_i$'s converge, these
 series may not converge as far as one might expect (\eg in a unit disk
 under the Dwork condition). This subtle phenomenon occurs in the
 normalized Bessel case $B_{\nu,\pi}$ for $p=2$ (${}_2F_0(\frac{\nu
 +1}{2},\frac{-\nu +1}{2}; \frac{\pm 1}{4\pi z})$ do not converge up
 to the unit circle): there is an annulus $\mathcal{C}(]1,r[)$ such that
 neither $y_\pm$ nor $\hat y_\pm$ provide any information about the
 solutions of $B_{\nu,\pi}$ in $\mathcal{C}(]1,r[)$.

 In particular, the question arises whether there is a basis of solutions 

 \begin{equation*}
  \sum{u_i}\Bigl(\frac{1}{z}\Bigr) z^{-{e}_i}(\log^{k_i} z),
 \end{equation*}
 where the ${u_i}$ are Laurent series in $\frac{1}{z}$ which converge in
 such an annulus $\mathcal{C}(]1,r[)$. More generally, concerning the
 existence of a `$p$-adic Fuchs normal form': 

 \medskip\noindent 
 $\bullet$ Does there exist a basis of solutions of $(\ast)$ of the form 
 \begin{equation*}
  y=\sum {u_i}\Bigl(\frac{1}{z}\Bigr) z^{-{e}_i} (\log^{k_i} z),
 \end{equation*}
 where the ${u_i}$ are Laurent series which converge in
 a given annulus $\mathcal{C}(]r',r[)$?

 \medskip\noindent A necessary condition has been put forward by
 P. Robba: {\it for any $\rho\in ]r',r[$, the radius of convergence of
 the solutions at the generic point $t_\rho$ of modulus $\rho$ is $\geq
 \rho$} (in the variable $z$, not $1/z$).

 The problems are: to which extent is the converse true? What can be
 said if Robba's condition is not fulfilled?

 \medskip\noindent These questions have recently become almost complete
 answers thanks to the Christol-Mebkhout theory of $p$-adic slopes of
 differential modules over annuli.  
\end{para}

\subsection{Local monodromy of $p$-adic differential equations. Outline
 of Christol-Mebkhout theory.}\label{sub::3.2}

\begin{para}\label{para::3.2.1} Let $K$ be a complete subfield of
 $\C_p$. We assume for simplicity that the restriction of the $p$-adic 
valuation to $K$ is discrete. 

The Robba $K$-algebra at `infinity' is 
 \begin{equation*}
\mathcal{R}=\mathcal{R}_{z,K}= \{K\text{-analytic functions on some }
 \mathcal{C}(]1,r[)\}.
\end{equation*}
 Endowed with $\partial = zd/dz$, this is a differential (integral)
 $K$-algebra. Although it is not simple as a differential algebra (there
 are many non-trivial differential ideals), it has the following useful
 property, which often allows one to replace $\mathcal{R}$ by its
 fraction field $Q(\mathcal{R})$ (\cf \cite{andre:_repres_bessel}.): 

\begin{quote}
 {\it any $\mathcal{R}[\partial]$-module $M$ which is of finite
 presentation over $\mathcal{R}$ is free over $\mathcal{R}$.} 
\end{quote}
 In the sequel, we simply refer to such an $M$ as a
 {\it differential module} over $\mathcal{R}$.
\end{para}

\begin{para}\label{para::3.2.2} For any $\rho >1$ sufficiently close to $1$,
 let us denote by $R(M,\rho)$ the infinimum of $\rho$ and the maximal
 radius $R$ such that $M$ is solvable (= has a basis of analytic
 solutions) in $\matheur{D}(t_\rho, R^-)$.

 One says that  $M$ is  {\it solvable} (at the inner boundary), if  
 $\displaystyle\lim_{\rho \to 1^+}R(M,\rho) = 1$ 

 \noindent(this is essentially Dwork's condition, except that $R(M,1)$
 itself is not defined.)
\end{para}

\begin{para}\label{para::3.2.3} A useful criterion (which goes back the
 Dwork's early studies) for $M$ being solvable is the existence 
 of a Frobenius structure. 

 Let $\sigma$ be an automorphism of $K$ inducing $t\mapsto
 t^{p^f}$ in characteristic $p$, and let $\varphi_f$ be the
 $\sigma$-linear endomorphism of $\mathcal{R}$ defined by $z\mapsto
 z^{p^f}$.

 One says that $M$ {\it has a Frobenius structure} if $\varphi_f^\ast
 M\cong M $. In what follows, the particular choice of such an
 isomorphism is irrelevant.

 \medskip\noindent
 {\it Example.} Let
 $M_{\nu,\pi}=\mathcal{R}/\mathcal{R}B_{\nu,\pi}$ be the differential
 module over $\mathcal{R}$ attached 
 to the normalized Bessel operator. It turns out that, up to isomorphism,
 $M_{\nu,\pi}$ does not depend on $\nu\in \Z_p$ 
 \cite{andre:_repres_bessel}. Choosing $\nu=0$ for instance, the existence of a
 Frobenius structure is known by the interpretation of $M_{\nu,\pi}$ in 
 terms of Kloosterman sums \cite{bfapfota} (in fact, for any $\nu \in \Q\cap
 \Z_p$, there is a relation to twisted Kloosterman sums, and a Frobenius
 structure has a geometric interpretation, \cf \cite{berthelot84:_cohom_dwork}).
\end{para}

\begin{para}\label{para::3.2.4} We know turn to $p$-adic slopes\footnote{which
 should not be confused with the so-called slopes of Frobenius which
 measure the $p$-adic size of a fixed Frobenius structure (when it
 exists); a unit-root $F$-crystal corresponding to the case of Frobenius
 slope zero.}, as defined by G. Christol and Z. Mebkhout.

\begin{itemize}
 \item A solvable $M$ is of {\it greatest slope} $\lambda \in
       [0,\infty[$ if $R(M,\rho) = \rho^{1-\lambda}$ for $\rho $
       close to $1$.
 \item $M$ is {\it purely of slope} $\lambda$ if any non-zero 
       solution at $t_\rho$ has radius of convergence $\rho^{1-\lambda}$ for
       $\rho $ close to $1$.
\end{itemize}
\end{para}

\begin{thm}\label{thm::3.2.5} Any solvable $M$ has a decreasing
 filtration $M_{>\lambda}$ such that $Gr_\lambda\; M$ is free and 
 purely of slope $\lambda$.  
\end{thm}
\par \noindent \cf \cite{sltdldedp3}. Following Christol-Mebkhout, one
defines the {\it $p$-adic Newton polygon} 
(at infinity) $NP_\infty$: the horizontal length of the segment of slope
$\lambda$ is  $rk \;Gr_\lambda\; M $.  In particular, the greatest slope
of $M$ in the above sense is the greatest slope of $NP_\infty$ in the
usual sense. The height of $NP_\infty$ is called the $p$-adic
irregularity of $M$. 

A solvable module $M$ is called {\it tame} (or $p$-adically regular) if
the only finite slope of $NP_\infty$ is zero. 
This amounts to say that $M$ satisfies the Robba condition of
\ref{para::3.1.5}: $R(M,\rho)=\rho$ when $\rho\to 1^+$.

\medskip One of the main results in the Christol-Mebkhout theory is
(\cite{sltdldedp3}): 
 
\begin{thm}\label{thm::3.2.6}
 $NP_\infty$ has integral vertices. 
\end{thm}
 
\par \noindent  A useful complement in the computation of $NP_\infty$
is given by (\cite{andre:_repres_bessel}):

\begin{lem}\label{lem::3.2.7} Assume $M$ comes from a differential module
 over the subring of $(K\cap \ovl{\Q})((\frac{1}{x}))$ consisting of
 convergent series in  $\mathcal{C}(]1,\infty[)$. Then the $p$-adic
 greatest slope $\lambda$ is $\leq$ than the greatest formal slope
 $\lambda_\infty$.
\end{lem}

 \noindent
 {\it Example.} It can be shown (using \ref{thm::3.2.8} below) that the
 Bessel module $M_{\nu,\pi}$ is not tame; \ref{thm::3.2.6} and
 \ref{lem::3.2.7} then leave only two possibilities for 
 $NP_\infty$, and it turns out that both occur  (\cf
 \cite{andre:_repres_bessel}).

 \begin{figure}[h]
  \begin{picture}(300,125)(0,-25)
   \put(0,0){\line(1,0){40}}
   \put(40,0){\line(2,1){40}}
   \put(40,0){\makebox(0,0)[c]{{\small $\bullet$}}}
   \put(60,0){\makebox(0,0)[c]{{\small $\bullet$}}}
   \put(80,0){\makebox(0,0)[c]{{\small $\bullet$}}}
   \put(80,20){\makebox(0,0)[c]{{\small $\bullet$}}}
   \put(80,20){\line(0,1){60}}
   \put(55,20){\makebox(0,0)[r]{{\footnotesize $\lambda=1/2$}}}
   \put(40,-20){\makebox(0,0)[c]{$p=2$}}

   \put(200,0){\line(1,0){40}}
   \put(240,0){\line(1,1){40}}
   \put(240,0){\makebox(0,0)[c]{{\small $\bullet$}}}
   \put(260,0){\makebox(0,0)[c]{{\small $\bullet$}}}
   \put(280,0){\makebox(0,0)[c]{{\small $\bullet$}}}
   \put(280,20){\makebox(0,0)[c]{{\small $\bullet$}}}
   \put(280,40){\makebox(0,0)[c]{{\small $\bullet$}}}
   \put(280,40){\line(0,1){40}}
   \put(255,25){\makebox(0,0)[r]{{\footnotesize $\lambda=1$}}}
   \put(240,-20){\makebox(0,0)[c]{$p>2$}}
  \end{picture}
  \caption{}
 \end{figure}
 \newpage
 \noindent
 As for the converse of Robba's criterion (\cf \ref{para::3.1.5}), one 
 has (\cf \cite{sltdldedp2}): 
\begin{thm}[$p$-adic local monodromy theorem, tame
 case]\label{thm::3.2.8}
 Assume $M$ has a Frobenius structure\footnote{there is a more 
 general version without Frobenius structure, in which one encounters 
 some trouble with the problem of small divisors. The existence of a
 Frobenius structure simplifies the statement, if not the proof, since
 it implies a priori that if $p$-adic exponents exist, they are
 rational, hence non-Liouville.}, and is purely of slope zero. Then 
 $M$ is an iterated extension of rank one differential modules of type
 $\partial-\alpha\;$, $\alpha\in \Q\cap \Z_p$ (called the $p$-{\rm adic
 exponents} of $M$).
 
 With respect to the conclusion about the rationality of the $p$-adic
 exponents, this result may be seen as a $p$-adic analogue of
 Grothendieck's $\ell$-adic local monodromy theorem for tame $\ell$-adic
 coefficients.
\end{thm}
The result, if not its proof, is closely related to the fact that a tame
(= of degree prime to $p$) Galois \'etale finite 
covering of any annulus $\mathcal{C}$ is Kummer (the tame algebraic
fundamental group of $\mathcal{C}$ is $\widehat{\Z}'(1)$, where the 
dash means as usual that the $p$-component is omitted,
\cf \cite{ecfnas}). The relationship will appear more clearly below.

\begin{para}\label{para::3.2.9} Let us return to our original differential
 equation $(\ast)$. Let us assume that the only singularities at finite
 distance are in $\matheur{D}(0,1^+)$.  Let us 
 also assume that $(\ast)$ satisfies Dwork's condition (so that the
 associated differential module over $\mathcal{R}$ is solvable), and that
 the formal exponents at $\infty$ are rational. Then the $p$-adic
 behavior of $(\ast)$ in a small punctured disk $\matheur{D}(\infty,
 \epsilon)^\ast$ centered at $\infty$ is dictated by the formal behavior
 (\cf \ref{para::3.1.5}), while the $p$-adic behavior in a thin annulus
 $\mathcal{C}(]1,r[)$ is predicted by Christol-Mebkhout theory. However,
 in the intermediate region $\mathcal{C}([r,1/\epsilon])$, nothing seems
 to be known.

 We shall return to this problem in 3.5. 
\end{para}

\subsection{Local monodromy of $p$-adic differential equations and
  Galois representations.}\label{sub::3.3}

\begin{para}\label{para::3.3.1} The differential sub-$K$-algebra 
 \begin{equation*}
  \mathcal{E}^{\dag} = \{\text{bounded $K$-analytic functions on some }
   \mathcal{C}(]1,r[)\}
 \end{equation*} 
 of bounded functions in $\mathcal{R}$ is actually a field. With the
 $p$-adic sup-norm, it is even a henselian field, with residue field
 $k((\frac{1}{z}))$ (where $k$ stands for the residue field of $K$).
 Hence any Galois extension $k((\frac{1}{z'}))/k((\frac{1}{z}))$ gives
 rise to an unramified Galois extension
 $\mathcal{E}^\dagger_{z',K'}/\mathcal{E}^\dagger_{z,K}$ (and, further,
 to a finite integral extension
 $\mathcal{R}_{z',K'}/\mathcal{R}_{z,K}$).
\end{para}

\begin{para}\label{para::3.3.2} In I.\ref{para-unit-root}, we have
 reviewed Katz' functor 
 from representations of fundamental groups of algebraic curves in
 characteristic $p$ to unit-root $F$-(iso)crystals. This construction has
 a local analogue, and one has the following theorem of N. Tsuzuki 
 \cite{tsuzuki98:_finit_f}:
\end{para}

\begin{thm}\label{thm::3.3.3} There is a natural equivalence of
 categories between the category of continuous $Q(W(k))$-linear 
 representations of
 $\Gal(k((\frac{1}{z}))^{\mathrm{sep}}/k((\frac{1}{z})))$ with 
 {\rm finite inertia} and the category of overconvergent unit-root
 $F$-isocrystals on $\Spec\;k((\frac{1}{z}))$.
\end{thm}

\par\noindent In particular, disregarding the precise Frobenius
structure, one gets a functor
\begin{equation*}
 \left\{\begin{matrix}
    \text{$p$-adic representations of $k((\frac{1}{z}))$}\\
    \text{with finite inertia}
   \end{matrix}
   \right\}
   \to
   \left\{\begin{matrix}
      \text{differential modules over $\mathcal{E}^\dagger$}\\
      \text{with Frobenius structure}
     \end{matrix}
     \right\}
\end{equation*}
hence, loosing even more structure, a functor

\begin{align*}
 &\bigl\{\begin{matrix}
	 \text{finite $p$-adic representations of $k((\frac{1}{z}))$}
	\end{matrix}\bigr\}\\
 &\qquad\qquad\qquad \to
 \left\{
 \begin{matrix}
  \text{differential modules over $ \mathcal{R}_{z,K}$ with Frobenius
  }\\
  \text{structure, which become trivial over some 
  }\\
  \text{finite extension $\mathcal{R}_{z',K'}$ of $\mathcal{R}_{z,K}$ as above}
 \end{matrix}
 \right\}
\end{align*}

\begin{dfn}\label{def::3.3.4} (R. Crew) A differential module $M$ over
$\mathcal{R}$ is {\rm quasi-unipotent} if there is a Galois extension
$k((\frac{1}{z'}))/k((\frac{1}{z}))$ such that after tensoring with the
corresponding finite extension $\mathcal{R}_{z',K'}$ of
$\mathcal{R}_{z,K} $, $M$ becomes an iterated extension of trivial
differential modules over $\mathcal{R}_{z',K'}$.

 There is an
 extensive work by Crew, S. Matsuda and Tsuzuki relating Galois
 representation theory and ramification theory of $k((\frac{1}{z})))$ and
 Christol-Mebkhout theory in the case of quasi-unipotent differential
 modules, \cf \cite{ceiatsc}.  In particular, they 
 establish a dictionary between Swan conductor and $p$-adic irregularity
 (in the above sense), break decomposition (in the sense of Katz) and
 $p$-adic slope decomposition... 
\end{dfn}

 \medskip\noindent {\it Examples.} (i) A {\it tame} module $M$
 with Frobenius structure is quasi-unipotent by \ref{thm::3.2.8}; indeed, it
 suffices to take a Kummer extension
 $k((\frac{1}{z^{1/n}}))/k((\frac{1}{z}))$. 

 (ii) $\mathcal{R}/\mathcal{R}(\partial-\pi z)$ is
 quasi-unipotent: take an Artin-Schreier extension defined by
 $z'-z^{\prime p}=z$.

 (iii) $p$-adic Bessel $M_{\nu,\pi}$ for odd $p$: a cyclic
 extension of $k((\frac{1}{z}))$ of degree $2p$ of 
 suffices. Indeed, it suffices to deal with the case $\nu = 1/2$ (\cf
\ref{para::3.2.3}). Then $(\frac{1}{z})^{\pm 1/2}.e^{\pm 2\pi z}$ is a
basis of solutions, and one reduces to the previous two examples.

\begin{coj}[Crew's $p$-adic local monodromy conjecture]\label{coj::3.3.5}
 Any differential module $M$ over $\mathcal{R}$, with Frobenius structure, is
 quasi-unipotent.
\end{coj}

This conjecture, which can be seen as a $p$-adic analogue of
Grothendieck's $\ell$-adic local monodromy theorem for wild $\ell$-adic
coefficients, and as a {\it non-abelian} extension of
\ref{thm::3.2.8} to the wild case, has now been solved
\cite{andre01:_filtr_hasse_arf}.
The proof relies heavily 
on theorems \ref{thm::3.2.5}, \ref{thm::3.2.6}, \ref{thm::3.2.8}.

\medskip\noindent A crucial test \cite{andre:_repres_bessel}, which is
at the origin 
of the solution, is the case of the {\it diadic} normalized Bessel
module $M_{\nu,2}$ (this test was suggested by Mebkhout as a potential
counterexample). This example provides a typical non-abelian instance of
\ref{coj::3.3.5}:

\begin{thm}\label{thm::3.3.6} For $K=\Q_4$, the normalized Bessel
 differential $\mathcal{E}^\dagger_{x,K}$-module $M$ is 
 quasi-unipotent.
\end{thm}

\medskip More precisely, let $E$ be the supersingular elliptic curve
 over $\F_4$ with affine equation 
  $Y^2+Y=X^3+\omega,\;$ ($\omega^3=1, \omega\neq 1$); then the covering
 $E\to E/\Aut(E)\cong \P^1\;$ is wildly totally ramified at the origin
 of $E$, and defines an unramified Galois extension of
 $\mathcal{E}^\dagger_{x,K}$ of degree $24$, with group $\Aut(E)\cong
 SL_2(\Z/3\Z)$, which contains all solutions of $M$. The underlying
 Galois representation of $G_{\F_4((1/z))}$ corresponds to the action of
 $\Aut(E)$ on $\mathrm{H}^1_{\mathrm{cris}}(E)$.

\subsection{A non-archimedean Riemann-Hilbert equivalence.}\label{sub::3.4} 

\begin{para}\label{para::3.4.1} Let us return to a general ordinary linear
 differential equation with polynomial or analytic coefficients 

 \begin{equation}
  \alpha_\mu(z) \partial^\mu y + \cdots + \alpha_0 (z) y = 0
  \;\text{ on $S\subset \A^1$, with $\partial=z\frac{d}{dz}$}.\tag{$\ast$}
 \end{equation}

 In the $p$-adic case, contrary to the complex situation, there
 does {\it not} exist, in general, any \'etale covering $S'\to S$ such
 that $(\ast)$ has a full set of solutions in $\mathcal{O}_{S'}$.
 
 \medskip\noindent {\it Example}: $\partial^2 y =0, \;S \supset
 \mathcal{C}(]1-\epsilon, 1[)$; there is no ``$p$-adic Riemann surface 
 of the logarithm'' on $S$.
 
 \medskip However, there are interesting global situations where this
 does happen. 

 \medskip\noindent {\it Example}: $\partial y - zy =0 $ on $S=\A^1$ :
 $S'=\matheur{D}(1,1^-)\xrightarrow{\log} S$, the infinite
 Galois \'etale covering with group $\mu_{p^\infty}\cong \Q_p/\Z_p$,
 trivializes the equation. 
 
 \medskip\noindent The object of this subsection is to investigate this
 phenomenon in full generality.
\end{para}

\begin{para}\label{para::3.4.2} Let $K$ be a complete non-archimedean
 field. We note that by definition, $K$-manifolds are `good' in the
 sense of \cite[1.2.6]{ecfnas}, \ie any point has an affinoid
 neighborhood. This property simplifies the theory of coherent sheaves
 (in the dictionary between rigid-analytic spaces and Berkovich analytic
 spaces, the theories of coherent sheaves and vector bundles correspond,
 \cf \cite[1.3.4, 1.6]{ecfnas})\footnote{Actually there are a number of
 technicalities involved in the analytic \'etale topology. First of all,
 since we restrict our attention to strictly analytic spaces in order to
 have the dictionary with rigid analytic varieties at disposal,
 Berkovich's definition of the \'etale site has to be slightly modified:
 this point is discussed in detail in R. Huber's book
 \cite[8.3]{ecoravaas}. Translated into rigid geometry, one gets the
 `partially proper \'etale' site of {\it loc.~cit.}. The associated
 sheaves correspond to the `overconvergent/conservative/constructible
 \'etale sheaves' of Van 
 der Put, Schneider, and De Jong \cite{jong96:_etale}, \cf
 \cite[8.2.12]{ecoravaas}. 
 We don't insist further on these subtleties.}.

 On the other hand, each coherent sheaf $M$ on $S$ gives rise to an
 \'etale sheaf $ M_{\mathrm{et}} $ on $S$ \cite[4.1.2]{ecfnas}.
 Moreover, one has \'etale descent for coherent sheaves: descent data are
 effective (\cf proof of \cite[4.1.10]{ecfnas}, and also
 \cite[3.2.3]{jong96:_etale} in the rigid setting). We shall use this only for
 locally free coherent sheaves (vector bundles) of rank $n$. These are
 in bijection with $\mathrm{H}^1_{\mathrm{et}}(S, GL_n)$.
  
\medskip In the sequel, $S$ will be a connected $K$-manifold, and $\ovl{s}$ a
geometric point of $S$. We assume that the point $s$ is $K$-rational,
 \ie is a classical point with residue field $K$.
\end{para}

\begin{para}\label{para::3.4.3} A continuous $K$-linear representation $\rho$
 of a topological group $\Gamma$ will be called {\it discrete} if its
 coimage $\Gamma/\ker \rho$ is discrete, or equivalently, if $\rho$
 factors through a discrete group.

 An {\it \'etale local system} on $S$ is a locally constant \'etale sheaf
 of $K$-vector spaces on $S$.
\end{para}

\begin{pro}\label{pro::3.4.4} There is a natural equivalence of categories   
\begin{equation*}
\{\text{discrete }\pi_1^{\mathrm{et}}(S,\ovl{s})\text{-representations }\}
 \xrightarrow{\sim}
 \{\text{\'etale local systems on $S$}\}.
\end{equation*}
\end{pro}

\begin{proof}
 Let $\rho $ be a discrete representation in a $K$-vector space
 $V$. This amounts to a continuous homomorphism $\pi_1^{\mathrm{et}}(S,
 \ovl{s}) \to \Aut_K V$, where $\Aut_K V$ denotes the discrete group of
 $K$-linear automorphisms of $V$.  On the other hand, any \'etale local
 system $\mathcal{V}$ on $S$ gives rise in the usual way to an \'etale 
 \v{C}ech cocycle with values in $\Aut_K V$, with $V=\mathcal{V}_s$, and
 its isomorphism class gives a well-defined element of 
 $\mathrm{H}^1_{\mathrm{et}}(S, \Aut_K V)$. Such a cohomology class may also be
 interpreted as usual as an $\Aut_K V$-principal homogeneous space
 $\mathcal{F}$ (as an \'etale sheaf) over $S$ up to isomorphism. What has
 to be shown is that such a principal homogeneous space is representable
 by an \'etale covering space. By lemma \ref{lem::1.2.6}. (ii), this is a local
 question on $S$. We thus may assume that $S$ is affinoid, and that for
 some finite \'etale covering $S'\to S$, $\mathcal{F}_{S'}$ is
 representable by an \'etale covering of $S'$. In that case, the
 representability of $\mathcal{F}$ follows from the descent theory of
 schemes.

 Starting from a representation $V$ of $\pi_1^{\mathrm{et}}(S,\ovl{s})$, we 
 have thus attached the $\Aut_K V$-principal homogeneous space and 
 etale sheaf $\FF$; the looked for \'etale local system is the 
 associated etale sheaf $\VV$ of $K$-vector spaces, with respect to 
 the canonical representation $V$ of $\Aut_K V$ (the pull-back of this 
 etale local system on $\FF$ is the constant sheaf with value $V$).
 
 Conversely, from an etale local system $\VV$, one attaches the $\Aut_K 
 V$-principal homogeneous space $\FF$ viewed as an \'etale covering 
 space (with $V=\VV_{\bar s}$); the structure of
 $\pi_1^{\mathrm{et}}(S,\ovl{s})$-set on its fiber at $\ovl{s}$ comes from
 a continuous homomorphism 
 $\pi_1^{\mathrm{et}}(S,\ovl{s})\to \Aut_K V$ (the cohomology class of
 $\FF$ in $\mathrm{H}^1(\pi_1^{\mathrm{et}}(S,\ovl{s}), \Aut_K V)$ is
 nothing but the latter homomorphism up to conjugation).
\end{proof}

In the sequel of this section, we assume that $K$ is of {\it
characteristic zero}.

\begin{para}\label{para::3.4.5} Let $M$ be a vector bundle (= locally free
 $\mathcal{O}_S$-module of finite rank) with a connection: $\nabla:
 M\to \Omega^1_S \otimes M$ (\ie a $K$-linear map satisfying Leibniz
 rule). This can also be interpreted $\mathcal{O}_S$-linearly, in terms
 of first-order jets, as a section $M\to \mathcal{P}^1(M)$ of the natural
 projection $\mathcal{P}^1(M)\to M$.

 \noindent The vector bundles with connection (\resp with integrable
 connection) form a tannakian category over $K$. The unit for $\otimes$
 is the trivial connection  $d: \mathcal{O}_S\to  \Omega^1_S $.

\medskip Because $s$ is a classical point with residue field $K$ (of
 characteristic zero), the stalk at $s$ of the sheaf $M^\nabla=\ker
 \nabla$ of germs of horizontal sections is a $K$-vector space of rank
 the rank of $M$, which can be identified with the fiber $M_{(s)}$ (note
 that since $S$ has a $K$-rational point, 
 $(\mathcal{O}_{S_{\mathrm{et}}})^{d_{S_{\mathrm{et}}}}=K$).

 \medskip To $M$ (\resp $\mathcal{P}^1(M)$), one attaches an \'etale sheaf
 $M_{\mathrm{et}}$ (\resp $\mathcal{P}^1(M)_{\mathrm{et}}$) --- locally free
 $\mathcal{O}_{S_{\mathrm{et}}}$-modules. Thus $\nabla$ also correspond
 to a section of the natural projection $\mathcal{P}^1(M)_{\mathrm{et}}\to 
 M_{\mathrm{et}}$, or equivalently, to a $K$-linear map
 $\nabla_{\mathrm{et}}: \;M_{\mathrm{et}}\to  \Omega^1_{S_{\mathrm{et}}}
 \otimes M_{\mathrm{et}}$. The kernel of the latter map is an \'etale
 sheaf of $K$-vector spaces denoted by $M_{\mathrm{et}}^\nabla$ (the
 \'etale sheaf of germs of horizontal sections). Let $t$ be point of
 $S$. The fiber of $M_{\mathrm{et}}^\nabla$ at $t$ (\cf
 \cite[4.2]{ecfnas} is a $\ovl{K}$-vector space endowed with a discrete
 $\Gal(\ovl{\mathscr{H}(t)}/{\mathscr{H}(t)})$-action. If
 $M_{\mathrm{et}}^\nabla$ is locally constant, the dimension of this
 vector space is independent of $t$; taking $t=s$, we find that it is
 the rank of $M$.  
\end{para}

 \medskip\noindent{\it Remark.} Any classical point admits a neighborhood
 over which $M^\nabla$ is constant of rank equal 
 to the rank of $M$ (by the non-archimedean analogue of Cauchy's theorem
 on solutions of analytic differential equations). 
 The constancy of $M_{\mathrm{et}}^\nabla$ imposes a constraint only at
 non-classical points of the Berkovich space $S$.  

\begin{thm}\label{thm::3.4.6} There is a natural equivalence of tannakian
 categories (the {\rm non-archimedean \'etale Riemann-Hilbert functor}):
 \begin{equation*}
  \left\{
   \begin{matrix}
    \text{discrete finite dimensional}\\
    \pi_1^{\mathrm{et}}(S, \ovl{s})\text{-representations}
   \end{matrix}
   \right\}
   \xrightarrow{RH^{\mathrm{et}}}
   \left\{
    \begin{matrix}
     \text{vector bundles with integrable}\\
     \text{connection $(M,\nabla)$ on $S$ such that}\\
     \text{$M_{\mathrm{et}}^\nabla$ is an \'etale local system}
    \end{matrix}
   \right\}.
 \end{equation*}
 In this correspondence, the subspace of 
 $\pi_1^{\mathrm{et}}(S, \ovl{s})$-invariants corresponds to the
 space of global sections of $M^\nabla$.
\end{thm}

\begin{proof}
 Using the previous proposition and \'etale descent of coherent
 modules, it suffices to show that the natural maps 
 \begin{equation*}
  \mathcal{V}\mapsto (\mathcal{V}\otimes_K \mathcal{O}_{S_{\mathrm{et}}},
   \id\otimes d_{S_{\mathrm{et}}}),\;\;\;(M,\nabla)\mapsto
   M_{\mathrm{et}}^\nabla,
 \end{equation*}
 induce quasi-inverse functors
 \begin{equation*}
  \left\{
   \begin{matrix}
    \text{locally constant \'etale sheaves}\\
    \text{of fin. dim. $K$-vector spaces}
   \end{matrix}
  \right\}
   \rightleftharpoons
   \left\{
    \begin{matrix}
     \text{vector bundles with integrable}\\
     \text{ connection $(M,\nabla)$ such that}\\
     \text{$M_{\mathrm{et}}^\nabla$ is an \'etale local system}
    \end{matrix}
   \right\}
 \end{equation*}
 It is clear that these maps are functorial (and the first one is clearly
 compatible with the tensor 
 structures). The point is to show that they are quasi-inverse. 
 
 On the one hand, it is clear that $\id\otimes d_{S_{\mathrm{et}}}$ is an
 integrable connection, and that $(\mathcal{V}\otimes_K
 \mathcal{O}_{S_{\mathrm{et}}} )^{\id\otimes
 d_{S_{\mathrm{et}}}}=\mathcal{V}$ (taking into account that
 $(\mathcal{O}_{S_{\mathrm{et}}} )^{d_{S_{\mathrm{et}}}}=K$).
 We have to show that if $(M,\nabla)  $ is a vector bundle with integrable
 connection such that $M_{\mathrm{et}}^\nabla$ is locally constant, then
 the natural horizontal map

 \begin{equation*}
 \phi: M_{\mathrm{et}}^{\nabla}\otimes_K \mathcal{O}_{S_{\mathrm{et}}}\to
 M_{\mathrm{et}}.
 \end{equation*}
 is an isomorphism. By the remark preceding the
 theorem, the stalks of $M_{\mathrm{et}}^\nabla$ are $\ovl{K}$-vector
 spaces of constant dimension, which is the rank of $M$. It follows that 
 $M_{\mathrm{et}}^\nabla\otimes_K \mathcal{O}_{S_{\mathrm{et}}}$ is the
 \'etale sheaf attached to a vector bundle of the same 
 rank as $M$. It is a general well-known fact about integrable
 connections\footnote{proven by passing to the completion around
 classical points.} that $\phi$ is a monomorphism of vector bundles, and
 that $\Im\phi$ is locally a direct summand of
 $M_{\mathrm{et}}$ (because the cokernel inherits an integrable
 connection, hence is locally free). To show that this it is an 
 epimorphism, it thus suffices to check that
 $M_{\mathrm{et}}^\nabla\otimes_K \mathcal{O}_{S_{\mathrm{et}}}$ and
 $M_{\mathrm{et}}$ have the same rank, which is the case by assumption.
\end{proof}

\begin{para}\label{para::3.4.7} The $ \pi_1^{\mathrm{et}}(S, \ovl{s})
 $-representation attached to $(M,\nabla)$ is called the {\it
mono\-dromy representation}. Note that the underlying vector space is
canonically identified with the fiber $M_s$. The image of the monodromy
representation is called the {\it monodromy group}. It may be
non-discrete (while the coimage is discrete by assumption); this
phenomenon is not specific to the non-archimedean situation, and is
indeed familiar in the complex case.
\end{para}

\begin{para}\label{para::3.4.8} The \'etale Riemann-Hilbert functor unifies
 several earlier constructions. Indeed, if the monodromy representation
 factors through $\pi_1^{\mathrm{top}}(S,s)$, we recover the
 `topological Riemann-Hilbert functor' considered in
 I.\ref{sub-connection}, whose image consists of connections with
 locally constant sheaves of solutions.

\noindent It also includes some cases of overconvergent unit-root
isocrystals with finite monodromy already considered in
I.\ref{sub-chiaroscuro} (in that case, the monodromy factors through a
discrete, hence finite, 
representation of $\pi_1^{\mathrm{alg}}(S, \ovl{s})$).

 Of course, connections in the image of the \'etale Riemann-Hilbert
 functor are still rather special\footnote{for instance, it does not
 account for the differential equation of the logarithm, as we 
 have said in \ref{para::3.4.1}; to see this, one can notice that otherwise the
 monodromy group would be contained in the additive 
 group of $K$, hence torsion-free, but this implies that the \'etale
 covering which trivializes the connection is a 
 topological covering (\ref{thm::2.1.8}), hence is trivial.}, but it
 seems difficult to extend its scope without invoking ``relative
 Fontaine rings''.
\end{para}

\begin{rems}
  (i) Whether $M_{\mathrm{et}}^\nabla$ is an
	\'etale local system or not can be checked locally on $S$. In 
	particular, the property for a connection to be in the essential image
	of $RH^{\mathrm{et}}$ is local. 

  (ii) The \'etale fundamental group $\pi_1^{\mathrm{et}}(\P^1_{\C_p},
	\ovl{s})$, although highly non-trivial, does not have any non-trivial
	discrete finite-dimensional representation. Indeed, 
	the corresponding connection on
	$\P^1_{\C_p}$ would be algebrizable, and necessarily trivial.
\end{rems}

\subsection{Examples and applications.}\label{sub::3.5} 

\begin{para}\label{para::3.5.1} For any full subcategory $\Cov^\bullet_S$ of
 $\Cov^{\mathrm{et}}_S$ which is stable under 
 taking connected components, fiber products  and quotients, we have seen
 that there is a natural continuous homomorphism 
 \begin{equation*}
  \pi_1^{\mathrm{et}}(S, \ovl{s})\to \pi_1^\bullet(S, \ovl{s}) 
 \end{equation*}
 with dense image (\ref{cor::1.4.8}). Any discrete representation
 of $ \pi_1^\bullet(S, \ovl{s})$ gives rise to a 
 discrete representation of $ \pi_1^{\mathrm{et}}(S, \ovl{s})$ with the
 same coimage. It follows that the \'etale Riemann-Hilbert 
 functor induces a fully faithful functor 
 \begin{equation*}
  \left\{
   \begin{matrix}
    \text{discrete fin. dim.}\\
    \text{$\pi_1^\bullet(S,\ovl{s})$-representations}
   \end{matrix}
   \right\}
   \xrightarrow{RH^\bullet}
   \left\{
    \begin{matrix}
     \text{vector bundles with integrable}\\
     \text{connection $(M,\nabla)$ on $S$}
    \end{matrix}
   \right\}.
 \end{equation*}
 Note that the monodromy group defined in \ref{para::3.4.7} coincides
 with the image 
 of $\pi_1^\bullet(S, \ovl{s})$.  We shall discuss the avatars
 $RH^\bullet$ in the cases of $\Cov^{\mathrm{top}}_S$,
 $\Cov^{\mathrm{alg}}_S$, $\Cov^{\mathrm{loc.alg}}_S$,
 $\Cov^{\mathrm{temp}}_S$.
\end{para}

\begin{para}\label{para::3.5.2} The functor $RH^{\mathrm{top}}$ was studied in
 I.\ref{sub-connection} (in the case $K=\C_p$). Its essential image consists 
 of connections which satisfy `Cauchy's theorem' at every (possibly
 non-classical) point of $S$, or equivalently, which have a full set of
 multivalued analytic solutions, or else, such that 
 $M^\nabla$ is a locally constant sheaf on the topological space $S$.
\end{para}

\begin{para}\label{para::3.5.3} The essential image of the functor
 $RH^{\mathrm{alg}}$ consists of objects $(M,\nabla)$ such that
 $M_{\mathrm{et}}^\nabla$ becomes constant on some finite \'etale
 covering of $S$.

 Let us assume for instance that $S$ is affinoid with good reduction
 $S^0$, and that $K$ is a finite extension of
 $\widehat{\Q_p^{\mathrm{ur}}}$ as in 
 I.\ref{sub-r-and-unit-root}. Then among the finite representations of $
 \pi_1^{\mathrm{alg}}(S, \ovl{s}) $, one can distinguish those which
 factor through $ \pi_1^{\mathrm{alg}}(S_0, \ovl{s}_0) $. The
 corresponding connections underlie overconvergent unit-root isocrystals
 with finite global monodromy. We have seen a standard example in
 I.\ref{exa-dwork-exp}, given by the differential equation
 $\frac{d}{dz} y = -\pi y$ 
 on $S=\matheur{D}(0,1^+)$, which is trivialized over an Artin-Schreier covering.
 
 Let us now consider a case of bad reduction, namely the case of an
 annulus $S_r=\mathcal{C}(]1,r[),\;r\in | K|, 
 r>1$ over any discretely valued $p$-adic field $K$. Let $\mathcal{R}$ be
 the Robba ring at infinity. Then $RH^{\mathrm{alg}}$ induces a 
 functor
 \begin{equation*}
  \{\text{finite representations of $\pi_1^{\mathrm{alg}}(S_r,
   \ovl{s})$}\}
   \to \{\text{$\mathcal{R}$-differential modules}\}.
 \end{equation*}
 In relation to the local monodromy theorem \ref{coj::3.3.5}, it would be very
 interesting to understand the essential image of this functor as $r\to
1^+$. Note however that there is a difficulty with base-points in
passing to the limit $r\to 1^+$ (the problem disappears if one considers
only abelian representations).

Nevertheless, any $\mathcal{R}$-differential module $M$ extends to a
vector bundle $(M_r,\nabla_r)$ with connection on some $S_r$. Let us
assume that $M$ is semisimple and has a Frobenius structure.

If $M$ is tame, then \ref{thm::3.2.8} ensures that $M$ is in the image
of the above functor
\begin{equation*}
 \{\text{finite $\pi_1^{\mathrm{alg}}(S_r, \ovl{s})$-representations}
  \}\to \{\text{$\mathcal{R}$-differential modules}\}
\end{equation*}
for $r$ close enough to
$1$. The local monodromy theorem \ref{coj::3.3.5} implies the same if
$M$ is wild, 
although the monodromy representation might be non-abelian in this case,
as \ref{thm::3.3.6} shows.
\end{para}

\begin{para}\label{para::3.5.4} The essential image of the functor
 $RH^{\mathrm{loc.alg}}$ consists of objects $(M,\nabla)$ such that for every
 compact subdomain $W\subset S$, $(M_{\mathrm{et}}^\nabla)_{| W}$ becomes
 constant on some finite \'etale covering of $W$.
 Although not explicitly introduced in \cite{ramero98}, it is clear that
 these connections played an important heuristic role in that paper,
 especially those coming from abelian representations of the locally
 algebraic fundamental group (1.5).

 The basic example is $(\mathcal{O}_{\A^1}, d-1)$ over $\C_p$ (with
 solution $e^z$). It is trivialized by the logarithmic \'etale covering
 of ${\A^1}$.
 
 Similarly, for any ramified polynomial $P(z^{1/m})$, the differential
 equation with basic solution $e^{P(z^{1/m})}$ is trivialized over a
 locally algebraic \'etale covering of a sufficiently small punctured
 disk $\matheur{D}_{\C_p}(\infty,\epsilon)^\ast$ around $\infty$ (the
 logarithmic 
 covering of a Kummer covering). By \ref{para::3.1.5}, we see that any
 differential 
 equation with polynomial coefficients --- or more generally meromorphic
 singularity at $\infty\;$ --- with {\it rational} formal exponents at
 $\infty$, and with no logarithm in the formal solutions $\hat y$ (\cf
 \ref{para::3.1.2}), corresponds to an object in the image of the
 locally algebraic Riemann-Hilbert functor
 \begin{equation*}
  \left\{
   \begin{matrix}
    \text{finite dimensional}\\
    \text{$\pi_1^{\mathrm{loc.alg}}(\matheur{D}(\infty,\epsilon)^\ast,
    \ovl{s})$-representations}
   \end{matrix}
  \right\}
  \to
  \{\text{$\mathcal{O}(\matheur{D}(\infty,\epsilon)^\ast)$-differential
  modules}\}
 \end{equation*}
 for $\epsilon$ small enough.

 There is a certain analogy between the logarithmic \'etale
 covering of $\matheur{D}_{\C_p}(\infty,\epsilon)^\ast$ and the {\it
 real blow-up} of $\matheur{D}_{\C}(\infty,\epsilon)^\ast$ which plays
 such an important role in 
 the study of irregular singularities of complex differential
 equations. In both cases, the point at infinity becomes a limiting
 circle\footnote{L. Ramero points out that the change from analytic to
 \'etale topology is essential in the $p$-adic case, while in the
 complex case, the real blow-up corresponds to completion for a
 different uniform structure of the punctured disk with no change of the
 topology.}.

\medskip Let us return to a general vector bundle $(M,\nabla)$ with
 connection on $\matheur{D}(\infty,1^-)^\ast$ which extends to a vector bundle
 with meromorphic connection at $\infty$. We assume that it is has a
 Frobenius structure near the unit circle, and that the formal exponents
 at infinity are rational. Regarding problem \ref{para::3.2.9}, a
 natural question is:
\end{para}

\begin{que}\label{que::3.5.5} Under which further conditions, if any,
 does there exist a locally algebraic covering of
 $\matheur{D}(\infty,1^-)^\ast$ 
 over which the connection becomes an iterated extension of trivial
 connections?
\end{que}
At any rate, in the cases where this holds, this would let Ramero's
locally algebraic fundamental group play the same role, 
over $\C_p$, as Ramis' wild fundamental group (over $\C$) in the local
study of meromorphic irregular singularities.

\begin{para}\label{para::3.5.6} The essential image of the functor
 $RH^{\mathrm{temp}}$ consists of objects $(M,\nabla)$ such that 
 $M_{\mathrm{et}}^\nabla$ becomes constant on some temperate \'etale
 covering of $S$.  This implies that there is a finite \'etale Galois
covering $f:\;T\to S$ such that the pull-back
$(M',\nabla')=f^\ast(M,\nabla)$ has a full set of multivalued analytic
solutions, \ie $M^{\prime \nabla'}$ is locally constant. Let $G$ be the group
$\Aut_S\,T$, and set $n=| G| $.
 
Assume that there are $n$ distinct $K$-rational points $t_i$ in $T$
 lying above $s$. Then the monodromy group admits the following {\it
 concrete interpretation}: 

\noindent let $\alpha_{ij}$ be a Berkovich path (up to homotopy)
connecting $t_i$ to $t_j$. Then $f_\ast(\alpha_{ij})\in
\pi_1^{\mathrm{temp}}(S,\ovl{s})$. When $i,j$ and $\alpha_{ij}$ vary,
they generate the image $\Aut_S\,\tilde T$ of
$\pi_1^{\mathrm{temp}}(S,\ovl{s})$ (note that, varying $j$, they
generate the finite quotient $G$). It follows that the monodromy group
consists of the images of the $f_\ast(\alpha_{ij})$'s in $GL(M_{(s)})$
by the monodromy representation.

Therefore, except for the occurrence of the finite \'etale covering, the 
picture is quite similar to the complex situation.

\end{para}

\begin{pro}\label{pro::3.5.7} Let $S=\A^1\setminus\{\zeta_1,\ldots,
 \zeta_m\}$. Assume\footnote{this very strong assumption is rarely met
 in practice.} that the residue characteristic $p$ of $K$ does not
 divide $n$. Then:
 \begin{enumerate}
  \renewcommand{\theenumi}{\alph{enumi}}
  \item for any open or closed disk $\matheur{D}\subset S$, the restriction of
	$\nabla$ to $\matheur{D}$ is trivial (\ie has a full set of solutions); in
	particular, $\nabla$ satisfies the Robba condition;
  \item for any open or closed annulus
	$\mathcal{C}\subset S$, the pull-back of $\nabla$ to the Kummer
	covering of $\mathcal{C}$ of degree $n$ is trivial.
 \end{enumerate}

 In particular,
 $p$-adic exponents in the sense of Christol-Mebkhout exist (in any
 annulus which does not contain any $\zeta_i$), and they are rational
 numbers with denominator $n$ (or a divisor of $n$).
\end{pro}

\begin{proof}
 One is immediately reduced to considering closed disks and
 annuli. We know the pull-back of $\nabla$ to $\tilde T$ is trivial. Let
 $\matheur{D}'$ (\resp $\mathcal{C}'$) be a connected component of $\matheur{D}\times_S
 \tilde T$ (\resp $\mathcal{C}\times_S \tilde T$). Because $n$ is assumed
 to be prime to $p$, the covering $\matheur{D}'\to \matheur{D}$ is trivial (\resp is a Kummer
 covering of degree dividing $n$), \cf \cite[6.3.3, 6.3.5]{ecfnas} (taking
 into account the fact that disks and annuli are simply-connected). The
 assertions follow.
\end{proof}

\medskip\noindent {\it Example.} (see also the example at the end of
\ref{lem::2.1.2}). Let $K$ be a $p$-adic field for any odd $p$. 
Let $T=E\setminus E[2]$ be a Tate curve deprived from its four
$2$-torsion points (assumed to be $K$-rational), written in 
Legendre form: 
\begin{equation*}
 y^2=x(x-1)(x-\lambda),\quad |\lambda(\lambda -1)|\neq 1.
\end{equation*}
Let $f$ be the projection from $T$ to $S=\A^1\setminus
\{0,1,\lambda \}$.
Here $\tilde T$ is of the form $\G_m\setminus \pm(\sqrt q)^\Z$, and
$T\cong \tilde T/ q^\Z $.

The sum of the canonical one-dimensional $K$-linear representation of
$q^\Z$ ($q$ acting by multiplication by a factor $q$) and its dual
extends to a representation of the semi-direct product of $q^\Z$ and
$\Z/2\Z$ (which is the Galois group of the covering $\tilde T\to S$). On
can choose the local monodromies at $0,\lambda, 1,\infty$ respectively
in such a way that they correspond to the matrices
\begin{equation*}
 \begin{pmatrix}0&1\\1&0\end{pmatrix},
 \begin{pmatrix}0&1\\1&0\end{pmatrix},
  \begin{pmatrix}0&q\\q^{-1}&0\end{pmatrix},
  \begin{pmatrix}0& q^{-1}\\q&0\end{pmatrix}
\end{equation*}
respectively. It is
then easy to check that the monodromy around an annulus surrounding $0$
and $\lambda$ alone ($|\lambda|$ being assumed $<1$ for convenience) is
trivial: the $p$-adic exponents are zero.
 
\begin{para}\label{para::3.5.8} Much more interesting examples arise in the
 context of $p$-adic period mappings. We place ourselves in the 
 setting of  II.\ref{sub:7.2}. Assuming the congruence subgroup sufficiently
 small, the picture was the following: 

\begin{center}
 \begin{picture}(150,90)(0,-5)
  \put(10,60){\makebox(0,0)[c]{$\til{S}$}}
  \put(25,60){\vector(1,0){20}}
  \put(25,45){\vector(1,-1){20}}
  \put(60,60){\makebox(0,0)[c]{$\mathcal{M}$}}
  \put(60,45){\vector(0,-1){20}}
  \put(60,10){\makebox(0,0)[c]{$S$}}
  \put(115,60){\makebox(0,0)[c]{$\mathcal{D}$}}
  \put(75,60){\vector(1,0){25}}
  \put(88,65){\makebox(0,0)[c]{${}^\mathcal{P}$}}
  \put(115,45){\vector(0,-1){20}}
  \put(75,10){\vector(1,0){20}}
  \put(115,10){\makebox(0,0)[c]{$\mathcal{D}/\Gamma$}}
 \end{picture}
\end{center}
where $S$ was a suitable connected `tube-like' domain
 in the Shimura variety $\Sh_{\C_p}$ ($\mathcal{M}^0$ was a deformation
 space for $p$-divisible groups quasi-isogenous to a fixed one, and
 $\mathcal{D}$ was the period domain).

 The period mapping turned out to be expressable in terms of quotients of
 solutions of the Gauss-Manin connection (on $S$), 
 and the Gauss-Manin connection was trivializable on $\mathcal{M}^0 $,
 hence on $\tilde S$ (II.\ref{cor:7.2.6}).

 \medskip Let us now abandon the assumption that the congruence subgroup
 is sufficiently small. Then a similar picture holds 
 true: $\Gamma$ is now only virtually torsion-free, and the covering
 $\mathcal{M}^0\to S $ is now ramified. By removing the 
 singularities, one obtains a {\it temperate \'etale covering}, which
 trivializes the Gauss-Manin 
 connection (using II.\ref{cor:7.2.6} after passing to a suitable
 finite \'etale covering). It follows that the 
 {\it Gauss-Manin connection (restricted to $S$ deprived from the
 singularities) is in the image of the temperate 
 Riemann-Hilbert functor}. 

 The situation is therefore similar the the complex one, except that $S$
 is only a small domain in the Shimura variety in 
 general. The most favorable situation occurs in the case of global
 uniformization (II.\ref{sub:7.3}), \ie when $S$ is a component of $\Sh$.
\end{para}

\begin{para}\label{para::3.5.9} Let us now give a criterion for a connection
 in the image of $RH^{\mathrm{et}}$ to be in the image of $RH^{\mathrm{temp}}$.
\end{para}

\begin{pro}\label{pro::3.5.10} Let $\rho:
\pi_1^{\mathrm{et}}(S,\ovl{s}) \to GL(V)$ be a finite-dimensional
 discrete representation (corresponding to an \'etale local system 
 $\mathcal{V}$). Let us consider the following conditions:
 \begin{enumerate}
  \renewcommand{\theenumi}{\alph{enumi}}
  \item $\rho$ factors through a discrete representation of
	$\pi_1^{\mathrm{temp}}(S,\ovl{s})$,
  \item $\Im \rho$ is virtually torsion-free,
  \item $\Im \rho$ is finitely generated.
 \end{enumerate}
 Then {\rm (c)} $\Rightarrow$ {\rm (b)} $\Rightarrow$ {\rm (a)}. If moreover
 $S$ is the 
 analytification of a smooth algebraic $K$-variety, 
 then {\rm (a)} $\Rightarrow$ {\rm (c)}.
\end{pro}

\begin{proof}
 (c) $\Rightarrow$ (b): this is Selberg's lemma: any finitely
 generated subgroup of $GL(V)$ is virtually torsion-free.

 (b) $\Rightarrow$ (a) follows from \ref{thm::2.1.8}.~(b). 

 (a) $\Rightarrow$ (c): under (a), $\Coim\,\rho$ is the Galois group of a
 Galois temperate covering of $S$. If $S$ is 
 algebraic, $\Coim\,\rho$ must be finitely generated by \ref{pro::2.1.6}. 
\end{proof}

\begin{para}\label{para::3.5.11} We finish the section by showing that the
 connections which lie in the image of the temperate 
Riemann-Hilbert functor are algebrizable (this is not automatic in the
 non-archimedean situation), and even fuchsian: 
\end{para}

\begin{thm}\label{thm::3.5.12} Assume that $S$ is the analytification of
 a smooth algebraic $K$-variety $S^{\mathrm{alg}}$. Let 
 $(M,\nabla)$ be a vector bundle with integrable connection on $S$ which
 comes from a discrete representation of 
 $\pi_1^{\mathrm{temp}}(S,\ovl{s})$. Then $(M,\nabla)$ is the
 analytification of a (unique) algebraic vector 
 bundle with integrable connection on $S^{\mathrm{alg}}$, which has
 regular singularities and exponents in $\Q/\Z$ at any 
divisorial valuation of $K(S)$.  
\end{thm}

\begin{proof}
 By assumption, there is a finite Galois \'etale covering $S'\to
 S$ (with group $G$) such that $(M,\nabla)$ becomes 
 trivial on the topological universal covering of $S'$. The pull-back
 $(M,\nabla)_{S'}$ corresponds to a (unique) 
 representation of 
 $\pi_1^{\mathrm{top}}(S', s')$ (for any point $s'$ above $s$). Let
 $\ovl{S}'$ be a smooth algebraic compactification of $S'$ 
 (Hironaka). We may replace $K$ by the completion of its algebraic
 closure. By \ref{pro::1.1.3}, the homomorphism 
 $\pi_1^{\mathrm{top}}(S', s')\to\pi_1^{\mathrm{top}}(\ovl{S}', s')$ is
 an isomorphism, and it follows that $(M,\nabla)_{S'}$ extends 
 (uniquely) to a vector bundle with integrable connection on 
 $\ovl{S}'$. By GAGA (K\"opf), the latter is the analytification of a
 (unique) algebraic vector 
 bundle with integrable connection on $\ovl{S}^{\prime \mathrm{alg}}$. It is
 now easy to conclude by Galois descent from $S^{\prime \mathrm{alg}}$ to 
 $S^{\mathrm{alg}}$ (we refer to 
 \cite{andre01:_de_rham} for the details on regular singularities,
 exponents, and GAGA functor for connections).
\end{proof}

\begin{para}\label{para::3.5.13} Besides discrete representations of
 $\pi_1^{\mathrm{temp}}(S,\ovl{s})$, it is also useful to 
 consider {\it projective discrete representations}, \ie discrete
 representations with value in $PGL(V)$.  They give rise to {\it
 projective connections} over $S^{\mathrm{alg}}$, at least when
 $S^{\mathrm{alg}}$ is an {\it affine curve}. 

 Indeed, replacing $S^{\mathrm{alg}}$ by a finite Galois \'etale
 covering $S^{\prime \mathrm{alg}}$ with Galois group $G$, we may
 replace $\pi_1^{\mathrm{temp}}(S')$ by $\pi_1^{\mathrm{top}}(S')$, 
 which is a free group, for which the lifting property is immediate. This
 provides a connection $\nabla'$ over $S^{\prime \mathrm{alg}}$,
 and the associated projective connection is endowed with a $G$-action
 compatible with the $G$-action on $S^{\prime \mathrm{alg}}$, hence
 descends to a projective connection on $S^{\mathrm{alg}}$.
\end{para}

One has more: 
\begin{pro}\label{pro::3.5.14} Then any projective discrete representation
 $\ovl{\rho}:$ $\pi_1^{\mathrm{temp}}(S,\ovl{s})$ $\to$ $PSL(V)$ lifts to
 a genuine discrete representation
 $\rho:\,\pi_1^{\mathrm{temp}}(S,\ovl{s})\to SL(V)$, after replacing $K$
 by a finite extension.
\end{pro}

\begin{proof}
 We have seen that $\ovl{\rho}$ corresponds to a projective connection on
 $S^{\mathrm{alg}}$. Now it is well known that any projective connection
 lifts to a genuine connection after replacing $K$ by a finite extension
 (because $\mathrm{H}^2(S^{\mathrm{alg}}_{\ovl{K}, et},
 \Z/2\Z)=0$, \cf \cite[2.2.2.1]{katz90:_expon}). It suffices to check that the
 pull-back $\nabla''$ of this connection on $S^{\prime \mathrm{alg}}$
 comes from a discrete representation
 $\rho:\,\pi_1^{\mathrm{temp}}(S',\ovl{s}')\to SL(V)$. But we have on the
 other hand our connection $\nabla'$ which is projectively equivalent to
 $\nabla''$, 
 and which we may assume to come from a discrete representation 
 $\pi_1^{\mathrm{temp}}(S',\ovl{s}')\to SL(V)$. Let $S''$ be a finite
 \'etale covering of $S'$ such that the pull-back of any class in
 $\mathrm{H}^1(S^{\mathrm{alg}}_{\ovl{K}, 
 et}, \Z/2\Z)$ vanishes in $\mathrm{H}^1(S^{\prime \prime\,
 \mathrm{alg}}_{\ovl{K}, et}, 
 \Z/2\Z)$. Since $\nabla'$ and $\nabla''$ are projectively equivalent,
 their pull-backs on $S''$ are isomorphic 
 ({\it loc.~cit.}). It follows that $\nabla''$ comes from a discrete
 representation $\pi_1^{\mathrm{temp}}(S',\ovl{s}')\to SL(V)$.        
\end{proof}

\medskip We leave it as an exercise to translate \ref{pro::3.5.10} in
the context of projective connections.

\newpage
\section{Non-archimedean orbifolds, uniformizing differential equations
and period mappings.}\label{sec::4}

\markboth{\thechapter. $p$-ADIC ORBIFOLDS AND MONODROMY.}%
{\thesection. NON-ARCHIMEDEAN ORBIFOLDS.}

\begin{abst}
We review the complex case, and then define and study non-archimedean
orbifolds, their associated orbifold fundamental groups, and
uniformizing differential equations.
\end{abst}

\subsection{Review of complex orbifolds.}\label{sub::4.1}

\begin{para}\label{para::4.1.1} There are several variants of the notion 
 of complex orbifold. We follow \cite{calip} and see orbifolds 
 as being locally quotients of a complex manifold by a finite group 
 acting faithfully. More precisely, a {\it complex orbifold} 
 $\mathcal{S}=(\ovl{S}, (Z_i,e_i))$ consists in 
 \begin{itemize}
  \renewcommand{\labelitemi}{\normalfont\bfseries\textendash}
  \item a normal analytic space $\ovl{S}$,
  \item a locally finite collection of irreducible (Weil) divisors
	$(Z_i)$,
  \item a positive integer\footnote{in other variants, one also allows
	$e_i=\infty$. This leads to `parabolic structures' which have
	indeed non-archimedean counterparts in positive characteristic
	(cusps of Drinfeld modular curves...). But we shall be mostly
	interested in non-archimedean orbifolds in characteristic zero,
	where such `parabolic structures' cannot occur. } $e_i$
	attached to each $Z_i$.
 \end{itemize}

\medskip\noindent It is assumed that $\ovl{S}$ is covered by so-called
 {\it orbifold charts}, \ie morphisms  $W 
\xrightarrow{\phi} V \subset \ovl{S} $ where 
 \begin{itemize}
  \renewcommand{\labelitemi}{\normalfont\bfseries\textendash}
  \item $W$ is a complex manifold,
  \item $\phi$ is a (ramified) covering in the sense of \ref{def::1.2.3}
	(\ie $V$ is covered by open subsets $U$ such that $\coprod V_j =
	\phi^{-1}U \to U$ and the restriction of $\phi$ to every $V_j$
	is finite),
  \item except for a closed subset $\Sigma$ of codimension $2$ in $V$,
	$\phi$ is \'etale above $V\setminus (V\cap\bigcup Z_i)$ (hence
	induces a topological covering of $V\setminus (V\cap\bigcup
	Z_i)$), and is {\it ramified with index $e_i$ above
	$Z_i$}.\footnote{it is easily seen that this definition is equivalent
	to that of \cite{calip}, using {\it loc.~cit.} 14.2.}
 \end{itemize}

\medskip Note that if $s$ is a non-singular point of $\ovl{S}$ not
 contained in $\bigcup Z_i$, then any smooth open neighborhood of $s$ in
 $\ovl{S}$ disjoint from $\bigcup Z_i$ provides an orbifold chart around
 $s$. Hence we may and shall assume that $\Sigma $ is contained in
 $\bigcup Z_i\cup \ovl{S}^{\mathrm{sing}}$ (where $\ovl{S}^{\mathrm{sing}}$ denotes the
 set of singular points of $\ovl{S}$).
 
 \noindent We shall often write
 $Z=\bigcup Z_i\cup
 \ovl{S}^{\mathrm{sing}}$, and
 $S=\ovl{S}\setminus Z$ (which is a manifold).  

\medskip\noindent {\it Remark.} By shrinking the $W$'s, one can assume
 that they are open subsets of $\C^n$; this matches with the conventional
 use of the word `chart'. However, we do not impose here this condition
 which would have no non-archimedean counterpart.

\end{para}

\begin{para}\label{para::4.1.2} A {\it morphism of orbifolds}
 $f:\mathcal{S}'\to \mathcal{S}$ is a morphism $\ovl{S}'\to \ovl{S}$, such
 that $\ovl{S}'$ is covered by orbifold charts  $\phi': W'\to V'$ with
 the property that $f\circ \phi'$ factors through some orbifold 
 chart $\phi:W\to V$ of $\mathcal{S}$: $f\circ \phi'=\phi\circ f'$ 
(for some morphism $f'\colon W'\to W$).

 \medskip\noindent {\it Example.}  Let $(\ovl{S}, (Z_i,e_i)_{i\in I})$ be
 an orbifold, and let $J\subset I$ be the subset of indices for which
 $e_i>1$. Then the identity of $\ovl{S}$ induces an isomorphism
 $(\ovl{S}, (Z_i,e_i)_{i\in I})\cong (\ovl{S}, (Z_i,e_i)_{i\in J})$. In
 other words, one can `forget' the divisors with index $e_i=1$.
 
 \noindent Considering complex manifolds as orbifolds with an empty
 collection of divisors ($S=\ovl{S}$), one can embed the 
 category of complex manifolds as a full subcategory of the category of
 complex orbifolds.

 \noindent It will be convenient to attach to any irreducible divisor
 $D\subset \ovl{S}$ an index $e_D\in \Z_{>0}$ which is $e_i$ if
 $D=Z_i$ and $1$ if $D$ is none of the $Z_i$'s.

\end{para}

\begin{para}\label{para::4.1.3} A morphism $f$ is called an {\it
 orbifold-covering} if it is a covering in the sense of \ref{def::1.2.3}
 and for any 
 orbifold chart $\phi': W'\to V'$ of $\mathcal{S}'$, $f\circ\phi': W'\to
 V=f(V')$ is an orbifold chart of $\mathcal{S}\;$ (in
 loose words: any lifting $f'$ as in \ref{para::4.1.2} is 
 \'etale).  
 
 \medskip\noindent The restriction of $f$ above $S$ is then an \'etale
 covering: indeed, this is clear except above
 a closed subset of codimension $\geq 2$ in $S$, thus above the whole of
 $S$ by purity since $S$ is smooth.
 
 \medskip\noindent {\it Example.} $(\A^1, (0,e'))
 \xrightarrow{z\mapsto z^n} (\A^1, (0,e)) $ is
 a morphism if and only if $e| ne'$, and is an orbifold-covering if
 and only if 
 $e=ne'$.   
  
\medskip \noindent {\it Remark.} Let $f:\mathcal{S}'\to \mathcal{S}$ be
 an orbifold-covering. As we just saw, we may assume 
 that all the indices $e'_i,e_j$ are $>1$. Then $f(Z')\subset Z$ and
 $f^{-1}(S)\subset S'$.
\end{para}

\begin{para}\label{para::4.1.4} Following Thurston, an orbifold is said to be
 good or {\it uniformizable} if it admits a {\it global} orbifold chart
 (\ie with $V=\ovl{S}$).  In that case, there exists a {\it universal} 
 global orbifold chart, called the {\it universal orbifold-covering} (it
 is indeed an orbifold-covering in the above sense); its underlying space
 is a simply-connected complex manifold if $\ovl{S}$ is connected.
 
 \medskip\noindent {\it Example.} By a classical result due to Fox,
 $(\P^1, (\zeta_i,e_i)_{i=0,\ldots, n})$ is uniformizable if and only if
 $n\geq 2$ or ($n=1$ and $e_0=e_1$). Under these conditions, 
 there exist algebraic global charts which are finite coverings of $\P^1$. 
\end{para}

\begin{para}\label{para::4.1.5} It is proved in {\it loc.~cit.} that orbifold
 charts are locally unique. Apart from a subset of codimension $2$,
 they are given in local coordinates around $\phi^{-1}(Z_i)$ by:
 $(z_1,z_2,\ldots, z_n)\mapsto (z_1^{e_i},z_2,\ldots, z_n)$.
 
 \noindent It follows that there exist orbifold charts of the form
 $W\to V=W/G$ ({\it loc.~cit.} 14.6), where $W$ is an
 open subset of $\C^n$ stable under a finite group $G$ of linear
 transformations.
\end{para}

\begin{para}\label{para::4.1.6} Let $\mathcal{S}=(\ovl{S}, (Z_i,e_i))$ be a
 connected orbifold, and let $\Gamma$ be a discrete group acting
 faithfully and properly by automorphisms on $\mathcal{S}$. The {\it orbifold
 quotient} $\mathcal{S}/\Gamma$ is defined as follows: the underlying
 analytic space is $\ovl{S}/\Gamma$. The ramifications divisors are the
 images of the $Z_i$'s, and the irreducible divisors $D$ with
 non-trivial fixer $\Gamma_D$; the
 latter are endowed with the multiplicity $e_D.| \Gamma_D|$.

 Note that the projection $\mathcal{S}\to \mathcal{S}/\Gamma$ is an
 orbifold-covering.
\end{para}

\begin{para}\label{para::4.1.7} One of the main applications of complex
 orbifolds is to represent moduli problems for isomorphism classes of
 objects with finite automorphism group. We have already alluded to this
 in II.\ref{sub:1.1}.  Here we shall content ourselves with examining
 the case of complex elliptic curves.

 \noindent Let $H$ be the subgroup of $SL_3(\R)$ consisting
 of matrices of the form
 $\begin{pmatrix}
 1&x_1&x_2\\ 0&a&b\\ 0&c&d
  \end{pmatrix}$ (this is a
 semi-direct product $SL_2(\R)\ltimes\R^2$).  It acts on the product
 $\mathfrak{h}$ 
 $\times\C$ of the Poincar\'e upper half plane $\mathfrak{h}$ and $\C$ by
 the rule 
 \begin{equation*}
  \begin{pmatrix}
   1&x_1&x_2\\0&a&b\\0&c&d
  \end{pmatrix}.
  (\tau, v)=\Bigl(\frac{a\tau+b}{c\tau +d},\frac{v+x_1\tau+x_2}{c\tau+d}\Bigr).
 \end{equation*}
 As is well-known, elliptic curves $E$ together with a symplectic isomorphism
 $\mathrm{H}_1(E,\Z)\cong \Z^2$ are classified by points of the
 Poincar\'e upper half plane $\mathfrak{h}$, the universal object being
 the quotient of $\mathfrak{h}\times\C$ by the action of the subgroup 
 $1.\Z^2\subset H$ (mapping to $\mathfrak{h}$ through the first projection). 

 Let $\Ker(SL_2(\Z)\to SL_2(\Z/N\Z))$ be the full congruence group of
 level $N$, $\Gamma(N)$ its image in $PGL(2)$, and $\tilde\Gamma(N)$ be
 the subgroup 
 \begin{equation*}
  \Ker(SL_2(\Z)\to SL_2(\Z/N\Z)).\Z^2\subset H.
 \end{equation*}
 For $N\geq 3$, the group $\tilde \Gamma(N)$ is
 torsion-free and acts freely on $\mathfrak{h}\times\C$; elliptic curves 
 together with a level $N$ structure (\cf II.\ref{para:1.1.2}) are
 classified by $\mathfrak{h}/\Gamma(N)$, the universal object being
 $(\mathfrak{h}\times\C)/\tilde \Gamma(N)\to \mathfrak{h}/\Gamma(N)$. 
 
 This still holds for $N=1$ or $2$, provided one interprets these
 quotients as orbifold quotients. For $N=2$, the element 
  $\begin{pmatrix}1&&0\\ &-1&\\ 0&&-1 \end{pmatrix} \in
 \tilde \Gamma(N)$ of order $2$ has fixed divisors in
 $\mathfrak{h}\times\C$: the quotient of
 $\{(\tau,x_1\tau+x_2) \mid x_1,x_2\in (\frac{1}{2}\Z)^2\}$ by $\tilde
 \Gamma(N)$, which is a disjoint union $\coprod_1^4 Z_i$ of four
 `horizontal' curves (horizontal with respect to the projection to 
 $\mathfrak{h}/\Gamma(2)=\P^1\setminus\{0,1,\infty\}$); they have the
 index $e_i=2$ (the order of the fixer).
 
This orbifold is uniformizable: the Legendre elliptic pencil provides a
global orbifold chart. Note however that, as an analytic variety,
$\mathfrak{h}\times\C/\tilde \Gamma(2)$ is a not an elliptic surface but
a $\P^1$-bundle.

The case $N=1$ is similar: $\mathfrak{h}/\Gamma(1)=\A^1$ is the $j$-line. It should be considered as an orbifold, with
divisor $\{j=0 \,{\rm(index \;3)} ,\;j=1728\,{\rm(index \;2)} \}$. In addition to the image of the horizontal divisors
$Z_i$ on $\mathfrak{h}\times\C/\tilde \Gamma(1)$, there are three vertical divisors with index $2$, and two vertical
divisors with index $3$.

\end{para}

\begin{para}\label{para::4.1.8} More generally, let us consider a Shimura
 variety of PEL type as in II.\ref{para:1.2.4}.  Recall from {\it
 loc.~cit.} the commutative diagram
\begin{equation*}
\begin{CD}
 \tilde S @>\mathcal{P}>> \mathcal{D}\\
 @V\mathcal{Q}VV @VV\mathcal{Q}V\\
 S @>\mathcal{P}>> \mathcal{D}/\Gamma
\end{CD}
\end{equation*}
where $S$ is a component of the Shimura variety, the horizontal maps are
 {\it isomorphisms}, and $\Gamma$ is a congruence subgroup. If $\Gamma$
 has torsion, it is appropriate to consider $S$ as an orbifold quotient 
 (defined by $\mathcal{Q}$ as in \ref{para::4.1.6}). 
\end{para}

\subsection{Orbifold fundamental groups.}\label{sub::4.2}

\begin{para}\label{para::4.2.1} Another important application of orbifolds is
 the Galois theory of {\it ramified} coverings.   

  Let $\mathcal{S}=(\ovl{S}, (Z_i,e_i))$ be a uniformizable connected
orbifold. We denote by $\tilde{\mathcal{S}}$ the universal
orbifold-covering (which is a simply-connected complex manifold,
together with a (ramified) covering morphism $\tilde{\mathcal{S}}\to
\ovl{S}$). We fix a base-point $s$ of $S =\ovl{S}\setminus \,(\bigcup
Z_i\cup\ovl{S}^{\mathrm{sing}})$, and a point $\tilde s$ of $\tilde{\mathcal{S}}$
above $s$. Then the pair $(\tilde{\mathcal{S}},\tilde s)$ is unique up to
unique isomorphism.
\end{para}

\begin{para}\label{para::4.2.2} Each divisor $Z_i$ defines a conjugacy class
of elements $\gamma_i$ in $\pi_1^{\mathrm{top}}(S,s)$: the local
monodromy at $Z_i$ based at $s$.

\medskip The {\it orbifold fundamental group}
$\pi_1^{\mathrm{orb}}(\mathcal{S},s)$ is the quotient of
$\pi_1^{\mathrm{top}}(S,s)$ by the {\it normal subgroup}
$\an{\gamma_i^{e_i}}$ generated by the $\gamma_i^{e_i}$'s.

 \medskip This definition is justified by the following fact
 (\cf \cite[14.10]{calip}): 
 one has a commutative square of functors 
 \begin{equation*}
  \begin{CD}
   \{\text{orbifold-coverings of $\mathcal{S}$}\} @>\subset>>
   \{\text{\'etale (= topological) coverings of $S$}\}\\
   @V\wr VV  @VV\wr V \\
   \{\text{$\pi_1^{\mathrm{orb}}(\mathcal{S},s)$-sets}\}
   @>\subset >>
   \{\text{$\pi_1^{\mathrm{top}}(S,s)$-sets}\}
  \end{CD}
 \end{equation*}
 where the top horizontal functor is the (fully faithful) restriction
 functor, and the vertical 
 functors are equivalences of categories. 

 The orbifold fundamental group may also be identified with the group of
 automorphisms of the (ramified) covering 
 $\tilde{\mathcal{S}}\to \ovl{S}$. Galois orbifold-coverings correspond to
 normal subgroups of $ \pi_1^{\mathrm{orb}}$.

 \medskip\noindent {\it Example.} Let $\mathcal{S}=(\P^1, (\zeta_i,
 e_i)_{i=0,\ldots, n})$ with $n\geq 2$ or ($n=1$ and $e_0=e_1$). Then
 $\pi_1^{\mathrm{orb}}(\mathcal{S},s)$ has generators $\gamma_0,\ldots, \gamma_n$
 and relations $\gamma_0^{e_0}=\ldots =\gamma_n^{e_n}= \Pi \gamma_i = 1$.
 In particular, it depends only, up to isomorphism, on the collection of
 integers $e_i$, and not on the $\zeta_i's$.

 \noindent The universal orbifold-covering is $\P^1, \A^1$ or the open
 unit disk $\matheur{D}$ according to whether $\sum (1-\frac{1}{e_i})$ is $< 2,\;=2,\;$
 or $>2$ respectively. 
 
\medskip\noindent {\it Remark.} The {\it finite} orbifolds coverings are
 classified by a profinite group
 $\pi_1^{\mathrm{alg.orb}}(\mathcal{S},s)$, which 
 is actually the profinite completion of
 $\pi_1^{\mathrm{orb}}(\mathcal{S},s)$.  When $\ovl{S}$ is 
 algebraic, this profinite group has a purely algebraic definition.

\end{para}

\begin{para}\label{para::4.2.3} Let $\phi:\,\ovl{S}'\to \ovl{S}$ be a global orbifold chart, such that the induced map $S'\to S$ is a
Galois (topological) covering with group $G$. Let $s'$ be a point of $S'$ above $s$. Then one has an exact sequence
\begin{equation*}
1\to \pi_1^{\mathrm{top}}(\ovl{S}',s')\to \pi_1^{\mathrm{orb}}(\mathcal{S},s) \to G \to 1.
\end{equation*}

\end{para}

\subsection{Uniformizing differential equations.}\label{sub::4.3}

\begin{para}\label{para::4.3.1} We follow the exposition in \cite[4,
 5.2]{fde} and \cite[4]{aeoaafg}. We now assume that $\ovl{S}$ is a 
complex smooth projective {\it curve}. Then $Z=\{\zeta_0,\ldots
 ,\zeta_n\}$ is just a finite set of points and $S=\ovl{S}\setminus Z$.

 We assume that $\mathcal{S}$ is uniformizable. By Riemann's
 uniformization theorem, the universal orbifold-covering
 $\tilde{\mathcal{S}}$ is then isomorphic to either
 $\P^1, \A^1$ or the open unit disk $\matheur{D}$.  We fix such an
 isomorphism\footnote{in some cases, it will be more 
 convenient to work with $\mathfrak{h}$ instead of $\matheur{D}$.}; the
 group $\Aut 
 \tilde{\mathcal{S}}\;$ is then a Lie subgroup of $PSL_2(\C)$. The
 orbifold fundamental group $\pi_1^{\mathrm{orb}}(\mathcal{S},s)$ is a discrete
 subgroup of $\Aut\, \tilde{\mathcal{S}}$ acting properly on
 $\tilde{\mathcal{S}}$, and the quotient space is
 $\tilde{\mathcal{S}}/\pi_1^{\mathrm{orb}}(\mathcal{S},s)=\ovl{S}$.      
\end{para}

\begin{para}\label{para::4.3.2} We fix a non-constant rational function $z\in
 \C(\ovl{S})$. We shall assume for simplicity that $z$ has no 
 pole on $S$. 

 \noindent Having identified $\tilde{\mathcal{S}}$ with a domain of
 $\P^1$ allows us to consider the canonical map 
 \begin{equation*}
  \tilde S\to \tilde{\mathcal{S}}\subset \P^1
 \end{equation*}
 as an analytic\footnote{in the case $\tilde{\mathcal{S}}=\P^1$,
 it is only meromorphic a priori. But in that case $\ovl{S}=\P^1$ and
 the covering $\P^1\to \P^1$ is finite. It will be convenient to
 normalize things so that $\infty$ lies above $\infty$ and $z$ is the
 standard coordinate on $\P^1$. With this choice, $\tau$ has no pole.} 
 function $\tau$ on $\tilde S$ (in loose terms, this is the multivalued
 inverse of the projection
 $\tilde{\mathcal{S}}\to
 \ovl{S}=\tilde{\mathcal{S}}/\pi_1^{\mathrm{orb}}(\mathcal{S},s)$,
 restricted to $S$). Notice that the derivation 
 $\frac{d}{dz}$ extends canonically to $\mathcal{O}(\tilde S)$. 
\end{para}

\begin{pro}\label{pro::4.3.3} There is a unique differential equation of
 the form
 \begin{equation*}
  y''- qy=0,\;\;\; q\in \C(\ovl{S}),
 \end{equation*}
 the {\rm uniformizing differential equation}, such that $\tau$ is the
 quotient of two solutions in $\mathcal{O}(\tilde S)$ of this
 differential equation. It is fuchsian and has singularities only at
 $Z=\{\zeta_0,\ldots ,\zeta_n\}$.  The monodromy group pointed at $s$ is
 a discrete finitely generated subgroup of $SL_2(\C)$ whose image in
 $PSL_2(\C)$ is $\pi_1^{\mathrm{orb}}(\mathcal{S},s)$.

 \hfill\break If $\ovl{S}=\P^1$ and $\zeta_n=\infty$, the exponents are
 rational numbers given by the Riemann scheme:
 \begin{equation*}
  \begin{pmatrix}
   \zeta_0 & \ldots & \zeta_{n-1} & \zeta_n=\infty \\
   &&&\\
   \frac{1}{2}(1+\frac{1}{e_0})& \ldots &\frac{1}{2}(1+\frac{1}{e_{n-1}})&
   \frac{-1}{2}(1+\frac{1}{e_{n}}) \\
   \\
   \frac{1}{2}(1-\frac{1}{e_0})& \ldots &\frac{1}{2}(1-\frac{1}{e_{n-1}})&
   \frac{-1}{2}(1-\frac{1}{e_{n}}) \\
  \end{pmatrix}.
 \end{equation*}
\cf {\it loc.~cit.}
\end{pro}

\begin{proof}
 Existence: $q$ is nothing but the Schwarz derivative
 \begin{equation*}
  \{\tau, z\}=
   \Bigl(\frac{d\tau}{dz}\Bigr)^{1/2}\frac{d^2}{dz^2}
   \Bigl[\Bigl(\frac{d\tau}{dz}\Bigr)^{-1/2}\Bigr],
 \end{equation*} 
 which is single-valued on $\ovl{S}$. A basis of solutions of the
 uniformizing differential equation is given by
 $(\frac{d\tau}{dz})^{-1/2}$ and $\tau.(\frac{d\tau}{dz})^{-1/2}$. The
 fact that $q\in \C(\ovl{S})$ follows immediately from the invariance of
 the Schwarz derivative under homographic transformations.

 Unicity: this is a purely algebraic fact: $q$ is determined by
 the projective connection attached to the differential equation $y''- qy=0$. 

 The determination of the Riemann scheme comes from this
 expression and the fact that around $\zeta_{i<n}$ (\resp around
 $\infty$) $\tau$ is given by $z-\zeta_i\mapsto c_i.(z-\zeta_i)^{1/e_i} +
 h.o.t.$ (\resp $1/z\mapsto c_n/z^{e_n} + h.o.t.$).
\end{proof}

\begin{rems}
  (i) The case of three singularities
 will be detailed in the next section. In the case of more than three
 singularities, it is extremely difficult in general to write down
 explicitly the uniformizing differential equations.

 (ii) In \cite{fde}, an extension of the theory of
 uniformizing differential equations to higher dimensions is
 developed. Its practical scope seems limited however to the case when
 $\tilde{\mathcal{S}}$ is a symmetric domain.

 (iii) The projective equivalence class of the uniformizing
 differential equation is a variant of what is called the Schwarz
 structure or canonical indigenous $\P^1$-bundle over $\ovl{S}$, \cf 
 \cite{gunning67:_special_rieman}, \cite{mochizuki96}.
\end{rems}

\subsection{Non-archimedean orbifolds.}\label{sub::4.4}

\begin{para}\label{para::4.4.1} We now work over a complete non-archimedean
 field $K$. The definition of orbifolds is parallel to the one 
 which we have selected in the complex situation. Being mainly interested
 in orbifolds which are quotients of 
 manifolds by finite groups, it is natural to replace everywhere \'etale
 (= topological) coverings in the complex case by {\it temperate}
 \'etale coverings in the non-archimedean context.

 \medskip\noindent A {\it $K$-orbifold} $\;\mathcal{S}=(\ovl{S},
 (Z_i,e_i))$ consists in 
 \begin{itemize}
  \renewcommand{\labelitemi}{\normalfont\bfseries\textendash}
  \item a normal $K$-analytic space $\ovl{S}$ (in the sense of
	Berkovich),
  \item a locally finite collection of irreducible divisors $(Z_i)$,
  \item a positive integer $e_i$ attached to each $Z_i$.
 \end{itemize}

 \noindent It is assumed that $\ovl{S}$ is covered by so-called
 {\it orbifold charts}, \ie morphisms
 $W \xrightarrow{\phi} V \subset \ovl{S} $ where 

\begin{itemize}
 \renewcommand{\labelitemi}{\normalfont\bfseries\textendash}
 \item $W$ is a $K$-manifold,
 \item $\phi$ is a (ramified) covering in the sense of \ref{def::1.2.3}
       (\ie $V$ is covered by open subsets $U$ such that 
       $\coprod V_j = \phi^{-1}U \to U$ and the restriction of $\phi$ to
       every $V_j$ is finite),
 \item except for a closed subset $\Sigma$ of codimension $2$ in $V$,
       $\phi$ restricts to a {\it temperate} \'etale covering 
       above $V\setminus (V\cap\bigcup Z_i)$, and is {\it ramified with index
       $e_i$ above $Z_i$}.
\end{itemize}  

 \medskip Here again, we may and shall assume that that $\Sigma $ is
 contained in $\bigcup Z_i\cup \ovl{S}^{\mathrm{sing}}$ and set
 $Z=\bigcup Z_i\cup \ovl{S}^{\mathrm{sing}},\;S=\ovl{S}\setminus Z$
 (which is a $K$-manifold).
\end{para}

\begin{rems}
  (i) At this point, because the notion of orbifold is local, it would
 not change if we had modified the definition of orbifold charts by
 replacing `temperate \'etale' by `finite \'etale' or by `\'etale'.  In
 contrast, we insist that orbifold charts are defined in terms of
 temperate coverings.
 
 (ii) Since we only prescribe the ramification indices and not
 the ramification groups, this notion of orbifold is 
 only reasonable in the tame case, \ie when $\ch(K) \nmid e_i$.  We 
 shall be mostly interested in the characteristic $0$ case.

 (iii) If $\dim \ovl{S} = 1$, then $S=\ovl{S}\setminus \bigcup Z_i$.
\end{rems}

\begin{para}\label{para::4.4.2} A {\it morphism of orbifolds}
 $f:\mathcal{S}'\to \mathcal{S}$ is a morphism $\ovl{S}'\to \ovl{S}$, such
 that $\ovl{S}'$ is covered by orbifold charts  $\phi': W'\to V'$ with
 the property that $f\circ \phi'$ factors through some orbifold chart
 $\phi:W\to V$ of $\mathcal{S}$: $f\circ \phi'=\phi\circ f'$.

 \medskip\noindent {\it Example.} Let $(\ovl{S}, (Z_i,e_i)_{i\in I})$ be
 an orbifold, and let $J\subset I$ be the subset of indices for which
 $e_i>1$. Then the identity of $\ovl{S}$ induces an isomorphism
 $(\ovl{S}, (Z_i,e_i)_{i\in I})\cong (\ovl{S}, (Z_i,e_i)_{i\in J})$. In
 other words, one can `forget' the divisors with index $e_i=1$.

 \noindent Considering $K$-manifolds as orbifolds with an empty
 collection of divisors ($S=\ovl{S}$), one can embed {\it the category of 
 $K$-manifolds as a full subcategory of the category of $K$-orbifolds}.

\noindent It will be convenient to attach to any irreducible divisor
 $D\subset \ovl{S}$ an index $e_D\in \Z_{>0}$ which is $e_i$ if
 $D=Z_i$ and $1$ if $D$ is none of the $Z_i$'s.
\end{para}

\begin{para}\label{para::4.4.3} A morphism $f$ is called an
{\it orbifold-covering} if it is a covering in the sense of
 \ref{def::1.2.3} and 
 for any orbifold chart $\phi': W'\to V'$ of $\mathcal{S}'$, $f\circ
 \phi': W'\to V=f(V')$ is an orbifold chart of $\mathcal{S}\;$ (note that if
 $\mathcal{S}'$ and $\mathcal{S}$ are just $K$-manifolds, an orbifold-covering
 amounts to a temperate \'etale covering, by our definition of orbifold
 charts).
 
 \noindent Orbifold-coverings of $\mathcal{S}$ form a category in an obvious
 way, denoted by $\Cov_\mathcal{S}$.
 
\medskip\noindent The restriction of $f$ above $S$ is then an \'etale
 covering: indeed, this is clear except above a closed subset of
 codimension $\geq 2$ in $S$, thus above the whole of $S$ by purity since
 $S$ is smooth.
 However, we do not know in general whether it
 is temperate (due to the non-local nature of the definition of
 temperate coverings).
 
\medskip\noindent {\it Example.} Assume that $\ch(K) \nmid n$. Then
 $(\A^1, (0, e')) \xrightarrow{z\mapsto z^n} (\A^1,
 (0, e)) $ is a morphism of orbifolds if and only if $e| ne'$, and is
 an orbifold-covering if and only if $e=ne'$.   

 \medskip \noindent {\it Remark.} Let $f:\mathcal{S}'\to \mathcal{S}$ be an
 orbifold-covering. As we just saw, we may assume that all the indices
 $e'_i,e_j$ are $>1$. Under this assumption, $f(Z')\subset Z$ and
$f^{-1}(S)\subset S'$.
\end{para}

\begin{para}\label{para::4.4.4} An orbifold is said to be {\it
 uniformizable}\footnote{we avoid the terminology `good' in 
 order to prevent confusion with good spaces in the sense of Berkovich,
 \cf last section.} if it admits a {\it global} orbifold chart.
 
 However, {\it universal global orbifold charts do not exist}; this is
 no surprise: there is no universal temperate covering of a
 $K$-manifold, except in trivial cases. For the same reason, {\it
 orbifold charts are not locally unique}.
 
 \medskip Let $f:\mathcal{S}'\to \mathcal{S}$ be an orbifold-covering of a
 uniformizable orbifold $\mathcal{S}$. Then the {\it restriction of $f$ to
 $f^{-1}(S)$ is a temperate \'etale covering} of $S$.

 Indeed, let $W\to \ovl{S}$ be a global orbifold chart. Then
 $W_{|S}\times_S f^{-1}(S)\to f'(W'_{|S}\times_S f^{-1}(S))= f^{-1}(S)$
 is an orbifold chart for $\mathcal{S}'$, hence $W_{| S}\times_S
 f^{-1}(S)\to S$ is an orbifold chart for $\mathcal{S}$. Since its image is
 $S$ itself, it is a temperate covering of $S$, of which $f^{-1}(S)\to S$ is 
 a quotient.

\medskip\noindent {\it Example.} Assume that $K$ is algebraically closed
of characteristic zero. Then $(\P^1, (\zeta_i, e_i)_{i=0,\ldots, n})$
is uniformizable if and only if $n\geq 2$ or ($n=1$ and $e_0=e_1$),
just as in the complex case. This actually follows from the complex
 case: one may assume that $K$ has the same cardinality as $\C$, choose
 an abstract isomorphism $K\cong \C$, and use the fact that in the
 complex case, one can find a global algebraic orbifold chart which is a
 finite ramified covering of $\ovl{S}$.
\end{para}

\begin{para}\label{para::4.4.5} Let us again assume that $K$ is of
 characteristic zero. Let $\mathcal{S}=(\ovl{S}, (Z_i,e_i))$ be a
 connected $K$-orbifold, and let $\Gamma$ be a virtually free discrete group acting
 faithfully and properly by automorphisms on $\mathcal{S}$.  The {\it
 orbifold quotient} $\mathcal{S}/\Gamma$ is defined as in the complex
 situation: the underlying analytic space is $\ovl{S}/\Gamma$. The
 ramification divisors are the images of the $Z_i$'s, and the
 irreducible divisors $D$ with non-trivial fixer $\Gamma_D$; the latter
 are endowed with the multiplicity $e_D.| \Gamma_D|$.
 Note that the projection $\mathcal{S}\to \mathcal{S}/\Gamma$ is an
 orbifold-covering.

 Interesting examples with $\ovl{S}$ open in $\P^1$ are discussed in
 \cite{van92:_discr_mumfor}. 
\end{para}

\begin{para}\label{para::4.4.6} One of the applications of $K$-orbifolds would
 be to represent moduli problems for isomorphism classes of objects with
 finite automorphism groups.

 Typical instances occur in the context of $p$-adic uniformization of
 Shimura varieties. Let us briefly mention how.  Recall from
 II.\ref{sub:7.2}, \ref{sub:7.3}
 the following basic commutative diagram

 \begin{equation*}
  \begin{CD}
   \mathcal{M}^0 @>\mathcal{P}>> \mathcal{D}\\
   @V\mathcal{Q}VV @VV\mathcal{Q}V\\
   S @>\mathcal{P}>> \mathcal{D}/\Gamma
  \end{CD}
 \end{equation*}
$S$ is a tube-like domain in the $p$-adic
 analytification of a Shimura variety, the horizontal maps are \'etale,
and $\Gamma$ is an arithmetic subgroup in the semi-simple group
$\ovl{G}^{\mathrm{ad}}$. Under the assumption that $\Gamma$ is small
enough, the maps $\mathcal{Q}$ are topological covering maps.

If one drops that assumption, one has to consider $S$ as an orbifold
quotient (defined by $\mathcal{Q}$ as in \ref{para::4.4.5}). This is
 especially relevant 
in the case of global uniformization, \ie when $S$ is a connected
component of the Shimura variety.

The fake projective planes of Mumford, Ishida, Kato, provide interesting
examples of this situation, \cf II.\ref{sub:7.5} and \cite{opuofpp}.
We shall see several other instances in the sequel.
\end{para}

\subsection{Non-archimedean orbifold fundamental groups.}\label{sub::4.5}

\medskip From now on, we assume that $K$ is {\it algebraically closed of
characteristic zero}. 

\begin{para}\label{para::4.5.1}  Let $\mathcal{S}=(\ovl{S}, (Z_i,e_i))$
 be a uniformizable connected orbifold.

As we have seen in \ref{para::4.4.4}, the restriction above $S$ of an
 orbifold-covering is a temperate \'etale covering. This gives 
rise to a restriction functor 
\begin{equation*}
   \Cov_\mathcal{S} \longrightarrow \Cov_S^{\mathrm{temp}}.
\end{equation*}
\end{para}

\begin{pro}\label{pro::4.5.2} This functor is fully faithful. Its image
 consists of those temperate \'etale coverings $S'\to S$ which extend to
 ramified coverings $\ovl{S}'\to \ovl{S}$, with $\ovl{S}'$ normal, such
 that the ramification index at each irreducible divisor above $Z_i$
 divides $e_i$ (for every $i$).

 \noindent This image is stable under
 taking connected components, fiber products (over $S$), quotients,
 finite disjoint unions, and Galois closure.
\end{pro}

\begin{proof}
 From Kiehl's theorem on the existence of tubular neighborhoods
 and the example at the end of \ref{para::4.4.3}, it follows that
 $\mathcal{S}'$ is 
 determined by $\ovl{S}'$ and by the ramification index $n_{ij}$ at each
 divisor $Z_{ij}$ above $Z_i$ (for every $i$): namely,
 $e'_{ij}=e_i/n_{ij}$. The first assertion then follows immediately from
 theorem \ref{thm::2.1.10}.  The second assertion is a straightforward
 consequence 
 of this characterization of the image.
\end{proof}

 \medskip 
 Let us fix a geometric base-point $\ovl{s}$ of $S
 =\ovl{S}\setminus
 \,(\bigcup Z_i\cup\ovl{S}^{\mathrm{sing}})$ (in dimension one, we could also
 choose a tangential base-point in the sense of 2.2).

\begin{dfn}\label{def::4.5.3} The orbifold fundamental group pointed at
 $\ovl{s}$ 
\begin{equation*}
 \pi_1^{\mathrm{orb}}(\mathcal{S},\ovl{s})
\end{equation*}
 is the topological group attached to $\Cov_\mathcal{S}$ viewed as a
 subcategory of $\Cov_S^{\mathrm{et}}$, \cf 1.4. 
\end{dfn}
This is a prodiscrete group (\ref{cor::1.4.7}). By \ref{thm::1.4.5}, one
has a commutative 
square of functors
\begin{equation*}
 \begin{CD}
  \left\{
  \begin{matrix}
   \text{disjoint unions of}\\
   \text{orbifold-coverings of $S$}
  \end{matrix}
  \right\}
    @>\subset >>
  \left\{
  \begin{matrix}
   \text{disjoint unions of}\\
   \text{temperate coverings of $S$}
  \end{matrix}
  \right\}\\
  @V\wr VV @VV\wr V\\
  \{\pi_1^{\mathrm{orb}}(\mathcal{S},s)\text{-sets}\}
  @>\subset >>
  \{\pi_1^{\mathrm{temp}}(S,s)\text{-sets}\}.
 \end{CD}
\end{equation*}
 Galois orbifold-coverings correspond to open normal subgroups of $
 \pi_1^{\mathrm{orb}}$.

\medskip\noindent{\it Example}. For $r\in | K^\times|$, we set $\matheur{D}_r=
\matheur{D}(0,r^+),\;\matheur{D}_r^\ast=\matheur{D}_r\setminus \{0\},$ and consider the orbifold
$\mathcal{S}_e = (\matheur{D},(0,e))$ (multiplicity $e$ at the origin). Then by
the same calculation as in \ref{pro::2.3.7}, one finds

\begin{equation*}
\pi_1^{\mathrm{orb}}(\mathcal{S}_e)=\pi_1^{\mathrm{alg}}(\matheur{D})\times \Z/e\Z 
\end{equation*} 
(we drop the base point in the notation since the group is abelian, \cf
\ref{pro::1.4.4}).

This can also be expressed in the following way: let $\gamma$ be a
topological generator of the $\widehat{\Z}(1)$-factor in
$\pi_1^{\mathrm{temp}}(\matheur{D}^\ast)$ (the factor corresponding to Kummer
coverings). Then \begin{equation*} \pi_1^{\mathrm{orb}}(\mathcal{S}_e
)=\pi_1^{\mathrm{temp}}(\matheur{D}^\ast,\ovl{s})/\an{\gamma^e}^- , \end{equation*}
where $\an{\gamma^e}^-$ denotes the closure of the (normal) subgroup
generated by $\gamma$.

\begin{para}\label{para::4.5.4} By Kiehl's theorem, $\ovl{S}$ is locally
 isomorphic to a product $\matheur{D}^\ast\times Z_i^{\mathrm{smooth}}$ in the
 neighborhood 
 of the non-singular part of $Z_i$; whence a homomorphism
 $\pi_1^{\mathrm{temp}}(\matheur{D}^\ast)\to \pi_1^{\mathrm{temp}}(S,\ovl{s})$,
 well-defined up to conjugation. Let us denote by $\gamma_i$ the image in
 $\pi_1^{\mathrm{temp}}(S,\ovl{s})$ of a topological generator of the
 $\widehat{\Z}(1)$-factor of $\pi_1^{\mathrm{temp}}(\matheur{D}^\ast)$ (``local
 monodromy'' at $Z_i$). The closure $\an{\gamma_i^{e_i}}^-$ of the normal
 subgroup of $\pi_1^{\mathrm{temp}}(S,\ovl{s})$ generated by the
 $\gamma_i^{e_i}$'s is a well-defined closed subgroup, independent of the
 choice of the $\gamma_i$'s.
\end{para}

\begin{thm}\label{thm::4.5.5}
 \begin{enumerate}
  \renewcommand{\theenumi}{\alph{enumi}}
  \item One has $\displaystyle
	\pi_1^{\mathrm{temp}}(S,\ovl{s}) \cong \limproj_{(e_i)}
	\pi_1^{\mathrm{orb}}((\ovl{S},(Z_i,e_i)),\ovl{s}),$ where the
	limit is taken according to the ordering by divisibility.
  \item The homomorphism $\pi_1^{\mathrm{temp}}(S,\ovl{s})\to
	\pi_1^{\mathrm{orb}}(\mathcal{S},\ovl{s})$ induces an injective strict
	homomorphism
	\begin{equation*}
	 \pi_1^{\mathrm{temp}}(S,\ovl{s})/\an{\gamma_i^{e_i}}^-\,\to \,
	  \pi_1^{\mathrm{orb}}(\mathcal{S},\ovl{s})
	\end{equation*}
	with dense image. This is an isomorphism if $S$ is the analytification
	of a smooth algebraic $K$-variety $S^{\mathrm{alg}}$; in particular
	$\pi_1^{\mathrm{temp}}(S,\ovl{s})/\an{\gamma_i }^-\,\cong \,
	\pi_1^{\mathrm{temp}}({\ovl{S}},\ovl{s}) $.
 \end{enumerate}
\end{thm}
 
\begin{proof}
  (a) reflects the fact that $\Cov_S^{\mathrm{temp}}= \liminj
 \Cov_{(\ovl{S},(Z_i,e_i)),\ovl{s})}$ which follows from the previous
 proposition.

  (b) The previous proposition shows more precisely that a
 $\pi_1^{\mathrm{temp}}(S,\ovl{s})$-set comes from a 
 $\pi_1^{\mathrm{orb}}(\mathcal{S},\ovl{s})$-set if and only if the
 $\gamma_i^{e_i}$'s act trivially. 
 
 \noindent We shall use lemma \ref{lem::1.4.9}. By construction
 of the topology, open subgroups of
 $\pi_1^{\mathrm{temp}}(S,\ovl{s})$ and
 $\pi_1^{\mathrm{orb}}(\mathcal{S},\ovl{s})$ respectively form a fundamental
 system of neighborhoods of 
 $1$. In fact, this is even true for normal open subgroups
 (existence of Galois closures). It follows that 
 $\pi_1^{\mathrm{temp}}(S,\ovl{s})/\an{\gamma_i^{e_i}}^-$ has the same
 property. The first part of b) then follows from
 \ref{lem::1.4.9}. (ii). The 
 second part follows from \ref{lem::1.4.9}. (iii) and (iv), since
 we have 
 shown that the quotient of $\pi_1^{\mathrm{temp}}(S,\ovl{s})$ by
 any closed normal subgroup (\eg $\an{\gamma_i^{e_i}}^-$) is
 complete when $S$ is algebraic (\ref{pro::2.1.6}).
\end{proof}

\begin{rems}
 (a) \ref{pro::2.3.9} is a special case of the
 last assertion of \ref{thm::4.5.5}. 

 (b) When $S$ is algebraic, the profinite completion of
 $\pi_1^{\mathrm{orb}}(\mathcal{S},\ovl{s})$ admits an 
 algebraic description. 

 (c) The situation is analogous to the classical one (over $\C$, \cf 4.2),
 with at least two noteworthy differences: 

 \begin{itemize}
 \item since there is no universal orbifold-covering, the geometrically
       meaningful object is not $\pi_1^{\mathrm{orb}}$ itself, but its
       discrete quotients,
 \item the non-archimedean $\pi_1^{\mathrm{orb}}$ is a finer invariant
       which takes into account the relative ``position'' of the $Z_i$. 
       For instance, if $\ovl{S}= \P^1_{\C_p}$, the classical
       $\pi_1^{\mathrm{orb}}$ 
       depends only on the numbers $e_i$, not on the branch 
       points $Z_i$ themselves. In the $p$-adic case, this is not the case. For
       instance, if $Z=\{0,1,\lambda,\infty\}$ and every 
       $e_i=2$, then the Legendre elliptic curve with parameter $\lambda$
       provides a global orbifold chart, 
       and it follows from \ref{para::2.3.2} that the orbifold fundamental
       group is either
       \begin{itemize}
	\item a (compact) semi-product of $\Z/2\Z$ by $\widehat{\Z}\times
	      \widehat{\Z}$ if $| \lambda(1-\lambda)| =1$,
	\item a (non-compact) semi-product of $\Z/2\Z$ by $ \Z\times
	      \widehat{\Z}$ if $| \lambda(1-\lambda)| \neq 1$.
       \end{itemize}
 \end{itemize}
\end{rems}

\begin{para}\label{para::4.5.6} Let $\phi:\,\ovl{S}'\to \ovl{S}$ be a global
 orbifold chart, such that the induced map $S'\to S$ is a 
finite Galois \'etale covering with group $G$. Let $\ovl{s}'$ be a
 geometric point of $S'$ above $\ovl{s}$.
\end{para}

\begin{pro}\label{pro::4.5.7} The universal topological covering
 $\widetilde{\ovl{S}'}$ of ${\ovl{S}'}$ is a Galois orbifold-covering of
 $\mathcal{S}$. Let $\Gamma= \Aut(\widetilde{\ovl{S}'}/\mathcal{S})$
 denote the corresponding (discrete) Galois group. Then one has a
 commutative diagram of strict homomorphisms of 
 topological groups, in which the rows are exact sequences
 \begin{equation*}
  \begin{CD}
   1 @>>> \pi_1^{\mathrm{temp}}(\ovl{S}',\ovl{s}') @>>>
   \pi_1^{\mathrm{orb}}(\mathcal{S},\ovl{s}) @>>> G @>>> 1\\
   @. @VVV @VVV @VVV @.\\
   1 @>>> \pi_1^{\mathrm{top}}(\ovl{S}',\ovl{s}') @>>> \Gamma @>>> G @>>> 1.
  \end{CD}
 \end{equation*}
\end{pro}

\begin{proof}
 Let us denote by $\mathcal{S}'$ the orbifold $\ovl{S}'$ with all 
 $\phi^{-1}(Z_i)$ endowed with multiplicity 
 $e_i=1$. The pull-back by $\phi$ provides the vertical arrows of a
 commutative diagram of functors
 \begin{equation*}
 \begin{CD}
  \Cov_\mathcal{S} @>\subset >> \Cov_S^{\mathrm{temp}}\\
  @VVV @VVV\\
  \Cov_{\mathcal{S}'}=\Cov_{\ovl{S}'}^{\mathrm{temp}}@>\subset>>
  \Cov_{S'}^{\mathrm{temp}}.
 \end{CD}
 \end{equation*}
 One concludes by \ref{cor::1.4.12} that there is an exact sequence

 \begin{equation*}
 1\to \pi_1^{\mathrm{temp}}(\ovl{S}',\ovl{s}')\to
  \pi_1^{\mathrm{orb}}(\mathcal{S},\ovl{s}) \to G \to 1.
 \end{equation*}
 On the other hand, it follows from \ref{pro::1.1.3} that the restriction
 of $\widetilde{\ovl{S}'}$ above $S'$ is nothing 
 but the universal topological covering $\tilde{  S'}$ of ${S'}$. By
 \ref{lem::2.1.2}, $\tilde S'\to S$ is a Galois temperate 
 covering, and it follows from that and \ref{pro::4.5.2} that
 $\widetilde{\ovl{S}'}\to \mathcal{S}$ is a Galois orbifold-covering. It
 is then clear that the corresponding epimorphism
 $\pi_1^{\mathrm{orb}}(\mathcal{S},\ovl{s})\to \Gamma$ gives rise to the above
 commutative diagram.
\end{proof} 

\begin{para}\label{para::4.5.8} Assume that $S$ is the analytification of a
 smooth algebraic $K$-variety $S^{\mathrm{alg}}$. According to
 \ref{thm::3.5.12}, 
 discrete representations of $\pi_1^{\mathrm{temp}}(S,\ovl{s})$ give rise
 to algebraic vector bundles with regular integrable 
 connection $(M,\nabla)$ on $S^{\mathrm{alg}}$. 
\end{para}

\begin{pro}\label{pro::4.5.9} An algebraic vector bundle with regular integrable
connection $(M,\nabla)$ on $S^{\mathrm{alg}}$ comes from a discrete representation $\rho$ of $\pi_1^{\mathrm{orb}}(\mathcal{S},\ovl{s})$ if and
only if there is a finite global orbifold chart $\phi:\,\ovl{S}'\to \ovl{S}$ such that $(\phi^\ast(M))^\nabla$ extends to a
locally constant sheaf on $\ovl{S}'$.   
\end{pro}

This follows from \ref{para::3.5.6} and \ref{thm::4.5.5}. 

\medskip\noindent {\it Remark.} If $S'/S$ is Galois, one can restrict
$\rho$ to the open normal subgroup $ 
\pi_1^{\mathrm{temp}}(\ovl{S}',s')$ of
$\pi_1^{\mathrm{orb}}(\mathcal{S},\ovl{s})\,$ (\ref{pro::3.5.7}), and this restriction
factors through $\pi_1^{\mathrm{top}}(\ovl{S}',s')$.

\subsection{$p$-adic uniformizing differential equations.}\label{sub::4.6}

\begin{para}\label{para::4.6.1} We now develop the non-archimedean counterpart
 of 4.3. We assume that $\ovl{S}$ is a 
 smooth projective {\it curve} over $K$. Then
 $Z=\{\zeta_0,\ldots ,\zeta_n\}$ is just a finite set of points and
 $S=\ovl{S}\setminus Z$.

 We assume moreover that $\mathcal{S}$ is
uniformizable. To imitate the complex theory,
we shall assume that {\it there exists a finite Galois global orbifold
 chart $\phi: \ovl{S}'\to \ovl{S}$ such that $\ovl{S}'$ 
is a Mumford curve} (of any genus $\geq 0$). 

The universal topological covering $\widetilde{\ovl{S}'}$ of
 ${\ovl{S}'}$ will replace the universal 
 orbifold-covering, which does not exist in the non-archimedean
 situation. By assumption, there is an embedding  
  \begin{equation*}
\widetilde{\ovl{S}'} \subset \P^1,
\end{equation*}
 which we fix, so that the group $\Aut\, \widetilde{\ovl{S}'}\;$
 becomes a Lie subgroup of $PGL_2(K)$.

The orbifold-covering $\widetilde{\ovl{S}'}\to \ovl{S}$ gives rise to a
 projective discrete representation 
\begin{equation*}
\ovl{\rho}:\,\pi_1^{\mathrm{orb}}(\mathcal{S},s) \to PGL_2(K).
\end{equation*}
 Its image $\Gamma$ is contained in the normalizer in $PGL_2(K)$ of the 
Schottky group $\pi_1^{\mathrm{top}}({\ovl{S}'},s')$, hence discrete; it
 acts properly on $\widetilde{\ovl{S}'}$, with quotient space 
$\widetilde{\ovl{S}'}/\Gamma=\ovl{S}$.  
\end{para}

\begin{para}\label{para::4.6.2} We fix a non-constant rational function
 $z\in K(\ovl{S})$, and assume for simplicity that $z$ has no 
 pole on $S$. 
 
 \noindent Having identified $\widetilde{\ovl{S}'}$ with a domain of
 $\P^1$ allows us to consider $z: \tau \mapsto z(\tau)$ as a meromorphic
 function defined on that domain (and analytic on 
 $\widetilde{  S'}\subset \widetilde{\ovl{S}'}$). Notice that the derivation
 $\frac{d}{dz}$ extends canonically to $\mathcal{O}(\widetilde{S'})$.  
\end{para}

\begin{pro}\label{pro::4.6.3} There is a unique differential equation of
 the form
 \begin{equation*}
  y''- qy=0,\;\;\; q\in K(\ovl{S}),
 \end{equation*}
 the {\rm uniformizing differential equation}, such that a local inverse
 $\tau(z)$ of the map $z(\tau)$ is given by the quotient of two solutions in
 $\mathcal{O}(\widetilde{ S'})$ of this differential equation. It is
 fuchsian and has singularities only at $Z=\{\zeta_0,\ldots ,\zeta_n\}$. 

 \hfill\break If $\ovl{S}=\P^1$ and $\zeta_n=\infty$, the exponents are
 rational numbers given 
 by the Riemann scheme:
\begin{equation*}
 \begin{pmatrix}
  \zeta_0 & \ldots & \zeta_{n-1} & \zeta_n=\infty \\
  &&&\\\frac{1}{2}(1+\frac{1}{e_0})&
  \ldots &\frac{1}{2}(1+\frac{1}{e_{n-1}})& \frac{-1}{2}(1+\frac{1}{e_{n}}) \\
  \\
  \frac{1}{2}(1-\frac{1}{e_0})&
  \ldots &\frac{1}{2}(1-\frac{1}{e_{n-1}})& \frac{-1}{2}(1-\frac{1}{e_{n}}) \\
 \end{pmatrix}.
\end{equation*}
\end{pro}

 \noindent The proof is the same as in the complex case: $q$ is nothing
 but the Schwarz derivative $\{\tau, z\}= 
 (\frac{d\tau}{dz})^{1/2}\frac{d^2}{dz^2}[(\frac{d\tau}{dz})^{-1/2}]$,
 which is single-valued on $\ovl{S}$. A basis of solutions of the
 uniformizing differential equation is given by
 $(\frac{d\tau}{dz})^{-1/2}$ and $\tau.(\frac{d\tau}{dz})^{-1/2}$.

\begin{pro}\label{pro::4.6.4} The projective connection attached to the
 uniformizing differential corresponds to the projective discrete
 representation $\ovl{\rho}$ in the sense of \ref{para::3.5.13}: it is
 in the image 
 of the temperate \'etale Riemann-Hilbert functor and its projective
 monodromy group is the discrete group $\Gamma$.
\end{pro}

\begin{proof}
 A basis of solutions of the symmetric square of the uniformizing 
 differential equation is given by $\frac{1}{{d\tau/
 dz}},\,\frac{\tau}{{d\tau/ dz}},\,\frac{\tau^2}{{d\tau/dz}} $. These 
 are elements of $\mathcal{O}(\widetilde{  S'})$ (in fact multivalued 
 meromorphic functions on $\ovl{S}'$). Criterion \ref{pro::4.5.9} then
 shows that 
 this symmetric square comes from a discrete representation of 
 $\pi_1^{\mathrm{orb}}(\mathcal{S},s) $. It is clear that this representation is 
 nothing but the composition of $\ovl{\rho}$ with the 
 standard homomorphism $PGL_2\to GL_3$. 
\end{proof}

\noindent [One could also argue by using an analogue of \ref{pro::4.5.9} for
projective connections.]
 
\medskip {\it Remark.} The connection itself attached to the
uniformizing differential equation corresponds to a discrete
representation of $\pi^{\mathrm{temp}}(S,\ovl{s})$ (\cf
\ref{pro::3.5.14}), which 
does not necessarily factor through $\pi_1^{\mathrm{orb}}(\mathcal{S},s) $ [it
factors through an orbifold fundamental group, up to replacing the
$e_i$'s by appropriate multiples].

\subsection{$p$-adic Shimura orbifolds and uniformizing differential
  equations.}\label{sub::4.7}

\begin{para}\label{para::4.7.1} We are now ready for a synthesis, relating the
 theory of $p$-adic period mappings surveyed in II with the 
theory of $p$-adic orbifolds. 

Let us consider a maximal order
$\mathcal{B}$ in an indefinite quaternion algebra
$B$ over
$\Q$, endowed with a positive involution $\ast$. ``Fake elliptic
 curves'' are $\ast$-polarized abelian surfaces with multiplication by
$\mathcal{B}$. They are classified by a smooth projective Shimura curve
$\mathcal{X}^+$ defined over $\Q$, \cf
 II.\ref{para:7.4.2}, which should actually be considered as an orbifold 
 $(\mathcal{X}^+, (\zeta_i,e_i))$, because the action of $\Gamma^+$ on
the Poincar\'e upper half-place $\mathfrak{h}$ has fixed points (the
$\zeta_i$'s are their images in $\mathcal{X}^+$).

 The Gauss-Manin connection $(\mathcal{H}, \nabla_{\mathrm{GM}})$
 attached to the 
 ``universal'' fake elliptic curve (in the sense of 
 orbifolds) is a connection with logarithmic poles at the $\zeta_i$'s. It
 has been analyzed in  II.\ref{sub:7.6}. The choice 
 of a quadratic splitting field
 $\Q(\sqrt{-d})/\Q$ for $B$ decomposes $(\mathcal{H},
 \nabla_{\mathrm{GM}})$ into 
 two pieces of rank two (the ``$+$ and $-$ 
 eigenspaces'' for the action of $\sqrt{-d}$).  It follows from the
 theory of the complex period mapping that a multivalued 
 section $\tau :\,\mathcal{X}^+(\C)\setminus 
 \{\zeta_i\}\to \mathfrak{h}$ of the projection $\mathfrak{h}\to
 \mathcal{X}^+(\C)$ is given by a quotient of two non-proportional 
solutions of $\nabla_{\mathrm{GM},+}$. One derives that {\it the uniformizing
 differential equation for $(\mathcal{X}^+,
(\zeta_i,e_i))$ is projectively equivalent to $\nabla_{\mathrm{GM},+}$}.

On the other hand, there is a canonical $\ovl{\Q}$-structure in
the $\C$-space of solutions at a given point $s$ (the Betti $\mathrm{H}^1( -
 /\ovl{\Q})$ of the fake elliptic curve with parameter $s$), which is
stable under (global) monodromy.
\end{para}
 
\begin{para}\label{para::4.7.2} According to the \v{C}erednik-Drinfeld theory
 surveyed in II.\ref{sub:7.4}, \ref{sub:7.6} (and the description of
 the $p$-adic period mapping given there as a quotient of solutions of
 the Gauss-Manin connection), there is an entirely parallel $p$-adic
 description of this geometric situation for any critical prime $p$, \ie
 any $p$ which divides the discriminant of $B$: $\mathcal{X}^+(\C_p)$ is
 a $p$-adic orbifold of the kind considered in 4.6, and {\it the
 $p$-adic uniformizing differential equation for $(\mathcal{X}^+,
 (\zeta_i,e_i))$ is projectively equivalent to
 $\nabla_{\mathrm{GM},+}$}, viewed 
 as a $p$-adic connection (\cf II.\ref{cor:7.2.6}, II.\ref{para:7.2.7}).

 By the theory of $p$-adic Betti lattices (\cf  II.\ref{para:7.6.2}),
 there is a canonical $\ovl{\Q}$-structure in the $\C_p$-space of
 solutions at a given point $s$, which is stable under (global) monodromy.
\end{para}

\begin{para}\label{para::4.7.3} Following  II.\ref{sub:7.5}, all this generalizes to
 Shimura curves defined by quaternion algebras $B$ over totally real
number fields $E$, and a ``congruence subgroup'' $\Gamma$ (here, we drop
the general assumption of II.\ref{sec:7} that $\Gamma$ is small enough). We
recall that these Shimura curves are not modular if $E\neq \Q$, but are
twisted forms of modular Shimura curves defined by quaternion algebras
$B^\bullet = B\otimes_E K$ over suitable totally imaginary quadratic
extensions $K/E$.  One can thus still attach a Gauss-Manin connection to
such a situation, and this connection is related to the uniformizing
differential equations both in the complex and in the $v$-adic
sense\footnote{if $p$ divides the level of the congruence group, the
orbifold is not quite of the kind considered in 4.6: the top space is
not the set of ordinary points for a Schottky group, but a finite
\'etale covering of such a space, \cf II.\ref{para:7.4.7}; however, the
construction of the uniformizing differential equation can still be done
in the same way, due to the fact that the finite \'etale coverings of
the Drinfeld space which occur are $PGL_2(E_v)$-equivariant.}, for {\it
critical places $v$ of $E$}, \cf II.\ref{sub:7.5}, II.\ref{sub:7.6} (we
fix embeddings $\C \leftarrow \ovl{\Q} \to \C_p$ for every critical
$v\mid p$, compatible with the privileged places $\infty_0$ and $v$ of
$E$, and denote by $\ovl{B}$ the quaternion algebra over $E$ obtained
from $B$ by interchanging the local invariants at $\infty_0$ and
$v$). We summarize these results as follows, using the terminology of
orbifolds:

\end{para}

\begin{thm}\label{thm::4.7.4} The Shimura curve $\mathcal{X}_\Gamma$ is
 defined over $\ovl{\Q}$, and may be considered, in a canonical way,
 both as a complex orbifold and as a $p$-adic orbifold (for any critical
 $v\mid p$).

 \noindent The complex uniformizing differential equation is defined over
 $\ovl{\Q}$ and coincides with the $p$-adic uniformizing differential
 equation. It is projectively equivalent to a 
 direct summand of the Gauss-Manin connection: the period mapping (\resp
 the $p$-adic period mapping) is given by a suitable quotient
 $\tau $ (\resp $\tau_p$) of two non-proportional solutions. The
 projective monodromy group of the uniformizing differential 
 equation over $\C$ (\resp over $\C_p$) is $\Gamma$ (\resp is a discrete
 arithmetic group $\Gamma_p$ in $PGL_2(\C_p)$ 
 attached to the quaternion algebra $\ovl{B}$).

 \noindent Moreover, there is a canonical $\ovl{\Q}$-structure in the
 space of solutions of the complex (\resp $p$-adic) uniformizing
 differential equation at any $\ovl{\Q}$-point, which is stable under
 global monodromy.
\end{thm}
Let us recall that the $p$-adic projective monodromy group, based at a
geometric point $\ovl{s}$, can be understood in a 
sense close to the usual intuitive one, \cf \ref{para::3.5.6}. In loose
words: there is a smooth projective curve $S$ and a finite 
orbifold-covering
$S\to \mathcal{X}_\Gamma$ such that the projective monodromy is given by
the values taken by $\tau_p$ when one passes 
from a point $ s'$ lying above $ s$ to another $ s''$ by all possible
paths on the Berkovich space $S_{\C_p}^{\mathrm{an}}$.

\begin{para}\label{para::4.7.5} Let us illustrate this section with a concrete
 example (\cf \cite{poaptg} for more details). We first 
define a quadrangle
$ABCD$ with circular edges in $\mathfrak{h}$, as follows 
\begin{figure}[h]
 \begin{picture}(250,250)(0,-20)
  \put(0,-20){\includegraphics[scale=0.75,keepaspectratio,clip]{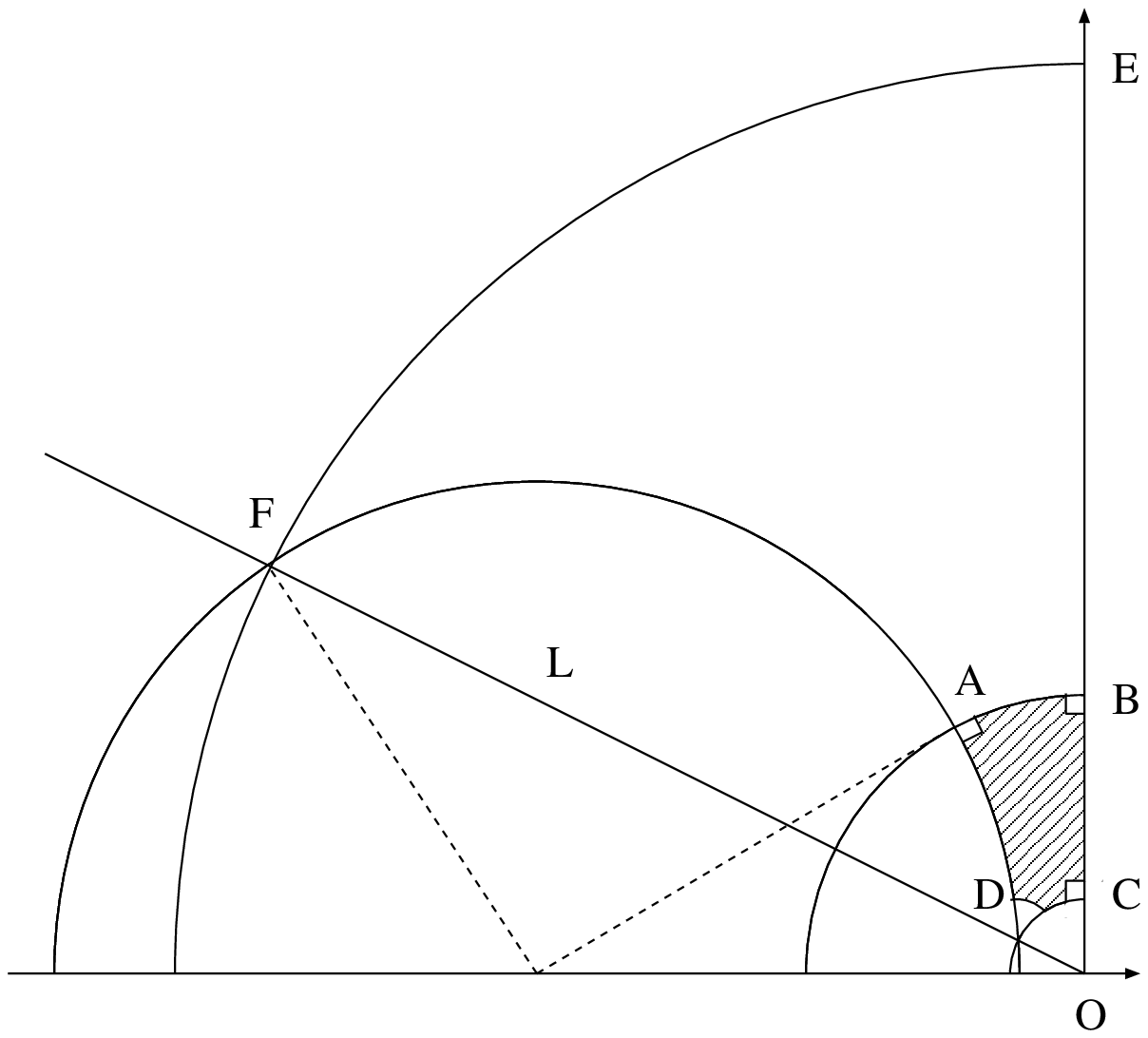}}
  \put(242,7){\makebox(0,0)[c]{{\tiny $\Sigma'$}}}
  \put(202,32){\makebox(0,0)[c]{{\footnotesize $\Sigma$}}}
  \put(180,80){\makebox(0,0)[c]{{\footnotesize $\Sigma''$}}}
  \put(182,-12){\makebox(0,0)[c]{{\small $-1$}}}
 \end{picture}
 \caption{}
\end{figure}
 
 $L$ is the half line of slope $-1/2$ through the origin,
 
 $\Sigma$ (\resp $\Sigma'$) is the unit half-circle centered at the
 origin (\resp at $2-\sqrt 3$),

 $D$ is the point $L\cap \Sigma'$,

 $\Sigma''$ is the half-circle centered on the horizontal axis,
 orthogonal to $\Sigma$ and passing through $D$, 

 $B$ (\resp $C$) is the intersection point of the vertical axis with
 $\Sigma$ (\resp $\Sigma'$), 

 $A$ is the intersection point of $\Sigma$ with $\Sigma'$ (specifically,
 $A=- \frac{1+2i}{\sqrt 5}$).

\medskip The quadrangle $ABCD$ is a hyperbolic quadrangle with angles
 $(\frac{\pi}{2},\frac{\pi}{2},\frac{\pi}{2},\frac{\pi}{6})$. The
 symmetries around the edges of this polygon generate a group of
 M\"obius transformations which contains a 
 subgroup $\Gamma^\ast$ of index $2$ of conformal transformations. This
 fuchsian group has already attracted the attention 
 of many a mathematician (J.F. Michon, M.F. Vigneras
 \cite[p. 123]{adadq}, F. Beukers, D. Krammer \cite{aeoaafg}, N. Elkies 
 \cite{elkies98:_shimur}, D. Kohel... ). It is $PSL_2(\R)$-conjugated to the 
 discrete group generated by the following unimodular 
matrices (up to sign)
\begin{equation*}
\begin{pmatrix}
 0&1\\-1&0
\end{pmatrix},\,
\begin{pmatrix}
 0&-\frac{\sqrt 5+1}{2}\\
 \frac{\sqrt 5-1}{2}&0
\end{pmatrix},\,
\begin{pmatrix}
 -\sqrt 3& 1+\sqrt 5 \\ 1-\sqrt 5&\sqrt 3
\end{pmatrix},\,
\begin{pmatrix}
 \frac{\sqrt 3+\sqrt {15}}{2}& -2\\2&\frac{\sqrt 3-\sqrt {15}}{2}
\end{pmatrix}
\end{equation*}
of respective orders $2,2,2,6$, whose product in this order is $1$. The
 orbifold quotient of $\mathfrak{h}$ by $\Gamma^\ast$ has 
been determined by Krammer --- and, later but independently, by Elkies: it is 
\begin{equation*}
\bigl(\P^1; (0,e_0=2), (1,e_1=2), (81, e_{81}=2), (\infty , e_\infty =
 6)\bigr),
\end{equation*}
 the points $A,B,C,D$ being mapped to $0,1,81,\infty$ respectively. 

Moreover, these authors have computed the uniformizing
equation\footnote{it was also found, earlier, by D. and G. Chudnovsky 
\cite[p.193]{chudnovsky89:_comput}, who did not give 
the proof.}. It is actually easier to write down the differential 
equation obtained from the uniformizing equation by 
multiplying all solutions by $P^{1/4}$, where $P=z(z-1)(z-81)$.  It is 
\begin{equation*}
 18P. y'' + 9 P' y' + (z-9). y =0.\tag{$\ast$}
\end{equation*}
\end{para}

\begin{para}\label{para::4.7.6} In general, it is quite hard to compute
 uniformizing differential equations (due to the presence of
 accessory parameters, when there are
 more than three singularities). The method of Krammer and Elkies takes
 advantage of the fact that the orbifold
 under consideration is actually a Shimura orbifold (the method, which
 seems to go back to H.P.F. Swinnerton-Dyer, exploits 
 the Hecke correspondences). 

 In fact, $\Gamma^\ast$ coincides with the group denoted by the same
 symbol in  II.\ref{para:7.4.2}, namely $\{g\in B^\times \mid 
 g\mathcal{B}=\mathcal{B}g, \Nr(g)>0\}$, for the quaternion algebra $B$ of
 discriminant $3.5$ over $\Q$  \cite[IV.3.B,C]{adadq}
 (\cf \cite[8]{poaptg} for more details on the coincidence of $\Gamma$
 with the group denoted by $W^+$ in \cite{aeoaafg}, up to 
 conjugation), and is the orbifold fundamental group of the Shimura
 orbifold $\mathcal{X}^\ast$:
 \begin{equation*}
  \mathcal{X}^\ast(\C)= \mathfrak{h}/\Gamma^\ast.
 \end{equation*} 
\end{para}

\begin{para}\label{para::4.7.7} Let us now consider the critical prime
 $5$ (similar phenomena occur for $p=3$). Note that there is a 
``confluence'' between the singularities $1$ and $81$ in characteristic
 $5$. Let $\ovl{B}$ be the (definite) quaternion 
algebra with discriminant $3$, and let $\ovl{\mathcal{B}}$ be a maximal
 order ($\cong \Z\bigl[1,i,\frac{i+j}{2},\frac{1+ij}{2}\bigr]$ with 
$i^2=-1,j^2=-3,ij=-ji\,$). Let $\Gamma_5^\ast$ be the image in
 $PGL_2(\Q_5)$ of the normalizer $\{g\in \ovl{B}^\times \mid 
 g\ovl{\mathcal{B}}[\frac{1}{5}]=\ovl{\mathcal{B}}[\frac{1}{5}]g\}$. 
 This is a discrete subgroup, and one has the \v{C}erednik uniformization 
\begin{equation*}
 \mathcal{X}^\ast(\C_p)= \Omega_{\C_5}/\Gamma_5^\ast
\end{equation*} 
 by the Drinfeld space (II.\ref{para:7.4.3}), which represents
 $\mathcal{X}^\ast$ as a $p$-adic orbifold.

\medskip It follows from what precedes that viewed as a $5$-adic
 differential equation, $(\ast)$ is obtained from the $p$-adic
differential uniformizing equation attached to $\mathcal{X}^\ast$ by
multiplying all solutions by $(z(z-1)(z-81))^{1/4}$. In particular, it
has projective global monodromy group $\Gamma_5^\ast$. This group is
generated by the {\it conjugacy classes} of local monodromies at
$0,1,81$, as follows from \ref{pro::2.3.9}\footnote{however, it is not
clear that 
it is generated by three such elements, as is stated erroneously in
loc.~cit.}.
\end{para}

\newpage
\section{$p$-adic triangle groups.}\label{sec::5}
\begin{abst}
In loose words, $p$-adic triangle groups are projective monodromy groups
of those hypergeometric differential equations which are in the image of
the \'etale Riemann-Hilbert functor. We give a purely geometric
description of them, and then construct all ``arithmetic'' $p$-adic
triangle groups using the \v{C}erednik-Drinfeld uniformization
of Shimura curves.
\end{abst}

\subsection{Review of Schwarz triangle groups.}\label{sub::5.1}

\begin{para}\label{para::5.1.1} A {\it Schwarz orbifold} is a complex orbifold
 of the form 
 \begin{equation*}
  \mathcal{S}(e_0,e_1,e_\infty)=(\P^1; (0,
   e_0), (1,e_1), (\infty, e_{\infty})).
 \end{equation*}
 It is uniformizable. The universal orbifold-covering is $\P^1, \A^1$ or 
 the open unit disk $\matheur{D}$ according to whether $\sum  \frac{1}{e_i}$ is
 $>1,\;=1,\;$ or $<1$ respectively.

 \noindent We fix a base-point $\ovl{s}$ (which can be a tangential
 base-point).
 
\medskip\noindent The {\it Schwarz triangle group\footnote{in conformity
 with the previous chapter, we do 
 not allow the $e_i$'s to be infinite. Thus we restrict our attention to
 the so-called cocompact triangle groups.}} $\Delta(e_0,e_1,e_\infty)$
 may be defined as $\pi_1^{\mathrm{orb}}(\mathcal{S}(e_0,e_1,e_\infty),s)$.

\noindent As an abstract group, it has generators 
$\gamma_0,\gamma_1, \gamma_\infty$ and relations
 $\gamma_0^{e_0}=\gamma_1^{e_1}=\gamma_\infty^{e_\infty}= \Pi \gamma_i = 
1$. A permutation of $e_0,e_1,e_\infty)$ does not change the group, only
 the orbifold by a homography interchanging 
 $0,1$ and $\infty$.  
\end{para}

\begin{para}\label{para::5.1.2} The uniformizing differential equation is
 fuchsian with singularities at $0,1,\infty$, and Riemann scheme 
\begin{equation*}
\begin{pmatrix}
 0 &  1 & \infty \\
 &&&\\
 \frac{1}{2}(1+\frac{1}{e_0})& \frac{1}{2}(1+\frac{1}{e_{1}})&
 \frac{-1}{2}(1+\frac{1}{e_{\infty}}) \\ 
 &&& \\
 \frac{1}{2}(1-\frac{1}{e_0})& \frac{1}{2}(1-\frac{1}{e_{1}})&
 \frac{-1}{2}(1-\frac{1}{e_{\infty}})
\end{pmatrix} .
\end{equation*}

 \medskip Multiplying solutions of this equation by  
\begin{equation*}
 z^{\frac{1}{2}(-1+\frac{1}{e_0})}(z-1)^{\frac{1}{2}(-1+\frac{1}{e_1})}
\end{equation*}
transforms this equation into a (projectively equivalent) differential equation of hypergeometric type
\begin{equation*}
 HGDE(a,b,c)\;\;\;\;\;\;\;\;\;\;\;\;\;\; z(z-1)y''+(c-(a+b+1)z)y'-ab y=0,  
\end{equation*}
with $1-c=\frac{1}{e_0},\;c-a-b=\frac{1}{e_1},\; \pm(a-b)=\frac{1}{e_\infty}$, whose Riemann scheme is    
 \begin{equation*}
  \begin{pmatrix}
   0 &  1 & \infty \\
   &&&\\
   0& 0& \frac{1}{2}(1-\frac{1}{e_0}-\frac{1}{e_1}+\frac{1}{e_{\infty}}) \\
   &&&\\
   \frac{1}{e_0}&\frac{1}{e_{1}}&
   \frac{1}{2}(1-\frac{1}{e_0}-\frac{1}{e_1}-\frac{1}{e_{\infty}})
  \end{pmatrix} .
 \end{equation*}
 The triangle group $\Delta(e_0,e_1,e_\infty)$ may thus be interpreted
 as the {\it projective monodromy group pointed at $\ovl{s}$ of this
 hypergeometric differential equation}.
\end{para}

\begin{para}\label{para::5.1.3} The quotient $\tau$ of a suitable basis of
 solutions of $ HGDE(a,b,c)$ sends the upper half-plane bijectively 
onto a curvilinear triangle $\mathcal{T}$ with vertices
 $\tau(0),\tau(1),\tau(\infty)$ and respective angles 
$\pi/e_0,\pi/e_1,\pi/e_\infty$. By the Schwarz reflection principle, the
 map $\tau$ can be extended to the lower half-plane 
through each of the ``segments'' $]0,1[$, $]1,\infty[$, $]\infty,0[$ by
 conformal reflections $s_\infty,s_0,s_1$ respectively. The products 
$\gamma_0=s_1s_\infty,\,\gamma_1=s_\infty s_0,\,\gamma_\infty=s_0s_1 $
 give the standard generators of $\Delta(e_0,e_1,e_\infty)$. By
 repeating these reflections, the occuring curvilinear triangles do not
 overlap (except on their edges) and provide a tiling of the universal
orbifold-covering ($\P^1, \A^1$ or $\matheur{D}$).
\end{para}

\begin{para}\label{para::5.1.4} A similar situation arises for any
 hypergeometric differential equation with rational parameters
 $(a,b,c)$.  However, if some of the differences of the exponents
$1-c,c-a-b,a-b$ is not the inverse of an integer, the triangles which
occur by iterating Schwarz reflections may overlap. Moreover, if the
local projective monodromies are of respective order $e_0, e_1,
e_\infty$ at $0,1,\infty$, the global projective monodromy group need
not be $\Delta(e_0,e_1,e_\infty)$. This is related to the fact that a
homomorphic image of $\Delta(e_0,e_1,e_\infty)$ which preserves the
order of $\gamma_0,\gamma_1,\gamma_\infty$ need not be isomorphic to
$\Delta(e_0,e_1,e_\infty)$ (\cf \cite{greenberg81:_homom_psl} for a
detailed discussion 
of homomorphic images of triangle groups\footnote{we are indebted to
M. Kapovich for this reference.}).
\end{para}

\begin{para}\label{para::5.1.5} The case $\sum  \frac{1}{e_i}>1$ (universal
 covering $\P^1$) corresponds to the Schwarz list of $HGDE$ with finite
 (projective) monodromy groups: the list of triples
 $(e_0,e_1,e_\infty)$ is, up to permutation:

 \medskip
 $(2,2,n)$ corresponding to a dihedral group,

 $(2,3,3)$ corresponding to the tetrahedral group,

$(2,3,4)$ corresponding to the octahedral group,

 $(2,3,5)$ corresponding to the icosahedral group.
 
 \medskip In the case $\sum  \frac{1}{e_i}=1$ (universal covering $\A^1$), there is a global orbifold chart given by an elliptic curve with complex
 multiplication. There are three possibilities for $(e_0,e_1,e_\infty)$, up to permutation:
 $(2,3,6)$,
 $(2,4,4)$,
 $(3,3,3)$.
 
 \medskip All the other triples belong to the ``hyperbolic case'':
 $\Delta(e_0,e_1,e_\infty)$ is a fuchsian group.
 
 For more detail on these topics, we refer to the  good treatises
 \cite{f2}, \cite{ntatg}, \cite{hfml}.

\end{para}

\begin{para}\label{para::5.1.6} As was remarked by J. Wolfart, Belyi's
 criterion for algebraic curves to be defined over $\ovl{\Q}$ can be
 expressed in terms of Schwarz orbifolds as follows:

 a complex projective curve $\ovl{S}$ is defined over $\ovl{\Q}$ if and
 only if it is the base of an orbifold-covering $(\ovl{S}, (Z_i,e_i))$
 of some Schwarz orbifold.
\end{para}

\subsection{$p$-adic Schwarz orbifolds and triangle groups.}\label{sub::5.2}

\begin{para}\label{para::5.2.1} A {\it $p$-adic Schwarz orbifold} is an
 orbifold over some complete subfield $K$ of $\C_p$  of the form 
\begin{equation*}
\mathcal{S}(e_0,e_1,e_\infty)=(\P^1, (0, e_0), (1,e_1), (\infty, e_{\infty})).
\end{equation*}
\noindent Up to replacing $K$ by a finite extension, it is uniformizable
 (\cf example in \ref{para::4.4.4}). However, there is no universal
 orbifold-covering (\ref{para::4.4.4}).
 
We want our $p$-adic triangle groups to be {\it discrete}
groups of fractional linear $p$-adic transformations. Thus it is not
 possible in general to define them as orbifold 
fundamental groups, but only as appropriate discrete quotients. This
 leads us to imitate the alternative definition \ref{para::5.1.2} 
of triangle groups as projective hypergeometric monodromy groups. 
\end{para}

\begin{para}\label{para::5.2.2} Let $e_0, e_1, e_\infty$ be integers $>1$, and
 define  
\begin{equation*}
 a=
 \frac{1}{2}\Bigl(1-\frac{1}{e_0}-\frac{1}{e_1}+\frac{1}{e_{\infty}}\Bigr),\,
 b=
 \frac{1}{2}\Bigl(1-\frac{1}{e_0}-\frac{1}{e_1}-\frac{1}{e_{\infty}}\Bigr),\, 
 c= 1- \frac{1}{e_0}.
\end{equation*} 
 We now consider the hypergeometric differential equation 
$HGDE(a,b,c)$ as a $p$-adic differential equation.
\end{para}

\begin{dfn}\label{def::5.2.3} Assume that $ HGDE(a,b,c) $ is in the
 image of the \'etale Riemann-Hilbert functor, 
\ie the \'etale sheaf of germs
 of solutions is locally constant on $\P^1_{\C_p}\setminus\{0,1,\infty\}$
 (\cf \ref{thm::3.4.6}).  If the corresponding projective monodromy group
\begin{equation*}
\Im\bigl[\pi_1^{\mathrm{et}}(\P^1\setminus\{0,1,\infty\},\ovl{s})\to
 PGL_2(\C_p)\bigr]
\end{equation*} is discrete and finitely generated, it is then called
a {\rm $p$-adic triangle group}, and denoted by
$\Delta_p(e_0,e_1,e_\infty)$.

\medskip
In the sequel, we say that ``{\it $\Delta_p(e_0,e_1,e_\infty)$ exists}''
to express that the assumptions of \ref{def::5.2.3} are 
fulfilled (this depends on the four numbers $p,e_0,e_1,e_\infty$).

\end{dfn}

\begin{pro}\label{pro::5.2.4}
 \begin{enumerate}
  \renewcommand{\theenumi}{\alph{enumi}}
  \item If $\Delta_p(e_0,e_1,e_\infty)$ exists,
	there exists a finite Galois covering 
	$\phi:\,\ovl{S}' \to \P^1$ ramified above $\{0,1,\infty\}$, such that
	the pull-back of $ HGDE(a,b,c) $ by $\phi$ has a full 
	set of multivalued analytic solutions over
	$S'=\phi^{-1}(\P^1\setminus\{0,1,\infty\})$.
  \item If $\Delta_p(e_0,e_1,e_\infty)$ exists, the homomorphism
	\begin{equation*}
	 \Im[\pi_1^{\mathrm{et}}(\P^1\setminus\{0,1,\infty\},\ovl{s})\to
	 PGL_2(\C_p)]
	\end{equation*}
	factors through a surjective homomorphism 
	\begin{equation*}
	 \pi_1^{\mathrm{orb}}(\mathcal{S}(e_0,e_1,e_\infty)) \to
	  \Delta_p(e_0,e_1,e_\infty). 
	\end{equation*}
  \item Conversely, if the projective connection attached to $HGDE(a,b,c)$
	comes from a projective representation
	\begin{equation*}
	 \pi_1^{\mathrm{orb}}(\mathcal{S}) \to  PGL_2(\C_p)
	\end{equation*}
	with discrete image, this image is a $p$-adic triangle group
	$\Delta_p(e_0,e_1,e_\infty)$.
  \item Every $p$-adic triangle group
	$\Delta_p(e_0,e_1,e_\infty)$ is a discontinuous subgroup of
	$PGL_2(K)$ for some finite extension $K$ of $\Q_p$.
 \end{enumerate}
\end{pro}

\begin{proof}
  (a) Because the Wronskian of $ HGDE(a,b,c) $ is algebraic, the
 determinant of the representation  
 $\rho:\,\pi_1^{\mathrm{et}}(\P^1_{\C_p} \setminus
 \{0,1,\infty\},\ovl{s})\to 
 GL_2(\C_p)$ factors through a finite representation of
 $\pi_1^{\mathrm{temp}}(\P^1_{\C_p}\setminus\{0,1,\infty\},\ovl{s})$.
 Under the assumptions of \ref{def::5.2.3}, $\Im\,\ovl{\rho}$ is
 finitely 
 generated, and it follows that $\Im \,\rho$ is finitely generated. By
 \ref{pro::3.5.10}, one derives that $\rho$ factors through
 $\pi_1^{\mathrm{temp}}(\P^1\setminus\{0,1,\infty\},\ovl{s})$. This
 implies a), \cf \ref{para::3.5.6} (where a concrete description
 of the monodromy is 
 also given).

   (b) We have just seen that $\rho$ factors through
 $\pi_1^{\mathrm{temp}}(\P^1_{\C_p}\setminus\{0,1,\infty\},\ovl{s})$. To
 check that the associated projective representation factors through 
 $\pi_1^{\mathrm{orb}}(\mathcal{S}(e_0,e_1,e_\infty))$ is a matter of
 projective local monodromy --- or more plainly, of difference of
 exponents --- at $0,1,\infty$ (\cf 4.5); this is clear.

  (c) Again because the Wronskian of $HGDE(a,b,c)$ is algebraic, the
 fact that the projective connection attached to $HGDE(a,b,c)$ comes
 from a projective representation of
 $\pi_1^{\mathrm{temp}}(\P^1\setminus\{0,1,\infty\},\ovl{s})$ implies
 that $ HGDE(a,b,c) $ itself comes from a representation of
 $\pi_1^{\mathrm{temp}}(\P^1\setminus\{0,1,\infty\},\ovl{s})$ (\cf the
 discussion in \ref{para::3.5.13}). By \ref{pro::3.5.10}, the
 image is finitely generated, and 
 the assertion follows.

  (d) If $\Delta_p(e_0,e_1,e_\infty)$ exists, there exists a finite
 Galois global orbifold chart $\ovl{S}'\to \P^1$ such that pull-back of
 $ HGDE(a,b,c) $ on the universal topological covering
 $\widetilde{\ovl{S}'}$ becomes trivial (as a connection). Notice that
 $\ovl{S}'$ is a projective curve defined over some number field, and a
 fortiori over some finite extension $K/\Q_p$.  Up to replacing $K$ by a
 finite extension, we may assume that $\widetilde{\ovl{S}'}$ is also
 defined over $K$ (and that the image in $\ovl{S}'$ of some geometric
 point above $\ovl{s}$ is $K$-rational). Thus $\rho$ factors through a
 representation of
 $\pi_1^{\mathrm{temp}}(\P^1_K\setminus\{0,1,\infty\},\ovl{s})\to
 GL_2(K)$, which makes clear that $\Delta_p(e_0,e_1,e_\infty)\subset
 PGL_2(K)$. Since $K$ is locally compact, any discrete subgroup of
 $PGL_2(K)$ is discontinuous \cite[1.6.4]{sgamc}.
\end{proof}

\medskip\noindent{\it Remark.} By (b), if $\Delta_p(e_0,e_1,e_\infty)$
exists, there {\it exists} a finite Galois global 
orbifold chart $\phi:\,\ovl{S}' \to \P^1$ such that the quotient of two
non-proportional solutions of $ HGDE(a,b,c) $ 
defines a multivalued meromorphic function on $\ovl{S}'$. However, one
should not believe that {\it any} finite Galois 
global orbifold chart $\phi:\,\ovl{S}' \to \P^1$ has this property.

\noindent Indeed, we shall see below that the diadic triangle group
$\Delta_2(4,4,4)$ exists. The Fermat curve $F_4$ over 
$\C_2$ provides a finite Galois global orbifold chart, with
$\pi_1^{\mathrm{top}}(F_4)=1$ (because $F_4$ has tree-like reduction). 
But the corresponding projective hypergeometric connection does not
become trivial over $F_4$, since $\Delta(4,4,4)$ is 
infinite.  

\begin{para}\label{para::5.2.5} The case $\sum \frac{1}{e_i}>1$ corresponds
 to the list of $p$-adic $HGDE$ with finite (projective) monodromy
 groups. This is the same as the Schwarz list: $(2,2,n)$, $(2,3,3)$,
$(2,3,4)$, $(2,3,5)$, up to permutation, and yields finite $p$-adic
triangle groups for every $p$ which are isomorphic to the corresponding
finite Schwarz triangle groups.

\medskip\noindent On the other hand, there is {\it no $p$-adic triangle
 group with $\sum  \frac{1}{e_i}=1$}. Indeed, there is a global
 algebraic orbifold chart given by an elliptic curve with complex
 multiplication.  Viewed over $\C_p$, these have good reduction, and it
 follows from \ref{para::4.5.6} and \ref{para::2.3.2} that the orbifold
 fundamental group is 
 compact. Any discrete representation would then be finite, and it would
 follow that the corresponding $HGDE$ has a basis of algebraic
 solutions. This is an algebraic property, which can be translated back
 over $\C$. One finds that such a $HGDE$ would belong to the Schwarz
 list, a contradiction.

\medskip In the sequel, we restrict our attention to the ``hyperbolic
 case'' $\sum  \frac{1}{e_i}<1$. It is not clear a 
 priori that there exists any hyperbolic $p$-adic triangle group.
\end{para}

\subsection{Mumford-Schwarz orbifolds.}\label{sub::5.3}

\begin{para}\label{para::5.3.1} Our aim is to give a more geometric
 description of $p$-adic triangle groups, independent of the theory of
 monodromy of $p$-adic differential equations. This involves the notion
 of Mumford-Schwarz orbifold.
\end{para}

\begin{dfn}\label{def::5.3.2} A Mumford-Schwarz orbifold is a $p$-adic
 Schwarz orbifold such that there exists a finite 
Galois global orbifold chart $\ovl{S}' \to \P^1$, $\ovl{S}'$ being a
 Mumford curve.  
\hfill\break In other
words, it is the orbifold quotient of a Mumford curve $\ovl{S}'$ by a
 finite group of automorphisms $G$ such that $\ovl{S}'/G= \P^1$ and the
Galois covering $\ovl{S}' \to \P^1$ is ramified exactly above $0,
1,\infty$ with ramification index $e_0,e_1,e_\infty$ respectively.
\end{dfn}

\medskip\noindent According to \ref{pro::4.5.7}, one has an exact sequence
\begin{equation*}
1\to \pi_1^{\mathrm{temp}}(\ovl{S}',\ovl{s}')\to
 \pi_1^{\mathrm{orb}}(\mathcal{S}(e_0,e_1,e_\infty),\ovl{s}) \to G \to 1. 
\end{equation*}

\noindent On the other hand, because $\ovl{S}'$ is a Mumford curve,
$\pi_1(\ovl{S}',\ovl{s}')$ is a Schottky group embedded in 
$PGL_2(K) $, and the universal topological covering
$\widetilde{\ovl{S}'}$ can be 
identified with the set $\Omega\subset \P^1_{\C_p}$ of ordinary points
for $\pi_1(\ovl{S}',\ovl{s}')$. The subgroup  
$\Gamma$ of all liftings of the elements of $G \subset \Aut (\ovl{S}')$
in $PGL_2(K)$ is a subgroup of finite index inside 
the normalizer of
$\pi_1(\ovl{S}',\ovl{s}')$ in $PGL_2(K)$; it sits in an extension 

\begin{equation*}
 1\to \pi_1^{\mathrm{top}}(\ovl{S}',\ovl{s}')\to \Gamma \to G \to 1.
\end{equation*}
{\it Remark.} A Mumford-Schwarz orbifold is hyperbolic if and only if some (\resp every) curve $\ovl{S}'$
as above has genus $\geq 2$. 

\begin{pro}\label{pro::5.3.3} Let $ \mathcal{S}=(\P^1; (0,
 e_0), (1,e_1), (\infty, e_{\infty})) $ be a Mumford-Schwarz orbifold. Then the
 $p$-adic triangle group $\Delta_p(e_0,e_1,e_\infty)$ exists, and
 coincides with the discrete group $\Gamma\;$ (in particular, $\Gamma$
 does not depend on the choice of $\ovl{S}'/ \P^1$). Moreover, one has a
 commutative diagram
 
 \begin{equation*}
  \begin{CD}
   1 @>>> \pi_1^{\mathrm{temp}}(\ovl{S}',\ovl{s}') @>>>
  \pi_1^{\mathrm{orb}}(\mathcal{S}(e_0,e_1,e_\infty),\ovl{s}) @>>> G @>>> 1\\
  @. @VVV @VVV @VVV @.\\
  1 @>>>\pi_1^{\mathrm{top}}(\ovl{S}',\ovl{s}')@>>>
  \Delta_p(e_0,e_1,e_\infty)@>>> G @>>> 1
  \end{CD}
 \end{equation*}
 with surjective vertical maps and exact arrows.
\end{pro}

\begin{proof}
  Let us consider the $p$-adic uniformizing differential equation
 attached to the orbifold-covering $\widetilde{\ovl{S}'}\to
 \mathcal{S}(e_0,e_1,e_\infty)$ (3.6). It is fuchsian with singularities
 at $0,1,\infty$, and Riemann scheme
 \begin{equation*}
 \begin{pmatrix}
 0 &  1 & \infty \\
 &&&\\
 \frac{1}{2}(1+\frac{1}{e_0})&  \frac{1}{2}(1+\frac{1}{e_{1}})&
 \frac{-1}{2}(1+\frac{1}{e_{\infty}}) \\
 &&& \\
 \frac{1}{2}(1-\frac{1}{e_0})&  \frac{1}{2}(1-\frac{1}{e_{1}})&
 \frac{-1}{2}(1-\frac{1}{e_{\infty}})     
 \end{pmatrix} .
 \end{equation*}

 \medskip Multiplying solutions of this equation by  
 \begin{equation*}
 z^{\frac{1}{2}(-1+\frac{1}{e_0})}(z-1)^{\frac{1}{2}(-1+\frac{1}{e_1})}
 \end{equation*}
 transforms the equation into a (projectively equivalent) $p$-adic
 differential equation of hypergeometric type $ HGDE(a,b,c) $ with
 $1-c=\frac{1}{e_0},\;c-a-b=\frac{1}{e_1},\;
 \pm(a-b)=\frac{1}{e_\infty}$.  By \ref{pro::4.6.4}, the projective connection
 attached to $ HGDE(a,b,c) $ is the projective connection attached to to
 the projective representation given by the composed homomorphism
 $\pi_1^{\mathrm{orb}}(\mathcal{S}(e_0,e_1,e_\infty),\ovl{s})\to \Gamma \to
 PGL_2(K)$, which has discrete image. By \ref{pro::5.2.4}. (b), we
 conclude that 
 $\Delta_p(e_0,e_1,e_\infty)$ exists and coincides with $\Gamma$. The
 commutative diagram in the statement comes from \ref{pro::4.5.7}.
\end{proof}

\begin{para}\label{para::5.3.4} F. Kato calls the $p$-adic triangle groups
 which arise in this way ``$p$-adic triangle groups of Mumford type''
\cite{kato:_schwar_mumfor}. Inspired by F. Herrlich's study of automorphisms of
Mumford curves (using graphs of groups), he has developed
geometrico-combinatorial techniques to construct Mumford-Schwarz
orbifolds, and even a complete classification scheme, including a
description of $p$-adic triangle groups ``of Mumford type'' as explicit
amalgamated sums. As a by-product, he has proved the following
\end{para}

\begin{thm}[F. Kato]\label{thm::5.3.5} For $p= 2, 3$ and $5$, there are
 infinitely many hyperbolic Mumford-Schwarz orbifolds.  For $p>5$, there
 is none.
\end{thm}

\begin{para}\label{para::5.3.6} The following theorem explains how to obtains
 all $p$-adic triangle groups from $p$-adic triangle groups ``of Mumford
 type''.
\end{para}

\begin{thm}\label{thm::5.3.7} $\Delta_p(e_0,e_1,e_\infty)$ exists if and
 only if $\mathcal{S}(e_0,e_1,e_\infty)$ is a finite orbifold-covering
 of some Mumford-Schwarz orbifold
 $\mathcal{S}(e'_0,e'_1,e'_\infty)$. Moreover, the latter can be chosen
 in such a way that
 $\Delta_p(e_0,e_1,e_\infty)=\Delta_p(e'_0,e'_1,e'_\infty)$ as a subgroup
 of $PGL_2(\C_p)$.
\end{thm}

\begin{proof}
 (a) Assume that we have a finite orbifold-covering $h:
 \mathcal{S}\to \mathcal{S}'$ of a Mumford-Schwarz orbifold by a 
 Schwarz orbifold. The pull-back by $h$ of the hypergeometric projective
 connection $ HGPC(a',b',c') $ related to 
 $e'_0,e'_1,e'_\infty$ as above is the hypergeometric projective
 connection $HGPC(a,b,c) $ related to $e_0,e_1,e_\infty$. 
 Because $\mathcal{S}'$ is a Mumford-Schwarz orbifold, its uniformizing
 differential equation, hence $ HGPC(a',b',c') $, comes 
 from a representation of the temperate fundamental group, and it follows
 that so does $ HGPC(a,b,c) $. Moreover, the monodromy group of
 $HGPC(a,b,c) $ is a subgroup of finite index of the monodromy group of
 $HGPC(a',b',c') $, which is the discrete finitely generated group
 $\Delta_p(e'_0,e'_1,e'_\infty)$. Therefore the monodromy 
 group of $HGPC(a,b,c) $ is itself discrete finitely generated, and is
 $\Delta_p(e_0,e_1,e_\infty)$.  

 (b) Conversely, assume that $\Delta_p(e_0,e_1,e_\infty)$ exists. Let 
 $\Sigma \to \P^1$ be the ramified covering corresponding to the quotient 
 $\Delta_p(e_0,e_1,e_\infty)$ of
 $\pi_1^{\mathrm{orb}}(\mathcal{S},\ovl{s}) $. By 
 definition of $\Delta_p(e_0,e_1,e_\infty)$ as projective monodromy
 group, any quotient of two non-proportional solutions of 
 $HGDE(a,b,c)$ defines an \'etale map
 \begin{equation*}
 \tau: \;\Sigma \to \P^1_{\C_p}
 \end{equation*} which is equivariant under
 $\Delta_p(e_0,e_1,e_\infty)$.

 On the other hand, $\Delta_p(e_0,e_1,e_\infty)$ being
 a finitely generated discontinuous subgroup of $PGL_2(K)$ for some
 finite extension $K/\Q_p$, it admits a normal 
 subgroup of finite index which is a Schottky group $\Gamma$. Let
 $\Omega\subset \P^1_{\C_p}$ denote the set of ordinary 
 points for $\Gamma$. Since $\Delta_p(e_0,e_1,e_\infty)$ normalizes
 $\Gamma$, it preserves $\Omega$ and the quotient 
 $H=\Delta_p(e_0,e_1,e_\infty)/\Gamma$ is a group of automorphisms of the
 Mumford curve $\Omega/\Gamma$. We denote by $T$ 
 the quotient of this Mumford curve by $H$ (a smooth projective curve).

 We claim that
 \begin{equation*}
  \tau(\Sigma) \subset \Omega.
 \end{equation*}
 Indeed, otherwise $\tau$ would have a value $\tau(s)$ which is a limit
 point for $\Gamma$, \ie an accumulation point of an infinite sequence
 of points $\gamma_n(x)$, $x\in \P^1$. Since $\tau$ is \'etale, hence
 open, we may assume that $x\in \Im\,\tau$, say $x=\tau(t)$. By
 $\Gamma$-equivariance of $\tau$, we see that $\tau(s)$ is an
 accumulation point of the infinite sequence $\tau(\gamma_n(t))$. This
 contradicts the fact that the set $\{\gamma_n(t)\}$ is discrete in
 $\Sigma$ and that $\tau$ is \'etale (hence locally invertible around
 $s$).

 Next, being $\Delta_p(e_0,e_1,e_\infty)$-equivariant, $\tau$ induces an
 analytic map
 \begin{equation*}
  \ovl{\tau}: \, \Sigma/\Delta_p(e_0,e_1,e_\infty)=\P^1
   \to \Omega/\Delta_p(e_0,e_1,e_\infty)=T.
 \end{equation*}
 This map is necessarily algebraic and finite, and $T\cong \P^1$ (we
 normalize the latter isomorphism by imposing that
 $\ovl{\tau}(\{0,1,\infty\}\subset \{0,1,\infty\}$). Notice that because
 $\tau$ is $\Delta_p(e_0,e_1,e_\infty)$-equivariant, 
 it follows that
 \begin{equation*}
  \tau(\Sigma)= \Omega 
 \end{equation*} 
 and that $\tau$ is a finite \'etale covering (not a topological
 covering, unless $\Sigma=\Omega$). We have a cartesian 
 square 
 \begin{equation*}
 \begin{CD}
  \Sigma @>>> \Omega \\
  @VVV @VVV\\
  \P^1 @>>> \P^1
 \end{CD}
 \end{equation*}
 in which the vertical maps are Galois
 orbifold-coverings with group $\Delta_p(e_0,e_1,e_\infty)$, and define
 orbifold quotients $\mathcal{S}$ and $\mathcal{S}'$ of $\Sigma$ and
 $\Omega$ respectively. Since $\tau$ is \'etale, we see that
 $\mathcal{S}$ is actually a finite orbifold-covering of
 $\mathcal{S}'$. This forces $\mathcal{S}'$ to be a Schwarz orbifold
 (like $\mathcal{S}$), hence a Mumford-Schwarz orbifold by construction.
\end{proof}

\begin{cor}\label{cor::5.3.7}
 \begin{enumerate}
  \renewcommand{\theenumi}{\alph{enumi}}  
  \item Assume that $\Delta_p(e'_0,e'_1,e'_\infty)$ is a $p$-adic
	triangle group ``of Mumford type''. Assume that
	$\Delta(e_0,e_1,e_\infty)$ is a subgroup of finite index of the
	corresponding ``complex'' triangle group
	$\Delta(e'_0,e'_1,e'_\infty)$. Then $\Delta_p(e_0,e_1,e_\infty)$ 
	exists, and can be realized as a subgroup of finite index of
	$\Delta_p(e'_0,e'_1,e'_\infty)$.
  \item If $\Delta_p(e_0,e_1,e_\infty)$ exists, and if the
	corresponding ``complex'' triangle group $\Delta(e_0,e_1,e_\infty)$ is
	maximal in its commensurability class of triangle groups, then
	$\Delta_p(e_0,e_1,e_\infty)$ is ``of Mumford type''.
 \end{enumerate}
\end{cor}
\begin{proof}
 (a) If $\Delta(e_0,e_1,e_\infty)$ is a subgroup of finite index of
 $\Delta(e'_0,e'_1,e'_\infty)$, there is a finite morphism
 $\P^1\to \P^1$, defined over a number field, which underlies a
 finite orbifold-covering $\mathcal{S}(e_0,e_1,e_\infty)\to
 \mathcal{S}(e'_0,e'_1,e'_\infty)$ in the complex sense, but also
 in the $p$-adic sense (for every embedding of the number field
 in $\C_p$). The assertion comes from the ``if'' part of the
 theorem.

 (b) If $\Delta_p(e_0,e_1,e_\infty)$ exists, then, similarly, by the
 ``only if'' part of \ref{para::5.3.6}, there is a finite
 morphism $\P^1\to 
 \P^1$ defined over a number field, which underlies a finite
 orbifold-covering $\mathcal{S}(e_0,e_1,e_\infty)\to
 \mathcal{S}(e'_0,e'_1,e'_\infty)$ in the $p$-adic sense, and
 also in the complex sense. Hence
 $\Delta(e_0,e_1,e_\infty)\subset \Delta(e'_0,e'_1,e'_\infty)$
 with equality if and only if
 $\mathcal{S}(e_0,e_1,e_\infty)=\mathcal{S}(e'_0,e'_1,e'_\infty)$.
\end{proof}

\begin{cor}\label{cor::5.3.8} There is no infinite $p$-adic triangle
 group for $p>5$.
\end{cor}
\begin{proof}
 This follows from \ref{thm::5.3.5} and \ref{para::5.3.6}.  
\end{proof}

\subsection{Arithmetic triangle groups. Review of Takeuchi's
  list.}\label{sub::5.4}

\begin{para}\label{para::5.4.1} (Cocompact) arithmetic triangle groups
 are hyperbolic triangle groups which are {\it commensurable with the
 unit groups of quaternion algebras $B$ over number fields}.

 In order to make this definition more precise, let us first place
 ourselves in greater generality. Let $E$ be a number field, and let $B$
 be a quaternion algebra over $E$. Let $V$ be a finite set of places of
 $E$ containing the places at infinity and at least one place (possibly
 infinite) where $B$ is unramified. Let us define
 \begin{equation*}
  G_{(V)} = \prod \,B_v^\times/E_v^\times = \prod \, PGL_2(E_v)
 \end{equation*} where the product runs over all places $v\in V$ which
 are unramified for
 $B$. Let $\mathcal{O}_{E}(V)$ be the ring of so-called $V$-integers of
 $E$ (elements which are integers at the places not in $V$),
 let $\mathcal{B}(V)$ be a maximal $\mathcal{O}_{E}(V)$-order in
 $B$, and let $PSL_1(\mathcal{B}(V))$ be the image in $G_{(V)}$ of the
 group of elements of $\mathcal{B}(V)$ of reduced norm $1$. 
 This is a cocompact discrete subgroup of $ G_{(V)}$, but the projection
 of $PSL_1(\mathcal{B}(V))$ in any proper factor of 
 $G_{(V)}$ is dense \cite[IV.1.1]{adadq}. 

 Any subgroup of $G_{(V)}$ which is commensurable to
 $PSL_1(\mathcal{B}(V))$ is called an {\it arithmetic subgroup of $G_{(V)}$}. 

 \noindent In particular, arithmetic subgroups of $PGL_2(\R)$ are
 subgroups which are commensurable with the unit groups
 $PSL_1(\mathcal{B})$ of quaternion algebras over {\it totally real}
 number fields $E$, which are {\it split at exactly one place at
 infinity $\infty_0$} (and one can take $V=\{\infty_0\}$).
\end{para}

\begin{para}\label{para::5.4.2} The complete list of cocompact
 arithmetic triangle groups (there are $76$ such groups, up to
 permutation of $e_0,e_1,e_\infty$), and of the corresponding quaternion
 algebras, has been established by K. Takeuchi; moreover, he has shown
 that two arithmetic triangle groups are commensurable with each other
 if and only if they come from the same quaternion algebra. Here is the
 list \cite[Table 1]{ccoatg}:

\renewcommand{\arraystretch}{1.5}
\begin{longtable}{lcp{0.63\textwidth}}
 $E$ & $\mathrm{disc}(B)$ & $ (e_0,e_1,e_\infty)\text{ up to permutation}$ \\
 \hline
 $\Q$ & $2.3$ & $(2,4,6)\, (2,6,6)\, (3,4,4)\, (3,6,6)$ \\
 $\Q(\sqrt 2)$ & $v_2$ &
 $(2,3,8)\, (2,4,8)\, (2,6,8)\, (2,8,8)\, (3,3,4)\, (3,8,8)\,$ 
 $(4,4,4)\, (4,6,6)\, (4,8,8)$ \\
 $\Q(\sqrt 3)$ & $v_2$ &
 $(2,3,12)\, (2,6,12)\, (3,3,6)\, (3,4,12)\, (3,12,12)\, (6,6,6)\,$ \\
 $\Q(\sqrt 3)$ & $v_3$ & $(2,4,12)\, (2,12,12)\, (4,4,6)\, (6,12,12)\,$ \\
 $\Q(\sqrt 5)$ & $v_2$ & $(2,4,5)\, (2,4,10)\, (2,5,5)\, (2,10,10)\,
 (4,4,5)\, (5,10,10)$ \\
 $\Q(\sqrt 5)$ & $v_3$ & $(2,5,6)\, (3,5,5)$ \\
 $\Q(\sqrt 5)$ & $v_5$ & $(2,3,10)\, (2,5,10)\, (3,3,5)\, (5,5,5)$ \\
 $\Q(\sqrt 6)$ & $v_2$ & $(3,4,6)$ \\
 $\Q(\cos \pi/7)$ & $1$ &
 $(2,3,7)\, (2,3,14)\, (2,4,7)\, (2,7,7)\, (2,7,14)\, (3,3,7)\,$ 
 $(7,7,7)$ \\ 
 $\Q(\cos \pi/9)$ & $1$ & $(2,3,9)\, (2,3,18)\, (2,9,18)\, (3,3,9)\,
 (3,6,18)\, (9,9,9)$ \\
 $\Q(\cos \pi/9)$ & $v_2.v_3$ &
 $(2,4,18)\, (2,18,18)\, (4,4,9)\, (9,18,18)$ \\
 $\Q(\cos \pi/8)$ & $v_2$ &
 $(2,3,16)\, (2,8,16)\, (3,3,8)\, (4,16,16)\, (8,8,8)$ \\
 $\Q(\cos \pi/10)$ & $v_2$ & $(2,5,20)\, (5,5,10)$ \\
 $\Q(\cos \pi/12)$ & $v_2$ & 
 $(2,3,24)\, (2,12,24)\, (3,3,12)\, (3,8,24)\, (6,24,24)\,$ $(12,12,12)$ \\
 $\Q(\cos \pi/15)$ & $v_3$ & $(2,5,30)\, (5,5,15)$ \\
 $\Q(\cos \pi/15)$ & $v_5$ &
 $(2,3,30)\, (2,15,30)\, (3,3,15)\, (3,10,30)\, (15,15,15)$ \\
 $\Q(\sqrt 2, \sqrt 5)$ & $v_2$ & $(2,5,8)\, (4,5,5)$ \\
 $\Q(\cos \pi/11)$ & $1$ & $(2,3,11).$\\
\end{longtable}

\medskip\noindent ($v_p$ denotes the unique place of $E$ above
 $p$. Notice that all fields 
 $E$ appearing in this table are Galois --- even abelian --- over $\Q$).
\end{para}

\begin{para}\label{para::5.4.3} We do not know whether all arithmetic
 triangle groups are 
 congruent subgroups (actually, unit groups of quaternion algebras seem
 to be unsolved hard cases of the congruence subgroup problem, \cf 
 \cite{rapinchuk92:_congr}). Nevertheless, in each commensurability
 class of arithmetic 
 triangle groups, there is a (unique) triangle group of the form
 \begin{equation*} \Gamma^\ast= \{g\in B^\times\mid 
 g\mathcal{B}=\mathcal{B}g, \Nr(g)>0\}: \end{equation*} this is the one
 which appears first in the each of the above rows, \cf
 \cite[Table 3]{ccoatg}. The corresponding Schwarz orbifolds are Shimura
 orbifolds $\mathcal{X}^\ast$. More generally, any Schwarz orbifold whose
 orbifold fundamental group is a congruence group, is a Shimura orbifold.
\end{para}

\subsection{Arithmetic $p$-adic triangle groups. A $p$-adic analogue of
Takeuchi's list.}\label{sub::5.5} 

\begin{para}\label{para::5.5.1} Arithmetic $p$-adic triangle groups are
 hyperbolic $p$-adic triangle groups which are {\it commensurable with
 the $p$-unit groups of quaternion algebras $\ovl{B}$ over number
 fields}.\footnote{In fact, these were the first $p$-adic triangle groups
 to be constructed \cite{poaptg}, and the existence of non-arithmetic
 hyperbolic $p$-adic triangle groups remained unclear for a while, before
 Kato's construction.}
 
 \medskip Let now $V$ be the set of $p$-adic and infinite places of
 $E$. It follows from the discussion of \ref{para::5.4.1} that
 arithmetic subgroups of 
 $G_{(V)}= PGL_2(E_v)$ are subgroups which are commensurable with the
 $p$-unit groups $PSL_1(\ovl{\mathcal{B}}[\frac{1}{p}])$ of quaternion
 algebras over {\it totally real} number fields $E$, which are {\it
 totally definite at infinity and split at exactly one place $v$ above
 $p$ }.
\end{para}

\begin{para}\label{para::5.5.2} Let $\Delta_p(e_0,e_1,e_\infty)\subset
 PGL_2(\C_p)$ be an arithmetic $p$-adic triangle group. Thus, there is a
 totally real number field $E \subset \C_p$, and a totally definite
 quaternion algebra $\ovl{B}$ over $E$, split at one place $v$ above $p$,
 and a maximal order $\ovl{\mathcal{B}}$, such that
 $\Delta_p(e_0,e_1,e_\infty)\subset PGL_2(E_v)$ is commensurable, up to
 conjugation, to
 $PSL_1\bigl(\ovl{\mathcal{B}}[\frac{1}{p}]\bigr)$. Instead of the 
 latter, it is equivalent and more convenient to work with the image
 $\Gamma_p^\ast$ in $PGL_2(E_v)$ of the normalizer $\bigl\{g\in
 \ovl{B}^\times\mid 
 \;g\ovl{\mathcal{B}}[\frac{1}{p}]=\ovl{\mathcal{B}}[\frac{1}{p}]g\bigr\}$.

 Let us show that the corresponding complex triangle group
 $\Delta_p(e_0,e_1,e_\infty)$ $\subset$ $PGL_2(\C)$ is arithmetic, attached
 to the same totally real number field $E$ (embedded into $\C$ via some
 real place $\infty_0$), and to the quaternion algebra $B$ over $E$
 obtained from $\ovl{B}$ by interchanging the local invariants at
 $\infty_0$ and $v$ (in particular $B$ is definite at every place at
 infinity except $\infty_0$, and {\it ramified at $v$} --- which turns out
 to be the unique place $v_p$ of $E$ above $p$ according to Takeuchi's
 classification).
 
 In the case $\Delta_p(e_0,e_1,e_\infty)=\Gamma_p^\ast$, this follows
 from \v{C}erednik's uniformization explained in  II.\ref{sub:7.4},
 II.\ref{sub:7.5}: in this 
 case, the Schwarz orbifold is a Shimura orbifold $\mathcal{X}^\ast$ and
 $\Delta(e_0,e_1,e_\infty)=\Gamma^\ast$.

 In general, the Schwarz
 orbifold $\mathcal{S}(e_0,e_1,e_\infty)$ is close to be a Shimura
 orbifold.  If fact, according to \ref{thm::5.3.7}, the $p$-adic Schwarz
 orbifold 
 $\mathcal{S}(e_0,e_1,e_\infty)$ is a finite orbifold-covering of a
 Mumford-Schwarz orbifold $\mathcal{S}(e'_0,e'_1,e'_\infty)$ with
 $\Delta_p(e_0,e_1,e_\infty)=\Delta_p(e'_0,e'_1,e'_\infty)$. Let $\Gamma$
 be a torsion-free subgroup of finite index of
 $\Delta_p(e'_0,e'_1,e'_\infty)\cap \Gamma_p^\ast$. Let $\Omega \subset
 \P^1_{\C_p}$ be the set of ordinary points for $\Gamma$. Then the
 Mumford curve $S=\Omega/\Gamma$ is simultaneously a finite orbifold
 covering of $\mathcal{S}(e'_0,e'_1,e'_\infty)$ and of $\mathcal{X}^\ast$
 (via the embeddings $\Gamma \subset \Delta_p(e'_0,e'_1,e'_\infty)$ and
 $\Gamma \subset \Gamma_p^\ast$ respectively).

 All these finite orbifold-coverings are defined over number
 fields $\subset \C$. Looking at the corresponding complex orbifolds, it
 is now clear that $\pi_1^{\mathrm{orb}}(\mathcal{S}(e_0,e_1,e_\infty))$
 $=\Delta(e_0,e_1,e_\infty)$ is commensurable with
 $\pi_1^{\mathrm{orb}}(\mathcal{X}^\ast)= \Gamma ^\ast$, whence the claim.
\end{para}

\begin{para}\label{para::5.5.3} Let us show, conversely, that if
 $\Delta(e_0,e_1,e_\infty)$ appears in Takeuchi's list, and if $B/E$ is
ramified at the place $v\mid p$, then $\Delta_p(e_0,e_1,e_\infty)$
exists and is an arithmetic $p$-adic triangle group; more precisely, it
is commensurable with the $p$-unit group of the quaternion algebras
$\ovl{B}/E$ obtained from $B$ by interchanging the local invariants at
$\infty_0$ (the indefinite place at infinity) and $v_p$ (in particular
$\ovl{B}$ is totally definite, and split at $v_p$).
 
 In the case $\Delta(e_0,e_1,e_\infty)=\Gamma^\ast$, \ie if
 $(e_0,e_1,e_\infty)$ is the first triple occurring on a row of 
 Takeuchi's table, this follows from \ref{thm::4.7.4}.     

 In general, $\Delta(e_0,e_1,e_\infty)$ is commensurable to an arithmetic
 triangle group $\Delta(e'_0,e'_1,e'_\infty)$ of the form $\Gamma^\ast$,
 and it is not difficult to deduce that $\Delta_p(e_0,e_1,e_\infty)$
 exists and is commensurable to
 $\Delta_p(e'_0,e'_1,e'_\infty)=\Gamma_p^\ast$. The point is that the
 corresponding respective projective connections $HGPC(a,b,c)$ and
 $HGPC(a',b',c')$ become isomorphic after suitable finite pull-backs.
\end{para}

\begin{para}\label{para::5.5.4} Putting \ref{para::5.5.2} and
 \ref{para::5.5.3} together, we conclude 
 that the arithmetic $p$-adic triangle groups
 $\Delta_p(e_0,e_1,e_\infty)$ are exactly those corresponding to the
 triples $(e_0,e_1,e_\infty)$ in Takeuchi's list when $v_p$ divides
 $\mathrm{disc}(B)$.  From that list, one gets:
\end{para}

\begin{thm}\label{thm::5.5.5} \cite{poaptg}. Up to permutation of the
 indices, there are $45$ arithmetic $2$-adic triangle groups, 
 $16$ arithmetic
 $3$-adic triangle groups, $9$ arithmetic $5$-adic triangle groups (and
 no arithmetic $p$-adic triangle group for $p>5$).

 More precisely, they are given by the following tables:

\begin{center}
 \renewcommand{\arraystretch}{1.5}
 \begin{longtable}{lcp{0.66\textwidth}}
 \multicolumn{3}{l}{\it Arithmetic $2$-adic triangle groups:}\\ \hline
 $E$ &$\mathrm{disc}(\ovl{B})$& $(e_0,e_1,e_\infty)
 \text{ \rm up to permutation}$ \\ \hline
 $\Q$ & $3$ & $(2,4,6)\,(2,6,6)\,(3,4,4)\,(3,6,6)$ \\
 $\Q(\sqrt 2)$ & $1$ &
 $(2,3,8)\, (2,4,8)\, (2,6,8)\, (2,8,8)\, (3,3,4)$\, $(3,8,8)\, 
 (4,4,4)\, (4,6,6)\, (4,8,8)$ \\
 $\Q(\sqrt 3)$ & $1$ &
 $(2,3,12)\, (2,6,12)\, (3,3,6)\, (3,4,12)\, (3,12,12)$\, $(6,6,6)$ \\
 $\Q(\sqrt 5)$ & $1$ &
 $(2,4,5)\, (2,4,10)\, (2,5,5)\, (2,10,10)\, (4,4,5)$\, $(5,10,10)$ \\
 $\Q(\sqrt 6)$ & $1$ & $(3,4,6)$ \\
 $\Q(\cos \pi/9)$ & $v_3$ & $(2,4,18)\, (2,18,18)\, (4,4,9)\, (9,18,18)$ \\
 $\Q(\cos \pi/8)$ & $1$ & $(2,3,16)\, (2,8,16)\, (3,3,8)\, (4,16,16)\,
 (8,8,8)$ \\ 
 $\Q(\cos \pi/10)$ & $1$ & $(2,5,20)\, (5,5,10)$  \\
 $\Q(\cos \pi/12)$ & $1$ &
 $(2,3,24)\, (2,12,24)\, (3,3,12)\, (3,8,24)\, (6,24,24)$\, $(12,12,12)$ \\
 $\Q(\sqrt 2, \sqrt 5)$ & $1$ & $(2,5,8)\, (4,5,5)$  \\ &&\\
 \multicolumn{3}{l}{\it Arithmetic $3$-adic triangle groups:} \\ \hline
 $\Q$ & $2$ & $(2,4,6)\, (2,6,6)\, (3,4,4)\, (3,6,6)$  \\
 $\Q(\sqrt 3)$ & $1$ & $(2,4,12)\, (2,12,12)\, (4,4,6)\, (6,12,12)$ \\
 $\Q(\sqrt 5)$ & $1$ & $(2,5,6)\, (3,5,5)$\\
 $\Q(\cos \pi/9)$ & $v_2$ & $(2,4,18)\, (2,18,18)\, (4,4,9)\, (9,18,18)$ \\
 $\Q(\cos \pi/15)$ & $1$ & $(2,5,30)\, (5,5,15)$ \\ && \\
 \multicolumn{3}{l}{\it Arithmetic $5$-adic triangle groups:}\\ \hline
 $\Q(\sqrt 5)$ & $1$ & $(2,3,10)\, (2,5,10)\, (3,3,5)\, (5,5,5)$ \\
 $\Q(\cos \pi/15)$ & $1$ &
 $(2,3,30)\, (2,15,30)\, (3,3,15)\, (3,10,30)\, (15,15,15)$ \\
 \end{longtable} 
\end{center}
\end{thm}

\begin{rems}
 (a) In these lists, one notice that $p$ divides
 at least one of the numbers $e_0,e_1,e_\infty$. This means that the
 parameters $a,b,c$ of the corresponding hypergeometric differential
 equation are not $p$-adic integers. In the literature on $p$-adic
 hypergeometric functions, on the contrary, $a,b,c$ are always assumed to
 be $p$-adic integers, for obvious convergence reasons.

 It is well-known that a fuchsian $p$-adic differential equations with
 Frobenius structure has exponents in $\Q\cap \Z_p$. Since this is not
 the case here, there is no Frobenius structure in the standard
 sense. However, there is some hidden Frobenius structure due to the
 ``geometric origin'' of the differential equation. We leave to the
 reader the intriguing question of finding out the right formulation of
 this hidden Frobenius structure (hint: the Frobenius structure is
 explicit at the level of the period mapping, \cf \cite{tpuosc}).

 \medskip
 (b) Takeuchi has also given a table with the various
 inclusions of triangle groups in each commensurability class. When
 $\Delta(e_0,e_1,e_\infty)$ is maximal (and $p$ ramified in the
 corresponding quaternion algebra), the corresponding $p$-adic triangle
 group $\Delta_p(e_0,e_1,e_\infty)$ is ``of Mumford type'', according to
 \ref{cor::5.3.7}.~(b).

 Arithmetic $p$-adic triangle groups not ``of Mumford type'' do exist,
 and are related to Shimura curves with {\it level divisible by $p$}.
 One can construct such $p$-adic triangle groups using Drinfeld's finite
 \'etale coverings $\Upsilon^n$ of the Drinfeld ``half space'' $\Omega$,
 which uniformize $p$-adically Shimura curves with level $p^n.M\,$ ($p$
 prime to $M$), \cf II.\ref{para:7.4.7}.  One has
 $\Gamma_p^+(p^nM)=\Gamma_p^+( M)$, 
 and the orbifold-quotient of a suitable intermediate covering $\Sigma$
 between $\Omega$ and $\Upsilon^n$ by $\Gamma_p^+( M)$ may be a Schwarz
 orbifold but not a Mumford-Schwarz orbifold. It seems that
 $\Delta_2(4,4,4)$ is of this type, for $p=2, \,n=M=1$ (we have not
 checked the details); the corresponding $\Gamma_2^+(1)$ is
 $\Delta_2(3,3,4)$.

 \medskip
 (c) Let $\ovl{S}$ be a smooth projective curve defined
 over $\ovl{\Q}$, and $\zeta_i$ be $\ovl{\Q}$-points on $\ovl{S}$
 endowed with indices $e_i$. We then say that the orbifold
 ${\mathcal{S}}=(\ovl{S}, (\zeta_i,e_i))$ is defined over $\ovl{\Q}$. In
 general, there is no reason for the uniformizing differential equation
 to be defined over $\ovl{\Q}$.

 However, this occurs if either ${\mathcal{S}}$ is a Schwarz orbifold
 (the uniformizing differential equation being essentially of
 hypergeometric type with rational parameters), or a Shimura orbifold
 (the uniformizing differential equation being essentially of Gauss-Manin
 type, \cf \ref{thm::4.7.4}). On the other hand, it is not difficult to
 see that the 
 property for the uniformizing differential equation to be defined over
 $\ovl{\Q}$ depends only on the commensurability class of the orbifold
 fundamental group $\pi_1^{\mathrm{orb}}(\mathcal{S})$.

 A transcendence conjecture of Chudnovsky, modified by Krammer 
 \cite[12]{aeoaafg}, predicts that the Schwarz and the Shimura 
 cases are the only exceptions, up to commensurability. Namely: 
 
 \begin{coj}\label{coj::5.5.6} (Chudnovsky-Krammer). Let ${\mathcal{S} }$
  be a one-dimensional hyperbolic orbifold defined over
  $\ovl{\Q}$. Assume that the uniformizing differential equation is
  defined over $\ovl{\Q}$. Then $\pi_1^{\mathrm{orb}}(\mathcal{S})$ is either
  commensurable to a triangle group, or is arithmetic.
 \end{coj}

 We {\it conjecture} the (straightforward) $p$-adic analogue of
 \ref{coj::5.5.6}.
 
 \medskip
 (d) By the theory of $p$-adic Betti lattices, there is
 a canonical $\ovl{\Q}$-structure in the $\C_p$-space of solutions at
 any point of $\P^1(\ovl{\Q})$ of the hypergeometric differential
 equation $HDGE(a,b,c)$ attached to an arithmetic $p$-adic triangle group
 (\cf \ref{thm::4.7.4}. If the base point is chosen to be one of the
 singularities $0,1,\infty$, one applies \ref{thm::4.7.4}
 to a finite orbifold-covering of the Schwarz orbifold which is a
 Shimura curve attached to a torsion-free congruence subgroup). 
 
 This raises the following (open) questions:

 {\it Is there a $\ovl{\Q}$-structure in the $\C_p$-space
 of solutions of the hypergeometric differential equation $HDGE(a,b,c)$
 attached to a non-arithmetic $p$-adic triangle group? Are
 non-arithmetic $p$-adic triangle groups related to some $p$-adic period
 mappings and suitable higher-dimensional Shimura orbifolds?}
\end{rems}

\subsection{On special values of hypergeometric functions attached to 
arithmetic $p$-adic triangle groups.}\label{sub::5.6}  

\begin{para}\label{para::5.6.1} Let $a,b,c$ be rational numbers such that
 $c>a+b$. Then the hypergeometric series ${}_2F_1(a,b;c;z)$ 
 converges at $1$ and takes the value (Gauss' formula) 
 \begin{equation*}
  {}_2F_1(a,b;c;1)= \frac{\Gamma(c)\Gamma(c-a-b)}{\Gamma(c-a)\Gamma(c-b)}\;.
 \end{equation*}   

\noindent If one considers ${}_2F_1(a,b;c;z)$ as a multivalued analytic
function on $\P^1\setminus\{0,1,\infty\}$, the 
values at $1$ form a finitely generated subgroup of $\C$; they generate
a $\ovl{\Q}$-subspace of $\C$ of dimension $\leq 2$,
namely:  
\begin{equation*}
 \frac{\Gamma(c)\Gamma(c-a-b)}{\Gamma(c-a)\Gamma(c-b)}\ovl{\Q} \, +\,
  \frac{\Gamma(c)\Gamma(a+b-c)}{\Gamma(a)\Gamma(b)}\ovl{\Q}. 
\end{equation*}
Similarly, its
values at $\infty$ generate the $\ovl{\Q}$-subspace of $\C$  
\begin{equation*}
 \frac{\Gamma(c)\Gamma(b-a)}{\Gamma(c-a)\Gamma(b)}\ovl{\Q} \, +\,
  \frac{\Gamma(c)\Gamma(a-b)}{\Gamma(a)\Gamma(c-b)}\ovl{\Q}. 
\end{equation*}   
\end{para}

\begin{para}\label{para::5.6.2} If $HGDE(a,b,c)$ is the differential equation
 attached as above to an {\it arithmetic} triangle group
$\Delta(e_0,e_1,e_\infty)$, these formulas are special cases of similar
formulas expressing the special values of ${}_2F_1(a,b;c;z)$ at $CM $
points, in terms of Gamma values, up to factors in $\ovl{\Q}$.

The point is that there is a $\ovl{\Q}$-structure in
the complex space of solutions which is stable under monodromy: the
 $\mathrm{H}_{1,\mathrm{B}}({\;\,},\ovl{\Q})$ of the fake elliptic curves $A_z$
 parametrized by $z\in\P^1$ ($\P^1$ being viewed as a Shimura curve) if
the quaternion algebra attached to $\Delta(e_0,e_1,e_\infty)$ is over
$\Q$, \cf II.\ref{para:7.6.3} for the general case. For any point $z\in
\P^1(\ovl{\Q})$, let $V_z$ be the $\ovl{\Q}$-structure in the complex
space of solutions at $z$ consisting of solutions whose Taylor
coefficients are in $\ovl{\Q}$. Then $V_z$ is related to
$\mathrm{H}_{1,\mathrm{B}}(A_{z,\C},\ovl{\Q})$ by the period matrix of
 $A_z$. At any CM 
point $z$ (according to an old observation of Shimura, the branch points
of a Shimura orbifold --- here, $0,1,\infty$ --- are actually CM points),
the Chowla-Selberg formula (I.\ref{thm-ogus-formula}) thus relates the two
$\ovl{\Q}$-structures in terms of $\Gamma$-values. In particular, $V_0$
is generated by ${}_2F_1(a,b,c;z)$ and
$z^{1-c}{}_2F_1(a+1-c,b+1-c,2-c;z)$, and the branch point $0$
corresponds to a CM point (for elliptic curves rather than fake elliptic
curves, this was worked out in I.\ref{coj-no-complex-multi}).
\end{para}

\begin{para}\label{para::5.6.3} There is an entirely analogous phenomenon for
 {\it arithmetic $p$-adic} triangle 
group $\Delta_p(e_0,e_1,e_\infty)$: thanks to II.\ref{thm:7.6.4}, there is a
 $\ovl{\Q}$-structure in 
the $\C_p$-space of solutions which is stable under $p$-adic monodromy:
 the ``$p$-adic Betti lattice''
 $\mathrm{H}_{1,\mathrm{B}}({\;\,},\ovl{\Q})$ of the fake 
 elliptic 
 curves $A_z$ if the quaternion algebra attached to
 $\Delta_p(e_0,e_1,e_\infty)$ is over $\Q$, \cf II.\ref{para:7.6.3} for the
 general case. ``Periods'' should be understood here in the sense of
 I.\ref{sub-p-adic-betti-lattice-supersing}; $p$-adic CM periods
 were computed in terms of $\Gamma_p$-values 
 in I.\ref{thm-ogus-formula}, I.\ref{para-ramif} (for $p\neq 2$), 
 I.\ref{coj-complex-mult}.

 Therefore, if $HGDE(a,b,c)$ is the differential equation attached as
 above to an {\it arithmetic $p$-adic} triangle group
 $\Delta_p(e_0,e_1,e_\infty)$, there is, at least for odd $p$, a
 expression in terms of $\Gamma_p$-values of the $\ovl{\Q}$-space of
 dimension $\leq 2$ generated by the values taken at any given CM point, \eg
 the point $1$, by ${}_2F_1(a,b,c;z)$ viewed as a multivalued $p$-adic
 analytic function on a suitable finite \'etale covering of
 $\P^1\setminus\{0,1,\infty\}$.

 We shall see examples in the next section.
\end{para}

\newpage

\section{The tale of $F(\frac{1}{24},\frac{7}{24},\frac{5}{6};z)$
(Escher's triangle group and its diadic and triadic twins).}\label{sec::6}
 
\markboth{\thechapter. $p$-ADIC ORBIFOLDS AND MONODROMY.}%
{\thesection. THE TALE OF $F(\frac{1}{24},\frac{7}{24},\frac{5}{6};z)$.}

\begin{abst}
 We study the rich geometry related to the simplest cocompact arithmetic
 triangle group: associated Shimura curves, period mapping and monodromy
 --- in the archimedean, diadic and triadic cases. At other primes, the
 story would be similar to that of $F(\frac{1}{2},\frac{1}{2},{1};z)$;
 we do not repeat it.
\end{abst}

\subsection{Escher's triangle group $\Delta(2,4,6)$.}\label{sub::6.1}

\begin{para}\label{para::6.1.1} This is the symmetry group of the
 black-and-white tesselation of the unit disk by hyperbolic 
 triangles with angles $(\frac{\pi}{2}, \frac{\pi}{4}, \frac{\pi}{6})$

\begin{figure}[h]
 \begin{picture}(300,170)(0,0)
  \put(0,0){\includegraphics[scale=0.43,clip]{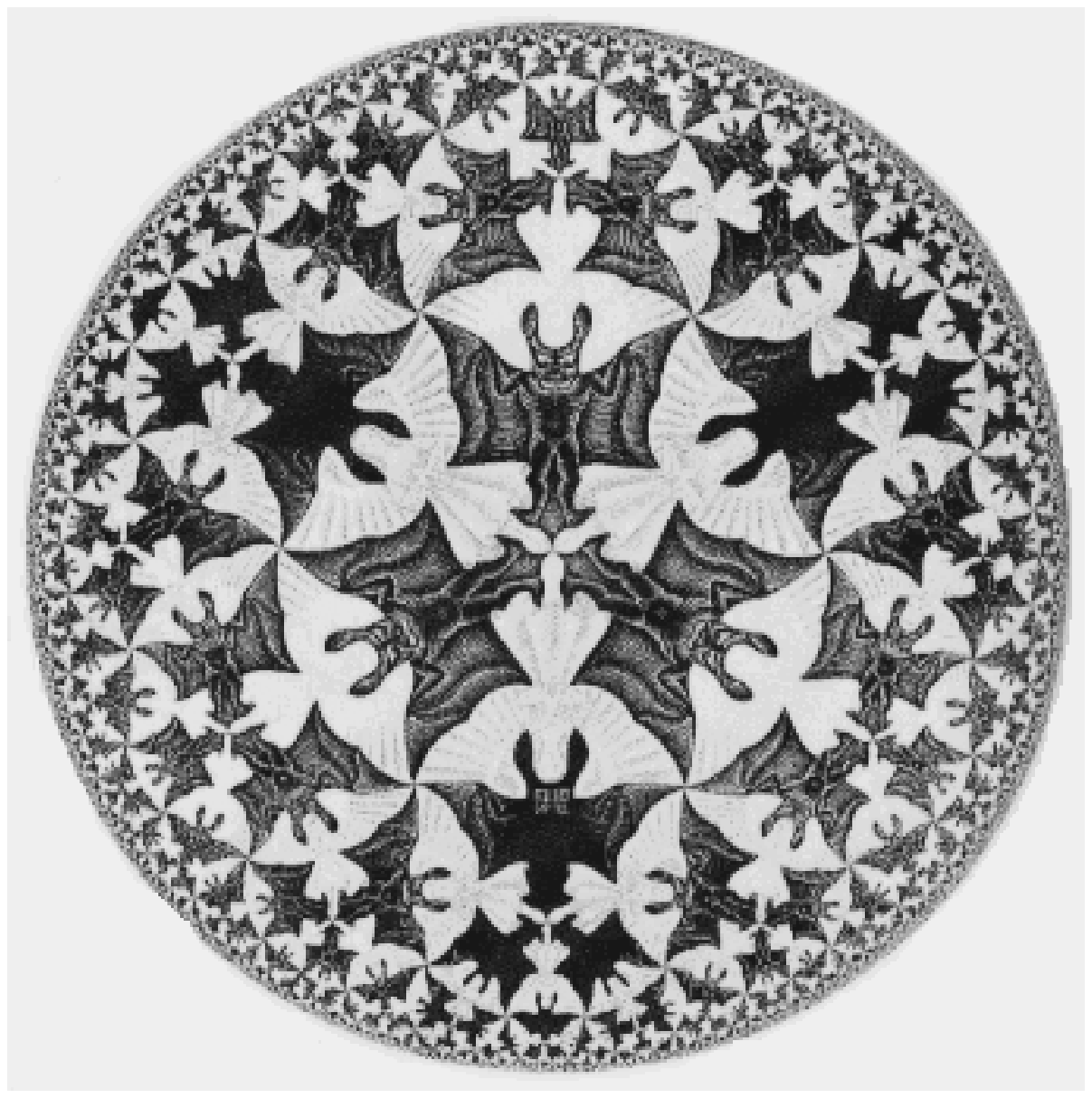}}
  \put(170,0){\includegraphics[scale=0.6,clip]{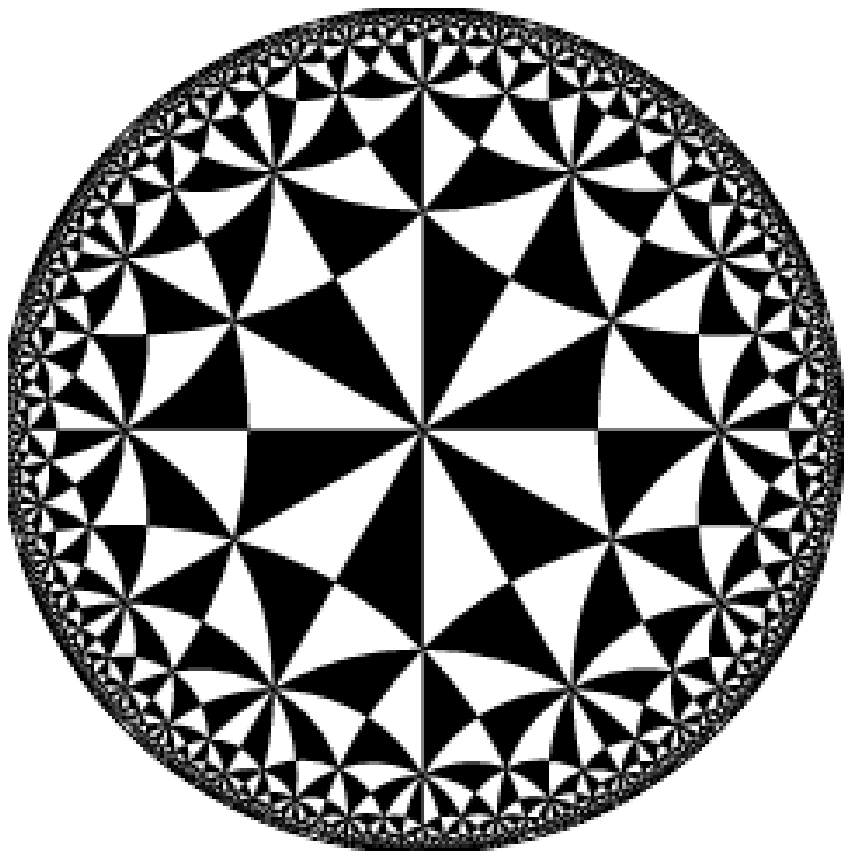}}
 \end{picture}
 \caption{}
\end{figure}

 M.C. Escher has engraved several celebrated variants of this Schwarz
 tesselation (``cirkel-limiet'', angels and 
 demons), {\it cf}  \cite{escher71:_m} \footnote{he has also used another
 Schwartz tesselation in a picture --- ``cirkel-limiet III'' ---
 brilliantly analyzed by H.S.M. Coxeter \cite{coxeter79:_escher_iii}}.

\medskip By gluing together two basic triangles, we get hyperbolic
 triangles with angles $(\frac{\pi}{3}, \frac{\pi}{4}, 
 \frac{\pi}{4})$ or $(\frac{\pi}{2}, \frac{\pi}{6}, \frac{\pi}{6})$

\begin{figure}[h]
 \begin{picture}(300,100)(0,-50)
  \put(0,-32){\includegraphics[scale=0.8,clip]{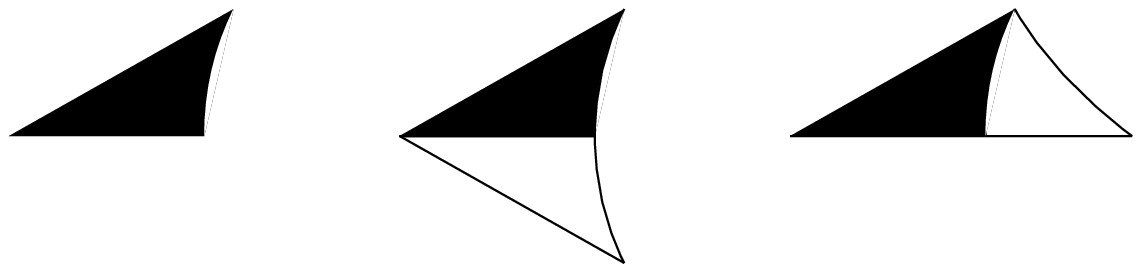}}
  {\small
  \put(0,-10){\makebox(0,0)[c]{$O$}}
  \put(47,-10){\makebox(0,0)[c]{$Q$}}
  \put(65,34){\makebox(0,0)[c]{$P$}}
  \put(-6,0){\makebox(0,0)[b]{$6$}}
  \put(55,0){\makebox(0,0)[b]{$2$}}
  \put(54,37){\makebox(0,0)[c]{$4$}}

  \put(85,0){\makebox(0,0)[b]{$3$}}
  \put(145,37){\makebox(0,0)[c]{$4$}}
  \put(152,-30){\makebox(0,0)[c]{$4$}}

  \put(237,37){\makebox(0,0)[c]{$2$}}
  \put(174,0){\makebox(0,0)[b]{$6$}}
  \put(269,0){\makebox(0,0)[b]{$6$}}
  }
 \end{picture}
 \caption{}
\end{figure}

\noindent This exhibits $\Delta(3,4,4)$ and $\Delta(2,6,6)$ as subgroups
 of index $2$ of $\Delta(2,4,6)$. 

\noindent Furthermore, by gluing two triangles with angles
 $(\frac{\pi}{2}, \frac{\pi}{6}, \frac{\pi}{6})$, one can get a 
triangle with angles $(\frac{\pi}{3}, \frac{\pi}{6}, \frac{\pi}{6})$,
 which exhibits $\Delta(3,6,6)$ as a subgroup 
 of index $2$ of $\Delta(2,6,6)$. 

\noindent Similarly, the intersection of $\Delta(2,6,6)$ and
 $\Delta(3,4,4)$ in $\Delta(2,4,6)$ is a quadrangle 
 group $\diamondsuit (2,2,3,3)$.
\end{para}

\begin{para}\label{para::6.1.2} The vertices of one of the basic triangles are

\begin{equation*}
 O=(0,0),\;Q=(\sqrt 2 -1, 0),\;P=\Bigl(\sqrt 3 \cos
 \frac{5\pi}{12}, \cos \frac{5\pi}{12}\Bigr).
\end{equation*}

\medskip\noindent
The matrix of the composition of the reflection along side $QP$, \resp
$PO$, followed by the reflection 
along side $OQ$ is $A=\sqrt {-1}\begin{pmatrix}\sqrt 2 & -1 \\1 &-\sqrt 2
\end{pmatrix}$, \resp  $B=\dfrac{1}{2}\begin{pmatrix}\sqrt 3 +\sqrt
{-1} & 0 \\0 & \sqrt 3 -\sqrt {-1}
\end{pmatrix}$. 

\noindent These matrices satisfy $A^2=B^6=(AB)^4=-1$ and, taken up to
 sign, generate the triangle group $\Delta(2,4,6)$, 
 (\cf \cite[II5]{ntatg}.)
\end{para}

\begin{para}\label{para::6.1.3} The group $\Delta(2,4,6)$ enjoys the following
 ``universal property'', discovered by Mennike in 1968: 

\begin{quote}
 {\it for all but a finite number of indefinite ternary quadratic forms
 $q$ with integral coefficients, the 
 (integral) associated unitary group of $q$ is contained in $\Delta(2,4,6)$.}
\end{quote}

 \noindent It coincides with $\Delta(2,4,6)$ in the case of
 $q=3x^2-y^2-z^2$, (\cf \cite[III3]{ntatg}).  These unitary groups
 have been thoroughly investigated by Fricke and Klein.
\end{para}

\subsection{Some Shimura curves attached to $D_{2.3}$.}\label{sub::6.2}

\begin{para}\label{para::6.2.1} The $\Q$-subalgebra of $M_2(\Q (\sqrt
 {-1}, \sqrt 2, \sqrt 3))$ generated by the matrices $A$ and $B$ 
 is the quaternion algebra $B=B_{2.3}$ over $\Q$ with discriminant
 $6$. This quaternion algebra is unique up to 
 isomorphism and admits standard generators $i,j$ and relations 
 $i^2=-1,\;j^2=3,\;ij=-ji$.  

 \noindent It follows that $\Delta(2,4,6)$ is a subgroup of
 $B_{2.3}^\times/\Q^\times$. According to \cite{ccoatg}, 
 $\Delta(2,4,6)$, $\Delta(3,4,4)$, $\Delta(2,6,6)$ and $\Delta(3,6,6)$
 are the only cocompact triangle 
 groups which give rise to a quaternion algebra over $\Q$. 

 \noindent As any indefinite quaternion algebra over $\Q$, $B_{2.3}$ has
 a unique maximal order $\mathcal{B}=\mathcal{B}_{2.3}$ 
 modulo conjugation, \cf \cite[p.99]{adadq}. It can be described
 as $\Z  + \Z i + \Z j + \Z \rho$, where 
 $\rho =\frac{1+i+j+ij}{2}$. According to \cite{ccoatg}, it can also be
 described as {\it the subring of 
 $B_{2.3}$ generated by the elements of $\diamondsuit (2,2,3,3)$} (the
 subgroup of $\Delta(2,4,6)$ of index $4$ 
 considered above). 

 Following the notation of  II.\ref{sub:7.4}, we introduce the following
 arithmetic subgroups of $B^\times/\Q^\times$:

\begin{itemize}
 \item the image $\Gamma^+(N)$ of the group $\{g\in
       (1+N\mathcal{B})^\times \mid \Nr(g)=1\}$ in 
       $PSU(1,1)(\R)$ $=\Aut(\matheur{D}(0,1^-))$, where $\Nr$ stands for the
       reduced norm ($\Gamma^+(1)=PSL_1(\mathcal{B})$ is abbreviated
       into $\Gamma^+$),
 \item the image $\Gamma^\ast$ in $PSU(1,1)(\R)$ of the group 
       $N(\mathcal{B})^+=\{g\in B^\times\mid
       g\mathcal{B}=\mathcal{B}g \mid \Nr(g)>0\}$, \ie 
       $N(\mathcal{B})^+/\Q^\times$. 
\end{itemize}
       
\medskip According to \cite{ccoatg}, one has: 
 
\begin{itemize}
 \item $\Delta(2,4,6)=\Gamma^\ast$, 
 \item $\diamondsuit (2,2,3,3)=\Gamma^+$.
\end{itemize}

These groups also appear in \cite[p.122]{adadq}, together
with the orbifolds

\noindent ${\mathcal{X}^+}_\C^{\mathrm{an}}= \matheur{D}(0,1^-)/\diamondsuit 
 (2,2,3,3)$ and ${\mathcal{X}^\ast}_\C^{\mathrm{an}}=\matheur{D}(0,1^-)/\Delta(2,4,6)$.
 The orbifold Euler characteristics are found to be $-\frac{1}{3}$ and 
 $-\frac{1}{12}$ respectively.  
\end{para}

\begin{para}\label{para::6.2.2} Any element of $N(\mathcal{B})^+$ can be
 written uniquely as a product 
\begin{equation*}
 q.g.(1+i)^{\epsilon_2}(3i+j+ij)^{\epsilon_3},
\end{equation*}
where $q\in \Q, \; g\in SL_1(\mathcal{B}),\;\epsilon_m=0 $
or $1$, \cf \cite{michon84:_courb_shimur}. The quotient group 

 \begin{align*}
  &\Delta(2,4,6)/\diamondsuit(2,2,3,3)\\
  = &\{1,w_2=[1+i],\;w_3=[3i+j+ij],\; w_6=w_2w_3=[-3+3i+2ij]\}\\
  \cong &(\Z/2\Z)^2
 \end{align*}
 is called the Atkin-Lehner group. The fixpoints of the involutions $w_m$
 have been the object of thorough 
 investigations, \cf {\it loc.~cit.} and \cite{shimura67:_const} (\cf also
 II.\ref{para:7.4.4}). In our case, the results are:

 $\bullet$ the two fixpoints of $w_2$ on $\matheur{D}(0,1^-)/\diamondsuit
 (2,2,3,3)$ have the same image in $\matheur{D}(0,1^-)/\Delta(2,4,6)$, which is the
 image of the fixpoint in $\matheur{D}(0,1^-)$ of any maximal 
 embedding of $\Z[\sqrt {-1}]$ in $\mathcal{B}$, 

$\bullet$ the two fixpoints of $w_3$ on $\matheur{D}(0,1^-)/\diamondsuit
(2,2,3,3)$ have the same image in $\matheur{D}(0,1^-)/\Delta(2,4,6)$, which is the
image of the fixpoint in $\matheur{D}(0,1^-)$ of any maximal 
embedding of $\Z[^3\sqrt{1}]$ in $\mathcal{B}$, 

$\bullet$ the two fixpoints of $w_2$ on $\matheur{D}(0,1^-)/\diamondsuit
(2,2,3,3)$ have the same image in $\matheur{D}(0,1^-)/\Delta(2,4,6)$, which is the
image of the fixpoint in $\matheur{D}(0,1^-)$ of any maximal 
embedding of $\Z[\sqrt {-6}]$ in $\mathcal{B}$, 

 From easy considerations of ramification, one draws that 

\begin{itemize}
 \item $\an{\Gamma^+,(1+i)}=\Delta (3,4,4)$
 \item $\an{\Gamma^+,(3i+j+ij)}=\Delta (2,6,6)$
 \item $\an{\Gamma^+,(1+i)}$ is the third intermediate group, which is
       a quadrangle group $\diamondsuit 
       (2,2,2,3)$ (not to be confused with $\diamondsuit (2,2,3,3)$).
\end{itemize}
\end{para}

\begin{lem}\label{lem::6.2.3} The congruence group $\Gamma^+(2)$ is
 torsion-free. The quotient $\Gamma^+/\Gamma^+(2)$ is the 
 tetrahedral group $\Delta(2,3,3)\cong \mathfrak{A}_4$ (of order $12$).
\end{lem}

\begin{proof}
 Recall that $\mathcal{B}=\Z + \Z i + \Z j + \Z \rho$. Any element
 $a+bi+cj+dij$ of finite order in the congruence subgroup
 $SL_1(\mathcal{B})(2)$ of level two must have reduced trace $2a=\pm 2$
 and reduced norm $a^2+b^2-3c^2-3d^2=1$; it is easy to see that the only
 such elements are $\pm 1$, which implies the first assertion.

 The quotient group $\Gamma^+/\Gamma^+(2)$ is represented by $\{1,i,j,ij,
 \frac{1\pm i\pm j \pm ij }{2}\}$ modulo $2$. It admits $\{1,\ovl{i},
 \ovl{j},\ovl{i}\ovl{j}\}$ as a normal subgroup, and
 $\{1,\ovl{\rho},\ovl{\rho}^2\}$ as a non-normal subgroup. its center is
 trivial.  It is thus isomorphic to the tetrahedral group, the group of
 oriented symmetries of the tetrahedron.
\begin{figure}[h]
 \begin{picture}(100,100)(0,0)
  \put(0,0){\includegraphics[scale=0.3,clip]{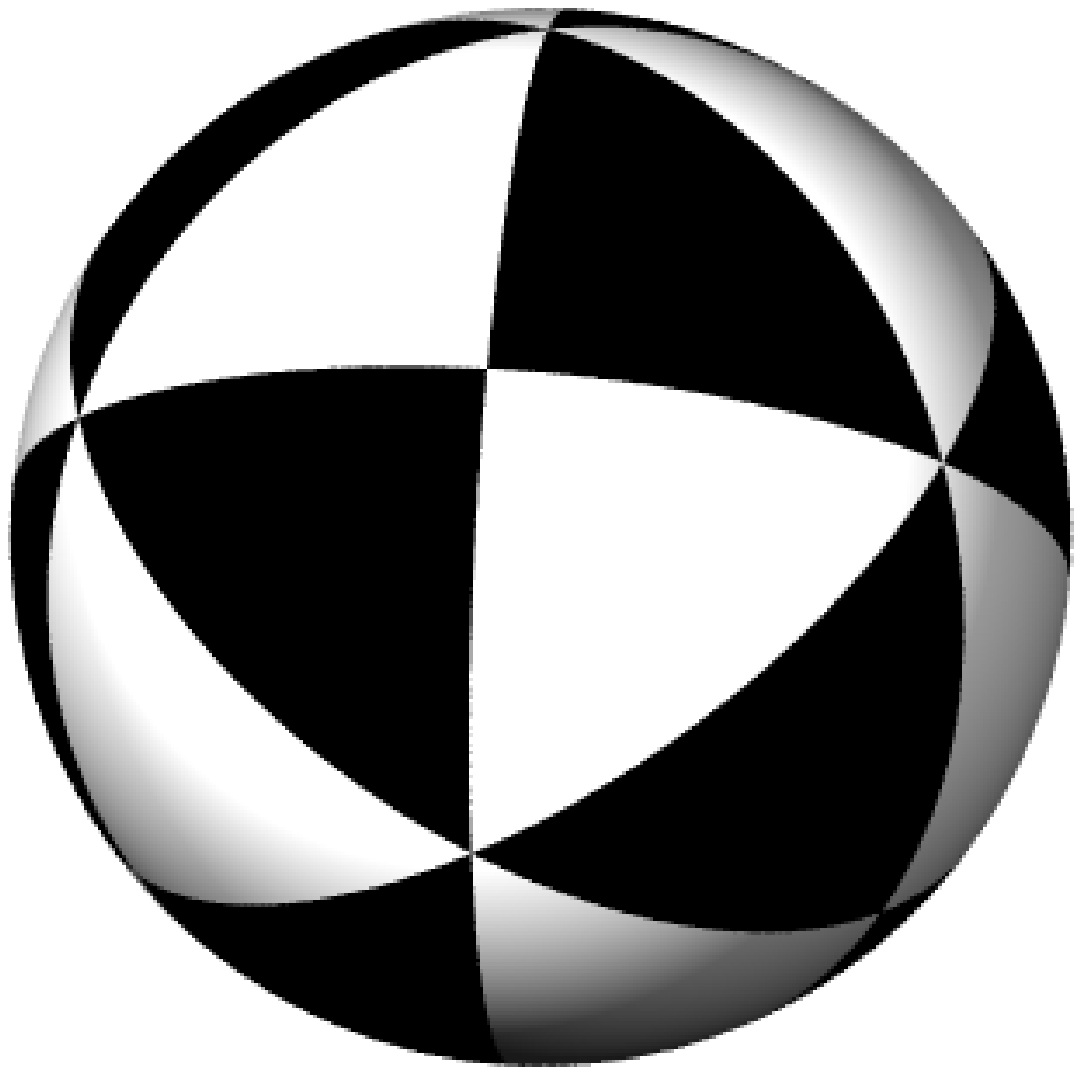}}
 \end{picture}
 \caption{}
\end{figure}
\end{proof}

\begin{cor}\label{cor::6.2.4}
 The associated Shimura curve $\mathcal{X}^+(2)$ has genus $g=3$. 
\end{cor}
\begin{proof}
 By definition,
 $\mathcal{X}^+(2)_\C^{\mathrm{an}}=\matheur{D}(0,1^-)/\Gamma^+(2)$ (here we
 prefer to identify the universal 
 covering with the unit disk rather than the Poincar\'e upper half plane
 as in  II.\ref{sub:7.4}). The assertion follows from 
 the computation
 $2-2g=\chi(\matheur{D}(0,1^-)/\Gamma^+(2))
 =12.\chi_{orb}(\matheur{D}(0,1^-)/\Gamma^+)=-4$.
\end{proof}

\begin{lem}\label{lem::6.2.5} The congruence group $\Gamma^+(2)$ is a
normal subgroup of $\Delta(2,4,6)$ $=$ $\Gamma^\ast$.  The quotient 
$\Gamma^\ast/\Gamma^+(2)$ is the extended
octahedral group $\Delta^\ast(2,3,4)$ $\cong$ $\mathfrak{S}_4\times \Z/2\Z$ 
(of order $48$). The generator of the factor $\Z/2\Z$ is 
given by the image $\tilde w_3$ of $(3i+j+ij)$.
\end{lem}

\par We refer to \cite[N. 58]{ford29:_autom} and
 \cite[p.27]{coxeter74:_regul} for a discussion 
 of $\Delta^\ast(2,3,4)$, the group of all symmetries of the octahedron.

\begin{figure}[h]
 \begin{picture}(100,100)(0,0)
  \put(0,0){\includegraphics[scale=0.3,clip]{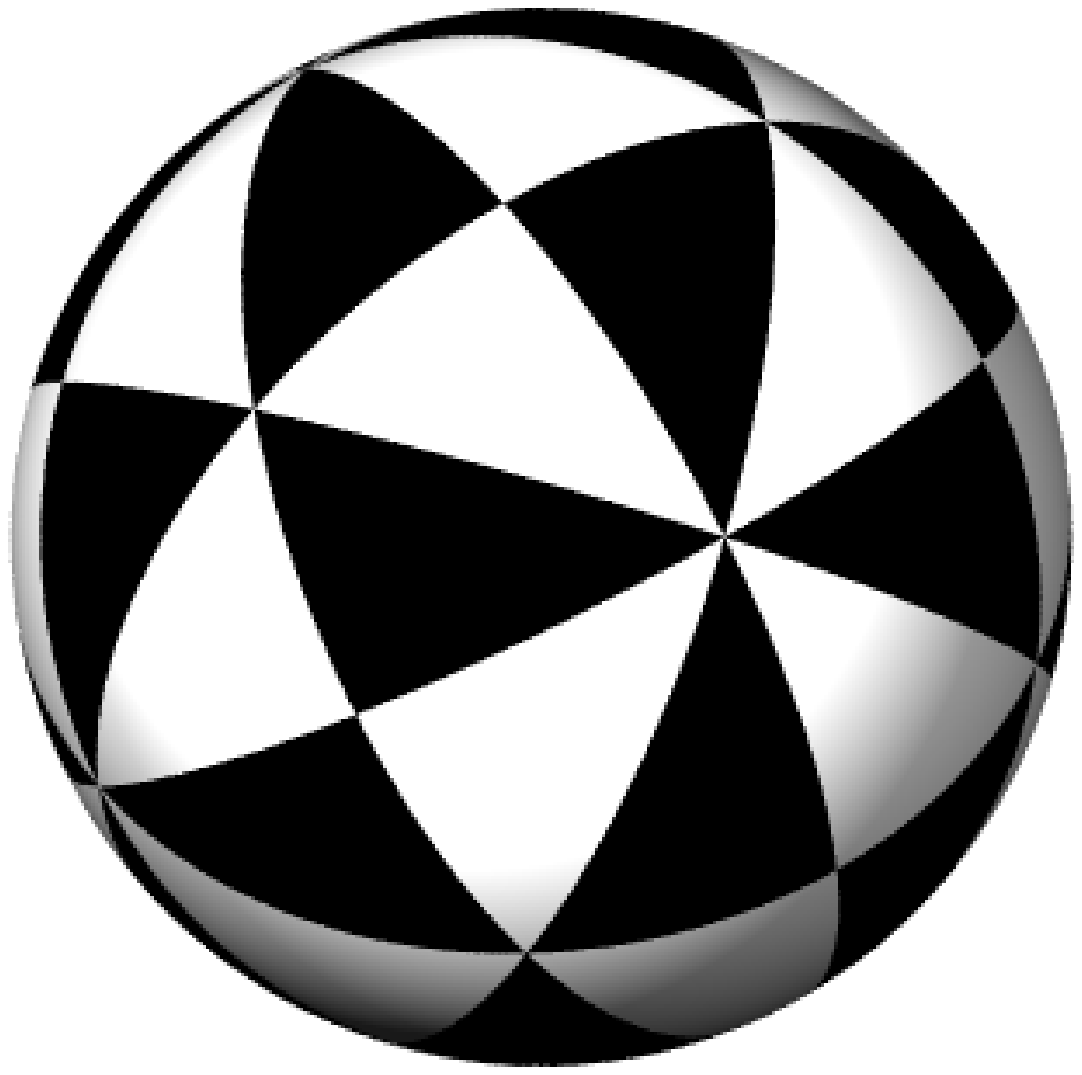}}
 \end{picture}
 \caption{}
\end{figure}

\begin{proof}
 Let us first show that $\Gamma^+(2)$ is a normal
 subgroup of $\Delta(2,6,6)=\an{\Gamma^+,(3i+j+ij)}$. It makes 
 sense to reduce the latter modulo $2$, $(3i+j+ij)$ being sent to $1 \in
 \Delta(2,3,3)$; on the other hand, $\Nr$ induces a 
 surjective homomorphism $\Delta(2,6,6)\to 3^\Z/9^\Z\cong \Z/2\Z$, and
 $(3i+j+ij)$ is sent to a generator of $\Z/2\Z$. It follows that
 $\Gamma^+(2)$ is a normal subgroup of $\Delta(2,6,6)$, with quotient
 group $\mathfrak{A}_4\times \Z/2\Z$.

 \noindent We next remark that every element of $\Gamma^+(2)$ can be written 

\begin{align*}
 &a+bi+cj+dij,\text{ with }a,b,c,d \in \Z,\;b \equiv c\equiv d
 \not\equiv a \;(2),\\
 &a^2+b^2-3c^2-3d^2=1.
\end{align*} 
 Since $(1+i)(a+bi+cj+dij)(1+i)^{-1}=a+bi-dj+cij,$ we see that $(1+i)$
 normalizes  $\Gamma^+(2)$, and conclude that $\Gamma^+(2)$ is normal in
 $\Gamma^\ast$. The quotient group sits in an extension 
 \begin{equation*}
  1\to \Delta(2,3,3)\times \Z/2\Z\to G \to \Z/2\Z\to 1,
 \end{equation*}
 the quotient $\Z/2\Z$ being generated by $(1+i)$. Since the image of
 $(1+i)$ in $G$ is not central, we conclude that $G\cong
 \mathfrak{S}_4\times \Z/2\Z$.
\end{proof} 

\begin{cor}\label{cor::6.2.6} The automorphism group of the curve
 $\mathcal{X}^+(2)_\C$ is $\Delta^\ast(2,3,4)$.
\end{cor}

\begin{proof}
 It follows from \ref{cor::6.2.4} and \ref{lem::6.2.5} that $\Aut
 \mathcal{X}^+(2)_\C$ contains $\Delta^\ast(2,3,4)$. On the other hand,
 any automorphism of $\mathcal{X}^+(2)_\C$ lifts to an element of
 $SU(1,1)(\R)$ which normalizes $\mathcal{B}$. The image in
 $PSU(1,1)(\R)$ of such an element lies in $\Gamma^\ast$, whence
 $\Aut\;\mathcal{X}^+(2)_\C\subset
 \Gamma^\ast/\Gamma^+(2)=\Delta^\ast(2,3,4)$.
\end{proof}

\medskip By similar arguments, one shows that the groups
$\Delta(3,4,4)/\Gamma^+(2)$ and
\noindent $\diamondsuit(2,2,2,3)/\Gamma^+(2)$ are the octahedral group
$\Delta(2,3,4)\cong \mathfrak{S}_4$.

\newpage
One reads from the following table of ramification that $\tilde w_3$ is
a hyperelliptic involution of $\mathcal{X}^+(2)_\C$.

\begin{figure}[h]
 \begin{picture}(300,480)(-10,-20)
  \put(5,-14){\includegraphics[scale=0.8,clip]{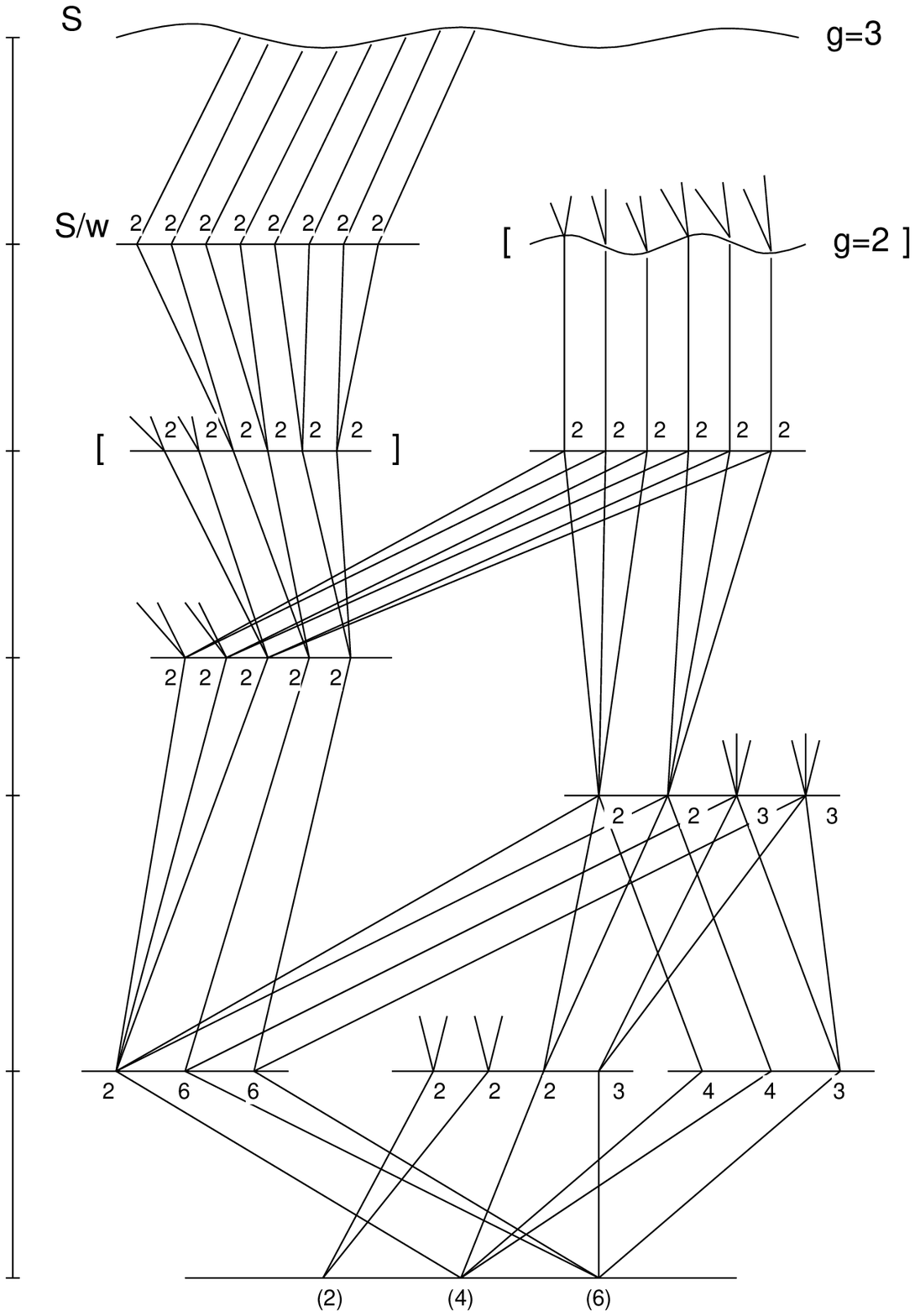}}
  \put(0,-20){\makebox(0,0)[c]{$\chi_{\mathrm{orb}}$}}
  \put(0,0){\makebox(0,0)[r]{$-\frac{1}{12}$}}
  \put(0,67){\makebox(0,0)[r]{$-\frac{1}{6}$}}
  \put(0,158){\makebox(0,0)[r]{$-\frac{1}{3}$}}
  \put(0,203){\makebox(0,0)[r]{$-\frac{1}{2}$}}
  \put(0,271){\makebox(0,0)[r]{{\footnotesize $-1$}}}
  \put(0,338){\makebox(0,0)[r]{{\footnotesize $-2$}}}
  \put(0,407){\makebox(0,0)[r]{{\footnotesize $-4$}}}
 \end{picture}
 \caption{}
\end{figure}
\vfill

\newpage
\begin{para}\label{para::6.2.7} Similar considerations apply to the congruence
 subgroup of level three $\Gamma^+(3)$. The Shimura curve 
$\mathcal{X}^+(3)_\C$ has genus $4$ and its automorphism group is of
 order $72$. We omit the details. 
\end{para}

\begin{para}\label{para::6.2.8} By the general theory of Shimura curves,
 $\mathcal{X}^+(2)$ is a moduli space for polarized abelian surfaces
 with quaternionic multiplication by $\mathcal{B}_{2.3}$ and level-two
 structure, \cf  II.\ref{sub:7.4}. It is defined over $\Q$, as well as
 $\mathcal{X}^\ast$ and the intermediate Shimura curves attached to
 $\Delta(3,4,4),\;\Delta(2,6,6),\;\diamondsuit(2,2,2,3)$ and
 $\diamondsuit(2,2,3,3)$.
 According to Y. Ihara, A. Kurihara \cite{oseoedscatmu}, A. Ogg and
 J.F. Michon, in each of the cases 
 $\Delta(2,4,6),\;\Delta(3,4,4),\;\Delta(2,6,6),\;\diamondsuit(2,2,2,3)$,
 the $\Q$-curve itself is $\P^1_\Q$, whereas for
 $\diamondsuit(2,2,3,3)$, it is the conic $x^2+y^2+3z^2=0$.\footnote{on
 the other hand, we do not know any hyperelliptic equation for
 $\mathcal{X}^+(2)$, even over $\C$.} The Atkin-Lehner involutions are
 defined over $\Q$.
\end{para}

\begin{para}\label{para::6.2.9} Any CM point on $\mathcal{X}^+(2)$
 parametrizes an abelian surface isogenous to the square of an elliptic 
 curve with complex multiplication. It follows from \ref{para::6.2.2}
 that the points 
 of ramification of the ``projection'' $\mathcal{X}^+(2)\to
 \mathcal{X}^\ast$ are CM points. More precisely, the corresponding
 imaginary quadratic fields are  

\begin{itemize}
 \item $\Q[\sqrt {-1}]$ for the points of ramification index $4$ (\eg 
       the image of the point $P\in \matheur{D}(0,1^-)$, \cf \ref{para::6.1.2})
 \item $\Q[ \sqrt{-3}]$ for the points of ramification index $6$ (\eg 
       the image of the point $O\in \matheur{D}(0,1^-)$) 
 \item $\Q[\sqrt {-6}]$ for the points of ramification index $2$ (\eg 
       the image of the point $Q\in \matheur{D}(0,1^-)$). 
\end{itemize}
\end{para}

\subsection{Triadic uniformization.}\label{sub::6.3}

\begin{para}\label{para::6.3.1} The Shimura curve $\mathcal{X}^+(2)$ has good
 reduction at any prime $p\neq 2,3$. As for the critical primes $2$ and
 $3$, the situation for $p=3$ is easier to describe than for $p=2$,
 since $3$ does not divide the level (II.\ref{sub:7.4}). 

 For the triadic \v{C}erednik-Drinfeld uniformization, the relevant
 quaternion algebra is the Hamilton definite quaternion algebra
 $\ovl{B}=B_{2.\infty}$ over $\Q$ with basis $(1,i,j,ij)$ and relations
 $i^2=j^2=-1, ij=-ji$. A maximal order (unique up to conjugation, {\it
 cf.} \cite[5.11]{adadq}) consists of the Hurwitz quaternions 

\begin{equation*}
\ovl{\mathcal{B}}=\ovl{\mathcal{B}}_{2.\infty} = \Z  + \Z i + \Z j + \Z
 \rho,\;\rho =\frac{1+i+j+ij}{2}. 
\end{equation*}
  
\noindent
The group $PGL_1(\ovl{\mathcal{B}})=\ovl{\mathcal{B}}^\times/\pm 1$ is
 finite, isomorphic to the tetrahedral group $\Delta(2,3,3)\cong
 \mathfrak{A}_4$.

\medskip 
Following the notation of  II.\ref{para:7.4.3}, we introduce the following
 $\{3,\infty\}$-arithmetic subgroups of $\ovl{B}^\times/\Q^\times$: 

\medskip $\bullet$ the image $\Gamma_3^+(N)$ of the group $\{g \in
 (1+N\ovl{\mathcal{B}}[\frac{1}{3}])^\times\mid \ord_3 \Nr(g) \;{\rm
 {even}}\}$ in $PGL_2(\Q_3)$ $=\Aut(\Omega_{\C_3})$  ($\Gamma_3^+(1)$
 is abbreviated into
 $\Gamma_3^+$),

 \medskip $\bullet$ the image $\Gamma_3^\ast$ in $PGL_2(\Q_3)$ of the group 
 
 \noindent $N(\mathcal{B} [\frac{1}{3}])=\{g\in \ovl{B}^\times \mid 
 g\ovl{\mathcal{B}}[\frac{1}{3}]=\ovl{\mathcal{B}}[\frac{1}{3}]g\}$,
 \ie $N(\mathcal{B}[\frac{1}{3}])/\Q^\times$.
 
\medskip According to \v{C}erednik-Drinfeld, one has:

\medskip
 $\bullet$
 ${\mathcal{X}^+(2)}_{\C_3}^{\mathrm{an}}=\Omega_{\C_3}/\Gamma_3^+(2)$ 
\medskip

$\bullet$ ${\mathcal{X}^+}_{\C_3}^{\mathrm{an}}=\Omega_{\C_3}/\Gamma_3^+$. 

\medskip

$\bullet$
 ${\mathcal{X}^\ast}_{\C_3}^{\mathrm{an}}=\Omega_{\C_3}/\Gamma_3^\ast$.

 In terms of triadic triangle (\resp quadrangle) groups, one can
 thus write 
 
\medskip
$\bullet$ $\Delta_3(2,4,6)=\Gamma_3^\ast$, 

$\bullet$ $\diamondsuit_3 (2,2,3,3)=\Gamma_3^+$.
\end{para}

\begin{para}\label{para::6.3.2} These uniformizations actually hold over
 $\Q_9$, the quadratic unramified extension of $\Q_3$. The corresponding
 Mumford curves over $\Q_9$ appear in \cite[IX]{sgamc}, without 
 referring to the \v{C}erednik theorem. The results of {\it loc.~cit.}
 allow to give a concrete description of the 
 arithmetic groups in terms of the generators $1,i,j,ij$ of $\ovl{B}$. 

 Let us fix a square root of $-1$ in $\Q_9$. A convenient matrix
 representation 
 for the elements of $\ovl{B}\otimes\Q_9$ is 

\begin{equation*}
a+bi+cj+dij \mapsto
 \begin{pmatrix}
  a+b\sqrt{-1} & c+d\sqrt{-1}
  \\-c+d\sqrt{-1} & a-b\sqrt{-1}  
 \end{pmatrix}
\end{equation*}  
 \noindent Of course, $\ovl{D}\otimes\Q_3\cong M_2(\Q_3)$, and any such
 matrix representation is conjugate to the previous after
 $\otimes\Q_9$ one by some $u\in PGL_2(\Q_9)$; $u$ is a kind of triadic
 Cayley transform which sends the Drinfeld space 
 $\Omega$ to the complement $u\Omega$ of the ``circle''
 $\{x+y\sqrt{-1}\mid x,y\in \Q_3, \;x^2+y^2=-1\}$ in the projective line  
 (here, we follow the usual convention that $PGL_2$ acts on the left on
 $\P^1$).
 
 After {\it loc.~cit.}, the congruence level subgroup of level two
 $PGL_1(\ovl{\mathcal{B}}[\frac{1}{3}])(2)$ is generated by 

 \begin{equation*}
  z_{\epsilon\epsilon'}=i+\epsilon j +\epsilon' ij,\;\;\epsilon=\pm 1,
   \epsilon' =\pm 1
 \end{equation*}

\noindent with the only relations $z_{\epsilon\epsilon'}^2=1$. The
 subgroup of index two
$\Gamma_3^+(2)$ consisting of words of even length (which is also the
 kernel of the homomorphism $PGL_1(\ovl{\mathcal{B}}[\frac{1}{3}])(2)\to
 3^\Z/9^\Z\cong \Z/2\Z$ induced by $\Nr$) is a {\it Schottky group of
 rank three} in $PGL_2(\Q_9)$, freely generated by
\begin{equation*}
\gamma_1=z_{++}z_{+-}=
 \begin{pmatrix}
  -1-2\sqrt{-1} & 2 \\-2 & -1+2\sqrt{-1} 
 \end{pmatrix},
\end{equation*}
\begin{equation*}
\gamma_2=z_{++}z_{-+}=
 \begin{pmatrix}
  -1-2\sqrt{-1} & -2\sqrt{-1}
  \\-2\sqrt{-1} & -1+2\sqrt{-1}  
 \end{pmatrix},
\end{equation*}
 \begin{equation*}
  \gamma_3=z_{++}z_{+-}=
   \begin{pmatrix}
    1 & 2-2\sqrt{-1} \\-2-2\sqrt{-1} & 1 
   \end{pmatrix}.
 \end{equation*}
Moreover $PGL_1(\ovl{\mathcal{B}}[\frac{1}{3}])$ is
 the semi-direct product of $PGL_1(\ovl{\mathcal{B}}[\frac{1}{3}])(2)$ 
 and the finite group $PGL_1(\ovl{\mathcal{B}})$, and $\Gamma_3^\ast$ is
 generated by $PGL_1(\ovl{\mathcal{B}}[\frac{1}{3}])$ and the 
 image of $(1+i)$.

At last, the dual graph of the stable reduction modulo $3$ of
 ${\mathcal{X}^+(2)}_{\Q_9}^{\mathrm{an}}$ is 

 \begin{figure}[h]
  \begin{picture}(100,40)(-50,-20)
   \put(-30,0){\makebox(0,0)[c]{$\bullet$}}
   \put(30,0){\makebox(0,0)[c]{$\bullet$}}
   \qbezier(-30,0)(-20,12)(0,12)
   \qbezier(30,0)(20,12)(0,12)
   \qbezier(-30,0)(0,10)(30,0)
   \qbezier(-30,0)(0,-10)(30,0)
   \qbezier(-30,0)(-20,-12)(0,-12)
   \qbezier(30,0)(20,-12)(0,-12)
  \end{picture}
  \caption{}
 \end{figure}
\end{para}

\begin{para}\label{para::6.3.3} The case of the prime $2$ is similar if one
 replaces $\mathcal{X}^+(2)$ by $\mathcal{X}^+(3)$\footnote{to deal with
 the diadic uniformization of $\mathcal{X}^+(2)$ itself would involve a
 certain finite \'etale covering of degree $12$ of $\Omega$.}. This is a
 Mumford curve over $\Q_4$, the unramified quadratic extension of
 $\Q_2$. The relevant quaternion algebra here is $B_{3.\infty}$. A
 maximal order (again unique up to conjugation) is given by $\Z+\Z i+\Z
 \frac{i+j}{2}+\Z \frac{1+ij}{2},\;\;i^2=-1,\;j^2=-3,\;ij=-ji$.
 
\begin{itemize}
 \item $\Delta_2(2,4,6)=\Gamma_2^\ast$, 
 \item $\diamondsuit_2 (2,2,3,3)=\Gamma_2^+$.
\end{itemize}

\end{para}

\subsection{The non-archimedean triangle groups $\Delta_3(2,4,6),
\Delta_2(2,4,6)$ and Herrlich's tree.}\label{sub::6.4}

\begin{para}\label{para::6.4.1} Our aim here is to explore the combinatorics
 of the finitely  
 generated groups $\Delta_3(2,4,6)$ and $ \Delta_2(2,4,6)$. 

\noindent It follows from the above considerations that
 $\Delta_3(2,4,6)$ sits in an exact sequence 
\begin{equation*}
1\to \an{\gamma_1,\gamma_2,\gamma_3}=\Gamma_3^+(2)\to \Delta_3(2,4,6)\to
 \mathfrak{S}_4\times \Z/2\Z =\Aut(\mathcal{X}^+(2))\to 1.
\end{equation*} 
 \noindent This extension is described by the outer action of
 $\mathfrak{S}_4\times \Z/2\Z$ on  $\an{\gamma_1,\gamma_2,\gamma_3}$. 
 Similarly, $\Delta_2(2,4,6)$ is an extension of a group of order $72$,
 $\Aut(\mathcal{X}^+(3))$, by a free group of rank four. In order to
 describe these extensions, we shall use F. Herrlich's work on 
 automorphisms of Mumford curves. 
\end{para}

\begin{para}\label{para::6.4.2} Let $S$ be a Mumford curve over some finite
 extension $K$ of $\Q_p$, and let $\Gamma 
 \subset PGL_2(K)$ be the corresponding Schottky group. It is known that
 any automorphism of $S$ lifts to $PGL_2(K)$: in fact, 
 the normalizer $N(\Gamma)$ of $\Gamma$ in sits $PGL_2(K)$ in an exact
 sequence \cite[IV]{sgamc}: 
\begin{equation*}
 1\to \Gamma\to N(\Gamma)\to  \Aut S\to 1.
\end{equation*} 
 \noindent Herrlich constructs a tree $\mathcal{T}$ endowed with an
 action of $N(\Gamma)$ as follows. One introduces a partition of the set
 of affinoid disks in $\P^1(K)$ into classes of the following form: a
 class contains a unique disk in $\A^1(K)$, of radius denoted by $r$ and
 all affinoid disks of the form $| z-a|\geq r$. The vertices of
 $\mathcal{T}$ are the classes of those disks which contain two fixpoints
 of some non-trivial element of $N(\Gamma)$, whose distance is the radius
 of the disk.  Two vertices are related by an edge if they can be
 represented by disks which are either both maximal, or nested in such a
 way that one is a maximal proper subdisk of the other. Herrlich shows
 that the natural action of $N(\Gamma)$ on affinoid disks induce an
 action without inversion on $\mathcal{T}$. This allows to apply the
 Bass-Serre theory of graphs of groups \cite{aas}, exploiting the fact
 that the stabilizer of vertices or edges belong to the short list of
 finite subgroups of $PGL_2(K)$
 \cite{herrlich80:_ordnun_autom_schot} \cite{herrlich80:_endlic_grupp}.
 When $S/\Aut(S)$ has genus $0$, the quotient $\mathcal{T}/N(\Gamma)$ is
 a finite tree of groups. The group $N(\Gamma)$ can then 
 be identified with its fundamental group, which is an amalgamated sum of
 the stabilizers of the vertices along the stabilizers of the edges.
\end{para}

\begin{para}\label{para::6.4.3} Let us apply this to
 $S=\mathcal{X}^+(2)_{\Q_9}$ with $N(\Gamma)=\Delta_3(2,4,6)$.  One has
 $(g-1)/\sharp(\Aut(S))=\frac{1}{24}$, and an inspection of the list of
 \cite{herrlich80:_ordnun_autom_schot} leaves only one possibility for
 the tree of groups $\mathcal{T}_3/\Delta_3(2,4,6)$, namely:

 \begin{center}
  \begin{picture}(100,35)(-50,-15)
   \put(-30,0){\makebox(0,0)[c]{$\mathfrak{D}_6$}}
   \put(-14.5,0){\makebox(0,0)[r]{$\bullet$}}
   \put(-15,0){\line(1,0){30}}
   \put(14.5,0){\makebox(0,0)[l]{$\circ$}}
   \put(30,0){\makebox(0,0)[c]{$\mathfrak{S}_4$}}
   \put(0,-10){\makebox(0,0)[c]{$\mathfrak{D}_3$}}
  \end{picture}
 \end{center}

 \medskip\noindent where $\mathfrak{D}_n$ are the dihedral groups.  
 In this case, the Herrlich tree $\mathcal{T}_3$ may be considered as a
 triadic analogue of the Schwarz-Escher tesselation:
 
\begin{figure}[h]
 \begin{picture}(180,180)(0,0)
  \put(0,0){\includegraphics[scale=0.5,clip]{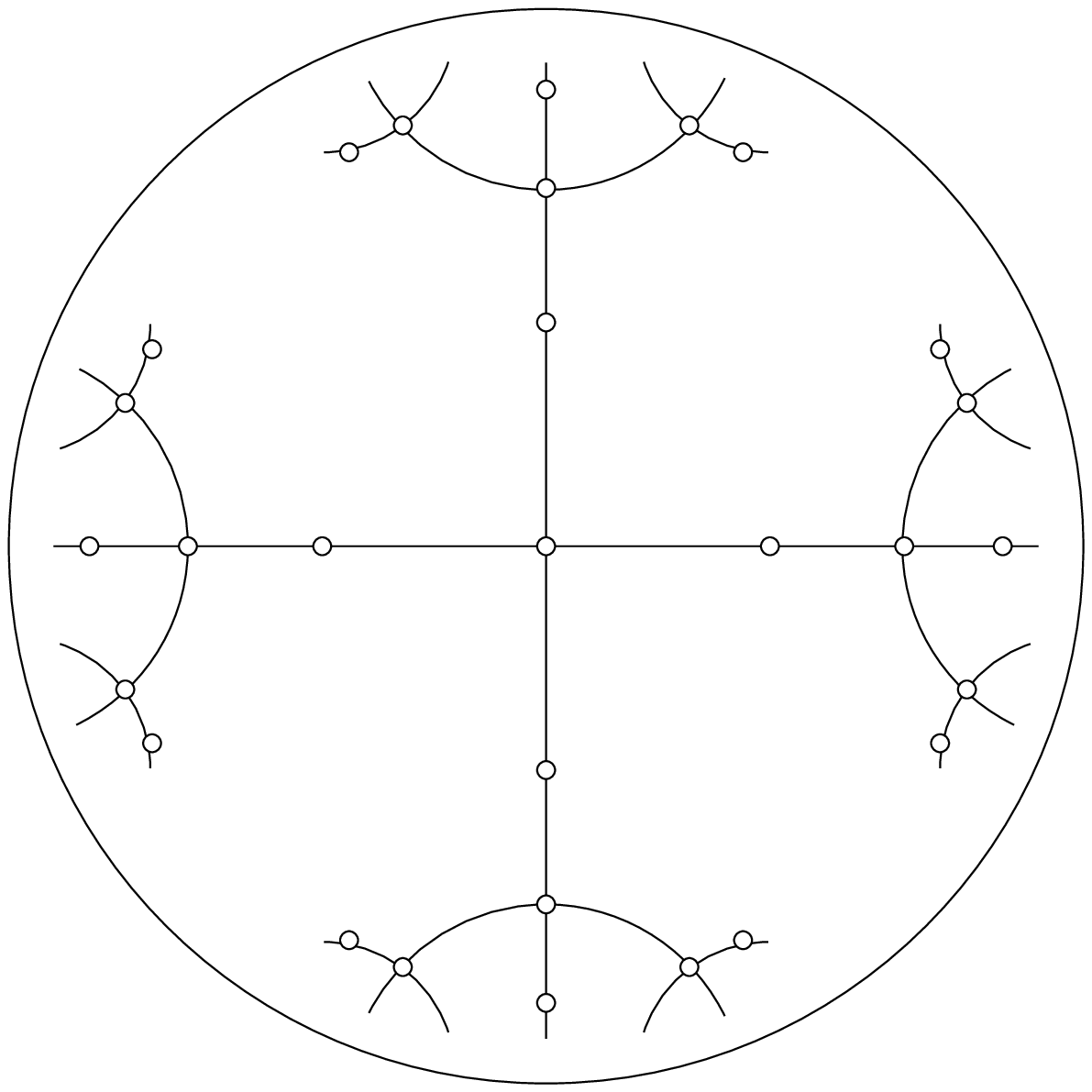}}
 \end{picture}
 \caption{}
\end{figure}

\end{para}
 
\begin{para}\label{para::6.4.4} In the case of $S=\mathcal{X}^+(3)_{\Q_4}$
 with $N(\Gamma)=\Delta_2(2,4,6)$, one still has
 $(g-1)/\sharp(\Aut(S))=\frac{1}{24}$, and an inspection of the list of
 \cite{herrlich80:_ordnun_autom_schot} leaves only two possibilities for
 the tree of groups $\mathcal{T}_2/\Delta_2(2,4,6)$, namely: 

 \begin{center}
  \begin{picture}(100,35)(-50,-15)
   \put(-30,0){\makebox(0,0)[c]{$\mathfrak{D}_6$}}
   \put(-14.5,0){\makebox(0,0)[r]{$\bullet$}}
   \put(-15,0){\line(1,0){30}}
   \put(14.5,0){\makebox(0,0)[l]{$\circ$}}
   \put(30,0){\makebox(0,0)[c]{$\mathfrak{D}_4$}}
   \put(0,-10){\makebox(0,0)[c]{$\mathfrak{D}_2$}}
  \end{picture}
 \end{center}

\noindent and 

 \begin{center}
  \begin{picture}(100,35)(-50,-15)
   \put(-30,0){\makebox(0,0)[c]{$\mathfrak{S}_4$}}
   \put(-14.5,0){\makebox(0,0)[r]{$\bullet$}}
   \put(-15,0){\line(1,0){30}}
   \put(14.5,0){\makebox(0,0)[l]{$\circ$}}
   \put(30,0){\makebox(0,0)[c]{$\mathfrak{S}_4$}}
   \put(0,-10){\makebox(0,0)[c]{$\mathfrak{D}_4$}}
  \end{picture}
 \end{center}

\medskip\noindent In order to select the right one, we remark that
 $\Delta_2(2,4,6)$ contains the finite group generated by $ij$ and the
 quotient of $(\mathcal{B}_{3.\infty})^\times=
 \bigl\{\pm1,\;\pm i,\; \frac{\pm i\pm j}{2},\; \frac{\pm 1\pm ij}{2}\bigr\}$
 by $\pm 1$. This group is
 isomorphic to $\mathfrak{D}_6\cong \Z/2\Z\;\times \mathfrak{D}_3$,
 which is not contained in $\mathfrak{S}_4$.  Since any finite subgroup
 of the amalgamated sum of two groups is conjugated to a subgroup of one
 of these groups \cite[4.3]{aas}, we conclude that
 $\Delta_2(2,4,6)\cong \mathfrak{D}_6
 \ast_{\mathfrak{D}_2}\mathfrak{D}_4$,
 not $\mathfrak{S}_4 \ast_{\mathfrak{D}_4}\mathfrak{S}_4$

\medskip\noindent
The Herrlich tree $\mathcal{T}_2$ may be considered as a diadic analogue of
the Schwarz-Escher tesselation:

\begin{figure}[h]
 \begin{picture}(180,180)(0,0)
  \put(0,0){\includegraphics[scale=0.5,clip]{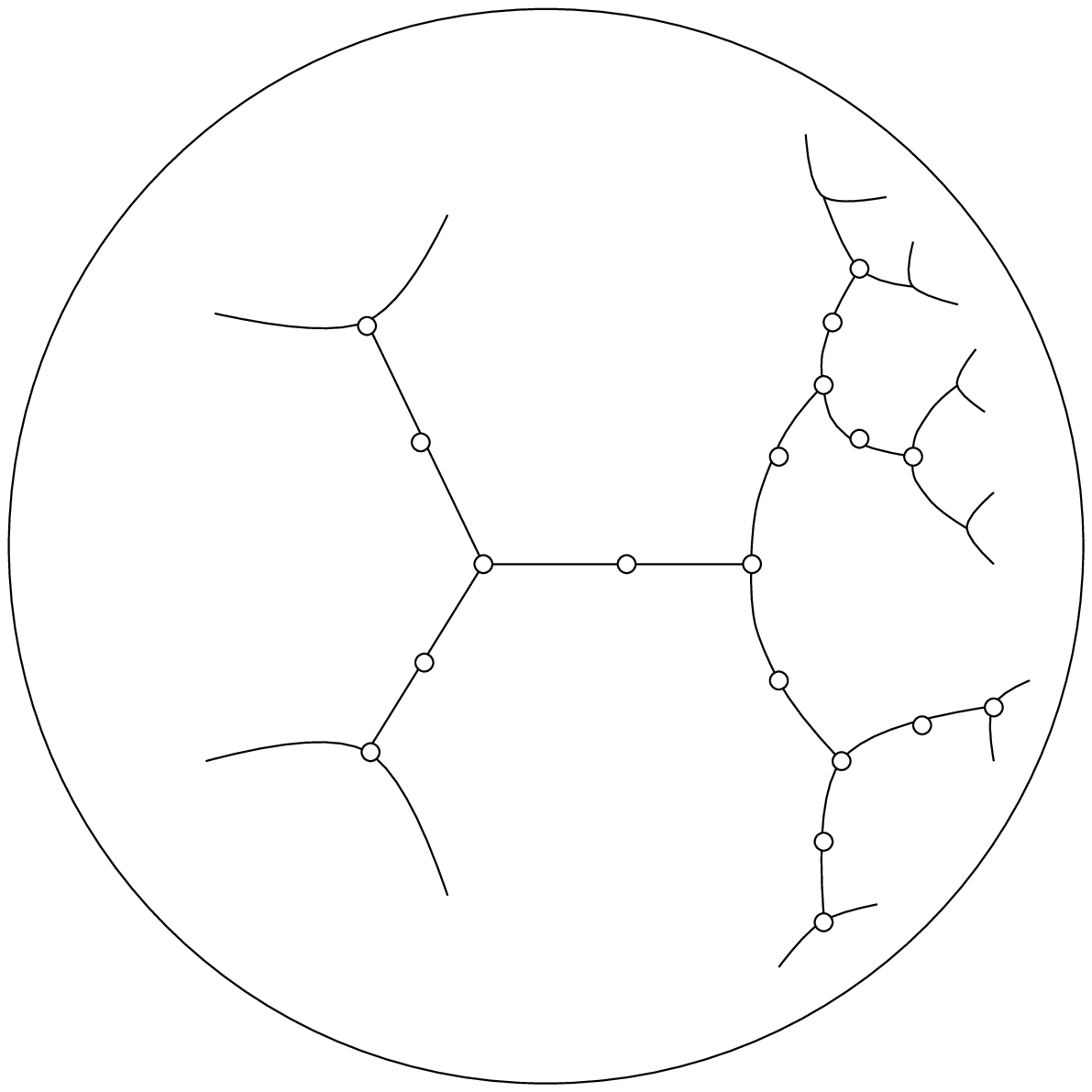}}
 \end{picture}
 \caption{}
\end{figure}

\bigskip\noindent
{\small (this is the set of points $u$ such that 
$j(\frac{e^{\pi i/3}-u}{1-e^{\pi i/3}u})\in [0,1728]$, 
where $j$ denotes the modular invariant).}

\medskip
It is important to pay attention to the fact that the symbol 
$\mathfrak{D}_6 \ast_{\mathfrak{D}_2}\mathfrak{D}_4$ is not free from
ambiguity: namely, there are essentially two non-equivalent ways of
embedding $\mathfrak{D}_2$ into $\mathfrak{D}_6\times \mathfrak{D}_4$:

{\it $1$st way$)$} the centers of $\mathfrak{D}_6$ and $\mathfrak{D}_4$ 
restrict to the same subgroup $\Z/2\Z$
of $\mathfrak{D}_2$. On identifying $\mathfrak{D}_6$ and $\mathfrak{D}_4$ with 
the matrix groups
\begin{equation*}
 \Bigl\langle\begin{pmatrix} 0& 1 \\ 1& 0
\end{pmatrix}, \begin{pmatrix} 0& -1 \\ 1& 1
\end{pmatrix}\Bigr\rangle\text{ and }
 \Bigl\langle\begin{pmatrix} 0& 1 \\ \pm1& 0
\end{pmatrix}\Bigr\rangle
\end{equation*}
respectively, the corresponding amalgam 
$\mathfrak{D}_6 \ast_{\mathfrak{D}_2}\mathfrak{D}_4$ is isomorphic 
to $GL_2(\Z)$;

{\it $2$nd way $)$} the trace of the center of $\mathfrak{D}_6\times 
\mathfrak{D}_4$ on $\mathfrak{D}_2$  is $\mathfrak{D}_2$ itself. The
corresponding amalgam $\mathfrak{D}_6
\ast_{\mathfrak{D}_2}\mathfrak{D}_4$ maps onto  
a semidirect product of $(\Z/3\Z)^2$ and ${\mathfrak{D}_4}$ (the 
latter acting on $(\Z/3\Z)^2$ as a matrix group
$\Bigl\langle
\begin{pmatrix}
 0& 1 \\ \pm1& 0
\end{pmatrix}
\Bigr\rangle$).
\end{para}

The following theorem summaries our results and settles the remaining 
ambiguity:

\begin{thm}\label{thm::6.4.5}
 \begin{enumerate}
  \item $\Delta_3(2,4,6)\cong \mathfrak{D}_6 
        \ast_{\mathfrak{D}_3} \mathfrak{S}_4$.
  \item $\Delta_2(2,4,6)\cong\mathfrak{D}_6\ast_{\mathfrak{D}_2}
	\mathfrak{D}_4$, where the embedding of $\mathfrak{D}_2$ into
	$\mathfrak{D}_6\times \mathfrak{D}_4$ is such that the trace of
	the center of $\mathfrak{D}_6\times \mathfrak{D}_4$ on
	$\mathfrak{D}_2$ is $\mathfrak{D}_2$ itself.
 \end{enumerate} 
\end{thm}

\begin{proof}
 Only the last assertion remains to be established. Let $\mathcal{S}$ be 
 the diadic Schwarz orbifold
 $(\P^1;(0,2),(1,4),(\infty,6))$.  There are surjective homomorphisms 
 (3.5.5, \ref{pro::4.3.3})
 \begin{equation*}
  \pi_1^{\mathrm{temp}}(\P^1\setminus
   \{0,1,\infty\})\to
   \pi_1^{\mathrm{orb}}(\mathcal{S})\to \Delta_2(2,4,6).
 \end{equation*} 
 The images 
 $\bar\gamma_0,\bar\gamma_1,\bar\gamma_\infty\in  \Delta_2(2,4,6)$
 of local monodromy elements $\gamma_0,\gamma_1,\gamma_\infty \in 
 \pi_1^{\mathrm{temp}}(\P^1\setminus
 \{0,1,\infty\})$ are elements of order $2,4,6$ respectively. By 
 \ref{pro::2.3.9}, the smallest closed normal
 subgroup of $\pi_1^{\mathrm{temp}}(\P^1\setminus
 \{0,1,\infty\})$ containing $\gamma_1,\gamma_\infty$ is 
 $\pi_1^{\mathrm{temp}}(\P^1\setminus
 \{0,1,\infty\})$ itself. It follows that the smallest normal
 subgroup of the discrete group $\Delta_2(2,4,6)$ containing 
 $\bar\gamma_1,\bar\gamma_\infty$ is $\Delta_2(2,4,6)$ itself.
 It thus suffices to show that if we where in the case $\mathfrak{D}_6 
 \ast_{\mathfrak{D}_2}\mathfrak{D}_4\cong
 GL_2(\Z)$, then $\bar\gamma_1,\bar\gamma_\infty\in SL_2(\Z)$. This is 
 clear, since any element of order $6$ (\resp $4$) in $ GL_2(\Z)$ is
 conjugate to $\pm
 \begin{pmatrix}
  0 & -1 \\ 1& 1
 \end{pmatrix}$
 (\resp to
 $\pm\begin{pmatrix}
      0& -1 \\ 1& 0
     \end{pmatrix}$).
\end{proof}

\medskip Part (ii) of this theorem allows us to answer at last 
question \ref{que::2.3.11} concerning
local monodromies $\gamma_0,\gamma_1,\gamma_\infty$ at $0,1,\infty$ 
in the temperate fundamental group
$\pi_1^{\mathrm{temp}}(\P^1_{\C_p}\setminus \{0,1,\infty\})$:

\begin{cor}\label{cor::6.4.6} Question \ref{que::2.3.11} has a negative
 answer for $p=2$: one cannot find local monodromies
 $\gamma_0,\gamma_1,\gamma_\infty$ such that 
 $\gamma_0.\gamma_1.\gamma_\infty = 1$.
\end{cor}

\begin{proof}
 Assume that such local monodromies exist.  By means of their images
 $\bar\gamma_0,\bar\gamma_1,\bar\gamma_\infty \in \Delta_2(2,4,6)\cong 
 \mathfrak{D}_6 \ast_{\mathfrak{D}_2}\mathfrak{D}_4 $ (notation of 
 the previous proof), it would be possible to construct
 a homomorphism
 \begin{equation*}
  \phi: \Delta (2,4,6)\to\Delta_2(2,4,6),
 \end{equation*}
 where $\Delta (2,4,6)$ 
 denotes the standard (archimedean) triangle group.
 Notice that the image of $\phi$ cannot be contained in any conjugate 
 of $\mathfrak{D}_6$ or $\mathfrak{D}_4$, since none of
 these finite groups contains elements of order $6$ and of order $4$ 
 simultaneously.

  According to \cite[6.3.5]{aas}, any (cocompact) triangle group 
 $\Delta$ has the
 fixpoint property $(FA)$: any action (without inversion) of $\Delta$ 
 on any tree has a fixpoint. It follows that the
 quotient
 $\Im\phi$ of $\Delta (2,4,6)$ would have property
 $(FA)$ as well. Since $Im\,\phi$ is a subgroup of the amalgam
 $\mathfrak{D}_6 \ast_{\mathfrak{D}_2}\mathfrak{D}_4$, this would 
 imply, according to \cite[6.2, prop. 21]{aas}, that
 $Im\,\phi$ is contained in a conjugate of $\mathfrak{D}_6$ or 
 $\mathfrak{D}_4$, a contradiction. 
\end{proof}

\begin{rem}
 The arithmetic $p$-adic triangle groups 
 $\Delta (e_0,e_1,e_\infty)$ do {\it not} have property
 $(FA)$. Indeed, their description as $p$-unit groups in quaternion 
 algebras shows that they admit a two-dimensional
 representation for which the eigenvalues are non-integral in general. 
 This property contradicts $(FA)$, \cf \cite[6.2, prop. 22]{aas}.
   Actually, Kato's description of $p$-adic triangle groups of Mumford 
 type as amalgams, together with \ref{thm::5.3.7}, implies that
 {\it no hyperbolic $p$-adic triangle group has property $(FA)$} (in 
 contrast to archimedean triangle groups).  Hence no hyperbolic $p$-adic
 triangle group can be generated by three elements $g_0,g_1,g_\infty$ of
 finite order satisfying $g_0g_1g_\infty=1$.
\end{rem}

\subsection{$F(\frac{1}{24}, \frac{7}{24},\frac{5}{6}; z)$ and the period 
mapping. Global analytic triadic continuation.}\label{sub::6.5}

\begin{para}\label{para::6.5.1} Rather than $\Delta(2,4,8)$, we now
 prefer to work with $\Delta(6,2,4)$ (this amounts to permute
 $0,1,\infty$ by $t\mapsto z=\frac{1}{1-t}$, and is more faithful to the
 first picture of this section which places the vertex 
 $O$ at the center of the disk). 

 The triangle group
$\Delta(6,2,4)$ is the projective monodromy group of the hypergeometric
 differential operator 
$L_{\frac{1}{24}, \frac{7}{24},\frac{5}{6}}$. More precisely,
 $L_{\frac{1}{24}, \frac{7}{24},\frac{5}{6}}(y)=0$ 
 is a {\it uniformizing differential equation} for the orbifold
 $\mathcal{X}^\ast$ viewed as the Schwarz orbifold 
\begin{equation*}
\left(\P^1,\;
 \begin{matrix}
  0& 1& \infty \\e_0=6& e_1=2& e_\infty = 4
 \end{matrix}
\right)
\end{equation*}
More explicitly, $F(\frac{1}{24}, \frac{7}{24},\frac{5}{6}; z)$ and
 $z^{1/6}F(\frac{5}{24}, \frac{11}{24},\frac{7}{6}; z)$ are 
independent solutions of $L_{\frac{1}{24}, \frac{7}{24},\frac{5}{6}}$
 and for a suitable constant $\kappa$, the mapping 

\begin{equation*}
\tau(z)= \kappa^{-1}. \frac{z^{1/6}F(\frac{5}{24},
 \frac{11}{24},\frac{7}{6}; z) }{F(\frac{1}{24},
 \frac{7}{24},\frac{5}{6}; z)}
\end{equation*} 

\noindent sends the upper half plane bijectively to the basic triangle
 $OQP$ of the Schwarz tesselation, and 
 $\tau(0)=O,\;\tau(1)=Q,\;\tau(\infty)=P$. 

\begin{figure}[h]
 \begin{picture}(320,100)(0,-15)
  \put(0,0){\includegraphics[scale=0.8,clip]{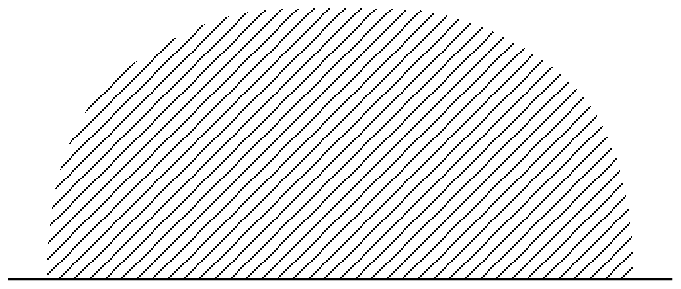}}
  \put(57,-1){\makebox(0,0)[b]{\small $\bullet$}}
  \put(57,-8){\makebox(0,0)[c]{\small $0$}}
  \put(93,-1){\makebox(0,0)[b]{\small $\bullet$}}
  \put(93,-8){\makebox(0,0)[c]{\small $1$}}
  \put(180,25){\vector(1,0){50}}
  \put(205,35){\makebox(0,0)[c]{$\tau$}}
  \put(250,0){\includegraphics[scale=0.8,clip]{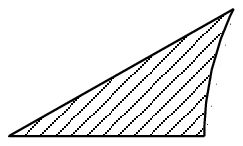}}
  \put(250,-8){\makebox(0,0)[c]{\small $O$}}
  \put(302,-8){\makebox(0,0)[c]{\small $Q$}}
  \put(312,38){\makebox(0,0)[c]{\small $P$}}
 \end{picture}
 \caption{}
\end{figure}

\noindent It extends to a multivalued analytic function on $\P^1
 \setminus \{0,1,\infty\}$ with value in $\matheur{D}(0,1^-)$. The 
inverse mapping
$z(\tau): \matheur{D}(0,1^-)\to \P^1$ is the (single-valued)
$\Delta(2,4,6)$-automorphic uniformizing mapping of the Shimura orbifold
 $\mathcal{X}^\ast $.
\end{para}

\begin{para}\label{para::6.5.2}  Actually, $\tau$ expresses the period mapping
 $\mathcal{P}$ for ``fake elliptic curves'' with multiplication by 
 $\mathcal{B}_{2.3}$, \cf  II.\ref{sub:1.2} and
 II.\ref{sub:7.4}: more precisely, 
 the multivalued function $\tau(z)$ is related to 
 the commutative diagram  
 \begin{equation*}
  \begin{CD}
   \widetilde{{\mathcal{X}^\ast_\C}^{\mathrm{an}}}
    @>{}^{\begin{subarray}\mathcal{P}\\\sim\end{subarray}}>>
   \mathcal{D}=\matheur{D}(0,1^-) @. \quad \subset \;\mathcal{D}^\vee=\P^1 \\
   @V\mathcal{Q}VV @VVV\\
   \left(\P^1,\;
   \begin{matrix}
    0& 1& \infty \\
    6& 2& 4
   \end{matrix}
   \right)
   ={\mathcal{X}^\ast_\C}^{\mathrm{an}}
   @>>> \mathcal{D}/\Delta(2,4,6)
  \end{CD}
 \end{equation*}
by
 \begin{equation*}
  \tau(z)= \mathcal{P}\circ {\mathcal{Q}}^{-1},\; z(\tau )= {\mathcal{Q}}\circ
   \mathcal{P}^{-1}.
 \end{equation*} 
\end{para}

\begin{para}\label{para::6.5.3} Using Gauss' formula for the value of
 hypergeometric functions at $1$, one can compute the constant $\kappa$, 
 \cf \cite[392]{f2}:
 \begin{equation*}
 \kappa=\frac{\Gamma (\frac{7}{6})}{\Gamma
  (\frac{5}{6})}\sqrt\frac{\Gamma (\frac{1}{24})\Gamma
  (\frac{7}{24})\Gamma (\frac{13}{24})\Gamma (\frac{19}{24}) 
 }{\Gamma (\frac{5}{24}) \Gamma (\frac{11}{24})\Gamma
 (\frac{17}{24})\Gamma (\frac{23}{24})}. 
 \end{equation*}

 Using the functional equations of $\Gamma$, one checks that $\kappa
 \sim \Bigl(\frac{\Gamma(\frac{1}{3})}{\Gamma(\frac{2}{3})}\Bigr)^3$
 $\pmod{\ovl{\Q}^\times}$, in conformity with II.\ref{thm:7.6.4} (taking into
 account the fact that $z=0$ corresponds to a fake elliptic curve with
 complex multiplication by an order in $M_2(\Q(\sqrt{-3}))$). 
\end{para}

\begin{para}\label{para::6.5.4} For $p\neq 2,3$, the $p$-adic radius of
 convergence of $F(\frac{1}{24}, \frac{7}{24},\frac{5}{6}; z)$ and 
 $F(\frac{5}{24}, \frac{11}{24},\frac{7}{6}; z)$ is $1$. The unit disk
 $|z|_p<1$ is ordinary or supersingular according 
 to whether $p\equiv 1$ $\pmod{3}$ or not. The ratio 
 $\frac{z^{1/6}F(\frac{5}{24}, \frac{11}{24},\frac{7}{6}; z)
 }{F(\frac{1}{24}, \frac{7}{24},\frac{5}{6}; z)}$ expresses a local 
 $p$-adic period mapping. 

 The $p$-divisible group of the fake elliptic curve parametrized by $z$
 is isogenous to the square of the $p$-divisible group of an elliptic
 curve; the period 
 mapping is thus completely analogous to one of those studied 
 in II.\ref{para:6.3.1}, \ref{para:6.3.3} (Dwork, 
 Gross-Hopkins). From the viewpoint of analytic continuation, this is
 essentially the situation studied in  I.\ref{sec-tale} 
 (for $F(\frac{1}{2},\frac{1}{2},{1}; z)$).
\end{para}

\begin{para}\label{para::6.5.5} For $p=3$ (\resp $p=2$), the
 situation is radically different; the $p$-adic radius of convergence of 
 $F(\frac{1}{24}, \frac{7}{24},\frac{5}{6}; z)$ and
 $F(\frac{5}{24}, \frac{11}{24},\frac{7}{6}; z)$ is $3^{-3/2}$
 (\resp $2^{-6}$).

 The ratio
 $ \frac{z^{1/6}F(\frac{5}{24},\frac{11}{24},\frac{7}{6}; z)}
 {F({\frac{1}{24},{\frac{7}{24}},\frac{5}{6}}; z)}$
 expresses the Drinfeld $3$-adic period mapping: 

 \begin{center}
  \begin{picture}(120,125)(0,-10)
   \put(45,105){\makebox(0,0)[c]{${}^\mathcal{P}$}}
   \put(45,95){\makebox(0,0)[c]{$\sim$}}
   \put(0,90){\makebox(0,0)[c]{$\til{\mathcal{X}^+(2)_{\C_3}^{\mathrm{an}}}$}}
   \put(0,75){\vector(0,-1){15}}
   \put(0,45){\makebox(0,0)[c]{${\mathcal{X}^+(2)}_{\C_3}^{\mathrm{an}}$}}
   \put(0,30){\vector(0,-1){15}}
   \put(0,0){\makebox(0,0)[c]{${\mathcal{X}^{\ast\,\mathrm{an}}_{\C_3}}2$}}
   \put(30,90){\vector(1,0){30}}
   \put(80,90){\makebox(0,0)[dc]{$\mathcal{D}$}}
   \put(97,90){\makebox(0,0)[c]{$\subset$}}
   \put(110,90){\makebox(0,0)[dl]{$\mathcal{D}^{\vee}=\P^1$}}
   \put(80,75){\vector(0,-1){60}}
   \put(80,0){\makebox(0,0)[c]{$\mathcal{D}/\Delta_3(2,4,6)$}}
   \put(25,0){\vector(1,0){15}}
  \end{picture}
 \end{center}

 \noindent 
 Here $\mathcal{D}$ is the image of the Drinfeld space $\Omega_{\C_3}$
 by a suitable homography (it would be interesting to 
 determine this homography explicitly). That this homography is actually
 defined over $\ovl{\Q}$ 
 comes from II.\ref{thm:7.6.4} (ramified variant), I.\ref{para-ramif},
 \ref{thm-cm-period} (the point is
 that the triadic periods, in the sense of
 I.\ref{sub-p-adic-betti-lattice-supersing}, of the abelian surface 
 with complex multiplication by $M_2(\Q[\sqrt{-3}])$ parametrized by
 $z=0$, belong to $\ovl{\Q}$).

 \bigskip Although $F(\frac{5}{24}, \frac{11}{24},\frac{7}{6};
 z)/F(\frac{1}{24}, \frac{7}{24},\frac{5}{6}; z)$ has very small triadic
 radius of convergence, its pull-back on the Shimura curve
 $\mathcal{X}^+(2)$ extend to a global multivalued meromorphic
 function, with poles only above $0,1,\infty$ (in other words, its
 pull-back on $\widetilde{\mathcal{X}^+(2)}_{\C_3} \cong \Omega_{\C_3}$
 is single-valued).

 Note however that the same is not true for $F(\frac{1}{24},
 \frac{7}{24},\frac{5}{6}; z)$ and $F(\frac{5}{24},
 \frac{11}{24},\frac{7}{6}; z)$ individually, for reasons of exponents
 at $\infty$.  Nevertheless, they enjoy a similar property, on replacing
 $\mathcal{X}^+(2)$ by a suitable ramified finite covering.

 \medskip Similarly, in the diadic case: the pull-back of
 $F(\frac{5}{24}, \frac{11}{24},\frac{7}{6}; z)/F(\frac{1}{24},
 \frac{7}{24},\frac{5}{6}; z)$ on the Shimura curve $\mathcal{X}^+(3)$
 extends to a global multivalued meromorphic functions (with poles only
 above $0,1,\infty$).
\end{para}

\subsection{A triadic analogue of Gauss' formula, and special values of
$F(\frac{1}{24}, \frac{7}{24},\frac{5}{6}; z)$ at CM points.}\label{sub::6.6}

\begin{para}\label{para::6.6.1} For our hypergeometric function, Gauss'
 formula reads
 \begin{equation*}
  F\Bigl(\frac{1}{24}, \frac{7}{24},\frac{5}{6};1\Bigr)=
 \frac{\Gamma(\frac{5}{6})\Gamma(\frac{1}{2})}
 {\Gamma(\frac{19}{24})\Gamma(\frac{13}{24})}.
 \end{equation*}

\noindent If one considers $F(\frac{1}{24}, \frac{7}{24},\frac{5}{6};
 z)$ as a multivalued analytic function on
 $\P^1\setminus\{0,1,\infty\}$, the values at $1$ generate a
 $\ovl{\Q}$-subspace of $\C$ of dimension $\leq 2$, namely:  
\begin{equation*}
 \frac{\Gamma(\frac{5}{6})\Gamma(\frac{1}{2})}
  {\Gamma(\frac{19}{24})\Gamma(\frac{13}{24})}\ovl{\Q} \,+ \,
  \frac{\Gamma(\frac{5}{6})\Gamma(-\frac{1}{2})}
  {\Gamma(\frac{1}{24})\Gamma(\frac{7}{24})}\ovl{\Q}.
\end{equation*}
Similarly, its values at $\infty$ generate the $\ovl{\Q}$-subspace 

\begin{equation*}
 \frac{\Gamma(\frac{5}{6})\Gamma(\frac{1}{4})}
  {\Gamma(\frac{19}{24})\Gamma(\frac{7}{24})}\ovl{\Q} \, +\,
  \frac{\Gamma(\frac{5}{6})\Gamma(-\frac{1}{4})}
  {\Gamma(\frac{1}{24})\Gamma(\frac{13}{24})}\ovl{\Q}.
\end{equation*}   
\end{para}

\begin{para}\label{para::6.6.2} Let us consider our hypergeometric
 function as a triadic function $F(\frac{1}{24},
 \frac{7}{24},\frac{5}{6}; 1)_3$, and more precisely as a multivalued
 $3$-adic analytic function on a suitable finite \'etale covering of
 $\P^1\setminus\{0,1,\infty\}$. 
 
 The values taken by $F(\frac{1}{24}, \frac{7}{24},\frac{5}{6}; z)_3$ at
 $1$ belong to $\ovl{\Q}$: this follows from \ref{para::5.6.3}, together
 with the fact that the ``periods'' (in the sense of
 I.\ref{sub-p-adic-betti-lattice-supersing}) attached to the 
 abelian surface with complex multiplication by $M_2(\Q[\sqrt{-24}])$ 
 parametrized by $z=1$, belong to $\ovl{\Q}$, by I.\ref{para-ramif}, 
 I.\ref{thm-cm-period}. 

 \medskip In the same way, one shows that the values taken by
 $F(\frac{1}{24}, \frac{7}{24},\frac{5}{6}; z)_3$ at $\infty$ belong to 
 $(\Gamma_3(\frac{3}{4})/\Gamma_3(\frac{1}{4}))^{-1/2}
 \ovl{\Q}=\Gamma_3(\frac{1}{4})\ovl{\Q}$.
 Conjecture I.\ref{coj-no-complex-multi} predicts that
 $\Gamma_3(\frac{1}{4})$ is a transcendental number.
\end{para}

\begin{para}\label{para::6.6.3} In \cite[Table 2]{elkies98:_shimur}, one finds a
 table of rational CM points $z$ (rational points which correspond to
 abelian surfaces with complex multiplication by an order in
 $M_2(\Q(\sqrt{-d}))$ for some fundamental discriminant $-d$); with the
 notation of {\it loc.~cit.}, they are given by $z=B/(B-A)= \pm B /|
 C|$. Some of them are: 

\vspace{5mm}
\begin{center}
 \renewcommand{\arraystretch}{1.2}
 \begin{tabular}{c @{\quad\qquad}l}
 $ z            $ & $ d=-\mathrm{disc}   $ \\ \hline 
 $ -3^3/2^2 7^2 $ & $ d=2^2.3.7$ \\
 $ -5^3/3^7     $ & $ d=2^3.5  $ \\
 $ 2^{10}/7^4   $ & $ d=3.17   $ \\
 $ -2^{10}/3^7  $ & $ d=19.    $ 
 \end{tabular} 
\end{center}
\vspace{3mm}

 Notice that for the first value of $z$, $F(\frac{1}{24},
 \frac{7}{24},\frac{5}{6}; \frac{-3^3}{2^2 7^2})_3$ converges. Moreover,
 since $(\frac{2^2.7}{3})=+1$, it follows from I.\ref{para-ramif},
 I.\ref{thm-cm-period}, that 
 the triadic ``periods'' attached to the abelian surface $A_{-3^3/ 2^2 7^2}$ 
 with complex multiplication by $M_2(\Q[\sqrt{-d}])$
 parametrized by $z=1$, is algebraic (like the periods of $A_0$). By 
 \ref{para::5.6.3}, it follows that the triadic number
 $F(\frac{1}{24}, \frac{7}{24},\frac{5}{6}; \frac{-3^3}{2^2 7^2})_3 $ is
 algebraic (can one compute this algebraic number?). 
 
 \medskip 
 On the other hand since
 $\Q(\sqrt{-d})\neq \Q(\sqrt{-3})$ in that case ($d=2^2.3.7 $), the real
 number $F(\frac{1}{24}, \frac{7}{24},\frac{5}{6}; \frac{-3^3}{2^2 7^2})$
 is transcendental; this follows from Wolfart's theorem 
 \cite{wolfart88:_werteh_funkt}
 (whose proof has been recently completed by work of Cohen-W\"ustholz
 \cite{cohen01:_number_theor} and Edixhoven-Yafaev).

 \medskip It can be proved, a contrario, that if the complex evaluation
 $F(\frac{1}{24}, \frac{7}{24},\frac{5}{6}; \xi)$ is algebraic for some
 algebraic value $\xi$ in the disc of convergence, then any $p$-adic
 evaluation $F(\frac{1}{24}, \frac{7}{24},\frac{5}{6}; \xi)_p$ which makes
 sense is also algebraic: indeed, by 
 \cite{wolfart88:_werteh_funkt}, \cite{cohen01:_number_theor}, such a
 $\xi$ parametrizes a CM abelian surface by some order in
 $M_2(\Q[\sqrt{-3}])$, and one concludes again using $p$-adic Betti
 lattices as in I.\ref{sub-p-adic-betti-lattice-supersing}.
 
 This is actually a very special case of a general, 
 but conjectural, ``principle of global relations'' for special values of
 solutions of Gauss-Manin connections, \cf \cite{tdmeigdvpdg}.
\end{para}

\appendix
\addtocontents{toc}{\protect\par\vskip10mm}

\chapter{Rapid Course in $p$-adic Analysis.}\label{app:rapid}
\renewcommand{\thechapter}{\Alph{chapter}}
\addtocontents{toc}{\protect\par\vskip2mm\hskip80mm by F. Kato \par\vskip5mm}

\vspace{1ex}
\begin{center}
By {\sc F.\ Kato}
\end{center}

\section{Introduction.}
%

In this appendix, $K$ always denotes a complete field 
with respect to a non-archimedean valuation $|\cdot|\colon K\rightarrow
\R_{\geq 0}$.
The norm $|\cdot|$ is almost always assumed to be non-trivial, unless
otherwise stated.

Let us assume for a while that $K$ is algebraically closed.
We view the field $K$ as an affine space, as we do in complex analysis; 
$K$ is a metrized space with the metric induced from the
valuation $|\cdot|$.
We can then follow the lines of classical complex analysis and define
convergence of power series, Taylor and Laurent
expansions, etc...
%

For instance let us define ``holomorphic'' functions to be $K$-valued
functions which are locally expressed by convergent power
series\footnote{The reader may wonder whether there is another approach
by means of Cauchy-Riemann type differential equations.  This is unlikely
because the differential calculus does not go well, see
\cite[Chap.\ 4]{pn}.}  (cf.\ \cite[Chap.\ 4]{pn}). 
This approach drives us however to several
problems, which come mainly from the fact that $K$ with the
metric topology is totally disconnected; for example:
\begin{enumerate}
 \item for an open set $U\subseteq K$ the ring of all such holomorphic
       functions on $U$ is huge. So is already the subring of all locally 
       constant functions.
 \item It can be shown that the sheaf of germs of such functions satisfies
       the principle of unique continuation; but not in a satisfactory way.
       In particular, \lemref{lem-unique-cont} of \chapref{chap-ana-p-period}
       holds but is of no use,
       because there is no non-empty connected open subset.
\end{enumerate}

The trouble becomes more apparent when we think the \emph{local}
representability by convergent power series in terms of
\emph{coverings}.  In the complex analytic situation, when we speak of
the holomorphy of a function $f$ defined over a connected open set $U$,
we tacitly take a open covering $\{U_i\}_{i\in I}$ of $U$ consisting of
sufficiently small open neighborhoods such that $f$ restricted on each
member $U_i$ can be seen as evaluation of a power series convergent in
$U_i$.  The point is that, since $U$ is connected, members of $\{U_i\}$
must have overlaps so that the local properties can be transmitted to
whole $U$.  But, in the non-archimedean situation, such coverings may be
refined in such a way that there are no overlap.  This is why
the analytic continuation does not work well.  So, for the sake of a
reasonable analytic theory, we have to ``limit'' coverings so that we
cannot take arbitrarily fine refinements.

A first method is due to Krasner.
He introduced the so-called {\it quasi-connected sets} and, on such
subsets, defined analytic functions as uniform convergent limits of
rational functions.

A more modern and systematic treatment was introduced by Tate, 
with the so-called {\it Rigid Analysis}.
Let us briefly view the main idea:
we first set
\begin{align*}
   & K\{t_1,\ldots,t_n\}\\
 = & \Biggl\{
 \sum_{\nu_1,\ldots,\nu_n\geq 0}
 a_{\nu_1,\ldots,\nu_n}t_1^{\nu_1}\cdots t_n^{\nu_n}
  \in K[[t_1,\ldots,t_n]]\,
 \Bigg|
 \begin{split}
  & a_{\nu_1,\ldots,\nu_n}\rightarrow 0\ \mathrm{for}\\
  & \nu_1+\cdots+\nu_n\rightarrow\infty
 \end{split} 
\Biggr\},
\end{align*}
which is now called the {\it Tate algebra} over $K$.
This is the ring of all functions expressed by power series convergent
on the closed disk $\matheur{D}(0,1^{+})=\{z\in K^n\ |\ |z_i|\leq 1\}$.
The algebra $K\{t_1,\ldots,t_n\}$ is endowed with the sup-norm 
$\|\cdot\|$, also called the {\it Gauss norm}.
We list up some known properties (cf.\ \cite[5.2]{na},
\cite[II.3]{garea}):

\begin{enumerate}
\item The ring $K\{t_1,\ldots,t_n\}$ satisfies the Weierstrass
Preparation Theorem.
\item The ring $K\{t_1,\ldots,t_n\}$ is Noetherian and factorial.
\item Every ideal of $K\{t_1,\ldots,t_n\}$ is closed.
\item For any maximal ideal $\mathfrak{m}$ of $K\{t_1,\ldots,t_n\}$, the
field $K\{t_1,\ldots,t_n\}/\mathfrak{m}$ is a finite extension of $K$.
\end{enumerate}

If, moreover, $K$ is algebraically closed, the set of all maximal
ideals, endowed with the natural topology, coincides with 
$\matheur{D}(0,1^{+})$.
So the situation is analogous to classical algebraic geometry; 
$\matheur{D}(0,1^{+})$ to the affine $n$-space, and $K\{t_1,\ldots,t_n\}$
to the coordinate ring.
In fact, quotients $\matheur{A}=K\{t_1,\ldots,t_n\}/I$ by ideals $I$ 
(called {\it affinoid algebras}), and their associated maximal spectra,
$\mathrm{Spm}(\matheur{A}$), called {\it affinoids}, forms the 
fundamental patches from which the non-archimedean function theory will be
developed.
Although the natural metric topology on $\mathrm{Spm}(\matheur{A})$ is 
terrible for the same reason as above, the situation becomes much
better as far as we deal only with those functions coming from 
$\matheur{A}$.
The trick is that, if the algebra $\matheur{A}$ is integral, then we can
pretend that the space $\mathrm{Spm}(\matheur{A})$ is ``connected''; for
example, the closed disk $\matheur{D}(0,1^{+})$ is ``connected'' in this
sense.
So, basically according to this point of view, we can globalize 
the situation by gluing affinoid patches to obtain a reasonable
theory of analytic functions.
The actual recipe to do it is furnished by the notion of Grothendieck
topology, which censors the plethora of coverings.
We will see this a little more precisely in what follows.

\begin{apprem}
Although the rigid analysis seems to provide the reasonable topological
and analytical framework, one must not expect it to be as nice as the
complex case.  Firstly, in classical complex analysis, when we expand
a holomorphic function centered at a point inside its region of
convergence, the resulting Taylor expansion may have a different region
of convergence than the original one. This is a important method of
analytic continuation.  But in the non-archimedean case,
this never happens (cf.\ \cite[4.4]{pn}).

Secondly, the topological framework
given as above is still not powerful enough to handle paths.
So it is hopeless to treat, for example, fundamental groups and 
monodromy in an intuitive way as in the classical complex case.
This difficulty is, in fact, remedied by Berkovich' theory 
which we shall see later on.
\end{apprem}

\section{Rigid analytic spaces.}

Let $K$ be as in the previous section; we do not assume in general 
that $K$ is algebraically closed.

As mentioned above, an {\it affinoid algebra}
$\matheur{A}$ over $K$ is a quotient $K\{t_1,\ldots,t_n\}/I$ for some
$n$ by an ideal $I$.  Since $I$ is closed, $\matheur{A}$ has a norm
$|\cdot |$ induced from the sup-norm $\|\cdot \|$ of
$K\{t_1,\ldots,t_n\}$ ({\it ``residue norm''}).  With the residue norm,
$\matheur{A}$ is a Banach $K$-algebra.  The norm $|\cdot |$ itself
depends on the presentation $\alpha\colon
K\{t_1,\ldots,t_n\}/I\stackrel{\sim}{\rightarrow}\matheur{A}$, while the
induced topology does not.  So, strictly speaking, we should write it
like $|\cdot |_{\alpha}$.  Clearly, $\matheur{A}$ is Noetherian and, for
any maximal ideal $\mathfrak{m}$ of $\matheur{A}$, the residue field
$\matheur{A}/\mathfrak{m}$ is a finite extension of $K$, and hence, the
valuation $|\cdot |$ of $K$ naturally extends to that of it, denoted
again by $|\cdot |$.

The associated {\it affinoid} $\mathrm{Spm}({\matheur{A}})$ is the set of all
maximal ideals of $\matheur{A}$.
For $x\in\mathrm{Spm}(\matheur{A})$ and $f\in\matheur{A}$ the value of $f$ at
$x$, denoted by $f(x)$, is the class of $f$ in $\matheur{A}/x$.
The set $\mathrm{Spm}({\matheur{A}})$ has the topology generated
by the subsets of form $\{x\in\mathrm{Spm}({\matheur{A}})\ |\ |f(x)|
\leq 1\}$ for $f\in\matheur{A}$.
But, as we pointed out above, this topology is not very interesting since 
it makes $\mathrm{Spm}({\matheur{A}})$ totally disconnected.
So we should specify the reasonable family of ``admissible'' open
sets (and coverings), on which the function theory will be built.

Let $\matheur{A}$ be an affinoid algebra over $K$, and $f_i$ 
($i=0,\ldots,n$) a collection of elements in $\matheur{A}$ which have no
common zeros on $\mathrm{Spm}(\matheur{A})$. 
The subspace
\begin{equation*}
 R=\{x\in\mathrm{Spm}(\matheur{A})\ |\ |f_i(x)|\leq |f_0(x)|,\ i=1,\ldots,n\}
\end{equation*}
can be identified with the affinoid $\mathrm{Spm}(\matheur{A}_R)$, where
\begin{equation*}
 \matheur{A}_R=\matheur{A}\widehat{\otimes}_KK\{t_1,\ldots,t_n\}/
 (f_1-t_1f_0,\ldots,f_n-t_nf_0).
\end{equation*}
A subset of this form is called a {\it rational subdomain}.
The identification comes as follows:
We first note that the morphism of affinoids $\mathrm{Spm}(\matheur{A}_R)
\rightarrow\mathrm{Spm}(\matheur{A})$ induced by
$\matheur{A}\rightarrow\matheur{A}_R$ maps $\mathrm{Spm}(\matheur{A}_R)$ to 
$R$.
Then $\matheur{A}_R$ is the unique solution of the following 
universal property (hence, it can be determined up to canonical isomorphism):
for any morphism of affinoids
$\phi\colon\mathrm{Spm}(\matheur{B})\rightarrow\mathrm{Spm}(\matheur{A})$ such
that $\phi(\mathrm{Spm}(\matheur{B}))\subset R$, there exists a unique 
$K$-homomorphism $\matheur{A}_R\rightarrow\matheur{B}$ such that
the resulting diagram
\vspace{1.5ex}
$$
\setlength{\unitlength}{1pt}
\begin{picture}(100,60)(10,0)
\put(0,45){$\mathrm{Spm}(\matheur{B})$}
\put(40,48){\vector(1,0){20}}
\put(46,53){$\scriptstyle{\phi}$}
\put(66,45){$\mathrm{Spm}(\matheur{A})$}
\put(29,6){$\mathrm{Spm}(\matheur{A}_R)$}
\multiput(33,38)(1.5,-3){6}{\circle*{.3}}\put(42,20){\vector(1,-2){2}}
\put(57,18){\vector(1,2){10}}
\end{picture}\vspace{-1ex}
$$
is commutative.

\begin{applem}\label{lem-rational-intersection}

$\mathrm{(1)}$
If $R$ and $S$ are rational subdomains in $\mathrm{Spm}(\matheur{A})$, then 
so is $R\cap S$.

$\mathrm{(2)}$
If $R_1$ is a rational subdomain in $\mathrm{Spm}(\matheur{A})$, and $R_2$ 
is a rational subdomain in $R_1$, then $R_2$ is a rational subdomain 
in $\mathrm{Spm}(\matheur{A})$.
\end{applem}

\begin{proof}
(1)
If $R=\{|f_i(x)|\leq |f_0(x)|,\ i=1,\ldots,n\}$ and
$S=\{|g_j(x)|\leq |g_0(x)|,\ j=1,\ldots,m\}$, then
one sees easily that 
$R\cap S=\{|f_i(x)g_j(x)|\leq |f_0(x)g_0(x)|,\ i=1,\ldots,n,
\ j=1,\ldots,m\}$.

(2) Exercise.
\end{proof}

\begin{appcor}\label{cor-rational-intersection}
Any subspace of the form
\begin{equation*}
 \{|f_i(x)|\leq 1,|g_j(x)|\geq 1,\ i=1,\ldots,n,\ j=1,\ldots,m\}
\end{equation*}
is a rational subdomain.
\end{appcor}

\begin{proof}
It is the intersection of rational subdomains
$\{|f_i(x)|\leq 1\}$ for $i=1,\ldots,n$ and
$\{1\leq |g_j(x)|\}$ for $j=1,\ldots,m$.
\end{proof}

\begin{appexa}[Rational subdomains of the unit polydisk]
Let us assume that $K$ is algebraically closed.
We consider rational subdomains of the unit polydisk
$\matheur{D}(0,1^{+})=\mathrm{Spm}(K\{t_1,\ldots,t_n\})$.

(1) {\it Closed polydisk}:\  
$\matheur{D}(0,|\pi|^{+})=\{x\in K^n\ |\ |x_i|\leq |\pi_i|,\ 
i=1,\ldots,n\}$,
for $0\neq\pi_i\in K$, $|\pi_i|\leq 1$, is a rational subdomain. 
The corresponding affinoid algebra is 
$$
K\left\{\frac{t_1}{\pi_1},\ldots,\frac{t_n}{\pi_n}\right\}.
$$

(2) {\it Annulus}:\ 
$\matheur{C}(0,|\pi^{(1)}|^{+},|\pi^{(2)}|^{+})=
\{x\in K^n\ |\ |\pi^{(1)}_i|\leq |x_i|\leq
|\pi^{(2)}_i|,\ i=1,\ldots,n\}$, 
for $0\neq\pi^{(j)}_i\in K$, $|\pi^{(j)}_i|\leq 1$, is a rational 
subdomain. The corresponding affinoid algebra is 
$$
K\left\{\frac{\pi^{(1)}_1}{t_1},\ldots,\frac{\pi^{(1)}_n}{t_n},
\frac{t_1}{\pi^{(2)}_1},\ldots,\frac{t_n}{\pi^{(2)}_n}\right\}.
$$
\end{appexa}

Let us next review the definition of Grothendieck topology
(in a narrow sense). Let $X$ be a topological space.
A {\it Grothendieck topology} (G-topology, in short) on $X$ is a pair 
$(\mathscr{T},\mathscr{C}\matheur{ov})$ consisting of
\begin{itemize}
\item a collection $\mathscr{T}$ of open subsets in $X$,
\item an assignment $U\mapsto\mathscr{C}\matheur{ov}(U)$ for any $U\in\mathscr{T}$,
where $\mathscr{C}\matheur{ov}(U)$ is a collection of coverings by elements in
$\mathscr{T}$,
\end{itemize}

\noindent
such that the following conditions are satisfied:

\vspace{1ex}
(1) $\emptyset\in\mathscr{T}$; $U,V\in\mathscr{T}\ \Rightarrow\ 
U\cap V\in\mathscr{T}$.

(2) $U\in\mathscr{T}\ \Rightarrow\ \{U\}\in\mathscr{C}\matheur{ov}(U)$.

(3) $\{U_i\}_{i\in I}\in\mathscr{C}\matheur{ov}(U),\ V\subseteq U,\ V\in
\mathscr{T}\ \Rightarrow\ \{U_i\cap V\}_{i\in I}\in
\mathscr{C}\matheur{ov}(V)$.

(4) $\{U_i\}_{i\in I}\in\mathscr{C}\matheur{ov}(U),\ \{U_{i,j}\}_{j\in J_i}
\in\mathscr{C}\matheur{ov}(U_i)\ \Rightarrow\ \{U_{i,j}\}_{i\in I,j\in J}
\in\mathscr{C}\matheur{ov}(U)$.

\vspace{1ex}
Elements in $\mathscr{T}$ are called {\it admissible open sets}, and elements
in $\mathscr{C}\matheur{ov}(U)$ are called {\it admissible coverings} of $U$.

Let $X=\mathrm{Spm}(\matheur{A})$ be an affinoid over $K$.
We introduce a Grothendieck topology to $X$
by the following recipe:

\begin{itemize}
\item An admissible open set is a rational subdomain.
\item For a rational subdomain, an admissible covering is a finite 
covering consisting of rational subdomains.
\end{itemize}
This is valid, due to Lemma \ref{lem-rational-intersection}.

We define the structure presheaf $\OO_X$, with respect to the Grothendieck
topology, by assigning for each rational subdomain $R$ the corresponding
affinoid algebra $\matheur{A}_R$.
Due to the following weak form of {\it Tate's acyclicity theorem}, the
presheaf $\OO_X$ is, in fact, a sheaf (a sheaf of local rings).

\begin{appthm}
Let $U_1,\ldots,U_m$ be rational subdomains of
$X=\mathrm{Spm}(\matheur{A})$. Set $U=U_1\cup\cdots\cup U_m$.
Then the sequence
$$
0\longrightarrow\OO_X(U)\longrightarrow\prod^m_{i=1}\OO_X(U_i)
\longrightarrow\prod^m_{i,j=1}\OO_X(U_i\cap U_j)
$$
is exact, where the last arrow is the difference of the two possible
restriction morphisms.
\end{appthm}

For the proof, see, for example, \cite[8.2]{na} or \cite[III.2.2]{garea}.

\begin{appdfn}[Rigid Analytic Space]
A {\it rigid analytic space} is a locally ringed space $(X,\OO_X)$, with 
a Grothendieck topology in the above sense, locally isomorphic to an
affinoid;
more precisely, there exists a covering $\{X_i\}_{i\in I}$ 
(possibly infinite) of $X$ by 
admissible open sets such that $(X_i,\OO_X|_{X_i})$ is isomorphic to 
a certain affinoid for each $i$.
\end{appdfn}

Let us see some examples of rigid analytic spaces:

\begin{appexa}[Projective Space]
The projective space $\mathbb{P}^n(K)$ has the natural rigid 
analytic structure.
Here we limit ourselves to demonstrate it only in the case $n=1$, and 
leave the general case to the reader.
Let $(X\colon Y)$ be the homogeneous coordinate in $\mathbb{P}^1(K)$
and set $z=X/Y$.
Set 
$$
\matheur{U}^{+}=\left\{(X\colon Y)\ |\ |X|\leq |Y|\right\}
\ \mathrm{and}\ 
\matheur{U}^{-}=\left\{(X\colon Y)\ |\ |X|\geq |Y|\right\}.
$$
Then $\matheur{U}^{+}$ is isomorphic to the closed disk
$\{z\in K\ |\ |z|\leq 1\}$ (the Southern Hemisphere with the equator)
, and $\matheur{U}^{-}$ to
$\{z\in K\cup\{\infty\}\ |\ |z|\geq 1\}$
(the Northern Hemisphere with the equator).
The intersection $\matheur{U}^{+}\cap\matheur{U}^{-}$ is thus isomorphic to
the circle $\{z\in K\ |\ |z|=1\}$.
Hence $\{\matheur{U}^{+},\matheur{U}^{-}\}$ gives the admissible covering of
$\mathbb{P}^1(K)$ which induces, in the obvious way, the structure
of a rigid analytic space.
\end{appexa}

More generally, an algebraic variety, separated and of finite type
over $K$, carries the canonical structure of a rigid analytic space.

Any open set (not necessarily admissible) $U$ in a rigid analytic 
space $X$ will be again a rigid analytic space ({\it ``open subanalytic
space''}); more precisely, $U$ has the induced analytic structure from 
that of $X$. The admissible open sets of $U$, for example, are those
of $X$ contained in $U$.

\section{Relation with Formal Geometry.}

In this section we need to consider the valuation ring 
$R$ of $K$. Let us denote by $\mathfrak{m}$ the maximal ideal of $R$, and
fix a topological generator $0\neq\pi\in\mathfrak{m}$.
The residue field $R/\mathfrak{m}$ is denoted by $k$.

Let $\matheur{A}$ be an affinoid algebra over $K$.
The {\it spectral semi-norm} on $\matheur{A}$ is the function
$|\cdot |_{\mathrm{sup}}\colon\matheur{A}\rightarrow\R_{\geq 0}$
defined, for any $f\in\matheur{A}$, by
\begin{equation*}
 |f|_{\mathrm{sup}}=\sup_{x\in\mathrm{Spm}(\matheur{A})}|f(x)|.
\end{equation*}
It is known that for any representation 
$\alpha\colon K\{t_1,\ldots,t_n\}/I\stackrel{\sim}{\rightarrow}\matheur{A}$,
we have $|f|_{\mathrm{sup}}\leq |f|_{\alpha}$ for every $f\in\matheur{A}$
(cf.\ \cite[6.2.1]{na}).

\begin{appthm}[Maximal modulus principle]
For an affinoid algebra 
$\matheur{A}$ and an element $f\in\matheur{A}$, there exists a point
$x\in\mathrm{Spm}(\matheur{A})$ such that $|f(x)|=|f|_{\mathrm{sup}}$.
\end{appthm}

We refer \cite[6.2.1]{na} for the proof.

Given an affinoid algebra $\matheur{A}$ over $K$, we introduce the following
notation:
\begin{equation*}
 \begin{array}{lcl}
 \matheur{A}^{\circ}&=&\{f\in\matheur{A}\ |\ |f|_{\mathrm{sup}}\leq 1\},\\
 \matheur{A}^{\circ\circ}&=&\{f\in\matheur{A}\ |\ |f|_{\mathrm{sup}}<1\},\\
 \ovl{\matheur{A}}&=&\matheur{A}^{\circ}/\matheur{A}^{\circ\circ};
 \end{array}
\end{equation*}
$\matheur{A}^{\circ}$ is a $R$-subalgebra of $\matheur{A}$, and 
$\matheur{A}^{\circ\circ}$ is an ideal of it.
For instance, $K^{\circ}=R$, $K^{\circ\circ}=\mathfrak{m}$, and $\ovl{K}=k$.

\begin{appexa}\label{exa-red-tate}
We can immediately calculate:
$$
\begin{array}{lrl}
{\matheur{T}^n_K}^{\circ}&=&R\{t_1,\ldots,t_n\}\\
&\colon=&\left\{
\begin{array}{r|l}
\sum_{\nu_1,\ldots,\nu_n\geq 0}
a_{\nu_1,\ldots,\nu_n}t_1^{\nu_1}\cdots t_n^{\nu_n}&
|a_{\nu_1,\ldots,\nu_n}|\rightarrow 0\ \mathrm{for}\\
\in R[[t_1,\ldots,t_n]]
&\nu_1+\cdots+\nu_n\rightarrow\infty
\end{array}
\right\}\\
\ovl{\matheur{T}^n_K}&=&k[t_1,\ldots,t_n].
\end{array}
$$
Note that ${\matheur{T}^n_K}^{\circ}$ is the $(\pi)$-adic completion of
the polynomial ring $R[t_1,\ldots,t_n]$ over $R$.
It should be noticed that $\ovl{\matheur{T}^n_K}$ is a 
{\it polynomial ring} over $k$; indeed, each element in
$R\{t_1,\ldots,t_n\}$ lies, by the convergence condition, in 
${\matheur{T}^n_K}^{\circ\circ}$ modulo finitely many terms.
\end{appexa}

\begin{apppro}
Let $\matheur{A}$ be an affinoid algebra over $K$. Then,
\begin{enumerate}
 \item $\matheur{A}^{\circ}$ is a model of $\matheur{A}$$;$
       i.e.\ $\matheur{A}^{\circ}\otimes_RK\cong\matheur{A}$,
 \item $\matheur{A}^{\circ}$ is $(\pi)$-adically completed, 
 \item $\matheur{A}^{\circ}$ is topologically of finite
       type over $R$$;$ i.e.\
       $R\{t_1,\ldots,t_n\}/\mathfrak{a}\cong\matheur{A}$ 
       for some ideal $\mathfrak{a}$ in $R\{t_1,\ldots,t_n\}$,
 \item $\matheur{A}^{\circ}$ is flat over $R$,
 \item $\ovl{\matheur{A}}$ is a $k$-algebra of finite type.
\end{enumerate}
\end{apppro}

\begin{proof}
(i) and (ii) are easy to see.
The flatness of $R$ is now equivalent to the lack of $R$-torsion.
So (iv) can be seen immediately.

Before proving (iii) and (v), we need some terminology:
let $\|\cdot\|_{\alpha}$ be a Banach norm of $\matheur{A}$ ($\alpha\colon
\matheur{T}^n_K\rightarrow\matheur{A}$: surjective).
We say an element $f$ of $\matheur{A}$ is {\it power-bounded} if 
$\{\|f^n\|_{\alpha}\ |\ n\in\N\}$ is bounded.
This does not depend on the choice of the presentation $\alpha$ since
the equivalence class of the resulting norms do not change.
An element $f$ of $\matheur{A}$ is said to be {\it topologically nilpotent}
if $\lim f^n=0$; here the limit is taken with respect to the Banach
norm $\|\cdot\|_{\alpha}$, and is not dependent of the choice.

Let us prove (iii).
By \cite[6.2.3]{na}, $\matheur{A}^{\circ}$ is the set of all power-bounded
elements of $\matheur{A}$, and $\matheur{A}^{\circ\circ}$ is the set of all
topologically nilpotent elements.
Hence a choice of a presentation 
$\alpha\colon\matheur{T}^n_K\rightarrow\matheur{A}$ induces surjections
$$
\begin{array}{lcl}
{\matheur{T}^n_K}^{\circ}&\longrightarrow&\matheur{A}^{\circ}\\
{\matheur{T}^n_K}^{\circ\circ}&\longrightarrow&\matheur{A}^{\circ\circ}.
\end{array}
$$
(iii) is due to the surjectivity of the first arrow.
By the surjectivity of these arrows and \ref{exa-red-tate}, we have
(v).
\end{proof}

\begin{appdfn}\label{dfn-admissiblealg}
An $R$-algebra $A$ is said to be {\it admissible} if $A$ is 
$(\pi)$-adically completed, flat over $R$, and topologically 
of finite type over $R$.
\end{appdfn}

\begin{apprem}
By \cite[1.1\ (c)]{farg1rs}, any admissible $R$-algebra $A$ is 
topologically of finite presentation, i.e., there exists a finitely
generated ideal $\mathfrak{a}$ in $R\{t_1,\ldots,t_n\}$ such that
$R\{t_1,\ldots,t_n\}/\mathfrak{a}\cong\matheur{A}$.
\end{apprem}

\begin{appdfn}[Formal model and Analytic reduction]
Let $X=\mathrm{Spm}\matheur{A}$ be an affinoid over $K$.
Then the formal scheme $\mathscr{X}=\Spf\matheur{A}^{\circ}$ over $R$ is called
the {\it formal model} of $X$, and the algebraic scheme 
$\ovl{X}=\Spec\ovl{\matheur{A}}$ over $k$ is called the {\it analytic
reduction} of $X$.
\end{appdfn}

We define the so-called {\it reduction map}
$$
\mathrm{Red}_X\colon
X=\mathrm{Spm}\matheur{A}\longrightarrow\ovl{X}=\Spec\ovl{\matheur{A}}
$$
by $\mathfrak{a}\mapsto\mathfrak{a}\cap\matheur{A}^{\circ}/\mathfrak{a}\cap
\matheur{A}^{\circ\circ}$.

\begin{apppro}
The map $\mathrm{Red}_X$ is continuous with respect to the 
Zariski topology of $\ovl{X}$ and the G-topology of $X$; 
it maps $X$ surjectively onto the set of all
closed points of $\ovl{X}$.
\end{apppro}

\begin{proof}
Let us consider the affine open set 
$U_{\ovl{f}}=\Spec\ovl{\matheur{A}}_{\ovl{f}}$ of $\Spec\ovl{\matheur{A}}$
relevant to $\ovl{f}\in\ovl{\matheur{A}}$.
Take $f\in\matheur{A}^{\circ}$ from the residue class $\ovl{f}$.
The function $f$ takes values less than equal to $1$ on 
$\mathrm{Spm}\matheur{A}$.
So $U_f=\{x\in\mathrm{Spm}\matheur{A}\ |\ |f(x)|=1\}$ is a rational
subdomain whose associated affinoid is $\mathrm{Spm}A\{f,f^{-1}\}$.
Since the condition $|f(x)|=1$ is equivalent to 
$\ovl{f}(\mathrm{Red}_X(x))\neq 0$, we have 
$\mathrm{Red}_X^{-1}(U_{\ovl{f}})=U_f$.
Then the first assertion follows.
The last one is almost trivial.
\end{proof}

\begin{appexa}[Analytic reduction of the unit disk]
Let us consider the analytic reduction of the unit disk 
$\matheur{D}^1=\matheur{D}(0,1^{+})=\mathrm{Spm}K\{T\}$
(here we assume that $K$ is algebraically closed for simplicity).
The analytic reduction is given by the affine line 
$\mathbb{A}^1_k=\Spec k[T]$.
The reduction map $\mathrm{Red}_{\matheur{D}^1}$ maps all ``interior''
points, i.e.\ those points lying in $\matheur{D}(0,1^{-})$ to the origin
in $\mathbb{A}^1_k$, and the other points to the non-zero points.
If we identify $\matheur{D}^1$ and $R$, $\mathbb{A}^1_k$and $k$, then
this is nothing but the mapping obtained by the reduction modulo
$\mathfrak{m}$.
\end{appexa}

It is very important to see how the formal models behave under 
localization of affinoids.
Let $X=\mathrm{Spm}\matheur{A}$ be an affinoid over $K$ and $f_0,\ldots,f_n\in
\matheur{A}$ a sequence of elements which does not have common zeros
on $\mathrm{Spm}\matheur{A}$.
We consider the admissible covering $X=\bigcup_iU_i$ consisting of 
rational subdomains
$$
U_i=\{x\in X \mid |f_j(x)|\leq|f_i(x)|\ \mathrm{for}\ j\neq i\};
$$
by multiplying a suitable power of $\pi$, we may assume that each
$f_i$ belongs to $\matheur{A}^{\circ}$.
For each $i$, 
\begin{equation*}
 U_i=\mathrm{Spm}\matheur{A}\left\{\frac{f_0}{f_i},\ldots,
			     \frac{f_n}{f_i}\right\},
\end{equation*}
and then the corresponding formal model is given by 
\begin{equation*}
 \mathscr{U}_i=\Spf\matheur{A}^{\circ}
 \left\{\frac{f_0}{f_i},\ldots,\frac{f_n}{f_i}\right\}/
 \left((\pi)\mathrm{-torsions}\right).
\end{equation*}
These formal schemes glue and the resulting morphism 
\begin{equation*}
 \bigcup_{0\leq i\leq n}\mathscr{U}_i\longrightarrow\mathscr{X}
  =\Spf\matheur{A}^{\circ}
\end{equation*}
is the formal blow-up along the ideal 
$\mathfrak{a}=(f_0,\ldots,f_n)\subset\matheur{A}^{\circ}$.
We note here that the ideal $\mathfrak{a}$ contains a power of $\pi$, or
what amounts to the same, $\mathfrak{a}$ is an open ideal, for 
$\mathfrak{a}\matheur{A}=\matheur{A}$.

This simple observation indicates:

\vspace{1ex}\noindent
{\bf Slogan (vague)}: {\it Refinements of admissible coverings corresponds
to formal blow-up's along open coherent ideals.}

\vspace{1ex}

\begin{appexa}
Let us assume in this example that $K$ is algebraically closed.
We consider the admissible covering $\matheur{D}^1=U_1\cup U_2$ of the
unit disk $\matheur{D}^1=\mathrm{Spm}K\{T\}$ given by
\begin{equation*}
 \begin{array}{cclcl}
 U_1&=&\{z\in K\ |\ |z|\leq|\pi|\}&=&\mathrm{Spm}K\left\{\frac{T}{\pi}
 \right\},\\
 U_2&=&\{z\in K\ |\ |\pi|\leq |z|\leq 1\}&=&
 \mathrm{Spm}K\left\{T,\frac{\pi}{T}\right\}.
 \end{array}
\end{equation*}
The corresponding formal models are
\begin{equation*}
 \begin{array}{cclcl}
 \mathscr{U}_1&=&\Spf R\left\{\frac{T}{\pi}\right\}
 &=&\Spf R\left\{T,U\right\}/(\pi U-T),\\
 \mathscr{U}_2&=&\Spf R\left\{T,\frac{\pi}{T}\right\}
 &=&\Spf R\left\{T,V\right\}/(TV-\pi)
 \end{array}
\end{equation*}
respectively. We glue them as $U=V^{-1}$ and get 
$\mathscr{U}_1\cup\mathscr{U}_2\rightarrow\mathscr{D}^1=\Spf R\{T\}$,
which is the formal blow-up along the ideal $(\pi,T)$.
Meanwhile, the analytic reduction is obtained by gluing the following 
two affine sets:
\begin{equation*}
 \begin{array}{ccl}
 \ovl{U}_1&=&\Spec k\left[T\right],\\
 \ovl{U}_2&=&\Spec k\left[T,V\right]/(TV),
 \end{array}
\end{equation*}
this is the normal crossing of the affine line with coordinate $T$ and 
the projective line with the inhomogeneous coordinate $V$.
\end{appexa}

We are going to state the Raynaud's theorem, which provides a close
relationship between rigid geometry and formal geometry; we need some
terminology:

\begin{appdfn}[Admissible formal scheme]
An $R$-formal scheme\footnote{We follow 
\cite[\S 10]{ega1} for generalities of formal schemes.}
$\mathscr{X}$ is said to
be {\it admissible} if it is Zariski locally isomorphic to admissible
$\Spf A$ with $A$ admissible in the sense of Definition 
\ref{dfn-admissiblealg}.
\end{appdfn}

\begin{appdfn}[Admissible formal blow-up]
Let $\mathscr{X}$ be an admissible formal scheme.
An {\it admissible formal blow-up} of $\mathscr{X}$ is the formal blow-up
along a coherent open ideal $\mathscr{I}\subset\OO_{\mathscr{X}}$, i.e., 
\begin{equation*}
 \mathscr{X}'=\lim_{\stackrel{\longrightarrow}{\lambda}}\mathrm{Proj}
 \bigoplus^{\infty}_{n=0}\left(\mathscr{I}^n\otimes_{\mathscr{O}_{\mathscr{X}}}
 \OO_{\mathscr{X}}/\pi^{\lambda}\right)\longrightarrow\mathscr{X}.
\end{equation*}
\end{appdfn}

We are going to construct the functor from the category of admissible 
formal schemes to the category of rigid analytic spaces; let us start 
by observing the affine case.
From an affine admissible formal scheme $X=\Spf A$, it is easy to
get the corresponding rigid analytic space: we just set 
\begin{equation*}
 A_{\mathrm{rig}}\colon=A\otimes_RK
\end{equation*}
which is an affinoid algebra over $K$.
If $A$ is presented as $A=R\{t_1,\ldots,t_n\}/\mathfrak{a}$, then
$K\{t_1,\ldots,t_n\}/\mathfrak{a}K\{t_1,\ldots,t_n\}$ gives a presentation
of $A_{\mathrm{rig}}$.

To globalize, we need to see the compatibility with localization.
Let us consider the completed localization $A\{f^{-1}\}$ of $A$ with
respect to $f\in A$.
The corresponding affinoid algebra is
\begin{equation*}
 \begin{array}{ccccl}
 A\{f^{-1}\}\otimes_RK&=&A\{T\}/(1-Tf)\otimes_RK\\
 &=&A_{\mathrm{rig}}\{T\}/(1-Tf)&=&A_{\mathrm{rig}}\{f^{-1}\},
 \end{array}
\end{equation*}
and the latter is the affinoid algebra attached to the
rational subdomain
\begin{equation*}
 \{x\in\mathrm{Spm}A_{\mathrm{rig}}\ |\ |f|\geq 1\};
\end{equation*}
this means that the completed localization of formal schemes just 
corresponds to the localization of rigid spaces with respect to the
G-topology.
So we globalize the recipe as above to obtain the functor
\begin{equation*}
 \mathrm{Rig}\colon\mathscr{X}\mapsto\mathscr{X}_{\mathrm{rig}}
\end{equation*}
from admissible formal $R$-schemes to $K$-rigid analytic spaces.
The rigid analytic space $\mathscr{X}_{\mathrm{rig}}$ is called the
{\it Raynaud generic fiber} of $\mathscr{X}$.

\begin{apppro}\label{pro-blow-up}
The functor $\mathrm{Rig}$ maps admissible formal blow-up's of 
admissible formal schemes to isomorphisms of rigid analytic spaces.
\end{apppro}

\begin{proof}
This has already been essentially shown just before the Slogan.
Indeed, we may limit ourselves to the affine case, for the centers
of blow-up's are coherent.
We have seen in this case that the Raynaud generic fiber of 
admissible formal blow-up's is the decomposition into rational
subdomains, which does not change the analytic structure.
\end{proof}

\vspace{1ex}
Now we state the theorem of Raynaud:
\begin{appthm}[Raynaud\ 1972]
The functor gives the equivalence of the following categories,
\begin{enumerate}
\item the category of quasi-compact admissible formal $R$-schemes,
      localized by admissible formal blow-up's,
\item the category of quasi-compact and quasi-separated rigid analytic
      spaces over $K$.
\end{enumerate}
\end{appthm}

The reader finds a nice comprehensive account of the proof of this
theorem in \cite{farg1rs}.

\section{Topology of Rigid Analytic Space.}

In this section, we assume, for simplicity, that $K$ is algebraically
closed.

\begin{appdfn}
 \begin{enumerate}
  \renewcommand{\labelenumi}{(\arabic{enumi})}
  \item A rigid analytic space $X$ is said to be {\it quasi-compact} if
	every admissible covering of $X$ has a finite admissible 
	refinement.
  \item A rigid analytic space $X$ is said to be {\it connected} if there
	is no admissible covering $\{U_i\}_{i\in I}$ of $X$
	such that 
	\begin{equation*}
	 \bigcup_{i\in I_1}U_i\cap\bigcup_{i\in I_2}U_i=\emptyset\ 
	  \mathrm{and}\ 
	  \bigcup_{i\in I_1}U_i\neq\emptyset\neq\bigcup_{i\in I_2}U_i
	\end{equation*}	
	for some non-empty subsets $I_1,I_2\subseteq I$ with $I=I_1\sqcup I_2$.
 \end{enumerate}
\end{appdfn}

Since admissible coverings of affinoids are fixed as finite, 
every affinoid is quasi-compact.
So, we see that a rigid analytic space is quasi-compact if and only if it is
a finite union of affinoids.
It is easy to see that a rigid analytic space $X$ is connected if 
and only if the ring $\Gamma(X,\OO_X)$ has no other idempotent than 
$0$ and $1$.

\begin{appdfn}
 A morphism $\pi\colon Y\rightarrow X$ of rigid analytic spaces over $K$
 is said to be an {\it analytic covering} if there exists an admissible
 covering\footnote{or else a topological covering, referring to the
 Grothendieck topology; cf. \chapref{chap-ana-p-period},
 \subsecref{sub-top-cov-et-cov}} $\{X_i\}_{i\in I}$ of $X$ such that for
 each $i\in I$ $\pi^{-1}(X_i)$ is isomorphic to the disjoint union of
 copies of $X_i$.  A connect rigid analytic space $X$ is said to be {\it
 simply connected} if there is no other connected analytic covering over
 $X$ than the trivial one $\mathrm{id}\colon X\rightarrow X$ (cf.\
 \cite{anopu}).
\end{appdfn}

\begin{applem}\label{lem-simply-connected}
Let $X$ be a connect rigid analytic space $X$.

$\mathrm{(1)}$
The space $X$ is simply connected if and only if every locally
constant sheaf of sets, with respect to the Grothendieck topology, 
is constant.

$\mathrm{(2)}$
If there exists an admissible covering $\{X_i\}_{i\in\N}$ of $X$
such that $(a)$ every $X_i$ is simply connected, and $(b)$
$X_i\subseteq X_{i+1}$ for all $i$, then $X$ is simply connected.
\end{applem}

\begin{proof}
Both are straightforward.
\end{proof}

\begin{appexa}[Topology of the closed disk]\label{para-closed-disk}
Let $\matheur{D}=\matheur{D}(0,1^{+})$ be the $1$-dimensional closed disk of
radius $1$.
A rational subdomain in $\matheur{D}$ of form
\begin{equation*}
 \{x\in\matheur{D}\ |\ |x-a|\leq\rho\ \mathrm{and}\ |x-a_i|\geq\rho_i\ 
 \mathrm{for}\ i=1,\ldots,s\},
\end{equation*}
where $\rho,\rho_i\in |K^{\ast}|$ and $a,a_i\in\matheur{D}$, is called a
{\it standard domain}.
This is the complement of finite union of open disks.
The following properties are easy to verify:
\begin{itemize}
\item If $S_1$ and $S_2$ are standard domains such that $S_1\cap S_2
\neq\emptyset$, then so are $S_1\cup S_2$ and $S_1\cap S_2$.
\item Every finite union of standard domains is uniquely decomposed 
as a disjoint union of standard domains.
\end{itemize}

What is more interesting is the following proposition:
\end{appexa}

\begin{apppro}\label{pro-standart-domain}
Every rational subdomain in $\matheur{D}$ is a finite union of standard
domains.
$($Hence standard domains generate the Grothendieck topology of 
$\matheur{D}$.$)$
\end{apppro}

\noindent
\begin{proof}[Sketch of Proof] (cf.\ \cite[III.1.18]{sgamc}).\ 
Let $R$ be the rational subdomain given by $|f_i(x)|\leq |f_0(x)|$
for $i=1,\ldots,n$.
Deforming those functions slightly, if necessary, we may assume 
$f_i$ and $f_j$ have no common zero for any $i\neq j$, and thus,
we may concentrate to the single inequality $|f_1(x)|\leq |f_0(x)|$.
By Weierstrass Preparation Theorem we may assume $f_0$ and $f_1$ are
polynomials having all their roots in $\matheur{D}$.
Then the proposition follows from an easy calculation.
\end{proof}

\begin{appcor}\label{cor-projective-line}
The Grothendieck topology of the projective line $\mathbb{P}^{1,\mathrm{an}}_K$
is generated by standard domains.
\end{appcor}

\begin{apppro}\label{pro-standard}
Standard domains are connected.
Moreover, every connected rational subdomain in
$\mathbb{P}^{1,\mathrm{an}}_K$ is simply connected.
\end{apppro}

\begin{proof}
The first assertion is obvious since the corresponding affinoid algebra
\begin{equation*}
 K\{t,t_0,\ldots,t_s\}/((t-a)-t_0\pi,\pi_i-t_i(t-a_i)\ \mathrm{for}\ 
 i=1,\ldots,s),
\end{equation*}
where $a,a_i\in K$ and $\pi,\pi_i\in K$, is integral.
(Note that $K\{t,t_0,\ldots,t_s\}$ is Noetherian and factorial.)
The second one is proved by using \ref{lem-simply-connected} (1).
Let $S$ be a standard domain and $\mathscr{F}$ a constant sheaf of
sets on $S$.
We can take an admissible covering $\{S_i\}$ of $S$ which trivializes 
$\mathscr{F}$.
The point is that any non-disjoint union of standard domains is 
standard, and hence, connected.
So the restriction maps $\mathscr{F}(S_i)\rightarrow
\mathscr{F}(S_i\cap S_j)$ must be bijective as far as $S_i\cap S_j
\neq\emptyset$.
It follows that $\mathscr{F}(S_i\cup S_j)
\rightarrow\mathscr{F}(S_i)$ are bijective, and repeating this argument,
we concludes $\mathscr{F}$ is constant.
\end{proof}

\vspace{1ex}
The following corollary may seen at first:

\begin{appcor}\label{cor-simply-connected}
Every connected open subanalytic space\footnote{We refer \cite[9.3.1]{na} 
for the definition of open subanalytic spaces.} in 
$\mathbb{P}^{1,\mathrm{an}}_K$ is simply connected.
\end{appcor}

\begin{proof}
A connected open subanalytic space $U$ has an admissible covering 
consisting of standard domains.
By \cite[III,2.6]{sgamc} and the easy fact that a non-disjoint union of 
connected sets is again connected, we see that a finite non-disjoint
union of standard domain is a connected rational subdomain.
Hence the corollary follows from Lemma \ref{lem-simply-connected} (2)
and Proposition \ref{pro-standard}.
\end{proof}

\vspace{1ex}
Finally we quote, without proof, a theorem by van der Put 
which indicates a close relation between the topology 
of a rigid analytic space and that of its analytic reduction:

\begin{appthm}\label{thm-vanderput}
Let $X$ be a quasi-compact rigid analytic space over $K$ which has 
the irreducible and smooth analytic reduction. Then $X$ is simply
connected.
\end{appthm}

For the proof, see \cite{anopu}.
\begin{apprem}
By the corollary, we know, for example, that the space
$\mathbb{G}^{\mathrm{an}}_K=K^{\times}$ is simply connected.  But, on
the other hand, it is easy to see that the map $K^{\times} \rightarrow
K^{\times}$ by $x\mapsto x^n$ is an analytic morphism.  As far as the
characteristic of $K$ does not divide $n$, we are tempted to think of
this morphisms also as an analytic covering.  But it is not.  Certainly,
this morphism induces isomorphisms between stalks of the structure
sheaf, but there is no {\it admissible} covering which trivializes this
morphism! 
This kind of morphisms is said to be {\it \'etale}; analytic coverings
are \'etale, but not vice versa.  Simply connected analytic spaces
may have many \'etale coverings.
\end{apprem}

\section{Berkovich' approach to non-archimedean analysis}

Berkovich' viewpoint provides an innovative approach to topological
problems of $p$-adic geometry and analysis, which is much closer to the
familiar intuition derived from the complex context
\cite{staagonf}\footnote{Already in \cite{coas} the {\it esquisse} of
Berkovich' idea can be seen.}.  His analysis, in fact, generalizes, or
properly speaking, ``completes'' the rigid analysis.  This situation is,
according to what he says in the introduction of his book, somehow
analogous to that of $\R$ completing $\Q$.  The standard archimedean
metric space $\Q$ is totally disconnected.  The rigid analytic viewpoint
corresponds, to some extent, to regarding every rational interval
$\{r\in\Q\ |\ a\leq r\leq b\}$ as ``connected'', and considering only those
functions which come from analytic functions on the corresponding real
interval.  Berkovich' non-archimedean analysis can be compared, in this
context, to analysis on $\R$ itself. Namely, his presentation of
analytic spaces complements those of rigid analysis by putting some
``generic points'', and consequently, simplifies their topology.

Let us be a little more precise about his theory.
Let $\matheur{A}$ be a Banach ring, i.e., a ring together with a norm 
$\|\cdot\|$ with respect to which $\matheur{A}$ is complete.

\begin{appdfn}\label{dfn-seminorm}

(1) A semi-norm $|\cdot |$ on $\matheur{A}$ is said to be {\it bounded} 
if there exists $C>0$ such that $|f|\leq C\|f\|$ for any $f\in\matheur{A}$.

(2) A semi-norm $|\cdot |$ on $\matheur{A}$ is said to be 
{\it multiplicative} if $|fg|=|f||g|$ for all $f,g\in\matheur{A}$.
\end{appdfn}

Replacing the maximal spectra $\Spm(\matheur{A})$ of rigid geometry, we have:
\begin{appdfn}\label{dfn-Berkovich-spectrum}
The Berkovich' {\it spectrum} $\matheur{M}(\matheur{A})$ for a Banach ring
$\matheur{A}$ is the set of all bounded multiplicative semi-norms on
$\matheur{A}$, endowed with the weakest topology so that the real
valued functions of form $|\cdot |\mapsto |f|$ for $f\in\matheur{A}$
are continuous.
\end{appdfn}

In case $\matheur{A}$ is an affinoid algebra (a {\it strict} affinoid
algebra, in Berkovich' term) over $K$, the affinoid space 
$\mathrm{Spm}(\matheur{A})$ is naturally viewed as a subspace, with 
the relative topology, of the Berkovich' spectrum $\matheur{M}(\matheur{A})$.
This is done by identifying an element $\mathfrak{m}$ in
$\mathrm{Spm}(\matheur{A})$ with the unique semi-norm $|\cdot |$ on 
$\matheur{A}$ such that $\mathfrak{m}=\Ker |\cdot |$.
The space $\mathrm{Spm}(\matheur{A})$ is identified, in terms of this 
correspondence, with the subspace of $\matheur{M}(\matheur{A})$ consisting of
semi-norms $|\cdot |$ such that $\dim_K\matheur{A}/\Ker |\cdot |$ is 
finite (the ``classical points'').

\begin{appexa}
Let us consider the usual absolute-value norm $|\cdot|_{\infty}$
on $\Z$. The spectrum
$\matheur{M}(\Z)$ consists of the following points:
\begin{enumerate}
\item $|\cdot|_{\infty,\epsilon}\colon=|\cdot|_{\infty}^{\epsilon}$
($0<\epsilon\leq 1$).
\item The $p$-adic norm $|\cdot|_{p,\epsilon}\colon$ with 
$|p|_{p,\epsilon}=\epsilon$ ($0<\epsilon\leq 1$).
\item The semi-norm $|\cdot|_p$ induced from the trivial norm on $\Z/p\Z$.
\item The trivial norm $|\cdot|_0$.
\end{enumerate}
As a topological space, $\matheur{M}(\Z)$ is a tree with end points
$|\cdot|_p$ and $|\cdot|_{\infty}$. Each of these points is connected
by a single edge with $|\cdot|_0$.
\end{appexa}

On the Berkovich' spectrum, the {\it affinoid subdomains} (incl.\ {\it 
rational subdomains}) are defined in the similar way to rigid 
analysis. 
Thus one can define the structure sheaves on spectra, which 
one calls affinoid spaces.
The local data consisting of Berkovich' affinoids and the structure
sheaves can glue to locally ringed spaces; thus we obtain {\it analytic
spaces} in Berkovich' sense.

\begin{appexa}
(Cf.\ \cite[1.5]{staagonf}.)
Let $K$ be a non archimedean complete field. An affine space over $K$ is
defined by
\begin{equation*}
 \mathbb{A}^n_K=\{\textrm{multiplicative seminorm}\ |\cdot|\ \textrm{on}\ 
 K[t_1,\ldots,t_n]\ |\ |\cdot|_{|_K}\ \textrm{is bounded}\}
\end{equation*}
together with the weakest topology so that $|\cdot|\mapsto |f|$
($f\in K[t_1,\ldots,t_n]$) is continuous.
This has a structure of analytic space so that the analytic functions
over it is characterized by limits of rational functions.
When $K=\C$ the space $\mathbb{A}^n_K$ is nothing but the usual affine
space in the complex analytic sense.
If $K$ is a non-archimedean field, then $\mathbb{A}^n_K$ is a 
union of balls; i.e., the spectrum of Tate algebra.
\end{appexa}

Here are some general properties of Berkovich' analytic spaces:
\begin{enumerate}
\item Every connected $K$-analytic space is arcwise connected.
\item If $X=\matheur{M}(\matheur{A})$ with $\matheur{A}$ an affinoid algebra over
$K$, then $\textrm{Krull-dim}\matheur{A}=\mathrm{dim}|X|$, where $|X|$
is the underlying topological space of $X$.
\item If a $K$-analytic space $X$ is smooth, then it is locally
contractible.
\end{enumerate}

Relation between Berkovich' analytic spaces and rigid analytic spaces is
as follows.
If a Berkovich' $K$-analytic space $X$ is formed only by spectra of 
(strict) affinoid algebras, we say $X$ is strict.
There exists a functor from the category of separated strict 
$K$-analytic spaces to that of separated rigid analytic spaces over $K$
by 
\begin{equation*}
 X\mapsto X_0=\{x\in X\ |\ [K(x)\colon K]<\infty\},
\end{equation*}
where $K(x)$ is the residue field at $x$.

\begin{appthm}
This functor is fully-faithful, and preserves fiber products.
\end{appthm}

\chapter{Overview of theory of $p$-adic uniformization.}
\renewcommand{\thechapter}{\Alph{chapter}}
\addtocontents{toc}{\protect\par\vskip2mm\hskip80mm by F. Kato \par\vskip5mm}

\vspace{1ex}
\begin{center}
By {\sc F.\ Kato}
\end{center}

\vspace{2ex}
In this appendix we always denote by $K$, $R$, and $\pi$ a complete 
discrete valuation field, its valuation ring, and a prime element of
$R$, respectively. 
We assume that the residue field $k=R/\pi R$ is finite and consists of 
$q$ elements.
The Mumford-Kurihara-Mustafin uniformization is a procedure to
construct nice analytic (and in many cases algebraic) varieties by
taking discrete Schottky-type quotients of a certain $p$-adic analogue
of symmetric space, so-called {\it Drinfeld symmetric space} $\Omega$, 
or its variants.

\section{Bruhat-Tits building.}

First we survey the construction of a certain simplicial complex, called 
Bruhat-Tits building (attached to $\mathrm{PGL}(n+1,K)$), which will be 
closely related with the Drinfeld symmetric space.
Let $n$ be a positive integer and $V$ an $n+1$ dimensional vector space
over $K$.
A {\it lattice} in $V$ is a finitely generated 
$R$-submodule of $V$ which spans $V$ over $K$. 
Every lattice is therefore a free $R$-module of rank $n+1$.
Let $\til{\Delta}_0$ be the set of all lattices in $V$. 
We say that two lattices $M_1$ and $M_2$ are {\it similar} if there
exists $\lambda\in K^{\times}$ such that $M_1=\lambda M_2$.
The similarity is obviously an equivalence relation.
We denote by
$\Delta_0$ the set of all similarity classes of lattices in $V$.

\begin{appdfn}
The {\it Bruhat-Tits building} (attached to $\mathrm{PGL}(n+1,K)$) is
the finite dimensional simplicial complex $\Delta$ with the vertex set
$\Delta_0$ defined as follows: A finite subset $\{\Lambda_0,\ldots,
\Lambda_l\}$ of $\Delta_0$ forms an $l$-simplex if and only if, after 
permuting indices if necessary, one can choose $M_i\in\Lambda_i$ for 
$0\leq i\leq n+1$ such that
$$
M_0\supsetneq M_1\supsetneq\cdots\supsetneq M_l\supsetneq\pi M_0.
$$
\end{appdfn}

The role of $M_0$ is by no means important; for example, one can shift
the indices like
$M_1\supsetneq\cdots\supsetneq M_l\supsetneq\pi M_0\supsetneq\pi M_1$.

To understand the structure of $\Delta$, let us fix one vertex
$\Lambda=[M]$. 
Suppose $\{\Lambda_0,\ldots,\Lambda_l\}$ is an $l$-simplex having
$\Lambda$ as a vertex; we may assume $\Lambda=\Lambda_0$, and can take
$M_i\in\Lambda_i$ for $1\leq i\leq l$ such that 
$M\supsetneq M_1\supsetneq\cdots\supsetneq M_l\supsetneq\pi M$.
Set $\ovl{M}_i\colon=M_i/\pi M$.
Then we get a flag 
$$
\ovl{M}\supsetneq\ovl{M}_1\supsetneq\cdots\supsetneq\ovl{M}_l\supsetneq 0
$$
of length $l+1$ consisting of subspaces of the $n+1$-dimensional
$k$-vector space $\ovl{M}=M/\pi M$.
Moreover, by the elementary algebra one finds that this yields a
one-to-one correspondence between $l$-simplices containing $\Lambda$ and 
length $l+1$ frags in $\ovl{M}=M/\pi M$.
By this quite easy observation one knows
\begin{enumerate}
\item $\Delta$ is an $n$-dimensional locally finite simplicial complex,
\item every simplex in $\Delta$ is a face of a chamber, a simplex having 
      maximal ($=n$) dimension.
\end{enumerate}
Furthermore, by the construction,
\begin{enumerate}
\setcounter{enumi}{2}
\item the group $\Aut_K(V)/K^{\times}(\cong\mathrm{PGL}(n+1,K))$
      naturally acts on $\Delta$ simplicially.
\end{enumerate}

\begin{appexa}[$n=1$ case]
In $n=1$ case, the vertices adjacent to a fix vertex $\Lambda$ are in 
one-to-one correspondence with the set of all $k$-rational points of
the projective line $\mathbb{P}^1_k$; in particular, each vertex has 
exactly $q+1$ adjacent vertices.
Moreover, it can be seen that the simplicial complex $\Delta=\Delta^1$
is a {\it tree}.
The following picture shows how it looks like in case $q=2$:
$$
\setlength{\unitlength}{1pt}
\begin{picture}(200,100)(0,0)
\put(100,50){\circle*{3}}
\put(100,50){\line(0,1){20}}\put(100,70){\circle*{3}}
\put(100,50){\line(2,-1){20}}\put(120,40){\circle*{3}}
\put(100,50){\line(-2,-1){20}}\put(80,40){\circle*{3}}
\put(100,70){\line(2,1){12}}\put(112,76){\circle*{3}}
\put(100,70){\line(-2,1){12}}\put(88,76){\circle*{3}}
\put(120,40){\line(0,-1){12}}\put(120,28){\circle*{3}}
\put(120,40){\line(2,1){12}}\put(132,46){\circle*{3}}
\put(80,40){\line(0,-1){12}}\put(80,28){\circle*{3}}
\put(80,40){\line(-2,1){12}}\put(68,46){\circle*{3}}
\multiput(112,76)(0,4){3}{\circle*{.3}}
\multiput(112,76)(4,-2){3}{\circle*{.3}}
\multiput(88,76)(0,4){3}{\circle*{.3}}
\multiput(88,76)(-4,-2){3}{\circle*{.3}}
\multiput(120,28)(4,-2){3}{\circle*{.3}}
\multiput(120,28)(-4,-2){3}{\circle*{.3}}
\multiput(132,46)(0,4){3}{\circle*{.3}}
\multiput(132,46)(4,-2){3}{\circle*{.3}}
\multiput(80,28)(4,-2){3}{\circle*{.3}}
\multiput(80,28)(-4,-2){3}{\circle*{.3}}
\multiput(68,46)(0,4){3}{\circle*{.3}}
\multiput(68,46)(-4,-2){3}{\circle*{.3}}
\end{picture}\vspace{-1ex}
$$
\end{appexa}

In general, the geometrical realization $|\Delta|$ of $\Delta$ is known
to be a contractible topological space.

\begin{appexa}[$n=2$ case]\label{appexa-B}
Already in case $n=2$, the structure of $\Delta=\Delta^2$ is fairly
complicated.
For a fixed vertex $\Lambda$, the set of all vertices adjacent of 
$\Lambda$ is divided in two; one for those corresponding to $q^2+q+1$
$k$-rational points in $\mathbb{P}^2_k$ and the other for those
corresponding to $q^2+q+1$ $k$-rational lines in $\mathbb{P}^2_k$.
So it is convenient to consider the surface $B$ which is obtained by
blowing-up of $\mathbb{P}^2_k$ at every $k$-rational points.
Two vertices adjacent to $\Lambda$ forms with $\Lambda$ a $2$-simplex
if and only if the corresponding lines in $B$ intersect.
\end{appexa}

\section{Drinfeld symmetric space.} 

To each lattice $M$ in $V$ we associate an $R$-scheme
$$
\mathbb{P}(M)\colon=\mathrm{Proj}(\mathrm{Sym}_RM).
$$
Obviously the scheme $\mathbb{P}(M)$ depends only on the similarity
class; so we should write as $\mathbb{P}(\Lambda)$, where $\Lambda=[M]$.
The scheme $\mathbb{P}(M)$ is isomorphic to $\mathbb{P}^n_R$; but the
isomorphism cannot be taken canonically, while there exists the
{\it canonical} isomorphism $\mathbb{V}\cong\mathbb{P}(M)\otimes_RK$ 
due to the canonical isomorphism $M\otimes_RK\cong V$.
So we get a bijection:
$$
\Delta_0\stackrel{\sim}{\longrightarrow}\left\{
\begin{array}{c|l}
(\mathbb{P},\phi) &\mathbb{P}\stackrel{\sim}{\rightarrow}\mathbb{P}^n_R,\\
&\phi\colon\mathbb{P}\otimes_RK\stackrel{\sim}{\rightarrow}\mathbb{P}(V)
\end{array}
\right\}/\sim,
$$
where $(\mathbb{P}_1,\phi_1)\sim(\mathbb{P}_2,\phi_2)$ if and only if
there exists an $R$-isomorphism
$\Phi\colon\mathbb{P}_1\stackrel{\sim}{\rightarrow}\mathbb{P}_2$
which makes the diagram 
$$
\begin{array}{ccc}
\mathbb{P}_1\otimes_RK&\stackrel{\Phi\otimes_R K}
{\longrightarrow}&\mathbb{P}_2\otimes_RK\\
\llap{$\phi_1$}\downarrow&&\downarrow\rlap{$\phi_2$}\\
\mathbb{P}(V)&=&\mathbb{P}(V)
\end{array}
$$
commutative.

Next we consider two vertices $\Lambda_1$ and $\Lambda_2$ which are 
adjacent with each other.
The canonical isomorphism between the generic fibers induce a natural
birational map between $\mathbb{P}(\Lambda_1)$ and 
$\mathbb{P}(\Lambda_2)$.
More explicitly, taking representatives $M_i\in\Lambda_i$ ($i=1,2$) so
that $M_1\supsetneq M_2\supsetneq\pi M_1$, we have the birational map
from $\mathbb{P}(\Lambda_1)$ to $\mathbb{P}(\Lambda_2)$ induced by
$M_2\hookrightarrow M_1$.
We define the {\it join}
$\mathbb{P}(\Lambda_1)\vee\mathbb{P}(\Lambda_2)$ to be the closure of
the graph in $\mathbb{P}(\Lambda_1)\times_R\mathbb{P}(\Lambda_2)$ of
this rational map.
It is easily seen that the join 
$\mathbb{P}(\Lambda_1)\vee\mathbb{P}(\Lambda_2)$ is $R$-isomorphic to
the blow-up of $\mathbb{P}(\Lambda_1)$ at the closed subscheme 
$\mathrm{Proj}(\mathrm{Sym}_kM_1/M_2)$.
In particular, the closed fiber of 
$\mathbb{P}(\Lambda_1)\vee\mathbb{P}(\Lambda_2)$ consists of two
rational varieties intersecting transversally.

\begin{applem}
The operation $\vee$ is associative and commutative.
\end{applem}

The proof is straightforward, and left to the reader. 
By the lemma, we can extend the construction to any finite subcomplex of
$\Delta$.
More precisely, for a finite subcomplex $\matheur{S}$, we take the joins 
of all vertices in $\matheur{S}$; the resulting $R$-scheme is denoted by 
$\mathbb{P}(\matheur{S})$.
The scheme $\mathbb{P}(\matheur{S})$ is a regular projective scheme over
$R$ having the following properties:
\begin{enumerate}
\item The generic fiber is canonically isomorphic to $\mathbb{P}(V)$.
\item The closed fiber is a reduced normal crossing of non-singular
      rational varieties of which the dual graph coincides with
      $\matheur{S}$.
\end{enumerate}

We want to do this for a general convex subcomplex of $\Delta$; this cannot 
be done in terms of schemes because we have to consider a limit of
blow-up's of $\mathbb{P}^n_R$ centered in the closed fiber.
But we can do this in terms of formal schemes, or rigid analytic spaces.
Let $\Delta_{\ast}$ be a convex subcomplex of $\Delta$.
For any finite subcomplex $\matheur{S}$ of $\Delta_{\ast}$ we define 
$\widehat{\Omega}(\matheur{S})$ to be the completion of
$\mathbb{P}(\matheur{S})$ along the closed fiber.
Consider the maximal Zariski open subset $\widehat{\Omega}(\matheur{S})'$
of $\widehat{\Omega}(\matheur{S})$ such that, for any finite subcomplex
$\matheur{T}$ of $\Delta_{\ast}$ containing $\matheur{S}$, the induced
morphism $\rho^{\matheur{T}}_{\matheur{S}}\colon\widehat{\Omega}(\matheur{T})
\rightarrow\widehat{\Omega}(\matheur{S})$ gives an isomorphism restricted to 
the pull-back of $\widehat{\Omega}(\matheur{S})'$.

\begin{appdfn}
For a convex subcomplex $\Delta_{\ast}$ of $\Delta$ we define
$$
\widehat{\Omega}(\Delta_{\ast})\colon=\bigcup_{\matheur{S}\colon
\textrm{\scriptsize{finite}}}
\widehat{\Omega}(\matheur{S})',
$$
where $\matheur{S}$ runs through all finite subcomplex of
$\Delta_{\ast}$.
The corresponding rigid analytic space (Raynaud generic fiber) is
denoted by $\Omega(\Delta_{\ast})$.
In case $\Delta_{\ast}=\Delta$ we simply write $\widehat{\Omega}=
\widehat{\Omega}(\Delta)$ and $\Omega=\Omega(\Delta)$.
The rigid analytic space $\Omega$ is called the {\it Drinfeld symmetric
space}.
\end{appdfn}

The formal scheme $\widehat{\Omega}(\Delta_{\ast})$ is formally locally of
finite type, of which the closed fiber is the normal crossings of
non-singular rational varieties having the dual graph isomorphic to 
$\Delta_{\ast}$.
By definition, the group $\mathrm{PGL}(n+1,K)$ naturally acts on 
$\widehat{\Omega}$, and hence $\Omega$.

\begin{appexa}
Let us consider $\widehat{\Omega}$ in $n=2$. 
The closed fiber $\widehat{\Omega}_0$ consists of countably many 
components which are all isomorphic to the surface $B$ considered 
in \ref{appexa-B}. 
In each component every $k$-rational line is a double curve in 
$\widehat{\Omega}_0$, and every $k$-rational point is a triple point.
\end{appexa}

{\it Structure of $\Omega$.}\ 
First we note that, as a point set, we have an equality
$$
\Omega=\mathbb{P}^{n,\mathrm{an}}_K\setminus\bigcup_{\matheur{H}\in
\mathcal{H}}\matheur{H},
$$
where $\mathcal{H}$ is the set of all $K$-rational hyperplanes.

Here we do not attempt to give a proof of this fact, but rather give a 
sketchy explanation in order to get a feeling.
We limit ourselves to $n=1$ case just for simplicity.
The formal scheme $\widehat{\Omega}$ is, roughly speaking, the {\it limit}
of successive blow-up's 
$$
\cdots\longrightarrow U_i\stackrel{\rho_i}{\longrightarrow}U_{i-1}
\longrightarrow\cdots\longrightarrow U_0
$$
with $U_0=\mathbb{P}^1_R$ and $\rho_i$ the blow-up centered at all the 
$k$-rational points in the closed fiber of $U_{i-1}$.
Let $z$ be a point in the generic fiber $\mathbb{P}^1_K$ of each $U_i$.
We consider the closure $\ovl{z}$ of $z$ in each $U_i$.
Since $U_i$ is a regular scheme, the section $\ovl{z}$ intersects the
closed fiber at a smooth point transversally. 
If $z$ is {\it not} a $K$-rational point, then for sufficiently large 
$N$, $\ovl{z}$ intersects the closed fiber of $U_i$ at a non-$k$-rational
smooth point for $i>N$.
But, if $z$ is $K$-rational, $\ovl{z}$ intersects the closed
fiber always at a $k$-rational point which is in the next blow-up center.
So if we consider the limit, such formal sections determined by $K$-rational 
points does not exist any more, while those by non-$K$-rational points 
certainly exist.
This explains the above equality in $n=1$.
For higher $n$ we can apply the similar argument.

Next we discuss the analytic structure of $\Omega$; viz.\ we show how
to introduce an admissible covering to $\Omega$; since we have
constructed the rigid analytic space $\Omega$ from the formal model, we 
must have one such covering such that the associated formal model 
recovers $\widehat{\Omega}$.
This can be actually described in an elegent way in terms of norms: 
A non-negative real valued function $\alpha\colon V\rightarrow
\mathbb{R}_{\geq 0}$ is said to be a {\it norm} over $K$ if it
satisfies:
\begin{enumerate}
\item $\alpha(x)>0$ whenever $x\neq 0$.
\item For $a\in K$, $\alpha(ax)=|a|\alpha(x)$.
\item $\alpha(x+y)\leq\mathrm{max}\{\alpha(x),\alpha(y)\}$.
\end{enumerate}
Two norms $\alpha$ and $\alpha'$ over $K$ is said to be {\it similar} if 
$\alpha'=\lambda\alpha$ for some $\lambda\in\mathbb{R}_{>0}$.
The set of all similarity classes of norms on $V$ over $K$ forms in an
obvious way a topological space. 
The following fact is well-known (\cite{tsopn}):

\begin{apppro}
The topological space of similarlity classes of norms on $V$ over $K$ is
$\mathrm{PGL}(n+1,K)$-equivaliantly homeomorphic to the topological 
realization $|\Delta|$ of the Bruhat-Tits building $\Delta$.
\end{apppro}

The homeomorphism is given as follows: To any lattice $M$ we need to 
associate a norm $\alpha_M$.
For a non-zero $x\in V$, the set $\{a\in K\ |\ ax\in M\}$ is a
fractional ideal of $R$, and hence is of form $(\pi^m)$.
We then define $\alpha_M(x)=q^m$ for non-zero $x\in V$, and
$\alpha_M(0)=0$.
It is easily verified that $\alpha_M$ is a norm.
The similarity class of $\alpha_M$ depends only on the similarity class
of $M$.
Thus we get a mapping from $\Delta_0$ to the set of similarity classes
of norms.
Let $\Lambda_0$ and $\Lambda_1$ be vertices adjacent with each other.
We can choose a basis $e_0,\ldots,e_n$ of $V$ such that 
$M_0=\bigoplus^n_{i=0}Re_i$ and $M_1=\bigoplus^{n-1}_{i=0}Re_i\oplus
R\pi e_n$ give representatives of $\Lambda_0$ and $\Lambda_1$,
respectively.
Then the norms $\alpha_{M_0}$ and $\alpha_{M_1}$ are given by
\begin{eqnarray*}
\alpha_{M_0}(\sum^n_{i=0}a_ie_i)&=&\mathrm{max}\{|a_0|,\ldots,|a_n|\},\\
\alpha_{M_1}(\sum^n_{i=0}a_ie_i)&=&\mathrm{max}\{|a_0|,\ldots,|a_{n-1}|,
|\pi|^{-1}|a_n|\},
\end{eqnarray*}
respectively.
For $0<t<1$, the class of the norm
$$
\alpha_t(\sum^n_{i=0}a_ie_i)=\mathrm{max}\{|a_0|,\ldots,|a_{n-1}|,
|\pi|^{-t}|a_n|\}
$$
corresponds to the point $t\Lambda_0+(1-t)\Lambda_1$ on the edge
connecting $\Lambda_0$ and $\Lambda_1$.

Now, to any point $z=(z_0\colon\cdots\colon z_n)\in\mathbb{P}(V)$
(subject to some $K$-basis $\{e_0,\ldots,e_n\}$of $V$) 
we associate, up to invertible factor, a non-negative real
valued function
$$
\alpha_z(\sum^n_{i=0}a_ie_i)=\biggm|
\sum_{i=0}^n a_iz_i\biggm|
$$
on $V$.
Note that 
{\it this determines a class of norm if and only if $z$ is not lying in any 
$K$-rational hyperplane.}
Thus we get a mapping 
$$
\rho\colon\Omega\longrightarrow |\Delta|.
$$
The admissible covering in question is actually given by the collection
of subsets of $\Omega$ consisting of pull-backs of each simplices of
$|\Delta|$ by $\rho$.
For example, notation being as above, 
$$
\rho^{-1}(\Lambda_0)=\{|z_0|=\cdots=|z_n|=1\}\setminus
\bigcup_{H}\{(z_0\colon\cdots\colon z_n)\ \mathrm{mod}\ \pi\in H\},
$$
where $H$ runs through all hyperplanes in $\mathbb{P}^n_k$, and 
$$
\rho^{-1}(t\Lambda_0+(1-t)\Lambda_1)=\{|z_0|=\cdots=|z_{n-1}|=1,
|z_n|=|\pi|^{-t}\}
$$
for $0<t<1$.

\section{Uniformization.}

Recall that, once we fix a $K$-basis of $V$, the group 
$\mathrm{PGL}(n+1,K)$ canonically acts on the Bruhat-Tits building 
$\Delta$ and the Drinfeld symmetric space $\Omega$.
We shall consider discrete groups in $\mathrm{PGL}(n+1,K)$ and the
quotients of $\Delta$ and $\Omega$ by them.
A subgroup $\Gamma\in\mathrm{PGL}(n+1,K)$ is said to be {\it
hyperbolic} if it acts on $\Delta$ discretely and freely.

To a hyperbolic subgroup $\Gamma$ we associate a convex subcomplex 
$\Delta_{\Gamma}$ of $\Delta$ by the following recipe:
In general, a convex subcomplex generated by vertices of form
$$
\left[\sum^n_{i=0}R\cdot\pi^{\sigma_i}\matheur{Y}_i\right],
\qquad
\sigma_i\in\Z,\qquad 0\leq i\leq n,
$$
subject to a basis $\matheur{Y}_0,\ldots,\matheur{Y}_n$ of $V$, is said
to be an {\it apartment}.
The geometrical realization of each apartment is a triangulation of 
$\mathbb{R}^n$ of $\mathrm{A}_n$-type, and $\Delta$ is the union of all
apartments.
There exists an obvious bijection between the set of all apartments and
the set of all $K$-split maximal tori in $\mathrm{PGL}(n+1,K)$.
For any $K$-split maximal tori $T$ the untersection $\Gamma\cap T$ is a
commutative discrete group of rank at most $n$.
If this rank is $n$ we call the apartment corresponding to $T$ a 
{\it $\Gamma$-apartment}, and we define $\Delta_{\Gamma}$ to be the
union of all $\Gamma$-apartments.
Clearly the group $\Gamma$ acts on $\Delta_{\Gamma}$; if the quotient
$\Delta_{\Gamma}/\Gamma$ is a finite simplicial complex, we say 
$\Gamma$ is {\it normal hyperbolic}.

Now let us be given a normal hyperbolic subgroup $\Gamma$.
For each $\gamma\in\Gamma$ and each finite subcomplex $\matheur{S}$ of
$\Delta_{\Gamma}$, $\gamma\matheur{S}$ is again a finite subcomplex of
$\Delta_{\Gamma}$, and then, $\gamma$ induces an $R$-isomorphism 
$\mathbb{P}(\matheur{S})\cong\mathbb{P}(\gamma\matheur{S})$ which
induces on generic fibers ($=\mathbb{P}^n_K$) exactly the linear 
automorphism $\gamma$.
Considering all $\matheur{S}$, gluing, and taking Raynaud generic fiber
give a rigid analytic automorphism of $\Omega(\Delta_{\Gamma})$ which 
we denote again by $\gamma$.
Thus we get an homomorphism 
$$
\Gamma\longrightarrow\mathrm{Aut}_{K\textrm{\scriptsize{-rig.}}}
\Omega(\Delta_{\Gamma}).
$$

\begin{appthm}[Mustafin\ \cite{nu}, Kurihara\ \cite{copubathp}]
Let $\Gamma$ be a normal hyperbolic subgroup of $\mathrm{PGL}(n+1,K)$
and $\Delta_{\Gamma}$ as above.
Then there exists, unique up to isomorphisms,
a rigid analytic space $X_{\Gamma}$, smooth, proper and of finite type
over $K$, and a topological covering map $\matheur{p}\colon
\Omega(\Delta_{\Gamma})\rightarrow X_{\Gamma}$ such that
\begin{enumerate}
\item $\matheur{p}\circ\gamma=\matheur{p}$ for any $\gamma\in\Gamma$,
\item for $x,y\in\Omega(\Delta_{\Gamma})$, $\matheur{p}(x)=\matheur{p}
(y)$ if and only if $x=\gamma(y)$ for some $\gamma\in\Gamma$.
\end{enumerate}
\end{appthm}

In the references, this theorem has been stated in terms of formal
schemes; viz.\ they took quotients of
$\widehat{\Omega}(\Delta_{\Gamma})$ to get a formal scheme 
$\widehat{X}_{\Gamma}$ flat and formally of finite type over $R$.
We have an admissible covering of $X_{\Gamma}$ induced from that of 
$\Omega(\Delta_{\Gamma})$; in terms of formal scheme, the closed fiber
of $\widehat{X}_{\Gamma}$ is a divisor of normal crossings with rational
components having the dual group isomorphic to $\Delta_{\Gamma}/\Gamma$.

\vspace{1ex}
{\it Algebraicity.}\ In $n=1$ case the rigid analytic curve $X_{\Gamma}$
is algebraizable for any $\Gamma$ (\cite{aacodcoclr}).
In higher dimensions, we know the following sufficient condition:
Let $\Gamma$ be a torsion-free {\it uniform lattice} in 
$\mathrm{PGL}(n+1,K)$, viz.\ a
discrete co-compact subgroup of finite co-volume without elements of 
finite order.
Then, since it is finite co-volume, we have $\Delta_{\Gamma}=\Delta$.
Hence in this case, we are going to take a quotient of the Drinfeld
symmetric space $\Omega$.

\begin{appthm}[loc.\ cit.]
If $\Gamma$ is a torsion-free uniform lattice, then the quotient rigid
analytic variety $X_{\Gamma}$ can be algebraizable to a non-singular
projective variety having the ample canonical class.
\end{appthm}

In terms of formal scheme, this means that the relative canonical sheaf
of $\widehat{X}_{\Gamma}/\mathrm{Spf}R$ is ample.

\section{Examples.}

{\it Tate curve.}\ 
Let $T$ be the maxial torus consisting of all invertible diagonal
matrices. We take $n$ elements $q_i\in K^{\times}$ ($1\leq i\leq n$)
with $|q_i|<1$ and set $q_0=1$.
Define $\Gamma\subset T$ to be the set of all diagonal matrices with
$(i,i)$-entry an integer power of $q_i$ for $0\leq i\leq n$.
Then $\Gamma$ is a normal hyperbolic subgroup of $\mathrm{PGL}(n+1,K)$
with $\Delta_{\Gamma}$ the apartment corresponding to $T$.
If $n=1$, the curve $\Omega(\Delta_{\Gamma})/\Gamma$ is known to be 
the so-called {\it Tate curve}; in this case
$\Omega(\Delta_{\Gamma})$ is simply $\mathbb{P}^1_K$ minus two points.

\vspace{1ex}
{\it Mumford curve.}\ 
Consider again the $n=1$ case.
The rigid analytic curve $X_{\Gamma}$ with $\Gamma$ generated by at
least two elements is called a {\it Mumford curve}.
A Mumford curve is, in fact, algebraized to a non-singular projective
curve of genus greater than $1$.
Moreover, the formal model provided by the admissible covering inhereted 
from $\Omega(\Delta_{\Gamma})$ gives a stable $k$-split multiplicative
reduction of that curve.
The converse is also known:

\begin{appthm}[Mumford\ \cite{aacodcoclr}]
Let $X$ be a smooth and proper rigid analytic curve over $K$, and
suppose $X$ admits a formal model with $k$-split multiplicative
reduction. Then $X$ is isomorphic to $X_{\Gamma}$ for a unique 
$\Gamma$.
\end{appthm}

{\it Surface case.}\ 
In higher dimension it becomes extremely difficult to find normal
hyperbolic groups; but there are several known examples of torsion-free
uniform lattices.
Let $\Gamma$ is a torsion-free uniform lattice, and $N$ the number of
$\Gamma$-orbits in the vertex set $\Delta_0$.
Mumford \cite{aaswkak9pq0} calculated the several numerical invariants of the
quotient surface $X_{\Gamma}$:
\begin{enumerate}
\item $\chi(\OO_{X_{\Gamma}})=N(q-1)^2(q+1)/3$.
\item $(\matheur{c}_{1,X_{\Gamma}})^2=3\matheur{c}_{2,X_{\Gamma}}=
3N(q-1)^2(q+1)$.
\item $\matheur{q}(X_{\Gamma})=0$,
\end{enumerate}
where $\matheur{q}$ denotes the irregularity.
It is very interesting that the Mumford-Kurihara-Mustafin uniformization 
in two dimension always produces the surfaces which satisfies
the equality in Miyaoka-Yau inequality.
Mumford \cite{aaswkak9pq0} constructed one example of $\Gamma$ in $q=2$ such
that $N=1$; the resulting surface is one of the so-called ``fake
projective planes''.
Recently, Cartwright, Mantero, Steger, and Zappa constructed several
examples of lattices in $q=2$ and $q=3$. 
Making use of this groups, in \cite{tsrtfnu}, Ishida and the author
discussed other possible fake projective planes, and showed that there
are at least three fake projective planes (incl.\ Mumford's one) which
are not isomorphic to among others.
We will see in more detail the construction of Mumford's fake projective
plane in the next section.

\vspace{1ex}
{\it Shimura variety case.}
Let $D$ be a quaternion algebra over a totally real field $F$.
We assume that $D$ ramifies at every infinite place except one $\infty_0$.
Let $\Gamma_{\infty}$ be a subgroup of $D^{\times}/F^{\times}$ defined
by some congruence condition; we shall assume that it is ``sufficiently
small''.
Then it is classically known that $\Gamma_{\infty}$ acts discontinuously 
on the upper $\frac{1}{2}$ plane $\matheur{H}$ such that the quotient
$\matheur{H}/\Gamma_{\infty}$ is a {\it projective} curve (Shimura curve
of PEL-type).
The canonical model $\matheur{Sh}$ is defined over some ray-class field
of $F$.

There must be a finite place $\nu$ at which $B$ ramifies.
Let $p$ be the residue characteristic of $\nu$.
We assume that $\Gamma_{\infty}$ is maximal at $\nu$.
Then the Shimura curve $\matheur{Sh}$ has a bad reduction 
at $\nu$, and due to Cherednik \cite{uoacbdsopwcq} and Drinfeld
\cite{copsr}, the associated rigid analytic space over 
$F^{\mathrm{ur}}_{\nu}$ has a $p$-adic uniformization, where
$F^{\mathrm{ur}}_{\nu}$ is the maximal unramified extension of
$F_{\nu}$.
Let $\ovl{D}$ be the quaternion algebra which has local invariants
obtained precisely by switching those of $D$ at $\nu$ and $\infty_0$
Then there exists a subgroup $\Gamma_{\nu}$ in 
$\ovl{D}^{\times}/F^{\times}$ defined basically by the same congruence 
relation as that of $\Gamma_{\infty}$ such that
$\Omega\otimes F^{\mathrm{ur}}/\Gamma_{\nu}$ is isomorphic to the
assocaied analytic space of $\matheur{Sh}$ at $\nu$.
(This isomorphism actually decends to some finite unramified extension
of $F_{\nu}$, depending on the data of connected components of
$\matheur{Sh}$.)

This story has been generalized to higher dimension by Rapoport-Zink
\cite{psfpg} and Boutot-Zink \cite{tpuosc}.
The next section is devoted to explain one special example of the 
$p$-adic uniformization of a Shimura surface.

\section{Mumford's fake projective plane.}
Our example of $p$-adically uniformizable Shimura surface is constructed 
by means of some division algebra $D$. Th following construction is
inspired by the Mumford's oroginal paper \cite{aaswkak9pq0}.

Let $\zeta=\zeta_7$ be a primitive $7$th-root of unity, and
set $L\colon=\Q(\zeta)$; $L$ is a cyclic extension of $\Q$ of degree $6$
having an intermediate quadratic extension $K\colon=\Q(\lambda)\ (\cong
\Q(\sqrt{-7}))$, where $\lambda=\zeta+\zeta^2+\zeta^4$.
The Galois group of $L/K$ is generated by the Frobenius map $\sigma\colon
\zeta\mapsto\zeta^2$, and that of $K/\Q$ is generated by the complex 
conjugation $z\mapsto\ovl{z}$.
Note that the prime $2$ decomposes on $K$, more explicitly, 
$2=\lambda\ovl{\lambda}$ gives the prime factorization, and the prime
$7$ ramifies.
We fix an infinite place $\varepsilon\colon K\hookrightarrow\C$ by
$\lambda\mapsto\frac{-1+\sqrt{-7}}{2}$.

We set $\mu\colon=\lambda/\ovl{\lambda}$, and define the central
division algebra $D$ over $K$ of dimension $9$ by
$$
D\colon=\bigoplus^2_{i=0}L\Pi^i
$$
with $\Pi^3=\mu$ and $\Pi z=z^{\sigma}\Pi$ for $z\in L$.
The involution $\ast$ on $D$ is defined by $z^{\ast}=\ovl{z}$ for $z\in L$
and $\Pi^{\ast}=\ovl{\mu}\Pi^2$.
It is easily verified that $\ast$ is positive.

We introduce a non-degenerate anti-symmetric $\Q$-bilinear
form $\psi=\psi_b$ on $V=D$ by putting 
$b=(\lambda-\ovl{\lambda})-\ovl{\lambda}\Pi+\ovl{\lambda}\Pi^2$.
Then define an algebraic group $G$ over $\Q$ as follows: 
Given a $\Q$-algebra $R$, we set
$$
G(R)\colon=\{g\in((D\otimes R)^{\mathrm{opp}})^{\times}\ |\ \psi(xg,yg)=
c(g)\psi(x,y),\ c(g)\in R^{\times},\ x,y\in V\otimes R\}.
$$
The algebraic group $G_{\R}$ is isomorphic to $\mathrm{GU}(2,1)$, which 
acts on the two dimensional complex unit-ball $\matheur{B}$.
In the following we will construct an open compact subgroup $C$ in
$G(\mathbb{A}_{\mathrm{f}})$.
By these data, we can define a Shimura variety
$$
\mathcal{S}_C\colon=G(\Q)\setminus\matheur{B}\times 
G(\mathbb{A}_{\mathrm{f}})/C.
$$
The canonical model of $\mathcal{S}_C$ is defined over $E\colon=
\varepsilon(K)$, which we denote by $\matheur{Sh}_C$.

\vspace{1ex}
{\it The congruence condition.}\ 
First we define an algebraic group $I$ over $\Q$ by 
$$
I(R)\colon=\{g\in(\mathrm{M}_d(K)\otimes R)^{\times}\ |\
H^{-1}\,^t\ovl{g}H\in R^{\times}\}
$$
for any $\Q$-algebra $R$, where 
$$
H=\left[
\begin{array}{lll}
3&\ovl{\lambda}&\ovl{\lambda}\\
\lambda&3&\ovl{\lambda}\\
\lambda&\lambda&3
\end{array}
\right].
$$
Then $I$ is an inner form of $G$ so that 
$G(\mathbb{A}^2_{\mathrm{f}})\cong I(\mathbb{A}^2_{\mathrm{f}})$.

To define the level structure, it suffices to define the 
prime-$2$-part $C^2$.
We define
$C^2$ in the form $C^2=C^{2,7}C_7$ with $C_7\in G(\Q_7)\cong I(\Q_7)$
and $C^{2,7}$ being the maximal one:
$$
C^{2,7}\colon=\{g\in G(\mathbb{A}^{2,7}_{\mathrm{f}})\ |\ \Gamma^{2,7}g
\subseteq\Gamma^{2,7}\},
$$
with $\Gamma^{2,7}=\Gamma\otimes_{\Z}\widehat{\Z}^{2,7}$, where
$\widehat{\Z}^{2,7}=\lim\Z/m\Z$ with the projective limit taken over $m$ 
such that $(2,m)=(7,m)=1$.

Now the component $C_7$ is defined as follows:
we have an isomorphism
$$
G(\Q_7)\stackrel{\sim}{\longrightarrow}\{g\in\mathrm{GL}_3(\til{\Q}_7)
\ |\ g^{\ast}Hg=c(g)H,\ c(g)\in\Q^{\times}_7\},
$$
where $\til{\Q}_7$ is the ramified quadratic extension of $\Q_7$, and
$\ast$ is the matrix transposition followed by the Galois action.
Let $\til{C}_7$ be the maximal open compact subgroup of $G(\Q_7)$ 
consisiting of $g$ which maps $\Gamma\otimes_{\Z}\Z_7$ to itself.
We consider modulo $\sqrt{-7}$ reduction of $\til{C}_7$, which is a
subgroup $G_0$ in $\mathrm{GL}_3(\mathbb{F}_7)$.
The matrix $H$ mod $\sqrt{-7}$ is of rank $1$, and has $2$ dimensional
null space $N_0$.
Restricting elements in $G_0$ to $N_0$, we obtain a homomorphism
$$
\pi\colon\til{C}_7\longrightarrow\mathrm{GL}_2(\mathbb{F}_7).
$$
The last group has $2^63^27$ elements. Let $S$ be a $2$-Sylow subgroup
of $\mathrm{GL}_2(\mathbb{F}_7)\cap\{g\ |\ \mathrm{det}\,g=\pm 1\}$.
Now we set $C_7\colon=\pi^{-1}(S)$.

\begin{appthm}\label{thm-uniformization}

The canonical model $\matheur{Sh}_C$ of the Shimura variety obtained by
the data as above has exactly three geometrically connected components defined
over $\Q(\zeta_7)$. Moreover, for any connected component $S$, its
base change $S\otimes_{\Q(\zeta_7)}\Q_2(\zeta_7)$ is isomorphic to the
Mumford's fake projective plane \cite{aaswkak9pq0} tensored by $\Q_2(\zeta_7)$.
\end{appthm}

More precisely, the associated rigid analytic space to 
$\matheur{Sh}_C\otimes_EE_{\lambda}$ is isomorphic to
$\Omega/\Gamma$ with
$$
\Gamma=I(\Q)\bigcap(J(\Q_2)\times C^2)\subset J(\Q_2)\cong
\mathrm{PGL}(3,\Q_2),
$$
where $I(\Q)$ is considered as a subgroup of 
$I(\mathbb{A}_{\mathrm{f}})\cong J(\Q_2)\times 
G(\mathbb{A}^2_{\mathrm{f}})$.
The proof of the theorem is in \cite{opuofpp}, where two more 
fake projective planes are
discussed as examples of Shimura varieties.

\chapter{$p$-adic symmetric domains and Totaro's theorem}
\renewcommand{\thechapter}{\Alph{chapter}}
\addtocontents{toc}{\protect\par\vskip2mm\hskip80mm by N. Tsuzuki\par\vskip5mm}

\vspace{1ex}
\begin{center}
By {\sc TSUZUKI Nobuo}
\end{center}
\vspace{10mm}

This appendix is a short exposition of 
M. Rapoport and T. Zink's construction of $p$-adic symmetric 
domains \cite{psfpg} and of B. Totaro's theorem \cite{tpipht}. 
Let $G$ be a connected reductive algebraic group over $\Q_p$. 
The set ${\mathcal F}$ of filtrations on an $F$-isocrystal with $G$-structure 
has a structure of a homogeneous space. 
Rapoport and Zink introduced a $p$-adic rigid analytic structure on 
the set ${\mathcal F}^{wa}$ of weakly admissible points in ${\mathcal F}$. 
They conjectured that the point in ${\mathcal F}^{wa}$ 
is characterized by the semistability 
in the sense of the geometric invariant theory \cite{git} 
and Totaro proved this conjecture. 

\section{Weakly admissible filtered isocrystals.}

We recall J.-M. Fontaine's definition of weakly admissible 
filtered $F$-isocrystals \cite{mgmfeadb}. 

\begin{apppara}
 Let $p$ be a prime number, $k$ a perfect field of 
 characteristic $p$, 
 $K_0$ an absolutely unramified discrete valuation field 
 of mixed characteristics $(0, p)$ with residue field $k$, 
 $\overline{K}_0$ an algebraic closure of $K_0$, 
 and $\sigma$ the Frobenius automorphism on $K_0$. 
\end{apppara}

\begin{appdfn}
 (1) An $F$-isocrystal over $k$, 
 (we simply say ``isocrystal''), is 
 a finite dimensional $K_0$-vector space $V$ 
 with a bijective $\sigma$-linear endomorphism $\Phi : V \rightarrow V$. 
 We denote the category of isocrystals over $k$ by $\mathrm{Isoc}(K_0)$. 

 (2) For a totally ramified finite extension $K$ of $K_0$ in $\overline{K}_0$, 
 a filtered isocrystal $(V, \Phi, F^{\bp})$ over $K$ 
 is an isocrystal $(V, \Phi)$ with a decreasing filtration $F^{\bp}$ 
 on the $K$-vector space $V \otimes_{K_0} K$ 
 such that $F^r = V \otimes_{K_0} K$ 
 for $r \ll 0$ and $F^s = 0$ for $s \gg 0$. 
 We denote the category of filtered isocrystals over $K$ by $MF(K)$. 
\end{appdfn}

Fontaine also introduced a filtered isocrystal with nilpotent operator $N$ \cite{fontaine94:_repres}. 
In this appendix we restrict our attension to filtered isocrystals with $N = 0$. 

The category $MF(K)$ is a $\Q_p$-linear 
additive category with $\otimes$ and internal $Hom${}'s, but not abelian. 
A subobject $(V', \Phi', {F'}^{\bp})$ of a filtered isocrystal 
$(V, \Phi, F^{\bp})$ is a $\Phi$-stable $K_0$-subspace $V'$ 
such that $\Phi' = \Phi|_{V'}$ and ${F'}^i = (V' \otimes_{K_0} K) \cap F^i$ 
for all $i$. 

\begin{appdfn}
 \label{weakad}
 Let $K$ be a totally ramified finite extension of $K_0$ 
 in $\overline{K}_0$. 
 A filtered isocrystal $(V, \Phi, F^{\bp})$ over $K$ 
 is weakly admissible if, for any subobject 
 $(V', \Phi', {F'}^{\bp}) \ne 0$, we have 
 $$
      \displaystyle{\mathop{\sum}_i}\ i \dim_{F'} \mathrm{gr}^i_{F'}(V'
 \otimes_{K_0} K) 
                            \leq \mathrm{ord}_p(\det(\Phi'))
 $$
 and the equality holds for $(V', \Phi', {F'}^{\bp}) = (V, \Phi, F^{\bp})$. 
 Here $\mathrm{ord}_p$ is an additive valuation of $K\sb 0$ normalized by 
 $\mathrm{ord}_p(p) = 1$. 
\end{appdfn}

The category of weakly admissible filtered isocrystals 
is an abelian category 
which is closed under duals in the category of filtered 
isocrystals. 
Fontaine proved that an admissible filtered isocrystal over $K$, that means 
a filtered isocrystal arising from a crystalline representation 
of the absolute Galois group of $K$ via Fontaine's functor, is weakly admissible 
and conjectured that a weakly admissible filtered isocrystal is admissible 
in \cite{mgmfeadb}. 
The category of admissible filtered isocrystals is 
a $\Q_p$-linear abelian category with $\otimes$ and duals. 
Hence, he also conjectured that the category of weakly admissible filtered 
isocrystals is closed under $\otimes$, 
and this was proved by G. Faltings in \cite{midag}. (See also \cite{tpipht}.)

In \cite{colmez00:_const} P. Colmez and Fontaine proved a weakly admissible filtered 
isocrystal is admissible.

\section{Filtered isocrystals with $G$-structure.}

\begin{apppara}
 Let $G$ be a linear algebraic group over $\Q_p$ and
 denote by
 $\mathrm{Rep}_{\Q_p}(G)$ the category of
 finite dimensional $\Q_p$-rational representations of $G$.
 An exact faithful $\otimes$-functor
 $\mathrm{Rep}_{\Q_p}(G) \rightarrow \mathrm{Isoc}(K_0)$
 is called an isocrystal with $G$-structure over $K_0$.

 Let $b \in G(K_0)$. Then, the functor
 $$
   \mathrm{Rep}_{\Q_p}(G)
                                        \rightarrow \mathrm{Isoc}(K_0)
 $$
 associated to $b$,
 defined by $V \mapsto (V \otimes K_0, b(\mathrm{id} \otimes \sigma))$,
 is an isocrystal with $G$-structure over $K_0$. \cite{iwas}
 Two elements $b$ and $b'$ in $G(K_0)$ are conjugate
 if and only if there is an element $g \in G(K_0)$ such that
 $gb\sigma(g)^{-1} = b'$.
 In this case, $g$ defines an isomorphism between the isocrystals with
 $G$-structure associated to $b$ and $b'$.

 If $G$ is connected and $k$ is algebraically closed,
 then any  isocrystal with $G$-structure over $K_0$ is associated to
 an element $b \in G(K_0)$ as above. \cite{otcasofwas}
\end{apppara}

\begin{apppara}
 \label{sfil}
 Let $\mathbb{D} = \displaystyle{\mathop{\lim}_\leftarrow}\ \G_m$
 be the pro-algebraic group over $\Q$
 whose character group is $\Q_p$. For an element $b \in G(K_0)$,
 R.E. Kottwitz defined a morphism
 $$
      \nu : \mathbb{D} \rightarrow G_{K_0}
 $$
 of algebraic groups over $K_0$ which is characterized by the property
 that, for any object $V$ in $\mathrm{Rep}_{\Q_p}(G)$,
 the $\Q$-grading of $V \otimes K_0$ associated to $- \nu$ is the
 slope grading of the isocrystal
 $(V \otimes K_0, b(\mathrm{id} \otimes \sigma))$.
 (The sign of our $\nu$ is different from the one in \cite{iwas}.)
 For a suitable positive integer $s$, $s\nu$ is regarded
 as a one-parameter subgroup of $G$ over $K_0$.
\end{apppara}

\begin{appdfn}
 \label{decent} A $\sigma$-conjugacy class $\bar{b}$
 of $G(K_0)$ is decent if there is an element $b \in \bar{b}$ such that
 $$
          (b\sigma)^s = s\nu(p)\sigma^s
 $$
 for some positive integer $s$.
\end{appdfn}

One knows that, for a decent $\sigma$-conjugacy class $\bar{b}$, $b$ and $\nu$
as above are defined over $\Q_{p^s}$.
If $G$ is connected and $k$ is algebraically closed, then any
$\sigma$-conjugacy
class is decent. \cite{iwas}

\begin{apppara}
 Let $K$ be a totally ramified finite extension
 of $K_0$ in $\overline{K}_0$.
 For a one-parameter subgroup $\lambda : \G_m \rightarrow G$
 over $K$ and an element $b \in G(K_0)$,
 we have an exact $\otimes$-functor
 $$
    {\mathcal I} : \mathrm{Rep}_{\Q_p}(G)
               \rightarrow MF(K)
 $$
 which is defined by $V \mapsto
 (V \otimes K_0, b(\mathrm{id} \otimes \sigma), F^{\bp}_\lambda)$.
 Here $V_{K, \lambda, j}$ is the subspace of $V \otimes K$
 of weight $j$ with respect to $\lambda$ and
 $$
    F^i_\lambda = \displaystyle{\mathop{\oplus}_{j \geq i}}\ V_{K, \lambda, j}
 $$
 is the weight filtration associated to $\lambda$.
\end{apppara}

\begin{appdfn}
 A pair $(\lambda, b)$ as above is weakly admissible
 if and only if
 the filtered isocrystal ${\mathcal I}(V)$ over $K$ is so for any
 object $V$ in $\mathrm{Rep}_{\Q_p}(G)$.
\end{appdfn}
To see the weak admissibility for $(\lambda, b)$, it is enough to check
the weak admissibility of ${\mathcal I}(V)$ for a faithful
representation $V$ of $G$.  Indeed, any representation of $G$ appears as
a direct summand of $V^{\otimes m} \otimes (V^\vee)^{\otimes n}$ and
${\mathcal I}(V)^{\otimes m} \otimes ({\mathcal I}(V)^\vee)^{\otimes n}$
is weakly admissible by Faltings (see \ref{weakad}).  Here $V^\vee$
(resp. ${\mathcal I}(V)^\vee$) is the dual of $V$ (resp.  ${\mathcal
I}(V)$).

\section{Totaro's theorem.}

In this section we assume that $k$ is algebraically closed.

\begin{apppara}
 Let $G$ be a reductive algebraic group over $\Q_p$.
 We fix a conjugacy class of a one-parameter subgroup
 $\lambda : \G_m \rightarrow G$ over $\overline{K}_0$.
 Here two one-parameter subgroups $\lambda, \lambda'$ are conjugate
 if and only if $g \lambda g^{-1} = \lambda'$ for some element $g \in
 G(\overline{K}_0)$.
 Then, there is a finite extension $E$ of $\Q_p$ in $\overline{K}_0$
 such that the conjugacy class of $\lambda$ is defined over $E$.
 Let us suppose that $\lambda$ is defined over $E$ and denote by
 $\breve{E}$ the composite field $EK_0$ in $\overline{K}_0$
\end{apppara}

\begin{apppara}
 Two one-parameter subgroups of $G$
 over $\overline{K}_0$ are equivalent
 if and only if they define the same weight filtration for any object
 in $\mathrm{Rep}_{\Q_p}(G)$.
 Note that, if two one-parameter subgroups are equivalent,
 then they belong to the same conjugacy class.

 Consider the functor
 $$
     R \mapsto
 \{\text{the equivalence classes in the conjugacy class of $\lambda$
 defined over $R$} \}
 $$
 on the category of $E$-algebras. If one defines an algebraic
 subgroup of $G$ over $E$ by
 $$
      P(\lambda)(\overline{K}_0) = \{ g \in G(\overline{K}_0) ~|~
             g \lambda g^{-1}~ \text{is equivalent to}~ \lambda \},
 $$
 then $P(\lambda)$ is parabolic and the functor above is represented by
 the projective variety $G_E/P(\lambda)$.
 We denote this homogeneous space over $E$ by ${\mathcal F}_\lambda$.
 If $V$ is a faithful representation
 in $\mathrm{Rep}_{\Q_p}(G)$
 and if we denote by $\mathrm{Flag}_\lambda(V)$ the flag variety over
 $\Q_p$
 which represents the functor
 $$
     R \mapsto \left\{
            \begin{array}{l}
             \text{the filtrations $F^{\bp}$ of $V \otimes R$
             as $R$-modules such that} \\
             \text{$F^i$ is a direct summand and
             $\mathrm{rank}_R F^i
                = \dim_{\overline{K}_0} F_\lambda^i(V \otimes \overline{K}_0)$}
            \end{array}
         \right\}
 $$
 on the category of $\Q_p$-algebras, then there is
 a natural $E$-closed immersion
 $$
      {\mathcal F}_\lambda
          \rightarrow \mathrm{Flag}_\lambda(V)
 \otimes_{\Q_p} E.
 $$
\end{apppara}

\begin{apppara}
 Let $b \in G(K_0)$.
 For a finite extension $K$ of $\breve{E}$, a point $\xi$
 in ${\mathcal F}_\lambda(K)$ is called weakly admissible if and only if
 the pair $(\xi, b)$ is weakly admissible. This condition is independent
 of the choice of the representative in the equivalence class $\xi$.
 We denote by ${\mathcal F}^{wa}_{\lambda, b}(K)$ the subset of weakly
 admissible points.
 Totaro gave a characterization of ${\mathcal F}^{wa}_{\lambda, b}$
 in the sense of geometric invariant theory. \cite{tpipht}
 We explain his theory in the rest of this section.
\end{apppara}

\begin{apppara}
 For a maximal torus $T$ in $G_{\overline{K}_0}$,
 let $X^*(T)$ be the free abelian group of characters,
 $X_*(T)$ the free abelian group of one-parameter subgroups,
 and $<\ , \ > : X^*(T) \times X_*(T) \rightarrow \Z$
 the perfect pairing with $\chi(\xi(t)) = t^{<\chi, \xi>}$.
 If $N(T)$ is the normalizer of $T$ in $G$, the Weyl group $W(T) = N(T)/T$
 acts $X_*(T)$ via inner automorphisms.

 Now we fix an invariant norm $||\ ||$ on $G$, a non-negative
 real valued function on the set of one-parameter subgroups
 of $G_{\overline{K}_0}$, such that

 a) $||g \xi g^{-1}|| = ||\xi||$ for any $g \in G(\overline{K}_0)$,

 b) for any maximal torus $T$, there is a positive definite rational valued
 bilinear form $(\ , \ )$ on $X_*(T) \otimes \Q$
 with $(\xi, \xi) = ||\xi||^2$,

 c) $||\gamma(\xi)|| = ||\xi||$
 for $\gamma \in \mathrm{Gal}(\overline{K}_0/K_0)$,
 where $\gamma(\xi)(t) = \gamma(\xi(t))$.

 \noindent
 The bilinear form on $X_*(T) \otimes \Q$
 as above is invariant under the action of
 the Weyl group by (a). For any maximal torus $T$,
 invariant norms are in one-to-one correspondence
 with $(\mathrm{Gal}(\overline{K}_0/K_0), W(T))$-invariant positive definite
 rational valued bilinear forms on $X_*(T) \otimes \Q$ since
 all maximal tori are conjugate and, if $g \xi g^{-1} \in X_*(T)$ for
 $\xi \in X_*(T)$ and $g \in G(\overline{K}_0)$, then there is $h \in W(T)$
 with $g \xi g^{-1} = h \xi h^{-1}$ by \cite{git}.
 (See also \cite{iiit} and \cite{tpipht}.) Hence, such an invariant norm exists.
\end{apppara}

\begin{apppara}
 Now we assume that $G$ is connected. Let $U(\lambda)$ be the unipotent
 radical of $P(\lambda)$, whose elements act on the graded space
 $\mathrm{gr}\ F_{\lambda}^{\bp}$ trivially. Then there is a bijection between
 the set of maximal tori of $G$ in $P(\lambda)$
 and the set of maximal tori
 of $P(\lambda)/U(\lambda)$ by the natural projection $T \mapsto \overline{T}$.
 Hence, the invariant norm on $G$ induces the one on $P(\lambda)/U(\lambda)$.
 Fix a maximal torus $T$ of $G$.
 Since the image of $\lambda$ is contained in the center of
 $P(\lambda)/U(\lambda)$,
 the perfect pairing associated to the invariant norm determines
 the dual of $\lambda$ in $X^*(\overline{T})$.
 This dual can extend to a character $\otimes \Q$ of
 $P(\lambda)/U(\lambda)$.
 Now we define a $G$-line bundle $\otimes \Q$, $L_\lambda$,
 on ${\mathcal F}_\lambda$
 by the associated one to the negative of the dual character $\otimes
 \Q$
 of $\lambda$. By construction, the line bundle ${\otimes \Q}$,
 $L_\lambda$, depends only on the conjugacy class of $\lambda$
 and is ample.

 Let $J$ be a smooth affine group scheme over $\Q_p$ such that
 $$
       J(\Q_p) = \{ g \in G(K_0) ~|~ g (b \sigma) = (b \sigma) g \}
 $$
 (which is introduced in \cite{psfpg}). 
 Since $J_{K_0} \subset G_{K_0}$, the pull back $L_{\lambda \breve{E}}$ 
 of $L_\lambda$ on ${\mathcal F}_{\lambda \breve{E}}$ 
 is an ample $J_{\breve{E}}$-line bundle. 

 By the same construction as above, $\nu$ in \ref{sfil} gives 
 a character $\otimes \Q$ of $P(\nu)$. 
 The opposite of this character $\otimes \Q$ 
 determines a $J_{\breve{E}}$-action $\otimes \Q$ 
 on the trivial line bundle on ${\mathcal F}_{\lambda \breve{E}}$ 
 since $J_{K_0} \subset P(\nu)$. We denote it by $L_\nu^0$. 

 We put a $J_{\breve{E}}$-line bundle ${\otimes \Q}$,
 $L = L_{\lambda \breve{E}} \otimes L_\nu^0$, 
on ${\mathcal F}_{\lambda \breve{E}}$. 
 Then it is ample and depends only on 
 $b$ and the conjugacy class of $\lambda$. 
 We denote by ${\mathcal F}_\lambda^{ss}(L)$ the set of semistable points 
 in ${\mathcal F}_\lambda$ 
 with respect to $L$ in the sense of D. Mumford \cite{git}. 
\end{apppara}

\begin{appthm}
 \cite{tpipht}
 {\it Suppose that $G$ is connected and reductive. 
 For any finite extension $K$ of $\breve{E}$, we have 
 $$
      {\mathcal F}_{\lambda, b}^{wa}(K) = {\mathcal F}_{\lambda 
 \breve{E}}^{ss}(L)(K). 
 $$
 }
\end{appthm}

We shall sketch Totaro's proof. 
First, let $G = GL(n)$ and let us consider the invariant norm 
induced by the pairing 
$$
    (\alpha, \beta) = \displaystyle{\mathop{\sum}_{i, j}}\ 
 ij \dim_{\overline{K}_0} 
\mathrm{gr}_{F_\alpha^{\bp}}^i\mathrm{gr}_{F_\beta^{\bp}}^j 
             (K^n) 
$$
for one-parameter subgroups $\otimes \Q$, $\alpha, \beta$ of 
$GL(n)$ over $K$. 
If one puts $\mu_\alpha(V) 
= (\displaystyle{\mathop{\sum}_i}\ i\dim_K 
\mathrm{gr}_{F_\alpha^{\bp}}^i(V \otimes K)) /\dim\ V$, 
then one has 
$$
     (\alpha, \beta) 
              = \int\ (\mu_\alpha(F_\beta^j) - \mu_\alpha(V))\dim_K 
F_\beta^j dj 
                      + \mu_\alpha(V)\mu_\beta(V)\dim\ V.
$$
So, $(\xi, b)$ is weakly admissible 
if and only if $(\xi, \alpha) + (\nu, \alpha) \leq 0$ 
for any one-parameter subgroup $\alpha$ of $GL(n)$ over $K_0$ 
with the filtration $F_\alpha^{\bp}$ as subisocrystals. 
In other words, $(\xi, b)$ is weakly admissible 
if and only if $(\xi, \alpha) + (\nu, \alpha) \leq 0$ 
for any one-parameter subgroup $\alpha$ of $J$ over $\Q_p$. 
Hence, the assertion follows from the calculation of Mumford's numerical 
invariant below. 

\begin{applem}
 {\it If $\mu(\xi, \alpha, L)$ is Mumford's numerical 
 invariant of $\xi \in {\mathcal F}_\lambda$ 
 for a one-parameter subgroup $\alpha$ of $J$ over $\Q_p$, then 
 $$
    \mu(\xi, \alpha, L) = - (\xi, \alpha) - (\nu, \alpha).
 $$
 }
\end{applem}

Next, let $G$ be arbitrary, $V$ a faithfully representation of $G$, and 
consider the invariant norm on $G$ induced from the above norm 
by the natural immersion $G \rightarrow GL(V)$. 
The Mumford's numerical invariant of weakly admissible points 
is non-negative for any one-parameter subgroup of $J(GL(V))$. 
Hence it is so for any one-parameter subgroup of $J$, 
and the weak admissibility implies the semistability. 
To see the converse, one needs to show that, if $\xi \in {\mathcal 
F}_\lambda^{ss}(L)(K)$, 
$(\xi, \alpha) + (\nu, \alpha) \leq 0$ 
for any one-parameter subgroup $\alpha$ of $J(GL(V))$ over $\Q_p$. 
If $\alpha$ is semistable for the $G_K$-line bundle $L_\alpha$ 
on ${\mathcal F}_\alpha$, the assertion follows from the first part. 
In the case where $\alpha$ is not semistable, one can use 
Kempf's filtration \cite{iiit} and Ramanan and Ramanathan's work 
\cite{srotif}, and obtains the required inequality. 

Finally one needs to prove the independence of the choice of the norm.
Suppose that the identity is valid for the particular norm.
Since $G$ is a quotient of a product of a torus and some simple algebraic
groups
by a finite central subgroup \cite{gr},
one can reduce the assertion in the case of tori and simple groups.
In the case of tori it was proved in \cite{psfpg},
and in the case of simple groups it is true
since the norm of the simple group
comes from the Killing form up to a positive rational multiple.
\hspace{\fill} $\square$

\section{$p$-adic symmetric domains.}

Let $k$ be the algebraic closure of the prime field $\F_p$,
$\C_p$ the $p$-adic completion of a fixed algebraic closure
$\overline\Q_p$ of $\Q_p$, and
$K_0 = \widehat\Q_p^{ur}$ the $p$-adic completion
of the unramified extension in $\C_p$.

\begin{apppara}
 Let $G$ be a reductive group over $\Q_p$, $b \in
 G(K_0)$,
 and fix a conjugacy class $\{ \lambda \}$ of a one-parameter subgroup $\lambda$
 of $G$ over $\overline\Q_p$.

 Rapoport and Zink gave a rigid analytic structure
 on ${\mathcal F}_{\lambda, b}^{wa}$ as an admissible open subset
 in ${\mathcal F}_{\lambda \breve{E}}$ and call it the $p$-adic
 symmetric domain 
 associated to the triple $(G, \{ \lambda \}, b)$ in \cite{psfpg}.
 This notion of $p$-adic
 symmetric domains is different from that of M. van der Put
 and H. Voskuil in \cite{ssatsagoalf}. Indeed, for any discrete co-compact subgroup
 $\Gamma$ of $G(\breve{E})$,
 the quotient ${\mathcal F}_{\lambda, b}^{wa}/\Gamma$
 is not always a proper analytic space over $\breve{E}$.
\end{apppara}

\begin{appthm}
 \cite{psfpg} {\it The set ${\mathcal F}_{\lambda,
 b}^{wa}$
 of weakly admissible points with respect to $b$ in ${\mathcal
 F}_\lambda(\C_p)$
 is an admissible open subset of ${\mathcal F}_{\lambda \breve{E}}$ as a
 rigid analytic
 space.}
\end{appthm}

Now we sketch the proof of the theorem in \cite{psfpg}.
By \cite{iwas} one may assume that the $\sigma$-conjugacy class of $b$
is decent with the decent equation $(b\sigma)^s = s\nu(p)\sigma^s$ as in
\ref{decent}.
Let $V$ be a faithful representation in
$\mathrm{Rep}_{\Q_p}(G)$, 
$V_s = V \otimes \Q_{p^s}$, and $\Phi_s = b(\mathrm{id} \otimes
\sigma)$.
Then $(V \otimes K_0, b(\mathrm{id} \otimes \sigma))
= (V_s, \Phi_s) \otimes_{\Q_{p^s}} K_0$.
Put $V_s = \displaystyle{\mathop{\oplus}_\lambda}\ V_{s, j}$
to be the isotypical decomposition for $\Phi_s$. The functor
$$
    R \mapsto \{ V' \subset V_s \otimes_{\Q_{p^s}}
R ~|~
                 V'~ \text{is a direct summand with}~
                 V' = \displaystyle{\mathop{\oplus}_j}\
                       V' \cap (V_{s, j}
                       \otimes_{\Q_{p^s}} R) \}
$$
on the category of $\Q_{p^s}$-algebras is represented
by a disjoint sum $T'$ of closed
subschemes of Grassmannians of $V_s$. $T'$ descends to a
$\Q_p$-variety
$T$ and one has
$$
     T(\Q_p) = \{ \Phi_s\text{-stable subspaces of}~ V_s \}.
$$
Indeed, $\Phi_s$ gives a descent datum $\alpha : T' \rightarrow {T'}^\sigma$,
where ${T'}^\sigma(R)$ is a set of direct summands of
$V_{s, \lambda} \otimes_{\Q_{p^s}, \sigma} R$ with
the isotypical decomposition,
and $\alpha^{s-1} \circ \cdots \circ \alpha : T' \rightarrow T'$
is the identity by the decent equation.

Consider the closed subscheme over $\Q_{p^s}$
$$
     {\mathcal H} \subset
           (\mathrm{Flag}_\lambda(V) \times T)
                    \otimes_{\Q_p} \Q_{p^s}
$$
which consists of pairs $(F^{\bp}, V')$ such that
$$
     \displaystyle{\mathop{\sum}_i}\ i\,\mathrm{rank}\
               \mathrm{gr}_{F^{\bp} \cap V'}^i\ (V')
                > \mathrm{ord}_p(\det(\Phi_s|_{V'_j})).
$$
Then, by the definition of weak admissibility, one has
$$
    \mathcal{F}_{\lambda, b}^{wa}(\C_p)
       = \mathcal{F}_\lambda(\C_p) \cap
\biggl(\mathrm{Flag}_\lambda(V)(\C_p)
        - \bigcup_{t\in T(\Q_p)}{\mathcal H}_t\biggr),
$$
where ${\mathcal F}_\lambda(\C_p)$ is identified with the image of
the immersion
${\mathcal F}_\lambda(\C_p) \subset \mathrm{Flag}_{\lambda}(\C_p)$.

Fix embeddings of $\mathrm{Flag}_\lambda(V)$ and $T$
in projective spaces over $\Q_p$
and a finite set $\{ f_j \}$ of bi-homogeneous polynomials
of definition of $\mathrm{Flag}_\lambda(V) \times T$
with integral coefficients.
For $\epsilon > 0$, consider a tubular neighbourhood
$$
    {\mathcal H}_t(\epsilon)
        = \{ x \in \mathrm{Flag}_\lambda(V)(\C_p)
                      ~|~ |f_j(x, t)| < \epsilon~ \text{for all}~ j \}
$$
of ${\mathcal H}_t$.
Here we choose unimodular representatives for $x$ and $t${}'s
and $|\ |$ is an absolute value on $\C_p$.
Then there is a finite set $S \subset T(\Q_p)$ such that
$ \bigcup_{t \in T(\Q_p)}{\mathcal H}_t(\epsilon)
= \bigcup_{t \in S}{\mathcal H}_t(\epsilon)$
for the local compactness.
${\mathcal F}_\epsilon = {\mathcal F}_\lambda(\C_p)
\cap (\mathrm{Flag}_\lambda(V)(\C_p) -
\bigcup_{t \in T(\Q_p)}\mathcal{H}_t(\epsilon))$
is an admissible open subset of ${\mathcal F}_{\lambda, b \breve{E}}^{wa}$,
hence
$$
    {\mathcal F}_1 \subset {\mathcal F}_{\frac{1}{2}}
          \subset {\mathcal F}_{\frac{1}{3}} \subset \cdots
$$
is an admissible covering of ${\mathcal F}_{\lambda, b \breve{E}}^{wa}$.
\hspace{\fill} $\square$

\nocite{ip}
\nocite{arotts}
\nocite{taotrc}
\nocite{aronef}
\nocite{gdefdlch}
\nocite{cdgplcdd1}
\nocite{craswsrts}
\nocite{farg}
\nocite{gardtk}
\nocite{grecdvadcp}
\nocite{iwas}
\nocite{tautbidnf}
\nocite{gray86:_linear_rieman_poinc}
\nocite{michon81:_courb_shimur}
\nocite{ogg83:_real_shimur}
\nocite{serre68:_corps}
\nocite{put00:_discon_pgl_k}

\newcommand{\noop}[1]{}\newcommand{\mathcalbigd}{{\mathcal
  D}}\newcommand{\bfbigq}{{\bf Q}}\newcommand{\bfbigz}{{\bf Z}}
\providecommand{\bysame}{\leavevmode\hbox to3em{\hrulefill}\thinspace}

\end{document}